\documentclass[10pt]{amsart}
\usepackage[english,russian]{babel}
\usepackage[cp1251]{inputenc}
\usepackage{amsmath}
\usepackage{amssymb}
\usepackage{amsfonts}

\usepackage[linktocpage=true, colorlinks=true, linkcolor=blue, citecolor=blue, urlcolor=blue]{hyperref}

\begin{document}
\renewcommand{\refname}{References}
\renewcommand\contentsname{Contents}

\thispagestyle{empty}

\title[The Hypotheses on Expansion of Iterated Stratonovich Stochastic
Integrals]
{The Hypotheses on Expansion of Iterated Stratonovich Stochastic
Integrals of Arbitrary Multiplicity and Their Partial Proof}
\author[D.F. Kuznetsov]{Dmitriy F. Kuznetsov}
\address{Dmitriy Feliksovich Kuznetsov
\newline\hphantom{iii} Peter the Great Saint-Petersburg Polytechnic University,
\newline\hphantom{iii} Polytechnicheskaya ul., 29,
\newline\hphantom{iii} 195251, Saint-Petersburg, Russia}
\email{sde\_kuznetsov@inbox.ru}

\linespread{1.25}

\maketitle {\small
\begin{quote}

\noindent{\sc Abstract. } 
In this article, we collected more than thirty theorems on
expansions of iterated Ito and Stratonovich 
stochastic integrals which have been formulated and proved by the author
in the period from 1997 to 2026.
These theorems open up a new direction for study of iterated Ito and Stratonovich stochastic integrals.
We consider two
theorems on expansion of iterated Ito stochastic integrals of arbitrary multiplicity 
$k$ $(k\in\mathbb{N})$
based on generalized multiple Fourier series converging in the sense of norm in Hilbert
space $L_2([t, T]^k)$.
We adapt these theorems on expansion of iterated Ito stochastic integrals
of arbitrary multiplicity $k$ $(k\in\mathbb{N})$
for iterated Stratonovich stochastic integrals
of multiplicities 1 to 8 (the case of continuously differentiable 
weight functions and a complete orthonormal system of Legendre polynomials or 
trigonometric functions in $L_2([t, T])$) and for iterated Stratonovich stochastic
integrals of multiplicities 1 to 6 (the case of 
an arbitrary complete orthonormal system of functions in $L_2([t, T])$).
On the base of the presented theorems we formulate several hypotheses
on expansions of iterated Stratonovich stochastic integrals of arbitrary 
multiplicity $k$ $(k\in\mathbb{N})$.
Recently, Hypothesis~8 (Sect.~27) has been proved for the case
of an arbitrary complete orthonormal system of functions in $L_2([t, T])$
and $p_1=\ldots=p_k=p,$ $k\in\mathbb{N}$ 
but under one additional condition (Theorems~57, 59) and without
this condition (Theorems~64, 65).
The considered expansions converge in the mean-square sense 
and contain only one operation of the limit transition
in contrast to its existing analogues.
The mentioned iterated Stratonovich stochastic integrals
are part of the 
Taylor--Stratonovich expansion. 
Therefore, the results of the article
can be applied to numerical
integration of Ito SDEs.
The results of the article were reformulated in the form of
theorems of the Wong--Zakai
type for iterated Stratonovich stochastic integrals.

\medskip

\end{quote}
}

\linespread{1.40}

\tableofcontents

\linespread{1.0}

\section{Introduction}

\vspace{5mm}

Let $(\Omega,$ ${\rm F},$ ${\sf P})$ be a complete probability space, let 
$\{{\rm F}_t, t\in[0,T]\}$ be a nondecreasing 
right-continous family of $\sigma$-algebras of ${\rm F},$
and let ${\bf f}_t$ be a standard $m$-dimensional Wiener 
stochastic process, which is
${\rm F}_t$-measurable for any $t\in[0, T].$ We assume that the components
${\bf f}_{t}^{(i)}$ $(i=1,\ldots,m)$ of this process are independent. 
Consider
an Ito stochastic differential equation (SDE) in the integral form

\vspace{-1mm}
\begin{equation}
\label{1.5.2}
{\bf x}_t={\bf x}_0+\int\limits_0^t {\bf a}({\bf x}_{\tau},\tau)d\tau+
\int\limits_0^t B({\bf x}_{\tau},\tau)d{\bf f}_{\tau},\ \ \
{\bf x}_0={\bf x}(0,\omega).
\end{equation}

\vspace{3mm}
\noindent
Here ${\bf x}_t$ is the $n$-dimensional stochastic process 
satisfying the equation (\ref{1.5.2}). 
The nonrandom functions ${\bf a}: \mathbb{R}^n\times[0, T]\to\mathbb{R}^n$,
$B: \mathbb{R}^n\times[0, T]\to\mathbb{R}^{n\times m}$
guarantee the existence and uniqueness up to stochastic equivalence 
of a strong solution
of the equation (\ref{1.5.2}) \cite{1}. The second integral on the 
right-hand side of (\ref{1.5.2}) is 
interpreted as the Ito stochastic integral.
Let ${\bf x}_0$ be an $n$-dimensional random variable, which is 
${\rm F}_0$-measurable and 
${\sf M}\bigl\{\left|{\bf x}_0\right|^2\bigr\}<\infty$ 
(${\sf M}$ denotes a mathematical expectation).
We assume that
${\bf x}_0$ and ${\bf f}_t-{\bf f}_0$ are independent when $t>0.$

It is well known \cite{KlPl2}-\cite{KPS}
that Ito SDEs are 
adequate mathematical models of dynamic systems of various 
physical nature under 
the influence of random disturbances. One of the effective approaches 
to the numerical integration of 
Ito SDEs is an approach based on 
the Taylor--Ito and 
Taylor--Stratonovich expansions
\cite{KlPl2}-\cite{12aa-afterxxx}. 
The most important feature of such 
expansions is a presence in them of the so-called iterated
Ito and Stratonovich stochastic integrals, which play the key 
role for solving the 
problem of numerical integration of Ito SDEs
and have the 
following form 

\vspace{-1mm}
\begin{equation}
\label{ito}
J[\psi^{(k)}]_{T,t}=\int\limits_t^T\psi_k(t_k) \ldots 
\int\limits_t^{t_{2}}
\psi_1(t_1) d{\bf w}_{t_1}^{(i_1)}\ldots
d{\bf w}_{t_k}^{(i_k)},
\end{equation}

\vspace{1mm}
\begin{equation}
\label{str}
J^{*}[\psi^{(k)}]_{T,t}=
{\int\limits_t^{*}}^T
\psi_k(t_k) \ldots 
{\int\limits_t^{*}}^{t_2}
\psi_1(t_1) d{\bf w}_{t_1}^{(i_1)}\ldots
d{\bf w}_{t_k}^{(i_k)},
\end{equation}

\vspace{3mm}
\noindent
where $\psi_1(\tau),\ldots,\psi_k(\tau): [t, T]\to\mathbb{R}$ are nonrandom
functions, ${\bf w}_{\tau}^{(i)}={\bf f}_{\tau}^{(i)}$
for $i=1,\ldots,m$ and
${\bf w}_{\tau}^{(0)}=\tau,$\ \
$i_1,\ldots,i_k = 0,1,\ldots,m,$
$$
\int\limits\ \hbox{and}\ \int\limits^{*}
$$ 

\vspace{3mm}
\noindent
denote Ito and 
Stratonovich stochastic integrals,
respectively.
In this paper we mainly use the definition of the Stratonovich 
stochastic integral from \cite{KlPl2} (also see \cite{2018a}, Sect.~2.1.1).

Note that $\psi_l(\tau)\equiv 1$ $(l=1,\ldots,k)$ and
$i_1,\ldots,i_k = 0,1,\ldots,m$ in  
\cite{KlPl2}-\cite{KlPl1}. At the same time
$\psi_l(\tau)\equiv (t-\tau)^{q_l}$ ($l=1,\ldots,k$,\ \
$q_1,\ldots,q_k=0, 1,\ldots $) and  
$i_1,\ldots,i_k = 1,\ldots,m$ in
\cite{kk5}-\cite{12aa-afterxxx}.

Effective solution 
of the problem of
combined mean-square approximation for collections 
of iterated Ito and Stratonovich stochastic integrals
(\ref{ito}) and (\ref{str}) 
composes the subject of the article.

We want to mention in short that there are 
two main criteria of numerical methods convergence 
for Ito SDEs \cite{KlPl2}-\cite{Mi3}:  
a strong or mean-square
criterion and a 
weak criterion where the subject of approximation is not the solution 
of Ito SDE, simply stated, but the 
distribution of Ito SDE solution.

Using the strong numerical methods, we may build
sample pathes
of Ito SDEs numerically. 
These methods require the combined mean-square approximation for collections 
of iterated Ito and Stratonovich stochastic integrals
(\ref{ito}) and (\ref{str}). 
The strong numerical methods are used when building new mathematical 
models on the basis of Ito SDEs
and solving
various mathematical problems connected with Ito SDEs.
Among these problems we mention signal 
filtering in the background of random
noise, stochastic optimal control, stochastic stability, evaluating 
the parameters of stochastic systems,
etc. \cite{KlPl2}-\cite{KPS}.

The problem of effective jointly numerical modeling 
(in accordance to the mean-square convergence criterion) of iterated 
Ito and Stratonovich stochastic integrals 
(\ref{ito}) and (\ref{str}) is 
difficult from 
theoretical and computing point of view \cite{KlPl2}-\cite{KPS},
\cite{2006}-\cite{rr}.

The only exception is connected with a narrow particular case, when 
$i_1=\ldots=i_k\ne 0$ and
$\psi_1(\tau),\ldots,\psi_k(\tau)\equiv \psi(\tau)$.
This case allows 
the investigation with using of the Ito formula 
\cite{KlPl2}-\cite{Mi3}.

Note that even for the mentioned coincidence ($i_1=\ldots=i_k\ne 0$), 
but for different 
functions $\psi_1(\tau),\ldots,\psi_k(\tau)$ the mentioned 
difficulties persist, and 
relatively simple families of 
iterated Ito and Stratonovich stochastic integrals, 
which can be often 
met in the applications, cannot be represented effectively in a finite 
form (within the framework of the mean-square approximation) 
using the system of standard 
Gaussian random variables.

Note that for a number of special types of Ito SDEs 
the problem of approximation of iterated Ito and Stratonovich
stochastic integrals may be simplified but cannot be solved. The equations
with additive scalar noise, with additive vector noise,
with non-additive  
scalar noise, with a small parameter are related to such 
types of equations \cite{KlPl2}-\cite{Mi3}. 
For the mentioned types of equations, simplifications 
are connected to the fact 
that some members
from stochastic Taylor expansions (Taylor--Ito and Taylor--Stratonovich 
expansions) are equal to zero 
or we may neglect some members (which include difficult for approximation 
iterated stochastic integrals) from these expansions
due to the presence of a small 
parameter \cite{KlPl2}-\cite{Mi3}.
In this article, we consider Ito SDEs 
with multidimentional and non-additive noise (non-commutative case).

Consider a brief overview of existing methods 
of the mean-square approximation of
iterated Ito and Stratonovich stochastic integrals.

Seems that iterated stochastic integrals may be approximated by multiple 
integral sums \cite{Mi2}, \cite{Mi3}, \cite{allen}. 
However, this approach implies the partitioning of the interval 
of integration $[t, T]$ for iterated stochastic integrals. The length
$T-t$ of this interval is already fairly
small (because it is a step of integration of numerical methods for 
Ito SDEs) and does not need to be partitioned.
Com\-pu\-ta\-ti\-onal experiments show that the application
of numerical simulation for iterated stochastic integrals 
(in which the interval of integration is
partitioned) leads to unacceptably high com\-pu\-ta\-ti\-onal cost and 
accumulation of computation errors \cite{2006}.

In \cite{Mi2} (also see \cite{KlPl2}, \cite{Mi3})
Milstein G.N. proposed to expand the integral (\ref{ito}) of multiplicity 2
($\psi_1(\tau),$ $\psi_2(\tau)\equiv 1$ and
$i_1,i_2=1,\ldots,m$)
into the iterated series of products
of standard Gaussian random variables by representing the Brownian 
bridge process
as the trigonometric Fourier series with random coefficients 
(the version of the so-called Karhunen--Loeve expansion
for the Brownian bridge process).
To obtain the Milstein expansion of (\ref{ito}) or (\ref{str}), 
the truncated Fourier
expansions of components of the Wiener process ${\bf f}_s$ must be
iteratively substituted in the single integrals, and the integrals
must be calculated, starting from the innermost integral.
This is a complicated procedure that does not lead to general
expansions of the integrals (\ref{ito}), (\ref{str}) 
of arbitrary multiplicity $k.$
For this reason, only expansions of single, double, and triple
stochastic integrals were presented 
in \cite{KlPl2} ($k=1, 2, 3$)
and in \cite{Mi2}, \cite{Mi3} ($k=1, 2$) 
for the simplest case $\psi_1(\tau), \psi_2(\tau), \psi_3(\tau)\equiv 1$ and
$i_1, i_2, i_3=0,1,\ldots,m.$
Moreover, the authors of the works
\cite{KlPl2}
(Sect.~5.8, pp.~202--204), \cite{KPS} (pp.~82-84),
\cite{Zapad-2} (pp.~438-439),  
\cite{Zapad-9} (pp.~263-264) use 
the Wong--Zakai approximation 
\cite{W-Z-1}, \cite{W-Z-2}, 
\cite{Watanabe} (without rigorous proof) within the frames
of the Milstein approach 
\cite{Mi2} based on the series expansion 
of the Brownian bridge process. See discussion in Sect.~15 of 
this paper for details.

Note that in \cite{rr} the method of approximation
of the double Ito stochastic integral (\ref{ito}) 
($\psi_1(\tau),$ $\psi_2(\tau)$ $\equiv 1$ and $i_1, i_2 =1,\ldots,m$) 
based on expansion
of the Wiener process using Haar functions and 
trigonometric functions has been considered.
The restrictions of the method \cite{rr} as well as 
the Milstein approach \cite{Mi2} are connected
with the iterated application of the operation
of limit transition at least starting from the second (in general case) 
and third 
multiplicity of iterated stochastic integrals.

It is necessary to note that the Milstein approach \cite{Mi2} 
excelled in several times or even in several orders
the methods of multiple integral sums \cite{Mi2}, \cite{Mi3}, \cite{allen}
(we mean here the diminishing of com\-pu\-ta\-ti\-onal costs).

An alternative strong approximation method (see Theorems 1 and 2 below) was 
proposed for (\ref{str}) in 
\cite{2018a}, Sect.~2.4 (also see \cite{2018}, \cite{2018aa}-\cite{12aa-afterxxx}, 
\cite{2010-1}-\cite{2013}, 
\cite{300a} (1997), \cite{400a} (1998),
\cite{arxiv-23}),
where $J^{*}[\psi^{(k)}]_{T,t}$ was represented as the multiple stochastic 
integral
from the certain discontinuous nonrandom function of $k$ variables, and the 
function
was then expressed as the iterated trigonometric Fourier series.
As a result,
an iterated series expansion of the integral 
(\ref{str}) in terms of products
of standard Gaussian random variables was obtained in 
\cite{2018a}, Sect.~2.4 (also see \cite{2018}, \cite{2018aa}-\cite{12aa-afterxxx}, 
\cite{2010-1}-\cite{2013}, 
\cite{300a} (1997), \cite{400a} (1998),
\cite{arxiv-23}) for an arbitrary multiplicity $k.$
This method is also valid for iterated Fourier--Legendre series
(at least for the case $k=2$) \cite{2018a} (Sect.~2.4.1).
Hereinafter, this method is referred to as the method of generalized
iterated Fourier series.

It was shown in \cite{2018a} (also see \cite{2018}, \cite{2018aa}-\cite{12aa-afterxxx}, 
\cite{2010-1}-\cite{2013}, 
\cite{300a} (1997), \cite{400a} (1998),
\cite{arxiv-23}) that the method of 
generalized iterated Fourier series leads to the 
Milstein expansion \cite{Mi2} of the integral (\ref{str})
in the case of trigonometric system 
and to the substantially simpler expansion 
of the integral (\ref{str}) in the case
of Legendre polynomials system (at least for the case
of multiplicity $k=2$ of the integral (\ref{str}), $i_1\ne i_2$).

Note that the method of generalized 
iterated Fourier series as well as the Milstein approach
\cite{Mi2}
leads to iterated application of the operation of limit transition. 
As mentioned above, this problem appears for triple stochastic integrals
($i_1, i_2, i_3=1,\ldots,m$)
or even for some double stochastic integrals 
in the case, when $\psi_1(\tau),$ 
$\psi_2(\tau)\not\equiv 1$ ($i_1, i_2=1,\ldots,m$)
\cite{2006}.

The mentioned problem (iterated application of the operation 
of limit transition) not appears 
in the method, which 
is considered for the integrals (\ref{ito}) in Theorems 3, 4 (see below)
\cite{2006} (2006) \cite{2011-2}-\cite{200a},
\cite{301a}-\cite{arxiv-8},
\cite{arxiv-12}, 
\cite{arxiv-24}-\cite{arxiv-9}.
The idea of this method is as follows: 
the iterated Ito stochastic 
integral (\ref{ito}) of multiplicity $k$ is represented as 
the multiple stochastic 
integral from the certain discontinuous nonrandom function of $k$ variables
defined on the hypercube $[t, T]^k$, where $[t, T]$ is the interval of 
integration of the iterated Ito stochastic 
integral (\ref{ito}). Then, 
the indicated 
nonrandom function is expanded in the hypercube $[t, T]^k$
into the generalized 
multiple Fourier series converging 
in the mean-square sense
in the space 
$L_2([t,T]^k)$. After a number of nontrivial transformations we come 
(see Theorems 3, 4 below) to the 
mean-square convergening expansion of the
iterated Ito stochastic 
integral (\ref{ito})
into the multiple 
series of products
of standard  Gaussian random 
variables. The coefficients of this 
series are the coefficients of 
generalized multiple Fourier series for the mentioned nonrandom function 
of $k$ variables, which can be calculated using the explicit formula 
regardless of the multiplicity $k$ of the
iterated Ito stochastic 
integral (\ref{ito}).
Hereinafter, this method is referred to as the method of generalized
multiple Fourier series.

Thus, we obtain the following useful possibilities
of the method of generalized multiple Fourier series.

1. There is the explicit formula (see (\ref{ppppa})) for calculation 
of expansion coefficients 
of the iterated Ito stochastic integral (\ref{ito}) with any
fixed multiplicity $k$. 

2. We have new possibilities for exact calculation of the mean-square 
error of approximation 
of the iterated Ito stochastic integral (\ref{ito})
\cite{2017}-\cite{12aa-afterxxx}, \cite{17a}, \cite{arxiv-2}. 

3. Since the used
multiple Fourier series is generalized in the sense
that it is built using various complete orthonormal
systems of functions in the space $L_2([t, T])$, then we 
have new possibilities 
for approximation --- we can 
use not only the trigonometric functions as in \cite{KlPl2}-\cite{Mi3}
but the Legendre polynomials as well as the systems of Haar
and Rademacher--Walsh functions.

4. As it turned out \cite{2006}-\cite{200a},
\cite{301a}-\cite{arxiv-8},
\cite{arxiv-12}, 
\cite{arxiv-24}-\cite{arxiv-9} it is more convenient to work 
with the Legendre polynomials for approximation
of the iterated Ito stochastic integrals (\ref{ito}). 
Approximations based on the Legendre polynomials essentially simpler 
than their analogues based on the trigonometric functions.
Another advantages of the application of Legendre polynomials 
in the framework of the mentioned direction are considered
in \cite{29a}, \cite{301a} (also see \cite{2018a}-\cite{12aa-afterxxx}).

5. The approach based on the Karhunen--Loeve expansion
of the Brownian bridge process (also see \cite{rr})
leads to 
iterated application of the operation of limit
transition (the operation of limit transition 
is implemented only once in Theorems 3, 4 (see below))
starting from  
the second or
third multiplicity
of the iterated Ito stochastic integrals (\ref{ito}).
Multiple series (the operation of limit transition 
is implemented only once) are more convenient 
for approximation than the iterated ones
(iterated application of the operation of limit
transition), 
since partial sums of multiple series converge for any possible case of  
convergence to infinity of their upper limits of summation 
(let us denote them as $p_1,\ldots, p_k$). 
For example, when $p_1=\ldots=p_k=p\to\infty$. 
For iterated series, the condition $p_1=\ldots=p_k=p\to\infty$
obviously does not guarantee the convergence of this series.
However, 
in \cite{KlPl2}
(Sect.~5.8, pp.~202--204), \cite{KPS} (pp.~82-84),
\cite{Zapad-2} (pp.~438-439),  
\cite{Zapad-9} (pp.~263-264) 
the authors use (without rigorous proof)
the condition $p_1=p_2=p_3=p\to\infty$
within the frames of the mentioned approach
based on the Karhunen--Loeve expansion of the Brownian bridge
process \cite{Mi2} together with the Wong--Zakai approximation
\cite{W-Z-1}, \cite{W-Z-2}, \cite{Watanabe} (see Sect.~15 for details).

6. As it turned out, the 
method of generalized
multiple Fourier series
can be adapted for the iterated
Stratonovich stochastic integrals (\ref{str}) at least
for multiplicities 1 to 6 \cite{2011-2}-\cite{12aa-afterxxx}, 
\cite{2010-2}-\cite{2013}, \cite{30a}, \cite{300a},
\cite{400a}, \cite{271a},
\cite{arxiv-5}-\cite{arxiv-8}, \cite{arxiv-23}, \cite{arxiv-9}. 
Expansions of these iterated Stratonovich 
stochastic integrals turned out
to be simpler (see Theorems 6--13, 15--18, 23, 24, 42--46, 48, 49, 51, 53, 55, 60, 62, 63 below) 
than the appropriate expansions
of the iterated Ito stochastic integrals (\ref{ito}) from Theorems 3, 4.

7. The method of generalized multiple Fourier series
has been applied for some other types of iterated stochastic integrals
(iterated stochastic integrals with respect to martingale 
Poisson random measures
and iterated stochastic integrals with respect to martingales)
as well as for approximation of iterated stochastic 
integrals with respect to the infinite-dimensional
$Q$-Wiener process \cite{2018a} (Chapters 1 and 7).

In this article, we collect more than thirty
theorems formulated and proved by the author
that develop the mentioned direction of investigations.
Moreover, on the base of the presented theorems, we formulate
several hypotheses (Hypotheses 1--8) on expansions of iterated Stratonovich 
stochastic integrals of arbitrary multiplicity $k.$
The results of the article prove (the cases of Legendre
polynomials and trigonometric functions)
Hypothesis 1 for 
$k=1,\ldots,8,$ Hypothesis 2 for $k=1,\ldots,5$ and Hypothesis 3 for $k=1,\ldots,3.$
Moreover, the proof of Hypothesis~8 
is given for the case 
of an arbitrary complete orthonormal system of functions in $L_2([t, T])$
($p_1=\ldots=p_k=p,$ $k\in\mathbb{N}$) but under one additional condition 
(Theorems~57, 59). Also, Hypothesis~8 is proved for $k=1,\ldots,6$
and various restrictions on the weight functions. Recently, 
Hypothesis~8 is proved for 
the case of an arbitrary complete orthonormal system of functions in $L_2([t, T])$
($p_1=\ldots=p_k=p,$ $k\in\mathbb{N}$) and without the aforementioned additional condition 
(Theorems~64, 65).

\vspace{5mm}

\section{Hypotheses on Expansions of Iterated Stratonovich Stochastic
Integrals of Arbitrary Multiplicity $k$}

\vspace{5mm}

Taking into account Theorems 1--13, 15--18, 23 (see below), let us 
formulate the following hypotheses on expansions of iterated 
Stratonovich stochastic
integrals of arbitrary multiplicity $k$.

\vspace{2mm}

{\bf Hypothesis 1}\ \cite{2011-2}-\cite{12aa-afterxxx}. {\it Assume that
$\{\phi_j(x)\}_{j=0}^{\infty}$ is a complete orthonormal
system of Legendre polynomials or trigonometric functions
in the space $L_2([t, T])$.
Then, for the iterated Stratonovich stochastic 
integral of $k$th multiplicity

\begin{equation}
\label{otitgo10}
I_{T,t}^{*(i_1\ldots i_k)}=
{\int\limits_t^{*}}^T
\ldots
{\int\limits_t^{*}}^{t_2}
d{\bf w}_{t_1}^{(i_1)}
\ldots d{\bf w}_{t_k}^{(i_k)}\ \ \ (i_1,\ldots, i_k=0, 1,\ldots,m)
\end{equation}

\vspace{4mm}
\noindent
the following 
expansion

\vspace{-2mm}
\begin{equation}
\label{feto1900otita}
I_{T,t}^{*(i_1\ldots i_k)}=
\hbox{\vtop{\offinterlineskip\halign{
\hfil#\hfil\cr
{\rm l.i.m.}\cr
$\stackrel{}{{}_{p\to \infty}}$\cr
}} }
\sum\limits_{j_1,\ldots j_k=0}^{p}
C_{j_k \ldots j_1}\zeta_{j_1}^{(i_1)}
\ldots
\zeta_{j_k}^{(i_k)}
\end{equation}

\vspace{4mm}
\noindent
converging in the mean-square sense is valid, where the Fourier coefficient 
$C_{j_k \ldots j_1}$ has the form

\vspace{-1mm}
$$
C_{j_k \ldots j_1}=\int\limits_t^T\phi_{j_k}(t_k)\ldots
\int\limits_t^{t_2}
\phi_{j_1}(t_1)
dt_1\ldots dt_k,
$$

\vspace{4mm}
\noindent
${\rm l.i.m.}$ is a limit in the mean-square sense,

\vspace{-1mm}
$$
\zeta_{j}^{(i)}=
\int\limits_t^T \phi_{j}(s) d{\bf w}_s^{(i)}
$$ 

\vspace{2mm}
\noindent
are independent standard Gaussian random variables for various 
$i$ or $j$ {\rm (}if $i\ne 0${\rm )},
${\bf w}_{\tau}^{(i)}={\bf f}_{\tau}^{(i)}$ are independent 
standard Wiener processes
$(i=1,\ldots,m)$ and 
${\bf w}_{\tau}^{(0)}=\tau.$}

\vspace{2mm}

Hypothesis 1 allows to approximate the iterated
Stratonovich stochastic integral 
$I_{T,t}^{*(i_1\ldots i_k)}$
by the sum

\vspace{-1mm}
\begin{equation}
\label{otit567}
I_{T,t}^{*(i_1\ldots i_k)p}=
\sum\limits_{j_1,\ldots j_k=0}^{p}
C_{j_k \ldots j_1}\zeta_{j_1}^{(i_1)}
\ldots
\zeta_{j_k}^{(i_k)},
\end{equation}

\vspace{4mm}
\noindent
where
$$
\lim_{p\to\infty}{\sf M}\left\{\Biggl(
I_{T,t}^{*(i_1 \ldots i_k)}-
I_{T,t}
^{*(i_1\ldots i_k)p}\Biggr)^2\right\}=0.
$$

\vspace{4mm}

The integrals (\ref{otitgo10}) are integrals from the 
Taylor--Stratonovich expansion \cite{KlPl2}. It means that
the approximations (\ref{otit567})
can be very useful for the numerical integration
of Ito SDEs.
The expansion (\ref{feto1900otita}) contains only one operation of the limit
transition and by this reason is convenient for approximation
of iterated Stratonovich stochastic integrals.
Moreover, the author supposes that the analogue of Hypothesis 1 will be
valid for the iterated Stratonovich stochastic integrals (\ref{str}).

\vspace{2mm}

{\bf Hypothesis 2}\ \cite{2018}-\cite{12aa-afterxxx}. {\it Assume that
$\{\phi_j(x)\}_{j=0}^{\infty}$ is a complete orthonormal
system of Legendre polynomials or trigonometric functions
in the space $L_2([t, T])$. Moreover,
every $\psi_l(\tau)$ $(l=1,\ldots,k)$ is 
an enough smooth nonrandom function
on $[t,T].$
Then, for the iterated Stratonovich stochastic integral {\rm (\ref{str})}
of $k$th multiplicity

\vspace{-1mm}
$$
J^{*}[\psi^{(k)}]_{T,t}=
{\int\limits_t^{*}}^T
\psi_k(t_k) \ldots 
{\int\limits_t^{*}}^{t_2}
\psi_1(t_1) d{\bf w}_{t_1}^{(i_1)}\ldots
d{\bf w}_{t_k}^{(i_k)}\ \ \ (i_1,\ldots, i_k=0, 1,\ldots,m)
$$

\vspace{4mm}
\noindent
the following 
expansion 

\vspace{-2mm}
\begin{equation}
\label{feto1900otitarrr}
J^{*}[\psi^{(k)}]_{T,t}=
\hbox{\vtop{\offinterlineskip\halign{
\hfil#\hfil\cr
{\rm l.i.m.}\cr
$\stackrel{}{{}_{p\to \infty}}$\cr
}} }
\sum\limits_{j_1,\ldots j_k=0}^{p}
C_{j_k \ldots j_1}\zeta_{j_1}^{(i_1)}
\ldots
\zeta_{j_k}^{(i_k)}
\end{equation}

\vspace{4mm}
\noindent
converging in the mean-square sense 
is valid, where the Fourier coefficient 
$C_{j_k \ldots j_1}$ has the form

\vspace{-1mm}
$$
C_{j_k \ldots j_1}=\int\limits_t^T\psi_k(t_k)\phi_{j_k}(t_k)\ldots
\int\limits_t^{t_2}
\psi_1(t_1)\phi_{j_1}(t_1)
dt_1\ldots dt_k,
$$

\vspace{4mm}
\noindent
${\rm l.i.m.}$ is a limit in the mean-square sense,

\vspace{-1mm}
$$
\zeta_{j}^{(i)}=
\int\limits_t^T \phi_{j}(s) d{\bf w}_s^{(i)}
$$ 

\vspace{2mm}
\noindent
are independent standard Gaussian random variables for various 
$i$ or $j$ {\rm (}if $i\ne 0${\rm )},
${\bf w}_{\tau}^{(i)}={\bf f}_{\tau}^{(i)}$ are independent 
standard Wiener processes
$(i=1,\ldots,m)$ and 
${\bf w}_{\tau}^{(0)}=\tau.$}

\vspace{2mm}

Hypothesis 2 allows to approximate the iterated
Stratonovich stochastic integral 
$J^{*}[\psi^{(k)}]_{T,t}$
by the sum

\vspace{-1mm}
\begin{equation}
\label{otit567r}
J^{*}[\psi^{(k)}]_{T,t}^p=
\sum\limits_{j_1,\ldots j_k=0}^{p}
C_{j_k \ldots j_1}\zeta_{j_1}^{(i_1)}
\ldots
\zeta_{j_k}^{(i_k)},
\end{equation}

\vspace{4mm}
\noindent
where
$$
\lim_{p\to\infty}{\sf M}\left\{\Biggl(
J^{*}[\psi^{(k)}]_{T,t}-
J^{*}[\psi^{(k)}]_{T,t}^p\Biggl)^2\right\}=0.
$$

\vspace{4mm}

Let us consider the more general statement, then 
Hypotheses 1 and 2.

\vspace{2mm}                                  

{\bf Hypothesis 3}\ \cite{2018a}-\cite{12aa-afterxxx}.\ {\it Assume that
$\{\phi_j(x)\}_{j=0}^{\infty}$ is a complete orthonormal
system of Legendre polynomials or trigonometric functions
in the space $L_2([t, T])$. Moreover,
every $\psi_l(\tau)$ $(l=1,\ldots,k)$ is 
an enough smooth nonrandom function
on $[t,T].$
Then, for the iterated Stratonovich stochastic integral {\rm (\ref{str})}
of $k$th multiplicity

$$
J^{*}[\psi^{(k)}]_{T,t}=
{\int\limits_t^{*}}^T
\psi_k(t_k) \ldots 
{\int\limits_t^{*}}^{t_2}
\psi_1(t_1) d{\bf w}_{t_1}^{(i_1)}\ldots
d{\bf w}_{t_k}^{(i_k)}\ \ \ (i_1,\ldots, i_k=0, 1,\ldots,m)
$$

\vspace{3mm}
\noindent
the following 
expansion 

\vspace{-2mm}
\begin{equation}
\label{feto1900otitarrri}
J^{*}[\psi^{(k)}]_{T,t}=
\hbox{\vtop{\offinterlineskip\halign{
\hfil#\hfil\cr
{\rm l.i.m.}\cr
$\stackrel{}{{}_{p_1,\ldots,p_k\to \infty}}$\cr
}} }
\sum\limits_{j_1=0}^{p_1}
\ldots\sum\limits_{j_k=0}^{p_k}
C_{j_k \ldots j_1}\zeta_{j_1}^{(i_1)}
\ldots
\zeta_{j_k}^{(i_k)}
\end{equation}

\vspace{4mm}
\noindent
converging in the mean-square sense 
is valid, where the Fourier coefficient 
$C_{j_k \ldots j_1}$ has the form

\vspace{-1mm}
$$
C_{j_k \ldots j_1}=\int\limits_t^T\psi_k(t_k)\phi_{j_k}(t_k)\ldots
\int\limits_t^{t_2}
\psi_1(t_1)\phi_{j_1}(t_1)
dt_1\ldots dt_k,
$$

\vspace{3mm}
\noindent
${\rm l.i.m.}$ is a limit in the mean-square sense,

\vspace{-1mm}
$$
\zeta_{j}^{(i)}=
\int\limits_t^T \phi_{j}(s) d{\bf w}_s^{(i)}
$$ 

\vspace{2mm}
\noindent
are independent standard Gaussian random variables for various 
$i$ or $j$ {\rm (}if $i\ne 0${\rm ),}
${\bf w}_{\tau}^{(i)}={\bf f}_{\tau}^{(i)}$ are independent 
standard Wiener processes
$(i=1,\ldots,m)$ and 
${\bf w}_{\tau}^{(0)}=\tau.$}

\vspace{2mm}

In the next section, we consider two theorems on expansions
of iterated Stratonovich stochastic integrals of arbitrary multiplicity $k$.
Expansions from Theorems 1 and 2 (see below) contain an iterated
operation of the limit transition in comparison with Hypotheses 1--3.
This feature creates some difficulties when estimating
the mean-square approximation error of iterated
Stratonovich stochastic integrals. On the other hand, Theorems 1 and 2
contain the same expansion terms of iterated 
Stratonovich stochastic integrals as in Hypotheses 1--3.

\vspace{5mm}

\section{Expansions of Iterated Stratonovich Stochastic Integrals of 
Arbitrary Multiplicity Based on 
Iterated Fourier Series
Converging Pointwise}

\vspace{5mm}

Let us formulate the following theorem. 

\vspace{2mm}

{\bf Theorem 1} \cite{2018a} (Sect.~2.4) (also see \cite{2018}, \cite{2018aa}-\cite{12aa-afterxxx}, 
\cite{2010-1}-\cite{2013}, 
\cite{300a} (1997), \cite{400a} (1998),
\cite{arxiv-23}).\ 
{\it Suppose that the functions $\psi_1(\tau),\ldots,\psi_k(\tau)$ are twice continuously
differentiable at the interval
$[t, T]$ and
$\{\phi_j(x)\}_{j=0}^{\infty}$ is a complete
orthonormal system
of trigonometric functions in the space $L_2([t, T])$. 
Then, the iterated Stratonovich stochastic integral

$$
J^{*}[\psi^{(k)}]_{T,t}=
{\int\limits_t^{*}}^T
\psi_k(t_k) \ldots 
{\int\limits_t^{*}}^{t_2}
\psi_1(t_1) d{\bf w}_{t_1}^{(i_1)}\ldots
d{\bf w}_{t_k}^{(i_k)}\ \ \ (i_1,\ldots, i_k=0, 1,\ldots,m)
$$

\vspace{4mm}
\noindent
is expanded into the 
converging 
in the mean of degree $2n$ $(n\in \mathbb{N})$
iterated series

\begin{equation}
\label{1500}
J^{*}[\psi^{(k)}]_{T,t}=
\sum_{j_1=0}^{\infty}\ldots\sum_{j_k=0}^{\infty}
C_{j_k\ldots j_1}
\zeta^{(i_1)}_{j_1}\ldots \zeta^{(i_k)}_{j_k},
\end{equation}

\vspace{3mm}
\noindent
i.e.

\vspace{-3mm}
$$
\lim\limits_{p_1\to\infty}\varlimsup\limits_{p_2\to\infty}
\ldots\varlimsup\limits_{p_k\to\infty}
{\sf M}\left\{\Biggl(J^{*}[\psi^{(k)}]_{T,t}-
\sum_{j_1=0}^{p_1}\ldots\sum_{j_k=0}^{p_k}
C_{j_k\ldots j_1}
\zeta^{(i_1)}_{j_1}\ldots 
\zeta^{(i_k)}_{j_k}\Biggr)^{2n}\right\}=0,
$$

\vspace{5mm}
\noindent
where $\varlimsup$ means $\limsup$,

\vspace{-1mm}
$$
\zeta_{j}^{(i)}=
\int\limits_t^T \phi_{j}(s) d{\bf w}_s^{(i)}
$$ 

\vspace{2mm}
\noindent
are independent standard Gaussian random variables
for 
various
$i$ or $j$ {\rm(}if $i\ne 0${\rm),}
the Fourier coefficient $C_{j_k\ldots j_1}$ has the form

\vspace{-1mm}
$$
C_{j_k\ldots j_1}=\int\limits_{[t,T]^k}
K^{*}(t_1,\ldots,t_k)\prod_{l=1}^{k} \phi_{j_l}(t_l)
dt_1\ldots dt_k=
$$

\vspace{2mm}
$$
=
\int\limits_t^T\psi_k(t_k)\phi_{j_k}(t_k)\ldots
\int\limits_t^{t_2}
\psi_1(t_1)\phi_{j_1}(t_1)
dt_1\ldots dt_k,
$$

\vspace{3mm}
\noindent
where

\vspace{-2mm}
$$
K^{*}(t_1,\ldots,t_k)=\prod\limits_{l=1}^k\psi_l(t_l)
\prod_{l=1}^{k-1}\Biggl({\bf 1}_{\{t_l<t_{l+1}\}}+
\frac{1}{2}{\bf 1}_{\{t_l=t_{l+1}\}}\Biggr),\ \ \ 
t_1,\ldots,t_k\in[t, T]
$$

\vspace{5mm}
\noindent
for $k\ge 2$ and 
$K^{*}(t_1)\equiv\psi_1(t_1)$ for $t_1\in[t, T].$ 
Here ${\bf 1}_A$ is the indicator of the set $A$. 
}

\vspace{2mm}

Let us consider the following theorem.

\vspace{2mm}

{\bf Theorem 2}\ \cite{2018a} (Sect.~2.4) (also see \cite{2018}, \cite{2018aa}-\cite{12aa-afterxxx}, 
\cite{2013} (2013), \cite{arxiv-23}).\
{\it Suppose that the functions $\psi_1(\tau),\ldots,\psi_k(\tau)$ are twice continuously
differentiable at the interval
$[t, T]$ and
$\{\phi_j(x)\}_{j=0}^{\infty}$ is a complete
orthonormal system
of trigonometric functions in the space $L_2([t, T])$. 
Then, the iterated Stratonovich stochastic integral 

$$
J^{*}[\psi^{(k)}]_{T,t}=
{\int\limits_t^{*}}^T
\psi_k(t_k) \ldots 
{\int\limits_t^{*}}^{t_2}
\psi_1(t_1) d{\bf w}_{t_1}^{(i_1)}\ldots
d{\bf w}_{t_k}^{(i_k)}\ \ \ (i_1,\ldots, i_k=0, 1,\ldots,m)
$$

\vspace{4mm}
\noindent
is expanded into the 
converging 
in the mean of degree $2n$ $(n\in\mathbb{N})$
iterated series

\begin{equation}
\label{1500e}
J^{*}[\psi^{(k)}]_{T,t}=
\sum_{j_k=0}^{\infty}\ldots\sum_{j_1=0}^{\infty}
C_{j_k\ldots j_1}
\zeta^{(i_1)}_{j_1}\ldots\zeta^{(i_k)}_{j_k},
\end{equation}

\vspace{4mm}
\noindent
i.e.

\vspace{-4mm}
$$
\lim\limits_{p_k\to\infty}\varlimsup\limits_{p_{k-1}\to\infty}
\ldots
\varlimsup\limits_{p_1\to\infty}
{\sf M}\left\{\Biggl(J^{*}[\psi^{(k)}]_{T,t}-
\sum_{j_k=0}^{p_k}\ldots\sum_{j_1=0}^{p_1}
C_{j_k\ldots j_1}
\zeta^{(i_1)}_{j_1}\ldots 
\zeta^{(i_k)}_{j_k}\Biggr)^{2n}\right\}=0;
$$

\vspace{4mm}
\noindent
another notations are the same as in Theorem {\rm 1}.}

\vspace{2mm}

It is not difficult to see that 
the members of expansions (\ref{feto1900otitarrr}), (\ref{1500}),
and (\ref{1500e}) 
are the same.
However, as mentioned before, the expansion (\ref{feto1900otitarrr})
contains only one operation of the limit transition. At the same time 
the expansions (\ref{1500}), (\ref{1500e}) contain an iterated operation
of the limit transiton.

In \cite{2018a} (Sect.~2.4.1) it is shown that Theorems~1 and 2 
will remain valid for the case when 
$\{\phi_j(x)\}_{j=0}^{\infty}$ is a complete
orthonormal system
of Legendre polynomials in the space $L_2([t, T])$, $k=2$, 
and $2n=2$ (the case of mean-square convergence).
In this case the function $\psi_2(\tau)$ is continuously 
differentiable at the interval $[t, T]$ and the
function $\psi_1(\tau)$ is twice continuously 
differentiable at the interval $[t, T]$.

As it tured out, the approach considered in the next section gives the key
to the proof of Hypotheses 1--3.

\vspace{5mm}

\section{Expansion of Iterated Ito Stochastic 
Integrals of Arbitrary 
Multiplicity $k$ Based 
on Generalized Multiple Fourier Series Converging in the Mean}

\vspace{5mm}

Suppose that every $\psi_l(\tau)$ $(l=1,\ldots,k)$ is a continuous 
nonrandom function on $[t, T]$
(the case $\psi_1(\tau),\ldots,\psi_k(\tau)\in L_2([t, T])$
will be considered in Theorem~4 (see below)).
Define the following function on the hypercube $[t, T]^k$

\vspace{-1mm}
\begin{equation}
\label{ppp}
K(t_1,\ldots,t_k)=
\begin{cases}
\psi_1(t_1)\ldots \psi_k(t_k),\ &\hbox{for}\ \ t_1<\ldots<t_k\\
~\\
~\\
0,\ &\hbox{otherwise}
\end{cases},\ \ \ \ t_1,\ldots,t_k\in[t, T],\ \ \ \ k\ge 2,
\end{equation}

\vspace{4mm}
\noindent
and 
$K(t_1)\equiv\psi_1(t_1)$ for $t_1\in[t, T].$

Suppose that $\{\phi_j(x)\}_{j=0}^{\infty}$
is a complete orthonormal system of functions in the space
$L_2([t, T])$. 
The function $K(t_1,\ldots,t_k)$ is piecewise continuous in the 
hypercube $[t, T]^k.$
At this situation it is well known that the generalized 
multiple Fourier series 
of $K(t_1,\ldots,t_k)\in L_2([t, T]^k)$ is converging 
to $K(t_1,\ldots,t_k)$ in the hypercube $[t, T]^k$ in 
the mean-square sense, i.e.

\vspace{1mm}
$$
\hbox{\vtop{\offinterlineskip\halign{
\hfil#\hfil\cr
{\rm lim}\cr
$\stackrel{}{{}_{p_1,\ldots,p_k\to \infty}}$\cr
}} }\Biggl\Vert
K(t_1,\ldots,t_k)-
\sum_{j_1=0}^{p_1}\ldots \sum_{j_k=0}^{p_k}
C_{j_k\ldots j_1}\prod_{l=1}^{k} \phi_{j_l}(t_l)
\Biggr\Vert_{L_2([t,T]^k)}=0,
$$

\vspace{5mm}
\noindent
where
\begin{equation}
\label{ppppa}
C_{j_k\ldots j_1}=\int\limits_{[t,T]^k}
K(t_1,\ldots,t_k)\prod_{l=1}^{k}\phi_{j_l}(t_l)dt_1\ldots dt_k
\end{equation}

\vspace{4mm}
\noindent
is the Fourier coefficient,

\vspace{-1mm}
$$
\left\Vert f\right\Vert_{L_2([t,T]^k)}=\left(\int\limits_{[t,T]^k}
f^2(t_1,\ldots,t_k)dt_1\ldots dt_k\right)^{1/2}.
$$

\vspace{5mm}

Consider the partition $\{\tau_j\}_{j=0}^N$ of $[t,T]$ such that

\begin{equation}
\label{1111}
t=\tau_0<\ldots <\tau_N=T,\ \ \
\Delta_N=
\hbox{\vtop{\offinterlineskip\halign{
\hfil#\hfil\cr
{\rm max}\cr
$\stackrel{}{{}_{0\le j\le N-1}}$\cr
}} }\Delta\tau_j\to 0\ \ \hbox{if}\ \ N\to \infty,\ \ \
\Delta\tau_j=\tau_{j+1}-\tau_j.
\end{equation}

\vspace{4mm}

{\bf Theorem 3}\ \cite{2006} (2006), \cite{2011-2}-\cite{200a}, 
\cite{301a}-\cite{arxiv-8},
\cite{arxiv-12}, 
\cite{arxiv-24}-\cite{arxiv-9}. 
{\it Suppose that
every $\psi_l(\tau)$ $(l=1,\ldots, k)$ is a continuous nonrandom function on 
$[t, T]$ and
$\{\phi_j(x)\}_{j=0}^{\infty}$ is a complete orthonormal system  
of continuous functions in the space $L_2([t,T]).$ Then

\vspace{1mm}
$$
J[\psi^{(k)}]_{T,t}\  =\ 
\hbox{\vtop{\offinterlineskip\halign{
\hfil#\hfil\cr
{\rm l.i.m.}\cr
$\stackrel{}{{}_{p_1,\ldots,p_k\to \infty}}$\cr
}} }\sum_{j_1=0}^{p_1}\ldots\sum_{j_k=0}^{p_k}
C_{j_k\ldots j_1}\Biggl(
\prod_{l=1}^k\zeta_{j_l}^{(i_l)}\ -
\Biggr.
$$

\vspace{2mm}
\begin{equation}
\label{tyyy}
-\ \Biggl.
\hbox{\vtop{\offinterlineskip\halign{
\hfil#\hfil\cr
{\rm l.i.m.}\cr
$\stackrel{}{{}_{N\to \infty}}$\cr
}} }\sum_{(l_1,\ldots,l_k)\in {\rm G}_k}
\phi_{j_{1}}(\tau_{l_1})
\Delta{\bf w}_{\tau_{l_1}}^{(i_1)}\ldots
\phi_{j_{k}}(\tau_{l_k})
\Delta{\bf w}_{\tau_{l_k}}^{(i_k)}\Biggr),
\end{equation}

\vspace{6mm}
\noindent
where $J[\psi^{(k)}]_{T,t}$ is defined by {\rm (\ref{ito}),}

\vspace{-1mm}
$$
{\rm G}_k={\rm H}_k\backslash{\rm L}_k,\ \ \
{\rm H}_k=\{(l_1,\ldots,l_k):\ l_1,\ldots,l_k=0,\ 1,\ldots,N-1\},
$$

$$
{\rm L}_k=\{(l_1,\ldots,l_k):\ l_1,\ldots,l_k=0,\ 1,\ldots,N-1;\
l_g\ne l_r\ (g\ne r);\ g, r=1,\ldots,k\},
$$

\vspace{4mm}
\noindent
${\rm l.i.m.}$ is a limit in the mean-square sense$,$
$i_1,\ldots,i_k=0,1,\ldots,m,$

\vspace{-2mm}
$$
\zeta_{j}^{(i)}=
\int\limits_t^T \phi_{j}(\tau) d{\bf w}_{\tau}^{(i)}
$$

\vspace{2mm}
\noindent
are independent standard Gaussian random variables
for various
$i$ or $j$ {\rm(}if $i\ne 0${\rm),}
$C_{j_k\ldots j_1}$ is the Fourier coefficient {\rm(\ref{ppppa}),}
$\Delta{\bf w}_{\tau_{j}}^{(i)}=
{\bf w}_{\tau_{j+1}}^{(i)}-{\bf w}_{\tau_{j}}^{(i)}$
$(i=0, 1,\ldots,m),$
$\left\{\tau_{j}\right\}_{j=0}^{N}$ is a partition of
the interval $[t, T],$ which satisfies the condition {\rm (\ref{1111})}.
}

\vspace{2mm}

Let us consider the transformed particular cases of Theorem 3 
(see (\ref{tyyy})) for
$k=1$--$6$ \cite{2006}-\cite{200a}, \cite{301a}-\cite{arxiv-8},
\cite{arxiv-12}, 
\cite{arxiv-24}-\cite{arxiv-9} (the case $k=7$ can be found
in \cite{2011-2}-\cite{12aa-afterxxx}, \cite{arxiv-1})

\vspace{1mm}
\begin{equation}
\label{a1}
J[\psi^{(1)}]_{T,t}
=\hbox{\vtop{\offinterlineskip\halign{
\hfil#\hfil\cr
{\rm l.i.m.}\cr
$\stackrel{}{{}_{p_1\to \infty}}$\cr
}} }\sum_{j_1=0}^{p_1}
C_{j_1}\zeta_{j_1}^{(i_1)},
\end{equation}

\vspace{4mm}
\begin{equation}
\label{a2}
J[\psi^{(2)}]_{T,t}
=\hbox{\vtop{\offinterlineskip\halign{
\hfil#\hfil\cr
{\rm l.i.m.}\cr
$\stackrel{}{{}_{p_1,p_2\to \infty}}$\cr
}} }\sum_{j_1=0}^{p_1}\sum_{j_2=0}^{p_2}
C_{j_2j_1}\Biggl(\zeta_{j_1}^{(i_1)}\zeta_{j_2}^{(i_2)}
-{\bf 1}_{\{i_1=i_2\ne 0\}}
{\bf 1}_{\{j_1=j_2\}}\Biggr),
\end{equation}

\vspace{6mm}
$$
J[\psi^{(3)}]_{T,t}=
\hbox{\vtop{\offinterlineskip\halign{
\hfil#\hfil\cr
{\rm l.i.m.}\cr
$\stackrel{}{{}_{p_1,\ldots,p_3\to \infty}}$\cr
}} }\sum_{j_1=0}^{p_1}\sum_{j_2=0}^{p_2}\sum_{j_3=0}^{p_3}
C_{j_3j_2j_1}\Biggl(
\zeta_{j_1}^{(i_1)}\zeta_{j_2}^{(i_2)}\zeta_{j_3}^{(i_3)}
-\Biggr.
$$
\begin{equation}
\label{a3}
\Biggl.-{\bf 1}_{\{i_1=i_2\ne 0\}}
{\bf 1}_{\{j_1=j_2\}}
\zeta_{j_3}^{(i_3)}
-{\bf 1}_{\{i_2=i_3\ne 0\}}
{\bf 1}_{\{j_2=j_3\}}
\zeta_{j_1}^{(i_1)}-
{\bf 1}_{\{i_1=i_3\ne 0\}}
{\bf 1}_{\{j_1=j_3\}}
\zeta_{j_2}^{(i_2)}\Biggr),
\end{equation}

\vspace{8mm}
$$
J[\psi^{(4)}]_{T,t}
=
\hbox{\vtop{\offinterlineskip\halign{
\hfil#\hfil\cr
{\rm l.i.m.}\cr
$\stackrel{}{{}_{p_1,\ldots,p_4\to \infty}}$\cr
}} }\sum_{j_1=0}^{p_1}\ldots\sum_{j_4=0}^{p_4}
C_{j_4\ldots j_1}\Biggl(
\prod_{l=1}^4\zeta_{j_l}^{(i_l)}
\Biggr.
-
$$
$$
-
{\bf 1}_{\{i_1=i_2\ne 0\}}
{\bf 1}_{\{j_1=j_2\}}
\zeta_{j_3}^{(i_3)}
\zeta_{j_4}^{(i_4)}
-
{\bf 1}_{\{i_1=i_3\ne 0\}}
{\bf 1}_{\{j_1=j_3\}}
\zeta_{j_2}^{(i_2)}
\zeta_{j_4}^{(i_4)}-
$$
$$
-
{\bf 1}_{\{i_1=i_4\ne 0\}}
{\bf 1}_{\{j_1=j_4\}}
\zeta_{j_2}^{(i_2)}
\zeta_{j_3}^{(i_3)}
-
{\bf 1}_{\{i_2=i_3\ne 0\}}
{\bf 1}_{\{j_2=j_3\}}
\zeta_{j_1}^{(i_1)}
\zeta_{j_4}^{(i_4)}-
$$
$$
-
{\bf 1}_{\{i_2=i_4\ne 0\}}
{\bf 1}_{\{j_2=j_4\}}
\zeta_{j_1}^{(i_1)}
\zeta_{j_3}^{(i_3)}
-
{\bf 1}_{\{i_3=i_4\ne 0\}}
{\bf 1}_{\{j_3=j_4\}}
\zeta_{j_1}^{(i_1)}
\zeta_{j_2}^{(i_2)}+
$$
$$
+
{\bf 1}_{\{i_1=i_2\ne 0\}}
{\bf 1}_{\{j_1=j_2\}}
{\bf 1}_{\{i_3=i_4\ne 0\}}
{\bf 1}_{\{j_3=j_4\}}
+
{\bf 1}_{\{i_1=i_3\ne 0\}}
{\bf 1}_{\{j_1=j_3\}}
{\bf 1}_{\{i_2=i_4\ne 0\}}
{\bf 1}_{\{j_2=j_4\}}+
$$
\begin{equation}
\label{a4}
+\Biggl.
{\bf 1}_{\{i_1=i_4\ne 0\}}
{\bf 1}_{\{j_1=j_4\}}
{\bf 1}_{\{i_2=i_3\ne 0\}}
{\bf 1}_{\{j_2=j_3\}}\Biggr),
\end{equation}

\vspace{8mm}
$$
J[\psi^{(5)}]_{T,t}
=\hbox{\vtop{\offinterlineskip\halign{
\hfil#\hfil\cr
{\rm l.i.m.}\cr
$\stackrel{}{{}_{p_1,\ldots,p_5\to \infty}}$\cr
}} }\sum_{j_1=0}^{p_1}\ldots\sum_{j_5=0}^{p_5}
C_{j_5\ldots j_1}\Biggl(
\prod_{l=1}^5\zeta_{j_l}^{(i_l)}
-\Biggr.
$$
$$
-
{\bf 1}_{\{i_1=i_2\ne 0\}}
{\bf 1}_{\{j_1=j_2\}}
\zeta_{j_3}^{(i_3)}
\zeta_{j_4}^{(i_4)}
\zeta_{j_5}^{(i_5)}-
{\bf 1}_{\{i_1=i_3\ne 0\}}
{\bf 1}_{\{j_1=j_3\}}
\zeta_{j_2}^{(i_2)}
\zeta_{j_4}^{(i_4)}
\zeta_{j_5}^{(i_5)}-
$$
$$
-
{\bf 1}_{\{i_1=i_4\ne 0\}}
{\bf 1}_{\{j_1=j_4\}}
\zeta_{j_2}^{(i_2)}
\zeta_{j_3}^{(i_3)}
\zeta_{j_5}^{(i_5)}-
{\bf 1}_{\{i_1=i_5\ne 0\}}
{\bf 1}_{\{j_1=j_5\}}
\zeta_{j_2}^{(i_2)}
\zeta_{j_3}^{(i_3)}
\zeta_{j_4}^{(i_4)}-
$$
$$
-
{\bf 1}_{\{i_2=i_3\ne 0\}}
{\bf 1}_{\{j_2=j_3\}}
\zeta_{j_1}^{(i_1)}
\zeta_{j_4}^{(i_4)}
\zeta_{j_5}^{(i_5)}-
{\bf 1}_{\{i_2=i_4\ne 0\}}
{\bf 1}_{\{j_2=j_4\}}
\zeta_{j_1}^{(i_1)}
\zeta_{j_3}^{(i_3)}
\zeta_{j_5}^{(i_5)}-
$$
$$
-
{\bf 1}_{\{i_2=i_5\ne 0\}}
{\bf 1}_{\{j_2=j_5\}}
\zeta_{j_1}^{(i_1)}
\zeta_{j_3}^{(i_3)}
\zeta_{j_4}^{(i_4)}
-{\bf 1}_{\{i_3=i_4\ne 0\}}
{\bf 1}_{\{j_3=j_4\}}
\zeta_{j_1}^{(i_1)}
\zeta_{j_2}^{(i_2)}
\zeta_{j_5}^{(i_5)}-
$$
$$
-
{\bf 1}_{\{i_3=i_5\ne 0\}}
{\bf 1}_{\{j_3=j_5\}}
\zeta_{j_1}^{(i_1)}
\zeta_{j_2}^{(i_2)}
\zeta_{j_4}^{(i_4)}
-{\bf 1}_{\{i_4=i_5\ne 0\}}
{\bf 1}_{\{j_4=j_5\}}
\zeta_{j_1}^{(i_1)}
\zeta_{j_2}^{(i_2)}
\zeta_{j_3}^{(i_3)}+
$$
$$
+
{\bf 1}_{\{i_1=i_2\ne 0\}}
{\bf 1}_{\{j_1=j_2\}}
{\bf 1}_{\{i_3=i_4\ne 0\}}
{\bf 1}_{\{j_3=j_4\}}\zeta_{j_5}^{(i_5)}+
{\bf 1}_{\{i_1=i_2\ne 0\}}
{\bf 1}_{\{j_1=j_2\}}
{\bf 1}_{\{i_3=i_5\ne 0\}}
{\bf 1}_{\{j_3=j_5\}}\zeta_{j_4}^{(i_4)}+
$$
$$
+
{\bf 1}_{\{i_1=i_2\ne 0\}}
{\bf 1}_{\{j_1=j_2\}}
{\bf 1}_{\{i_4=i_5\ne 0\}}
{\bf 1}_{\{j_4=j_5\}}\zeta_{j_3}^{(i_3)}+
{\bf 1}_{\{i_1=i_3\ne 0\}}
{\bf 1}_{\{j_1=j_3\}}
{\bf 1}_{\{i_2=i_4\ne 0\}}
{\bf 1}_{\{j_2=j_4\}}\zeta_{j_5}^{(i_5)}+
$$
$$
+
{\bf 1}_{\{i_1=i_3\ne 0\}}
{\bf 1}_{\{j_1=j_3\}}
{\bf 1}_{\{i_2=i_5\ne 0\}}
{\bf 1}_{\{j_2=j_5\}}\zeta_{j_4}^{(i_4)}+
{\bf 1}_{\{i_1=i_3\ne 0\}}
{\bf 1}_{\{j_1=j_3\}}
{\bf 1}_{\{i_4=i_5\ne 0\}}
{\bf 1}_{\{j_4=j_5\}}\zeta_{j_2}^{(i_2)}+
$$
$$
+
{\bf 1}_{\{i_1=i_4\ne 0\}}
{\bf 1}_{\{j_1=j_4\}}
{\bf 1}_{\{i_2=i_3\ne 0\}}
{\bf 1}_{\{j_2=j_3\}}\zeta_{j_5}^{(i_5)}+
{\bf 1}_{\{i_1=i_4\ne 0\}}
{\bf 1}_{\{j_1=j_4\}}
{\bf 1}_{\{i_2=i_5\ne 0\}}
{\bf 1}_{\{j_2=j_5\}}\zeta_{j_3}^{(i_3)}+
$$
$$
+
{\bf 1}_{\{i_1=i_4\ne 0\}}
{\bf 1}_{\{j_1=j_4\}}
{\bf 1}_{\{i_3=i_5\ne 0\}}
{\bf 1}_{\{j_3=j_5\}}\zeta_{j_2}^{(i_2)}+
{\bf 1}_{\{i_1=i_5\ne 0\}}
{\bf 1}_{\{j_1=j_5\}}
{\bf 1}_{\{i_2=i_3\ne 0\}}
{\bf 1}_{\{j_2=j_3\}}\zeta_{j_4}^{(i_4)}+
$$
$$
+
{\bf 1}_{\{i_1=i_5\ne 0\}}
{\bf 1}_{\{j_1=j_5\}}
{\bf 1}_{\{i_2=i_4\ne 0\}}
{\bf 1}_{\{j_2=j_4\}}\zeta_{j_3}^{(i_3)}+
{\bf 1}_{\{i_1=i_5\ne 0\}}
{\bf 1}_{\{j_1=j_5\}}
{\bf 1}_{\{i_3=i_4\ne 0\}}
{\bf 1}_{\{j_3=j_4\}}\zeta_{j_2}^{(i_2)}+
$$
$$
+
{\bf 1}_{\{i_2=i_3\ne 0\}}
{\bf 1}_{\{j_2=j_3\}}
{\bf 1}_{\{i_4=i_5\ne 0\}}
{\bf 1}_{\{j_4=j_5\}}\zeta_{j_1}^{(i_1)}+
{\bf 1}_{\{i_2=i_4\ne 0\}}
{\bf 1}_{\{j_2=j_4\}}
{\bf 1}_{\{i_3=i_5\ne 0\}}
{\bf 1}_{\{j_3=j_5\}}\zeta_{j_1}^{(i_1)}+
$$
\begin{equation}
\label{a5}
+\Biggl.
{\bf 1}_{\{i_2=i_5\ne 0\}}
{\bf 1}_{\{j_2=j_5\}}
{\bf 1}_{\{i_3=i_4\ne 0\}}
{\bf 1}_{\{j_3=j_4\}}\zeta_{j_1}^{(i_1)}\Biggr),
\end{equation}

\vspace{8mm}

$$
J[\psi^{(6)}]_{T,t}
=\hbox{\vtop{\offinterlineskip\halign{
\hfil#\hfil\cr
{\rm l.i.m.}\cr
$\stackrel{}{{}_{p_1,\ldots,p_6\to \infty}}$\cr
}} }\sum_{j_1=0}^{p_1}\ldots\sum_{j_6=0}^{p_6}
C_{j_6\ldots j_1}\Biggl(
\prod_{l=1}^6
\zeta_{j_l}^{(i_l)}
-\Biggr.
$$
$$
-
{\bf 1}_{\{i_1=i_6\ne 0\}}
{\bf 1}_{\{j_1=j_6\}}
\zeta_{j_2}^{(i_2)}
\zeta_{j_3}^{(i_3)}
\zeta_{j_4}^{(i_4)}
\zeta_{j_5}^{(i_5)}-
{\bf 1}_{\{i_2=i_6\ne 0\}}
{\bf 1}_{\{j_2=j_6\}}
\zeta_{j_1}^{(i_1)}
\zeta_{j_3}^{(i_3)}
\zeta_{j_4}^{(i_4)}
\zeta_{j_5}^{(i_5)}-
$$
$$
-
{\bf 1}_{\{i_3=i_6\ne 0\}}
{\bf 1}_{\{j_3=j_6\}}
\zeta_{j_1}^{(i_1)}
\zeta_{j_2}^{(i_2)}
\zeta_{j_4}^{(i_4)}
\zeta_{j_5}^{(i_5)}-
{\bf 1}_{\{i_4=i_6\ne 0\}}
{\bf 1}_{\{j_4=j_6\}}
\zeta_{j_1}^{(i_1)}
\zeta_{j_2}^{(i_2)}
\zeta_{j_3}^{(i_3)}
\zeta_{j_5}^{(i_5)}-
$$
$$
-
{\bf 1}_{\{i_5=i_6\ne 0\}}
{\bf 1}_{\{j_5=j_6\}}
\zeta_{j_1}^{(i_1)}
\zeta_{j_2}^{(i_2)}
\zeta_{j_3}^{(i_3)}
\zeta_{j_4}^{(i_4)}-
{\bf 1}_{\{i_1=i_2\ne 0\}}
{\bf 1}_{\{j_1=j_2\}}
\zeta_{j_3}^{(i_3)}
\zeta_{j_4}^{(i_4)}
\zeta_{j_5}^{(i_5)}
\zeta_{j_6}^{(i_6)}-
$$
$$
-
{\bf 1}_{\{i_1=i_3\ne 0\}}
{\bf 1}_{\{j_1=j_3\}}
\zeta_{j_2}^{(i_2)}
\zeta_{j_4}^{(i_4)}
\zeta_{j_5}^{(i_5)}
\zeta_{j_6}^{(i_6)}-
{\bf 1}_{\{i_1=i_4\ne 0\}}
{\bf 1}_{\{j_1=j_4\}}
\zeta_{j_2}^{(i_2)}
\zeta_{j_3}^{(i_3)}
\zeta_{j_5}^{(i_5)}
\zeta_{j_6}^{(i_6)}-
$$
$$
-
{\bf 1}_{\{i_1=i_5\ne 0\}}
{\bf 1}_{\{j_1=j_5\}}
\zeta_{j_2}^{(i_2)}
\zeta_{j_3}^{(i_3)}
\zeta_{j_4}^{(i_4)}
\zeta_{j_6}^{(i_6)}-
{\bf 1}_{\{i_2=i_3\ne 0\}}
{\bf 1}_{\{j_2=j_3\}}
\zeta_{j_1}^{(i_1)}
\zeta_{j_4}^{(i_4)}
\zeta_{j_5}^{(i_5)}
\zeta_{j_6}^{(i_6)}-
$$
$$
-
{\bf 1}_{\{i_2=i_4\ne 0\}}
{\bf 1}_{\{j_2=j_4\}}
\zeta_{j_1}^{(i_1)}
\zeta_{j_3}^{(i_3)}
\zeta_{j_5}^{(i_5)}
\zeta_{j_6}^{(i_6)}-
{\bf 1}_{\{i_2=i_5\ne 0\}}
{\bf 1}_{\{j_2=j_5\}}
\zeta_{j_1}^{(i_1)}
\zeta_{j_3}^{(i_3)}
\zeta_{j_4}^{(i_4)}
\zeta_{j_6}^{(i_6)}-
$$
$$
-
{\bf 1}_{\{i_3=i_4\ne 0\}}
{\bf 1}_{\{j_3=j_4\}}
\zeta_{j_1}^{(i_1)}
\zeta_{j_2}^{(i_2)}
\zeta_{j_5}^{(i_5)}
\zeta_{j_6}^{(i_6)}-
{\bf 1}_{\{i_3=i_5\ne 0\}}
{\bf 1}_{\{j_3=j_5\}}
\zeta_{j_1}^{(i_1)}
\zeta_{j_2}^{(i_2)}
\zeta_{j_4}^{(i_4)}
\zeta_{j_6}^{(i_6)}-
$$
$$
-
{\bf 1}_{\{i_4=i_5\ne 0\}}
{\bf 1}_{\{j_4=j_5\}}
\zeta_{j_1}^{(i_1)}
\zeta_{j_2}^{(i_2)}
\zeta_{j_3}^{(i_3)}
\zeta_{j_6}^{(i_6)}+
$$
$$
+
{\bf 1}_{\{i_1=i_2\ne 0\}}
{\bf 1}_{\{j_1=j_2\}}
{\bf 1}_{\{i_3=i_4\ne 0\}}
{\bf 1}_{\{j_3=j_4\}}
\zeta_{j_5}^{(i_5)}
\zeta_{j_6}^{(i_6)}+
{\bf 1}_{\{i_1=i_2\ne 0\}}
{\bf 1}_{\{j_1=j_2\}}
{\bf 1}_{\{i_3=i_5\ne 0\}}
{\bf 1}_{\{j_3=j_5\}}
\zeta_{j_4}^{(i_4)}
\zeta_{j_6}^{(i_6)}+
$$
$$
+
{\bf 1}_{\{i_1=i_2\ne 0\}}
{\bf 1}_{\{j_1=j_2\}}
{\bf 1}_{\{i_4=i_5\ne 0\}}
{\bf 1}_{\{j_4=j_5\}}
\zeta_{j_3}^{(i_3)}
\zeta_{j_6}^{(i_6)}
+
{\bf 1}_{\{i_1=i_3\ne 0\}}
{\bf 1}_{\{j_1=j_3\}}
{\bf 1}_{\{i_2=i_4\ne 0\}}
{\bf 1}_{\{j_2=j_4\}}
\zeta_{j_5}^{(i_5)}
\zeta_{j_6}^{(i_6)}+
$$
$$
+
{\bf 1}_{\{i_1=i_3\ne 0\}}
{\bf 1}_{\{j_1=j_3\}}
{\bf 1}_{\{i_2=i_5\ne 0\}}
{\bf 1}_{\{j_2=j_5\}}
\zeta_{j_4}^{(i_4)}
\zeta_{j_6}^{(i_6)}
+{\bf 1}_{\{i_1=i_3\ne 0\}}
{\bf 1}_{\{j_1=j_3\}}
{\bf 1}_{\{i_4=i_5\ne 0\}}
{\bf 1}_{\{j_4=j_5\}}
\zeta_{j_2}^{(i_2)}
\zeta_{j_6}^{(i_6)}+
$$
$$
+
{\bf 1}_{\{i_1=i_4\ne 0\}}
{\bf 1}_{\{j_1=j_4\}}
{\bf 1}_{\{i_2=i_3\ne 0\}}
{\bf 1}_{\{j_2=j_3\}}
\zeta_{j_5}^{(i_5)}
\zeta_{j_6}^{(i_6)}
+
{\bf 1}_{\{i_1=i_4\ne 0\}}
{\bf 1}_{\{j_1=j_4\}}
{\bf 1}_{\{i_2=i_5\ne 0\}}
{\bf 1}_{\{j_2=j_5\}}
\zeta_{j_3}^{(i_3)}
\zeta_{j_6}^{(i_6)}+
$$
$$
+
{\bf 1}_{\{i_1=i_4\ne 0\}}
{\bf 1}_{\{j_1=j_4\}}
{\bf 1}_{\{i_3=i_5\ne 0\}}
{\bf 1}_{\{j_3=j_5\}}
\zeta_{j_2}^{(i_2)}
\zeta_{j_6}^{(i_6)}
+
{\bf 1}_{\{i_1=i_5\ne 0\}}
{\bf 1}_{\{j_1=j_5\}}
{\bf 1}_{\{i_2=i_3\ne 0\}}
{\bf 1}_{\{j_2=j_3\}}
\zeta_{j_4}^{(i_4)}
\zeta_{j_6}^{(i_6)}+
$$
$$
+
{\bf 1}_{\{i_1=i_5\ne 0\}}
{\bf 1}_{\{j_1=j_5\}}
{\bf 1}_{\{i_2=i_4\ne 0\}}
{\bf 1}_{\{j_2=j_4\}}
\zeta_{j_3}^{(i_3)}
\zeta_{j_6}^{(i_6)}
+
{\bf 1}_{\{i_1=i_5\ne 0\}}
{\bf 1}_{\{j_1=j_5\}}
{\bf 1}_{\{i_3=i_4\ne 0\}}
{\bf 1}_{\{j_3=j_4\}}
\zeta_{j_2}^{(i_2)}
\zeta_{j_6}^{(i_6)}+
$$
$$
+
{\bf 1}_{\{i_2=i_3\ne 0\}}
{\bf 1}_{\{j_2=j_3\}}
{\bf 1}_{\{i_4=i_5\ne 0\}}
{\bf 1}_{\{j_4=j_5\}}
\zeta_{j_1}^{(i_1)}
\zeta_{j_6}^{(i_6)}
+
{\bf 1}_{\{i_2=i_4\ne 0\}}
{\bf 1}_{\{j_2=j_4\}}
{\bf 1}_{\{i_3=i_5\ne 0\}}
{\bf 1}_{\{j_3=j_5\}}
\zeta_{j_1}^{(i_1)}
\zeta_{j_6}^{(i_6)}+
$$
$$
+
{\bf 1}_{\{i_2=i_5\ne 0\}}
{\bf 1}_{\{j_2=j_5\}}
{\bf 1}_{\{i_3=i_4\ne 0\}}
{\bf 1}_{\{j_3=j_4\}}
\zeta_{j_1}^{(i_1)}
\zeta_{j_6}^{(i_6)}
+
{\bf 1}_{\{i_6=i_1\ne 0\}}
{\bf 1}_{\{j_6=j_1\}}
{\bf 1}_{\{i_3=i_4\ne 0\}}
{\bf 1}_{\{j_3=j_4\}}
\zeta_{j_2}^{(i_2)}
\zeta_{j_5}^{(i_5)}+
$$
$$
+
{\bf 1}_{\{i_6=i_1\ne 0\}}
{\bf 1}_{\{j_6=j_1\}}
{\bf 1}_{\{i_3=i_5\ne 0\}}
{\bf 1}_{\{j_3=j_5\}}
\zeta_{j_2}^{(i_2)}
\zeta_{j_4}^{(i_4)}
+
{\bf 1}_{\{i_6=i_1\ne 0\}}
{\bf 1}_{\{j_6=j_1\}}
{\bf 1}_{\{i_2=i_5\ne 0\}}
{\bf 1}_{\{j_2=j_5\}}
\zeta_{j_3}^{(i_3)}
\zeta_{j_4}^{(i_4)}+
$$
$$
+
{\bf 1}_{\{i_6=i_1\ne 0\}}
{\bf 1}_{\{j_6=j_1\}}
{\bf 1}_{\{i_2=i_4\ne 0\}}
{\bf 1}_{\{j_2=j_4\}}
\zeta_{j_3}^{(i_3)}
\zeta_{j_5}^{(i_5)}
+
{\bf 1}_{\{i_6=i_1\ne 0\}}
{\bf 1}_{\{j_6=j_1\}}
{\bf 1}_{\{i_4=i_5\ne 0\}}
{\bf 1}_{\{j_4=j_5\}}
\zeta_{j_2}^{(i_2)}
\zeta_{j_3}^{(i_3)}+
$$
$$
+
{\bf 1}_{\{i_6=i_1\ne 0\}}
{\bf 1}_{\{j_6=j_1\}}
{\bf 1}_{\{i_2=i_3\ne 0\}}
{\bf 1}_{\{j_2=j_3\}}
\zeta_{j_4}^{(i_4)}
\zeta_{j_5}^{(i_5)}
+
{\bf 1}_{\{i_6=i_2\ne 0\}}
{\bf 1}_{\{j_6=j_2\}}
{\bf 1}_{\{i_3=i_5\ne 0\}}
{\bf 1}_{\{j_3=j_5\}}
\zeta_{j_1}^{(i_1)}
\zeta_{j_4}^{(i_4)}+
$$
$$
+
{\bf 1}_{\{i_6=i_2\ne 0\}}
{\bf 1}_{\{j_6=j_2\}}
{\bf 1}_{\{i_4=i_5\ne 0\}}
{\bf 1}_{\{j_4=j_5\}}
\zeta_{j_1}^{(i_1)}
\zeta_{j_3}^{(i_3)}
+
{\bf 1}_{\{i_6=i_2\ne 0\}}
{\bf 1}_{\{j_6=j_2\}}
{\bf 1}_{\{i_3=i_4\ne 0\}}
{\bf 1}_{\{j_3=j_4\}}
\zeta_{j_1}^{(i_1)}
\zeta_{j_5}^{(i_5)}+
$$
$$
+
{\bf 1}_{\{i_6=i_2\ne 0\}}
{\bf 1}_{\{j_6=j_2\}}
{\bf 1}_{\{i_1=i_5\ne 0\}}
{\bf 1}_{\{j_1=j_5\}}
\zeta_{j_3}^{(i_3)}
\zeta_{j_4}^{(i_4)}
+
{\bf 1}_{\{i_6=i_2\ne 0\}}
{\bf 1}_{\{j_6=j_2\}}
{\bf 1}_{\{i_1=i_4\ne 0\}}
{\bf 1}_{\{j_1=j_4\}}
\zeta_{j_3}^{(i_3)}
\zeta_{j_5}^{(i_5)}+
$$
$$
+
{\bf 1}_{\{i_6=i_2\ne 0\}}
{\bf 1}_{\{j_6=j_2\}}
{\bf 1}_{\{i_1=i_3\ne 0\}}
{\bf 1}_{\{j_1=j_3\}}
\zeta_{j_4}^{(i_4)}
\zeta_{j_5}^{(i_5)}
+
{\bf 1}_{\{i_6=i_3\ne 0\}}
{\bf 1}_{\{j_6=j_3\}}
{\bf 1}_{\{i_2=i_5\ne 0\}}
{\bf 1}_{\{j_2=j_5\}}
\zeta_{j_1}^{(i_1)}
\zeta_{j_4}^{(i_4)}+
$$
$$
+
{\bf 1}_{\{i_6=i_3\ne 0\}}
{\bf 1}_{\{j_6=j_3\}}
{\bf 1}_{\{i_4=i_5\ne 0\}}
{\bf 1}_{\{j_4=j_5\}}
\zeta_{j_1}^{(i_1)}
\zeta_{j_2}^{(i_2)}
+
{\bf 1}_{\{i_6=i_3\ne 0\}}
{\bf 1}_{\{j_6=j_3\}}
{\bf 1}_{\{i_2=i_4\ne 0\}}
{\bf 1}_{\{j_2=j_4\}}
\zeta_{j_1}^{(i_1)}
\zeta_{j_5}^{(i_5)}+
$$
$$
+
{\bf 1}_{\{i_6=i_3\ne 0\}}
{\bf 1}_{\{j_6=j_3\}}
{\bf 1}_{\{i_1=i_5\ne 0\}}
{\bf 1}_{\{j_1=j_5\}}
\zeta_{j_2}^{(i_2)}
\zeta_{j_4}^{(i_4)}
+
{\bf 1}_{\{i_6=i_3\ne 0\}}
{\bf 1}_{\{j_6=j_3\}}
{\bf 1}_{\{i_1=i_4\ne 0\}}
{\bf 1}_{\{j_1=j_4\}}
\zeta_{j_2}^{(i_2)}
\zeta_{j_5}^{(i_5)}+
$$
$$
+
{\bf 1}_{\{i_6=i_3\ne 0\}}
{\bf 1}_{\{j_6=j_3\}}
{\bf 1}_{\{i_1=i_2\ne 0\}}
{\bf 1}_{\{j_1=j_2\}}
\zeta_{j_4}^{(i_4)}
\zeta_{j_5}^{(i_5)}
+
{\bf 1}_{\{i_6=i_4\ne 0\}}
{\bf 1}_{\{j_6=j_4\}}
{\bf 1}_{\{i_3=i_5\ne 0\}}
{\bf 1}_{\{j_3=j_5\}}
\zeta_{j_1}^{(i_1)}
\zeta_{j_2}^{(i_2)}+
$$
$$
+
{\bf 1}_{\{i_6=i_4\ne 0\}}
{\bf 1}_{\{j_6=j_4\}}
{\bf 1}_{\{i_2=i_5\ne 0\}}
{\bf 1}_{\{j_2=j_5\}}
\zeta_{j_1}^{(i_1)}
\zeta_{j_3}^{(i_3)}
+
{\bf 1}_{\{i_6=i_4\ne 0\}}
{\bf 1}_{\{j_6=j_4\}}
{\bf 1}_{\{i_2=i_3\ne 0\}}
{\bf 1}_{\{j_2=j_3\}}
\zeta_{j_1}^{(i_1)}
\zeta_{j_5}^{(i_5)}+
$$
$$
+
{\bf 1}_{\{i_6=i_4\ne 0\}}
{\bf 1}_{\{j_6=j_4\}}
{\bf 1}_{\{i_1=i_5\ne 0\}}
{\bf 1}_{\{j_1=j_5\}}
\zeta_{j_2}^{(i_2)}
\zeta_{j_3}^{(i_3)}
+
{\bf 1}_{\{i_6=i_4\ne 0\}}
{\bf 1}_{\{j_6=j_4\}}
{\bf 1}_{\{i_1=i_3\ne 0\}}
{\bf 1}_{\{j_1=j_3\}}
\zeta_{j_2}^{(i_2)}
\zeta_{j_5}^{(i_5)}+
$$
$$
+
{\bf 1}_{\{i_6=i_4\ne 0\}}
{\bf 1}_{\{j_6=j_4\}}
{\bf 1}_{\{i_1=i_2\ne 0\}}
{\bf 1}_{\{j_1=j_2\}}
\zeta_{j_3}^{(i_3)}
\zeta_{j_5}^{(i_5)}
+
{\bf 1}_{\{i_6=i_5\ne 0\}}
{\bf 1}_{\{j_6=j_5\}}
{\bf 1}_{\{i_3=i_4\ne 0\}}
{\bf 1}_{\{j_3=j_4\}}
\zeta_{j_1}^{(i_1)}
\zeta_{j_2}^{(i_2)}+
$$
$$
+
{\bf 1}_{\{i_6=i_5\ne 0\}}
{\bf 1}_{\{j_6=j_5\}}
{\bf 1}_{\{i_2=i_4\ne 0\}}
{\bf 1}_{\{j_2=j_4\}}
\zeta_{j_1}^{(i_1)}
\zeta_{j_3}^{(i_3)}
+
{\bf 1}_{\{i_6=i_5\ne 0\}}
{\bf 1}_{\{j_6=j_5\}}
{\bf 1}_{\{i_2=i_3\ne 0\}}
{\bf 1}_{\{j_2=j_3\}}
\zeta_{j_1}^{(i_1)}
\zeta_{j_4}^{(i_4)}+
$$
$$
+
{\bf 1}_{\{i_6=i_5\ne 0\}}
{\bf 1}_{\{j_6=j_5\}}
{\bf 1}_{\{i_1=i_4\ne 0\}}
{\bf 1}_{\{j_1=j_4\}}
\zeta_{j_2}^{(i_2)}
\zeta_{j_3}^{(i_3)}
+
{\bf 1}_{\{i_6=i_5\ne 0\}}
{\bf 1}_{\{j_6=j_5\}}
{\bf 1}_{\{i_1=i_3\ne 0\}}
{\bf 1}_{\{j_1=j_3\}}
\zeta_{j_2}^{(i_2)}
\zeta_{j_4}^{(i_4)}+
$$
$$
+
{\bf 1}_{\{i_6=i_5\ne 0\}}
{\bf 1}_{\{j_6=j_5\}}
{\bf 1}_{\{i_1=i_2\ne 0\}}
{\bf 1}_{\{j_1=j_2\}}
\zeta_{j_3}^{(i_3)}
\zeta_{j_4}^{(i_4)}-
$$
$$
-
{\bf 1}_{\{i_6=i_1\ne 0\}}
{\bf 1}_{\{j_6=j_1\}}
{\bf 1}_{\{i_2=i_5\ne 0\}}
{\bf 1}_{\{j_2=j_5\}}
{\bf 1}_{\{i_3=i_4\ne 0\}}
{\bf 1}_{\{j_3=j_4\}}-
$$
$$
-
{\bf 1}_{\{i_6=i_1\ne 0\}}
{\bf 1}_{\{j_6=j_1\}}
{\bf 1}_{\{i_2=i_4\ne 0\}}
{\bf 1}_{\{j_2=j_4\}}
{\bf 1}_{\{i_3=i_5\ne 0\}}
{\bf 1}_{\{j_3=j_5\}}-
$$
$$
-
{\bf 1}_{\{i_6=i_1\ne 0\}}
{\bf 1}_{\{j_6=j_1\}}
{\bf 1}_{\{i_2=i_3\ne 0\}}
{\bf 1}_{\{j_2=j_3\}}
{\bf 1}_{\{i_4=i_5\ne 0\}}
{\bf 1}_{\{j_4=j_5\}}-
$$
$$
-
{\bf 1}_{\{i_6=i_2\ne 0\}}
{\bf 1}_{\{j_6=j_2\}}
{\bf 1}_{\{i_1=i_5\ne 0\}}
{\bf 1}_{\{j_1=j_5\}}
{\bf 1}_{\{i_3=i_4\ne 0\}}
{\bf 1}_{\{j_3=j_4\}}-
$$
$$
-
{\bf 1}_{\{i_6=i_2\ne 0\}}
{\bf 1}_{\{j_6=j_2\}}
{\bf 1}_{\{i_1=i_4\ne 0\}}
{\bf 1}_{\{j_1=j_4\}}
{\bf 1}_{\{i_3=i_5\ne 0\}}
{\bf 1}_{\{j_3=j_5\}}-
$$
$$
-
{\bf 1}_{\{i_6=i_2\ne 0\}}
{\bf 1}_{\{j_6=j_2\}}
{\bf 1}_{\{i_1=i_3\ne 0\}}
{\bf 1}_{\{j_1=j_3\}}
{\bf 1}_{\{i_4=i_5\ne 0\}}
{\bf 1}_{\{j_4=j_5\}}-
$$
$$
-
{\bf 1}_{\{i_6=i_3\ne 0\}}
{\bf 1}_{\{j_6=j_3\}}
{\bf 1}_{\{i_1=i_5\ne 0\}}
{\bf 1}_{\{j_1=j_5\}}
{\bf 1}_{\{i_2=i_4\ne 0\}}
{\bf 1}_{\{j_2=j_4\}}-
$$
$$
-
{\bf 1}_{\{i_6=i_3\ne 0\}}
{\bf 1}_{\{j_6=j_3\}}
{\bf 1}_{\{i_1=i_4\ne 0\}}
{\bf 1}_{\{j_1=j_4\}}
{\bf 1}_{\{i_2=i_5\ne 0\}}
{\bf 1}_{\{j_2=j_5\}}-
$$
$$
-
{\bf 1}_{\{i_3=i_6\ne 0\}}
{\bf 1}_{\{j_3=j_6\}}
{\bf 1}_{\{i_1=i_2\ne 0\}}
{\bf 1}_{\{j_1=j_2\}}
{\bf 1}_{\{i_4=i_5\ne 0\}}
{\bf 1}_{\{j_4=j_5\}}-
$$
$$
-
{\bf 1}_{\{i_6=i_4\ne 0\}}
{\bf 1}_{\{j_6=j_4\}}
{\bf 1}_{\{i_1=i_5\ne 0\}}
{\bf 1}_{\{j_1=j_5\}}
{\bf 1}_{\{i_2=i_3\ne 0\}}
{\bf 1}_{\{j_2=j_3\}}-
$$
$$
-
{\bf 1}_{\{i_6=i_4\ne 0\}}
{\bf 1}_{\{j_6=j_4\}}
{\bf 1}_{\{i_1=i_3\ne 0\}}
{\bf 1}_{\{j_1=j_3\}}
{\bf 1}_{\{i_2=i_5\ne 0\}}
{\bf 1}_{\{j_2=j_5\}}-
$$
$$
-
{\bf 1}_{\{i_6=i_4\ne 0\}}
{\bf 1}_{\{j_6=j_4\}}
{\bf 1}_{\{i_1=i_2\ne 0\}}
{\bf 1}_{\{j_1=j_2\}}
{\bf 1}_{\{i_3=i_5\ne 0\}}
{\bf 1}_{\{j_3=j_5\}}-
$$
$$
-
{\bf 1}_{\{i_6=i_5\ne 0\}}
{\bf 1}_{\{j_6=j_5\}}
{\bf 1}_{\{i_1=i_4\ne 0\}}
{\bf 1}_{\{j_1=j_4\}}
{\bf 1}_{\{i_2=i_3\ne 0\}}
{\bf 1}_{\{j_2=j_3\}}-
$$
$$
-
{\bf 1}_{\{i_6=i_5\ne 0\}}
{\bf 1}_{\{j_6=j_5\}}
{\bf 1}_{\{i_1=i_2\ne 0\}}
{\bf 1}_{\{j_1=j_2\}}
{\bf 1}_{\{i_3=i_4\ne 0\}}
{\bf 1}_{\{j_3=j_4\}}-
$$
\begin{equation}
\label{a6}
\Biggl.-
{\bf 1}_{\{i_6=i_5\ne 0\}}
{\bf 1}_{\{j_6=j_5\}}
{\bf 1}_{\{i_1=i_3\ne 0\}}
{\bf 1}_{\{j_1=j_3\}}
{\bf 1}_{\{i_2=i_4\ne 0\}}
{\bf 1}_{\{j_2=j_4\}}\Biggr),
\end{equation}

\vspace{6mm}
\noindent
where ${\bf 1}_A$ is the indicator of the set $A$.

It was shown that Theorem 3 is valid for convergence 
in the mean of degree 
$2n$ $(n\in \mathbb{N})$ \cite{2011-2}-\cite{2013}, \cite{arxiv-1}.
Moreover, the convergence with probability 1 (further w.~p.~1) 
is proved in Theorem 3 
for iterated Ito stochastic integrals of multiplicity $k$
for the cases of Legendre polynomials and trigonometric functions
\cite{arxiv-1}-\cite{arxiv-3},
\cite{100-000-2}, \cite{100-000-3} (also see \cite{2018a}-\cite{12aa-afterxxx}).

As it turned out, Theorem 3 remains valid for some discontinuous
complete
orthonormal systems  
of functions in the space $L_2([t,T]).$ For example, Theorem 3 is true
for the system of Haar functions 
as well as for the system of Rademacher--Walsh functions
\cite{2006}-\cite{2013}, \cite{arxiv-1}.

In \cite{2011-2}-\cite{2013}, \cite{arxiv-26b} we demonstrate that approach
to expansion of iterated Ito stochastic integrals
considered in Theorem 3 is essentially
general and allows some modifications for
other types of iterated stochastic integrals.
Versions of 
Theorem 3 for iterated stochastic integrals
with respect to martingale Poisson measures and 
for iterated stochastic integrals with respect to martingales are obtained
in \cite{2011-2}-\cite{2013}, \cite{arxiv-26b}.
The mentioned theorems are 
sufficiently natural according to
general properties of martingales.
Another modification of Theorem 3 can be 
found in \cite{2018a}-\cite{12aa-afterxxx}, \cite{arxiv-1}, \cite{arxiv-26b}, 
where complete orthonormal with weight $r(x)\ge 0$
systems of functions 
in the space $L_2([t, T])$ were considered.

A generalization of Theorem 3 (see Theorem 4 below)
for the case
of an arbitrary complete orthonormal system
of functions in the space $L_2([t,T])$ 
and $\psi_1(\tau),\ldots,\psi_k(\tau)\in L_2([t, T])$
is given in \cite{2018a} (Sect.~1.11), \cite{arxiv-1} (Sect.~15), \cite{diffjournal}, \cite{new-2023a}.
Moreover, Theorems 3 and 4 allow us to calculate exactly
the mean-square approximation error for the iterated Ito stochastic
integral (\ref{ito}) of arbitrary multiplicity $k$
(see \cite{2018}, 
\cite{2018a}-\cite{12aa-afterxxx}, \cite{arxiv-2}).
Here we consider an approxination as the expression before passing to the limit
in (\ref{tyyy}) or (\ref{leto6000}) (see below).

Application of Theorem 3 and Theorem 4 (see below) for the mean-square
approximation of iterated stochastic integrals 
with respect to the 
infinite-dimensional $Q$-Wiener process can be found
in the monographs \cite{2018a}-\cite{12aa-afterxxx} (Chapter 7) and in \cite{31a}-\cite{31ddd}.

Consider the generalization of formulas (\ref{a1})--(\ref{a6}) 
for the case 
of an arbitrary multiplicity $k$ of the stochastic integral
$J[\psi^{(k)}]_{T,t}$ as well as 
for the case
of an arbitrary complete orthonormal systems  
of functions in the space $L_2([t,T])$ 
and $\psi_1(\tau),\ldots,\psi_k(\tau)\in L_2([t, T]).$
In order to do this, let us
consider the unordered
set $\{1, 2, \ldots, k\}$ 
and separate it into two parts:
the first part consists of $r$ unordered 
pairs (sequence order of these pairs is also unimportant) and the 
second one consists of the 
remaining $k-2r$ numbers.
So, we have

\vspace{-2mm}
\begin{equation}
\label{leto5007}
(\{
\underbrace{\{g_1, g_2\}, \ldots, 
\{g_{2r-1}, g_{2r}\}}_{\small{\hbox{part 1}}}
\},
\{\underbrace{q_1, \ldots, q_{k-2r}}_{\small{\hbox{part 2}}}
\}),
\end{equation}

\vspace{3mm}
\noindent
where 
$\{g_1, g_2, \ldots, 
g_{2r-1}, g_{2r}, q_1, \ldots, q_{k-2r}\}=\{1, 2, \ldots, k\},$
braces   
mean an unordered 
set, and pa\-ren\-the\-ses mean an ordered set.

We will say that (\ref{leto5007}) is a partition 
and consider the sum with respect to all possible
partitions

\vspace{-1mm}
\begin{equation}
\label{leto5008}
\sum_{\stackrel{(\{\{g_1, g_2\}, \ldots, 
\{g_{2r-1}, g_{2r}\}\}, \{q_1, \ldots, q_{k-2r}\})}
{{}_{\{g_1, g_2, \ldots, 
g_{2r-1}, g_{2r}, q_1, \ldots, q_{k-2r}\}=\{1, 2, \ldots, k\}}}}
a_{g_1 g_2, \ldots, 
g_{2r-1} g_{2r}, q_1 \ldots q_{k-2r}},
\end{equation}

\vspace{3mm}
\noindent
where $a_{g_1 g_2, \ldots, 
g_{2r-1} g_{2r}, q_1 \ldots q_{k-2r}}\in \mathbb{R}.$

Below there are several examples of sums in the form (\ref{leto5008})

$$
\sum_{\stackrel{(\{g_1, g_2\})}{{}_{\{g_1, g_2\}=\{1, 2\}}}}
a_{g_1 g_2}=a_{12},
$$

$$
\sum_{\stackrel{(\{\{g_1, g_2\}, \{g_3, g_4\}\})}
{{}_{\{g_1, g_2, g_3, g_4\}=\{1, 2, 3, 4\}}}}
a_{g_1 g_2, g_3 g_4}=a_{12,34} + a_{13,24} + a_{23,14},
$$

\vspace{5mm}
$$
\sum_{\stackrel{(\{g_1, g_2\}, \{q_1, q_{2}\})}
{{}_{\{g_1, g_2, q_1, q_{2}\}=\{1, 2, 3, 4\}}}}
a_{g_1 g_2, q_1 q_{2}}=
$$

$$
=a_{12,34}+a_{13,24}+a_{14,23}
+a_{23,14}+a_{24,13}+a_{34,12},
$$

\vspace{6mm}
$$
\sum_{\stackrel{(\{g_1, g_2\}, \{q_1, q_{2}, q_3\})}
{{}_{\{g_1, g_2, q_1, q_{2}, q_3\}=\{1, 2, 3, 4, 5\}}}}
a_{g_1 g_2, q_1 q_{2}q_3}
=
$$

$$
=a_{12,345}+a_{13,245}+a_{14,235}
+a_{15,234}+a_{23,145}+a_{24,135}+
$$

$$
+a_{25,134}+a_{34,125}+a_{35,124}+a_{45,123},
$$

\vspace{7mm}
$$
\sum_{\stackrel{(\{\{g_1, g_2\}, \{g_3, g_{4}\}\}, \{q_1\})}
{{}_{\{g_1, g_2, g_3, g_{4}, q_1\}=\{1, 2, 3, 4, 5\}}}}
a_{g_1 g_2, g_3 g_{4},q_1}
=
$$

$$
=
a_{12,34,5}+a_{13,24,5}+a_{14,23,5}+
a_{12,35,4}+a_{13,25,4}+
$$

$$
+a_{15,23,4}+
a_{12,54,3}+a_{15,24,3}+
a_{14,25,3}+a_{15,34,2}+a_{13,54,2}+
$$

$$
+
a_{14,53,2}+
a_{52,34,1}+a_{53,24,1}+a_{54,23,1}.
$$

\vspace{9mm}

Now we can generalize Theorem 3.

\vspace{2mm}

{\bf Theorem 4}\ \cite{2018a} (Sect.~1.11), \cite{arxiv-1} (Sect.~15), \cite{diffjournal}, \cite{new-2023a}.
{\it Suppose that
$\psi_1(\tau),\ldots,\psi_k(\tau)\in L_2([t, T])$ and
$\{\phi_j(x)\}_{j=0}^{\infty}$ is an arbitrary complete orthonormal system  
of functions in the space $L_2([t,T]).$
Then the following expansion

\vspace{2mm}
$$
J[\psi^{(k)}]_{T,t}=
\hbox{\vtop{\offinterlineskip\halign{
\hfil#\hfil\cr
{\rm l.i.m.}\cr
$\stackrel{}{{}_{p_1,\ldots,p_k\to \infty}}$\cr
}} }
\sum\limits_{j_1=0}^{p_1}\ldots
\sum\limits_{j_k=0}^{p_k}
C_{j_k\ldots j_1}\Biggl(
\prod_{l=1}^k\zeta_{j_l}^{(i_l)}+\sum\limits_{r=1}^{[k/2]}
(-1)^r \times \Biggr.
$$

\vspace{3mm}
\begin{equation}
\label{leto6000}
\times
\sum_{\stackrel{(\{\{g_1, g_2\}, \ldots, 
\{g_{2r-1}, g_{2r}\}\}, \{q_1, \ldots, q_{k-2r}\})}
{{}_{\{g_1, g_2, \ldots, 
g_{2r-1}, g_{2r}, q_1, \ldots, q_{k-2r}\}=\{1,2, \ldots, k\}}}}
\prod\limits_{s=1}^r
{\bf 1}_{\{i_{g_{{}_{2s-1}}}=~i_{g_{{}_{2s}}}\ne 0\}}
\Biggl.{\bf 1}_{\{j_{g_{{}_{2s-1}}}=~j_{g_{{}_{2s}}}\}}
\prod_{l=1}^{k-2r}\zeta_{j_{q_l}}^{(i_{q_l})}\Biggr)
\end{equation}

\vspace{6mm}
\noindent 
that converges in the mean-square sense is valid$,$ where 
$[x]$ is an integer part of a real number $x$ and 
$\prod\limits_{\emptyset}
\stackrel{\sf def}{=}1,$ $\sum\limits_{\emptyset}
\stackrel{\sf def}{=}0;$
another notations are the same as in Theorem~{\rm 3.}
}

\vspace{4mm}

In particular from (\ref{leto6000}) for $k=5$ we obtain

\vspace{2mm}

$$
J[\psi^{(5)}]_{T,t}=
\hbox{\vtop{\offinterlineskip\halign{
\hfil#\hfil\cr
{\rm l.i.m.}\cr
$\stackrel{}{{}_{p_1,\ldots,p_5\to \infty}}$\cr
}} }
\sum\limits_{j_1=0}^{p_1}\ldots
\sum\limits_{j_5=0}^{p_5}
C_{j_5 \ldots j_1}
\Biggl(
\zeta_{j_1}^{(i_1)}\ldots \zeta_{j_5}^{(i_5)}
-\Biggr.
$$

\vspace{2mm}
$$
-\sum\limits_{\stackrel{(\{g_1, g_2\}, \{q_1, q_{2}, q_3\})}
{{}_{\{g_1, g_2, q_{1}, q_{2}, q_3\}=\{1, 2, 3, 4, 5\}}}}
{\bf 1}_{\{i_{g_{{}_{1}}}=~i_{g_{{}_{2}}}\ne 0\}}
{\bf 1}_{\{j_{g_{{}_{1}}}=~j_{g_{{}_{2}}}\}}
\prod_{l=1}^{3}\zeta_{j_{q_l}}^{(i_{q_l})}+
$$

\vspace{2mm}
$$
\Biggl.+
\sum_{\stackrel{(\{\{g_1, g_2\}, 
\{g_{3}, g_{4}\}\}, \{q_1\})}
{{}_{\{g_1, g_2, g_{3}, g_{4}, q_1\}=\{1, 2, 3, 4, 5\}}}}
{\bf 1}_{\{i_{g_{{}_{1}}}=~i_{g_{{}_{2}}}\ne 0\}}
{\bf 1}_{\{j_{g_{{}_{1}}}=~j_{g_{{}_{2}}}\}}
\Biggl.{\bf 1}_{\{i_{g_{{}_{3}}}=~i_{g_{{}_{4}}}\ne 0\}}
{\bf 1}_{\{j_{g_{{}_{3}}}=~j_{g_{{}_{4}}}\}}
\zeta_{j_{q_1}}^{(i_{q_1})}\Biggr).
$$

\vspace{6mm}

\noindent
The last equality obviously agrees with
(\ref{a5}).

It should be noted that an analogue of Theorem 4 for multiple Ito 
stochastic integrals was considered 
in \cite{Rybakov1000}. 
Note that we use another notations in comparison with \cite{Rybakov1000}.
Moreover, the proof of an analogue of Theorem 4
from \cite{Rybakov1000} is different from the proof given in 
\cite{2018a} (Sect.~1.11), \cite{arxiv-1} (Sect.~15), \cite{diffjournal}, \cite{new-2023a}.

\vspace{5mm}

\section{The Idea of the Proof of Hypotheses 1, 2, and 3}

\vspace{5mm}

Let us consider the idea of the proof of Hypotheses 1--3.
Introduce the following notations

\vspace{1mm}
$$
J[\psi^{(k)}]_{T,t}^{s_l,\ldots,s_1}\ \stackrel{\rm def}{=}\ 
\prod_{q=1}^l {\bf 1}_{\{i_{s_q}=
i_{s_{q}+1}\ne 0\}}\ \times
$$

\vspace{1mm}
$$
\times
\int\limits_t^T\psi_k(t_k)\ldots \int\limits_t^{
t_{s_l+3}}\psi_{s_l+2}(t_{s_l+2})
\int\limits_t^{t_{s_l+2}}\psi_{s_l}(t_{s_l+1})
\psi_{s_l+1}(t_{s_l+1})\times
$$

\vspace{1mm}
$$
\times
\int\limits_t^{t_{s_l+1}}\psi_{s_l-1}(t_{s_l-1})
\ldots
\int\limits_t^{t_{s_1+3}}\psi_{s_1+2}(t_{s_1+2})
\int\limits_t^{t_{s_1+2}}\psi_{s_1}(t_{s_1+1})
\psi_{s_1+1}(t_{s_1+1})\times
$$

\vspace{1mm}
$$
\times
\int\limits_t^{t_{s_1+1}}\psi_{s_1-1}(t_{s_1-1})
\ldots \int\limits_t^{t_2}\psi_1(t_1)
d{\bf w}_{t_1}^{(i_1)}\ldots d{\bf w}_{t_{s_1-1}}^{(i_{s_1-1})}
dt_{s_1+1}d{\bf w}_{t_{s_1+2}}^{(i_{s_1+2})}\ldots
$$

\vspace{1mm}
\begin{equation}
\label{30.1}
\ldots
d{\bf w}_{t_{s_l-1}}^{(i_{s_l-1})}
dt_{s_l+1}d{\bf w}_{t_{s_l+2}}^{(i_{s_l+2})}\ldots d{\bf w}_{t_k}^{(i_k)},
\end{equation}

\vspace{4mm}
\noindent
where $(s_l,\ldots,s_1)\in{\rm A}_{k,l},$

\begin{equation}
\label{30.5550001}
{\rm A}_{k,l}=\left\{(s_l,\ldots,s_1):\
s_l>s_{l-1}+1,\ldots,s_2>s_1+1;\ s_l,\ldots,s_1=1,\ldots,k-1\right\},
\end{equation}

\vspace{4mm}
\noindent
$l=1,2,\ldots,\left[k/2\right],$\ \
$i_s=0, 1,\ldots,m,$\ \
$s=1,\ldots,k,$\ \ $[x]$ is an
integer part of a real number $x,$\ \
${\bf 1}_A$ is the indicator of the set $A$.

Let us formulate the statement on connection 
between
iterated 
Ito and Stratonovich stochastic integrals 
(\ref{ito}) and (\ref{str})
of arbitrary multiplicity $k$.

\vspace{2mm}

{\bf Theorem 5} \cite{300a} (1997) (also see \cite{2006}-\cite{12aa-afterxxx}). 
{\it Suppose that
every $\psi_l(\tau)$ $(l=1,\ldots,k)$ is a continuous
nonrandom function at the interval $[t, T]$.
Then, the following relation between iterated
Ito and Stra\-to\-no\-vich stochastic integrals {\rm (\ref{ito})}
and {\rm (\ref{str})} is correct

\begin{equation}
\label{30.4}
J^{*}[\psi^{(k)}]_{T,t}=J[\psi^{(k)}]_{T,t}+
\sum_{r=1}^{\left[k/2\right]}\frac{1}{2^r}
\sum_{(s_r,\ldots,s_1)\in {\rm A}_{k,r}}
J[\psi^{(k)}]_{T,t}^{s_r,\ldots,s_1}\ \ \ \hbox{{\rm w.\ p.\ 1}},
\end{equation}

\vspace{4mm}
\noindent
where $\sum\limits_{\emptyset}$ is supposed to be equal to zero.}

\vspace{2mm}

Note that the condition of continuity of the functions
$\psi_1(\tau),\ldots, \psi_k(\tau)$ 
is related to the definition \cite{KlPl2} 
of the Stratonovich stochastic integral that we use.

According to (\ref{tyyy}), we have

\vspace{1mm}
\begin{equation}
\label{tyyy1}
\hbox{\vtop{\offinterlineskip\halign{
\hfil#\hfil\cr
{\rm l.i.m.}\cr
$\stackrel{}{{}_{p_1,\ldots,p_k\to \infty}}$\cr
}} }\sum_{j_1=0}^{p_1}\ldots\sum_{j_k=0}^{p_k}
C_{j_k\ldots j_1}\prod_{g=1}^k\zeta_{j_g}^{(i_g)}=
J[\psi^{(k)}]_{T,t}+
$$

\vspace{1mm}
$$
+
\hbox{\vtop{\offinterlineskip\halign{
\hfil#\hfil\cr
{\rm l.i.m.}\cr
$\stackrel{}{{}_{p_1,\ldots,p_k\to \infty}}$\cr
}} }\sum_{j_1=0}^{p_1}\ldots\sum_{j_k=0}^{p_k}
C_{j_k\ldots j_1}
~\hbox{\vtop{\offinterlineskip\halign{
\hfil#\hfil\cr
{\rm l.i.m.}\cr
$\stackrel{}{{}_{N\to \infty}}$\cr
}} }\sum_{(l_1,\ldots,l_k)\in {G}_k}
\prod_{g=1}^k
\phi_{j_{g}}(\tau_{l_g})
\Delta{\bf w}_{\tau_{l_g}}^{(i_g)}.
\end{equation}

\vspace{6mm}

From (\ref{tyyy1}) and (\ref{30.4}) it follows that

\begin{equation}
\label{tt1}
J^{*}[\psi^{(k)}]_{T,t}=\hbox{\vtop{\offinterlineskip\halign{
\hfil#\hfil\cr
{\rm l.i.m.}\cr
$\stackrel{}{{}_{p_1,\ldots,p_k\to \infty}}$\cr
}} }\sum_{j_1=0}^{p_1}\ldots\sum_{j_k=0}^{p_k}
C_{j_k\ldots j_1}\prod_{g=1}^k\zeta_{j_g}^{(i_g)}
\end{equation}

\vspace{3mm}
\noindent
if 

\vspace{-2mm}
$$
\sum_{r=1}^{\left[k/2\right]}\frac{1}{2^r}
\sum_{(s_r,\ldots,s_1)\in {\rm A}_{k,r}}
J[\psi^{(k)}]_{T,t}^{s_r,\ldots,s_1}=
$$

\vspace{3mm}
$$
=
\hbox{\vtop{\offinterlineskip\halign{
\hfil#\hfil\cr
{\rm l.i.m.}\cr
$\stackrel{}{{}_{p_1,\ldots,p_k\to \infty}}$\cr
}} }\sum_{j_1=0}^{p_1}\ldots\sum_{j_k=0}^{p_k}
C_{j_k\ldots j_1}
~\hbox{\vtop{\offinterlineskip\halign{
\hfil#\hfil\cr
{\rm l.i.m.}\cr
$\stackrel{}{{}_{N\to \infty}}$\cr
}} }\sum_{(l_1,\ldots,l_k)\in {G}_k}
\prod_{g=1}^k
\phi_{j_{g}}(\tau_{l_g})
\Delta{\bf w}_{\tau_{l_g}}^{(i_g)}\ \ \ \hbox{w.\ p.\ 1.}
$$

\vspace{6mm}

In the case $p_1=\ldots=p_k=p$ and 
$\psi_l(\tau)\equiv 1$ $(l=1,\ldots,k)$
from (\ref{tt1}) we obtain the statement of Hypothesis 1 
(see (\ref{feto1900otita})).

If 
$p_1=\ldots=p_k=p$ and every $\psi_l(\tau)$ $(l=1,\ldots,k)$
is 
an enough smooth nonrandom function
on $[t,T],$ then
from (\ref{tt1}) we obtain the statement of Hypothesis 2 
(see (\ref{feto1900otitarrr})).

In the case when every $\psi_l(\tau)$ $(l=1,\ldots,k)$
is 
an enough smooth nonrandom function
on $[t,T]$ 
from (\ref{tt1}) we obtain the statement of Hypothesis 3 
(see (\ref{feto1900otitarrri})).

In the following sections we consider some theorems proving 
Hypothesis 1 for 
$k=1,\ldots,6,$ Hypothesis 2 for $k=1,\ldots,5$ and Hypothesis 3 for $k=1,\ldots,3.$
Moreover, the proof of Hypotheses 2 and 3
is given for an arbitrary $k$ under the condition of convergence of trace series.
The case $k=1$ obviously directly follows 
from Theorem 3 (see (\ref{a1})).

\vspace{5mm}

\section{Expansions of Iterated Stratonovich 
Stochastc Integrals of Multiplicity 2--4. Some Old Results}

\vspace{5mm}

As it turned out, approximations of the iterated Stratonovich 
stochastic integrals (\ref{str}) 
are essentially simpler
than the appropriate approximations of the iterated Ito
stochastic integrals (\ref{ito}) 
based on Theorems 3 and 4.
For the first time 
this fact was mentioned
in \cite{2006} (2006). 

According to the standard connection between Ito
and Stratonovich stochastic integrals, 
the iterated Ito and Stratonovich stochastic integrals 
(\ref{ito}) and (\ref{str}) of first
multiplicity are equal to each other w.~p.~1.
So, we begin the consideration from the multiplicity $k=2$.

The following theorems adapt Theorems 3, 4 for the integrals
(\ref{str}) of multiplicities 2--4.

\vspace{2mm}

{\bf Theorem 6} \cite{2011-2}-\cite{12aa-afterxxx}, \cite{2010-2}-\cite{2013},
\cite{arxiv-5}.
{\it Suppose that the following conditions are fulfilled{\rm :}

{\rm 1}. The function $\psi_2(\tau)$ is continuously 
differentiable at the interval $[t, T]$ and the
function $\psi_1(\tau)$ is twice continuously 
differentiable at the interval $[t, T]$.

{\rm 2}. $\{\phi_j(x)\}_{j=0}^{\infty}$ is a complete orthonormal system 
of Legendre polynomials or trigonometric functions
in the space $L_2([t, T]).$

Then, the iterated Stratonovich stochastic integral of 
second multiplicity

$$
J^{*}[\psi^{(2)}]_{T,t}={\int\limits_t^{*}}^T\psi_2(t_2)
{\int\limits_t^{*}}^{t_2}\psi_1(t_1)d{\bf w}_{t_1}^{(i_1)}
d{\bf w}_{t_2}^{(i_2)}\ \ \ (i_1, i_2=0, 1,\ldots,m)
$$

\vspace{4mm}
\noindent
is expanded into the 
converging in the mean-square sense double series

$$
J^{*}[\psi^{(2)}]_{T,t}=
\hbox{\vtop{\offinterlineskip\halign{
\hfil#\hfil\cr
{\rm l.i.m.}\cr
$\stackrel{}{{}_{p_1,p_2\to \infty}}$\cr
}} }\sum_{j_1=0}^{p_1}\sum_{j_2=0}^{p_2}
C_{j_2j_1}\zeta_{j_1}^{(i_1)}\zeta_{j_2}^{(i_2)},
$$

\vspace{4mm}
\noindent
where notations are the same as in 
Theorems {\rm 3, 4}.}

\vspace{2mm}

Proving Theorem 6 \cite{2011-2}-\cite{12aa-afterxxx}, \cite{2010-2}-\cite{2013},
\cite{arxiv-5},
we used Theorem 3 and double integration by parts. This
procedure leads to the condition of double 
continuously 
differentiability of the
function $\psi_1(\tau)$ at the interval $[t, T]$.
The mentioned condition can be weakened, but the 
proof becomes more complicated. As a result, we
have the following theorem.

\vspace{2mm}

{\bf Theorem 7} \cite{2018a}-\cite{12aa-afterxxx}, \cite{30a}, \cite{arxiv-8}.\
{\it Suppose that the following conditions are fulfilled{\rm :}

{\rm 1}. Every $\psi_l(\tau)$ $(l=1, 2)$ is a continuously 
differentiable func\-ti\-on at the interval $[t, T]$.

{\rm 2}. $\{\phi_j(x)\}_{j=0}^{\infty}$ is a complete orthonormal system 
of Legendre polynomials or trigonometric functions
in the space $L_2([t, T]).$

Then, the iterated 
Stratonovich stochastic integral of second multiplicity

$$
J^{*}[\psi^{(2)}]_{T,t}={\int\limits_t^{*}}^T\psi_2(t_2)
{\int\limits_t^{*}}^{t_2}\psi_1(t_1)d{\bf w}_{t_1}^{(i_1)}
d{\bf w}_{t_2}^{(i_2)}\ \ \ (i_1, i_2=0, 1,\ldots,m)
$$

\vspace{4mm}
\noindent
is expanded into the 
converging in the mean-square sense double series

$$
J^{*}[\psi^{(2)}]_{T,t}=
\hbox{\vtop{\offinterlineskip\halign{
\hfil#\hfil\cr
{\rm l.i.m.}\cr
$\stackrel{}{{}_{p_1,p_2\to \infty}}$\cr
}} }\sum_{j_1=0}^{p_1}\sum_{j_2=0}^{p_2}
C_{j_2j_1}\zeta_{j_1}^{(i_1)}\zeta_{j_2}^{(i_2)},
$$

\vspace{4mm}
\noindent
where notations are the same as in
Theorems {\rm 3, 4}.}

\vspace{2mm}

The proof of Theorem 7 \cite{2018a}-\cite{12aa-afterxxx}, \cite{30a}, \cite{arxiv-8} is based on 
Theorem 3 and double 
Fourier--Legendre series as well as double trigonometric Fourier series
summarized by Pringsheim method
at the square $[t, T]^2.$ 

Recently, Theorem~7 has been generalized to the case 
of an arbitrary complete orthonormal system
of functions in the space $L_2([t,T])$ 
and $\psi_1(\tau),\psi_2(\tau)\in L_2([t, T])$
(see \cite{2018a}, Sect.~2.18) or Theorem~42 below.

The following 4 theorems (Theorems 8--11) adapt Theorems 3, 4 for the integrals
(\ref{str}) of multiplicity 3.

\vspace{2mm}

{\bf Theorem 8}\ \cite{2011-2}-\cite{12aa-afterxxx}, \cite{2010-2}-\cite{2013},
\cite{arxiv-7}.\
{\it Suppose that
$\{\phi_j(x)\}_{j=0}^{\infty}$ is a complete orthonormal
system of Legendre polynomials or trigonometric func\-ti\-ons
in the space $L_2([t, T])$.
Then, for the iterated Stratonovich stochastic integral of 
third multiplicity

$$
{\int\limits_t^{*}}^T
{\int\limits_t^{*}}^{t_3}
{\int\limits_t^{*}}^{t_2}
d{\bf f}_{t_1}^{(i_1)}
d{\bf f}_{t_2}^{(i_2)}d{\bf f}_{t_3}^{(i_3)}\ \ \ (i_1, i_2, i_3=1,\ldots,m)
$$

\vspace{4mm}
\noindent
the following 
expansion 

\vspace{-1mm}
\begin{equation}
\label{feto19001}
{\int\limits_t^{*}}^T
{\int\limits_t^{*}}^{t_3}
{\int\limits_t^{*}}^{t_2}
d{\bf f}_{t_1}^{(i_1)}
d{\bf f}_{t_2}^{(i_2)}d{\bf f}_{t_3}^{(i_3)}\ 
=
\hbox{\vtop{\offinterlineskip\halign{
\hfil#\hfil\cr
{\rm l.i.m.}\cr
$\stackrel{}{{}_{p_1,p_2,p_3\to \infty}}$\cr
}} }\sum_{j_1=0}^{p_1}\sum_{j_2=0}^{p_2}\sum_{j_3=0}^{p_3}
C_{j_3 j_2 j_1}\zeta_{j_1}^{(i_1)}\zeta_{j_2}^{(i_2)}\zeta_{j_3}^{(i_3)}
\stackrel{\rm def}{=}
$$
$$
\stackrel{\rm def}{=}
\sum\limits_{j_1, j_2, j_3=0}^{\infty}
C_{j_3 j_2 j_1}\zeta_{j_1}^{(i_1)}\zeta_{j_2}^{(i_2)}\zeta_{j_3}^{(i_3)}
\end{equation}

\vspace{3mm}
\noindent
that converges in the mean-square sense is valid, where

$$
C_{j_3 j_2 j_1}=\int\limits_t^T
\phi_{j_3}(t_3)\int\limits_t^{t_3}
\phi_{j_2}(t_2)
\int\limits_t^{t_2}
\phi_{j_1}(t_1)dt_1dt_2dt_3;
$$

\vspace{4mm}
\noindent
another notations are the same as in Theorems {\rm 3, 4}.}

\vspace{2mm}

Obviously, that Theorem 8 proves Hypothesis 3 
for the case $k=3$ and $\psi_l(\tau)\equiv 1$ ($l=1, 2, 3$).

Let us consider the generalization of Theorem 8 (the case of Legendre
polynomials) for the binomial weight functions.

\vspace{2mm}

{\bf Theorem 9}\ \cite{2011-2}-\cite{12aa-afterxxx}, \cite{2010-2}-\cite{2013},
\cite{arxiv-7}.\ {\it Suppose that
$\{\phi_j(x)\}_{j=0}^{\infty}$ is a complete orthonormal
system of Legendre polynomials
in the space $L_2([t, T])$.
Then, for the iterated Stratonovich stochastic integral of 
third multiplicity

$$
I_{{l_1l_2l_3}_{T,t}}^{*(i_1i_2i_3)}={\int\limits_t^{*}}^T(t-t_3)^{l_3}
{\int\limits_t^{*}}^{t_3}(t-t_2)^{l_2}
{\int\limits_t^{*}}^{t_2}(t-t_1)^{l_1}
d{\bf f}_{t_1}^{(i_1)}
d{\bf f}_{t_2}^{(i_2)}d{\bf f}_{t_3}^{(i_3)}\ \ \ (i_1, i_2, i_3=1,\ldots,m)
$$

\vspace{4mm}
\noindent
the following 
expansion 

\vspace{1mm}
$$
I_{{l_1l_2l_3}_{T,t}}^{*(i_1i_2i_3)}=
\hbox{\vtop{\offinterlineskip\halign{
\hfil#\hfil\cr
{\rm l.i.m.}\cr
$\stackrel{}{{}_{p_1,p_2,p_3\to \infty}}$\cr
}} }\sum_{j_1=0}^{p_1}\sum_{j_2=0}^{p_2}\sum_{j_3=0}^{p_3}
C_{j_3 j_2 j_1}\zeta_{j_1}^{(i_1)}\zeta_{j_2}^{(i_2)}\zeta_{j_3}^{(i_3)}
\stackrel{\rm def}{=}
$$

\vspace{2mm}
\begin{equation}
\label{feto1900}
\stackrel{\rm def}{=}
\sum\limits_{j_1, j_2, j_3=0}^{\infty}
C_{j_3 j_2 j_1}\zeta_{j_1}^{(i_1)}\zeta_{j_2}^{(i_2)}\zeta_{j_3}^{(i_3)}
\end{equation}

\vspace{4mm}
\noindent
that converges in the mean-square sense 
is valid for each of the following cases

\vspace{2mm}
\noindent
{\rm 1}.\ $i_1\ne i_2,\ i_2\ne i_3,\ i_1\ne i_3$\ and
$l_1, l_2, l_3=0, 1, 2,\ldots $\\
{\rm 2}.\ $i_1=i_2\ne i_3$ and $l_1=l_2\ne l_3$\ and
$l_1, l_2, l_3=0, 1, 2,\ldots $\\
{\rm 3}.\ $i_1\ne i_2=i_3$ and $l_1\ne l_2=l_3$\ and
$l_1, l_2, l_3=0, 1, 2,\ldots $\\
{\rm 4}.\ $i_1, i_2, i_3=1,\ldots,m;$ $l_1=l_2=l_3=l$\ and $l=0, 1, 
2,\ldots,$\\

\noindent
where

\vspace{-3mm}
$$
C_{j_3 j_2 j_1}=\int\limits_t^T(t-t_3)^{l_3}\phi_{j_3}(t_3)
\int\limits_t^{t_3}(t-t_2)^{l_2}\phi_{j_2}(t_2)
\int\limits_t^{t_2}(t-t_1)^{l_1}\phi_{j_1}(t_1)dt_1dt_2dt_3;
$$

\vspace{2mm}
\noindent
another notations are the same as in Theorems {\rm 3, 4}.}

\vspace{2mm}

We can introduce the weight functions $\psi_l(\tau)$ $(l=1, 2, 3)$
with some properties of smoothness. However, we consider in this case
the more specific method of series summation ($p_1=p_2=p_3=p\to\infty$).

\vspace{2mm}

{\bf Theorem 10}\ \cite{2011-2}-\cite{12aa-afterxxx}, \cite{2010-2}-\cite{2013}.\
{\it Suppose that
$\{\phi_j(x)\}_{j=0}^{\infty}$ is a complete orthonormal
system of Legendre polynomials or trigonometric func\-ti\-ons
in the space $L_2([t, T])$
and $\psi_l(\tau)$ $(l=1, 2, 3)$ are continuously
differentiable functions at the interval $[t, T]$.
Then, for the iterated Stratonovich stochastic 
integral of third multiplicity

$$
J^{*}[\psi^{(3)}]_{T,t}={\int\limits_t^{*}}^T\psi_3(t_3)
{\int\limits_t^{*}}^{t_3}\psi_2(t_2)
{\int\limits_t^{*}}^{t_2}\psi_1(t_1)
d{\bf f}_{t_1}^{(i_1)}
d{\bf f}_{t_2}^{(i_2)}d{\bf f}_{t_3}^{(i_3)}\ \ \ (i_1, i_2, i_3=1,\ldots,m)
$$

\vspace{4mm}
\noindent
the following 
expansion 

\begin{equation}
\label{feto19000}
J^{*}[\psi^{(3)}]_{T,t}=
\hbox{\vtop{\offinterlineskip\halign{
\hfil#\hfil\cr
{\rm l.i.m.}\cr
$\stackrel{}{{}_{p\to \infty}}$\cr
}} }
\sum\limits_{j_1, j_2, j_3=0}^{p}
C_{j_3 j_2 j_1}\zeta_{j_1}^{(i_1)}\zeta_{j_2}^{(i_2)}\zeta_{j_3}^{(i_3)}
\end{equation}

\vspace{4mm}
\noindent
that converges in the mean-square sense is valid 
for each of the following cases

\vspace{2mm}
\noindent
{\rm 1}.\ $i_1\ne i_2,\ i_2\ne i_3,\ i_1\ne i_3,$\\
{\rm 2}.\ $i_1=i_2\ne i_3$ and
$\psi_1(\tau)\equiv\psi_2(\tau),$\\
{\rm 3}.\ $i_1\ne i_2=i_3$ and
$\psi_2(\tau)\equiv\psi_3(\tau),$\\
{\rm 4}.\ $i_1, i_2, i_3=1,\ldots,m$
and $\psi_1(\tau)\equiv\psi_2(\tau)\equiv\psi_3(\tau),$\

\vspace{2mm}
\noindent
where

\vspace{-1mm}
$$
C_{j_3 j_2 j_1}=\int\limits_t^T\psi_3(t_3)\phi_{j_3}(t_3)
\int\limits_t^{t_3}\psi_2(t_2)\phi_{j_2}(t_2)
\int\limits_t^{t_2}\psi_1(t_1)\phi_{j_1}(t_1)dt_1dt_2dt_3;
$$

\vspace{4mm}
\noindent
another notations are the same as in Theorems {\rm 3, 4}.}

\vspace{2mm}

We can omit Cases 1--4 in Theorem 10 in the case
when the functions 
$\psi_1(\tau)$ and $\psi_3(\tau)$
are
twice
continuously differentiable.

\vspace{2mm}

{\bf Theorem 11}\ \cite{2017}-\cite{12aa-afterxxx}, \cite{2013},
\cite{arxiv-5}.\
{\it Suppose that
$\{\phi_j(x)\}_{j=0}^{\infty}$ is a complete orthonormal
system of Legendre polynomials or trigonomertic func\-ti\-ons
in the space $L_2([t, T])$. Furthermore, let
the function $\psi_2(\tau)$ is continuously
differentiable at the interval $[t, T]$ and
the functions $\psi_1(\tau),$ $\psi_3(\tau)$ are twice continuously
differentiable at the interval $[t, T]$.
Then, for the iterated Stra\-to\-no\-vich stochastic integral 
of third multiplicity

$$
J^{*}[\psi^{(3)}]_{T,t}={\int\limits_t^{*}}^T\psi_3(t_3)
{\int\limits_t^{*}}^{t_3}\psi_2(t_2)
{\int\limits_t^{*}}^{t_2}\psi_1(t_1)
d{\bf f}_{t_1}^{(i_1)}
d{\bf f}_{t_2}^{(i_2)}d{\bf f}_{t_3}^{(i_3)}\ \ \ (i_1, i_2, i_3=1,\ldots,m)
$$

\vspace{4mm}
\noindent
the following 
expansion

\begin{equation}
\label{feto19000a}
J^{*}[\psi^{(3)}]_{T,t}=
\hbox{\vtop{\offinterlineskip\halign{
\hfil#\hfil\cr
{\rm l.i.m.}\cr
$\stackrel{}{{}_{p\to \infty}}$\cr
}} }
\sum\limits_{j_1, j_2, j_3=0}^{p}
C_{j_3 j_2 j_1}\zeta_{j_1}^{(i_1)}\zeta_{j_2}^{(i_2)}\zeta_{j_3}^{(i_3)}
\end{equation}

\vspace{4mm}
\noindent
that converges in the mean-square sense is valid, where

\vspace{-1mm}
$$
C_{j_3 j_2 j_1}=\int\limits_t^T\psi_3(t_3)\phi_{j_3}(t_3)
\int\limits_t^{t_3}\psi_2(t_2)\phi_{j_2}(t_2)
\int\limits_t^{t_2}\psi_1(t_1)\phi_{j_1}(t_1)dt_1dt_2dt_3;
$$

\vspace{4mm}
\noindent
another notations are the same as in Theorems {\rm 3, 4}.}

The following theorem adapts Theorems 3, 4 for the integrals
(\ref{str}) ($\psi_l(\tau)\equiv 1,\ l=1,\ldots,4$) 
of multiplicity 4.

\vspace{2mm}

{\bf Theorem 12} \cite{2017}-\cite{12aa-afterxxx},
\cite{2013}, \cite{arxiv-5}.\ 
{\it Suppose that
$\{\phi_j(x)\}_{j=0}^{\infty}$ is a complete orthonormal
system of Legendre polynomials or trigonometric func\-ti\-ons
in the space $L_2([t, T])$.
Then, for the iterated Stra\-to\-no\-vich stochastic integral 
of fourth multiplicity

$$
I_{T,t}^{*(i_1 i_2 i_3 i_4)}=
{\int\limits_t^{*}}^T
{\int\limits_t^{*}}^{t_4}
{\int\limits_t^{*}}^{t_3}
{\int\limits_t^{*}}^{t_2}
d{\bf w}_{t_1}^{(i_1)}
d{\bf w}_{t_2}^{(i_2)}d{\bf w}_{t_3}^{(i_3)}d{\bf w}_{t_4}^{(i_4)}\ \ \ 
(i_1, i_2, i_3, i_4=0, 1,\ldots,m)
$$

\vspace{4mm}
\noindent
the following 
expansion 

\vspace{-1mm}
$$
I_{T,t}^{*(i_1 i_2 i_3 i_4)}=
\hbox{\vtop{\offinterlineskip\halign{
\hfil#\hfil\cr
{\rm l.i.m.}\cr
$\stackrel{}{{}_{p\to \infty}}$\cr
}} }
\sum\limits_{j_1, j_2, j_3, j_4=0}^{p}
C_{j_4 j_3 j_2 j_1}\zeta_{j_1}^{(i_1)}\zeta_{j_2}^{(i_2)}\zeta_{j_3}^{(i_3)}
\zeta_{j_4}^{(i_4)}
$$

\vspace{4mm}
\noindent
that converges in the mean-square sense 
is valid, where

$$
C_{j_4 j_3 j_2 j_1}=\int\limits_t^T\phi_{j_4}(t_4)\int\limits_t^{t_4}
\phi_{j_3}(t_3)
\int\limits_t^{t_3}\phi_{j_2}(t_2)\int\limits_t^{t_2}\phi_{j_1}(t_1)
dt_1dt_2dt_3dt_4,
$$

\vspace{4mm}
\noindent
${\bf w}_{\tau}^{(i)}={\bf f}_{\tau}^{(i)}$ $(i=1,\ldots,m)$ 
are independent 
standard Wiener processes and 
${\bf w}_{\tau}^{(0)}=\tau.$}

\vspace{2mm}

It is obvious that Theorem 12 prove Hypothesis 1 
for the case $k=4$.

\vspace{5mm}

\section{Expansion of Iterated 
Stratonovich Stochastic Integrals of Arbitrary 
Multiplicity $k$ $(k\in \mathbb{N})$. Proof of Hypothesis 2 Under the Condition of 
Convergence of Trace Series}

\vspace{5mm}

In this section, we 
prove the expansion of iterated 
Stratonovich stochastic integrals of arbitrary 
multiplicity $k$ ($k\in \mathbb{N}$)
under the condition of 
convergence of trace series.

Let us introduce some notations
and formulate some auxiliary results.

Consider the unordered
set $\{1, 2, \ldots, k\}$ 
and separate it into two parts:
the first part consists of $r$ unordered 
pairs (sequence order of these pairs is also unimportant) and the 
second one consists of the 
remaining $k-2r$ numbers.
So, we have 

\vspace{-2mm}
\begin{equation}
\label{leto5007xxx}
(\{
\underbrace{\{g_1, g_2\}, \ldots, 
\{g_{2r-1}, g_{2r}\}}_{\small{\hbox{part 1}}}
\},
\{\underbrace{q_1, \ldots, q_{k-2r}}_{\small{\hbox{part 2}}}
\}),
\end{equation}

\vspace{3mm}
\noindent
where 

\vspace{-4mm}
$$
\{g_1, g_2, \ldots, 
g_{2r-1}, g_{2r}, q_1, \ldots, q_{k-2r}\}=\{1, 2, \ldots, k\},
$$

\vspace{3mm}
\noindent
braces   
mean an unordered 
set, and pa\-ren\-the\-ses mean an ordered set.

Recall that the expression (\ref{leto5007xxx}) is called the partition. 
Let us consider the sum with respect to all possible
partitions

\vspace{-1mm}
$$
\sum_{\stackrel{(\{\{g_1, g_2\}, \ldots, 
\{g_{2r-1}, g_{2r}\}\}, \{q_1, \ldots, q_{k-2r}\})}
{{}_{\{g_1, g_2, \ldots, 
g_{2r-1}, g_{2r}, q_1, \ldots, q_{k-2r}\}=\{1, 2, \ldots, k\}}}}
a_{g_1 g_2, \ldots, 
g_{2r-1} g_{2r}, q_1 \ldots q_{k-2r}},
$$

\vspace{6mm}
\noindent
where $a_{g_1 g_2, \ldots, 
g_{2r-1} g_{2r}, q_1 \ldots q_{k-2r}}\in \mathbb{R}.$

Consider the Fourier coefficient

\vspace{-1mm}
\begin{equation}
\label{after3000}
C_{j_k \ldots j_1}=\int\limits_t^T\psi_k(t_k)\phi_{j_k}(t_k)\ldots
\int\limits_t^{t_2}
\psi_1(t_1)\phi_{j_1}(t_1)
dt_1\ldots dt_k
\end{equation}

\vspace{3mm}
\noindent
corresponding to the function {\rm (\ref{ppp}), where 
$\{\phi_j(x)\}_{j=0}^{\infty}$ is a complete orthonormal system
of functions 
in the space $L_2([t, T])$. At that we suppose $\phi_0(x)=1/\sqrt{T-t}.$ 

Denote

\vspace{-1mm}
$$
C_{j_k \ldots j_{l+1}j_l j_l j_{l-2} \ldots j_1}\biggl|_{(j_l j_l)\curvearrowright (\cdot) }\biggr.
\stackrel{\sf def}{=}
$$

\vspace{2mm}
$$
\stackrel{\sf def}{=}\int\limits_t^T\psi_k(t_k)\phi_{j_k}(t_k)\ldots
\int\limits_t^{t_{l+2}}\psi_{l+1}(t_{l+1})\phi_{j_{l+1}}(t_{l+1})
\int\limits_t^{t_{l+1}}\psi_{l}(t_{l})\psi_{l-1}(t_{l})\times
$$

\begin{equation}
\label{after900}
\times
\int\limits_t^{t_{l}}\psi_{l-2}(t_{l-2})\phi_{j_{l-2}}(t_{l-2})\ldots 
\int\limits_t^{t_2}
\psi_1(t_1)\phi_{j_1}(t_1)
dt_1\ldots dt_{l-2}dt_{l}t_{l+1}\ldots dt_k=
\end{equation}

\vspace{2mm}
$$
=\sqrt{T-t}\int\limits_t^T\psi_k(t_k)\phi_{j_k}(t_k)\ldots
\int\limits_t^{t_{l+2}}\psi_{l+1}(t_{l+1})\phi_{j_{l+1}}(t_{l+1})
\int\limits_t^{t_{l+1}}\psi_{l}(t_{l})\psi_{l-1}(t_{l})\phi_0(t_l)\times
$$

\vspace{2mm}
$$
\times
\int\limits_t^{t_{l}}\psi_{l-2}(t_{l-2})\phi_{j_{l-2}}(t_{l-2})\ldots 
\int\limits_t^{t_2}
\psi_1(t_1)\phi_{j_1}(t_1)
dt_1\ldots dt_{l-2}dt_{l}t_{l+1}\ldots dt_k=
$$

\vspace{2mm}
$$
=\sqrt{T-t}\hat C_{j_k \ldots j_{l+1} 0 j_{l-2} \ldots j_1},
$$

\vspace{5mm}
\noindent
i.e. $\sqrt{T-t}\hat C_{j_k \ldots j_{l+1} 0 j_{l-2} \ldots j_1}$
is again the Fourier coefficient of type $C_{j_k \ldots j_1}$
but with a new shorter multi-index
$j_k \ldots j_{l+1} 0 j_{l-2} \ldots j_1$ 
and new weight functions $\psi_1(\tau),$ $\ldots,$ $\psi_{l-2}(\tau),$
$\sqrt{T-t}\psi_{l-1}(\tau)\psi_{l}(\tau),$ $\psi_{l+1}(\tau)$, $\ldots,$
$\psi_{k}(\tau)$ (also we suppose that $\{l, l-1\}$ is one of the pairs
$\{g_1, g_2\}, \ldots,\{g_{2r-1}, g_{2r}\}$).

Let

\vspace{-3mm}
$$
C_{j_k \ldots j_{l+1}j_l j_l j_{l-2} \ldots j_1}\biggl|_{(j_l j_l)\curvearrowright j_m}\biggr.
\stackrel{\sf def}{=}
$$

\vspace{3mm}
$$
\stackrel{\sf def}{=}\int\limits_t^T\psi_k(t_k)\phi_{j_k}(t_k)\ldots
\int\limits_t^{t_{l+2}}\psi_{l+1}(t_{l+1})\phi_{j_{l+1}}(t_{l+1})
\int\limits_t^{t_{l+1}}\psi_{l}(t_{l})\psi_{l-1}(t_{l})\phi_{j_m}(t_l)\times
$$

\vspace{1mm}
\begin{equation}
\label{after2000}
\times
\int\limits_t^{t_{l}}\psi_{l-2}(t_{l-2})\phi_{j_{l-2}}(t_{l-2})\ldots 
\int\limits_t^{t_2}
\psi_1(t_1)\phi_{j_1}(t_1)
dt_1\ldots dt_{l-2}dt_{l}t_{l+1}\ldots dt_k=
\end{equation}

\vspace{2mm}
$$
= \bar C_{j_k \ldots j_{l+1} j_m j_{l-2} \ldots j_1},
$$

\vspace{5mm}
\noindent
i.e. $\bar C_{j_k \ldots j_{l+1} j_m j_{l-2} \ldots j_1}$
is again the Fourier coefficient of type $C_{j_k \ldots j_1}$
but with a new shorter multi-index
$j_k \ldots j_{l+1} j_m j_{l-2} \ldots j_1$ 
and new weight functions $\psi_1(\tau),$ $\ldots,$ $\psi_{l-2}(\tau),$
$\psi_{l-1}(\tau)\psi_{l}(\tau),$ $\psi_{l+1}(\tau)$, $\ldots,$
$\psi_{k}(\tau)$ (also we suppose that $\{l-1, l\}$ is one of the pairs
$\{g_1, g_2\}, \ldots,\{g_{2r-1}, g_{2r}\}$).

Denote
$$
\bar C^{(p)}_{j_k\ldots j_q \ldots j_1}\biggl|_{q\ne g_1,g_2,\ldots,g_{2r-1}, g_{2r}}
\stackrel{\sf def}{=}
$$

\vspace{2mm}
\begin{equation}
\label{nov100}
\stackrel{\sf def}{=}
\sum\limits_{j_{g_{2r-1}}=p+1}^{\infty}
\sum\limits_{j_{g_{2r-3}}=p+1}^{\infty}
\ldots \sum\limits_{j_{g_{3}}=p+1}^{\infty}
\sum\limits_{j_{g_{1}}=p+1}^{\infty}
C_{j_k \ldots j_1}\biggl|_{j_{g_1}=j_{g_2},\ldots, j_{g_{2r-1}}=j_{g_{2r}}}\biggl..
\end{equation}

\vspace{6mm}

Introduce the following notation

$$
S_l \left\{\bar C^{(p)}_{j_k\ldots j_q \ldots j_1}\biggl|_{q\ne g_1,g_2,\ldots,g_{2r-1}, g_{2r}}
\right\}
\stackrel{\sf def}{=}
\frac{1}{2}{\bf 1}_{\{g_{2l}=g_{2l-1}+1\}}
\sum\limits_{j_{g_{2r-1}}=p+1}^{\infty}
\sum\limits_{j_{g_{2r-3}}=p+1}^{\infty}
\ldots 
$$

\vspace{2mm}
\begin{equation}
\label{nov101}
\ldots
\sum\limits_{j_{g_{2l+1}}=p+1}^{\infty}
\sum\limits_{j_{g_{2l-3}}=p+1}^{\infty}
\ldots
\sum\limits_{j_{g_{3}}=p+1}^{\infty}
\sum\limits_{j_{g_{1}}=p+1}^{\infty}
C_{j_k \ldots j_1}\biggl|_{(j_{g_{2l}} j_{g_{2l-1}})\curvearrowright (\cdot),
j_{g_1}=j_{g_2},\ldots, j_{g_{2r-1}}=j_{g_{2r}}}\biggr..
\end{equation}

\vspace{5mm}

Note that the operation $S_l$ $(l=1,2,\ldots,r)$ acts on the value

\begin{equation}
\label{after301}
\bar C^{(p)}_{j_k\ldots j_q \ldots j_1}\biggl|_{q\ne g_1,g_2,\ldots,g_{2r-1}, g_{2r}}
\end{equation}

\vspace{4mm}
\noindent
as follows: $S_l$ multiplies (\ref{after301})
by ${\bf 1}_{\{g_{2l}=g_{2l-1}+1\}}/2,$ removes the summation

\vspace{1mm}
$$
\sum\limits_{j_{g_{2l-1}}=p+1}^{\infty},
$$

\vspace{3mm}
\noindent
and replaces 
$$
C_{j_k \ldots j_1}\biggl|_{j_{g_1}=j_{g_2},\ldots, j_{g_{2r-1}}=j_{g_{2r}}}\biggl.
$$

\vspace{3mm}
\noindent
with
\begin{equation}
\label{after300}
C_{j_k \ldots j_1}\biggl|_{(j_{g_{2l}} j_{g_{2l-1}})\curvearrowright (\cdot),
j_{g_1}=j_{g_2},\ldots, j_{g_{2r-1}}=j_{g_{2r}}}\biggr..
\end{equation}

\vspace{3mm}

Note that we write

$$
C_{j_k \ldots j_1}\biggl|_{(j_{g_{1}} j_{g_{2}})\curvearrowright (\cdot),
j_{g_{1}}=j_{g_{2}}}\biggr.
=
C_{j_k \ldots j_1}\biggl|_{(j_{g_{1}} j_{g_{1}})\curvearrowright (\cdot),
j_{g_{1}}=j_{g_{2}}},\biggr.
$$

\vspace{2mm}
$$                                                                   
C_{j_k \ldots j_1}\biggl|_{(j_{g_{1}} j_{g_{2}})\curvearrowright j_{m},
j_{g_{1}}=j_{g_{2}}}\biggr.
=
C_{j_k \ldots j_1}\biggl|_{(j_{g_{1}} j_{g_{1}})\curvearrowright j_{m},
j_{g_{1}}=j_{g_{2}}}\biggr.,\ \ \ 
$$

\vspace{2mm}
$$
C_{j_k \ldots j_1}\biggl|_{(j_{g_{1}} j_{g_{2}})\curvearrowright (\cdot),
(j_{g_{3}} j_{g_{4}})\curvearrowright (\cdot),
j_{g_{1}}=j_{g_{2}}, j_{g_{3}}=j_{g_{4}}}\biggr.
=
C_{j_k \ldots j_1}\biggl|_{(j_{g_{1}} j_{g_{1}})\curvearrowright (\cdot)
(j_{g_{3}} j_{g_{3}})\curvearrowright (\cdot),
j_{g_{1}}=j_{g_{2}}, j_{g_{3}}=j_{g_{4}}}\biggr.,\ \ \ \ldots
$$

\vspace{6mm}
        
Since (\ref{after300}) is again the Fourier coefficient, 
then the action 
of superposition $S_l S_m$ on 
(\ref{after300}) is obvious. For example,
for $r=3$

\vspace{-1mm}
$$
S_3S_2S_1
\left\{\bar C^{(p)}_{j_k\ldots j_q \ldots j_1}\biggl|_{q\ne g_1,g_2,\ldots,g_{5}, g_{6}}
\right\}
=
$$

\vspace{2mm}
$$
=\frac{1}{2^3}
\prod_{s=1}^3
{\bf 1}_{\{g_{2s}=g_{2s-1}+1\}}
C_{j_k \ldots j_1}\Biggl|_{(j_{g_2} j_{g_1})\curvearrowright (\cdot)
(j_{g_{4}} j_{g_{3}})\curvearrowright (\cdot)
(j_{g_{6}} j_{g_{5}})\curvearrowright (\cdot),
j_{g_{1}}=j_{g_{2}}, j_{g_{3}}=j_{g_{4}}, j_{g_{5}}=j_{g_{6}}}\Biggr.,
$$

\vspace{5mm}
$$
S_3S_1
\left\{\bar C^{(p)}_{j_k\ldots j_q \ldots j_1}\biggl|_{q\ne g_1,g_2,\ldots,g_{5}, g_{6}}
\right\}
=
$$

\vspace{2mm}
$$
=\frac{1}{2^2}
{\bf 1}_{\{g_{6}=g_{5}+1\}}{\bf 1}_{\{g_{2}=g_{1}+1\}}
\sum_{j_{g_3}=p+1}^{\infty}
C_{j_k \ldots j_1}\Biggl|_{(j_{g_{2}} j_{g_{1}})\curvearrowright (\cdot)
(j_{g_{6}} j_{g_{5}})\curvearrowright (\cdot),
j_{g_{1}}=j_{g_{2}}, j_{g_{3}}=j_{g_{4}}, j_{g_{5}}=j_{g_{6}}}\Biggr.,
$$

\vspace{5mm}
$$
S_2
\left\{\bar C^{(p)}_{j_k\ldots j_q \ldots j_1}\biggl|_{q\ne g_1,g_2,\ldots,g_{5}, g_{6}}
\right\}
=
$$

\vspace{2mm}
$$
=\frac{1}{2}
{\bf 1}_{\{g_{4}=g_{3}+1\}}
\sum_{j_{g_1}=p+1}^{\infty}\sum_{j_{g_5}=p+1}^{\infty}
C_{j_k \ldots j_1}\Biggl|_{(j_{g_{4}} j_{g_{3}})\curvearrowright (\cdot),
j_{g_{1}}=j_{g_{2}}, j_{g_{3}}=j_{g_{4}}, j_{g_{5}}=j_{g_{6}}}\Biggr..
$$

\vspace{6mm}

{\bf Theorem~13}\ \cite{2018a}, \cite{arxiv-4}, \cite{arxiv-5}, \cite{arxiv-9},
\cite{new-art-1-xxy}.\ 
{\it Assume that
the continuously differentiable func\-ti\-ons 
$\psi_l(\tau)$ $(l=1,\ldots,k)$ and 
the complete orthonormal system $\{\phi_j(x)\}_{j=0}^{\infty}$
of continuous functions $(\phi_0(x)=1/\sqrt{T-t})$ 
in the space $L_2([t, T])$ are such that the following 
conditions are satisfied{\rm :}

\vspace{2mm}

{\rm 1.}\ The equality 

\vspace{-2mm}
\begin{equation}
\label{after200}
\frac{1}{2}
\int\limits_t^s \Phi_1(t_1)\Phi_2(t_1)dt_1
=\sum_{j_1=0}^{\infty}
\int\limits_t^s
\Phi_2(t_2)\phi_{j_1}(t_2)\int\limits_t^{t_2}
\Phi_1(t_1)\phi_{j_1}(t_1)dt_1 dt_2
\end{equation}

\vspace{3mm}
\noindent
holds for all $s\in (t, T],$ where the nonrandom functions 
$\Phi_1(\tau),$ $\Phi_2(\tau)$
are continuously differentiable on $[t, T]$
and the series on the right-hand side of {\rm (\ref{after200})}
converges absolutely.

{\rm 2.}\ The estimates

\vspace{-1mm}
$$
\left|\int\limits_t^s \phi_{j}(\tau)\Phi_1(\tau)d\tau\right|
\le \frac{\Psi_1(s)}{j^{1/2+\alpha}},\ \ \ 
\left|\int\limits_s^T \phi_{j}(\tau)\Phi_2(\tau)d\tau\right|\le
\frac{\Psi_1(s)}{j^{1/2+\alpha}},
$$

\vspace{2mm}
$$
\left|\sum_{j=p+1}^{\infty}\int\limits_t^s
\Phi_2(\tau)\phi_{j}(\tau)\int\limits_t^{\tau}
\Phi_1(\theta)\phi_{j}(\theta)d\theta d\tau\right|\le \frac{\Psi_2(s)}{p^{\beta}}
$$

\vspace{4mm}
\noindent
hold for all $s\in (t, T)$ and for some $\alpha, \beta >0,$ where 
$\Phi_1(\tau),$ $\Phi_2(\tau)$
are continuously differentiable nonrandom functions on $[t, T],$\ $j, p\in \mathbb{N},$
and

\vspace{-1mm}
$$
\int\limits_t^T \Psi_1^2(\tau) d\tau<\infty,\ \ \ 
\int\limits_t^T \left|\Psi_2(\tau)\right| d\tau<\infty.
$$

\vspace{2mm}

{\rm 3.}\ The condition

\vspace{1mm}
$$
\lim\limits_{p\to\infty}
\sum\limits_{\stackrel{j_1,\ldots,j_q,\ldots,j_k=0}{{}_{q\ne g_1, g_2, \ldots, g_{2r-1},
g_{2r}}}}^p
\left(S_{l_1}S_{l_2}\ldots S_{l_{d}}
\left\{\bar C^{(p)}_{j_k\ldots j_q \ldots j_1}\biggl|_{q\ne g_1,g_2,\ldots,g_{2r-1}, g_{2r}}
\right\}\right)^2=0
$$

\vspace{4mm}
\noindent
holds for all possible $g_1,g_2,\ldots,g_{2r-1},g_{2r}$ {\rm (}see {\rm (\ref{leto5007xxx}))}
and $l_1, l_2, \ldots, l_{d}$ such that
$l_1, l_2, \ldots, l_{d}\in \{1,2,$ $\ldots,$ $r\},$\
$l_1>l_2>\ldots >l_{d},$\ $d=0, 1, 2,\ldots, r-1,$\ 
where $r=1, 2,\ldots,[k/2]$ and

\vspace{1mm}
$$
S_{l_1}S_{l_2}\ldots S_{l_{d}}
\left\{\bar C^{(p)}_{j_k\ldots j_q \ldots j_1}\biggl|_{q\ne g_1,g_2,\ldots,g_{2r-1}, g_{2r}}
\right\}\stackrel{\sf def}{=}
\bar C^{(p)}_{j_k\ldots j_q \ldots j_1}\biggl|_{q\ne g_1,g_2,\ldots,g_{2r-1}, g_{2r}}
$$

\vspace{3mm}
\noindent
for $d=0.$

Then, for the iterated Stratonovich stochastic integral 
of arbitrary multiplicity $k$

\vspace{-1mm}
\begin{equation}
\label{afterstr}
J^{*}[\psi^{(k)}]_{T,t}^{(i_1\ldots i_k)}=
{\int\limits_t^{*}}^T
\psi_k(t_k) \ldots 
{\int\limits_t^{*}}^{t_2}
\psi_1(t_1) d{\bf w}_{t_1}^{(i_1)}\ldots
d{\bf w}_{t_k}^{(i_k)}
\end{equation}

\vspace{2mm}
\noindent
the following 
expansion 

\vspace{-4mm}

\begin{equation}
\label{after1}
J^{*}[\psi^{(k)}]_{T,t}^{(i_1\ldots i_k)}=
\hbox{\vtop{\offinterlineskip\halign{
\hfil#\hfil\cr
{\rm l.i.m.}\cr
$\stackrel{}{{}_{p\to \infty}}$\cr
}} }
\sum\limits_{j_1,\ldots,j_k=0}^{p}
C_{j_k \ldots j_1}\prod\limits_{l=1}^k \zeta_{j_l}^{(i_l)}
\end{equation}

\vspace{3mm}
\noindent
that converges in the mean-square sense is valid, where 

\vspace{-3mm}
\begin{equation}
\label{after1000}
C_{j_k \ldots j_1}=\int\limits_t^T\psi_k(t_k)\phi_{j_k}(t_k)\ldots
\int\limits_t^{t_2}
\psi_1(t_1)\phi_{j_1}(t_1)
dt_1\ldots dt_k
\end{equation}

\vspace{2mm}
\noindent
is the Fourier coefficient, 
${\rm l.i.m.}$ is a limit in the mean-square sense,
$i_1, \ldots, i_k=0, 1,\ldots,m,$

\vspace{-2mm}
$$
\zeta_{j}^{(i)}=
\int\limits_t^T \phi_{j}(s) d{\bf w}_s^{(i)}
$$ 

\vspace{2mm}
\noindent
are independent standard Gaussian random variables for various 
$i$ or $j$ {\rm (}in the case when $i\ne 0${\rm )},
${\bf w}_{\tau}^{(i)}={\bf f}_{\tau}^{(i)}$ 
for $i=1,\ldots,m$ and 
${\bf w}_{\tau}^{(0)}=\tau.$}

\vspace{2mm}

{\bf Proof.} First note that (\ref{after200}) is fulfilled 
(see \cite{2018a}, Sect.~2.1.4 or \cite{rybakov5000}). The proof of Theorem~13 will consist 
of several steps.

{\bf Step~1.}\ Let us find a representation of the quantity

$$
\sum\limits_{j_1,\ldots,j_k=0}^{p}
C_{j_k \ldots j_1}\prod\limits_{l=1}^k \zeta_{j_l}^{(i_l)}
$$

\vspace{3mm}
\noindent
that will be convenient for further consideration.

Let us consider the following
multiple stochastic integral 

\begin{equation}
\label{mult11www}
\hbox{\vtop{\offinterlineskip\halign{
\hfil#\hfil\cr
{\rm l.i.m.}\cr
$\stackrel{}{{}_{N\to \infty}}$\cr
}} }
\sum\limits_{\stackrel{j_1,\ldots,j_k=0}{{}_{j_q\ne j_r;\ q\ne r;\ 
q, r=1,\ldots, k}}}^{N-1}
\Phi\left(\tau_{j_1},\ldots,\tau_{j_k}\right)
\prod\limits_{l=1}^k
\Delta{\bf w}_{\tau_{j_l}}^{(i_l)}
\stackrel{\rm def}{=}J'[\Phi]_{T,t}^{(i_1\ldots i_k)},
\end{equation}

\vspace{3mm}
\noindent
where for simplicity we assume that
$\Phi(t_1,\ldots,t_k):\ [t, T]^k\to\mathbb{R}$ is a 
continuous nonrandom
function on $[t, T]^k.$ Moreover, 
$\Delta{\bf w}_{\tau_{j}}^{(i)}=
{\bf w}_{\tau_{j+1}}^{(i)}-{\bf w}_{\tau_{j}}^{(i)}$
$(i=0, 1,\ldots,m),$
$\left\{\tau_{j}\right\}_{j=0}^{N}$ is a partition of
$[t,T],$ which satisfies the condition {\rm (\ref{1111}),}
$i_1,\ldots,i_k=0, 1,\ldots, m.$

The stochastic integral with respect to the scalar standard Wiener process
$(i_1=\ldots=i_k\ne 0)$
and similar to {\rm (\ref{mult11www})} was considered in {\rm \cite{ito1951} (1951)}
and is called the multiple Wiener stochastic integral {\rm \cite{ito1951}.}

Note that the following well known estimate

\begin{equation}
\label{wiener1}
{\sf M}\left\{\left(J'[\Phi]_{T,t}^{(i_1\ldots i_k)}\right)^2\right\}\le
C_k
\int\limits_{[t,T]^k}
\Phi^2(t_1,\ldots,t_k)dt_1\ldots dt_k
\end{equation}

\vspace{2mm}
\noindent
is true for the multiple Wiener stochastic integral,
where $J'[\Phi]_{T,t}^{(i_1\ldots i_k)}$ is defined by {\rm (\ref{mult11www})}
and $C_k$ is a constant.

From the proof of Theorem~3 (see the proof of Theorem~5.1 in 
the original paper \cite{2006} (2006) in Russian
or proof of Theorem~1.1 in the monographs \cite{2018a}-\cite{12aa-afterxxx}
in English)
it follows that (\ref{tyyy}) can be written as 

\vspace{-2mm}
\begin{equation}
\label{tyyyarg}
J[\psi^{(k)}]_{T,t}^{(i_1\ldots i_k)}\  =\ 
\hbox{\vtop{\offinterlineskip\halign{
\hfil#\hfil\cr
{\rm l.i.m.}\cr
$\stackrel{}{{}_{p_1,\ldots,p_k\to \infty}}$\cr
}} }\sum_{j_1=0}^{p_1}\ldots\sum_{j_k=0}^{p_k}
C_{j_k\ldots j_1}J'[\phi_{j_1}\ldots \phi_{j_k}]^{(i_1\ldots i_k)}_{T,t},
\end{equation}

\vspace{3mm}
\noindent 
where $J'[\phi_{j_1}\ldots \phi_{j_k}]^{(i_1\ldots i_k)}_{T,t}$
is the multiple Wiener stochatic integral defined by (\ref{mult11www}) and
$J[\psi^{(k)}]_{T,t}^{(i_1\ldots i_k)}$ is the iterated Ito stochastic integral
(\ref{ito}), i.e.

\vspace{-1mm}
\begin{equation}
\label{afterxgxg1}
J[\psi^{(k)}]_{T,t}^{(i_1\ldots i_k)}=\int\limits_t^T\psi_k(t_k) \ldots 
\int\limits_t^{t_{2}}
\psi_1(t_1) d{\bf w}_{t_1}^{(i_1)}\ldots
d{\bf w}_{t_k}^{(i_k)}.
\end{equation}

\vspace{2mm}

Let us consider the following
multiple stochastic integral 

\vspace{-1mm}
\begin{equation}
\label{30.34ququ}
\hbox{\vtop{\offinterlineskip\halign{
\hfil#\hfil\cr
{\rm l.i.m.}\cr
$\stackrel{}{{}_{N\to \infty}}$\cr
}} }\sum_{j_1,\ldots,j_k=0}^{N-1}
\Phi\left(\tau_{j_1},\ldots,\tau_{j_k}\right)
\prod\limits_{l=1}^k\Delta{\bf w}_{\tau_{j_l}}^{(i_l)}
\stackrel{\rm def}{=}J[\Phi]_{T,t}^{(i_1\ldots i_k)},
\end{equation}

\vspace{2mm}
\noindent
where we assume that
$\Phi(t_1,\ldots,t_k):\ [t, T]^k\to\mathbb{R}$ is a 
continuous nonrandom
function on $[t, T]^k.$ 
Another notations are the same as in 
(\ref{mult11www}).

The stochastic integral with respect to the scalar standard Wiener process
($i_1=\ldots=i_k\ne 0$)
and similar to (\ref{30.34ququ}) 
(the function $\Phi(t_1,\ldots,t_k)$ is assumed to be symmetric
on the hypercube $[t, T]^k$)
has been considered in literature
(see, for example, Remark~1.5.7 \cite{bugh1}). 
The integral (\ref{30.34ququ})
is sometimes called the multiple Stratonovich stochastic integral.
This is due to the fact that the following rule
of the classical integral calculus holds for this integral

\vspace{1mm}
$$
J[\Phi]_{T,t}^{(i_1\ldots i_k)}=
J[\varphi_1]_{T,t}^{(i_1)}\ldots J[\varphi_k]_{T,t}^{(i_k)}\ \ \ \hbox{\rm w.~p.~1},
$$

\vspace{3mm}
\noindent
where $\Phi(t_1,\ldots,t_k)=\varphi_1(t_1)\ldots \varphi_k(t_k)$
and
$$
J[\varphi_l]_{T,t}^{(i_l)}
=\int\limits_t^T \varphi_l(s) d{\bf w}_{s}^{(i_l)}\ \ \ (l=1,\ldots,k).
$$

\vspace{2mm}

{\bf Theorem~14}\ \cite{2018a}, \cite{12aa-afterxxx}.\ 
{\it Suppose that $\Phi(t_1,\ldots,t_k):\ [t, T]^k\to\mathbb{R}$ is a 
continuous nonrandom
func\-ti\-on on $[t, T]^k.$ Furthermore,
$\{\phi_j(x)\}_{j=0}^{\infty}$ is a complete orthonormal system  
of functions in the space $L_2([t,T]),$ each function $\phi_j(x)$ of which 
for finite $j$ is continuous at the interval $[t, T]$ except may be
for the finite number of points 
of the finite discontinuity as well as $\phi_j(x)$
right-continuous at the interval  $[t, T].$
Then the following expansion

$$
J'[\Phi]_{T,t}^{(i_1\ldots i_k)}=
\hbox{\vtop{\offinterlineskip\halign{
\hfil#\hfil\cr
{\rm l.i.m.}\cr
$\stackrel{}{{}_{p_1,\ldots,p_k\to \infty}}$\cr
}} }
\sum\limits_{j_1=0}^{p_1}\ldots
\sum\limits_{j_k=0}^{p_k}
C_{j_k\ldots j_1}\Biggl(
\prod_{l=1}^k\zeta_{j_l}^{(i_l)}+\sum\limits_{r=1}^{[k/2]}
(-1)^r \times
\Biggr.
$$

\vspace{2mm}
\begin{equation}
\label{quq11}
\times
\sum_{\stackrel{(\{\{g_1, g_2\}, \ldots, 
\{g_{2r-1}, g_{2r}\}\}, \{q_1, \ldots, q_{k-2r}\})}
{{}_{\{g_1, g_2, \ldots, 
g_{2r-1}, g_{2r}, q_1, \ldots, q_{k-2r}\}=\{1, 2, \ldots, k\}}}}
\prod\limits_{s=1}^r
{\bf 1}_{\{i_{g_{{}_{2s-1}}}=~i_{g_{{}_{2s}}}\ne 0\}}
\Biggl.{\bf 1}_{\{j_{g_{{}_{2s-1}}}=~j_{g_{{}_{2s}}}\}}
\prod_{l=1}^{k-2r}\zeta_{j_{q_l}}^{(i_{q_l})}\Biggr)
\end{equation}

\vspace{4mm}
\noindent
con\-verg\-ing in the mean-square sense is valid$,$
where $J'[\Phi]^{(i_1\ldots i_k)}_{T,t}$
is the multiple Wiener stochatic integral defined by {\rm (\ref{mult11www}),}
\begin{equation}
\label{quq12}
C_{j_k\ldots j_1}=\int\limits_{[t,T]^k}
\Phi(t_1,\ldots,t_k)\prod_{l=1}^{k}\phi_{j_l}(t_l)dt_1\ldots dt_k
\end{equation}

\vspace{2mm}
\noindent
is the Fourier coefficient. Another notations are the same as in Theorems~{\rm 3, 4}.}

\vspace{2mm}

From (\ref{leto6000}) and (\ref{tyyyarg}) (also see Theorem~5 in \cite{diffjournal}
or Theorem~5 in \cite{new-2023a}) we conclude that

\vspace{1mm}
$$
J'[\phi_{j_1}\ldots \phi_{j_k}]_{T,t}^{(i_1\ldots i_k)}=\prod_{l=1}^k\zeta_{j_l}^{(i_l)}+
$$

\begin{equation}
\label{2023abc300}
+
\sum\limits_{r=1}^{[k/2]}
(-1)^r 
\hspace{-3mm}\sum_{\stackrel{(\{\{g_1, g_2\}, \ldots, 
\{g_{2r-1}, g_{2r}\}\}, \{q_1, \ldots, q_{k-2r}\})}
{{}_{\{g_1, g_2, \ldots, 
g_{2r-1}, g_{2r}, q_1, \ldots, q_{k-2r}\}=\{1, 2, \ldots, k\}}}}
\prod\limits_{s=1}^r
{\bf 1}_{\{i_{g_{{}_{2s-1}}}=~i_{g_{{}_{2s}}}\ne 0\}}
\Biggl.{\bf 1}_{\{j_{g_{{}_{2s-1}}}=~j_{g_{{}_{2s}}}\}}
\prod_{l=1}^{k-2r}\zeta_{j_{q_l}}^{(i_{q_l})}
\end{equation}

\vspace{5mm}
\noindent
w.~p.~1, where notations are the same as in Theorem 4
and
$J'[\phi_{j_1}\ldots\phi_{j_k}]_{T,t}^{(i_1\ldots i_k)}$
is the multiple Wiener stochastic
integral (\ref{mult11www}). For a more detailed derivation of (\ref{2023abc300}), 
see \cite{diffjournal}, \cite{new-2023a} (also see \cite{2018a})).

Using (\ref{2023abc300}), we obtain

\vspace{1mm}
$$
\prod_{l=1}^k\zeta_{j_l}^{(i_l)}=J'[\phi_{j_1}\ldots \phi_{j_k}]_{T,t}^{(i_1\ldots i_k)}-
$$

\begin{equation}
\label{2023abc300xz1}
-\sum\limits_{r=1}^{[k/2]}
(-1)^r 
\hspace{-3mm}\sum_{\stackrel{(\{\{g_1, g_2\}, \ldots, 
\{g_{2r-1}, g_{2r}\}\}, \{q_1, \ldots, q_{k-2r}\})}
{{}_{\{g_1, g_2, \ldots, 
g_{2r-1}, g_{2r}, q_1, \ldots, q_{k-2r}\}=\{1, 2, \ldots, k\}}}}
\prod\limits_{s=1}^r
{\bf 1}_{\{i_{g_{{}_{2s-1}}}=~i_{g_{{}_{2s}}}\ne 0\}}
\Biggl.{\bf 1}_{\{j_{g_{{}_{2s-1}}}=~j_{g_{{}_{2s}}}\}}
\prod_{l=1}^{k-2r}\zeta_{j_{q_l}}^{(i_{q_l})}
\end{equation}

\vspace{4mm}
\noindent
w.~p.~1.

By iteratively applying the formula (\ref{2023abc300xz1})
(also see (\ref{a2})--(\ref{a6})), we obtain the following
representation of the product
$$
\prod_{l=1}^k\zeta_{j_l}^{(i_l)}
$$

\vspace{2mm}
\noindent
as the sum of some constant value and multiple Wiener stochastic integrals 
of 
multiplicities not exceeding $k$ 

\vspace{-1mm}
$$
\prod_{l=1}^k\zeta_{j_l}^{(i_l)}=J'[\phi_{j_1}\ldots \phi_{j_k}]_{T,t}^{(i_1\ldots i_k)}
+\sum\limits_{r=1}^{[k/2]}
\sum_{\stackrel{(\{\{g_1, g_2\}, \ldots, 
\{g_{2r-1}, g_{2r}\}\}, \{q_1, \ldots, q_{k-2r}\})}
{{}_{\{g_1, g_2, \ldots, 
g_{2r-1}, g_{2r}, q_1, \ldots, q_{k-2r}\}=\{1, 2, \ldots, k\}}}}
\prod\limits_{s=1}^r
{\bf 1}_{\{i_{g_{{}_{2s-1}}}=~i_{g_{{}_{2s}}}\ne 0\}}\times
$$

\vspace{3mm}
\begin{equation}
\label{after8xx1}
\times
{\bf 1}_{\{j_{g_{{}_{2s-1}}}=~j_{g_{{}_{2s}}}\}}
J'[\phi_{j_{q_1}}\ldots \phi_{j_{q_{k-2r}}}]_{T,t}^{(i_{q_1}\ldots i_{q_{k-2r}})}\ \ \ \hbox{w.~p.~1,}
\end{equation}

\vspace{6mm}
\noindent
where
$J'[\phi_{j_{q_1}}\ldots \phi_{j_{q_{k-2r}}}]_{T,t}^{(i_{q_1}\ldots i_{q_{k-2r}})}
\stackrel{\sf def}{=}1$
for $k=2r$.

Multiplying both sides of the equality (\ref{after8xx1}) by $C_{j_k\ldots j_1}$
and summing over $j_1,\ldots,j_k,$ we get w.~p.~1

\vspace{1mm}
$$
\sum_{j_1=0}^{p_1}\ldots\sum_{j_k=0}^{p_k}
C_{j_k\ldots j_1}
\prod_{l=1}^k\zeta_{j_l}^{(i_l)}
=
\sum_{j_1=0}^{p_1}\ldots\sum_{j_k=0}^{p_k}
C_{j_k\ldots j_1}
J'[\phi_{j_1}\ldots \phi_{j_k}]_{T,t}^{(i_1\ldots i_k)}+
$$

\vspace{3mm}
$$
+\sum_{j_1=0}^{p_1}\ldots\sum_{j_k=0}^{p_k}
C_{j_k\ldots j_1}
\sum\limits_{r=1}^{[k/2]}
\sum_{\stackrel{(\{\{g_1, g_2\}, \ldots, 
\{g_{2r-1}, g_{2r}\}\}, \{q_1, \ldots, q_{k-2r}\})}
{{}_{\{g_1, g_2, \ldots, 
g_{2r-1}, g_{2r}, q_1, \ldots, q_{k-2r}\}=\{1, 2, \ldots, k\}}}}
\prod\limits_{s=1}^r
{\bf 1}_{\{i_{g_{{}_{2s-1}}}=~i_{g_{{}_{2s}}}\ne 0\}}\times
$$

\vspace{3mm}
\begin{equation}
\label{after8}
\times{\bf 1}_{\{j_{g_{{}_{2s-1}}}=~j_{g_{{}_{2s}}}\}}
J'[\phi_{j_{q_1}}\ldots \phi_{j_{q_{k-2r}}}]_{T,t}^{(i_{q_1}\ldots i_{q_{k-2r}})}\ \ \ \hbox{w.~p.~1.}
\end{equation}

\vspace{5mm}

Denote
\begin{equation}
\label{afterxx1}
K_{p_1\ldots p_k}(t_1,\ldots,t_k)=
\sum_{j_1=0}^{p_1}\ldots\sum_{j_k=0}^{p_k}
C_{j_k\ldots j_1}
\prod_{l=1}^k\phi_{j_l}(t_l),
\end{equation}

\vspace{3mm}
\begin{equation}
\label{afterxx2}
K_{p_1\ldots p_k}^{g_1\ldots g_{2r}, q_1\ldots q_{k-2r}}(t_{q_1},\ldots,t_{q_{k-2r}})=
\sum_{j_1=0}^{p_1}\ldots\sum_{j_k=0}^{p_k}
C_{j_k\ldots j_1}
\prod\limits_{s=1}^r
{\bf 1}_{\{j_{g_{{}_{2s-1}}}=~j_{g_{{}_{2s}}}\}}
\prod_{l=1}^{k-2r}\phi_{j_{q_l}}(t_{q_l}),
\end{equation}

\vspace{5mm}
\noindent
where $C_{j_k\ldots j_1}$ is defined by (\ref{after1000}) and
$\prod\limits_{\emptyset}\stackrel{\sf def}{=}1.$

\vspace{2mm}

The equality (\ref{after8}) can be written as 

\vspace{2mm}
$$
J[K_{p_1\ldots p_k}]_{T,t}^{(i_1\ldots i_k)}=
J'[K_{p_1\ldots p_k}]_{T,t}^{(i_1\ldots i_k)}+
$$

\vspace{1mm}
\begin{equation}
\label{after7}
+\sum\limits_{r=1}^{[k/2]}
\sum_{\stackrel{(\{\{g_1, g_2\}, \ldots, 
\{g_{2r-1}, g_{2r}\}\}, \{q_1, \ldots, q_{k-2r}\})}
{{}_{\{g_1, g_2, \ldots, 
g_{2r-1}, g_{2r}, q_1, \ldots, q_{k-2r}\}=\{1, 2, \ldots, k\}}}}
\prod\limits_{s=1}^r
{\bf 1}_{\{i_{g_{{}_{2s-1}}}=~i_{g_{{}_{2s}}}\ne 0\}}\
J'[K_{p_1\ldots p_k}^{g_1\ldots g_{2r}, q_1\ldots q_{k-2r}}]_{T,t}^{(i_{q_1}\ldots i_{q_{k-2r}})}
\end{equation}

\vspace{4mm}
\noindent
w.~p.~{\rm 1}, where
$K_{p_1\ldots p_k}(t_1,\ldots,t_k)$ and 
$K_{p_1\ldots p_k}^{g_1\ldots g_{2r}, q_1\ldots q_{k-2r}}(t_{q_1},\ldots,t_{q_{k-2r}})$
have the form (\ref{afterxx1}), (\ref{afterxx2}),
$J[K_{p_1\ldots p_k}]_{T,t}^{(i_1\ldots i_k)}$ 
is the multiple Stra\-to\-no\-vich stochastic integral defined by (\ref{30.34ququ}),
$J'[K_{p_1\ldots p_k}]_{T,t}^{(i_1\ldots i_k)}$ and 
$J'[K_{p_1\ldots p_k}^{g_1\ldots g_{2r}, q_1\ldots q_{k-2r}}]_{T,t}^{(i_{q_1}\ldots i_{q_{k-2r}})}$
are multiple Wiener stochastic
integrals defined by (\ref{mult11www}).             

Passing to the limit 
$\hbox{\vtop{\offinterlineskip\halign{
\hfil#\hfil\cr
{\rm l.i.m.}\cr
$\stackrel{}{{}_{p_1,\ldots,p_k\to \infty}}$\cr
}} }$ ($p_1=\ldots =p_k=p$) in (\ref{after8}) or (\ref{after7}), 
we get w.~p.~1 (see (\ref{tyyyarg})) 

\vspace{1mm}
$$
\hbox{\vtop{\offinterlineskip\halign{
\hfil#\hfil\cr
{\rm l.i.m.}\cr
$\stackrel{}{{}_{p\to \infty}}$\cr
}} }
\sum_{j_1,\ldots,j_k=0}^{p}
C_{j_k \ldots j_1}\prod\limits_{l=1}^k \zeta_{j_l}^{(i_l)}
=
J[\psi^{(k)}]_{T,t}^{(i_1\ldots i_k)}+
$$

\vspace{3mm}
$$
+
\hbox{\vtop{\offinterlineskip\halign{
\hfil#\hfil\cr
{\rm l.i.m.}\cr
$\stackrel{}{{}_{p\to \infty}}$\cr
}} }\sum_{j_1,\ldots,j_k=0}^{p}
C_{j_k \ldots j_1}
\sum\limits_{r=1}^{[k/2]}
\sum_{\stackrel{(\{\{g_1, g_2\}, \ldots, 
\{g_{2r-1}, g_{2r}\}\}, \{q_1, \ldots, q_{k-2r}\})}
{{}_{\{g_1, g_2, \ldots, 
g_{2r-1}, g_{2r}, q_1, \ldots, q_{k-2r}\}=\{1, 2, \ldots, k\}}}}
\prod\limits_{s=1}^r
{\bf 1}_{\{i_{g_{{}_{2s-1}}}=~i_{g_{{}_{2s}}}\ne 0\}}\times
$$

\vspace{3mm}
\begin{equation}
\label{after3}
\times
{\bf 1}_{\{j_{g_{{}_{2s-1}}}=~j_{g_{{}_{2s}}}\}}
J'[\phi_{j_{q_1}}\ldots \phi_{j_{q_{k-2r}}}]_{T,t}^{(i_{q_1}\ldots i_{q_{k-2r}})}=
\end{equation}

\vspace{3mm}
$$
=
J[\psi^{(k)}]_{T,t}^{(i_1\ldots i_k)}+
\hbox{\vtop{\offinterlineskip\halign{
\hfil#\hfil\cr
{\rm l.i.m.}\cr
$\stackrel{}{{}_{p\to \infty}}$\cr
}} }
\sum\limits_{r=1}^{[k/2]}
\sum_{\stackrel{(\{\{g_1, g_2\}, \ldots, 
\{g_{2r-1}, g_{2r}\}\}, \{q_1, \ldots, q_{k-2r}\})}
{{}_{\{g_1, g_2, \ldots, 
g_{2r-1}, g_{2r}, q_1, \ldots, q_{k-2r}\}=\{1, 2, \ldots, k\}}}}
\prod\limits_{s=1}^r
{\bf 1}_{\{i_{g_{{}_{2s-1}}}=~i_{g_{{}_{2s}}}\ne 0\}}\times
$$

\vspace{3mm}
\begin{equation}
\label{after7xx}
\times 
J'[K_{p\ldots p}^{g_1\ldots g_{2r}, q_1\ldots q_{k-2r}}]_{T,t}^{(i_{q_1}\ldots i_{q_{k-2r}})}
\end{equation}

\vspace{4mm}
\noindent
w.~p.~{\rm 1}, where
$J[\psi^{(k)}]_{T,t}^{(i_1\ldots i_k)}$ is the iterated Ito stochastic
integral (\ref{afterxgxg1}).

If we prove that w.~p.~1

$$
\sum_{r=1}^{\left[k/2\right]}\frac{1}{2^r}
\sum_{(s_r,\ldots,s_1)\in {\rm A}_{k,r}}
J[\psi^{(k)}]_{T,t}^{s_r,\ldots,s_1}=
$$

\vspace{3mm}
$$
=
\hbox{\vtop{\offinterlineskip\halign{
\hfil#\hfil\cr
{\rm l.i.m.}\cr
$\stackrel{}{{}_{p\to \infty}}$\cr
}} }\sum_{j_1,\ldots,j_k=0}^{p}
C_{j_k \ldots j_1}
\sum\limits_{r=1}^{[k/2]}
\sum_{\stackrel{(\{\{g_1, g_2\}, \ldots, 
\{g_{2r-1}, g_{2r}\}\}, \{q_1, \ldots, q_{k-2r}\})}
{{}_{\{g_1, g_2, \ldots, 
g_{2r-1}, g_{2r}, q_1, \ldots, q_{k-2r}\}=\{1, 2, \ldots, k\}}}}
\prod\limits_{s=1}^r
{\bf 1}_{\{i_{g_{{}_{2s-1}}}=~i_{g_{{}_{2s}}}\ne 0\}}\times
$$

\vspace{3mm}
\begin{equation}
\label{after4}
\times
{\bf 1}_{\{j_{g_{{}_{2s-1}}}=~j_{g_{{}_{2s}}}\}}
J'[\phi_{j_{q_1}}\ldots \phi_{j_{q_{k-2r}}}]_{T,t}^{(i_{q_1}\ldots i_{q_{k-2r}})},
\end{equation}

\vspace{4mm}
\noindent
then (see (\ref{after3}), (\ref{after4}), and Theorem~5)

$$
\hbox{\vtop{\offinterlineskip\halign{
\hfil#\hfil\cr
{\rm l.i.m.}\cr
$\stackrel{}{{}_{p\to \infty}}$\cr
}} }
\sum_{j_1,\ldots,j_k=0}^{p}
C_{j_k \ldots j_1}\prod\limits_{l=1}^k \zeta_{j_l}^{(i_l)}
=
$$

\vspace{1mm}
\begin{equation}
\label{after333}
=
J[\psi^{(k)}]_{T,t}^{(i_1\ldots i_k)}+
\sum_{r=1}^{\left[k/2\right]}\frac{1}{2^r}
\sum_{(s_r,\ldots,s_1)\in {\rm A}_{k,r}}
J[\psi^{(k)}]_{T,t}^{s_r,\ldots,s_1}=
J^{*}[\psi^{(k)}]_{T,t}^{(i_1\ldots i_k)}
\end{equation}

\vspace{3.5mm}
\noindent
w.~p.~1, where notations in (\ref{after333}) are the same as in Theorem~5. 
Thus Theorem~13 will be proved.

From (\ref{after7}) we have that 
the multiple Stratonovich stochastic integral $J[K_{p_1\ldots p_k}]_{T,t}^{(i_1\ldots i_k)}$ 
of mul\-ti\-pli\-ci\-ty $k$ is 
expressed as a sum of some constant value and
multiple Wiener stochastic 
integrals 
$J'[K_{p_1\ldots p_k}]_{T,t}^{(i_1\ldots i_k)}$
and 
$J'[K_{p_1\ldots p_k}^{g_1\ldots g_{2r}, q_1\ldots q_{k-2r}}]_{T,t}^{(i_{q_1}\ldots i_{q_{k-2r}})}$
of multiplicities $k,$ $k-2,$ $k-4$, $\ldots,$ $k-2[k/2]$  
($r=1,2,\ldots,[k/2]$).

The formulas (\ref{after8}), (\ref{after7}) can be considered
as new representations
of the Hu-Meyer formula
for the case of a multidimensional Wiener process 
\cite{Rybakov3000} (also see \cite{bugh1}, \cite{bugh3})
and kernel 
$K_{p_1\ldots p_k}(t_1,\ldots,t_k)$ (see (\ref{afterxx1})).

Note that the equality (\ref{after7}) can be obtained from 
(\ref{quq11}) if we consider (\ref{quq11}) for
$\Phi(t_1,\ldots,t_k)=K_{p_1\ldots p_k}(t_1,\ldots,t_k)$
and without passing to the limit
$\hbox{\vtop{\offinterlineskip\halign{
\hfil#\hfil\cr
{\rm l.i.m.}\cr
$\stackrel{}{{}_{p_1,\ldots,p_k\to \infty}}$\cr
}} }$

For example, for $k=2,3,4,5,6$ we have from (\ref{after8}) w.~p.~1

\vspace{1mm}
\begin{equation}
\label{after32}
\sum_{j_1=0}^{p_1}\sum_{j_2=0}^{p_2}
C_{j_2j_1}\zeta_{j_1}^{(i_1)}\zeta_{j_2}^{(i_2)}=
J'[K_{p_1p_2}]_{T,t}^{(i_1 i_2)}
+\sum_{j_1=0}^{p_1}\sum_{j_2=0}^{p_2}
C_{j_2j_1}
{\bf 1}_{\{i_1=i_2\ne 0\}}{\bf 1}_{\{j_1=j_2\}},
\end{equation}

\vspace{5mm}
$$
\sum_{j_1=0}^{p_1}\sum_{j_2=0}^{p_2}\sum_{j_3=0}^{p_3}
C_{j_3j_2j_1}
\zeta_{j_1}^{(i_1)}\zeta_{j_2}^{(i_2)}\zeta_{j_3}^{(i_3)}=
J'[K_{p_1p_2p_3}]_{T,t}^{(i_1 i_2 i_3)}+
$$
$$
+
\sum_{j_1=0}^{p_1}\sum_{j_2=0}^{p_2}\sum_{j_3=0}^{p_3}
C_{j_3j_2j_1}
\Biggl(
{\bf 1}_{\{i_1=i_2\ne 0\}}
{\bf 1}_{\{j_1=j_2\}}
J'[\phi_{j_3}]^{(i_3)}_{T,t}
+{\bf 1}_{\{i_2=i_3\ne 0\}}
{\bf 1}_{\{j_2=j_3\}}
J'[\phi_{j_1}]^{(i_1)}_{T,t}+\Biggr.
$$
\begin{equation}
\label{after33}
\Biggl.
+{\bf 1}_{\{i_1=i_3\ne 0\}}
{\bf 1}_{\{j_1=j_3\}}
J'[\phi_{j_2}]^{(i_2)}_{T,t}\Biggr),
\end{equation}

\vspace{6mm}

$$
\sum_{j_1=0}^{p_1}\ldots\sum_{j_4=0}^{p_4}
C_{j_4 j_3 j_2 j_1}
\zeta_{j_1}^{(i_1)}\zeta_{j_2}^{(i_2)}\zeta_{j_3}^{(i_3)}
\zeta_{j_4}^{(i_4)}=
J'[K_{p_1p_2p_3p_4}]_{T,t}^{(i_1 i_2 i_3 i_4)}+
$$
$$
+\sum_{j_1=0}^{p_1}\ldots\sum_{j_4=0}^{p_4}
C_{j_4 j_3 j_2 j_1}\Biggl(
\Biggr.
{\bf 1}_{\{i_1=i_2\ne 0\}}
{\bf 1}_{\{j_1=j_2\}}
J'[\phi_{j_3}\phi_{j_4}]^{(i_3 i_4)}_{T,t}
+
$$
$$
+
{\bf 1}_{\{i_1=i_3\ne 0\}}
{\bf 1}_{\{j_1=j_3\}}
J'[\phi_{j_2}\phi_{j_4}]^{(i_2 i_4)}_{T,t}
+
{\bf 1}_{\{i_1=i_4\ne 0\}}
{\bf 1}_{\{j_1=j_4\}}
J'[\phi_{j_2}\phi_{j_3}]^{(i_2 i_3)}_{T,t}
+
$$
$$
+
{\bf 1}_{\{i_2=i_3\ne 0\}}
{\bf 1}_{\{j_2=j_3\}}
J'[\phi_{j_1}\phi_{j_4}]^{(i_1 i_4)}_{T,t}
+
{\bf 1}_{\{i_2=i_4\ne 0\}}
{\bf 1}_{\{j_2=j_4\}}
J'[\phi_{j_1}\phi_{j_3}]^{(i_1 i_3)}_{T,t}
+
$$
$$
+
{\bf 1}_{\{i_3=i_4\ne 0\}}
{\bf 1}_{\{j_3=j_4\}}
J'[\phi_{j_1}\phi_{j_2}]^{(i_1 i_2)}_{T,t}
+
$$
$$
+
{\bf 1}_{\{i_1=i_2\ne 0\}}
{\bf 1}_{\{j_1=j_2\}}
{\bf 1}_{\{i_3=i_4\ne 0\}}
{\bf 1}_{\{j_3=j_4\}}
+
{\bf 1}_{\{i_1=i_3\ne 0\}}
{\bf 1}_{\{j_1=j_3\}}
{\bf 1}_{\{i_2=i_4\ne 0\}}
{\bf 1}_{\{j_2=j_4\}}+
$$
\begin{equation}
\label{after34}
+\Biggl.
{\bf 1}_{\{i_1=i_4\ne 0\}}
{\bf 1}_{\{j_1=j_4\}}
{\bf 1}_{\{i_2=i_3\ne 0\}}
{\bf 1}_{\{j_2=j_3\}}\Biggr),
\end{equation}

\vspace{7mm}
$$
\sum_{j_1=0}^{p_1}\ldots\sum_{j_5=0}^{p_5}
C_{j_5 j_4 j_3 j_2 j_1}
\zeta_{j_1}^{(i_1)}\zeta_{j_2}^{(i_2)}\zeta_{j_3}^{(i_3)}
\zeta_{j_4}^{(i_4)}\zeta_{j_5}^{(i_5)}
=
J'[K_{p_1p_2p_3p_4p_5}]_{T,t}^{(i_1 i_2 i_3 i_4 i_5)}
+
$$
$$
+\sum_{j_1=0}^{p_1}\ldots\sum_{j_5=0}^{p_5}
C_{j_5 j_4 j_3 j_2 j_1}\Biggl(
{\bf 1}_{\{i_1=i_2\ne 0\}}
{\bf 1}_{\{j_1=j_2\}}
J'[\phi_{j_3}\phi_{j_4}
\phi_{j_5}]^{(i_3i_4i_5)}_{T,t}+
$$
$$
+
{\bf 1}_{\{i_1=i_3\ne 0\}}
{\bf 1}_{\{j_1=j_3\}}
J'[\phi_{j_2}\phi_{j_4}\phi_{j_5}]^{(i_2i_4i_5)}_{T,t}+
{\bf 1}_{\{i_1=i_4\ne 0\}}
{\bf 1}_{\{j_1=j_4\}}
J'[\phi_{j_2}\phi_{j_3}\phi_{j_5}]^{(i_2i_3i_5)}_{T,t}
+
$$
$$
+
{\bf 1}_{\{i_1=i_5\ne 0\}}
{\bf 1}_{\{j_1=j_5\}}
J'[\phi_{j_2}\phi_{j_3}\phi_{j_4}]^{(i_2i_3i_4)}_{T,t}
+
{\bf 1}_{\{i_2=i_3\ne 0\}}
{\bf 1}_{\{j_2=j_3\}}
J'[\phi_{j_1}\phi_{j_4}\phi_{j_5}]^{(i_1i_4i_5)}_{T,t}
+
$$
$$
+
{\bf 1}_{\{i_2=i_4\ne 0\}}
{\bf 1}_{\{j_2=j_4\}}
J'[\phi_{j_1}\phi_{j_3}\phi_{j_5}]^{(i_1i_3i_5)}_{T,t}
+
{\bf 1}_{\{i_2=i_5\ne 0\}}
{\bf 1}_{\{j_2=j_5\}}
J'[\phi_{j_1}\phi_{j_3}\phi_{j_4}]^{(i_1i_3i_4)}_{T,t}
+
$$
$$
+
{\bf 1}_{\{i_3=i_4\ne 0\}}
{\bf 1}_{\{j_3=j_4\}}
J'[\phi_{j_1}\phi_{j_2}\phi_{j_5}]^{(i_1i_2i_5)}_{T,t}
+
{\bf 1}_{\{i_3=i_5\ne 0\}}
{\bf 1}_{\{j_3=j_5\}}
J'[\phi_{j_1}\phi_{j_2}\phi_{j_4}]^{(i_1i_2i_4)}_{T,t}
+
$$
$$
+{\bf 1}_{\{i_4=i_5\ne 0\}}
{\bf 1}_{\{j_4=j_5\}}
J'[\phi_{j_1}\phi_{j_2}\phi_{j_3}]^{(i_1i_2i_3)}_{T,t}
+
$$
$$
+
{\bf 1}_{\{i_1=i_2\ne 0\}}
{\bf 1}_{\{j_1=j_2\}}
{\bf 1}_{\{i_3=i_4\ne 0\}}
{\bf 1}_{\{j_3=j_4\}}J'[\phi_{j_5}]^{(i_5)}_{T,t}+
$$
$$
+
{\bf 1}_{\{i_1=i_2\ne 0\}}
{\bf 1}_{\{j_1=j_2\}}
{\bf 1}_{\{i_3=i_5\ne 0\}}
{\bf 1}_{\{j_3=j_5\}}
J'[\phi_{j_4}]^{(i_4)}_{T,t}
+
$$
$$
+
{\bf 1}_{\{i_1=i_2\ne 0\}}
{\bf 1}_{\{j_1=j_2\}}
{\bf 1}_{\{i_4=i_5\ne 0\}}
{\bf 1}_{\{j_4=j_5\}}
J'[\phi_{j_3}]^{(i_3)}_{T,t}+
$$
$$
+
{\bf 1}_{\{i_1=i_3\ne 0\}}
{\bf 1}_{\{j_1=j_3\}}
{\bf 1}_{\{i_2=i_4\ne 0\}}
{\bf 1}_{\{j_2=j_4\}}
J'[\phi_{j_5}]^{(i_5)}_{T,t}
+
$$
$$
+
{\bf 1}_{\{i_1=i_3\ne 0\}}
{\bf 1}_{\{j_1=j_3\}}
{\bf 1}_{\{i_2=i_5\ne 0\}}
{\bf 1}_{\{j_2=j_5\}}
J'[\phi_{j_4}]^{(i_4)}_{T,t}+
$$
$$
+
{\bf 1}_{\{i_1=i_3\ne 0\}}
{\bf 1}_{\{j_1=j_3\}}
{\bf 1}_{\{i_4=i_5\ne 0\}}
{\bf 1}_{\{j_4=j_5\}}
J'[\phi_{j_2}]^{(i_2)}_{T,t}+
$$
$$
+
{\bf 1}_{\{i_1=i_4\ne 0\}}
{\bf 1}_{\{j_1=j_4\}}
{\bf 1}_{\{i_2=i_3\ne 0\}}
{\bf 1}_{\{j_2=j_3\}}
J'[\phi_{j_5}]^{(i_5)}_{T,t}
+
$$
$$
+
{\bf 1}_{\{i_1=i_4\ne 0\}}
{\bf 1}_{\{j_1=j_4\}}
{\bf 1}_{\{i_2=i_5\ne 0\}}
{\bf 1}_{\{j_2=j_5\}}
J'[\phi_{j_3}]^{(i_3)}_{T,t}
+
$$
$$
+
{\bf 1}_{\{i_1=i_4\ne 0\}}
{\bf 1}_{\{j_1=j_4\}}
{\bf 1}_{\{i_3=i_5\ne 0\}}
{\bf 1}_{\{j_3=j_5\}}
J'[\phi_{j_2}]^{(i_2)}_{T,t}
+
$$
$$
+
{\bf 1}_{\{i_1=i_5\ne 0\}}
{\bf 1}_{\{j_1=j_5\}}
{\bf 1}_{\{i_2=i_3\ne 0\}}
{\bf 1}_{\{j_2=j_3\}}
J'[\phi_{j_4}]^{(i_4)}_{T,t}
+
$$
$$
+
{\bf 1}_{\{i_1=i_5\ne 0\}}
{\bf 1}_{\{j_1=j_5\}}
{\bf 1}_{\{i_2=i_4\ne 0\}}
{\bf 1}_{\{j_2=j_4\}}
J'[\phi_{j_3}]^{(i_3)}_{T,t}+
$$
$$
+
{\bf 1}_{\{i_1=i_5\ne 0\}}
{\bf 1}_{\{j_1=j_5\}}
{\bf 1}_{\{i_3=i_4\ne 0\}}
{\bf 1}_{\{j_3=j_4\}}
J'[\phi_{j_2}]^{(i_2)}_{T,t}+
$$
$$
+
{\bf 1}_{\{i_2=i_3\ne 0\}}
{\bf 1}_{\{j_2=j_3\}}
{\bf 1}_{\{i_4=i_5\ne 0\}}
{\bf 1}_{\{j_4=j_5\}}
J'[\phi_{j_1}]^{(i_1)}_{T,t}+
$$
$$
+
{\bf 1}_{\{i_2=i_4\ne 0\}}
{\bf 1}_{\{j_2=j_4\}}
{\bf 1}_{\{i_3=i_5\ne 0\}}
{\bf 1}_{\{j_3=j_5\}}
J'[\phi_{j_1}]^{(i_1)}_{T,t}+
$$
\begin{equation}
\label{after35}
+\Biggl.
{\bf 1}_{\{i_2=i_5\ne 0\}}
{\bf 1}_{\{j_2=j_5\}}
{\bf 1}_{\{i_3=i_4\ne 0\}}
{\bf 1}_{\{j_3=j_4\}}
J'[\phi_{j_1}]^{(i_1)}_{T,t}
\Biggr),
\end{equation}

\vspace{8mm}

$$
\sum_{j_1=0}^{p_1}\ldots\sum_{j_6=0}^{p_6}
C_{j_6 j_5 j_4 j_3 j_2 j_1}
\zeta_{j_1}^{(i_1)}\zeta_{j_2}^{(i_2)}\zeta_{j_3}^{(i_3)}
\zeta_{j_4}^{(i_4)}\zeta_{j_5}^{(i_5)}\zeta_{j_6}^{(i_6)}
=
J'[K_{p_1p_2p_3p_4p_5p_6}]_{T,t}^{(i_1 i_2 i_3 i_4 i_5 i_6)}+
$$
$$
+\sum_{j_1=0}^{p_1}\ldots\sum_{j_6=0}^{p_6}
C_{j_6j_5j_4j_3j_2j_1}\Biggl(
{\bf 1}_{\{i_1=i_6\ne 0\}}
{\bf 1}_{\{j_1=j_6\}}
J'[\phi_{j_2}\phi_{j_3}\phi_{j_4}\phi_{j_5}]^{(i_2i_3i_4i_5)}_{T,t}
+
$$
$$
+
{\bf 1}_{\{i_2=i_6\ne 0\}}
{\bf 1}_{\{j_2=j_6\}}
J'[\phi_{j_1}\phi_{j_3}\phi_{j_4}\phi_{j_5}]^{(i_1i_3i_4i_5)}_{T,t}
+
{\bf 1}_{\{i_3=i_6\ne 0\}}
{\bf 1}_{\{j_3=j_6\}}
J'[\phi_{j_1}\phi_{j_2}\phi_{j_4}\phi_{j_5}]^{(i_1i_2i_4i_5)}_{T,t}
+
$$
$$
+
{\bf 1}_{\{i_4=i_6\ne 0\}}
{\bf 1}_{\{j_4=j_6\}}
J'[\phi_{j_1}\phi_{j_2}\phi_{j_3}\phi_{j_5}]^{(i_1i_2i_3i_5)}_{T,t}
+
{\bf 1}_{\{i_5=i_6\ne 0\}}
{\bf 1}_{\{j_5=j_6\}}
J'[\phi_{j_1}\phi_{j_2}\phi_{j_3}\phi_{j_4}]^{(i_1i_2i_3i_4)}_{T,t}
+
$$
$$
+{\bf 1}_{\{i_1=i_2\ne 0\}}
{\bf 1}_{\{j_1=j_2\}}
J'[\phi_{j_3}\phi_{j_4}\phi_{j_5}\phi_{j_6}]^{(i_3i_4i_5i_6)}_{T,t}
+
{\bf 1}_{\{i_1=i_3\ne 0\}}
{\bf 1}_{\{j_1=j_3\}}
J'[\phi_{j_2}\phi_{j_4}\phi_{j_5}\phi_{j_6}]^{(i_2i_4i_5i_6)}_{T,t}
+
$$
$$
+{\bf 1}_{\{i_1=i_4\ne 0\}}
{\bf 1}_{\{j_1=j_4\}}
J'[\phi_{j_2}\phi_{j_3}\phi_{j_5}\phi_{j_6}]^{(i_2i_3i_5i_6)}_{T,t}
+
{\bf 1}_{\{i_1=i_5\ne 0\}}
{\bf 1}_{\{j_1=j_5\}}
J'[\phi_{j_2}\phi_{j_3}\phi_{j_4}\phi_{j_6}]^{(i_2i_3i_4i_6)}_{T,t}
+
$$
$$
+{\bf 1}_{\{i_2=i_3\ne 0\}}
{\bf 1}_{\{j_2=j_3\}}
J'[\phi_{j_1}\phi_{j_4}\phi_{j_5}\phi_{j_6}]^{(i_1i_4i_5i_6)}_{T,t}
+
{\bf 1}_{\{i_2=i_4\ne 0\}}
{\bf 1}_{\{j_2=j_4\}}
J'[\phi_{j_1}\phi_{j_3}\phi_{j_5}\phi_{j_6}]^{(i_1i_3i_5i_6)}_{T,t}
+
$$
$$
+
{\bf 1}_{\{i_2=i_5\ne 0\}}
{\bf 1}_{\{j_2=j_5\}}
J'[\phi_{j_1}\phi_{j_3}\phi_{j_4}\phi_{j_6}]^{(i_1i_3i_4i_6)}_{T,t}
+
{\bf 1}_{\{i_3=i_4\ne 0\}}
{\bf 1}_{\{j_3=j_4\}}
J'[\phi_{j_1}\phi_{j_2}\phi_{j_5}\phi_{j_6}]^{(i_1i_2i_5i_6)}_{T,t}
+
$$
$$
+
{\bf 1}_{\{i_3=i_5\ne 0\}}
{\bf 1}_{\{j_3=j_5\}}
J'[\phi_{j_1}\phi_{j_2}\phi_{j_4}\phi_{j_6}]^{(i_1i_2i_4i_6)}_{T,t}
+
{\bf 1}_{\{i_4=i_5\ne 0\}}
{\bf 1}_{\{j_4=j_5\}}
J'[\phi_{j_1}\phi_{j_2}\phi_{j_3}\phi_{j_6}]^{(i_1i_2i_3i_6)}_{T,t}
+
$$
$$
+
{\bf 1}_{\{i_1=i_2\ne 0\}}
{\bf 1}_{\{j_1=j_2\}}
{\bf 1}_{\{i_3=i_4\ne 0\}}
{\bf 1}_{\{j_3=j_4\}}
J'[\phi_{j_5}\phi_{j_6}]^{(i_5i_6)}_{T,t}+
$$
$$
+
{\bf 1}_{\{i_1=i_2\ne 0\}}
{\bf 1}_{\{j_1=j_2\}}
{\bf 1}_{\{i_3=i_5\ne 0\}}
{\bf 1}_{\{j_3=j_5\}}
J'[\phi_{j_4}\phi_{j_6}]^{(i_4i_6)}_{T,t}
+
$$
$$
+
{\bf 1}_{\{i_1=i_2\ne 0\}}
{\bf 1}_{\{j_1=j_2\}}
{\bf 1}_{\{i_4=i_5\ne 0\}}
{\bf 1}_{\{j_4=j_5\}}
J'[\phi_{j_3}\phi_{j_6}]^{(i_3i_6)}_{T,t}
+
$$
$$
+
{\bf 1}_{\{i_1=i_3\ne 0\}}
{\bf 1}_{\{j_1=j_3\}}
{\bf 1}_{\{i_2=i_4\ne 0\}}
{\bf 1}_{\{j_2=j_4\}}
J'[\phi_{j_5}\phi_{j_6}]^{(i_5i_6)}_{T,t}
+
$$
$$
+
{\bf 1}_{\{i_1=i_3\ne 0\}}
{\bf 1}_{\{j_1=j_3\}}
{\bf 1}_{\{i_2=i_5\ne 0\}}
{\bf 1}_{\{j_2=j_5\}}
J'[\phi_{j_4}\phi_{j_6}]^{(i_4i_6)}_{T,t}
+
$$
$$
+{\bf 1}_{\{i_1=i_3\ne 0\}}
{\bf 1}_{\{j_1=j_3\}}
{\bf 1}_{\{i_4=i_5\ne 0\}}
{\bf 1}_{\{j_4=j_5\}}
J'[\phi_{j_2}\phi_{j_6}]^{(i_2i_6)}_{T,t}
+
$$
$$
+
{\bf 1}_{\{i_1=i_4\ne 0\}}
{\bf 1}_{\{j_1=j_4\}}
{\bf 1}_{\{i_2=i_3\ne 0\}}
{\bf 1}_{\{j_2=j_3\}}
J'[\phi_{j_5}\phi_{j_6}]^{(i_5i_6)}_{T,t}
+
$$
$$
+
{\bf 1}_{\{i_1=i_4\ne 0\}}
{\bf 1}_{\{j_1=j_4\}}
{\bf 1}_{\{i_2=i_5\ne 0\}}
{\bf 1}_{\{j_2=j_5\}}
J'[\phi_{j_3}\phi_{j_6}]^{(i_3i_6)}_{T,t}
+
$$
$$
+
{\bf 1}_{\{i_1=i_4\ne 0\}}
{\bf 1}_{\{j_1=j_4\}}
{\bf 1}_{\{i_3=i_5\ne 0\}}
{\bf 1}_{\{j_3=j_5\}}
J'[\phi_{j_2}\phi_{j_6}]^{(i_2i_6)}_{T,t}
+
$$
$$
+
{\bf 1}_{\{i_1=i_5\ne 0\}}
{\bf 1}_{\{j_1=j_5\}}
{\bf 1}_{\{i_2=i_3\ne 0\}}
{\bf 1}_{\{j_2=j_3\}}
J'[\phi_{j_4}\phi_{j_6}]^{(i_4i_6)}_{T,t}
+
$$
$$
+
{\bf 1}_{\{i_1=i_5\ne 0\}}
{\bf 1}_{\{j_1=j_5\}}
{\bf 1}_{\{i_2=i_4\ne 0\}}
{\bf 1}_{\{j_2=j_4\}}
J'[\phi_{j_3}\phi_{j_6}]^{(i_3i_6)}_{T,t}
+
$$
$$
+
{\bf 1}_{\{i_1=i_5\ne 0\}}
{\bf 1}_{\{j_1=j_5\}}
{\bf 1}_{\{i_3=i_4\ne 0\}}
{\bf 1}_{\{j_3=j_4\}}
J'[\phi_{j_2}\phi_{j_6}]^{(i_2i_6)}_{T,t}
+
$$
$$
+
{\bf 1}_{\{i_2=i_3\ne 0\}}
{\bf 1}_{\{j_2=j_3\}}
{\bf 1}_{\{i_4=i_5\ne 0\}}
{\bf 1}_{\{j_4=j_5\}}
J'[\phi_{j_1}\phi_{j_6}]^{(i_1i_6)}_{T,t}
+
$$
$$
+{\bf 1}_{\{i_2=i_4\ne 0\}}
{\bf 1}_{\{j_2=j_4\}}
{\bf 1}_{\{i_3=i_5\ne 0\}}
{\bf 1}_{\{j_3=j_5\}}
J'[\phi_{j_1}\phi_{j_6}]^{(i_1i_6)}_{T,t}
+
$$
$$
+
{\bf 1}_{\{i_2=i_5\ne 0\}}
{\bf 1}_{\{j_2=j_5\}}
{\bf 1}_{\{i_3=i_4\ne 0\}}
{\bf 1}_{\{j_3=j_4\}}
J'[\phi_{j_1}\phi_{j_6}]^{(i_1i_6)}_{T,t}
+
$$
$$
+{\bf 1}_{\{i_6=i_1\ne 0\}}
{\bf 1}_{\{j_6=j_1\}}
{\bf 1}_{\{i_3=i_4\ne 0\}}
{\bf 1}_{\{j_3=j_4\}}
J'[\phi_{j_2}\phi_{j_5}]^{(i_2i_5)}_{T,t}
+
$$
$$
+
{\bf 1}_{\{i_6=i_1\ne 0\}}
{\bf 1}_{\{j_6=j_1\}}
{\bf 1}_{\{i_3=i_5\ne 0\}}
{\bf 1}_{\{j_3=j_5\}}
J'[\phi_{j_2}\phi_{j_4}]^{(i_2i_4)}_{T,t}
+
$$
$$
+{\bf 1}_{\{i_6=i_1\ne 0\}}
{\bf 1}_{\{j_6=j_1\}}
{\bf 1}_{\{i_2=i_5\ne 0\}}
{\bf 1}_{\{j_2=j_5\}}
J'[\phi_{j_3}\phi_{j_4}]^{(i_3i_4)}_{T,t}
+
$$
$$
+
{\bf 1}_{\{i_6=i_1\ne 0\}}
{\bf 1}_{\{j_6=j_1\}}
{\bf 1}_{\{i_2=i_4\ne 0\}}
{\bf 1}_{\{j_2=j_4\}}
J'[\phi_{j_3}\phi_{j_5}]^{(i_3i_5)}_{T,t}
+
$$
$$
+{\bf 1}_{\{i_6=i_1\ne 0\}}
{\bf 1}_{\{j_6=j_1\}}
{\bf 1}_{\{i_4=i_5\ne 0\}}
{\bf 1}_{\{j_4=j_5\}}
J'[\phi_{j_2}\phi_{j_3}]^{(i_2i_3)}_{T,t}
+
$$
$$
+
{\bf 1}_{\{i_6=i_1\ne 0\}}
{\bf 1}_{\{j_6=j_1\}}
{\bf 1}_{\{i_2=i_3\ne 0\}}
{\bf 1}_{\{j_2=j_3\}}
J'[\phi_{j_4}\phi_{j_5}]^{(i_4i_5)}_{T,t}
+
$$
$$
+{\bf 1}_{\{i_6=i_2\ne 0\}}
{\bf 1}_{\{j_6=j_2\}}
{\bf 1}_{\{i_3=i_5\ne 0\}}
{\bf 1}_{\{j_3=j_5\}}
J'[\phi_{j_1}\phi_{j_4}]^{(i_1i_4)}_{T,t}
+
$$
$$
+
{\bf 1}_{\{i_6=i_2\ne 0\}}
{\bf 1}_{\{j_6=j_2\}}
{\bf 1}_{\{i_4=i_5\ne 0\}}
{\bf 1}_{\{j_4=j_5\}}
J'[\phi_{j_1}\phi_{j_3}]^{(i_1i_3)}_{T,t}
+
$$
$$
+{\bf 1}_{\{i_6=i_2\ne 0\}}
{\bf 1}_{\{j_6=j_2\}}
{\bf 1}_{\{i_3=i_4\ne 0\}}
{\bf 1}_{\{j_3=j_4\}}
J'[\phi_{j_1}\phi_{j_5}]^{(i_1i_5)}_{T,t}
+
$$
$$
+
{\bf 1}_{\{i_6=i_2\ne 0\}}
{\bf 1}_{\{j_6=j_2\}}
{\bf 1}_{\{i_1=i_5\ne 0\}}
{\bf 1}_{\{j_1=j_5\}}
J'[\phi_{j_3}\phi_{j_4}]^{(i_3i_4)}_{T,t}
+
$$
$$
+{\bf 1}_{\{i_6=i_2\ne 0\}}
{\bf 1}_{\{j_6=j_2\}}
{\bf 1}_{\{i_1=i_4\ne 0\}}
{\bf 1}_{\{j_1=j_4\}}
J'[\phi_{j_3}\phi_{j_5}]^{(i_3i_5)}_{T,t}
+
$$
$$
+
{\bf 1}_{\{i_6=i_2\ne 0\}}
{\bf 1}_{\{j_6=j_2\}}
{\bf 1}_{\{i_1=i_3\ne 0\}}
{\bf 1}_{\{j_1=j_3\}}
J'[\phi_{j_4}\phi_{j_5}]^{(i_4i_5)}_{T,t}
+
$$
$$
+{\bf 1}_{\{i_6=i_3\ne 0\}}
{\bf 1}_{\{j_6=j_3\}}
{\bf 1}_{\{i_2=i_5\ne 0\}}
{\bf 1}_{\{j_2=j_5\}}
J'[\phi_{j_1}\phi_{j_4}]^{(i_1i_4)}_{T,t}
+
$$
$$
+
{\bf 1}_{\{i_6=i_3\ne 0\}}
{\bf 1}_{\{j_6=j_3\}}
{\bf 1}_{\{i_4=i_5\ne 0\}}
{\bf 1}_{\{j_4=j_5\}}
J'[\phi_{j_1}\phi_{j_2}]^{(i_1i_2)}_{T,t}
+
$$
$$
+{\bf 1}_{\{i_6=i_3\ne 0\}}
{\bf 1}_{\{j_6=j_3\}}
{\bf 1}_{\{i_2=i_4\ne 0\}}
{\bf 1}_{\{j_2=j_4\}}
J'[\phi_{j_1}\phi_{j_5}]^{(i_1i_5)}_{T,t}
+
$$
$$
+
{\bf 1}_{\{i_6=i_3\ne 0\}}
{\bf 1}_{\{j_6=j_3\}}
{\bf 1}_{\{i_1=i_5\ne 0\}}
{\bf 1}_{\{j_1=j_5\}}
J'[\phi_{j_2}\phi_{j_4}]^{(i_2i_4)}_{T,t}
+
$$
$$
+{\bf 1}_{\{i_6=i_3\ne 0\}}
{\bf 1}_{\{j_6=j_3\}}
{\bf 1}_{\{i_1=i_4\ne 0\}}
{\bf 1}_{\{j_1=j_4\}}
J'[\phi_{j_2}\phi_{j_5}]^{(i_2i_5)}_{T,t}
+
$$
$$
+
{\bf 1}_{\{i_6=i_3\ne 0\}}
{\bf 1}_{\{j_6=j_3\}}
{\bf 1}_{\{i_1=i_2\ne 0\}}
{\bf 1}_{\{j_1=j_2\}}
J'[\phi_{j_4}\phi_{j_5}]^{(i_4i_5)}_{T,t}
+
$$
$$
+{\bf 1}_{\{i_6=i_4\ne 0\}}
{\bf 1}_{\{j_6=j_4\}}
{\bf 1}_{\{i_3=i_5\ne 0\}}
{\bf 1}_{\{j_3=j_5\}}
J'[\phi_{j_1}\phi_{j_2}]^{(i_1i_2)}_{T,t}
+
$$
$$
+
{\bf 1}_{\{i_6=i_4\ne 0\}}
{\bf 1}_{\{j_6=j_4\}}
{\bf 1}_{\{i_2=i_5\ne 0\}}
{\bf 1}_{\{j_2=j_5\}}
J'[\phi_{j_1}\phi_{j_3}]^{(i_1i_3)}_{T,t}
+
$$
$$
+{\bf 1}_{\{i_6=i_4\ne 0\}}
{\bf 1}_{\{j_6=j_4\}}
{\bf 1}_{\{i_2=i_3\ne 0\}}
{\bf 1}_{\{j_2=j_3\}}
J'[\phi_{j_1}\phi_{j_5}]^{(i_1i_5)}_{T,t}
+
$$
$$
+
{\bf 1}_{\{i_6=i_4\ne 0\}}
{\bf 1}_{\{j_6=j_4\}}
{\bf 1}_{\{i_1=i_5\ne 0\}}
{\bf 1}_{\{j_1=j_5\}}
J'[\phi_{j_2}\phi_{j_3}]^{(i_2i_3)}_{T,t}
+
$$
$$
+{\bf 1}_{\{i_6=i_4\ne 0\}}
{\bf 1}_{\{j_6=j_4\}}
{\bf 1}_{\{i_1=i_3\ne 0\}}
{\bf 1}_{\{j_1=j_3\}}
J'[\phi_{j_2}\phi_{j_5}]^{(i_2i_5)}_{T,t}
+
$$
$$
+
{\bf 1}_{\{i_6=i_4\ne 0\}}
{\bf 1}_{\{j_6=j_4\}}
{\bf 1}_{\{i_1=i_2\ne 0\}}
{\bf 1}_{\{j_1=j_2\}}
J'[\phi_{j_3}\phi_{j_5}]^{(i_3i_5)}_{T,t}
+
$$
$$
+{\bf 1}_{\{i_6=i_5\ne 0\}}
{\bf 1}_{\{j_6=j_5\}}
{\bf 1}_{\{i_3=i_4\ne 0\}}
{\bf 1}_{\{j_3=j_4\}}
J'[\phi_{j_1}\phi_{j_2}]^{(i_1i_2)}_{T,t}
+
$$
$$
+
{\bf 1}_{\{i_6=i_5\ne 0\}}
{\bf 1}_{\{j_6=j_5\}}
{\bf 1}_{\{i_2=i_4\ne 0\}}
{\bf 1}_{\{j_2=j_4\}}
J'[\phi_{j_1}\phi_{j_3}]^{(i_1i_3)}_{T,t}
+
$$
$$
+{\bf 1}_{\{i_6=i_5\ne 0\}}
{\bf 1}_{\{j_6=j_5\}}
{\bf 1}_{\{i_2=i_3\ne 0\}}
{\bf 1}_{\{j_2=j_3\}}
J'[\phi_{j_1}\phi_{j_4}]^{(i_1i_4)}_{T,t}
+
$$
$$
+
{\bf 1}_{\{i_6=i_5\ne 0\}}
{\bf 1}_{\{j_6=j_5\}}
{\bf 1}_{\{i_1=i_4\ne 0\}}
{\bf 1}_{\{j_1=j_4\}}
J'[\phi_{j_2}\phi_{j_3}]^{(i_2i_3)}_{T,t}
+
$$
$$
+{\bf 1}_{\{i_6=i_5\ne 0\}}
{\bf 1}_{\{j_6=j_5\}}
{\bf 1}_{\{i_1=i_3\ne 0\}}
{\bf 1}_{\{j_1=j_3\}}
J'[\phi_{j_2}\phi_{j_4}]^{(i_2i_4)}_{T,t}
+
$$
$$
+
{\bf 1}_{\{i_6=i_5\ne 0\}}
{\bf 1}_{\{j_6=j_5\}}
{\bf 1}_{\{i_1=i_2\ne 0\}}
{\bf 1}_{\{j_1=j_2\}}
J'[\phi_{j_3}\phi_{j_4}]^{(i_3i_4)}_{T,t}
+
$$
$$
+
{\bf 1}_{\{i_6=i_1\ne 0\}}
{\bf 1}_{\{j_6=j_1\}}
{\bf 1}_{\{i_2=i_5\ne 0\}}
{\bf 1}_{\{j_2=j_5\}}
{\bf 1}_{\{i_3=i_4\ne 0\}}
{\bf 1}_{\{j_3=j_4\}}+
$$
$$
+
{\bf 1}_{\{i_6=i_1\ne 0\}}
{\bf 1}_{\{j_6=j_1\}}
{\bf 1}_{\{i_2=i_4\ne 0\}}
{\bf 1}_{\{j_2=j_4\}}
{\bf 1}_{\{i_3=i_5\ne 0\}}
{\bf 1}_{\{j_3=j_5\}}+
$$
$$
+
{\bf 1}_{\{i_6=i_1\ne 0\}}
{\bf 1}_{\{j_6=j_1\}}
{\bf 1}_{\{i_2=i_3\ne 0\}}
{\bf 1}_{\{j_2=j_3\}}
{\bf 1}_{\{i_4=i_5\ne 0\}}
{\bf 1}_{\{j_4=j_5\}}+
$$
$$
+
{\bf 1}_{\{i_6=i_2\ne 0\}}
{\bf 1}_{\{j_6=j_2\}}
{\bf 1}_{\{i_1=i_5\ne 0\}}
{\bf 1}_{\{j_1=j_5\}}
{\bf 1}_{\{i_3=i_4\ne 0\}}
{\bf 1}_{\{j_3=j_4\}}+
$$
$$
+
{\bf 1}_{\{i_6=i_2\ne 0\}}
{\bf 1}_{\{j_6=j_2\}}
{\bf 1}_{\{i_1=i_4\ne 0\}}
{\bf 1}_{\{j_1=j_4\}}
{\bf 1}_{\{i_3=i_5\ne 0\}}
{\bf 1}_{\{j_3=j_5\}}+
$$
$$
+
{\bf 1}_{\{i_6=i_2\ne 0\}}
{\bf 1}_{\{j_6=j_2\}}
{\bf 1}_{\{i_1=i_3\ne 0\}}
{\bf 1}_{\{j_1=j_3\}}
{\bf 1}_{\{i_4=i_5\ne 0\}}
{\bf 1}_{\{j_4=j_5\}}+
$$
$$
+
{\bf 1}_{\{i_6=i_3\ne 0\}}
{\bf 1}_{\{j_6=j_3\}}
{\bf 1}_{\{i_1=i_5\ne 0\}}
{\bf 1}_{\{j_1=j_5\}}
{\bf 1}_{\{i_2=i_4\ne 0\}}
{\bf 1}_{\{j_2=j_4\}}+
$$
$$
+
{\bf 1}_{\{i_6=i_3\ne 0\}}
{\bf 1}_{\{j_6=j_3\}}
{\bf 1}_{\{i_1=i_4\ne 0\}}
{\bf 1}_{\{j_1=j_4\}}
{\bf 1}_{\{i_2=i_5\ne 0\}}
{\bf 1}_{\{j_2=j_5\}}+
$$
$$
+
{\bf 1}_{\{i_3=i_6\ne 0\}}
{\bf 1}_{\{j_3=j_6\}}
{\bf 1}_{\{i_1=i_2\ne 0\}}
{\bf 1}_{\{j_1=j_2\}}
{\bf 1}_{\{i_4=i_5\ne 0\}}
{\bf 1}_{\{j_4=j_5\}}+
$$
$$
+
{\bf 1}_{\{i_6=i_4\ne 0\}}
{\bf 1}_{\{j_6=j_4\}}
{\bf 1}_{\{i_1=i_5\ne 0\}}
{\bf 1}_{\{j_1=j_5\}}
{\bf 1}_{\{i_2=i_3\ne 0\}}
{\bf 1}_{\{j_2=j_3\}}+
$$
$$
+
{\bf 1}_{\{i_6=i_4\ne 0\}}
{\bf 1}_{\{j_6=j_4\}}
{\bf 1}_{\{i_1=i_3\ne 0\}}
{\bf 1}_{\{j_1=j_3\}}
{\bf 1}_{\{i_2=i_5\ne 0\}}
{\bf 1}_{\{j_2=j_5\}}+
$$
$$
+
{\bf 1}_{\{i_6=i_4\ne 0\}}
{\bf 1}_{\{j_6=j_4\}}
{\bf 1}_{\{i_1=i_2\ne 0\}}
{\bf 1}_{\{j_1=j_2\}}
{\bf 1}_{\{i_3=i_5\ne 0\}}
{\bf 1}_{\{j_3=j_5\}}+
$$
$$
+
{\bf 1}_{\{i_6=i_5\ne 0\}}
{\bf 1}_{\{j_6=j_5\}}
{\bf 1}_{\{i_1=i_4\ne 0\}}
{\bf 1}_{\{j_1=j_4\}}
{\bf 1}_{\{i_2=i_3\ne 0\}}
{\bf 1}_{\{j_2=j_3\}}+
$$
$$
+
{\bf 1}_{\{i_6=i_5\ne 0\}}
{\bf 1}_{\{j_6=j_5\}}
{\bf 1}_{\{i_1=i_2\ne 0\}}
{\bf 1}_{\{j_1=j_2\}}
{\bf 1}_{\{i_3=i_4\ne 0\}}
{\bf 1}_{\{j_3=j_4\}}+
$$
\begin{equation}
\label{after36}
\Biggl.+
{\bf 1}_{\{i_6=i_5\ne 0\}}
{\bf 1}_{\{j_6=j_5\}}
{\bf 1}_{\{i_1=i_3\ne 0\}}
{\bf 1}_{\{j_1=j_3\}}
{\bf 1}_{\{i_2=i_4\ne 0\}}
{\bf 1}_{\{j_2=j_4\}}\Biggr).
\end{equation}

\vspace{5mm}

Note that the relation (\ref{after34})
can be written in the following form

\vspace{2mm}
$$
\sum_{j_1=0}^{p_1}\ldots \sum_{j_4=0}^{p_4}
C_{j_4 j_3 j_2 j_1}\zeta_{j_1}^{(i_1)}\zeta_{j_2}^{(i_2)}\zeta_{j_3}^{(i_3)}\zeta_{j_4}^{(i_4)}=
\sum_{j_1=0}^{p_1}\ldots \sum_{j_4=0}^{p_4}C_{j_4 j_3 j_2 j_1}
J'[\phi_{j_1}\phi_{j_2}\phi_{j_3}\phi_{j_4}]_{T,t}^{(i_1i_2i_3i_4)}+
$$

\vspace{2mm}
$$
+{\bf 1}_{\{i_1=i_2\ne 0\}}
\sum_{j_3=0}^{p_3}\sum_{j_4=0}^{p_4}\left(\sum_{j_1=0}^{\min\{p_1,p_2\}} C_{j_4 j_3 j_1 j_1}\right)
J'[\phi_{j_3}\phi_{j_4}]_{T,t}^{(i_3i_4)}+
$$

\vspace{2mm}
$$
+{\bf 1}_{\{i_1=i_3\ne 0\}}
\sum_{j_2=0}^{p_2}\sum_{j_4=0}^{p_4}\left(\sum_{j_3=0}^{\min\{p_1,p_3\}} C_{j_4 j_3 j_2 j_3}\right)
J'[\phi_{j_2}\phi_{j_4}]_{T,t}^{(i_2i_4)}+
$$

\vspace{2mm}
$$
+{\bf 1}_{\{i_1=i_4\ne 0\}}
\sum_{j_2=0}^{p_2}\sum_{j_3=0}^{p_3}\left(\sum_{j_4=0}^{\min\{p_1,p_4\}} C_{j_4 j_3 j_2 j_4}\right)
J'[\phi_{j_2}\phi_{j_3}]_{T,t}^{(i_2i_3)}+
$$

\vspace{2mm}
$$
+{\bf 1}_{\{i_2=i_3\ne 0\}}
\sum_{j_1=0}^{p_1}\sum_{j_4=0}^{p_4}\left(\sum_{j_3=0}^{\min\{p_2,p_3\}} C_{j_4 j_3 j_3 j_1}\right)
J'[\phi_{j_1}\phi_{j_4}]_{T,t}^{(i_1i_4)}+
$$

\vspace{2mm}
$$
+{\bf 1}_{\{i_2=i_4\ne 0\}}
\sum_{j_1=0}^{p_1}\sum_{j_3=0}^{p_3}\left(\sum_{j_4=0}^{\min\{p_2,p_4\}} C_{j_4 j_3 j_4 j_1}\right)
J'[\phi_{j_1}\phi_{j_3}]_{T,t}^{(i_1i_3)}+
$$

\vspace{2mm}
$$
+{\bf 1}_{\{i_3=i_4\ne 0\}}
\sum_{j_1=0}^{p_1}\sum_{j_2=0}^{p_2}\left(\sum_{j_4=0}^{\min\{p_3,p_4\}} C_{j_4 j_4 j_2 j_1}\right)
J'[\phi_{j_1}\phi_{j_2}]_{T,t}^{(i_1i_2)}+
$$

\vspace{2mm}
$$
+{\bf 1}_{\{i_2=i_3\ne 0\}}{\bf 1}_{\{i_1=i_4\ne 0\}}
\sum_{j_2=0}^{\min\{p_2,p_3\}}\sum_{j_4=0}^{\min\{p_1,p_4\}}C_{j_4 j_2 j_2 j_4}+
$$

\vspace{2mm}
$$
+{\bf 1}_{\{i_2=i_4\ne 0\}}{\bf 1}_{\{i_1=i_3\ne 0\}}
\sum_{j_3=0}^{\min\{p_1,p_3\}}\sum_{j_4=0}^{\min\{p_2,p_4\}}C_{j_4 j_3 j_4 j_3}+
$$

\vspace{2mm}
$$
+{\bf 1}_{\{i_3=i_4\ne 0\}}{\bf 1}_{\{i_1=i_2\ne 0\}}
\sum_{j_2=0}^{\min\{p_1,p_2\}}\sum_{j_4=0}^{\min\{p_3,p_4\}}C_{j_4 j_4 j_2 j_2}\ \ \ \hbox{w.~p.~1}.
$$

\vspace{6mm}

Further, we will use the representation (\ref{after8}) for $p_1=\ldots=p_k=p,$ i.e.

\vspace{1mm}
$$
\sum_{j_1,\ldots,j_k=0}^{p}
C_{j_k\ldots j_1}
\prod_{l=1}^k \zeta_{j_l}^{(i_l)}=
\sum_{j_1,\ldots,j_k=0}^{p}
C_{j_k\ldots j_1}
J'[\phi_{j_1}\ldots \phi_{j_k}]_{T,t}^{(i_1\ldots i_k)}+
$$

\vspace{2mm}
$$
+\sum_{j_1,\ldots,j_k=0}^{p}
C_{j_k\ldots j_1}
\sum\limits_{r=1}^{[k/2]}
\sum_{\stackrel{(\{\{g_1, g_2\}, \ldots, 
\{g_{2r-1}, g_{2r}\}\}, \{q_1, \ldots, q_{k-2r}\})}
{{}_{\{g_1, g_2, \ldots, 
g_{2r-1}, g_{2r}, q_1, \ldots, q_{k-2r}\}=\{1, 2, \ldots, k\}}}}
\prod\limits_{s=1}^r
{\bf 1}_{\{i_{g_{{}_{2s-1}}}=~i_{g_{{}_{2s}}}\ne 0\}}\times
$$

\vspace{3mm}
\begin{equation}
\label{after8xx}
\times{\bf 1}_{\{j_{g_{{}_{2s-1}}}=~j_{g_{{}_{2s}}}\}}
J'[\phi_{j_{q_1}}\ldots \phi_{j_{q_{k-2r}}}]_{T,t}^{(i_{q_1}\ldots i_{q_{k-2r}})}\ \ \ \hbox{w.~p.~1.}
\end{equation}

\vspace{7mm}

{\bf Step~2.}\ Let us prove that 

\begin{equation}
\label{after80}
\sum_{j_l=0}^{\infty} C_{j_k \ldots j_{l+1} j_l j_{l-1} \ldots j_{s+1} j_l j_{s-1} \ldots j_1}=0
\end{equation}
or

\begin{equation}
\label{after80xx}
\sum_{j_l=0}^{p} C_{j_k \ldots j_{l+1} j_l j_{l-1} \ldots j_{s+1} j_l j_{s-1} \ldots j_1}=
-\sum_{j_l=p+1}^{\infty} C_{j_k \ldots j_{l+1} j_l j_{l-1} \ldots j_{s+1} j_l j_{s-1} \ldots j_1},
\end{equation}

\vspace{4mm}
\noindent
where
$l-1\ge s+1.$

Our further proof will not fundamentally depend on the weight
functions $\psi_1(\tau),\ldots,\psi_k(\tau).$
There\-fo\-re, sometimes in subsequent consideration we assume
that $\psi_1(\tau),\ldots,\psi_k(\tau)\equiv 1.$

We have 

\vspace{-3mm}
$$
C_{j_k \ldots j_{l+1} j_l j_{l-1} \ldots j_{s+1} j_l j_{s-1} \ldots j_1}=
$$

\vspace{2mm}
$$
=\int\limits_t^T \phi_{j_k}(t_k)\ldots \int\limits_t^{t_{l+2}} \phi_{j_{l+1}}(t_{l+1})
\int\limits_t^{t_{l+1}} \phi_{j_{l}}(t_{l})
\int\limits_t^{t_{l}} \phi_{j_{l-1}}(t_{l-1})\ldots
$$

$$
\ldots
\int\limits_t^{t_{s+2}} \phi_{j_{s+1}}(t_{s+1})
\int\limits_t^{t_{s+1}} \phi_{j_{l}}(t_{s})
\int\limits_t^{t_{s}} \phi_{j_{s-1}}(t_{s-1})\ldots
$$

$$
\ldots \int\limits_t^{t_{2}} \phi_{j_{1}}(t_{1})dt_1\ldots dt_{s-1}dt_{s}dt_{s+1}\ldots
dt_{l-1}dt_{l}dt_{l+1}\ldots dt_k=
$$

$$
=\int\limits_t^{T} \phi_{j_{s+1}}(t_{s+1})
\int\limits_t^{t_{s+1}} \phi_{j_{l}}(t_{s})
\int\limits_t^{t_{s}} \phi_{j_{s-1}}(t_{s-1})\ldots
\int\limits_t^{t_{2}} \phi_{j_{1}}(t_{1})dt_1\ldots dt_{s-1}dt_{s}\times
$$

$$
\times \left(~
\int\limits_{t_{s+1}}^T \phi_{j_{s+2}}(t_{s+2})
\ldots \int\limits_{t_{l-2}}^T \phi_{j_{l-1}}(t_{l-1})
\int\limits_{t_{l-1}}^T \phi_{j_{l}}(t_{l})
\int\limits_{t_{l}}^T \phi_{j_{l+1}}(t_{l+1})\ldots \right.
$$

$$
\left.\ldots
\int\limits_{t_{k-1}}^T \phi_{j_k}(t_k)dt_k\ldots 
dt_{l+1}dt_{l}dt_{l-1}\ldots dt_{s+2}\right)dt_{s+1}=
$$

$$
=\int\limits_t^{T} \phi_{j_{s+1}}(t_{s+1})
\int\limits_t^{t_{s+1}} \phi_{j_{l}}(t_{s})
\underbrace{
\int\limits_t^{t_{s}} \phi_{j_{s-1}}(t_{s-1})\ldots
\int\limits_t^{t_{2}} \phi_{j_{1}}(t_{1})dt_1\ldots dt_{s-1}}
_{G_{j_{s-1}\ldots j_1}(t_s)}
dt_{s}\times
$$
$$
\times
\int\limits_{t_{s+1}}^T \phi_{j_{l}}(t_{l})
\underbrace{\int\limits_{t_{l}}^T \phi_{j_{l+1}}(t_{l+1})
\ldots
\int\limits_{t_{k-1}}^T \phi_{j_k}(t_k)dt_k\ldots 
dt_{l+1}}_{H_{j_{k}\ldots j_{l+1}}(t_l)}\times
$$
$$
\times \left(~
\underbrace{\int\limits_{t_{s+1}}^{t_l} \phi_{j_{l-1}}(t_{l-1})
\ldots \int\limits_{t_{s+1}}^{t_{s+3}} \phi_{j_{s+2}}(t_{s+2})
dt_{s+2}\ldots dt_{l-1}}_{Q_{j_{l-1}\ldots j_{s+2}}(t_l,t_{s+1})}
dt_{l}\right)dt_{s+1}=
$$

\vspace{4mm}
$$
=\int\limits_t^{T} \phi_{j_{s+1}}(t_{s+1})
\int\limits_t^{t_{s+1}} \phi_{j_{l}}(t_{s})
G_{j_{s-1}\ldots j_1}(t_s)
dt_{s}\times
$$
\begin{equation}
\label{after9}
\times
\int\limits_{t_{s+1}}^T \phi_{j_{l}}(t_{l})
H_{j_{k}\ldots j_{l+1}}(t_l)
Q_{j_{l-1}\ldots j_{s+2}}(t_l,t_{s+1})
dt_{l}dt_{s+1}.
\end{equation}

\vspace{5mm}

Using the additive property of the integral, we obtain

\vspace{2mm}
$$
Q_{j_{l-1}\ldots j_{s+2}}(t_l,t_{s+1})=
$$

\vspace{1mm}
$$
=
\int\limits_{t_{s+1}}^{t_l} \phi_{j_{l-1}}(t_{l-1})
\ldots \int\limits_{t_{s+1}}^{t_{s+3}} \phi_{j_{s+2}}(t_{s+2})
dt_{s+2}\ldots dt_{l-1}=
$$

\vspace{1mm}
$$
=
\int\limits_{t_{s+1}}^{t_l} \phi_{j_{l-1}}(t_{l-1})
\ldots \int\limits_{t_{s+1}}^{t_{s+4}}\phi_{j_{s+3}}(t_{s+3})
\int\limits_{t}^{t_{s+3}} \phi_{j_{s+2}}(t_{s+2})
dt_{s+2}dt_{s+3}\ldots dt_{l-1}-
$$

\vspace{1mm}
$$
-\int\limits_{t_{s+1}}^{t_l} \phi_{j_{l-1}}(t_{l-1})
\ldots \int\limits_{t_{s+1}}^{t_{s+4}}\phi_{j_{s+3}}(t_{s+3})
dt_{s+3}\ldots dt_{l-1}\int\limits_{t}^{t_{s+1}} \phi_{j_{s+2}}(t_{s+2})dt_{s+2}=
$$

\vspace{-2mm}
$$
\ldots 
$$

\vspace{-4mm}
\begin{equation}
\label{after10}
=\sum_{m=1}^d h^{(m)}_{j_{l-1}\ldots j_{s+2}}(t_l)
q^{(m)}_{j_{l-1}\ldots j_{s+2}}(t_{s+1}),\ \ \ d<\infty.
\end{equation}

\vspace{4mm}

Combining (\ref{after9}) and (\ref{after10}), we have

$$
\sum_{j_l=0}^p C_{j_k \ldots j_{l+1} j_l j_{l-1} \ldots j_{s+1} j_l j_{s-1} \ldots j_1}=
$$

$$
=\sum_{m=1}^d \left(
\int\limits_t^{T} \phi_{j_{s+1}}(t_{s+1})q^{(m)}_{j_{l-1}\ldots j_{s+2}}(t_{s+1})
\sum_{j_l=0}^p \int\limits_t^{t_{s+1}} \phi_{j_{l}}(t_{s})
G_{j_{s-1}\ldots j_1}(t_s)
dt_{s}\times\right.
$$

\begin{equation}
\label{after11}
\left.\times
\int\limits_{t_{s+1}}^T \phi_{j_{l}}(t_{l})
H_{j_{k}\ldots j_{l+1}}(t_l)
h^{(m)}_{j_{l-1}\ldots j_{s+2}}(t_l)
dt_{l}dt_{s+1}\right).
\end{equation}

\vspace{4mm}

Using the generalized Parseval equality, we obtain

$$
\sum_{j_l=0}^{\infty} \int\limits_t^{t_{s+1}} \phi_{j_{l}}(t_{s})
G_{j_{s-1}\ldots j_1}(t_s)
dt_{s}
\int\limits_{t_{s+1}}^T \phi_{j_{l}}(t_{l})
H_{j_{k}\ldots j_{l+1}}(t_l)
h^{(m)}_{j_{l-1}\ldots j_{s+2}}(t_l)
dt_{l}=
$$

\begin{equation}
\label{after400}
=\int\limits_t^T {\bf 1}_{\{\tau<t_{s+1}\}}
G_{j_{s-1}\ldots j_1}(\tau) \cdot
{\bf 1}_{\{\tau>t_{s+1}\}}
H_{j_{k}\ldots j_{l+1}}(\tau)
h^{(m)}_{j_{l-1}\ldots j_{s+2}}(\tau)
d\tau=0.
\end{equation}

\vspace{4mm}

From (\ref{after11}) and (\ref{after400}) we get

$$
\sum_{j_l=0}^p C_{j_k \ldots j_{l+1} j_l j_{l-1} \ldots j_{s+1} j_l j_{s-1} \ldots j_1}=
$$

$$
=-\sum_{m=1}^d \left(
\int\limits_t^{T} \phi_{j_{s+1}}(t_{s+1})q^{(m)}_{j_{l-1}\ldots j_{s+2}}(t_{s+1})
\sum_{j_l=p+1}^{\infty} \int\limits_t^{t_{s+1}} \phi_{j_{l}}(t_{s})
G_{j_{s-1}\ldots j_1}(t_s)
dt_{s}\times\right.
$$

\begin{equation}
\label{after450}
\left.\times
\int\limits_{t_{s+1}}^T \phi_{j_{l}}(t_{l})
H_{j_{k}\ldots j_{l+1}}(t_l)
h^{(m)}_{j_{l-1}\ldots j_{s+2}}(t_l)
dt_{l}dt_{s+1}\right).
\end{equation}

\vspace{4mm}

Combining Condition~2 of Theorem~13 and (\ref{after9})--(\ref{after11}), (\ref{after450}), we have

$$
\sum_{j_l=0}^p C_{j_k \ldots j_{l+1} j_l j_{l-1} \ldots j_{s+1} j_l j_{s-1} \ldots j_1}=
$$

\vspace{2mm}
$$
=-\sum_{j_l=p+1}^{\infty}\sum_{m=1}^d \left(
\int\limits_t^{T} \phi_{j_{s+1}}(t_{s+1})q^{(m)}_{j_{l-1}\ldots j_{s+2}}(t_{s+1})
\int\limits_t^{t_{s+1}} \phi_{j_{l}}(t_{s})
G_{j_{s-1}\ldots j_1}(t_s)
dt_{s}\times\right.
$$

\vspace{2mm}
$$
\left.\times
\int\limits_{t_{s+1}}^T \phi_{j_{l}}(t_{l})
H_{j_{k}\ldots j_{l+1}}(t_l)
h^{(m)}_{j_{l-1}\ldots j_{s+2}}(t_l)
dt_{l}dt_{s+1}\right)=
$$

\vspace{2mm}
$$
=-\sum_{j_l=p+1}^{\infty}
\int\limits_t^T \phi_{j_k}(t_k)\ldots \int\limits_t^{t_{l+2}} \phi_{j_{l+1}}(t_{l+1})
\int\limits_t^{t_{l+1}} \phi_{j_{l}}(t_{l})
\int\limits_t^{t_{l}} \phi_{j_{l-1}}(t_{l-1})\ldots
$$

\vspace{2mm}
$$
\ldots
\int\limits_t^{t_{s+2}} \phi_{j_{s+1}}(t_{s+1})
\int\limits_t^{t_{s+1}} \phi_{j_{l}}(t_{s})
\int\limits_t^{t_{s}} \phi_{j_{s-1}}(t_{s-1})\ldots
$$

\vspace{2mm}
$$
\ldots \int\limits_t^{t_{2}} \phi_{j_{1}}(t_{1})dt_1\ldots dt_{s-1}dt_{s}dt_{s+1}\ldots
dt_{l-1}dt_{l}dt_{l+1}\ldots dt_k=
$$

\vspace{1mm}
\begin{equation}
\label{after79}
=-\sum_{j_l=p+1}^{\infty} C_{j_k \ldots j_{l+1} j_l j_{l-1} \ldots j_{s+1} j_l j_{s-1} \ldots j_1}.
\end{equation}

\vspace{6mm}

The equality (\ref{after79}) implies (\ref{after80}), (\ref{after80xx}).

\vspace{2mm}

{\bf Step 3.}\ Under the conditions of Theorem~13 we prove that

\vspace{-1mm}
\begin{equation}
\label{after500}
\sum_{j_l=0}^{p} C_{j_k \ldots j_{l+1} j_l j_l j_{l-2} \ldots j_1}
=
\frac{1}{2}
C_{j_k \ldots j_1}\biggl|_{(j_l j_l)\curvearrowright (\cdot)}\biggr.
-\sum_{j_l=p+1}^{\infty} C_{j_k \ldots j_{l+1} j_l j_l j_{l-2} \ldots j_1}.
\end{equation}

\vspace{4mm}

Denote
$$
C_{j_{l-2}\ldots j_1}(t_{l-1})=
\int\limits_t^{t_{l-1}} \psi_{l-2}(t_{l-2})\phi_{j_{l-2}}(t_{l-2})\ldots
\int\limits_t^{t_{2}} \psi_1(t_1)\phi_{j_{1}}(t_{1})dt_1\ldots dt_{l-2}.
$$

\vspace{2mm}

Using the integration order replacement and 
Condition~1 of Theorem~13, we obtain

$$
\sum_{j_l=0}^{\infty} C_{j_k \ldots j_{l+1} j_l j_l j_{l-2} \ldots j_1}
=
\sum_{j_l=0}^{\infty} \int\limits_t^T \psi_k(t_k)\phi_{j_k}(t_k)\ldots 
\int\limits_t^{t_{l+2}} 
\psi_{l+1}(t_{l+1})\phi_{j_{l+1}}(t_{l+1})\times
$$

\vspace{1mm}
$$
\times
\int\limits_t^{t_{l+1}} 
\psi_l(t_l)\phi_{j_{l}}(t_{l})
\int\limits_t^{t_{l}} \psi_{l-1}(t_{l-1}) \phi_{j_{l}}(t_{l-1})
C_{j_{l-2}\ldots j_1}(t_{l-1})
dt_{l-1}
dt_{l}dt_{l+1}\ldots dt_k=
$$

\vspace{1mm}
$$
=
\sum_{j_l=0}^{\infty} \int\limits_t^T 
\psi_l(t_l)\phi_{j_{l}}(t_{l})
\int\limits_t^{t_{l}} \psi_{l-1}(t_{l-1}) \phi_{j_{l}}(t_{l-1})
C_{j_{l-2}\ldots j_1}(t_{l-1})
dt_{l-1}\times
$$

\vspace{1mm}
$$
\times
\int\limits_{t_l}^T
\psi_{l+1}(t_{l+1})\phi_{j_{l+1}}(t_{l+1})\ldots
\int\limits_{t_{k-1}}^T
\psi_{k}(t_{k})\phi_{j_{k}}(t_{k})dt_k\ldots dt_{l+1}dt_l=
$$

\vspace{1mm}

$$
=\frac{1}{2}
\sum_{j_l=0}^{\infty} \int\limits_t^T 
\psi_l(t_l)\psi_{l-1}(t_{l})
C_{j_{l-2}\ldots j_1}(t_{l})
\int\limits_{t_l}^T
\psi_{l+1}(t_{l+1})\phi_{j_{l+1}}(t_{l+1})\ldots
\int\limits_{t_{k-1}}^T
\psi_{k}(t_{k})\phi_{j_{k}}(t_{k})dt_k\ldots dt_{l+1}dt_l=
$$

\vspace{1mm}
$$
=
\frac{1}{2}\sum_{j_l=0}^{\infty} \int\limits_t^T \psi_k(t_k)\phi_{j_k}(t_k)\ldots 
\int\limits_t^{t_{l+2}} 
\psi_{l+1}(t_{l+1})\phi_{j_{l+1}}(t_{l+1})
\int\limits_t^{t_{l+1}} 
\psi_l(t_l)\psi_{l-1}(t_{l})
C_{j_{l-2}\ldots j_1}(t_{l})
dt_{l}dt_{l+1}\ldots dt_k=
$$

\vspace{1mm}
\begin{equation}
\label{r12345x}
=
\frac{1}{2}
C_{j_k \ldots j_1}\biggl|_{(j_l j_l)\curvearrowright (\cdot)}\biggr..
\end{equation}

\vspace{5mm}

The equality (\ref{after500}) is proved.

\vspace{4mm}

{\bf Step~4.}\ Passing to the limit 
$\hbox{\vtop{\offinterlineskip\halign{
\hfil#\hfil\cr
{\rm l.i.m.}\cr
$\stackrel{}{{}_{p\to \infty}}$\cr
}} }$ 
in (\ref{after8xx}), we have (see (\ref{tyyyarg}))

\vspace{2mm}
$$
\hbox{\vtop{\offinterlineskip\halign{
\hfil#\hfil\cr
{\rm l.i.m.}\cr
$\stackrel{}{{}_{p\to \infty}}$\cr
}} }\sum_{j_1,\ldots,j_k=0}^{p}
C_{j_k\ldots j_1}
\zeta_{j_1}^{(i_1)}\ldots \zeta_{j_k}^{(i_k)}
=J[\psi^{(k)}]_{T,t}^{(i_1\ldots i_k)}
+
$$

\vspace{2mm}
$$
+
\sum\limits_{r=1}^{[k/2]}
\sum_{\stackrel{(\{\{g_1, g_2\}, \ldots, 
\{g_{2r-1}, g_{2r}\}\}, \{q_1, \ldots, q_{k-2r}\})}
{{}_{\{g_1, g_2, \ldots, 
g_{2r-1}, g_{2r}, q_1, \ldots, q_{k-2r}\}=\{1, 2, \ldots, k\}}}}
\prod\limits_{s=1}^r
{\bf 1}_{\{i_{g_{{}_{2s-1}}}=~i_{g_{{}_{2s}}}\ne 0\}}\times
$$

\vspace{3mm}
\begin{equation}
\label{after501}
\times \hbox{\vtop{\offinterlineskip\halign{
\hfil#\hfil\cr
{\rm l.i.m.}\cr
$\stackrel{}{{}_{p\to \infty}}$\cr
}} }\sum_{j_1,\ldots,j_k=0}^{p}
C_{j_k\ldots j_1}
\prod\limits_{s=1}^r{\bf 1}_{\{j_{g_{{}_{2s-1}}}=~j_{g_{{}_{2s}}}\}}
J'[\phi_{j_{q_1}}\ldots \phi_{j_{q_{k-2r}}}]_{T,t}^{(i_{q_1}\ldots i_{q_{k-2r}})}\ \ \ 
\hbox{w.~p.~1.}
\end{equation}

\vspace{6mm}

Taking into account (\ref{after80xx}) and (\ref{after500}), we obtain for $r=1$

\vspace{1mm}
$$
{\bf 1}_{\{i_{g_{{}_{1}}}=~i_{g_{{}_{2}}}\ne 0\}}\hbox{\vtop{\offinterlineskip\halign{
\hfil#\hfil\cr
{\rm l.i.m.}\cr
$\stackrel{}{{}_{p\to \infty}}$\cr
}} }\sum_{j_1,\ldots,j_k=0}^{p}
C_{j_k\ldots j_1}
{\bf 1}_{\{j_{g_{{}_{1}}}=~j_{g_{{}_{2}}}\}}
J'[\phi_{j_{q_1}}\ldots \phi_{j_{q_{k-2}}}]_{T,t}^{(i_{q_1}\ldots i_{q_{k-2}})}=
$$

\vspace{3mm}
$$
=-{\bf 1}_{\{i_{g_{{}_{1}}}=~i_{g_{{}_{2}}}\ne 0\}}\hbox{\vtop{\offinterlineskip\halign{
\hfil#\hfil\cr
{\rm l.i.m.}\cr
$\stackrel{}{{}_{p\to \infty}}$\cr
}} }\sum\limits_{j_{g_1}=p+1}^{\infty}
\sum\limits_{\stackrel{j_1,\ldots,j_q,\ldots,j_k=0}{{}_{q\ne g_1, g_2}}}^p
C_{j_k\ldots j_1}\biggl|_{j_{g_{{}_{1}}}=~j_{g_{{}_{2}}}}\biggr. 
\hspace{-1mm}{\bf 1}_{\{g_2>g_1+1\}}
J'[\phi_{j_{q_1}}\ldots \phi_{j_{q_{k-2}}}]_{T,t}^{(i_{q_1}\ldots i_{q_{k-2}})}+
$$

\vspace{2mm}
$$
+
{\bf 1}_{\{i_{g_{{}_{1}}}=~i_{g_{{}_{2}}}\ne 0\}}\hbox{\vtop{\offinterlineskip\halign{
\hfil#\hfil\cr
{\rm l.i.m.}\cr
$\stackrel{}{{}_{p\to \infty}}$\cr
}} }\sum\limits_{\stackrel{j_1,\ldots,j_q,\ldots,j_k=0}{{}_{q\ne g_1, g_2}}}^p
\frac{1}{2}
C_{j_k \ldots j_1}\biggl|_{(j_{g_2} j_{g_1})\curvearrowright (\cdot),
j_{g_{{}_{1}}}=~j_{g_{{}_{2}}}}\biggr.
\hspace{-1mm}{\bf 1}_{\{g_2=g_1+1\}}
J'[\phi_{j_{q_1}}\ldots \phi_{j_{q_{k-2}}}]_{T,t}^{(i_{q_1}\ldots i_{q_{k-2}})}-
$$

\vspace{2mm}
$$
-{\bf 1}_{\{i_{g_{{}_{1}}}=~i_{g_{{}_{2}}}\ne 0\}}\hbox{\vtop{\offinterlineskip\halign{
\hfil#\hfil\cr
{\rm l.i.m.}\cr
$\stackrel{}{{}_{p\to \infty}}$\cr
}} }\sum\limits_{j_{g_1}=p+1}^{\infty}
\sum\limits_{\stackrel{j_1,\ldots,j_q,\ldots,j_k=0}{{}_{q\ne g_1, g_2}}}^p
C_{j_k \ldots j_1}\biggl|_{j_{g_{{}_{1}}}=~j_{g_{{}_{2}}}}\biggr. 
\hspace{-1mm}{\bf 1}_{\{g_2=g_1+1\}}
J'[\phi_{j_{q_1}}\ldots \phi_{j_{q_{k-2}}}]_{T,t}^{(i_{q_1}\ldots i_{q_{k-2}})}=
$$

\vspace{4mm}
$$
=-{\bf 1}_{\{i_{g_{{}_{1}}}=~i_{g_{{}_{2}}}\ne 0\}}\hbox{\vtop{\offinterlineskip\halign{
\hfil#\hfil\cr
{\rm l.i.m.}\cr
$\stackrel{}{{}_{p\to \infty}}$\cr
}} }\sum\limits_{j_{g_1}=p+1}^{\infty}
\sum\limits_{\stackrel{j_1,\ldots,j_q,\ldots,j_k=0}{{}_{q\ne g_1, g_2}}}^p
C_{j_k \ldots j_1}\biggl|_{j_{g_{{}_{1}}}=~j_{g_{{}_{2}}}}\biggr. 
J'[\phi_{j_{q_1}}\ldots \phi_{j_{q_{k-2}}}]_{T,t}^{(i_{q_1}\ldots i_{q_{k-2}})}+
$$

\begin{equation}
\label{after600}
+
{\bf 1}_{\{i_{g_{{}_{1}}}=~i_{g_{{}_{2}}}\ne 0\}}
\hbox{\vtop{\offinterlineskip\halign{
\hfil#\hfil\cr
{\rm l.i.m.}\cr
$\stackrel{}{{}_{p\to \infty}}$\cr
}} }\hspace{-1mm}\sum\limits_{\stackrel{j_1,\ldots,j_q,\ldots,j_k=0}{{}_{q\ne g_1, g_2}}}^p
\frac{1}{2}
C_{j_k \ldots j_1}\biggl|_{(j_{g_2} j_{g_1})\curvearrowright (\cdot),
j_{g_{{}_{1}}}=~j_{g_{{}_{2}}}}\biggr.
\hspace{-1mm}{\bf 1}_{\{g_2=g_1+1\}}
J'[\phi_{j_{q_1}}\ldots \phi_{j_{q_{k-2}}}]_{T,t}^{(i_{q_1}\ldots i_{q_{k-2}})}=
\end{equation}

\vspace{4mm}

\begin{equation}
\label{after607}
=\frac{1}{2}{\bf 1}_{\{g_2=g_1+1\}}
J[\psi^{(k)}]_{T,t}^{g_1}+{\bf 1}_{\{i_{g_{{}_{1}}}=~i_{g_{{}_{2}}}\ne 0\}}
\hbox{\vtop{\offinterlineskip\halign{
\hfil#\hfil\cr
{\rm l.i.m.}\cr
$\stackrel{}{{}_{p\to \infty}}$\cr
}} } R_{T,t}^{(p)1,g_1,g_2}\ \ \ \hbox{w.~p.~1},
\end{equation}

\vspace{6.5mm}
\noindent 
where $J[\psi^{(k)}]_{T,t}^{g_1}$ $(g_1=1,2,\ldots,k-1)$ is
defined by (\ref{30.1}), 

\vspace{1mm}

$$
R_{T,t}^{(p)1,g_1,g_2}=
-\sum\limits_{\stackrel{j_1,\ldots,j_q,\ldots,j_k=0}{{}_{q\ne g_1, g_2}}}^p
\bar C^{(p)}_{j_k\ldots j_q \ldots j_1}\biggl|_{q\ne g_1,g_2}
J'[\phi_{j_{q_1}}\ldots \phi_{j_{q_{k-2}}}]_{T,t}^{(i_{q_1}\ldots i_{q_{k-2}})}.
$$

\vspace{4mm}

Let us explain the transition from 
(\ref{after600}) to (\ref{after607}). We have for $g_2=g_1+1$

$$
{\bf 1}_{\{i_{g_{{}_{1}}}=~i_{g_{{}_{2}}}\ne 0\}}
\hbox{\vtop{\offinterlineskip\halign{
\hfil#\hfil\cr
{\rm l.i.m.}\cr
$\stackrel{}{{}_{p\to \infty}}$\cr
}} }\sum\limits_{\stackrel{j_1,\ldots,j_q,\ldots,j_k=0}{{}_{q\ne g_1, g_2}}}^p
\frac{1}{2}
C_{j_k \ldots j_1}\biggl|_{(j_{g_2} j_{g_1})\curvearrowright (\cdot),
j_{g_{{}_{1}}}=~j_{g_{{}_{2}}}}\biggr.
J'[\phi_{j_{q_1}}\ldots \phi_{j_{q_{k-2}}}]_{T,t}^{(i_{q_1}\ldots i_{q_{k-2}})}=
$$

\vspace{3mm}
$$
=\frac{1}{2}{\bf 1}_{\{i_{g_{{}_{1}}}=~i_{g_{{}_{2}}}\ne 0\}}
\hbox{\vtop{\offinterlineskip\halign{
\hfil#\hfil\cr
{\rm l.i.m.}\cr
$\stackrel{}{{}_{p\to \infty}}$\cr
}} }\sum\limits_{\stackrel{j_1,\ldots,j_q,\ldots,j_k=0}{{}_{q\ne g_1, g_2}}}^p
C_{j_k \ldots j_1}\biggl|_{(j_{g_2} j_{g_1})\curvearrowright 0,
j_{g_{{}_{1}}}=~j_{g_{{}_{2}}}}\biggr.
\zeta_{0}^{(0)} J'[\phi_{j_{q_1}}\ldots \phi_{j_{q_{k-2}}}]_{T,t}^{(i_{q_1}\ldots i_{q_{k-2}})}=
$$

\vspace{4.5mm}
$$
=\frac{1}{2}{\bf 1}_{\{i_{g_{{}_{1}}}=~i_{g_{{}_{2}}}\ne 0\}}
\hbox{\vtop{\offinterlineskip\halign{
\hfil#\hfil\cr
{\rm l.i.m.}\cr
$\stackrel{}{{}_{p\to \infty}}$\cr
}} }\sum\limits_{\stackrel{j_1,\ldots,j_q,\ldots,j_k=0}{{}_{q\ne g_1, g_2}}}^p
\sum_{j_{m_1}=0}^p
C_{j_k \ldots j_1}\biggl|_{(j_{g_2} j_{g_1})\curvearrowright  j_{m_1},
j_{g_{{}_{1}}}=~j_{g_{{}_{2}}}}\biggr.
\times
$$

\vspace{3.5mm}
$$
\times\zeta_{j_{m_1}}^{(0)} 
J'[\phi_{j_{q_1}}\ldots \phi_{j_{q_{k-2}}}]_{T,t}^{(i_{q_1}\ldots i_{q_{k-2}})}=
$$

\vspace{4mm}
$$
=\frac{1}{2}{\bf 1}_{\{i_{g_{{}_{1}}}=~i_{g_{{}_{2}}}\ne 0\}}
\hbox{\vtop{\offinterlineskip\halign{
\hfil#\hfil\cr
{\rm l.i.m.}\cr
$\stackrel{}{{}_{p\to \infty}}$\cr
}} }\sum\limits_{\stackrel{j_1,\ldots,j_q,\ldots,j_k=0}{{}_{q\ne g_1, g_2}}}^p
\sum_{j_{m_1}=0}^p
C_{j_k \ldots j_1}\biggl|_{(j_{g_2} j_{g_1})\curvearrowright  j_{m_1},
j_{g_{{}_{1}}}=~j_{g_{{}_{2}}}}\biggr.\times
$$

\vspace{3mm}
\begin{equation}
\label{after608}
\times J'[\phi_{ j_{m_1}} \phi_{j_{q_1}}\ldots 
\phi_{j_{q_{k-2}}}]_{T,t}^{(0 i_{q_1}\ldots i_{q_{k-2}})}=
\end{equation}

\vspace{4mm}
\begin{equation}
\label{after609}
=
\frac{1}{2}J[\psi^{(k)}]_{T,t}^{g_1}\ \ \ \hbox{w.~p.~1},
\end{equation}

\vspace{4mm}
\noindent
where 

\vspace{-3mm}
$$
C_{j_k \ldots j_1}\biggl|_{(j_{g_2} j_{g_1})\curvearrowright j_{m_1},
j_{g_{{}_{1}}}=~j_{g_{{}_{2}}}, g_2=g_1+1}\biggr.
=
$$

\vspace{3mm}
$$
=
\int\limits_t^T \psi_k(t_k)\phi_{j_k}(t_k)\ldots 
\int\limits_t^{t_{g_1+3}} 
\psi_l(t_{g_1+2})\phi_{j_{g_1+2}}(t_{g_1+2})
\int\limits_t^{t_{g_1+2}} 
\psi_{g_1+1}(t_{g_1})\psi_{g_1}(t_{g_1})\phi_{j_{m_1}}(t_{g_1})\times
$$

\vspace{1mm}
$$
\times
\int\limits_t^{t_{g_1}} \psi_l(t_{g_1-1})\phi_{j_{g_1-1}}(t_{g_1-1})\ldots
\int\limits_t^{t_{2}} \psi_1(t_1)\phi_{j_{1}}(t_{1})dt_1\ldots dt_{g_1-1}dt_{g_1}
dt_{g_1+2}\ldots dt_k,
$$

\vspace{3mm}
$$
\zeta_{j_{m_1}}^{(0)}=\int\limits_t^T\phi_{j_{m_1}}(\tau)d{\bf w}_{\tau}^{(0)}
=\int\limits_t^T\phi_{j_{m_1}}(\tau)d\tau=
\left\{
\begin{matrix}
\sqrt{T-t} &\hbox{if}\ j_{m_1}=0\cr\cr
0 &\hbox{if}\ j_{m_1}\ne 0
\end{matrix},\right.
$$ 

\vspace{3mm}
$$
\phi_0(\tau)=\frac{1}{\sqrt{T-t}}.
$$

\vspace{6mm}
\noindent
The transition from (\ref{after608}) to (\ref{after609}) is based
on (\ref{tyyyarg}).

By Condition~3 of Theorem~13 we have (also see the property (\ref{wiener1}) of
multiple Wiener stochastic integral)

\vspace{-2mm}
$$
\lim\limits_{p\to\infty}
{\sf M}\left\{\left(R_{T,t}^{(p)1,g_1,g_2}\right)^2\right\}
\le 
K \lim\limits_{p\to\infty}
\sum\limits_{\stackrel{j_1,\ldots,j_q,\ldots,j_k=0}{{}_{q\ne g_1, g_2}}}^p
\left(\bar C^{(p)}_{j_k\ldots j_q \ldots j_1}\biggl|_{q\ne g_1,g_2}\right)^2 =0,
$$
  
\vspace{3mm}
\noindent
where constant $K$ does not depend on $p$. 

Thus

\vspace{-2mm}
$$
{\bf 1}_{\{i_{g_{{}_{1}}}=~i_{g_{{}_{2}}}\ne 0\}}\hbox{\vtop{\offinterlineskip\halign{
\hfil#\hfil\cr
{\rm l.i.m.}\cr
$\stackrel{}{{}_{p\to \infty}}$\cr
}} }\sum_{j_1,\ldots,j_k=0}^{p}
C_{j_k\ldots j_1}
{\bf 1}_{\{j_{g_{{}_{1}}}=~j_{g_{{}_{2}}}\}}
J'[\phi_{j_{q_1}}\ldots \phi_{j_{q_{k-2}}}]_{T,t}^{(i_{q_1}\ldots i_{q_{k-2}})}=
$$

\vspace{2mm}
$$
=\frac{1}{2}{\bf 1}_{\{g_2=g_1+1\}}
J[\psi^{(k)}]_{T,t}^{g_1}\ \ \ \hbox{w.~p.~1}.
$$

\vspace{5mm}

Involving into consideration the second pair $\{g_3,g_4\}$
(the first pair is $\{g_1,g_2\}$), we obtain
from (\ref{after600}) for $r=2$

$$
\prod\limits_{s=1}^2
{\bf 1}_{\{i_{g_{{}_{2s-1}}}=~i_{g_{{}_{2s}}}\ne 0\}}\hbox{\vtop{\offinterlineskip\halign{
\hfil#\hfil\cr
{\rm l.i.m.}\cr
$\stackrel{}{{}_{p\to \infty}}$\cr
}} }\sum_{j_1,\ldots,j_k=0}^{p}
C_{j_k\ldots j_1}
\prod\limits_{s=1}^2{\bf 1}_{\{j_{g_{{}_{2s-1}}}=~j_{g_{{}_{2s}}}\}}\times
$$

\vspace{4mm}
$$
\times
J'[\phi_{j_{q_1}}\ldots \phi_{j_{q_{k-4}}}]_{T,t}^{(i_{q_1}\ldots i_{q_{k-4}})}=
$$

\vspace{3mm}
$$
=\prod\limits_{s=1}^2
{\bf 1}_{\{i_{g_{{}_{2s-1}}}=~i_{g_{{}_{2s}}}\ne 0\}}
\times
$$

$$
\times\hbox{\vtop{\offinterlineskip\halign{
\hfil#\hfil\cr
{\rm l.i.m.}\cr
$\stackrel{}{{}_{p\to \infty}}$\cr
}} }
\sum\limits_{\stackrel{j_1,\ldots,j_q,\ldots,j_k=0}{{}_{q\ne g_1, g_2, g_3, g_4}}}^p
\Biggl(\frac{1}{4}
C_{j_k \ldots j_1}\biggl|_{(j_{g_2} j_{g_1})\curvearrowright (\cdot)
(j_{g_4} j_{g_3})\curvearrowright (\cdot),
j_{g_{{}_{1}}}=~j_{g_{{}_{2}}}, j_{g_{{}_{3}}}=~j_{g_{{}_{4}}}}\biggr.
\prod\limits_{s=1}^2 {\bf 1}_{\{g_{2s}=g_{2s-1}+1\}}-\Biggr.
$$

\vspace{2mm}
$$
-\frac{1}{2}\sum\limits_{j_{g_1}=p+1}^{\infty}
C_{j_k \ldots j_1}\biggl|_{(j_{g_4} j_{g_3})\curvearrowright (\cdot),
j_{g_{{}_{1}}}=~j_{g_{{}_{2}}}, j_{g_{{}_{3}}}=~j_{g_{{}_{4}}}}
{\bf 1}_{\{g_{4}=g_{3}+1\}}-
$$

\vspace{2mm}
$$
-
\frac{1}{2}\sum\limits_{j_{g_3}=p+1}^{\infty}
C_{j_k \ldots j_1}\biggl|_{(j_{g_2} j_{g_1})\curvearrowright (\cdot),
j_{g_{{}_{1}}}=~j_{g_{{}_{2}}}, j_{g_{{}_{3}}}=~j_{g_{{}_{4}}}}\biggr.
{\bf 1}_{\{g_{2}=g_{1}+1\}}+
$$

\vspace{2mm}
\begin{equation}
\label{after610}
\Biggl.+
\sum\limits_{j_{g_3}=p+1}^{\infty}\sum\limits_{j_{g_1}=p+1}^{\infty}
C_{j_k \ldots j_1}\biggl|_{j_{g_{{}_{1}}}=~j_{g_{{}_{2}}}, j_{g_{{}_{3}}}=~j_{g_{{}_{4}}}} 
\biggr.
\Biggr)J'[\phi_{j_{q_1}}\ldots 
\phi_{j_{q_{k-4}}}]_{T,t}^{(i_{q_1}\ldots i_{q_{k-4}})}=
\end{equation}

\vspace{2mm}
\begin{equation}
\label{after610x}
=
\frac{1}{4}\prod\limits_{s=1}^2 {\bf 1}_{\{g_{2s}=g_{2s-1}+1\}}
J[\psi^{(k)}]_{T,t}^{s_2,s_1}+  
\prod\limits_{s=1}^2
{\bf 1}_{\{i_{g_{{}_{2s-1}}}=~i_{g_{{}_{2s}}}\ne 0\}}
\hbox{\vtop{\offinterlineskip\halign{
\hfil#\hfil\cr
{\rm l.i.m.}\cr
$\stackrel{}{{}_{p\to \infty}}$\cr
}} } R_{T,t}^{(p)2,g_1,g_2,g_3,g_4}
\end{equation}

\vspace{5mm}
\noindent
w.~p.~1, where $g_3\stackrel{\sf def}{=}s_2,$ $g_1\stackrel{\sf def}{=}s_1,$
$(s_2,s_1)\in {\rm A}_{k,2},$ $J[\psi^{(k)}]_{T,t}^{s_2,s_1}$ is
defined by (\ref{30.1}) and ${\rm A}_{k,2}$ is defined by (\ref{30.5550001}),

$$
R_{T,t}^{(p)2,g_1,g_2,g_3,g_4}=
\sum\limits_{\stackrel{j_1,\ldots,j_q,\ldots,j_k=0}{{}_{q\ne g_1, g_2, g_3, g_4}}}^p
\left(
\bar C^{(p)}_{j_k\ldots j_q \ldots j_1}\biggl|_{q\ne g_1,g_2,g_3,g_4}-\right.
$$

\vspace{3mm}
$$
-S_1\left\{
\bar C^{(p)}_{j_k\ldots j_q \ldots j_1}\biggl|_{q\ne g_1,g_2,g_3,g_4}\right\}
\left.-S_2\left\{
\bar C^{(p)}_{j_k\ldots j_q \ldots j_1}\biggl|_{q\ne g_1,g_2,g_3,g_4}\right\}\right)\times
$$

\vspace{5mm}
$$
\times
J'[\phi_{j_{q_1}}\ldots \phi_{j_{q_{k-4}}}]_{T,t}^{(i_{q_1}\ldots i_{q_{k-4}})}.
$$

\vspace{7mm}

Let us explain the transition from (\ref{after610}) to (\ref{after610x}).
We have for $g_2=g_1+1,$ $g_4=g_3+1$

\vspace{1mm}
$$
\hbox{\vtop{\offinterlineskip\halign{
\hfil#\hfil\cr
{\rm l.i.m.}\cr
$\stackrel{}{{}_{p\to \infty}}$\cr
}} }
\sum\limits_{\stackrel{j_1,\ldots,j_q,\ldots,j_k=0}{{}_{q\ne g_1, g_2, g_3, g_4}}}^p
\frac{1}{4}
C_{j_k \ldots j_1}\biggl|_{(j_{g_2} j_{g_1})\curvearrowright (\cdot)
(j_{g_4} j_{g_3})\curvearrowright (\cdot),
j_{g_{{}_{1}}}=~j_{g_{{}_{2}}}, j_{g_{{}_{3}}}=~j_{g_{{}_{4}}}}\biggr.
\times
$$

\vspace{3mm}
$$
\times 
\prod\limits_{s=1}^2
{\bf 1}_{\{i_{g_{{}_{2s-1}}}=~i_{g_{{}_{2s}}}\ne 0\}}
J'[\phi_{j_{q_1}}\ldots \phi_{j_{q_{k-4}}}]_{T,t}^{(i_{q_1}\ldots i_{q_{k-4}})}=
$$

\vspace{6mm}
$$
=\frac{1}{4}
\hbox{\vtop{\offinterlineskip\halign{
\hfil#\hfil\cr
{\rm l.i.m.}\cr
$\stackrel{}{{}_{p\to \infty}}$\cr
}} }
\sum\limits_{\stackrel{j_1,\ldots,j_q,\ldots,j_k=0}{{}_{q\ne g_1, g_2, g_3, g_4}}}^p
C_{j_k \ldots j_1}\biggl|_{(j_{g_2} j_{g_1})\curvearrowright 0
(j_{g_4} j_{g_3})\curvearrowright 0,
j_{g_{{}_{1}}}=~j_{g_{{}_{2}}}, j_{g_{{}_{3}}}=~j_{g_{{}_{4}}}}\biggr.
\times
$$

\vspace{4mm}
$$
\times
\prod\limits_{s=1}^2
{\bf 1}_{\{i_{g_{{}_{2s-1}}}=~i_{g_{{}_{2s}}}\ne 0\}}
\zeta_{0}^{(0)}\zeta_{0}^{(0)}
J'[\phi_{j_{q_1}}\ldots \phi_{j_{q_{k-4}}}]_{T,t}^{(i_{q_1}\ldots i_{q_{k-4}})}=
$$

\vspace{5mm}
$$
=\frac{1}{4}
\hbox{\vtop{\offinterlineskip\halign{
\hfil#\hfil\cr
{\rm l.i.m.}\cr
$\stackrel{}{{}_{p\to \infty}}$\cr
}} }
\sum\limits_{\stackrel{j_1,\ldots,j_q,\ldots,j_k=0}{{}_{q\ne g_1, g_2, g_3, g_4}}}^p
\sum\limits_{ j_{m_1}, j_{m_3}=0}^p
C_{j_k \ldots j_1}\biggl|_{(j_{g_2} j_{g_1})\curvearrowright j_{m_1}
(j_{g_4} j_{g_3})\curvearrowright j_{m_3},
j_{g_{{}_{1}}}=~j_{g_{{}_{2}}}, j_{g_{{}_{3}}}=~j_{g_{{}_{4}}}}\biggr.
\times
$$

\vspace{4mm}
$$
\times
\prod\limits_{s=1}^2
{\bf 1}_{\{i_{g_{{}_{2s-1}}}=~i_{g_{{}_{2s}}}\ne 0\}}
\zeta_{ j_{m_1}}^{(0)}\zeta_{j_{m_3}}^{(0)}
J'[\phi_{j_{q_1}}\ldots \phi_{j_{q_{k-4}}}]_{T,t}^{(i_{q_1}\ldots i_{q_{k-4}})}=
$$

\vspace{6mm}
$$
=\frac{1}{4}
\hbox{\vtop{\offinterlineskip\halign{
\hfil#\hfil\cr
{\rm l.i.m.}\cr
$\stackrel{}{{}_{p\to \infty}}$\cr
}} }
\sum\limits_{\stackrel{j_1,\ldots,j_q,\ldots,j_k=0}{{}_{q\ne g_1, g_2, g_3, g_4}}}^p
\sum\limits_{j_{m_1}, j_{m_3}=0}^p
C_{j_k \ldots j_1}\biggl|_{(j_{g_2} j_{g_1})\curvearrowright j_{m_1}
(j_{g_4} j_{g_3})\curvearrowright j_{m_3},
j_{g_{{}_{1}}}=~j_{g_{{}_{2}}}, j_{g_{{}_{3}}}=~j_{g_{{}_{4}}}}\biggr.
\times
$$

\vspace{4mm}
\begin{equation}
\label{after900x}
\times
\prod\limits_{s=1}^2
{\bf 1}_{\{i_{g_{{}_{2s-1}}}=~i_{g_{{}_{2s}}}\ne 0\}}
J'[\phi_{ j_{m_1}}\phi_{j_{m_3}}
\phi_{j_{q_1}}\ldots \phi_{j_{q_{k-4}}}]_{T,t}^{(0 0 i_{q_1}\ldots i_{q_{k-4}})}=
\end{equation}

\vspace{3mm}
\begin{equation}
\label{after901}
=
\frac{1}{4}
J[\psi^{(k)}]_{T,t}^{s_2,s_1}\ \ \ \hbox{w.~p.~1}.
\end{equation}

\vspace{6mm}
\noindent
The transition from (\ref{after900x}) to (\ref{after901}) is based
on (\ref{tyyyarg}).

Note that

\vspace{-2mm}
$$
C_{j_k \ldots j_1}\biggl|_{(j_{g_2} j_{g_1})\curvearrowright j_{m_1},
j_{g_{{}_{1}}}=~j_{g_{{}_{2}}}}\biggr.
=C_{j_k \ldots j_1}\biggl|_{(j_{g_1} j_{g_1})\curvearrowright j_{m_1},
j_{g_{{}_{1}}}=~j_{g_{{}_{2}}}}\biggr.
$$

\vspace{5mm}
\noindent
is the Fourier coefficient, where $g_{2}=g_{1}+1.$
Therefore, the value

$$
C_{j_k \ldots j_1}\biggl|_{(j_{g_2} j_{g_1})\curvearrowright j_{m_1}
(j_{g_4} j_{g_3})\curvearrowright  j_{m_3},
j_{g_{{}_{1}}}=~j_{g_{{}_{2}}},
j_{g_{{}_{3}}}=~j_{g_{{}_{4}}}
}\biggr.=
$$

\vspace{2mm}
$$
=C_{j_k \ldots j_1}\biggl|_{(j_{g_1} j_{g_1})\curvearrowright j_{m_1}
(j_{g_3} j_{g_3})\curvearrowright  j_{m_3},
j_{g_{{}_{1}}}=~j_{g_{{}_{2}}},
j_{g_{{}_{3}}}=~j_{g_{{}_{4}}}
}\biggr.
$$

\vspace{5mm}
\noindent
is determined recursively using (\ref{after2000}) in an obvious way
for $g_2=g_1+1$ and $g_4=g_3+1$.

By Condition~3 of Theorem~13 we have (also see the property (\ref{wiener1}) of
multiple Wiener stochastic integral)

$$
\lim\limits_{p\to\infty}
{\sf M}\left\{\left(R_{T,t}^{(p)2,g_1,g_2,g_3,g_4}\right)^2\right\}
\le
K \lim\limits_{p\to\infty}
\sum\limits_{\stackrel{j_1,\ldots,j_q,\ldots,j_k=0}{{}_{q\ne g_1, g_2, g_3, g_4}}}^p
\left(
\left(
\bar C^{(p)}_{j_k\ldots j_q \ldots j_1}\biggl|_{q\ne g_1,g_2,g_3,g_4}\biggr.\right)^2+
\right.
$$

\vspace{3mm}
$$
\left.+\left(S_1\left\{
\bar C^{(p)}_{j_k\ldots j_q \ldots j_1}\biggl|_{q\ne g_1,g_2,g_3,g_4}\biggr.
\right\}\right)^2+
\left(S_2\left\{
\bar C^{(p)}_{j_k\ldots j_q \ldots j_1}\biggl|_{q\ne g_1,g_2,g_3,g_4}\biggr.
\right\}\right)^2\right)=0,
$$

\vspace{6mm}
\noindent
where constant $K$ is independent of $p$.

Thus

\vspace{-2mm}
$$
\prod\limits_{s=1}^2
{\bf 1}_{\{i_{g_{{}_{2s-1}}}=~i_{g_{{}_{2s}}}\ne 0\}}\hbox{\vtop{\offinterlineskip\halign{
\hfil#\hfil\cr
{\rm l.i.m.}\cr
$\stackrel{}{{}_{p\to \infty}}$\cr
}} }\sum_{j_1,\ldots,j_k=0}^{p}
C_{j_k\ldots j_1}
\prod\limits_{s=1}^2{\bf 1}_{\{j_{g_{{}_{2s-1}}}=~j_{g_{{}_{2s}}}\}}
\times
$$

\vspace{4mm}
$$
\times J'[\phi_{j_{q_1}}\ldots \phi_{j_{q_{k-4}}}]_{T,t}^{(i_{q_1}\ldots i_{q_{k-4}})}=
\frac{1}{4}\prod\limits_{s=1}^2 {\bf 1}_{\{g_{2s}=g_{2s-1}+1\}}
J[\psi^{(k)}]_{T,t}^{s_2,s_1}\ \ \ \hbox{w.~p.~1},
$$

\vspace{4mm}
\noindent
where $g_3\stackrel{\sf def}{=}s_2,$ $g_1\stackrel{\sf def}{=}s_1,$
$(s_2,s_1)\in {\rm A}_{k,2},$ $J[\psi^{(k)}]_{T,t}^{s_2,s_1}$ is
defined by (\ref{30.1}) and ${\rm A}_{k,2}$ is defined by (\ref{30.5550001}).

Involving into consideration the third pair $\{g_6,g_5\}$
($\{g_1,g_2\}$ is the first pair and $\{g_4,g_3\}$ is the second pair), we obtain
from (\ref{after610}) for $r=3$

\vspace{1mm}
$$
\prod\limits_{s=1}^3
{\bf 1}_{\{i_{g_{{}_{2s-1}}}=~i_{g_{{}_{2s}}}\ne 0\}}\hbox{\vtop{\offinterlineskip\halign{
\hfil#\hfil\cr
{\rm l.i.m.}\cr
$\stackrel{}{{}_{p\to \infty}}$\cr
}} }\sum_{j_1,\ldots,j_k=0}^{p}
C_{j_k\ldots j_1}
\prod\limits_{s=1}^3{\bf 1}_{\{j_{g_{{}_{2s-1}}}=~j_{g_{{}_{2s}}}\}}
\times
$$

\vspace{4mm}
$$
\times
J'[\phi_{j_{q_1}}\ldots \phi_{j_{q_{k-6}}}]_{T,t}^{(i_{q_1}\ldots i_{q_{k-6}})}=
\prod\limits_{s=1}^3
{\bf 1}_{\{i_{g_{{}_{2s-1}}}=~i_{g_{{}_{2s}}}\ne 0\}}\times
$$

\vspace{1mm}
$$
\times\hbox{\vtop{\offinterlineskip\halign{
\hfil#\hfil\cr
{\rm l.i.m.}\cr
$\stackrel{}{{}_{p\to \infty}}$\cr
}} }\hspace{-3mm}
\sum\limits_{\stackrel{j_1,\ldots,j_q,\ldots,j_k=0}{{}_{q\ne g_1, g_2, g_3, g_4, g_5, g_6}}}^p
\hspace{-2mm}
\left(\frac{1}{2^3}
C_{j_k \ldots j_1}\biggl|_{(j_{g_2} j_{g_1})\curvearrowright (\cdot)
(j_{g_4} j_{g_3})\curvearrowright (\cdot)
(j_{g_6} j_{g_5})\curvearrowright (\cdot),
j_{g_{{}_{1}}}=~j_{g_{{}_{2}}}, j_{g_{{}_{3}}}=~j_{g_{{}_{4}}}, j_{g_{{}_{5}}}=~j_{g_{{}_{6}}}
}\biggr.\right.
\times
$$

\vspace{2mm}
$$
\times
\prod\limits_{s=1}^3
{\bf 1}_{\{g_{2s}=g_{2s-1}+1\}}-
$$

$$
-\frac{1}{2^2}\sum\limits_{j_{g_1}=p+1}^{\infty}
C_{j_k \ldots j_1}\biggl|_{(j_{g_4} j_{g_3})\curvearrowright (\cdot)
(j_{g_6} j_{g_5})\curvearrowright (\cdot),
j_{g_{{}_{1}}}=~j_{g_{{}_{2}}}, j_{g_{{}_{3}}}=~j_{g_{{}_{4}}}, j_{g_{{}_{5}}}=~j_{g_{{}_{6}}}}
{\bf 1}_{\{g_{4}=g_{3}+1\}}{\bf 1}_{\{g_{6}=g_{5}+1\}}-
$$

\vspace{2mm}
$$
-\frac{1}{2^2}\sum\limits_{j_{g_3}=p+1}^{\infty}
C_{j_k \ldots j_1}\biggl|_{(j_{g_2} j_{g_1})\curvearrowright (\cdot)
(j_{g_6} j_{g_5})\curvearrowright (\cdot),
j_{g_{{}_{1}}}=~j_{g_{{}_{2}}}, j_{g_{{}_{3}}}=~j_{g_{{}_{4}}}, j_{g_{{}_{5}}}=~j_{g_{{}_{6}}}}
{\bf 1}_{\{g_{2}=g_{1}+1\}}{\bf 1}_{\{g_{6}=g_{5}+1\}}-
$$

\vspace{2mm}
$$
-\frac{1}{2^2}\sum\limits_{j_{g_5}=p+1}^{\infty}
C_{j_k \ldots j_1}\biggl|_{(j_{g_2} j_{g_1})\curvearrowright (\cdot)
(j_{g_4} j_{g_3})\curvearrowright (\cdot),
j_{g_{{}_{1}}}=~j_{g_{{}_{2}}}, j_{g_{{}_{3}}}=~j_{g_{{}_{4}}}, j_{g_{{}_{5}}}=~j_{g_{{}_{6}}}}
{\bf 1}_{\{g_{2}=g_{1}+1\}}{\bf 1}_{\{g_{4}=g_{3}+1\}}+
$$

\vspace{4mm}
$$
+\frac{1}{2}\sum\limits_{j_{g_3}=p+1}^{\infty}\sum\limits_{j_{g_1}=p+1}^{\infty}
C_{j_k \ldots j_1}\biggl|_{(j_{g_6} j_{g_5})\curvearrowright (\cdot),
j_{g_{{}_{1}}}=~j_{g_{{}_{2}}}, j_{g_{{}_{3}}}=~j_{g_{{}_{4}}}, j_{g_{{}_{5}}}=~j_{g_{{}_{6}}}}
{\bf 1}_{\{g_{6}=g_{5}+1\}}+
$$

\vspace{4mm}
$$
+\frac{1}{2}\sum\limits_{j_{g_5}=p+1}^{\infty}\sum\limits_{j_{g_1}=p+1}^{\infty}
C_{j_k \ldots j_1}\biggl|_{(j_{g_4} j_{g_3})\curvearrowright (\cdot),
j_{g_{{}_{1}}}=~j_{g_{{}_{2}}}, j_{g_{{}_{3}}}=~j_{g_{{}_{4}}}, j_{g_{{}_{5}}}=~j_{g_{{}_{6}}}}
{\bf 1}_{\{g_{4}=g_{3}+1\}}+
$$

\vspace{4mm}
$$
+\frac{1}{2}\sum\limits_{j_{g_5}=p+1}^{\infty}\sum\limits_{j_{g_3}=p+1}^{\infty}
C_{j_k \ldots j_1}\biggl|_{(j_{g_2} j_{g_1})\curvearrowright (\cdot),
j_{g_{{}_{1}}}=~j_{g_{{}_{2}}}, j_{g_{{}_{3}}}=~j_{g_{{}_{4}}}, j_{g_{{}_{5}}}=~j_{g_{{}_{6}}}}
{\bf 1}_{\{g_{2}=g_{1}+1\}}-
$$

\vspace{4mm}
$$
\left.-
\sum\limits_{j_{g_5}=p+1}^{\infty}
\sum\limits_{j_{g_3}=p+1}^{\infty}\sum\limits_{j_{g_1}=p+1}^{\infty}
C_{j_k \ldots j_1}\biggl|_{j_{g_{{}_{1}}}=~j_{g_{{}_{2}}}, 
j_{g_{{}_{3}}}=~j_{g_{{}_{4}}}, j_{g_{{}_{5}}}=~j_{g_{{}_{6}}}}
\right)\times
$$

\vspace{6mm}
$$
\times
J'[\phi_{j_{q_1}}\ldots \phi_{j_{q_{k-6}}}]_{T,t}^{(i_{q_1}\ldots i_{q_{k-6}})}=
$$

\vspace{4mm}
$$
=
\frac{1}{2^3}\prod\limits_{s=1}^3
{\bf 1}_{\{g_{2s}=g_{2s-1}+1\}}
J[\psi^{(k)}]_{T,t}^{s_3,s_2,s_1}+
\prod\limits_{s=1}^3
{\bf 1}_{\{i_{g_{{}_{2s-1}}}=~i_{g_{{}_{2s}}}\ne 0\}}
\hbox{\vtop{\offinterlineskip\halign{
\hfil#\hfil\cr
{\rm l.i.m.}\cr
$\stackrel{}{{}_{p\to \infty}}$\cr
}} } R_{T,t}^{(p)3,g_1,g_2,\ldots,g_5,g_6}
$$

\vspace{5mm}
\noindent
w.~p.~1, where $g_{2i-1}\stackrel{\sf def}{=}s_i;$\ $i=1,2,3,$ 
$(s_3,s_2,s_1)\in {\rm A}_{k,3},$ $J[\psi^{(k)}]_{T,t}^{s_3,s_2,s_1}$ is
defined by (\ref{30.1}) and ${\rm A}_{k,3}$ is defined by (\ref{30.5550001}),

\vspace{1mm}

$$
R_{T,t}^{(p)3,g_1,g_2,\ldots,g_5,g_6}=
\sum\limits_{\stackrel{j_1,\ldots,j_q,\ldots,j_k=0}{{}_{q\ne g_1,g_2,\ldots,g_5,g_6}}}^p
\left(
-\bar C^{(p)}_{j_k\ldots j_q \ldots j_1}\biggl|_{q\ne g_1,g_2,\ldots,g_5,g_6}+\right.
$$

\vspace{4mm}
$$
+S_1\left\{
\bar C^{(p)}_{j_k\ldots j_q \ldots j_1}\biggl|_{q\ne g_1,g_2,\ldots,g_5,g_6}\right\}
+S_2\left\{
\bar C^{(p)}_{j_k\ldots j_q \ldots j_1}\biggl|_{q\ne g_1,g_2,\ldots,g_5,g_6}\right\}+
$$

\vspace{4mm}
$$       
+S_3\left\{
\bar C^{(p)}_{j_k\ldots j_q \ldots j_1}\biggl|_{q\ne g_1,g_2,\ldots,g_5,g_6}\right\}-
$$

\vspace{4mm}
$$
-S_3S_1\left\{
\bar C^{(p)}_{j_k\ldots j_q \ldots j_1}\biggl|_{q\ne g_1,g_2,\ldots,g_5,g_6}\right\}
-S_3S_2\left\{
\bar C^{(p)}_{j_k\ldots j_q \ldots j_1}\biggl|_{q\ne g_1,g_2,\ldots,g_5,g_6}\right\}-
$$

\vspace{6mm}
$$
\left.-S_2S_1\left\{
\bar C^{(p)}_{j_k\ldots j_q \ldots j_1}\biggl|_{q\ne g_1,g_2,\ldots,g_5,g_6}\right\}\right)
J'[\phi_{j_{q_1}}\ldots \phi_{j_{q_{k-6}}}]_{T,t}^{(i_{q_1}\ldots i_{q_{k-6}})}.
$$

\vspace{7mm}

By Condition~3 of Theorem~13 we have (also see the property (\ref{wiener1}) of
multiple Wiener stochastic integral)

\vspace{1mm}

$$
\lim\limits_{p\to\infty}
{\sf M}\left\{\left(R_{T,t}^{(p)3,g_1,g_2,\ldots,g_5,g_6}\right)^2\right\}
\le K \lim\limits_{p\to\infty}
\sum\limits_{\stackrel{j_1,\ldots,j_q,\ldots,j_k=0}{{}_{q\ne g_1, g_2,\ldots, g_5, g_6}}}^p
\left(\left(
\bar C^{(p)}_{j_k\ldots j_q \ldots j_1}\biggl|_{q\ne g_1,g_2,\ldots, g_5, g_6}\right)^2+\right.
$$

\vspace{4mm}
$$
+\left(S_1\left\{
\bar C^{(p)}_{j_k\ldots j_q \ldots j_1}\biggl|_{q\ne g_1,g_2,\ldots,g_5,g_6}\right\}\right)^2
+\left(S_2\left\{
\bar C^{(p)}_{j_k\ldots j_q \ldots j_1}\biggl|_{q\ne g_1,g_2,\ldots,g_5,g_6}\right\}\right)^2+
$$

\vspace{4mm}
$$
+\left(S_3\left\{
\bar C^{(p)}_{j_k\ldots j_q \ldots j_1}\biggl|_{q\ne g_1,g_2,\ldots,g_5,g_6}\right\}\right)^2+
$$

\vspace{4mm}
$$
+\left(S_3S_1\left\{
\bar C^{(p)}_{j_k\ldots j_q \ldots j_1}\biggl|_{q\ne g_1,g_2,\ldots,g_5,g_6}\right\}\right)^2
+\left(S_3S_2\left\{
\bar C^{(p)}_{j_k\ldots j_q \ldots j_1}\biggl|_{q\ne g_1,g_2,\ldots,g_5,g_6}\right\}\right)^2+
$$

\vspace{4mm}
$$
\left.+\left(S_2S_1\left\{
\bar C^{(p)}_{j_k\ldots j_q \ldots j_1}\biggl|_{q\ne g_1,g_2,\ldots,g_5,g_6}\right\}\right)^2
\right)=0,
$$

\vspace{5mm}
\noindent
where constant $K$ does not depend on $p$.

Thus

\vspace{-1mm}
$$
\hbox{\vtop{\offinterlineskip\halign{
\hfil#\hfil\cr
{\rm l.i.m.}\cr
$\stackrel{}{{}_{p\to \infty}}$\cr
}} }\prod\limits_{s=1}^3
{\bf 1}_{\{i_{g_{{}_{2s-1}}}=~i_{g_{{}_{2s}}}\ne 0\}}\hbox{\vtop{\offinterlineskip\halign{
\hfil#\hfil\cr
{\rm l.i.m.}\cr
$\stackrel{}{{}_{p\to \infty}}$\cr
}} }\sum_{j_1,\ldots,j_k=0}^{p}
C_{j_k\ldots j_1}
\prod\limits_{s=1}^3{\bf 1}_{\{j_{g_{{}_{2s-1}}}=~j_{g_{{}_{2s}}}\}}\times
$$

\vspace{5mm}
$$
\times
J'[\phi_{j_{q_1}}\ldots \phi_{j_{q_{k-6}}}]_{T,t}^{(i_{q_1}\ldots i_{q_{k-6}})}=
\frac{1}{2^3}\prod\limits_{s=1}^3
{\bf 1}_{\{g_{2s}=g_{2s-1}+1\}}
J[\psi^{(k)}]_{T,t}^{s_3,s_2,s_1}
\ \ \ \hbox{w.~p.~1},
$$

\vspace{6mm}
\noindent
where $g_{2i-1}\stackrel{\sf def}{=}s_i;$\ $i=1,2,3,$ 
$(s_3,s_2,s_1)\in {\rm A}_{k,3},$ $J[\psi^{(k)}]_{T,t}^{s_3,s_2,s_1}$ is
defined by (\ref{30.1}) and ${\rm A}_{k,3}$ is defined by (\ref{30.5550001}).

Repeating the previous steps, we obtain for an arbitrary $r$ 
($r=1,2,$ $\ldots,$ $[k/2]$)

\vspace{1mm}
$$
\prod\limits_{s=1}^r
{\bf 1}_{\{i_{g_{{}_{2s-1}}}=~i_{g_{{}_{2s}}}\ne 0\}}\hbox{\vtop{\offinterlineskip\halign{
\hfil#\hfil\cr
{\rm l.i.m.}\cr
$\stackrel{}{{}_{p\to \infty}}$\cr
}} }\sum_{j_1,\ldots,j_k=0}^{p}
C_{j_k\ldots j_1}
\prod\limits_{s=1}^r{\bf 1}_{\{j_{g_{{}_{2s-1}}}=~j_{g_{{}_{2s}}}\}}\times
$$

\vspace{6mm}
$$
\times
J'[\phi_{j_{q_1}}\ldots \phi_{j_{q_{k-2r}}}]_{T,t}^{(i_{q_1}\ldots i_{q_{k-2r}})}=
\prod\limits_{s=1}^r
{\bf 1}_{\{i_{g_{{}_{2s-1}}}=~i_{g_{{}_{2s}}}\ne 0\}}\times
$$

\vspace{2mm}
$$
\times\hbox{\vtop{\offinterlineskip\halign{
\hfil#\hfil\cr
{\rm l.i.m.}\cr
$\stackrel{}{{}_{p\to \infty}}$\cr
}} }\hspace{-2.5mm}
\sum\limits_{\stackrel{j_1,\ldots,j_q,\ldots,j_k=0}{{}_{q\ne g_1, g_2,\ldots, g_{2r-1}, g_{2r}}}}^p
\hspace{-2.5mm}
\frac{1}{2^r}
C_{j_k \ldots j_1}\biggl|_{(j_{g_2} j_{g_1})\curvearrowright (\cdot)
\ldots (j_{g_{2r}} j_{g_{2r-1}})\curvearrowright (\cdot),
j_{g_{{}_{1}}}=~j_{g_{{}_{2}}},\ldots, j_{g_{{}_{2r-1}}}=~j_{g_{{}_{2r}}}
}\biggr.
\times
$$

\vspace{5mm}
$$
\times\prod\limits_{s=1}^r
{\bf 1}_{\{g_{2s}=g_{2s-1}+1\}}
J'[\phi_{j_{q_1}}\ldots \phi_{j_{q_{k-2r}}}]_{T,t}^{(i_{q_1}\ldots i_{q_{k-2r}})}
+
$$

\vspace{3mm}
\begin{equation}
\label{after903}
+\prod\limits_{s=1}^r
{\bf 1}_{\{i_{g_{{}_{2s-1}}}=~i_{g_{{}_{2s}}}\ne 0\}}
\hbox{\vtop{\offinterlineskip\halign{
\hfil#\hfil\cr
{\rm l.i.m.}\cr
$\stackrel{}{{}_{p\to \infty}}$\cr
}} } R_{T,t}^{(p)r,g_1,g_2,\ldots,g_{2r-1},g_{2r}}=
\end{equation}

\vspace{3mm}
\begin{equation}
\label{after904}
=\frac{1}{2^r}\prod\limits_{s=1}^r
{\bf 1}_{\{g_{2s}=g_{2s-1}+1\}}
J[\psi^{(k)}]_{T,t}^{s_r, \ldots, s_1} + 
\prod\limits_{s=1}^r
{\bf 1}_{\{i_{g_{{}_{2s-1}}}=~i_{g_{{}_{2s}}}\ne 0\}}
\hbox{\vtop{\offinterlineskip\halign{
\hfil#\hfil\cr
{\rm l.i.m.}\cr
$\stackrel{}{{}_{p\to \infty}}$\cr
}} } R_{T,t}^{(p)r,g_1,g_2,\ldots,g_{2r-1},g_{2r}}
\end{equation}

\vspace{6mm}
\noindent
w.~p.~1, where $g_{2i-1}\stackrel{\sf def}{=}s_i;$\ $i=1,2,\ldots,r;$\
$r=1,2,\ldots,\left[k/2\right],$ 
$(s_r,\ldots,s_1)\in {\rm A}_{k,r},$ $J[\psi^{(k)}]_{T,t}^{s_r,\ldots,s_1}$ is
defined by (\ref{30.1}) and ${\rm A}_{k,r}$ is defined by (\ref{30.5550001}),

\vspace{1mm}
$$
R_{T,t}^{(p)r,g_1,g_2,\ldots,g_{2r-1},g_{2r}}=
\sum\limits_{\stackrel{j_1,\ldots,j_q,\ldots,j_k=0}{{}_{q\ne g_1, g_2, \ldots, g_{2r-1}, g_{2r}}}}^p
\left(
(-1)^r \bar C^{(p)}_{j_k\ldots j_q \ldots j_1}\biggl|_{q\ne g_1,g_2, \ldots, g_{2r-1}, g_{2r}}+\right.
$$

\vspace{4mm}
$$
+(-1)^{r-1}\sum\limits_{l_1=1}^r S_{l_1}\left\{
\bar C^{(p)}_{j_k\ldots j_q \ldots j_1}\biggl|_{q\ne g_1,g_2, \ldots, g_{2r-1}, g_{2r}}\right\}+
$$

\vspace{5mm}
$$
+(-1)^{r-2}\sum\limits_{\stackrel{l_1,l_2=1}{{}_{l_1>l_2}}}^r
S_{l_1}S_{l_2}\left\{
\bar C^{(p)}_{j_k\ldots j_q \ldots j_1}\biggl|_{q\ne g_1,g_2, \ldots, g_{2r-1}, g_{2r}}\right\}+
$$

\vspace{1mm}
$$
\ldots
$$

$$
\left.+(-1)^{1}\sum\limits_{\stackrel{l_1,l_2,\ldots, l_{r-1}=1}{{}_{l_1>l_2>\ldots > l_{r-1}}}}^r
S_{l_1}S_{l_2}\ldots S_{l_{r-1}}\left\{
\bar C^{(p)}_{j_k\ldots j_q \ldots j_1}\biggl|_{q\ne g_1,g_2, \ldots, g_{2r-1}, g_{2r}}\right\}
\right)\times
$$

\vspace{4mm}
\begin{equation}
\label{afterr1}
\times
J'[\phi_{j_{q_1}}\ldots \phi_{j_{q_{k-2r}}}]_{T,t}^{(i_{q_1}\ldots i_{q_{k-2r}})}.
\end{equation}

\vspace{8mm}

Let us explain the transition from (\ref{after903}) to (\ref{after904}).
We have for $g_2=g_1+1,\ldots,g_{2r}=g_{2r-1}+1$

\vspace{1mm}
$$
\hbox{\vtop{\offinterlineskip\halign{
\hfil#\hfil\cr
{\rm l.i.m.}\cr
$\stackrel{}{{}_{p\to \infty}}$\cr
}} }
\sum\limits_{\stackrel{j_1,\ldots,j_q,\ldots,j_k=0}{{}_{q\ne g_1, g_2,\ldots, g_{2r-1}, g_{2r}}}}^p
\frac{1}{2^r}
C_{j_k \ldots j_1}\biggl|_{(j_{g_2} j_{g_1})\curvearrowright (\cdot)
\ldots (j_{g_{2r}} j_{g_{2r-1}})\curvearrowright (\cdot),
j_{g_{{}_{1}}}=~j_{g_{{}_{2}}},\ldots, j_{g_{{}_{2r-1}}}=~j_{g_{{}_{2r}}}}\biggr.
\times 
$$

\vspace{4mm}
$$
\times
\prod\limits_{s=1}^r
{\bf 1}_{\{i_{g_{{}_{2s-1}}}=~i_{g_{{}_{2s}}}\ne 0\}}
J'[\phi_{j_{q_1}}\ldots \phi_{j_{q_{k-2r}}}]_{T,t}^{(i_{q_1}\ldots i_{q_{k-2r}})}=
$$

\vspace{5mm}
$$
=\frac{1}{2^r}\hbox{\vtop{\offinterlineskip\halign{
\hfil#\hfil\cr
{\rm l.i.m.}\cr
$\stackrel{}{{}_{p\to \infty}}$\cr
}} }
\sum\limits_{\stackrel{j_1,\ldots,j_q,\ldots,j_k=0}{{}_{q\ne g_1, g_2,\ldots, g_{2r-1}, g_{2r}}}}^p
C_{j_k \ldots j_1}\biggl|_{(j_{g_2} j_{g_1})\curvearrowright 0
\ldots (j_{g_{2r}} j_{g_{2r-1}})\curvearrowright 0,
j_{g_{{}_{1}}}=~j_{g_{{}_{2}}},\ldots, j_{g_{{}_{2r-1}}}=~j_{g_{{}_{2r}}}}\biggr.
\times 
$$

\vspace{4mm}
$$
\times\prod\limits_{s=1}^r
{\bf 1}_{\{i_{g_{{}_{2s-1}}}=~i_{g_{{}_{2s}}}\ne 0\}}
\left(\zeta_{0}^{(0)}\right)^r 
J'[\phi_{j_{q_1}}\ldots \phi_{j_{q_{k-2r}}}]_{T,t}^{(i_{q_1}\ldots i_{q_{k-2r}})}=
$$

\vspace{7mm}
$$
=\frac{1}{2^r}\hbox{\vtop{\offinterlineskip\halign{
\hfil#\hfil\cr
{\rm l.i.m.}\cr
$\stackrel{}{{}_{p\to \infty}}$\cr
}} }
\sum\limits_{\stackrel{j_1,\ldots,j_q,\ldots,j_k=0}{{}_{q\ne g_1, g_2,\ldots, g_{2r-1}, g_{2r}}}}^p
\sum_{j_{m_1}, j_{m_3}\ldots,j_{m_{2r-1}}=0}^{p}
\prod\limits_{s=1}^r
{\bf 1}_{\{i_{g_{{}_{2s-1}}}=~i_{g_{{}_{2s}}}\ne 0\}}
\times
$$

\vspace{3mm}
$$
\times
C_{j_k \ldots j_1}\biggl|_{(j_{g_2} j_{g_1})\curvearrowright j_{m_1}
\ldots (j_{g_{2r}} j_{g_{2r-1}})\curvearrowright  j_{m_{2r-1}},
j_{g_{{}_{1}}}=~j_{g_{{}_{2}}},\ldots, j_{g_{{}_{2r-1}}}=~j_{g_{{}_{2r}}}}\biggr.
\times 
$$

\vspace{7mm}
$$
\times\zeta_{ j_{m_{1}}}^{(0)} \zeta_{ j_{m_{3}}}^{(0)}\ldots \zeta_{ j_{m_{2r-1}}}^{(0)}
J'[\phi_{j_{q_1}}\ldots \phi_{j_{q_{k-2r}}}]_{T,t}^{(i_{q_1}\ldots i_{q_{k-2r}})}=
$$

\vspace{7mm}
$$
=\frac{1}{2^r}\hbox{\vtop{\offinterlineskip\halign{
\hfil#\hfil\cr
{\rm l.i.m.}\cr
$\stackrel{}{{}_{p\to \infty}}$\cr
}} }
\sum\limits_{\stackrel{j_1,\ldots,j_q,\ldots,j_k=0}{{}_{q\ne g_1, g_2,\ldots, g_{2r-1}, g_{2r}}}}^p
\sum_{j_{m_1}, j_{m_3}\ldots,j_{m_{2r-1}}=0}^{p}
\prod\limits_{s=1}^r
{\bf 1}_{\{i_{g_{{}_{2s-1}}}=~i_{g_{{}_{2s}}}\ne 0\}}\times
$$

\vspace{5mm}
$$
\times
C_{j_k \ldots j_1}\biggl|_{(j_{g_2} j_{g_1})\curvearrowright j_{m_1}
\ldots (j_{g_{2r}} j_{g_{2r-1}})\curvearrowright  j_{m_{2r-1}},
j_{g_{{}_{1}}}=~j_{g_{{}_{2}}},\ldots, j_{g_{{}_{2r-1}}}=~j_{g_{{}_{2r}}}}\biggr.
\times
$$

\vspace{6mm}
\begin{equation}
\label{after905}
\times
J'[\phi_{ j_{m_1}}\phi_{ j_{m_3}}
\ldots \phi_{j_{m_{2r-1}}}
\phi_{j_{q_1}}\ldots \phi_{j_{q_{k-2r}}}]_{T,t}^{(00\ldots 0 i_{q_1}\ldots i_{q_{k-2r}})}=
\end{equation}

\vspace{4mm}
\begin{equation}
\label{after906}
=\frac{1}{2^r}
J[\psi^{(k)}]_{T,t}^{s_r, \ldots, s_1}\ \ \ \hbox{w.~p.~1}.
\end{equation}

\vspace{6mm}
\noindent
The transition from (\ref{after905}) to (\ref{after906}) is based
on (\ref{tyyyarg}).

Note that

\vspace{-1mm}
$$
C_{j_k \ldots j_1}\biggl|_{(j_{g_2} j_{g_1})\curvearrowright j_{m_1},
j_{g_{{}_{1}}}=~j_{g_{{}_{2}}}}\biggr.
=
C_{j_k \ldots j_1}\biggl|_{(j_{g_1} j_{g_1})\curvearrowright j_{m_1},
j_{g_{{}_{1}}}=~j_{g_{{}_{2}}}}\biggr.
$$

\vspace{5mm}
\noindent
is the Fourier coefficient, where $g_{2}=g_{1}+1$. 
Therefore, the value

$$
C_{j_k \ldots j_1}\biggl|_{(j_{g_2} j_{g_1})\curvearrowright j_{m_1}
\ldots (j_{g_{2d}} j_{g_{2d-1}})\curvearrowright j_{m_{2d-1}},
j_{g_{{}_{1}}}=~j_{g_{{}_{2}}},\ldots,
j_{g_{{}_{2d-1}}}=~j_{g_{{}_{2d}}}}\biggr.=
$$

\vspace{1mm}
$$
=C_{j_k \ldots j_1}\biggl|_{(j_{g_1} j_{g_1})\curvearrowright j_{m_1}
\ldots (j_{g_{2d-1}} j_{g_{2d-1}})\curvearrowright j_{m_{2d-1}},
j_{g_{{}_{1}}}=~j_{g_{{}_{2}}},\ldots,
j_{g_{{}_{2d-1}}}=~j_{g_{{}_{2d}}}}\biggr.
$$

\vspace{6mm}
\noindent
is determined recursively using (\ref{after2000}) in an obvious way
for $g_2=g_1+1,$ $\ldots,$ $g_{2d}=g_{2d-1}+1$ and $d=2,\ldots,r$.

By Condition~3 of Theorem~13 we have (also see the property (\ref{wiener1}) of
multiple Wiener stochastic integral)

\vspace{1mm}
$$
\lim\limits_{p\to\infty}
{\sf M}\left\{\left(R_{T,t}^{(p)r,g_1,g_2,\ldots,g_{2r-1},g_{2r}}\right)^2\right\}
\le 
$$

\vspace{2mm}
$$
\le
K \lim\limits_{p\to\infty}
\sum\limits_{\stackrel{j_1,\ldots,j_q,\ldots,j_k=0}{{}_{q\ne g_1, g_2, \ldots, g_{2r-1}, g_{2r}}}}^p
\left(\left(
\bar C^{(p)}_{j_k\ldots j_q \ldots j_1}\biggl|_{q\ne g_1,g_2, \ldots, g_{2r-1}, g_{2r}}
\right)^2+\right.
$$

\vspace{3mm}
$$
+\sum\limits_{l_1=1}^r \left(S_{l_1}\left\{
\bar C^{(p)}_{j_k\ldots j_q \ldots j_1}\biggl|_{q\ne g_1,g_2, \ldots, g_{2r-1}, g_{2r}}\right\}
\right)^2+
$$

\vspace{4mm}
$$
+\sum\limits_{\stackrel{l_1,l_2=1}{{}_{l_1>l_2}}}^r
\left(S_{l_1}S_{l_2}\left\{
\bar C^{(p)}_{j_k\ldots j_q \ldots j_1}\biggl|_{q\ne g_1,g_2, \ldots, g_{2r-1}, g_{2r}}\right\}
\right)^2+
$$

\vspace{1mm}
$$
\ldots
$$

\vspace{-1mm}

$$
\left.+\sum\limits_{\stackrel{l_1,l_2,\ldots, l_{r-1}=1}{{}_{l_1>l_2>\ldots > l_{r-1}}}}^r
\left(S_{l_1}S_{l_2}\ldots S_{l_{r-1}}\left\{
\bar C^{(p)}_{j_k\ldots j_q \ldots j_1}\biggl|_{q\ne g_1,g_2, \ldots, g_{2r-1}, g_{2r}}\right\}
\right)^2
\right)=0,
$$

\vspace{7mm}
\noindent
where constant $K$ does not depend on $p$.

So we have

\vspace{-1mm}
$$
\prod\limits_{s=1}^r
{\bf 1}_{\{i_{g_{{}_{2s-1}}}=~i_{g_{{}_{2s}}}\ne 0\}}\hbox{\vtop{\offinterlineskip\halign{
\hfil#\hfil\cr
{\rm l.i.m.}\cr
$\stackrel{}{{}_{p\to \infty}}$\cr
}} }\sum_{j_1,\ldots,j_k=0}^{p}
C_{j_k\ldots j_1}
\prod\limits_{s=1}^r{\bf 1}_{\{j_{g_{{}_{2s-1}}}=~j_{g_{{}_{2s}}}\}}\times
$$

\vspace{4mm}
$$
\times
J'[\phi_{j_{q_1}}\ldots \phi_{j_{q_{k-2r}}}]_{T,t}^{(i_{q_1}\ldots i_{q_{k-2r}})}=
$$

\vspace{2mm}
\begin{equation}
\label{after801}
=\frac{1}{2^r}\prod\limits_{s=1}^r
{\bf 1}_{\{g_{2s}=g_{2s-1}+1\}}
J[\psi^{(k)}]_{T,t}^{s_r, \ldots, s_1}\ \ \ \hbox{w.~p.~1},
\end{equation}

\vspace{4mm}
\noindent
where $g_{2i-1}\stackrel{\sf def}{=}s_i;$\ $i=1,2,\ldots,r;$\
$r=1,2,\ldots,\left[k/2\right],$ 
$(s_r,\ldots,s_1)\in {\rm A}_{k,r},$ $J[\psi^{(k)}]_{T,t}^{s_r,\ldots,s_1}$ is
defined by (\ref{30.1}) and ${\rm A}_{k,r}$ is defined by (\ref{30.5550001}).

Note that

\vspace{1mm}
$$
\sum_{\stackrel{(\{\{g_1, g_2\}, \ldots, 
\{g_{2r-1}, g_{2r}\}\}, \{q_1, \ldots, q_{k-2r}\})}
{{}_{\{g_1, g_2, \ldots, 
g_{2r-1}, g_{2r}, q_1, \ldots, q_{k-2r}\}=\{1, 2, \ldots, k\}}}}
\Biggl|_{g_2=g_1+1, g_3=g_2+1,\ldots, g_{2r}=g_{2r-1}+1}\Biggr.
A_{g_1,g_3,\ldots,g_{2r-1}}=
$$

\vspace{2mm}
\begin{equation}
\label{after800}
=\sum\limits_{(s_r,\ldots,s_1)\in {\rm A}_{k,r}}
A_{s_1,s_2,\ldots,s_r},
\end{equation}

\vspace{4mm}
\noindent
where $A_{g_1,g_3,\ldots,g_{2r-1}},$ 
$A_{s_1,s_2,\ldots,s_r}$ are scalar values,
$g_{2i-1}=s_i;$\ $i=1,2,\ldots,r;$\ $r=1,2,\ldots,\left[k/2\right],$
${\rm A}_{k,r}$ is defined by (\ref{30.5550001}):

$$
{\rm A}_{k,r}
=\bigl\{(s_r,\ldots,s_1):\
s_r>s_{r-1}+1,\ldots,s_2>s_1+1,\ s_r,\ldots,s_1=1,\ldots,k-1\bigr\}.
$$

\vspace{5mm}

Using (\ref{after501}), (\ref{after801}), (\ref{after800}), and Theorem~5,
we finally get                                                           

\vspace{0.5mm}
$$
\hbox{\vtop{\offinterlineskip\halign{
\hfil#\hfil\cr
{\rm l.i.m.}\cr
$\stackrel{}{{}_{p\to \infty}}$\cr
}} }\sum_{j_1,\ldots,j_k=0}^{p}
C_{j_k\ldots j_1}
\prod\limits_{l=1}^k \zeta_{j_l}^{(i_l)}=
\hbox{\vtop{\offinterlineskip\halign{
\hfil#\hfil\cr
{\rm l.i.m.}\cr
$\stackrel{}{{}_{p\to \infty}}$\cr
}} }\sum_{j_1,\ldots,j_k=0}^{p}
C_{j_k\ldots j_1}
\zeta_{j_1}^{(i_1)}\ldots \zeta_{j_k}^{(i_k)}
=
$$

\vspace{3mm}
\begin{equation}
\label{afteru11}
=J[\psi^{(k)}]_{T,t}^{(i_1\ldots i_k)}+
\sum_{r=1}^{\left[k/2\right]}\frac{1}{2^r}
\sum_{(s_r,\ldots,s_1)\in {\rm A}_{k,r}}
J[\psi^{(k)}]_{T,t}^{s_r, \ldots, s_1}=J^{*}[\psi^{(k)}]_{T,t}^{(i_1\ldots i_k)}
\end{equation}

\vspace{5mm}
\noindent
w.~p.~1, where (see (\ref{30.1}))

\vspace{-1mm}
$$
J[\psi^{(k)}]_{T,t}^{s_r, \ldots, s_1} \stackrel{\rm def}{=}\ 
\prod_{q=1}^r {\bf 1}_{\{i_{s_q}=
i_{s_{q}+1}\ne 0\}}\ \times
$$

\vspace{2mm}
$$
\times
\int\limits_t^T\psi_k(t_k)\ldots \int\limits_t^{t_{s_r+3}}
\psi_{s_r+2}(t_{s_r+2})
\int\limits_t^{t_{s_r+2}}\psi_{s_r}(t_{s_r+1})
\psi_{s_r+1}(t_{s_r+1}) \times
$$

\vspace{2mm}
$$
\times
\int\limits_t^{t_{s_r+1}}\psi_{s_r-1}(t_{s_r-1})
\ldots
\int\limits_t^{t_{s_1+3}}\psi_{s_1+2}(t_{s_1+2})
\int\limits_t^{t_{s_1+2}}\psi_{s_1}(t_{s_1+1})
\psi_{s_1+1}(t_{s_1+1}) \times
$$

\vspace{2mm}
$$
\times
\int\limits_t^{t_{s_1+1}}\psi_{s_1-1}(t_{s_1-1})
\ldots \int\limits_t^{t_2}\psi_1(t_1)
d{\bf w}_{t_1}^{(i_1)}\ldots d{\bf w}_{t_{s_1-1}}^{(i_{s_1-1})}
dt_{s_1+1}d{\bf w}_{t_{s_1+2}}^{(i_{s_1+2})}\ldots
$$

\vspace{2mm}
\begin{equation}
\label{afterito1}
\ldots\
d{\bf w}_{t_{s_r-1}}^{(i_{s_r-1})}
dt_{s_r+1}d{\bf w}_{t_{s_r+2}}^{(i_{s_r+2})}\ldots d{\bf w}_{t_k}^{(i_k)}.
\end{equation}

\vspace{5mm}

Theorem~13 is proved.

Let us make a number of remarks about Theorem~{\rm 13}.
An expansion similar to {\rm (\ref{after1})}
was obtained in {\rm \cite{Rybakov3000}}, where
the author used the definition (\ref{january19})
of the Stratonovich stochastic integral, which differs 
from the definition we use in this article
{\rm \cite{KlPl2}}.
The proof from {\rm \cite{Rybakov3000}} is somewhat simpler than the
proof proposed in this work. 
However$,$ the results from {\rm \cite{Rybakov3000}} 
were obtained under 
the condition of convergence of 
trace series. 
The verification of this condition for the kernel {\rm (\ref{ppp})} is a separate 
problem. 
In our proof, 
we essentially use the structure of the Fourier coefficients {\rm (\ref{after1000})}
corresponding to the kernel $K(t_1,\ldots,t_k)$ of the form {\rm (\ref{ppp}).} 
This circumstance actually made it possible
to prove Theorem~{\rm 13} using not the condition of
finiteness of trace series, but using the condition 
of convergence to zero of explicit expressions for the remainders
of the mentioned series. This leaves hope 
that it is possible to estimate the rate of convergence in Theorem
{\rm 13}
(see Theorems~19--22 below).

Note that under the conditions of Theorem~{\rm 13}
the sequential order of the series {\rm (}also see {\rm (\ref{after80xx}), (\ref{after500}))}

\vspace{-2mm}
$$
\sum\limits_{j_{g_{2r-1}}=p+1}^{\infty}
\sum\limits_{j_{g_{2r-3}}=p+1}^{\infty}
\ldots \sum\limits_{j_{g_{3}}=p+1}^{\infty}
\sum\limits_{j_{g_{1}}=p+1}^{\infty}
$$

\vspace{3mm}
\noindent
is not important.

We also note that the first and second conditions 
of Theorem~{\rm 13} are satisfied for complete
orthonormal systems of Legendre polynomials 
and trigonometric functions in the space
$L_2([t, T])$ (see the proofs of 
Theorems 16--18 below).
Moreover$,$ the equality (\ref{after200}) is true for an arbitrary 
basis in $L_2([t, T])$ 
(see \cite{2018a}, Sect.~2.1.4 or \cite{rybakov5000}).
Note that in the proofs of 
Theorems~6--12, 16--18, 23
the conditions of Theorem~{\rm 13}
are verified for various special cases
of iterated Stratonovich stochastic integrals
of multiplicities {\rm 2--6} with respect to 
the components of the multidimensional Wiener process.

It should be noted that (see (\ref{afterr1}))

\vspace{-1mm}
$$
(-1)^r \bar C^{(p)}_{j_k\ldots j_q \ldots j_1}\biggl|_{q\ne g_1,g_2, \ldots, g_{2r-1}, g_{2r}}+
$$

\vspace{3mm}
$$
+(-1)^{r-1}\sum\limits_{l_1=1}^r S_{l_1}\left\{
\bar C^{(p)}_{j_k\ldots j_q \ldots j_1}\biggl|_{q\ne g_1,g_2, \ldots, g_{2r-1}, g_{2r}}\right\}+
$$

\vspace{4mm}
$$
+(-1)^{r-2}\sum\limits_{\stackrel{l_1,l_2=1}{{}_{l_1>l_2}}}^r
S_{l_1}S_{l_2}\left\{
\bar C^{(p)}_{j_k\ldots j_q \ldots j_1}\biggl|_{q\ne g_1,g_2, \ldots, g_{2r-1}, g_{2r}}\right\}+
$$

\vspace{1mm}
$$
\ldots
$$

\vspace{1mm}
$$
+(-1)^{1}\sum\limits_{\stackrel{l_1,l_2,\ldots, l_{r-1}=1}{{}_{l_1>l_2>\ldots > l_{r-1}}}}^r
S_{l_1}S_{l_2}\ldots S_{l_{r-1}}\left\{
\bar C^{(p)}_{j_k\ldots j_q \ldots j_1}\biggl|_{q\ne g_1,g_2, \ldots, g_{2r-1}, g_{2r}}\right\}
=
$$

\vspace{4mm}
$$
=\sum\limits_{j_{g_1}, j_{g_3},\ldots, j_{g_{2r-1}}=0}^p
C_{j_k\ldots j_1}\biggl|_{j_{g_1}=j_{g_2},\ldots, j_{g_{2r-1}}=j_{g_{2r}}}-
$$

\vspace{2mm}
\begin{equation}
\label{drdr1000}
-\frac{1}{2^r} \prod\limits_{l=1}^r {\bf 1}_{\{g_{2l}=g_{2l-1}+1\}}
C_{j_k \ldots j_1}\biggl|_{(j_{g_2} j_{g_1})\curvearrowright (\cdot)
\ldots (j_{g_{2r}} j_{g_{2r-1}})\curvearrowright (\cdot),
j_{g_{{}_{1}}}=~j_{g_{{}_{2}}},\ldots, j_{g_{{}_{2r-1}}}=~j_{g_{{}_{2r}}}
}\biggr.,
\end{equation}

\vspace{5mm}
\noindent
where the meaning of the notations used in (\ref{afterr1}) is preserved.

For example, from (\ref{drdr1000}) for the case $r=2$ we get

\vspace{1mm}
$$
\sum\limits_{j_{g_3}=p+1}^{\infty}\sum\limits_{j_{g_1}=p+1}^{\infty}
C_{j_k \ldots j_1}\biggl|_{j_{g_{{}_{1}}}=~j_{g_{{}_{2}}}, j_{g_{{}_{3}}}=~j_{g_{{}_{4}}}}-
$$

\vspace{3mm}
$$
-
\frac{1}{2}{\bf 1}_{\{g_{4}=g_{3}+1\}}\sum\limits_{j_{g_1}=p+1}^{\infty}
C_{j_k \ldots j_1}\biggl|_{(j_{g_4} j_{g_3})\curvearrowright (\cdot),
j_{g_{{}_{1}}}=~j_{g_{{}_{2}}}, j_{g_{{}_{3}}}=~j_{g_{{}_{4}}}}
-
$$

\vspace{3mm}
$$
-
\frac{1}{2}{\bf 1}_{\{g_{2}=g_{1}+1\}}\sum\limits_{j_{g_3}=p+1}^{\infty}
C_{j_k \ldots j_1}\biggl|_{(j_{g_2} j_{g_1})\curvearrowright (\cdot),
j_{g_{{}_{1}}}=~j_{g_{{}_{2}}}, j_{g_{{}_{3}}}=~j_{g_{{}_{4}}}}\biggr.
=
$$

\vspace{3mm}
$$
=\sum\limits_{j_{g_1}=0}^p \sum\limits_{j_{g_{3}}=0}^p
C_{j_k\ldots j_1}\biggl|_{j_{g_1}=j_{g_2},j_{g_{3}}=j_{g_{4}}}-
$$

\vspace{3mm}
$$
-\frac{1}{4}{\bf 1}_{\{g_{2}=g_{1}+1\}}{\bf 1}_{\{g_{4}=g_{3}+1\}}
C_{j_k \ldots j_1}\biggl|_{(j_{g_2} j_{g_1})\curvearrowright (\cdot)
(j_{g_4} j_{g_3})\curvearrowright (\cdot),
j_{g_{{}_{1}}}=~j_{g_{{}_{2}}}, j_{g_{{}_{3}}}=~j_{g_{{}_{4}}}}\biggr..
$$

\vspace{6mm}

As a result, Condition 3 of Theorem~13 can be replaced by a weaker con\-di\-ti\-on

\vspace{1mm}
$$
\lim\limits_{p\to\infty}
\sum\limits_{\stackrel{j_1,\ldots,j_q,\ldots,j_k=0}{{}_{q\ne g_1, g_2, \ldots, g_{2r-1},
g_{2r}}}}^p
\Biggl(\sum\limits_{j_{g_1}, j_{g_3},\ldots, j_{g_{2r-1}}=0}^p
C_{j_k\ldots j_1}\biggl|_{j_{g_1}=j_{g_2},\ldots, j_{g_{2r-1}}=j_{g_{2r}}}-\Biggr.
$$

\vspace{3mm}
\begin{equation}
\label{drdr1001}
\Biggl.-\frac{1}{2^r} \prod\limits_{l=1}^r {\bf 1}_{\{g_{2l}=g_{2l-1}+1\}}
C_{j_k \ldots j_1}\biggl|_{(j_{g_2} j_{g_1})\curvearrowright (\cdot)
\ldots (j_{g_{2r}} j_{g_{2r-1}})\curvearrowright (\cdot),
j_{g_{{}_{1}}}=~j_{g_{{}_{2}}},\ldots, j_{g_{{}_{2r-1}}}=~j_{g_{{}_{2r}}}
}\biggr.\Biggr)^2=0,
\end{equation}

\vspace{6mm}
\noindent
where $r=1, 2,\ldots,[k/2]$.

However, Condition~3 of Theorem~13 itself contains 
a way of proving of the condition (\ref{drdr1001}), which is partially 
realized in the proof of Theorems~16--18, 23 (see below). 

In fact, when proving Theorem~18 (the case $r=3$
is proved in Theorem~23 for $\psi_1(\tau),\ldots,\psi_6(\tau)\equiv 1$), 
we proved the following equality

\vspace{1mm}
$$
\lim\limits_{p\to\infty}
\sum\limits_{j_{g_1}=0}^p \sum\limits_{j_{g_{3}}=0}^p
C_{j_k\ldots j_1}\biggl|_{j_{g_1}=j_{g_2},j_{g_{3}}=j_{g_{4}}}=
$$

\vspace{2mm}
$$
=\frac{1}{4}{\bf 1}_{\{g_{2}=g_{1}+1\}}{\bf 1}_{\{g_{4}=g_{3}+1\}}
C_{j_k \ldots j_1}\biggl|_{(j_{g_2} j_{g_1})\curvearrowright (\cdot)
(j_{g_4} j_{g_3})\curvearrowright (\cdot),
j_{g_{{}_{1}}}=~j_{g_{{}_{2}}}, j_{g_{{}_{3}}}=~j_{g_{{}_{4}}}}\biggr..
$$

\vspace{6mm}

On the other hand, iterative application of (\ref{after500}) gives

\vspace{1mm}
$$
\sum\limits_{j_{g_1}=0}^{\infty} \sum\limits_{j_{g_3}=0}^{\infty}\ldots \sum\limits_{j_{g_{2r-1}}=0}^{\infty}
C_{j_k\ldots j_1}\biggl|_{j_{g_1}=j_{g_2},\ldots, j_{g_{2r-1}}=j_{g_{2r}}}=
$$

\vspace{2mm}
$$
=\frac{1}{2^r} \prod\limits_{l=1}^r {\bf 1}_{\{g_{2l}=g_{2l-1}+1\}}
C_{j_k \ldots j_1}\biggl|_{(j_{g_2} j_{g_1})\curvearrowright (\cdot)
\ldots (j_{g_{2r}} j_{g_{2r-1}})\curvearrowright (\cdot),
j_{g_{{}_{1}}}=~j_{g_{{}_{2}}},\ldots, j_{g_{{}_{2r-1}}}=~j_{g_{{}_{2r}}}
},
$$

\vspace{6mm}
\noindent
where $r=1, 2,\ldots,[k/2]$.

Taking into account 
the modification of Theorem 3 for the case
of in\-teg\-ra\-tion interval $[t, s]$ $(s\in (t, T])$ 
of iterated Ito sto\-chas\-tic integrals
(see Theorem~1.11 in \cite{2018a}, \cite{tu11} or Theorem~1.24 in \cite{2018a}), 
we can formulate an analogue of Theorem~13
for the case of integration interval $[t, s]$ $(s\in (t, T])$ 
of iterated Stratonovich stochastic integrals of multiplicity 
$k$ ($k\in \mathbb{N}$).

Denote

\vspace{-2mm}
$$
\bar C^{(p)}_{j_k\ldots j_q \ldots j_1}(s)\biggl|_{q\ne g_1,g_2,\ldots,g_{2r-1}, g_{2r}}
\stackrel{\sf def}{=}
$$

\vspace{3mm}
$$
\stackrel{\sf def}{=}
\sum\limits_{j_{g_{2r-1}}=p+1}^{\infty}
\sum\limits_{j_{g_{2r-3}}=p+1}^{\infty}
\ldots \sum\limits_{j_{g_{3}}=p+1}^{\infty}
\sum\limits_{j_{g_{1}}=p+1}^{\infty}
C_{j_k \ldots j_1}(s)\biggl|_{j_{g_1}=j_{g_2},\ldots, j_{g_{2r-1}}=j_{g_{2r}}}\biggl.
$$

\vspace{5mm}

\noindent
and introduce the following notation

\vspace{3mm}
$$
S_l \left\{\bar C^{(p)}_{j_k\ldots j_q \ldots j_1}(s)\biggl|_{q\ne g_1,g_2,\ldots,g_{2r-1}, g_{2r}}
\right\}
\stackrel{\sf def}{=}
\frac{1}{2}{\bf 1}_{\{g_{2l}=g_{2l-1}+1\}}
\sum\limits_{j_{g_{2r-1}}=p+1}^{\infty}
\sum\limits_{j_{g_{2r-3}}=p+1}^{\infty}
\ldots 
$$

\vspace{3mm}
$$
\ldots
\sum\limits_{j_{g_{2l+1}}=p+1}^{\infty}
\sum\limits_{j_{g_{2l-3}}=p+1}^{\infty}
\ldots
\sum\limits_{j_{g_{3}}=p+1}^{\infty}
\sum\limits_{j_{g_{1}}=p+1}^{\infty}
C_{j_k \ldots j_1}(s)\biggl|_{(j_{g_{2l}} j_{g_{2l-1}})\curvearrowright (\cdot),
j_{g_1}=j_{g_2},\ldots, j_{g_{2r-1}}=j_{g_{2r}}}\biggr.,
$$

\vspace{6mm}
\noindent
where $l=1,2,\ldots,r,$ 
$$
C_{j_k \ldots j_1}(s)\Biggl|_{(j_{g_{2l}} j_{g_{2l-1}})\curvearrowright (\cdot)}\Biggr.
$$

\vspace{2mm}
\noindent
is defined by analogy with (\ref{after900}), 

\vspace{-2mm}
\begin{equation}
\label{after1300}
C_{j_k \ldots j_1}(s)=\int\limits_t^s\psi_k(t_k)\phi_{j_k}(t_k)\ldots
\int\limits_t^{t_2}
\psi_1(t_1)\phi_{j_1}(t_1)
dt_1\ldots dt_k.
\end{equation}

\vspace{3mm}

{\bf Theorem~15}\ \cite{2018a}, \cite{arxiv-5}, \cite{arxiv-9}, \cite{new-art-1-xxy}.\  
{\it Assume that
the continuously differentiable func\-ti\-ons 
$\psi_l(\tau)$ $(l=1,\ldots,k)$ and 
the complete orthonormal system $\{\phi_j(x)\}_{j=0}^{\infty}$
of continuous functions $(\phi_0(x)=1/\sqrt{T-t})$ 
in the space $L_2([t, T])$ are such that the following 
conditions are satisfied{\rm :}

\vspace{2mm}

{\rm 1.}\ The equality 

\vspace{-2mm}
\begin{equation}
\label{after1301}
\frac{1}{2}
\int\limits_t^s \Phi_1(t_1)\Phi_2(t_1)dt_1
=\sum_{j_1=0}^{\infty}
\int\limits_t^s
\Phi_2(t_2)\phi_{j_1}(t_2)\int\limits_t^{t_2}
\Phi_1(t_1)\phi_{j_1}(t_1)dt_1 dt_2
\end{equation}

\vspace{3mm}
\noindent
holds for all $s\in (t, T],$ where the nonrandom functions 
$\Phi_1(\tau),$ $\Phi_2(\tau)$
are continuously differentiable on $[t, T]$
and the series on the right-hand side of {\rm (\ref{after1301})}
converges absolutely.

\vspace{2mm}

{\rm 2.}\ The estimates

\vspace{-2mm}
$$
\left|\int\limits_t^{s} \phi_{j}(\tau)\Phi_1(\tau)d\tau\right|
\le \frac{\Psi_1(s)}{j^{1/2+\alpha}},\ \ \ 
\left|\int\limits_{\tau}^s \phi_{j}(\theta)\Phi_2(\theta)d\theta\right|\le
\frac{\Psi_2(s,\tau)}{j^{1/2+\alpha}},
$$

\vspace{2mm}
$$
\left|\sum_{j=p+1}^{\infty}\int\limits_t^{s}
\Phi_2(\tau)\phi_{j}(\tau)\int\limits_t^{\tau}
\Phi_1(\theta)\phi_{j}(\theta)d\theta d\tau\right|\le \frac{\Psi_3(s)}{p^{\beta}}
$$

\vspace{4mm}
\noindent
hold for all $s, \tau$ such that $t < \tau < s < T$ and for some $\alpha, \beta >0,$ where 
$\Phi_1(\tau),$ $\Phi_2(\tau)$
are continuously differentiable nonrandom functions on $[t, T],$\ $j, p\in \mathbb{N},$
and

\vspace{-1mm}
$$
\int\limits_t^s \left|\Psi_1(\tau)
\Psi_2(s,\tau)\right| d \tau<\infty,\ \ \ 
\int\limits_t^s \left|\Psi_3(\tau)\right| d\tau<\infty
$$

\vspace{2mm}
\noindent
for all $s\in (t, T).$

\vspace{4mm}

{\rm 3.}\ The condition

$$
\lim\limits_{p\to\infty}
\sum\limits_{\stackrel{j_1,\ldots,j_q,\ldots,j_k=0}{{}_{q\ne g_1, g_2, \ldots, g_{2r-1},
g_{2r}}}}^p
\left(S_{l_1}S_{l_2}\ldots S_{l_{d}}
\left\{\bar C^{(p)}_{j_k\ldots j_q \ldots j_1}(s)\biggl|_{q\ne g_1,g_2,\ldots,g_{2r-1}, g_{2r}}
\right\}\right)^2=0
$$

\vspace{4mm}
\noindent
holds for all possible $g_1,g_2,\ldots,g_{2r-1},g_{2r}$ {\rm (}see {\rm (\ref{leto5007xxx}))}
and $l_1, l_2, \ldots, l_{d}$ such that
$l_1, l_2, \ldots, l_{d}\in \{1,2,$ $\ldots,$ $r\},$\
$l_1>l_2>\ldots >l_{d},$\ $d=0, 1, 2,\ldots, r-1,$
where $r=1, 2,\ldots,[k/2]$ and

\vspace{2mm}
$$
S_{l_1}S_{l_2}\ldots S_{l_{d}}
\left\{\bar C^{(p)}_{j_k\ldots j_q \ldots j_1}(s)\biggl|_{q\ne g_1,g_2,\ldots,g_{2r-1}, g_{2r}}
\right\}\stackrel{\sf def}{=}
\bar C^{(p)}_{j_k\ldots j_q \ldots j_1}(s)\biggl|_{q\ne g_1,g_2,\ldots,g_{2r-1}, g_{2r}}
$$

\vspace{2mm}
\noindent
for $d=0.$

\vspace{4mm}

Then, for the iterated Stratonovich stochastic integral 
of arbitrary multiplicity $k$

\vspace{-1mm}
\begin{equation}
\label{afterstr1}
J^{*}[\psi^{(k)}]_{s,t}^{(i_1\ldots i_k)}=
{\int\limits_t^{*}}^s
\psi_k(t_k) \ldots 
{\int\limits_t^{*}}^{t_2}
\psi_1(t_1) d{\bf w}_{t_1}^{(i_1)}\ldots
d{\bf w}_{t_k}^{(i_k)}
\end{equation}

\vspace{3mm}
\noindent
the following 
expansion 

\vspace{-1mm}

$$
J^{*}[\psi^{(k)}]_{s,t}^{(i_1\ldots i_k)}=
\hbox{\vtop{\offinterlineskip\halign{
\hfil#\hfil\cr
{\rm l.i.m.}\cr
$\stackrel{}{{}_{p\to \infty}}$\cr
}} }
\sum\limits_{j_1,\ldots,j_k=0}^{p}
C_{j_k \ldots j_1}(s)\prod\limits_{l=1}^k \zeta_{j_l}^{(i_l)}
$$

\vspace{5mm}
\noindent
that converges in the mean-square sense is valid, where $C_{j_k \ldots j_1}(s)$
is the Fourier coefficient {\rm (\ref{after1300}),} 
${\rm l.i.m.}$ is a limit in the mean-square sense,
$i_1, \ldots, i_k=0, 1,\ldots,m,$ $s\in (t, T),$

$$
\zeta_{j}^{(i)}=
\int\limits_t^T \phi_{j}(s) d{\bf w}_s^{(i)}
$$ 

\vspace{3mm}
\noindent
are independent standard Gaussian random variables for various 
$i$ or $j$ {\rm (}in the case when $i\ne 0${\rm )},
${\bf w}_{\tau}^{(i)}={\bf f}_{\tau}^{(i)}$ 
for $i=1,\ldots,m$ and 
${\bf w}_{\tau}^{(0)}=\tau.$}

\vspace{2mm}

The results presented below in this section show that Conditions 1 and 2 
of Theorem~13 are satisfied for complete
orthonormal systems of Legendre polynomials 
and trigonometric functions in the space
$L_2([t, T]).$ Note that Condition 1 of Theorem~13
is fulfilled for an arbitrary complete orthonormal system  
of functions in the space $L_2([t,T])$ and
$\psi_1(\tau),\ldots,\psi_k(\tau)\in L_2([t, T])$
(see recent publications \cite{2018a} (Sect.~2.1.4) or \cite{rybakov5000}).

In Sect.~2.1.2 of the monograpths \cite{2018a}--\cite{12aa-afterxxx}, the following formula is proved

\vspace{-2mm}
\begin{equation}
\label{5t}
\frac{1}{2}
\int\limits_t^T\psi_1(t_1)\psi_2(t_1)dt_1
=\sum_{j_1=0}^{\infty}
C_{j_1j_1},
\end{equation}

\vspace{1mm}
\noindent
where
$$
C_{j_1 j_1}=
\int\limits_t^T
\psi_2(t_2)\phi_{j_1}(t_2)\int\limits_t^{t_2}
\psi_1(t_1)\phi_{j_1}(t_1)dt_1 dt_2,
$$

\vspace{2mm}
\noindent
$\{\phi_j(x)\}_{j=0}^{\infty}$ is a complete orthonormal system 
of Legendre polynomials or trigonometric functions
in the space $L_2([t, T]),$ 
the functions $\psi_1(\tau)$, $\psi_2(\tau)$ are continuously 
differentiable at the interval $[t, T]$.

Moreover (see Sect.~2.1.2 of the monograpths \cite{2018a}--\cite{12aa-afterxxx}), 
the following estimate 

\vspace{-2mm}
\begin{equation}
\label{5tafter}
\left\vert \sum\limits_{j_1=p+1}^{\infty}C_{j_1 j_1}
\right\vert\le \frac{C}{p},
\end{equation}

\vspace{2mm}
\noindent
holds under the above assumptions, where constant $C$ does not depend on $p$.

The relations (\ref{5t}) and (\ref{5tafter}) have been modified  
for the Legendre polynomial system as follows
(see Sect.~2.8, 2.13 of the monograpth \cite{tu11})

\vspace{-2mm}
\begin{equation}
\label{5tafter1}
\frac{1}{2}
\int\limits_t^s\psi_1(t_1)\psi_2(t_1)dt_1
=\sum_{j_1=0}^{\infty}
C_{j_1j_1}(s),
\end{equation}

\vspace{-1mm}
\begin{equation}
\label{5tafter3}
\left\vert \sum\limits_{j_1=p+1}^{\infty}C_{j_1 j_1}(s)
\right\vert \le \frac{C}{p}\left(\frac{1}{\left(1-z^2(s)\right)^{1/4}}
+1\right),
\end{equation}

\vspace{2mm}
\noindent
where $s\in (t, T)$ ($s$ is fixed, the case $s=T$ corresponds 
to (\ref{5t}) and (\ref{5tafter})), constant $C$ does not depend $p$,
the functions $\psi_1(\tau)$, $\psi_2(\tau)$ are continuously 
differentiable at the interval $[t, T]$,

\vspace{-1mm}
$$
C_{j_1 j_1}(s)=
\int\limits_t^s
\psi_2(t_2)\phi_{j_1}(t_2)\int\limits_t^{t_2}
\psi_1(t_1)\phi_{j_1}(t_1)dt_1 dt_2,
$$

\vspace{-1mm}
\begin{equation}
\label{zz1}
z(s)=\left(s-\frac{T+t}{2}\right)\frac{2}{T-t}.
\end{equation}

\vspace{3mm}

For the trigonometric case, the estimate (\ref{5tafter3}) is replaced by \cite{2018a}, \cite{tu11}

\vspace{-1mm}
\begin{equation}
\label{5tafter4}
\left\vert \sum\limits_{j_1=p+1}^{\infty}C_{j_1 j_1}(s)
\right\vert \le \frac{C}{p},
\end{equation}

\vspace{4mm}
\noindent
where $s\in [t, T]$, constant $C$ does not depend on $p$.

Note the well known estimate for the Legendre polynomials

\vspace{-1mm}
\begin{equation}
\label{ogo23}
|P_j(y)| <\frac{K}{\sqrt{j+1}(1-y^2)^{1/4}},\ \ \ 
y\in (-1, 1),\ \ \ j\in \mathbb{N},
\end{equation}

\vspace{4mm}
\noindent
where $P_j(y)$ is the Legendre polynomial, 
constant $K$ does not depend on $y$ and $j$.

We also note the following useful estimates for the case of Legendre polynomials
(\cite{2018a}-\cite{12aa-afterxxx}, Chapters 1, 2)

\vspace{-3mm}
\begin{equation}
\label{101xx}
\left|
\int\limits_t^x\psi(\tau)\phi_{j}(\tau)d\tau
\right| <
\frac{C}{j}\Biggl(\frac{1}{(1-(z(x))^2)^{1/4}}+1\Biggr),
\end{equation}

\vspace{1mm}
\begin{equation}
\label{101xxx}
\left|
\int\limits_x^T\psi(\tau)\phi_{j}(\tau)d\tau
\right| <
\frac{C}{j}\Biggl(\frac{1}{(1-(z(x))^2)^{1/4}}+1\Biggr),
\end{equation}

\vspace{1mm}
\begin{equation}
\label{101oh}
\left|
\int\limits_v^x \psi(\tau)\phi_{j}(\tau)d\tau
\right| <
\frac{C}{j}\Biggl(\frac{1}{(1-(z(x))^2)^{1/4}}+
\frac{1}{(1-(z(v))^2)^{1/4}}+1\Biggr),
\end{equation}

\vspace{4mm}
\noindent
where $j\in\mathbb{N},$\ $z(x), z(v)\in (-1, 1),$\ $x, v\in (t, T),$ 
the function $\psi(\tau)$ is continuously 
differentiable at the interval $[t, T]$,
constant $C$ does not depend on $j$.

For the case of trigonometric functions we note the following obvious estimates

\vspace{-1mm}
\begin{equation}
\label{101xxyy}
\left|
\int\limits_t^x\psi(\tau)\phi_{j}(\tau)d\tau
\right| <
\frac{C}{j},
\end{equation}

\vspace{1mm}
\begin{equation}
\label{101xxyy1}
\left|
\int\limits_x^T\psi(\tau)\phi_{j}(\tau)d\tau
\right| <
\frac{C}{j},
\end{equation}

\vspace{1mm}
\begin{equation}
\label{101xxyy2}
\left|
\int\limits_v^x\psi(\tau)\phi_{j}(\tau)d\tau
\right| <
\frac{C}{j},
\end{equation}

\vspace{4mm}
\noindent
where $j\in\mathbb{N},$\ $x, v\in [t, T],$ 
the function $\psi(\tau)$ is continuously 
differentiable at the interval $[t, T]$,
constant $C$ does not depend on $j$.

\vspace{5mm}

\section{Expansion of Iterated Stratonovich Stochastic Integrals
of Multiplicity 3. The Case $p_1=p_2=p_3\to \infty$ and 
Continuously Differentiable 
Weight Functions $\psi_1(\tau),$ $\psi_2(\tau),$ $\psi_3(\tau)$ 
(The Cases of Legendre 
Polynomials and Trigonometric Functions)}

\vspace{5mm}

In this section, we present a simple proof of Theorem~11  
based on Theorem~13. In this case, the conditions of Theorem~11 will be weakened.

First, consider the following equalities 

\vspace{-1mm}
\begin{equation}
\label{after1400}
\frac{1}{2}
\int\limits_{t_1}^{t_2} \Phi_1(\tau)\Phi_2(\tau)d\tau
=\sum_{j=0}^{\infty}
\int\limits_{t_1}^{t_2}
\Phi_2(\tau)\phi_{j}(\tau)\int\limits_{t_1}^{\tau}
\Phi_1(\theta)\phi_{j}(\theta)d\theta d\tau,
\end{equation}

\vspace{1mm}
\begin{equation}
\label{after1401}
\frac{1}{2}
\int\limits_{t_1}^{t_2} \Phi_1(\tau)\Phi_2(\tau)d\tau
=\sum_{j=0}^{\infty}
\int\limits_{t_1}^{t_2}
\Phi_1(\theta)\phi_{j}(\theta)\int\limits_{\theta}^{t_2}
\Phi_2(\tau)\phi_{j}(\tau)d\tau d\theta
\end{equation}

\vspace{3mm}
\noindent
that will be used further, where $t\le t_1<t_2\le T,$  $\Phi_1(\tau), \Phi_2(\tau)\in L_2([t,T]),$
$\{\phi_j(x)\}_{j=0}^{\infty}$ is an arbitrary complete orthonormal system of funtions
in $L_2([t, T]).$
The equality (\ref{after1401}) has been proved in \cite{2018a} (Sect.~2.7.2).
Using (\ref{after1401}) and Fubini's Theorem, we get (\ref{after1400})
(also see \cite{rybakov5000}).

\vspace{2mm}

{\bf Theorem 16}\ \cite{2018a}, \cite{arxiv-5}, \cite{arxiv-9}, \cite{new-art-1-xxy}.\ 
{\it Suppose that 
$\{\phi_j(x)\}_{j=0}^{\infty}$ is a complete orthonormal system of 
Legendre polynomials or trigonometric func\-ti\-ons in the space $L_2([t, T]).$
Furthermore, let $\psi_1(\tau), \psi_2(\tau),$ $\psi_3(\tau)$ are continuously dif\-ferentiable 
nonrandom functions on $[t, T].$ 
Then, for the 
iterated Stra\-to\-no\-vich stochastic integral of third multiplicity

\vspace{-1mm}
\begin{equation}
\label{after1500}
J^{*}[\psi^{(3)}]_{T,t}={\int\limits_t^{*}}^T\psi_3(t_3)
{\int\limits_t^{*}}^{t_3}\psi_2(t_2)
{\int\limits_t^{*}}^{t_2}\psi_1(t_1)
d{\bf w}_{t_1}^{(i_1)}
d{\bf w}_{t_2}^{(i_2)}d{\bf w}_{t_3}^{(i_3)}
\end{equation}

\vspace{3mm}
\noindent
the following 
expansion 

\vspace{-2mm}
$$
J^{*}[\psi^{(3)}]_{T,t}
=\hbox{\vtop{\offinterlineskip\halign{
\hfil#\hfil\cr
{\rm l.i.m.}\cr
$\stackrel{}{{}_{p\to \infty}}$\cr
}} }
\sum\limits_{j_1, j_2, j_3=0}^{p}
C_{j_3 j_2 j_1}\zeta_{j_1}^{(i_1)}\zeta_{j_2}^{(i_2)}\zeta_{j_3}^{(i_3)}
$$

\vspace{5mm}
\noindent
that converges in the mean-square sense is valid, where
$i_1, i_2, i_3=0, 1,\ldots,m,$

\vspace{-1mm}
$$
C_{j_3 j_2 j_1}=\int\limits_t^T\psi_3(t_3)\phi_{j_3}(t_3)
\int\limits_t^{t_3}\psi_2(t_2)\phi_{j_2}(t_2)
\int\limits_t^{t_2}\psi_1(t_1)\phi_{j_1}(t_1)dt_1dt_2dt_3
$$

\vspace{1mm}
\noindent
and
$$
\zeta_{j}^{(i)}=
\int\limits_t^T \phi_{j}(s) d{\bf w}_s^{(i)}
$$ 

\vspace{2mm}
\noindent
are independent standard Gaussian random variables for various 
$i$ or $j$
{\rm (}in the case when $i\ne 0${\rm ),}
${\bf w}_{\tau}^{(i)}={\bf f}_{\tau}^{(i)}$ for
$i=1,\ldots,m$ and 
${\bf w}_{\tau}^{(0)}=\tau.$}

\vspace{2mm}

{\bf Proof.}\ As follows from Sect.~7, Conditions 1 and 2
of Theorem~13 are satisfied for complete
orthonormal systems of Legendre polynomials 
and trigonometric functions in the space
$L_2([t, T]).$ Let us verify Condition 3 of Theorem~13 for
the iterated Stratonovich stochastic integral (\ref{after1500}). 
Thus, we have to check the following conditions

\vspace{-2mm}
\begin{equation}
\label{after1600}
\lim\limits_{p\to \infty}
\sum_{j_3=0}^p\left(\sum_{j_1=p+1}^{\infty}
C_{j_3 j_1 j_1}\right)^2=0,
\end{equation}

\vspace{2mm}
\begin{equation}
\label{after1601}
\lim\limits_{p\to \infty}
\sum_{j_1=0}^p\left(\sum_{j_3=p+1}^{\infty}
C_{j_3 j_3 j_1}\right)^2=0,
\end{equation}

\vspace{2mm}
\begin{equation}
\label{after1602}
\lim\limits_{p\to \infty}
\sum_{j_2=0}^p\left(\sum_{j_1=p+1}^{\infty}
C_{j_1 j_2 j_1}\right)^2=0.
\end{equation}

\vspace{4mm}

We have
$$
\sum_{j_3=0}^p\left(\sum_{j_1=p+1}^{\infty}
C_{j_3 j_1 j_1}\right)^2=
$$

\vspace{2mm}
\begin{equation}
\label{after1800}
=
\sum_{j_3=0}^p\left(\sum_{j_1=p+1}^{\infty}
\int\limits_t^T\psi_3(t_3)\phi_{j_3}(t_3)
\int\limits_t^{t_3}\psi_2(t_2)\phi_{j_1}(t_2)
\int\limits_t^{t_2}\psi_1(t_1)\phi_{j_1}(t_1)dt_1dt_2dt_3\right)^2=
\end{equation}

\vspace{2mm}
\begin{equation}
\label{after1801}
=
\sum_{j_3=0}^p\left(
\int\limits_t^T\psi_3(t_3)\phi_{j_3}(t_3)
\sum_{j_1=p+1}^{\infty}\int\limits_t^{t_3}\psi_2(t_2)\phi_{j_1}(t_2)
\int\limits_t^{t_2}\psi_1(t_1)\phi_{j_1}(t_1)dt_1dt_2dt_3\right)^2\le
\end{equation}

\vspace{2mm}
\begin{equation}
\label{after1802}
\le
\sum_{j_3=0}^{\infty}\left(
\int\limits_t^T\psi_3(t_3)\phi_{j_3}(t_3)
\sum_{j_1=p+1}^{\infty}\int\limits_t^{t_3}\psi_2(t_2)\phi_{j_1}(t_2)
\int\limits_t^{t_2}\psi_1(t_1)\phi_{j_1}(t_1)dt_1dt_2dt_3\right)^2=
\end{equation}

\vspace{2mm}
\begin{equation}
\label{after1803}
=\int\limits_t^T\psi_3^2(t_3)
\left(\sum_{j_1=p+1}^{\infty}\int\limits_t^{t_3}\psi_2(t_2)\phi_{j_1}(t_2)
\int\limits_t^{t_2}\psi_1(t_1)\phi_{j_1}(t_1)dt_1dt_2\right)^2 dt_3\le
\end{equation}
\begin{equation}
\label{after1804}
\le \frac{K}{p^2}\ \to\  0
\end{equation}

\vspace{4mm}
\noindent
if $p\to\infty,$ where constant $K$ does not depend on $p.$

Note that the transition from (\ref{after1800}) to (\ref{after1801}) is based on 
the estimate (\ref{5tafter3}) for the polynomial case and its analogue
(\ref{5tafter4}) for the
trigonometric case, 
the transition from (\ref{after1802}) to (\ref{after1803}) is based on the Parseval equality, 
and the transition from (\ref{after1803}) to (\ref{after1804}) 
is also based on the 
estimate (\ref{5tafter3}) and its analogue
(\ref{5tafter4}) for the
trigonometric case.

By analogy with the previous case we have 

$$
\sum_{j_1=0}^p\left(\sum_{j_3=p+1}^{\infty}
C_{j_3 j_3 j_1}\right)^2=
$$

\vspace{2mm}
$$
=
\sum_{j_1=0}^p\left(\sum_{j_3=p+1}^{\infty}
\int\limits_t^T\psi_3(t_3)\phi_{j_3}(t_3)
\int\limits_t^{t_3}\psi_2(t_2)\phi_{j_3}(t_2)
\int\limits_t^{t_2}\psi_1(t_1)\phi_{j_1}(t_1)dt_1dt_2dt_3\right)^2=
$$

\vspace{2mm}
\begin{equation}
\label{after1900}
=
\sum_{j_1=0}^p\left(\sum_{j_3=p+1}^{\infty}
\int\limits_t^T \psi_1(t_1)\phi_{j_1}(t_1) \int\limits_{t_1}^{T} \psi_2(t_2)\phi_{j_3}(t_2)
\int\limits_{t_2}^T\psi_3(t_3)\phi_{j_3}(t_3)
dt_3dt_2dt_1\right)^2=
\end{equation}

\vspace{2mm}
\begin{equation}
\label{after1901}
=
\sum_{j_1=0}^p\left(
\int\limits_t^T \psi_1(t_1)\phi_{j_1}(t_1) 
\sum_{j_3=p+1}^{\infty}\int\limits_{t_1}^{T} \psi_2(t_2)\phi_{j_3}(t_2)
\int\limits_{t_2}^T\psi_3(t_3)\phi_{j_3}(t_3)
dt_3dt_2dt_1\right)^2\le
\end{equation}

\vspace{2mm}
$$
\le
\sum_{j_1=0}^{\infty}\left(
\int\limits_t^T \psi_1(t_1)\phi_{j_1}(t_1) 
\sum_{j_3=p+1}^{\infty}\int\limits_{t_1}^{T} \psi_2(t_2)\phi_{j_3}(t_2)
\int\limits_{t_2}^T\psi_3(t_3)\phi_{j_3}(t_3)
dt_3dt_2dt_1\right)^2=
$$

\vspace{2mm}
\begin{equation}
\label{after1902}
=
\int\limits_t^T \psi_1^2(t_1)\left(
\sum_{j_3=p+1}^{\infty}\int\limits_{t_1}^{T} \psi_2(t_2)\phi_{j_3}(t_2)
\int\limits_{t_2}^T\psi_1(t_3)\phi_{j_3}(t_3)
dt_3dt_2\right)^2 dt_1\le
\end{equation}

\vspace{2mm}
\begin{equation}
\label{after1903}
\le \frac{K}{p^2}\ \to\  0
\end{equation}

\vspace{4mm}
\noindent
if $p\to\infty,$ where constant $K$ is independent of $p.$

The transition from (\ref{after1900}) to (\ref{after1901}) is based on 
analogues of the estimates (\ref{5tafter3}), (\ref{5tafter4}) 
for the value

$$
\left|\sum_{j_3=p+1}^{\infty}\int\limits_{t_1}^{T} \psi_2(t_2)\phi_{j_3}(t_2)
\int\limits_{t_2}^T\psi_3(t_3)\phi_{j_3}(t_3)
dt_3dt_2\right|
$$

\vspace{3mm}
\noindent
for the polynomial and
trigonometric cases, 
the transition from (\ref{after1902}) to (\ref{after1903})
is also based on the mentioned analogues of the estimates 
(\ref{5tafter3}), (\ref{5tafter4}).

Further, we have 

$$
\sum_{j_2=0}^p\left(\sum_{j_1=p+1}^{\infty}
C_{j_1 j_2 j_1}\right)^2=
$$

\vspace{2mm}
$$
=
\sum_{j_2=0}^p\left(\sum_{j_1=p+1}^{\infty}
\int\limits_t^T\psi_3(t_3)\phi_{j_1}(t_3)
\int\limits_t^{t_3}\psi_2(t_2)\phi_{j_2}(t_2)
\int\limits_t^{t_2}\psi_1(t_1)\phi_{j_1}(t_1)dt_1dt_2dt_3\right)^2=
$$

\vspace{2mm}
\begin{equation}
\label{after1920}
=
\sum_{j_2=0}^p\left(\sum_{j_1=p+1}^{\infty}
\int\limits_t^T \psi_2(t_2)\phi_{j_2}(t_2) \int\limits_{t}^{t_2} \psi_1(t_1)\phi_{j_1}(t_1)
dt_1 \int\limits_{t_2}^T\psi_3(t_3)\phi_{j_1}(t_3)
dt_3dt_2\right)^2=
\end{equation}

\vspace{2mm}
\begin{equation}
\label{after1921}
=
\sum_{j_2=0}^p\left(
\int\limits_t^T \psi_2(t_2)\phi_{j_2}(t_2) 
\sum_{j_1=p+1}^{\infty}\int\limits_{t}^{t_2} \psi_1(t_1)\phi_{j_1}(t_1)
dt_1 \int\limits_{t_2}^T\psi_3(t_3)\phi_{j_1}(t_3)
dt_3dt_2\right)^2\le
\end{equation}

\vspace{2mm}
$$
\le
\sum_{j_2=0}^{\infty}\left(
\int\limits_t^T \psi_2(t_2)\phi_{j_2}(t_2) \sum_{j_1=p+1}^{\infty}
\int\limits_{t}^{t_2} \psi_1(t_1)\phi_{j_1}(t_1)
dt_1 \int\limits_{t_2}^T\psi_3(t_3)\phi_{j_1}(t_3)
dt_3dt_2\right)^2=
$$

\vspace{2mm}
\begin{equation}
\label{after1922}
=
\int\limits_t^T \psi_2^2(t_2)
\left(\sum_{j_1=p+1}^{\infty}\int\limits_{t}^{t_2} \psi_1(t_1)\phi_{j_1}(t_1)
dt_1 \int\limits_{t_2}^T\psi_3(t_3)\phi_{j_1}(t_3)
dt_3\right)^2 dt_2.
\end{equation}

\vspace{4mm}

The transition from (\ref{after1920}) to (\ref{after1921}) 
is based on the estimates (\ref{101xx}), (\ref{101xxx}) and its obvious analogues
(\ref{101xxyy}), (\ref{101xxyy1}) for the trigonometric case.
However, the estimates (\ref{101xx}), (\ref{101xxx}) cannot be used to estimate
the right-hand side of (\ref{after1922}), 
since we get the divergent integral.
For this reason, we will obtain a new estimate based on the relation \cite{2018a}-\cite{12aa-afterxxx}

$$
\int\limits_t^x\psi(s)\phi_{j_1}(s)ds=
\frac{\sqrt{T-t}\sqrt{2j_1+1}}{2}
\int\limits_{-1}^{z(x)}P_{j_1}(y)
\psi(u(y))dy=
$$

\vspace{2mm}
$$
=\frac{\sqrt{T-t}}{2\sqrt{2j_1+1}}\Biggl((P_{j_1+1}(z(x))-
P_{j_1-1}(z(x)))\psi(x)-\Biggr.
$$

\vspace{2mm}
\begin{equation}
\label{otit6000x}
\Biggl.-
\frac{T-t}{2}
\int\limits_{-1}^{z(x)}((P_{j_1+1}(y)-P_{j_1-1}(y))
{\psi}'(u(y))dy\Biggr),
\end{equation}

\vspace{4mm}
\noindent
where $x\in (t, T),$ $j_1\ge p+1,$ 
$z(x)$ is defined by (\ref{zz1}), $P_j(x)$ is the Legendre
polynomial,
${\psi}'$ is a derivative of the continuously differentiable function $\psi(s)$
with respect to the variable $u(y),$

$$
u(y)=\frac{T-t}{2}y+\frac{T+t}{2}.
$$

\vspace{4mm}

From (\ref{ogo23}) and the estimate $\left| P_j(y) \right|\le 1$, $y\in [-1, 1]$
we obtain

\vspace{-1mm}
\begin{equation}
\label{after5000}
\left|P_{j}(y)\right|=\left|P_{j}(y)\right|^{\varepsilon}
\cdot \left|P_{j}(y)\right|^{1-\varepsilon}\le \left|P_{j}(y)\right|^{1-\varepsilon}
< \frac{C}{j^{1/2-\varepsilon/2} (1-y^2)^{1/4-\varepsilon/4}},
\end{equation}

\vspace{3mm}
\noindent
where $y\in (-1, 1),$ $j\in\mathbb{N},$ and $\varepsilon$ is an arbitrary
small positive real number.

Combining (\ref{otit6000x}) and (\ref{after5000}), we have the following estimate

\vspace{-1mm}
\begin{equation}
\label{after1940}
\left|
\int\limits_t^s\psi_1(\tau)\phi_{j_1}(\tau)d\tau
\right| <
\frac{C}{(j_1)^{1-\varepsilon/2}}\Biggl(\frac{1}{(1-z^2(s))^{1/4-\varepsilon/4}}+1\Biggr),
\end{equation}

\vspace{3mm}
\noindent
where $s\in (t, T),$
$z(s)$ is defined by (\ref{zz1}), constant $C$ does not depend on $j_1.$

Similarly to (\ref{after1940}) we obtain 

\vspace{-1mm}
\begin{equation}
\label{after1941}
\left|\int\limits_s^T\psi_3(\tau)\phi_{j_1}(\tau)d\tau
\right| <
\frac{C}{(j_1)^{1-\varepsilon/2}}\Biggl(\frac{1}{(1-z^2(s))^{1/4-\varepsilon/4}}+1\Biggr),
\end{equation}

\vspace{3mm}
\noindent
where $s\in (t, T),$ constant $C$ is independent of $j_1.$

Combining (\ref{101xx}) and (\ref{after1941}), we have

\vspace{-1mm}
$$
\left|
\int\limits_t^s\psi_1(\tau)\phi_{j_1}(\tau)d\tau
\int\limits_s^T\psi_3(\tau)\phi_{j_1}(\tau)d\tau\right|
<
$$

\vspace{2mm}
\begin{equation}
\label{after4000}
<\frac{L}{(j_1)^{2-\varepsilon/2}}\Biggl(\frac{1}{(1-z^2(s))^{1/4-\varepsilon/4}}+1\Biggr)
\Biggl(\frac{1}{(1-z^2(s))^{1/4}}+1\Biggr),
\end{equation}

\vspace{4mm}
\noindent
where $s\in (t, T),$ 
$z(s)$ is defined by (\ref{zz1}), constant $L$ does not depend on $j_1.$

Observe that

\vspace{-2mm}
\begin{equation}
\label{after1944}
\sum\limits_{j_1=p+1}^{\infty}\frac{1}{(j_1)^{2-\varepsilon/2}}
\le \int\limits_{p}^{\infty}\frac{dx}{x^{2-\varepsilon/2}}=
\frac{1}{(1-\varepsilon/2)p^{1-\varepsilon/2}}.
\end{equation}

\vspace{4mm}

Applying (\ref{after4000}) and (\ref{after1944})
to estimate the right-hand side of (\ref{after1922}) gives

\vspace{-1mm}
\begin{equation}
\label{after5001}
\sum_{j_2=0}^p\left(\sum_{j_1=p+1}^{\infty}
C_{j_1 j_2 j_1}\right)^2\le \frac{K}{p^{2-\varepsilon}}\ \to\  0
\end{equation}

\vspace{3mm}
\noindent
if $p\to\infty,$ where 
$\varepsilon$ is an arbitrary
small positive real number,
constant $K$ is independent of $p$.

The estimation of the right-hand side of (\ref{after1922})
for the trigonometric case is carried out using the estimates
(\ref{101xxyy}), (\ref{101xxyy1}). At that we obtain the
estimate (\ref{after5001}) with $\varepsilon=0.$
Theorem~16 is proved.

\vspace{5mm}

\section{Expansion of Iterated Stratonovich Stochastic Integrals
of Multiplicity 4. The Case $p_1=\ldots =p_4\to \infty$ and 
Continuously Differentiable 
Weight Functions $\psi_1(\tau),$ $\ldots,$ $\psi_4(\tau)$ 
(The Cases of Legendre 
Polynomials and Trigonometric Functions)}

\vspace{5mm}  

{\bf Theorem 17}\ \cite{2018a}, \cite{arxiv-5}, \cite{arxiv-9}, \cite{new-art-1-xxy}.\ 
{\it Suppose that 
$\{\phi_j(x)\}_{j=0}^{\infty}$ is a complete orthonormal system of 
Legendre polynomials or trigonometric func\-ti\-ons in the space $L_2([t, T]).$
Furthermore, let $\psi_1(\tau), \ldots,$ $\psi_4(\tau)$ are continuously dif\-ferentiable 
nonrandom functions on $[t, T].$ 
Then, for the 
iterated Stra\-to\-no\-vich stochastic integral of fourth multiplicity

\vspace{-1mm}
\begin{equation}
\label{after2500}
J^{*}[\psi^{(4)}]_{T,t}={\int\limits_t^{*}}^T\psi_4(t_4)
{\int\limits_t^{*}}^{t_4}\psi_3(t_3)
{\int\limits_t^{*}}^{t_3}\psi_2(t_2)
{\int\limits_t^{*}}^{t_2}\psi_1(t_1)
d{\bf w}_{t_1}^{(i_1)}
d{\bf w}_{t_2}^{(i_2)}d{\bf w}_{t_3}^{(i_3)}d{\bf w}_{t_4}^{(i_4)}
\end{equation}

\vspace{3mm}
\noindent
the following 
expansion 

\vspace{-1mm}
$$
J^{*}[\psi^{(4)}]_{T,t}
=\hbox{\vtop{\offinterlineskip\halign{
\hfil#\hfil\cr
{\rm l.i.m.}\cr
$\stackrel{}{{}_{p\to \infty}}$\cr
}} }
\sum\limits_{j_1, j_2, j_3, j_4=0}^{p}
C_{j_4 j_3 j_2 j_1}\zeta_{j_1}^{(i_1)}\zeta_{j_2}^{(i_2)}\zeta_{j_3}^{(i_3)}
\zeta_{j_4}^{(i_4)}
$$

\vspace{4mm}
\noindent
that converges in the mean-square sense is valid, where
$i_1, i_2, i_3, i_4=0, 1,\ldots,m,$

\vspace{-1mm}
$$
C_{j_4 j_3 j_2 j_1}=
\int\limits_t^T\psi_4(t_4)\phi_{j_4}(t_4)
\int\limits_t^{t_4}\psi_3(t_3)\phi_{j_3}(t_3)
\int\limits_t^{t_3}\psi_2(t_2)\phi_{j_2}(t_2)
\int\limits_t^{t_2}\psi_1(t_1)\phi_{j_1}(t_1)dt_1\times
$$

$$
\times
dt_2dt_3dt_4
$$

\vspace{3mm}
\noindent
and
$$
\zeta_{j}^{(i)}=
\int\limits_t^T \phi_{j}(s) d{\bf w}_s^{(i)}
$$ 

\vspace{2mm}
\noindent
are independent standard Gaussian random variables for various 
$i$ or $j$
{\rm (}in the case when $i\ne 0${\rm ),}
${\bf w}_{\tau}^{(i)}={\bf f}_{\tau}^{(i)}$ for
$i=1,\ldots,m$ and 
${\bf w}_{\tau}^{(0)}=\tau.$}

\vspace{2mm}

{\bf Proof.}\ As follows from Sect.~7, Conditions 1 and 2
of Theorem~13 are satisfied for complete
orthonormal systems of Legendre polynomials 
and trigonometric functions in the space
$L_2([t, T]).$ Let us verify Condition 3 of Theorem~13 for
the iterated Stratonovich stochastic integral (\ref{after2500}). 
Thus, we have to check the following conditions

\vspace{-2mm}
\begin{equation}
\label{after2501}
\lim\limits_{p\to \infty}
\sum_{j_3,j_4=0}^p\left(\sum_{j_1=p+1}^{\infty}
C_{j_4 j_3 j_1 j_1}\right)^2=0,
\end{equation}

\begin{equation}
\label{after2502}
\lim\limits_{p\to \infty}
\sum_{j_2,j_4=0}^p\left(\sum_{j_1=p+1}^{\infty}
C_{j_4 j_1 j_2 j_1}\right)^2=0,
\end{equation}

\begin{equation}
\label{after2503}
\lim\limits_{p\to \infty}
\sum_{j_2,j_3=0}^p\left(\sum_{j_1=p+1}^{\infty}
C_{j_1 j_3 j_2 j_1}\right)^2=0,
\end{equation}

\begin{equation}
\label{after2504}
\lim\limits_{p\to \infty}
\sum_{j_1,j_4=0}^p\left(\sum_{j_2=p+1}^{\infty}
C_{j_4 j_2 j_2 j_1}\right)^2=0,
\end{equation}

\begin{equation}
\label{after2505}
\lim\limits_{p\to \infty}
\sum_{j_1,j_3=0}^p\left(\sum_{j_2=p+1}^{\infty}
C_{j_2 j_3 j_2 j_1}\right)^2=0,
\end{equation}

\begin{equation}
\label{after2506}
\lim\limits_{p\to \infty}
\sum_{j_1,j_2=0}^p\left(\sum_{j_3=p+1}^{\infty}
C_{j_3 j_3 j_2 j_1}\right)^2=0,
\end{equation}

\begin{equation}
\label{after2508}
\lim\limits_{p\to \infty}
\left(\sum_{j_2=p+1}^{\infty}\sum_{j_1=p+1}^{\infty}
C_{j_2 j_1 j_2 j_1}\right)^2=0,
\end{equation}

\begin{equation}
\label{after2509}
\lim\limits_{p\to \infty}
\left(\sum_{j_2=p+1}^{\infty}\sum_{j_1=p+1}^{\infty}
C_{j_1 j_2 j_2 j_1}\right)^2=0,
\end{equation}

\begin{equation}
\label{after2507}
\lim\limits_{p\to \infty}
\left(\sum_{j_3=p+1}^{\infty}\sum_{j_1=p+1}^{\infty}
C_{j_3 j_3 j_1 j_1}\right)^2=0,
\end{equation}

\begin{equation}
\label{after1602xxxx}
\lim\limits_{p\to \infty}
\left(\sum_{j_3=p+1}^{\infty}
C_{j_3 j_3 j_1 j_1}\biggl|_{(j_{1} j_{1})\curvearrowright (\cdot)}\right)^2=0,
\end{equation}

\begin{equation}
\label{after1602xxxy}
\lim\limits_{p\to \infty}
\left(\sum_{j_1=p+1}^{\infty}
C_{j_3 j_3 j_1 j_1}\biggl|_{(j_{3} j_{3})\curvearrowright (\cdot)}\right)^2=0,
\end{equation}

\begin{equation}
\label{after1602xxxz}
\lim\limits_{p\to \infty}
\left(\sum_{j_1=p+1}^{\infty}
C_{j_1 j_2 j_2 j_1}\biggl|_{(j_{2} j_{2})\curvearrowright (\cdot)}\biggr.\right)^2=0,
\end{equation}

\vspace{5mm}
\noindent
where we use the notation (\ref{after900}) in (\ref{after1602xxxx})--(\ref{after1602xxxz}).

Applying arguments similar to those we used in the proof of Theorem~16, we obtain
for (\ref{after2501})

$$
\sum_{j_3,j_4=0}^p\left(\sum_{j_1=p+1}^{\infty}
C_{j_4 j_3 j_1 j_1}\right)^2=
\sum_{j_3,j_4=0}^p\left(\sum_{j_1=p+1}^{\infty}
\int\limits_t^T\psi_4(t_4)\phi_{j_4}(t_4)
\int\limits_t^{t_4}\psi_3(t_3)\phi_{j_3}(t_3)\times\right.
$$

\begin{equation}
\label{after5010}
\left.\times
\int\limits_t^{t_3}\psi_2(t_2)\phi_{j_1}(t_2)
\int\limits_t^{t_2}\psi_1(t_1)\phi_{j_1}(t_1)dt_1dt_2dt_3dt_4\right)^2=
\end{equation}

$$
=
\sum_{j_3,j_4=0}^p\left(
\int\limits_t^T\psi_4(t_4)\phi_{j_4}(t_4)
\int\limits_t^{t_4}\psi_3(t_3)\phi_{j_3}(t_3)\times\right.
$$

\begin{equation}
\label{after5011}
\left.\times
\sum_{j_1=p+1}^{\infty}\int\limits_t^{t_3}\psi_2(t_2)\phi_{j_1}(t_2)
\int\limits_t^{t_2}\psi_1(t_1)\phi_{j_1}(t_1)dt_1dt_2dt_3dt_4\right)^2\le
\end{equation}

$$
\le
\sum_{j_3,j_4=0}^{\infty}\left(
\int\limits_t^T\psi_4(t_4)\phi_{j_4}(t_4)
\int\limits_t^{t_4}\psi_3(t_3)\phi_{j_3}(t_3)\times\right.
$$

\begin{equation}
\label{after5002}
\left.\times
\sum_{j_1=p+1}^{\infty}\int\limits_t^{t_3}\psi_2(t_2)\phi_{j_1}(t_2)
\int\limits_t^{t_2}\psi_1(t_1)\phi_{j_1}(t_1)dt_1dt_2dt_3dt_4\right)^2=
\end{equation}

\begin{equation}
\label{after5003}
=
\int\limits_{[t, T]^2}{\bf 1}_{\{t_3<t_4\}}
\psi_4^2(t_4)
\psi_3^2(t_3)\left(
\sum_{j_1=p+1}^{\infty}\int\limits_t^{t_3}\psi_2(t_2)\phi_{j_1}(t_2)
\int\limits_t^{t_2}\psi_1(t_1)\phi_{j_1}(t_1)dt_1dt_2\right)^2 dt_3dt_4\le
\end{equation}

\begin{equation}
\label{after5004}
\le \frac{K}{p^2}\ \to\  0
\end{equation}

\vspace{4mm}
\noindent
if $p\to\infty,$ where constant $K$ is independent of $p.$

Note that the transition from (\ref{after5010}) to (\ref{after5011}) is based on 
the estimate (\ref{5tafter3}) for the polynomial case and its analogue
for the
trigonometric case, 
the transition from (\ref{after5002}) to (\ref{after5003}) is based on the Parseval equality, 
and the transition from (\ref{after5003}) to (\ref{after5004}) 
is also based on the estimate (\ref{5tafter3})
and its analogue
for the
trigonometric case.

Further, we have for (\ref{after2502})

$$
\sum_{j_2,j_4=0}^p\left(\sum_{j_1=p+1}^{\infty}
C_{j_4 j_1 j_2 j_1}\right)^2=
\sum_{j_2,j_4=0}^p\left(\sum_{j_1=p+1}^{\infty}
\int\limits_t^T\psi_4(t_4)\phi_{j_4}(t_4)
\int\limits_t^{t_4}\psi_3(t_3)\phi_{j_1}(t_3)\times\right.
$$

\begin{equation}
\label{after98}
\left.\times
\int\limits_t^{t_3}\psi_2(t_2)\phi_{j_2}(t_2)
\int\limits_t^{t_2}\psi_1(t_1)\phi_{j_1}(t_1)dt_1dt_2dt_3dt_4\right)^2=
\end{equation}

$$
=
\sum_{j_2,j_4=0}^p\left(\sum_{j_1=p+1}^{\infty}
\int\limits_t^T\psi_4(t_4)\phi_{j_4}(t_4)
\int\limits_t^{t_4}\psi_2(t_2)\phi_{j_2}(t_2)\times\right.
$$

\begin{equation}
\label{after99}
\left.\times
\int\limits_t^{t_2}\psi_1(t_1)\phi_{j_1}(t_1)dt_1
\int\limits_{t_2}^{t_4}\psi_3(t_3)\phi_{j_1}(t_3)dt_3dt_2dt_4\right)^2=
\end{equation}

$$
=
\sum_{j_2,j_4=0}^p\left(
\int\limits_t^T\psi_4(t_4)\phi_{j_4}(t_4)
\int\limits_t^{t_4}\psi_2(t_2)\phi_{j_2}(t_2)\times\right.
$$

$$
\left.\times
\sum_{j_1=p+1}^{\infty}\int\limits_t^{t_2}\psi_1(t_1)\phi_{j_1}(t_1)dt_1
\int\limits_{t_2}^{t_4}\psi_3(t_3)\phi_{j_1}(t_3)dt_3dt_2dt_4\right)^2\le
$$

$$
\le
\sum_{j_2,j_4=0}^{\infty}\left(
\int\limits_t^T\psi_4(t_4)\phi_{j_4}(t_4)
\int\limits_t^{t_4}\psi_2(t_2)\phi_{j_2}(t_2)\times\right.
$$

$$
\left.\times
\sum_{j_1=p+1}^{\infty}\int\limits_t^{t_2}\psi_1(t_1)\phi_{j_1}(t_1)dt_1
\int\limits_{t_2}^{t_4}\psi_3(t_3)\phi_{j_1}(t_3)dt_3dt_2dt_4\right)^2=
$$

$$
=
\int\limits_{[t, T]^2}{\bf 1}_{\{t_2<t_4\}}
\psi_4^2(t_4)
\psi_2^2(t_2)
\left(\sum_{j_1=p+1}^{\infty}\int\limits_t^{t_2}\psi_1(t_1)\phi_{j_1}(t_1)dt_1
\int\limits_{t_2}^{t_4}\psi_3(t_3)\phi_{j_1}(t_3)dt_3\right)^2dt_2dt_4\le 
$$

\begin{equation}
\label{after6009}
\le \frac{K}{p^{2-\varepsilon}}\ \to\  0
\end{equation}

\vspace{4mm}
\noindent
if $p\to\infty,$ where $\varepsilon$ is an arbitrary small positive real number
for the polynomial case and $\varepsilon=0$ for the 
trigonometric case, constant $K$ does not depend on $p.$

The relation (\ref{after6009}) was obtained by the same method as (\ref{after5004}). 
Note that in obtaining (\ref{after6009}) 
we used the estimates (\ref{101oh}), (\ref{after1940}) for the polynomial 
case and (\ref{101xxyy}), (\ref{101xxyy2}) for the trigonometric case.
We also used the integration order replacement in the iterated Riemann integrals
(see (\ref{after98}), (\ref{after99})).

Repeating the previous steps for (\ref{after2503}) and (\ref{after2504}), we get

$$
\sum_{j_2,j_3=0}^p\left(\sum_{j_1=p+1}^{\infty}
C_{j_1 j_3 j_2 j_1}\right)^2=
\sum_{j_2,j_3=0}^p\left(\sum_{j_1=p+1}^{\infty}
\int\limits_t^T\psi_4(t_4)\phi_{j_1}(t_4)
\int\limits_t^{t_4}\psi_3(t_3)\phi_{j_3}(t_3)\times\right.
$$

$$
\left.\times
\int\limits_t^{t_3}\psi_2(t_2)\phi_{j_2}(t_2)
\int\limits_t^{t_2}\psi_1(t_1)\phi_{j_1}(t_1)dt_1dt_2dt_3dt_4\right)^2=
$$

$$
=
\sum_{j_2,j_3=0}^p\left(\sum_{j_1=p+1}^{\infty}
\int\limits_t^T\psi_3(t_3)\phi_{j_3}(t_3)
\int\limits_t^{t_3}\psi_2(t_2)\phi_{j_2}(t_2)\times\right.
$$

$$
\left.\times
\int\limits_t^{t_2}\psi_1(t_1)\phi_{j_1}(t_1)dt_1
\int\limits_{t_3}^{T}\psi_4(t_4)\phi_{j_1}(t_4)dt_4 dt_2 dt_3\right)^2=
$$

$$
=
\sum_{j_2,j_3=0}^p\left(
\int\limits_t^T\psi_3(t_3)\phi_{j_3}(t_3)
\int\limits_t^{t_3}\psi_2(t_2)\phi_{j_2}(t_2)\times\right.
$$

$$
\left.\times
\sum_{j_1=p+1}^{\infty}\int\limits_t^{t_2}\psi_1(t_1)\phi_{j_1}(t_1)dt_1
\int\limits_{t_3}^{T}\psi_4(t_4)\phi_{j_1}(t_4)dt_4 dt_2 dt_3\right)^2\le
$$

$$
\le
\sum_{j_2,j_3=0}^{\infty}\left(
\int\limits_t^T\psi_3(t_3)\phi_{j_3}(t_3)
\int\limits_t^{t_3}\psi_2(t_2)\phi_{j_2}(t_2)\times\right.
$$

$$
\left.\times
\sum_{j_1=p+1}^{\infty}\int\limits_t^{t_2}\psi_1(t_1)\phi_{j_1}(t_1)dt_1
\int\limits_{t_3}^{T}\psi_4(t_4)\phi_{j_1}(t_4)dt_4 dt_2 dt_3\right)^2=
$$

$$
=
\int\limits_{[t, T]^2}{\bf 1}_{\{t_2<t_3\}}
\psi_3^2(t_3)
\psi_2^2(t_2)
\left(\sum_{j_1=p+1}^{\infty}\int\limits_t^{t_2}\psi_1(t_1)\phi_{j_1}(t_1)dt_1
\int\limits_{t_3}^{T}\psi_4(t_4)\phi_{j_1}(t_4)dt_4\right)^2dt_2dt_3\le 
$$
\begin{equation}
\label{after59}
\le \frac{K}{p^2}\ \to\  0
\end{equation}

\vspace{4mm}
\noindent
if $p\to\infty,$ where constant $K$ does not depend on $p;$

$$
\sum_{j_1,j_4=0}^p\left(\sum_{j_2=p+1}^{\infty}
C_{j_4 j_2 j_2 j_1}\right)^2=
\sum_{j_1,j_4=0}^p\left(\sum_{j_2=p+1}^{\infty}
\int\limits_t^T\psi_4(t_4)\phi_{j_4}(t_4)
\int\limits_t^{t_4}\psi_3(t_3)\phi_{j_2}(t_3)\times\right.
$$

$$
\left.\times
\int\limits_t^{t_3}\psi_2(t_2)\phi_{j_2}(t_2)
\int\limits_t^{t_2}\psi_1(t_1)\phi_{j_1}(t_1)dt_1dt_2dt_3dt_4\right)^2=
$$

$$
=
\sum_{j_1,j_4=0}^p\left(\sum_{j_2=p+1}^{\infty}
\int\limits_t^T\psi_4(t_4)\phi_{j_4}(t_4)
\int\limits_t^{t_4}\psi_1(t_1)\phi_{j_1}(t_1)\times\right.
$$

$$
\left.\times
\int\limits_{t_1}^{t_4}\psi_2(t_2)\phi_{j_2}(t_2)
\int\limits_{t_2}^{t_4}\psi_3(t_3)\phi_{j_2}(t_3)dt_3dt_2 dt_1 dt_4\right)^2=
$$

$$
=
\sum_{j_1,j_4=0}^p\left(
\int\limits_t^T\psi_4(t_4)\phi_{j_4}(t_4)
\int\limits_t^{t_4}\psi_1(t_1)\phi_{j_1}(t_1)\times\right.
$$

$$
\left.\times
\sum_{j_2=p+1}^{\infty}\int\limits_{t_1}^{t_4}\psi_2(t_2)\phi_{j_2}(t_2)
\int\limits_{t_2}^{t_4}\psi_3(t_3)\phi_{j_2}(t_3)dt_3dt_2 dt_1 dt_4\right)^2\le
$$

$$
\le
\sum_{j_1,j_4=0}^{\infty}\left(
\int\limits_t^T\psi_4(t_4)\phi_{j_4}(t_4)
\int\limits_t^{t_4}\psi_1(t_1)\phi_{j_1}(t_1)\times\right.
$$

$$
\left.\times
\sum_{j_2=p+1}^{\infty}\int\limits_{t_1}^{t_4}\psi_2(t_2)\phi_{j_2}(t_2)
\int\limits_{t_2}^{t_4}\psi_3(t_3)\phi_{j_2}(t_3)dt_3dt_2 dt_1 dt_4\right)^2=
$$

\begin{equation}
\label{after72}
=
\int\limits_{[t, T]^2}{\bf 1}_{\{t_1<t_4\}}
\psi_4^2(t_4)
\psi_1^2(t_1)
\left(\sum_{j_2=p+1}^{\infty}\int\limits_{t_1}^{t_4}\psi_2(t_2)\phi_{j_2}(t_2)
\int\limits_{t_2}^{t_4}\psi_3(t_3)\phi_{j_2}(t_3)dt_3 dt_2\right)^2 dt_1dt_4.
\end{equation}

\vspace{4mm} 

Note that, by virtue of the additivity property of the integral, we have

\begin{equation}
\label{after8000}
\sum_{j_2=p+1}^{\infty}\int\limits_{t_1}^{t_4}\psi_2(t_2)\phi_{j_2}(t_2)
\int\limits_{t_2}^{t_4}\psi_3(t_3)\phi_{j_2}(t_3)dt_3 dt_2=
\end{equation}

$$
=\sum_{j_2=p+1}^{\infty}
\int\limits_{t}^{t_4}\psi_3(t_3)\phi_{j_2}(t_3)
\int\limits_{t}^{t_3}\psi_2(t_2)\phi_{j_2}(t_2)dt_2 dt_3-
$$

$$
-\sum_{j_2=p+1}^{\infty}\int\limits_{t}^{t_1}\psi_3(t_3)\phi_{j_2}(t_3)
\int\limits_{t}^{t_3}\psi_2(t_2)\phi_{j_2}(t_2)dt_2 dt_3-
$$

\begin{equation}
\label{after71}
-\sum_{j_2=p+1}^{\infty}
\int\limits_{t_1}^{t_4}\psi_3(t_3)\phi_{j_2}(t_3)dt_3
\int\limits_{t}^{t_1}\psi_2(t_2)\phi_{j_2}(t_2)dt_2.
\end{equation}

\vspace{4mm}

However, all three series on the right-hand side of (\ref{after71})
have already been evaluated in (\ref{after5004}) and (\ref{after6009}).
From (\ref{after72}) and (\ref{after71}) we finally obtain

\begin{equation}
\label{after8001}
\sum_{j_1,j_4=0}^p\left(\sum_{j_2=p+1}^{\infty}
C_{j_4 j_2 j_2 j_1}\right)^2\le \frac{K}{p^{2-\varepsilon}}\ \to\  0
\end{equation}

\vspace{3mm}
\noindent
if $p\to\infty,$ where $\varepsilon$ is an arbitrary small positive real number
for the polynomial case and $\varepsilon=0$ for the 
trigonometric case, constant $K$ does not depend on $p.$

In complete analogy with (\ref{after6009}), we have for (\ref{after2505})

$$
\sum_{j_1,j_3=0}^p\left(\sum_{j_2=p+1}^{\infty}
C_{j_2 j_3 j_2 j_1}\right)^2=
\sum_{j_1,j_3=0}^p\left(\sum_{j_2=p+1}^{\infty}
\int\limits_t^T\psi_4(t_4)\phi_{j_2}(t_4)
\int\limits_t^{t_4}\psi_3(t_3)\phi_{j_3}(t_3)\times\right.
$$

$$
\left.\times
\int\limits_t^{t_3}\psi_2(t_2)\phi_{j_2}(t_2)
\int\limits_t^{t_2}\psi_1(t_1)\phi_{j_1}(t_1)dt_1dt_2dt_3dt_4\right)^2=
$$

$$
=
\sum_{j_1,j_3=0}^p\left(\sum_{j_2=p+1}^{\infty}
\int\limits_t^T\psi_3(t_3)\phi_{j_3}(t_3)
\int\limits_t^{t_3}\psi_2(t_2)\phi_{j_2}(t_2)\times\right.
$$

$$
\left.\times
\int\limits_t^{t_2}\psi_1(t_1)\phi_{j_1}(t_1)dt_1 dt_2
\int\limits_{t_3}^{T}\psi_4(t_4)\phi_{j_2}(t_4)dt_4dt_3\right)^2=
$$

$$
=
\sum_{j_1,j_3=0}^p\left(\sum_{j_2=p+1}^{\infty}
\int\limits_t^T\psi_3(t_3)\phi_{j_3}(t_3)
\int\limits_t^{t_3}\psi_1(t_1)\phi_{j_1}(t_1)\times\right.
$$

$$
\left.\times
\int\limits_{t_1}^{t_3}\psi_2(t_2)\phi_{j_2}(t_2)dt_2 dt_1
\int\limits_{t_3}^{T}\psi_4(t_4)\phi_{j_2}(t_4)dt_4dt_3\right)^2=
$$

$$
=
\sum_{j_1,j_3=0}^p\left(
\int\limits_t^T\psi_3(t_3)\phi_{j_3}(t_3)
\int\limits_t^{t_3}\psi_1(t_1)\phi_{j_1}(t_1)\times\right.
$$

$$
\left.\times
\sum_{j_2=p+1}^{\infty}\int\limits_{t_1}^{t_3}\psi_2(t_2)\phi_{j_2}(t_2)dt_2 dt_1
\int\limits_{t_3}^{T}\psi_4(t_4)\phi_{j_2}(t_4)dt_4dt_3\right)^2\le
$$

$$
\le
\sum_{j_1,j_3=0}^{\infty}\left(
\int\limits_t^T\psi_3(t_3)\phi_{j_3}(t_3)
\int\limits_t^{t_3}\psi_1(t_1)\phi_{j_1}(t_1)\times\right.
$$

$$
\left.\times
\sum_{j_2=p+1}^{\infty}\int\limits_{t_1}^{t_3}\psi_2(t_2)\phi_{j_2}(t_2)dt_2 
\int\limits_{t_3}^{T}\psi_4(t_4)\phi_{j_2}(t_4)dt_4 dt_1 dt_3\right)^2=
$$

$$
=
\int\limits_{[t,T]^2}{\bf 1}_{\{t_1<t_3\}}\psi_3^2(t_3)
\psi_1^2(t_1)
\left(\sum_{j_2=p+1}^{\infty}\int\limits_{t_1}^{t_3}\psi_2(t_2)\phi_{j_2}(t_2)dt_2 
\int\limits_{t_3}^{T}\psi_4(t_4)\phi_{j_2}(t_4)dt_4\right)^2
dt_1 dt_3
\le
$$

\begin{equation}
\label{after9030}
\le \frac{K}{p^{2-\varepsilon}}\ \to\  0
\end{equation}

\vspace{4mm}
\noindent
if $p\to\infty,$ where $\varepsilon$ is an arbitrary small positive real number
for the polynomial case and $\varepsilon=0$ for the 
trigonometric case, constant $K$ does not depend on $p.$

We have for (\ref{after2506})

$$
\sum_{j_1,j_2=0}^p\left(\sum_{j_3=p+1}^{\infty}
C_{j_3 j_3 j_2 j_1}\right)^2=
\sum_{j_1,j_2=0}^p\left(\sum_{j_3=p+1}^{\infty}
\int\limits_t^T\psi_4(t_4)\phi_{j_3}(t_4)
\int\limits_t^{t_4}\psi_3(t_3)\phi_{j_3}(t_3)\times\right.
$$

$$
\left.\times
\int\limits_t^{t_3}\psi_2(t_2)\phi_{j_2}(t_2)
\int\limits_t^{t_2}\psi_1(t_1)\phi_{j_1}(t_1)dt_1dt_2dt_3dt_4\right)^2=
$$

$$
=
\sum_{j_1,j_2=0}^p\left(\sum_{j_3=p+1}^{\infty}
\int\limits_t^T\psi_1(t_1)\phi_{j_1}(t_1)
\int\limits_{t_1}^{T}\psi_2(t_2)\phi_{j_2}(t_2)\times\right.
$$

$$
\left.\times
\int\limits_{t_2}^{T}\psi_3(t_3)\phi_{j_3}(t_3)
\int\limits_{t_3}^{T}\psi_4(t_4)\phi_{j_3}(t_4)dt_4dt_3dt_2dt_1\right)^2=
$$

$$
=
\sum_{j_1,j_2=0}^p\left(
\int\limits_t^T\psi_1(t_1)\phi_{j_1}(t_1)
\int\limits_{t_1}^{T}\psi_2(t_2)\phi_{j_2}(t_2)\times\right.
$$

$$
\left.\times
\sum_{j_3=p+1}^{\infty}\int\limits_{t_2}^{T}\psi_3(t_3)\phi_{j_3}(t_3)
\int\limits_{t_3}^{T}\psi_4(t_4)\phi_{j_3}(t_4)dt_4dt_3dt_2dt_1\right)^2\le
$$

$$
\le
\sum_{j_1,j_2=0}^{\infty}\left(
\int\limits_t^T\psi_1(t_1)\phi_{j_1}(t_1)
\int\limits_{t_1}^{T}\psi_2(t_2)\phi_{j_2}(t_2)\times\right.
$$

$$
\left.\times
\sum_{j_3=p+1}^{\infty}\int\limits_{t_2}^{T}\psi_3(t_3)\phi_{j_3}(t_3)
\int\limits_{t_3}^{T}\psi_4(t_4)\phi_{j_3}(t_4)dt_4dt_3dt_2dt_1\right)^2=
$$

\begin{equation}
\label{after8002}
=
\int\limits_{[t,T]^2}{\bf 1}_{\{t_1<t_2\}}\psi_1^2(t_1)
\psi_2^2(t_2)
\left(\sum_{j_3=p+1}^{\infty}\int\limits_{t_2}^{T}\psi_3(t_3)\phi_{j_3}(t_3)
\int\limits_{t_3}^{T}\psi_4(t_4)\phi_{j_3}(t_4)dt_4dt_3\right)^2 dt_2dt_1.
\end{equation}

\vspace{4mm}

It is easy to see that the integral (see (\ref{after8002}))

$$
\int\limits_{t_2}^{T}\psi_3(t_3)\phi_{j_3}(t_3)
\int\limits_{t_3}^{T}\psi_4(t_4)\phi_{j_3}(t_4)dt_4dt_3
$$

\vspace{3mm}
\noindent
is similar to the integral from the formula
(\ref{after8000}) if in the last integral we substitute $t_4=T.$
Therefore, by analogy with (\ref{after8001}), we obtain

\begin{equation}
\label{after9040}
\sum_{j_1,j_2=0}^p\left(\sum_{j_3=p+1}^{\infty}
C_{j_3 j_3 j_2 j_1}\right)^2
\le \frac{K}{p^{2-\varepsilon}}\ \to\  0
\end{equation}

\vspace{3mm}
\noindent
if $p\to\infty,$ where $\varepsilon$ is an arbitrary small positive real number
for the polynomial case and $\varepsilon=0$ for the 
trigonometric case, constant $K$ does not depend on $p.$

Now consider (\ref{after2508})--(\ref{after2507}).
We have for (\ref{after2508}) (see {\bf Step~2} in the proof
of Theorem~13)

$$
\left(\sum_{j_2=p+1}^{\infty}\sum_{j_1=p+1}^{\infty}
C_{j_2 j_1 j_2 j_1}\right)^2=
\left(\sum_{j_1=0}^{p}\sum_{j_2=p+1}^{\infty}
C_{j_2 j_1 j_2 j_1}\right)^2\le 
$$

\begin{equation}
\label{after8003}
\le (p+1)
\sum_{j_1=0}^{p}\left(\sum_{j_2=p+1}^{\infty}
C_{j_2 j_1 j_2 j_1}\right)^2.
\end{equation}

\vspace{3mm}

Consider (\ref{after2505}) and (\ref{after9030}). We have

$$
\sum_{j_1=0}^p\left(\sum_{j_2=p+1}^{\infty}
C_{j_2 j_1 j_2 j_1}\right)^2=
\sum_{j_1,j_3=0}^p\left(\sum_{j_2=p+1}^{\infty}
C_{j_2 j_3 j_2 j_1}\right)^2\Biggl|_{j_1=j_3}\Biggr.\le
$$

\begin{equation}
\label{after8009}
\le\sum_{j_1,j_3=0}^p\left(\sum_{j_2=p+1}^{\infty}
C_{j_2 j_3 j_2 j_1}\right)^2
 \le \frac{K}{p^{2-\varepsilon}},
\end{equation}

\vspace{3mm}
\noindent
where $\varepsilon$ is an arbitrary small positive real number
for the polynomial case and $\varepsilon=0$ for the 
trigonometric case, constant $K$ does not depend on $p.$
Combining (\ref{after8003}) and (\ref{after8009}), we obtain

$$
\left(\sum_{j_2=p+1}^{\infty}\sum_{j_1=p+1}^{\infty}
C_{j_2 j_1 j_2 j_1}\right)^2\le \frac{(p+1)K}{p^{2-\varepsilon}}\le
\frac{K_1}{p^{1-\varepsilon}}\ \to\  0
$$

\vspace{3mm}
\noindent
if $p\to\infty,$ where constant $K_1$ does not depend on $p.$

Similarly for (\ref{after2509}) we have (see (\ref{after2504}), (\ref{after8001}))

$$
\left(\sum_{j_2=p+1}^{\infty}\sum_{j_1=p+1}^{\infty}
C_{j_1 j_2 j_2 j_1}\right)^2=
\left(\sum_{j_1=0}^{p}\sum_{j_2=p+1}^{\infty}
C_{j_1 j_2 j_2 j_1}\right)^2\le 
$$

\begin{equation}
\label{after9002}
\le (p+1)\sum_{j_1=0}^{p}
\left(\sum_{j_2=p+1}^{\infty}
C_{j_1 j_2 j_2 j_1}\right)^2,
\end{equation}

\vspace{2mm}
$$
\sum_{j_1=0}^p\left(\sum_{j_2=p+1}^{\infty}
C_{j_1 j_2 j_2 j_1}\right)^2=
\sum_{j_1,j_4=0}^p\left(\sum_{j_2=p+1}^{\infty}
C_{j_4 j_2 j_2 j_1}\right)^2\Biggl|_{j_1=j_4}\Biggr.\le
$$

\begin{equation}
\label{after9011}
\le
\sum_{j_1,j_4=0}^p\left(\sum_{j_2=p+1}^{\infty}
C_{j_4 j_2 j_2 j_1}\right)^2\le \frac{K}{p^{2-\varepsilon}},
\end{equation}

\vspace{4mm}
\noindent
where $\varepsilon$ is an arbitrary small positive real number
for the polynomial case and $\varepsilon=0$ for the 
trigonometric case, constant $K$ does not depend on $p.$
Combining (\ref{after9002}) and (\ref{after9011}), we obtain

$$
\left(\sum_{j_2=p+1}^{\infty}\sum_{j_1=p+1}^{\infty}
C_{j_1 j_2 j_2 j_1}\right)^2\le \frac{(p+1)K}{p^{2-\varepsilon}}\le
\frac{K_1}{p^{1-\varepsilon}}\ \to\  0
$$

\vspace{3mm}
\noindent
if $p\to\infty,$ where constant $K_1$ does not depend on $p.$

Consider (\ref{after2507}). Using (\ref{after500}), we get

$$
\sum_{j_3=p+1}^{\infty}\sum_{j_1=p+1}^{\infty}
C_{j_3 j_3 j_1 j_1}=\sum_{j_3=p+1}^{\infty}\sum_{j_1=0}^{\infty}
C_{j_3 j_3 j_1 j_1}-\sum_{j_3=p+1}^{\infty}\sum_{j_1=0}^{p}
C_{j_3 j_3 j_1 j_1}=
$$

\begin{equation}
\label{after9041}
=\frac{1}{2}\sum_{j_3=p+1}^{\infty}
C_{j_3 j_3 j_1 j_1}\biggl|_{(j_1 j_1)\curvearrowright (\cdot) }\biggr.
-\sum_{j_3=p+1}^{\infty}\sum_{j_1=0}^{p}
C_{j_3 j_3 j_1 j_1},
\end{equation}

\vspace{3mm}
\noindent
where (see (\ref{after900}))
$$
C_{j_3 j_3 j_1 j_1}\biggl|_{(j_1 j_1)\curvearrowright (\cdot) }\biggr.=
$$

$$
=
\int\limits_t^T\psi_4(t_4)\phi_{j_3}(t_4)
\int\limits_t^{t_4}\psi_3(t_3)\phi_{j_3}(t_3)
\int\limits_t^{t_3}\psi_2(t_2)\psi_1(t_2)dt_2 dt_3 dt_4.
$$

\vspace{3mm}

From the estimate (\ref{5tafter}) for the polynomial and
trigonometric cases we get

\vspace{-1mm}
\begin{equation}
\label{after9042}
\left|\sum_{j_3=p+1}^{\infty}
C_{j_3 j_3 j_1 j_1}\biggl|_{(j_1 j_1)\curvearrowright (\cdot) }\biggr.\right|\le
\frac{C}{p},
\end{equation}

\vspace{3mm}
\noindent
where constant $C$ is independent of $p.$ 

Further, we have (see (\ref{after9040}))

$$
\left(\sum_{j_1=0}^{p}\sum_{j_3=p+1}^{\infty}
C_{j_3 j_3 j_1 j_1}\right)^2\le (p+1)\sum_{j_1=0}^{p}
\left(\sum_{j_3=p+1}^{\infty}C_{j_3 j_3 j_1 j_1}\right)^2=
$$

$$
=(p+1)\sum_{j_1,j_2=0}^{p}
\left(\sum_{j_3=p+1}^{\infty}C_{j_3 j_3 j_2 j_1}\right)^2\Biggl|_{j_1=j_2}\Biggr.\le
$$

\begin{equation}
\label{after9043}
\le
(p+1)\sum_{j_1,j_2=0}^{p}
\left(\sum_{j_3=p+1}^{\infty}C_{j_3 j_3 j_2 j_1}\right)^2\le
\frac{(p+1)K}{p^{2-\varepsilon}}\le
\frac{K_1}{p^{1-\varepsilon}},
\end{equation}

\vspace{4mm}
\noindent
where constant $K_1$ does not depend on $p.$

Combining (\ref{after9041})--(\ref{after9043}), we obtain

$$
\left(\sum_{j_3=p+1}^{\infty}\sum_{j_1=p+1}^{\infty}
C_{j_3 j_3 j_1 j_1}\right)^2
\le \frac{K_2}{p^{1-\varepsilon}}\ \to\  0
$$

\vspace{4mm}
\noindent
if $p\to\infty,$ where constant $K_2$ does not depend on $p.$

Let us prove (\ref{after1602xxxx})--(\ref{after1602xxxz}).
It is not difficult to see that the estimate (\ref{after9042})
proves (\ref{after1602xxxx}).

Using the integration order replacement, we have

\vspace{-1mm}
$$
\sum_{j_1=p+1}^{\infty}
C_{j_3 j_3 j_1 j_1}\biggl|_{(j_{3} j_{3})\curvearrowright (\cdot)}=
$$

\vspace{2mm}
$$
=\sum_{j_1=p+1}^{\infty}\int\limits_t^T\psi_4(t_4)\psi_3(t_4)
\int\limits_t^{t_4}\psi_2(t_2)\phi_{j_1}(t_2)
\int\limits_t^{t_2}\psi_1(t_1)\phi_{j_1}(t_1)dt_1dt_2dt_4=
$$

\vspace{2mm}
\begin{equation}
\label{afterafter1}
=
\sum_{j_1=p+1}^{\infty}
\int\limits_t^{T}\left(\psi_2(t_2)
\int\limits_{t_2}^T\psi_4(t_4)\psi_3(t_4)dt_4\right)\phi_{j_1}(t_2)
\int\limits_t^{t_2}\psi_1(t_1)\phi_{j_1}(t_1)dt_1dt_2,
\end{equation}

\vspace{4mm}
$$
\sum_{j_1=p+1}^{\infty}
C_{j_1 j_2 j_2 j_1}\biggl|_{(j_{2} j_{2})\curvearrowright (\cdot)}\biggr.=
$$

\vspace{2mm}
$$
=
\sum_{j_1=p+1}^{\infty}\int\limits_t^T\psi_4(t_4)\phi_{j_1}(t_4)
\int\limits_t^{t_4}\psi_3(t_3)\psi_2(t_3)
\int\limits_t^{t_3}\psi_1(t_1)\phi_{j_1}(t_1)dt_1dt_3dt_4=
$$

\vspace{2mm}
$$
=
\sum_{j_1=p+1}^{\infty}\int\limits_t^T\psi_4(t_4)\phi_{j_1}(t_4)
\int\limits_t^{t_4}\psi_1(t_1)\phi_{j_1}(t_1)\int\limits_{t_1}^{t_4}         
\psi_3(t_3)\psi_2(t_3)dt_3dt_1dt_4=
$$

\vspace{2mm}
$$
=
\sum_{j_1=p+1}^{\infty}\int\limits_t^T\psi_4(t_4)\phi_{j_1}(t_4)
\int\limits_t^{t_4}\psi_1(t_1)\phi_{j_1}(t_1)
\left(\int\limits_{t}^{t_4}-\int\limits_{t}^{t_1}\right)
\psi_3(t_3)\psi_2(t_3)dt_3dt_1dt_4=
$$

\vspace{2mm}
\begin{equation}
\label{afterafter198}
=
\sum_{j_1=p+1}^{\infty}\int\limits_t^T\left(\psi_4(t_4)
\int\limits_{t}^{t_4}
\psi_3(t_3)\psi_2(t_3)dt_3\right)
\phi_{j_1}(t_4)
\int\limits_t^{t_4}\psi_1(t_1)\phi_{j_1}(t_1)
dt_1dt_4-
\end{equation}

\vspace{2mm}
\begin{equation}
\label{afterafter199}
-
\sum_{j_1=p+1}^{\infty}\int\limits_t^T\psi_4(t_4)\phi_{j_1}(t_4)
\int\limits_t^{t_4}\left(\psi_1(t_1)
\int\limits_{t}^{t_1}\
\psi_3(t_3)\psi_2(t_3)dt_3\right)\phi_{j_1}(t_1)dt_1dt_4.
\end{equation}

\vspace{5mm}

Applying the estimate (\ref{5tafter}) (polynomial 
and trigonometric cases)   
to the right-hand sides of (\ref{afterafter1})--(\ref{afterafter199}), we get

\vspace{-2mm}
\begin{equation}
\label{after9042xxx}
\left|\sum_{j_3=p+1}^{\infty}
C_{j_3 j_3 j_1 j_1}\biggl|_{(j_3 j_3)\curvearrowright (\cdot) }\biggr.\right|\le
\frac{C}{p},
\end{equation}

\vspace{1mm}
\begin{equation}
\label{after9042xxxe}
\left|\sum_{j_1=p+1}^{\infty}
C_{j_1 j_2 j_2 j_1}\biggl|_{(j_{2} j_{2})\curvearrowright (\cdot)}\biggr.\right|\le
\frac{C}{p},
\end{equation}

\vspace{4mm}
\noindent
where constant $C$ is independent of $p.$
The estimates (\ref{after9042xxx}), (\ref{after9042xxxe})
prove (\ref{after1602xxxy}), (\ref{after1602xxxz}).

The relations (\ref{after2501})--(\ref{after1602xxxz}) are proved. Theorem~17 is proved.

\vspace{5mm}

\section{Expansion of Iterated Stratonovich Stochastic Integrals
of Multiplicity 5. The Case $p_1=\ldots =p_5\to \infty$ and 
Continuously Differentiable 
Weight Functions $\psi_1(\tau),$ $\ldots,$ $\psi_5(\tau)$ 
(The Cases of Legendre 
Polynomials and Trigonometric Functions)}

\vspace{5mm}

{\bf Theorem 18}\ \cite{2018a}, \cite{arxiv-5}, \cite{arxiv-9}, \cite{new-art-1-xxy}.\ 
{\it Suppose that 
$\{\phi_j(x)\}_{j=0}^{\infty}$ is a complete orthonormal system of 
Legendre polynomials or trigonometric func\-ti\-ons in the space $L_2([t, T]).$
Furthermore, let $\psi_1(\tau), \ldots,$ $\psi_5(\tau)$ are continuously dif\-ferentiable 
nonrandom functions on $[t, T].$ 
Then, for the 
iterated Stra\-to\-no\-vich stochastic integral of fifth multiplicity

\vspace{-1mm}
\begin{equation}
\label{after10001}
J^{*}[\psi^{(5)}]_{T,t}={\int\limits_t^{*}}^T\psi_5(t_5)
\ldots
{\int\limits_t^{*}}^{t_2}\psi_1(t_1)
d{\bf w}_{t_1}^{(i_1)}
\ldots d{\bf w}_{t_5}^{(i_5)}
\end{equation}

\vspace{3mm}
\noindent
the following 
expansion 

\vspace{-2mm}
$$
J^{*}[\psi^{(5)}]_{T,t}
=\hbox{\vtop{\offinterlineskip\halign{
\hfil#\hfil\cr
{\rm l.i.m.}\cr
$\stackrel{}{{}_{p\to \infty}}$\cr
}} }
\sum\limits_{j_1, \ldots, j_5=0}^{p}
C_{j_5 \ldots j_1}\zeta_{j_1}^{(i_1)}\ldots
\zeta_{j_5}^{(i_5)}
$$

\vspace{5mm}
\noindent
that converges in the mean-square sense is valid, where
$i_1, \ldots, i_5=0, 1,\ldots,m,$

\vspace{-1mm}
$$
C_{j_5 \ldots j_1}=
\int\limits_t^T\psi_5(t_5)\phi_{j_5}(t_5)\ldots
\int\limits_t^{t_2}\psi_1(t_1)\phi_{j_1}(t_1)dt_1\ldots dt_5
$$

\vspace{3mm}
\noindent
and
$$
\zeta_{j}^{(i)}=
\int\limits_t^T \phi_{j}(s) d{\bf w}_s^{(i)}
$$ 

\vspace{2mm}
\noindent
are independent standard Gaussian random variables for various 
$i$ or $j$
{\rm (}in the case when $i\ne 0${\rm ),}
${\bf w}_{\tau}^{(i)}={\bf f}_{\tau}^{(i)}$ for
$i=1,\ldots,m$ and 
${\bf w}_{\tau}^{(0)}=\tau.$}

\vspace{2mm}

{\bf Proof.}\ Note that in this proof we write $k$ instead of 5 when this is true for 
an arbitrary $k$ $(k\in \mathbb{N}).$
As follows from Sect.~7, Conditions 1 and 2
of Theorem~{\rm 13} are satisfied for complete
orthonormal systems of Legendre polynomials 
and trigonometric functions in the space
$L_2([t, T]).$ Let us verify Condition 3 of Theorem~13 for
the iterated Stratonovich stochastic integral (\ref{after10001}). 
Thus, we have to check the following conditions

\vspace{-1mm}
\begin{equation}
\label{after14000}
\lim\limits_{p\to\infty}
\sum\limits_{j_{q_1},j_{q_2},j_{q_3}=0}^p
\left(\sum_{j_{g_1}=p+1}^{\infty}
C_{j_5\ldots j_1}\biggl|_{j_{g_1}=j_{g_2}}\biggr.\right)^2=0,
\end{equation}

\vspace{1mm}
\begin{equation}
\label{after14001}
\lim\limits_{p\to\infty}
\sum\limits_{j_{q_1}=0}^p
\left(\sum_{j_{g_1}=p+1}^{\infty}\sum_{j_{g_3}=p+1}^{\infty}
C_{j_5\ldots j_1}\biggl|_{j_{g_1}=j_{g_2},j_{g_3}=j_{g_4}}\biggr.\right)^2=0,
\end{equation}

\vspace{1mm}
\begin{equation}
\label{afterafter001}
\lim\limits_{p\to\infty}
\sum\limits_{j_{q_1}=0}^p
\left(\sum_{j_{g_3}=p+1}^{\infty}
C_{j_5\ldots j_1}\biggl|_{(j_{g_2}j_{g_1})\curvearrowright (\cdot),
j_{g_1}=j_{g_2},
j_{g_3}=j_{g_4}, g_2=g_1+1}\biggr.\right)^2=0,
\end{equation}

\vspace{5mm}
\noindent
where 
$\left(\{g_1,g_2\},\{g_3,g_4\}, \{q_1\}\right)$ and 
$\left(\{g_1,g_2\}, \{q_1, q_2,q_3\}\right)$ 
are partitions of the set $\{1,2,\ldots,5\}$ that is
$\{g_1,g_2,g_3,g_4,q_1\}=\{g_1,g_2,q_1,q_2,q_3\}=\{1,2,\ldots,5\};$
braces mean an unordered 
set, and pa\-ren\-the\-ses mean an ordered set.

Let us find a representation for
$C_{j_k\ldots j_1}\bigl|_{j_{g_1}=j_{g_2},\ g_2>g_1+1}\biggr.$ 
that will be convenient for further con\-si\-de\-ra\-ti\-on.

Using the integration order replacement in the Riemann integrals, we obtain

\vspace{-1mm}
$$
\int\limits_t^T h_{k}(t_k)\ldots \int\limits_t^{t_{l+2}} h_{l+1}(t_{l+1})
\int\limits_t^{t_{l+1}} h_{l}(t_{l})
\int\limits_t^{t_{l}} h_{l-1}(t_{l-1})\ldots
\int\limits_t^{t_2} h_{1}(t_1)
dt_1\ldots 
$$

$$
\ldots
dt_{l-1}dt_{l}dt_{l+1}\ldots dt_k=
$$

\vspace{1mm}
$$
=\int\limits_t^T h_{k}(t_k)\ldots \int\limits_t^{t_{l+2}} h_{l+1}(t_{l+1})
\int\limits_t^{t_{l+1}} h_{1}(t_{1})
\int\limits_{t_1}^{t_{l+1}} h_{2}(t_{2})\ldots
\int\limits_{t_{l-2}}^{t_{l+1}} h_{l-1}(t_{l-1})
\int\limits_{t_{l-1}}^{t_{l+1}} h_{l}(t_{l})dt_l\times
$$

$$
\times dt_{l-1}\ldots dt_2dt_{1}dt_{l+1}\ldots dt_k=
$$

\vspace{1mm}
$$
=\int\limits_t^T h_{k}(t_k)\ldots \int\limits_t^{t_{l+2}} h_{l+1}(t_{l+1})
\left(\int\limits_{t}^{t_{l+1}} h_{l}(t_{l})dt_l\right)\int\limits_t^{t_{l+1}} h_{1}(t_{1})
\int\limits_{t_1}^{t_{l+1}} h_{2}(t_{2})\ldots
\int\limits_{t_{l-2}}^{t_{l+1}} h_{l-1}(t_{l-1})
\times
$$

$$
\times dt_{l-1}\ldots dt_2dt_{1}dt_{l+1}\ldots dt_k-
$$

\vspace{1mm}
$$
-\int\limits_t^T h_{k}(t_k)\ldots \int\limits_t^{t_{l+2}} h_{l+1}(t_{l+1})
\int\limits_t^{t_{l+1}} h_{1}(t_{1})
\int\limits_{t_1}^{t_{l+1}} h_{2}(t_{2})\ldots
\int\limits_{t_{l-2}}^{t_{l+1}} h_{l-1}(t_{l-1})
\left(\int\limits_{t}^{t_{l-1}} h_{l}(t_{l})dt_l\right)\times
$$

$$
\times dt_{l-1}\ldots dt_2dt_{1}dt_{l+1}\ldots dt_k=
$$

\vspace{1mm}
$$
=\int\limits_t^T h_{k}(t_k)\ldots \int\limits_t^{t_{l+2}} h_{l+1}(t_{l+1})
\left(\int\limits_t^{t_{l+1}} h_{l}(t_{l})dt_l\right)
\int\limits_t^{t_{l+1}} h_{l-1}(t_{l-1})\ldots
$$

$$
\ldots 
\int\limits_t^{t_2} h_{1}(t_1)
dt_1\ldots dt_{l-1}dt_{l+1}\ldots dt_k-
$$

\vspace{1mm}
$$
-\int\limits_t^T h_{k}(t_k)\ldots \int\limits_t^{t_{l+2}} h_{l+1}(t_{l+1})
\int\limits_t^{t_{l+1}} h_{l-1}(t_{l-1})\left(\int\limits_t^{t_{l-1}} h_{l}(t_{l})dt_l\right)
\int\limits_t^{t_{l-1}} h_{l-2}(t_{l-2})
\ldots
$$

\begin{equation}
\label{after81}
\ldots
\int\limits_t^{t_2} h_{1}(t_1)
dt_1\ldots dt_{l-2}dt_{l-1}dt_{l+1}\ldots dt_k,
\end{equation}

\vspace{4mm}
\noindent
where $2<l<k-1$ and $h_1(\tau),\ldots,h_k(\tau)$ are continuous functions on the interval
$[t, T].$ The case $k=1$ is obvious. By analogy with (\ref{after81}) we have for $l=k$ 

\vspace{-1mm}
$$
\int\limits_t^{T} h_{l}(t_{l})
\int\limits_t^{t_{l}} h_{l-1}(t_{l-1})\ldots
\int\limits_t^{t_2} h_{1}(t_1)
dt_1\ldots 
dt_{l-1}dt_{l}=
$$

\vspace{2mm}
$$
=\int\limits_{t}^{T} h_{1}(t_{1})
\int\limits_{t_1}^{T} h_{2}(t_{2})\ldots
\int\limits_{t_{l-2}}^{T} h_{l-1}(t_{l-1})\int\limits_{t_{l-1}}^{T} h_{l}(t_{l})
dt_ldt_{l-1}\ldots dt_2dt_{1}=
$$

\vspace{2mm}
$$
=\left(\int\limits_{t}^{T} h_{l}(t_{l})
dt_l\right)\int\limits_{t}^{T} h_{1}(t_{1})
\int\limits_{t_1}^{T} h_{2}(t_{2})\ldots
\int\limits_{t_{l-2}}^{T} h_{l-1}(t_{l-1})
dt_{l-1}\ldots dt_2dt_{1}-
$$

\vspace{2mm}
$$
-\int\limits_{t}^{T}h_{1}(t_{1})
\int\limits_{t_1}^{T} h_{2}(t_{2})\ldots
\int\limits_{t_{l-2}}^{T} h_{l-1}(t_{l-1})\left(
\int\limits_{t}^{t_{l-1}} h_{l}(t_{l})dt_l\right)
dt_{l-1}\ldots dt_2dt_{1}=
$$

\vspace{2mm}
$$
=\left(\int\limits_{t}^{T} h_{l}(t_{l})
dt_l\right)\int\limits_{t}^{T} h_{l-1}(t_{l-1})
\ldots
\int\limits_{t}^{t_2} h_{1}(t_{1})
dt_{1}\ldots dt_{l-1}-
$$

\vspace{2mm}
\begin{equation}
\label{after82}
-\int\limits_{t}^{T}
h_{l-1}(t_{l-1})\left(
\int\limits_{t}^{t_{l-1}} h_{l}(t_{l})dt_l\right)
\int\limits_{t}^{t_{l-1}} h_{l-2}(t_{l-2})\ldots
\int\limits_{t}^{t_2} 
h_{1}(t_{1})
dt_{1}\ldots dt_{l-1}.
\end{equation}

\vspace{4mm}

The formulas (\ref{after81}), (\ref{after82})
will be used further.

Our further proof will not fundamentally depend on the weight
functions $\psi_1(\tau),\ldots,\psi_k(\tau).$
Therefore, sometimes in subsequent consideration we assume for simplicity
that $\psi_1(\tau),\ldots,\psi_k(\tau)\equiv 1.$

Let us continue the proof. Applying (\ref{after81}) 
to $C_{j_k \ldots j_{l+1} j_l j_{l-1} \ldots j_{s+1} j_l j_{s-1} \ldots j_1}$ 
(more precisely to $h_s(t_s)$ $=$ $\psi_s(t_s)\phi_{j_{l}}(t_{s})$), we obtain
for $l+1\le k,$ $s-1\ge 1,$ $l-1\ge s+1$

\begin{equation}
\label{after90a}
\sum_{j_l=p+1}^{\infty}
C_{j_k \ldots j_{l+1} j_l j_{l-1} \ldots j_{s+1} j_l j_{s-1} \ldots j_1}=
\end{equation}

\vspace{2mm}
$$
=
\sum_{j_l=p+1}^{\infty}
\int\limits_t^T \phi_{j_k}(t_k)\ldots 
\int\limits_t^{t_{l+2}} \phi_{j_{l+1}}(t_{l+1})
\int\limits_t^{t_{l+1}} \phi_{j_{l}}(t_{l})
\int\limits_t^{t_{l}}\phi_{j_{l-1}}(t_{l-1})\ldots
$$

\vspace{1mm}
$$
\ldots \int\limits_t^{t_{s+2}} \phi_{j_{s+1}}(t_{s+1})
\int\limits_t^{t_{s+1}} \phi_{j_{l}}(t_{s})
\int\limits_t^{t_{s}} \phi_{j_{s-1}}(t_{s-1})\ldots 
$$

\vspace{1mm}
$$
\ldots\int\limits_t^{t_{2}} \phi_{j_{1}}(t_{1})dt_1
\ldots dt_{s-1}dt_{s}dt_{s+1}\ldots
dt_{l-1}dt_{l}dt_{l+1}\ldots dt_k=
$$

\vspace{2mm}
$$
=
\sum_{j_l=p+1}^{\infty}
\int\limits_t^T \phi_{j_k}(t_k)\ldots \int\limits_t^{t_{l+2}} \phi_{j_{l+1}}(t_{l+1})
\int\limits_t^{t_{l+1}} \phi_{j_{l}}(t_{l})
\int\limits_t^{t_{l}} \phi_{j_{l-1}}(t_{l-1})\ldots
$$

\vspace{1mm}
$$
\ldots
\int\limits_t^{t_{s+2}} \phi_{j_{s+1}}(t_{s+1})
\left(\int\limits_t^{t_{s+1}} \phi_{j_{l}}(t_{s})dt_s\right)
\int\limits_t^{t_{s+1}} \phi_{j_{s-1}}(t_{s-1})\ldots
$$

\vspace{1mm}
$$
\ldots \int\limits_t^{t_{2}} \phi_{j_{1}}(t_{1})dt_1\ldots dt_{s-1}dt_{s+1}\ldots
dt_{l-1}dt_{l}dt_{l+1}\ldots dt_k-
$$

\vspace{2mm}
$$
-\sum_{j_l=p+1}^{\infty}
\int\limits_t^T \phi_{j_k}(t_k)\ldots \int\limits_t^{t_{l+2}} \phi_{j_{l+1}}(t_{l+1})
\int\limits_t^{t_{l+1}} \phi_{j_{l}}(t_{l})
\int\limits_t^{t_{l}} \phi_{j_{l-1}}(t_{l-1})\ldots
$$

\vspace{1mm}
$$
\ldots
\int\limits_t^{t_{s+2}} \phi_{j_{s+1}}(t_{s+1})
\int\limits_t^{t_{s+1}} \phi_{j_{s-1}}(t_{s-1})
\left(\int\limits_t^{t_{s-1}} \phi_{j_{l}}(t_{s})dt_s\right)
\int\limits_t^{t_{s-1}}\phi_{j_{s-2}}(t_{s-2})
\ldots 
$$

\vspace{1mm}
$$
\ldots \int\limits_t^{t_{2}} \phi_{j_{1}}(t_{1})dt_1\ldots dt_{s-2}dt_{s-1}dt_{s+1}\ldots
dt_{l-1}dt_{l}dt_{l+1}\ldots dt_k =
$$

\vspace{2mm}
$$
=\sum_{j_l=p+1}^{\infty} A_{j_k \ldots j_{l+1} j_l j_{l-1} \ldots j_{s+1} j_l j_{s-1} \ldots j_1}-
\sum_{j_l=p+1}^{\infty} B_{j_k \ldots j_{l+1} j_l j_{l-1} \ldots j_{s+1} j_l j_{s-1} \ldots j_1}.
$$

\vspace{5mm}

Now we apply the formula (\ref{after81}) 
to $A_{j_k \ldots j_{l+1} j_l j_{l-1} \ldots j_{s+1} j_l j_{s-1} \ldots j_1},$
$B_{j_k \ldots j_{l+1} j_l j_{l-1} \ldots j_{s+1} j_l j_{s-1} \ldots j_1}$
(more precisely to $h_l(t_l)=\psi_l(t_l)\phi_{j_{l}}(t_{l})$). Then we have
for $l+1\le k,$ $s-1\ge 1,$ $l-1\ge s+1$

$$
\sum_{j_l=p+1}^{\infty}
C_{j_k \ldots j_{l+1} j_l j_{l-1} \ldots j_{s+1} j_l j_{s-1} \ldots j_1}=
$$

\vspace{2mm}
$$
=\int\limits_{[t, T]^{k-2}} \sum\limits_{d=1}^4
F_p^{(d)}(t_1,\ldots,t_{s-1},t_{s+1},\ldots ,t_{l-1},t_{l+1},\ldots,t_k)\times
$$

$$
\times 
\prod_{\stackrel{g=1}{{}_{g\ne l, s}}}^{k}
\psi_g(t_g)\phi_{j_g}(t_g)
dt_1\ldots dt_{s-1}dt_{s+1}\ldots dt_{l-1}dt_{l+1}\ldots dt_k=
$$

\begin{equation}
\label{after85a}
=\sum\limits_{d=1}^4
C_{j_k \ldots j_{l+1} j_{l-1} \ldots j_{s+1} j_{s-1} \ldots j_1}^{*(d)}
=\sum\limits_{d=1}^4
C_{j_k \ldots j_q \ldots j_1}^{*(d)}\biggl|_{q\ne l, s}\biggr.,
\end{equation}

\vspace{7mm}
\noindent
where
$$
F_p^{(1)}(t_1,\ldots,t_{s-1},t_{s+1},\ldots ,t_{l-1},t_{l+1},\ldots,t_k)=
$$

\vspace{-2mm}
\begin{equation}
\label{after90}
=
{\bf 1}_{\{t_1< \ldots <t_{s-1}<t_{s+1}< \ldots <t_{l-1}<t_{l+1}< \ldots <t_k\}}
\sum_{j_l=p+1}^{\infty}~
\int\limits_{t}^{t_{s+1}} \psi_s(\tau) \phi_{j_{l}}(\tau)d\tau
\int\limits_{t}^{t_{l+1}} \psi_l(\tau) \phi_{j_{l}}(\tau)d\tau,
\end{equation}

\vspace{3mm}
$$
F_p^{(2)}(t_1,\ldots,t_{s-1},t_{s+1},\ldots ,t_{l-1},t_{l+1},\ldots,t_k)=
$$

\vspace{-2mm}
\begin{equation}
\label{after91}
=
{\bf 1}_{\{t_1< \ldots <t_{s-1}<t_{s+1}< \ldots <t_{l-1}<t_{l+1}< \ldots <t_k\}}
\sum_{j_l=p+1}^{\infty}~
\int\limits_{t}^{t_{s-1}} \psi_s(\tau) \phi_{j_{l}}(\tau)d\tau
\int\limits_{t}^{t_{l-1}} \psi_l(\tau) \phi_{j_{l}}(\tau)d\tau,
\end{equation}

\vspace{3mm}
$$
F_p^{(3)}(t_1,\ldots,t_{s-1},t_{s+1},\ldots ,t_{l-1},t_{l+1},\ldots,t_k)=
$$

\vspace{-2mm}
\begin{equation}
\label{after92}
=
-{\bf 1}_{\{t_1< \ldots <t_{s-1}<t_{s+1}< \ldots <t_{l-1}<t_{l+1}< \ldots <t_k\}}
\sum_{j_l=p+1}^{\infty}~
\int\limits_{t}^{t_{s-1}} \psi_s(\tau)\phi_{j_{l}}(\tau)d\tau
\int\limits_{t}^{t_{l+1}} \psi_l(\tau)\phi_{j_{l}}(\tau)d\tau,
\end{equation}

\vspace{3mm}
$$
F_p^{(4)}(t_1,\ldots,t_{s-1},t_{s+1},\ldots ,t_{l-1},t_{l+1},\ldots,t_k)=
$$

\vspace{-2mm}
\begin{equation}
\label{after93}
=
-{\bf 1}_{\{t_1< \ldots <t_{s-1}<t_{s+1}< \ldots <t_{l-1}<t_{l+1}< \ldots <t_k\}}
\sum_{j_l=p+1}^{\infty}~
\int\limits_{t}^{t_{s+1}} \psi_s(\tau)\phi_{j_{l}}(\tau)d\tau
\int\limits_{t}^{t_{l-1}} \psi_l(\tau)\phi_{j_{l}}(\tau)d\tau.
\end{equation}

\vspace{4mm}

By analogy with (\ref{after85a}) we can consider the expressions

\vspace{-1mm}
\begin{equation}
\label{afterx200}
\sum_{j_l=p+1}^{\infty}
C_{j_l j_{k-1}\ldots j_2 j_l},
\end{equation}

\begin{equation}
\label{afterx201}
\sum_{j_l=p+1}^{\infty}
C_{j_k \ldots j_{l+1} j_l j_{l-1} \ldots j_2 j_l}\ \ \ (l+1\le k),
\end{equation}

\begin{equation}
\label{afterx202}
\sum_{j_l=p+1}^{\infty}
C_{j_l j_{k-1} \ldots j_{s+1} j_l j_{s-1} \ldots j_1}\ \ \ (s-1\ge 1).
\end{equation}

\vspace{4mm}

Then we have for (\ref{afterx200})--(\ref{afterx202}) (see (\ref{after81}), (\ref{after82}))

\begin{equation}
\label{afterx204}
\sum_{j_l=p+1}^{\infty}
C_{j_l j_{k-1} \ldots j_2 j_l}=
\int\limits_{[t, T]^{k-2}} \sum\limits_{d=1}^2
G_p^{(d)}(t_2,\ldots,t_{k-1})
\prod_{g=2}^{k-1}
\psi_g(t_g)\phi_{j_g}(t_g)
dt_2\ldots dt_{k-1},
\end{equation}

\vspace{4mm}
$$
\sum_{j_l=p+1}^{\infty}
C_{j_k \ldots j_{l+1} j_l j_{l-1} \ldots j_2 j_l}=
\int\limits_{[t, T]^{k-2}} \sum\limits_{d=1}^2
E_p^{(d)}(t_2,\ldots,t_{l-1},t_{l+1}, \ldots, t_{k})\times
$$

\vspace{-2mm}
\begin{equation}
\label{afterx205}
\times \prod_{\stackrel{g=2}{{}_{g\ne l}}}^{k}
\psi_g(t_g)\phi_{j_g}(t_g)
dt_2\ldots dt_{l-1}dt_{l+1} \ldots dt_{k},
\end{equation}

\vspace{3mm}
$$
\sum_{j_l=p+1}^{\infty}
C_{j_l j_{k-1} \ldots j_{s+1} j_l j_{s-1} \ldots j_1}=
\int\limits_{[t, T]^{k-2}} \sum\limits_{d=1}^4
D_p^{(d)}(t_1,\ldots,t_{s-1},t_{s+1}, \ldots, t_{k-1})\times
$$

\vspace{-2mm}
\begin{equation}
\label{afterx206}
\times \prod_{\stackrel{g=1}{{}_{g\ne s}}}^{k-1}
\psi_g(t_g)\phi_{j_g}(t_g)
dt_1\ldots dt_{s-1}dt_{s+1} \ldots dt_{k-1},
\end{equation}

\vspace{4mm}
\noindent
where
$$
G_p^{(1)}(t_2,\ldots,t_{k-1})=
{\bf 1}_{\{t_2< \ldots <t_{k-1}\}}
\sum_{j_l=p+1}^{\infty}~
\int\limits_{t}^{T} \psi_k(\tau)\phi_{j_{l}}(\tau)d\tau
\int\limits_{t}^{t_{2}} \psi_1(\tau)\phi_{j_{l}}(\tau)d\tau,
$$

\vspace{3mm}
$$
G_p^{(2)}(t_2,\ldots,t_{k-1})=
-{\bf 1}_{\{t_2< \ldots <t_{k-1}\}}
\sum_{j_l=p+1}^{\infty}~
\int\limits_{t}^{t_{k-1}} \psi_k(\tau)\phi_{j_{l}}(\tau)d\tau
\int\limits_{t}^{t_{2}} \psi_1(\tau)\phi_{j_{l}}(\tau)d\tau,
$$

\vspace{3mm}
$$
E_p^{(1)}(t_2,\ldots,t_{l-1}, t_{l+1}, \ldots, t_k)=
$$

\vspace{-2mm}
$$
=
{\bf 1}_{\{t_2< \ldots <t_{l-1}<t_{l+1}<\ldots < t_k\}}
\sum_{j_l=p+1}^{\infty}~
\int\limits_{t}^{t_{l+1}} \psi_l(\tau)\phi_{j_{l}}(\tau)d\tau
\int\limits_{t}^{t_{2}} \psi_1(\tau)\phi_{j_{l}}(\tau)d\tau,
$$

\vspace{3mm}
$$
E_p^{(2)}(t_2,\ldots,t_{l-1}, t_{l+1}, \ldots, t_k)=
$$

\vspace{-2mm}
$$
=
-{\bf 1}_{\{t_2< \ldots <t_{l-1}<t_{l+1}<\ldots < t_k\}}
\sum_{j_l=p+1}^{\infty}~
\int\limits_{t}^{t_{l-1}} \psi_l(\tau)\phi_{j_{l}}(\tau)d\tau
\int\limits_{t}^{t_{2}} \psi_1(\tau)\phi_{j_{l}}(\tau)d\tau,
$$

\vspace{3mm}
$$
D_p^{(1)}(t_1,\ldots,t_{s-1}, t_{s+1}, \ldots, t_{k-1})=
$$

\vspace{-2mm}
$$
=
{\bf 1}_{\{t_1< \ldots <t_{s-1}<t_{s+1}<\ldots < t_{k-1}\}}
\sum_{j_l=p+1}^{\infty}~
\int\limits_{t}^{T} \psi_k(\tau)\phi_{j_{l}}(\tau)d\tau
\int\limits_{t}^{t_{s+1}} \psi_s(\tau)\phi_{j_{l}}(\tau)d\tau,
$$

\vspace{3mm}
$$
D_p^{(2)}(t_1,\ldots,t_{s-1}, t_{s+1}, \ldots, t_{k-1})=
$$

\vspace{-2mm}
$$
=
-{\bf 1}_{\{t_1< \ldots <t_{s-1}<t_{s+1}<\ldots < t_{k-1}\}}
\sum_{j_l=p+1}^{\infty}~
\int\limits_{t}^{T} \psi_k(\tau)\phi_{j_{l}}(\tau)d\tau
\int\limits_{t}^{t_{s-1}} \psi_s(\tau)\phi_{j_{l}}(\tau)d\tau,
$$

\vspace{3mm}
$$
D_p^{(3)}(t_1,\ldots,t_{s-1}, t_{s+1}, \ldots, t_{k-1})=
$$

\vspace{-2mm}
$$
=
-{\bf 1}_{\{t_1< \ldots <t_{s-1}<t_{s+1}<\ldots < t_{k-1}\}}
\sum_{j_l=p+1}^{\infty}~
\int\limits_{t}^{t_{k-1}} \psi_k(\tau)\phi_{j_{l}}(\tau)d\tau
\int\limits_{t}^{t_{s+1}} \psi_s(\tau)\phi_{j_{l}}(\tau)d\tau,
$$

\vspace{3mm}
$$
D_p^{(4)}(t_1,\ldots,t_{s-1}, t_{s+1}, \ldots, t_{k-1})=
$$

\vspace{-2mm}
$$
=
{\bf 1}_{\{t_1< \ldots <t_{s-1}<t_{s+1}<\ldots < t_{k-1}\}}
\sum_{j_l=p+1}^{\infty}~
\int\limits_{t}^{t_{k-1}} \psi_k(\tau)\phi_{j_{l}}(\tau)d\tau
\int\limits_{t}^{t_{s-1}} \psi_s(\tau)\phi_{j_{l}}(\tau)d\tau.
$$

\vspace{5mm}

Now let us consider the value
$C_{j_k\ldots j_1}\bigl|_{j_{g_1}=j_{g_2},\ g_2=g_1+1}\biggr.$. 
To do this, we will make the following transformations

\vspace{-1mm}
$$
\int\limits_t^T h_{k}(t_k)\ldots \int\limits_t^{t_{l+2}} h_{l+1}(t_{l+1})
\int\limits_t^{t_{l+1}} h_{l}(t_{l})
\int\limits_t^{t_{l}} h_{l}(t_{l-1})
\int\limits_t^{t_{l-1}} h_{l-2}(t_{l-2})\ldots
\int\limits_t^{t_2} h_{1}(t_1)
dt_1\ldots 
$$

\vspace{0.5mm}
$$
\ldots
dt_{l-2}dt_{l-1}dt_{l}dt_{l+1}\ldots dt_k=
$$

\vspace{0.5mm}
$$
=\int\limits_t^T h_{k}(t_k)\ldots \int\limits_t^{t_{l+2}} h_{l+1}(t_{l+1})
\int\limits_t^{t_{l+1}} h_{1}(t_{1})
\int\limits_{t_1}^{t_{l+1}} h_{2}(t_{2})\ldots
\int\limits_{t_{l-3}}^{t_{l+1}} h_{l-2}(t_{l-2})\times
$$

\vspace{0.5mm}
$$
\times
\left(\int\limits_{t}^{t_{l+1}}-
\int\limits_{t}^{t_{l-2}}~\right)h_{l}(t_{l-1})
\left(\int\limits_{t}^{t_{l+1}}-
\int\limits_{t}^{t_{l-1}}~\right)
h_{l}(t_{l})dt_l
dt_{l-1}dt_{l-2}\ldots dt_2dt_{1}dt_{l+1}\ldots dt_k=
$$

\vspace{2mm}
$$
=\int\limits_t^T h_{k}(t_k)\ldots \int\limits_t^{t_{l+2}} h_{l+1}(t_{l+1})
\left(\int\limits_t^{t_{l+1}} h_{l}(t_{l})dt_l
\int\limits_t^{t_{l+1}} h_{l}(t_{l-1})dt_{l-1}\right)
\int\limits_t^{t_{l+1}} h_{1}(t_{1})
\times
$$

\vspace{0.5mm}
$$
\times
\int\limits_{t_1}^{t_{l+1}} h_{2}(t_{2})\ldots \int\limits_{t_{l-3}}^{t_{l+1}} h_{l-2}(t_{l-2})
dt_{l-2}\ldots dt_2dt_{1}dt_{l+1}\ldots dt_k-
$$

\vspace{2mm}
$$
-\int\limits_t^T h_{k}(t_k)\ldots \int\limits_t^{t_{l+2}} h_{l+1}(t_{l+1})
\left(\int\limits_t^{t_{l+1}} h_{l}(t_{l})dt_l\right)
\int\limits_t^{t_{l+1}} h_{1}(t_{1})
\int\limits_{t_1}^{t_{l+1}} h_{2}(t_{2})\ldots 
$$

\vspace{0.5mm}
$$
\ldots \int\limits_{t_{l-3}}^{t_{l+1}} h_{l-2}(t_{l-2})
\left(\int\limits_t^{t_{l-2}} h_{l}(t_{l-1})dt_{l-1}\right)
dt_{l-2}\ldots dt_2dt_{1}dt_{l+1}\ldots dt_k-
$$

\vspace{2mm}
$$
-\int\limits_t^T h_{k}(t_k)\ldots \int\limits_t^{t_{l+2}} h_{l+1}(t_{l+1})
\left(\int\limits_t^{t_{l+1}} h_{l}(t_{l-1})
\int\limits_t^{t_{l-1}} h_{l}(t_{l})dt_l dt_{l-1}\right)
\int\limits_t^{t_{l+1}} h_{1}(t_{1})\times
$$

\vspace{0.5mm}
$$
\times
\int\limits_{t_1}^{t_{l+1}} h_{2}(t_{2})\ldots 
\int\limits_{t_{l-3}}^{t_{l+1}} h_{l-2}(t_{l-2})
dt_{l-2}\ldots dt_2dt_{1}dt_{l+1}\ldots dt_k+
$$

\vspace{2mm}
$$
+\int\limits_t^T h_{k}(t_k)\ldots \int\limits_t^{t_{l+2}} h_{l+1}(t_{l+1})
\int\limits_t^{t_{l+1}} h_{1}(t_{1})
\int\limits_{t_1}^{t_{l+1}} h_{2}(t_{2})\ldots \int\limits_{t_{l-3}}^{t_{l+1}} h_{l-2}(t_{l-2})\times
$$

\vspace{0.5mm}
$$
\times
\left(\int\limits_t^{t_{l-2}} h_{l}(t_{l-1})
\int\limits_t^{t_{l-1}} h_{l}(t_{l})dt_ldt_{l-1}\right)
dt_{l-2}\ldots dt_2dt_{1}dt_{l+1}\ldots dt_k=
$$

\vspace{2mm}
$$
=\int\limits_t^T h_{k}(t_k)\ldots \int\limits_t^{t_{l+2}} h_{l+1}(t_{l+1})
\left(\int\limits_t^{t_{l+1}} h_{l}(t_{l})dt_l
\int\limits_t^{t_{l+1}} h_{l}(t_{l-1})dt_{l-1}\right)
\int\limits_t^{t_{l+1}} h_{l-2}(t_{l-2})
\times
$$

\vspace{0.5mm}
$$
\times
\int\limits_{t}^{t_{l-2}} h_{l-3}(t_{l-3})\ldots \int\limits_{t}^{t_2} h_{1}(t_{1})dt_1
\ldots dt_{l-3}dt_{l-2}dt_{l+1}\ldots dt_k-
$$

\vspace{2mm}
$$
-\int\limits_t^T h_{k}(t_k)\ldots \int\limits_t^{t_{l+2}} h_{l+1}(t_{l+1})
\left(\int\limits_t^{t_{l+1}} h_{l}(t_{l})dt_l\right)
\int\limits_t^{t_{l+1}} h_{l-2}(t_{l-2})
\times
$$

\vspace{0.5mm}
$$
\times\left(\int\limits_t^{t_{l-2}} h_{l}(t_{l-1})dt_{l-1}\right)
\int\limits_{t}^{t_{l-2}} h_{l-3}(t_{l-3})\ldots \int\limits_{t}^{t_2} h_{1}(t_{1})dt_1
\ldots dt_{l-3}dt_{l-2}dt_{l+1}\ldots dt_k-
$$

\vspace{2mm}
$$
-\int\limits_t^T h_{k}(t_k)\ldots \int\limits_t^{t_{l+2}} h_{l+1}(t_{l+1})
\left(\int\limits_t^{t_{l+1}} h_{l}(t_{l-1})
\int\limits_t^{t_{l-1}} h_{l}(t_{l})dt_l dt_{l-1}\right)
\times
$$

\vspace{0.5mm}
$$
\times
\int\limits_t^{t_{l+1}} h_{l-2}(t_{l-2})
\int\limits_{t}^{t_{l-2}} h_{l-3}(t_{l-3})\ldots \int\limits_{t}^{t_2} h_{1}(t_{1})dt_1
\ldots dt_{l-3}dt_{l-2}dt_{l+1}\ldots dt_k+
$$

\vspace{2mm}
$$
+\int\limits_t^T h_{k}(t_k)\ldots \int\limits_t^{t_{l+2}} h_{l+1}(t_{l+1})
\int\limits_t^{t_{l+1}}
h_{l-2}(t_{l-2})
\left(\int\limits_t^{t_{l-2}} h_{l}(t_{l-1})
\int\limits_t^{t_{l-1}} h_{l}(t_{l})dt_ldt_{l-1}\right)
\times
$$

\vspace{0.5mm}
\begin{equation}
\label{after9031}
\times
\int\limits_{t}^{t_{l-2}} h_{l-3}(t_{l-3})\ldots \int\limits_{t}^{t_2} h_{1}(t_{1})dt_1
\ldots dt_{l-3}dt_{l-2}dt_{l+1}\ldots dt_k,
\end{equation}

\vspace{2mm}
\noindent
where $l+1\le k,$ $l-2\ge 1,$ and
$h_1(\tau),\ldots,h_k(\tau)$ are continuous functions on the interval
$[t, T].$

Applying (\ref{after9031}) 
to $C_{j_k \ldots j_{l+1} j_l j_l j_{l-2} \ldots \ldots j_1}$, we obtain
for $l+1\le k,$ $l-2\ge 1$

$$
\sum_{j_l=p+1}^{\infty}
C_{j_k \ldots j_{l+1} j_l j_{l} j_{l-2} \ldots \ldots j_1}=
$$

$$
=\int\limits_{[t, T]^{k-2}} \sum\limits_{d=1}^4
H_p^{(d)}(t_1,\ldots,t_{l-2},t_{l+1},\ldots,t_k)\times
$$

$$
\times
\prod_{\stackrel{g=1}{{}_{g\ne l-1, l}}}^{k}
\psi_g(t_g)\phi_{j_g}(t_g)
dt_1\ldots dt_{l-2}dt_{l+1}\ldots dt_k=
$$

\begin{equation}
\label{after85ayyy}
=\sum\limits_{d=1}^4
C_{j_k \ldots j_{l+1} j_{l-2} \ldots j_1}^{**(d)}
=\sum\limits_{d=1}^4
C_{j_k \ldots j_q \ldots j_1}^{**(d)}\biggl|_{q\ne l-1, l}\biggr.,
\end{equation}

\vspace{5mm}
\noindent
where

\vspace{-3mm}
$$
H_p^{(1)}(t_1,\ldots,t_{l-2},t_{l+1},\ldots,t_k)=
$$

\vspace{-2mm}
\begin{equation}
\label{after90yyy}
=
{\bf 1}_{\{t_1<\ldots<t_{l-2}<t_{l+1}<\ldots<t_k\}}
\sum_{j_l=p+1}^{\infty}~
\int\limits_{t}^{t_{l+1}} \psi_l(\tau) \phi_{j_{l}}(\tau)d\tau
\int\limits_{t}^{t_{l+1}} \psi_{l-1}(\tau) \phi_{j_{l}}(\tau)d\tau,
\end{equation}

\vspace{3mm}
$$
H_p^{(2)}(t_1,\ldots,t_{l-2},t_{l+1},\ldots,t_k)=
$$

\vspace{-2mm}
\begin{equation}
\label{after91yyy}
=
-{\bf 1}_{\{t_1<\ldots<t_{l-2}<t_{l+1}<\ldots<t_k\}}
\sum_{j_l=p+1}^{\infty}~
\int\limits_{t}^{t_{l+1}} \psi_l(\tau) \phi_{j_{l}}(\tau)d\tau
\int\limits_{t}^{t_{l-2}} \psi_{l-1}(\tau) \phi_{j_{l}}(\tau)d\tau,
\end{equation}

\vspace{3mm}
$$
H_p^{(3)}(t_1,\ldots,t_{l-2},t_{l+1},\ldots,t_k)=
$$

\vspace{-2mm}
\begin{equation}
\label{after92yyy}
=
-{\bf 1}_{\{t_1<\ldots<t_{l-2}<t_{l+1}<\ldots<t_k\}}
\sum_{j_l=p+1}^{\infty}~
\int\limits_{t}^{t_{l+1}} \psi_{l-1}(\tau)\phi_{j_{l}}(\tau)
\int\limits_{t}^{\tau} \psi_l(\theta)\phi_{j_{l}}(\theta)d\theta d\tau,
\end{equation}

\vspace{3mm}
$$
H_p^{(4)}(t_1,\ldots,t_{l-2},t_{l+1},\ldots,t_k)=
$$

\vspace{-2mm}
\begin{equation}
\label{after93yyy}
=
{\bf 1}_{\{t_1<\ldots<t_{l-2}<t_{l+1}<\ldots<t_k\}}
\sum_{j_l=p+1}^{\infty}~
\int\limits_{t}^{t_{l-2}} \psi_{l-1}(\tau)\phi_{j_{l}}(\tau)
\int\limits_{t}^{\tau} \psi_l(\theta)\phi_{j_{l}}(\theta)d\theta d\tau.
\end{equation}

\vspace{5mm}

By analogy with (\ref{after85ayyy}) we can consider the expressions

\vspace{-2mm}
\begin{equation}
\label{afterx300}
\sum_{j_l=p+1}^{\infty}
C_{j_k \ldots j_{l+1} j_l j_l},
\end{equation}

\vspace{2mm}
\begin{equation}
\label{afterx301}
\sum_{j_l=p+1}^{\infty}
C_{j_l j_l j_{k-2} \ldots j_1}.
\end{equation}

\vspace{3mm}

Then we have for (\ref{afterx300}), (\ref{afterx301}) 
(see (\ref{after9031}) and its analogue for $t_{l+1}=T$)

\begin{equation}
\label{afterx500}
\sum_{j_l=p+1}^{\infty}
C_{j_k \ldots j_{l+1} j_l j_{l}}=
\int\limits_{[t, T]^{k-2}} 
L_p(t_3,\ldots,t_k)\prod_{g=3}^{k}
\psi_g(t_g)\phi_{j_g}(t_g)
dt_3\ldots dt_k,
\end{equation}

\vspace{2mm}
\begin{equation}
\label{afterx501}
\sum_{j_l=p+1}^{\infty}
C_{j_l j_l j_{k-2} \ldots j_{1}}=
\int\limits_{[t, T]^{k-2}} 
\sum\limits_{d=1}^4 M_p^{(d)}(t_1,\ldots,t_{k-2})\prod_{g=1}^{k-2}
\psi_g(t_g)\phi_{j_g}(t_g)
dt_1\ldots dt_{k-2},
\end{equation}

\vspace{5mm}
\noindent
where
$$
L_p(t_3,\ldots,t_k)
=
{\bf 1}_{\{t_3<\ldots<t_k\}}
\sum_{j_l=p+1}^{\infty}~
\int\limits_{t}^{t_{3}} \psi_2(\tau) \phi_{j_{l}}(\tau)
\int\limits_{t}^{\tau} \psi_{1}(\theta) \phi_{j_{l}}(\theta)d\theta d\tau,
$$

\vspace{3mm}
$$
M_p^{(1)}(t_1,\ldots,t_{k-2})=
$$

\vspace{-2mm}
$$
=
{\bf 1}_{\{t_1<\ldots<t_{k-2}\}}
\sum_{j_l=p+1}^{\infty}~
\int\limits_{t}^{T} \psi_k(\tau) \phi_{j_{l}}(\tau)d\tau
\int\limits_{t}^{T} \psi_{k-1}(\tau) \phi_{j_{l}}(\tau)d\tau,
$$

\vspace{3mm}
$$
M_p^{(2)}(t_1,\ldots,t_{k-2})=
$$

\vspace{-2mm}
$$
=
-{\bf 1}_{\{t_1<\ldots<t_{k-2}\}}
\sum_{j_l=p+1}^{\infty}~
\int\limits_{t}^{T} \psi_k(\tau) \phi_{j_{l}}(\tau)d\tau
\int\limits_{t}^{t_{k-2}} \psi_{k-1}(\tau) \phi_{j_{l}}(\tau)d\tau,
$$

\vspace{3mm}
$$
M_p^{(3)}(t_1,\ldots,t_{k-2})=
$$

\vspace{-2mm}
$$
=
-{\bf 1}_{\{t_1<\ldots<t_{k-2}\}}
\sum_{j_l=p+1}^{\infty}~
\int\limits_{t}^{T} \psi_{k-1}(\tau) \phi_{j_{l}}(\tau)
\int\limits_{t}^{\tau} \psi_{k}(\theta) \phi_{j_{l}}(\theta)d\theta d\tau,
$$

\vspace{3mm}
$$
M_p^{(4)}(t_1,\ldots,t_{k-2})=
$$

\vspace{-2mm}
$$
=
{\bf 1}_{\{t_1<\ldots<t_{k-2}\}}
\sum_{j_l=p+1}^{\infty}~
\int\limits_{t}^{t_{k-2}} \psi_{k-1}(\tau) \phi_{j_{l}}(\tau)
\int\limits_{t}^{\tau} \psi_{k}(\theta) \phi_{j_{l}}(\theta)d\theta d\tau.
$$

\vspace{4mm}

It is important to note that 
$C_{j_k \ldots j_{l+1} j_{l-2} \ldots j_1}^{*(d)},$
$C_{j_k \ldots j_{l+1} j_{l-2} \ldots j_1}^{**(d)}$
$(d=1,\ldots,4)$ 
are Fourier coefficients (see (\ref{after85a}), (\ref{after85ayyy})), 
that is, we can use Parseval's equality 
in the further proof.

Combining the equalities (\ref{after85a})--(\ref{after93})
(the case $g_2>g_1+1$), using Parseval's equality and applying the 
estimates for integrals from basis functions that we 
used in the proof of Theorems 16, 17, we obtain for (\ref{after85a}) 

$$
\sum\limits_{j_{q_1},\ldots,j_{q_{k-2}}=0}^p
\left(\sum_{j_{g_1}=p+1}^{\infty}
C_{j_k\ldots j_1}\biggl|_{j_{g_1}=j_{g_2}, g_2>g_1+1}\biggr.\right)^2=
$$

\vspace{1mm}
$$
=\sum\limits_{\stackrel{j_{1},\ldots,j_q,\ldots,j_{k}=0}{{}_{q\ne g_1, g_2}}}^p
\left(\sum_{j_{g_1}=p+1}^{\infty}
C_{j_k\ldots j_1}\biggl|_{j_{g_1}=j_{g_2}, g_2>g_1+1}\biggr.\right)^2=
$$

\vspace{1mm}
$$
= 
\sum\limits_{\stackrel{j_{1},\ldots,j_q,\ldots,j_{k}=0}{{}_{q\ne g_1, g_2}}}^p
\left(\sum\limits_{d=1}^4 
C_{j_k \ldots j_q \ldots j_1}^{*(d)}\biggl|_{q\ne g_1, g_2}\biggr.\right)^2
\le
\sum\limits_{\stackrel{j_{1},\ldots,j_q,\ldots,j_{k}=0}{{}_{q\ne g_1, g_2}}}^{\infty}
\left(\sum\limits_{d=1}^4
C_{j_k \ldots j_q \ldots j_1}^{*(d)}\biggl|_{q\ne g_1, g_2}\biggr.\right)^2
=
$$

\vspace{3mm}
$$
=\sum\limits_{\stackrel{j_{1},\ldots,j_q,\ldots,j_{k}=0}{{}_{q\ne g_1, g_2}}}^{\infty}
\left(~
\int\limits_{[t, T]^{k-2}} \sum\limits_{d=1}^4
F_p^{(d)}(t_1,\ldots,t_{g_1-1},t_{g_1+1},\ldots ,t_{g_2-1},t_{g_2+1},\ldots,t_k)\times\right.
$$

\vspace{1mm}
$$
\left.\times 
\prod_{\stackrel{q=1}{{}_{q\ne g_1, g_2}}}^{k}
\psi_q(t_q)\phi_{j_q}(t_q)
dt_1\ldots dt_{g_1-1}dt_{g_1+1}\ldots dt_{g_2-1}dt_{g_2+1}\ldots dt_k\right)^2=
$$

\vspace{1mm}
$$
=\int\limits_{[t, T]^{k-2}} 
\left(~
\sum\limits_{d=1}^4
F_p^{(d)}(t_1,\ldots,t_{g_1-1},t_{g_1+1},\ldots, t_{g_2-1},t_{g_2+1},\ldots,t_k)
\prod_{\stackrel{q=1}{{}_{q\ne g_1, g_2}}}^{k}
\psi_q(t_q)\right)^2\times
$$

\vspace{1mm}
$$
\times
dt_1\ldots dt_{g_1-1}dt_{g_1+1}\ldots dt_{g_2-1}dt_{g_2+1}\ldots dt_k\le
$$

\vspace{1mm}
$$
\le 4\sum\limits_{d=1}^4 \int\limits_{[t, T]^{k-2}} 
\left(
F_p^{(d)}(t_1,\ldots,t_{g_1-1},t_{g_1+1},\ldots, t_{g_2-1},t_{g_2+1},\ldots,t_k)
\prod_{\stackrel{q=1}{{}_{q\ne g_1, g_2}}}^{k}
\psi_q(t_q)\right)^2\times
$$

\vspace{1mm}
$$
\times dt_1\ldots dt_{g_1-1}dt_{g_1+1}\ldots dt_{g_2-1}dt_{g_2+1}\ldots dt_k\le
$$

\vspace{1mm}
\begin{equation}
\label{after12001}
\le \frac{K}{p^{2-\varepsilon}}\ \to\  0
\end{equation}

\vspace{4mm}
\noindent
if $p\to\infty,$ where $\varepsilon$ is an arbitrary small positive real number
for the polynomial case and $\varepsilon=0$ for the 
trigonometric case, constant $K$ does not depend on $p.$
The cases (\ref{afterx200})--(\ref{afterx202}) are considered
analogously.

Absolutely similarly (see (\ref{after12001}))
combining the equalities (\ref{after85ayyy})--(\ref{after93yyy})
(the case $g_2=g_1+1$), using Parseval's equality and applying the 
estimates for integrals from basis functions that we 
used in the proof of Theorems 16, 17, we get for (\ref{after85ayyy})

$$
\sum\limits_{j_{q_1},\ldots,j_{q_{k-2}}=0}^p
\left(\sum_{j_{g_1}=p+1}^{\infty}
C_{j_k\ldots j_1}\biggl|_{j_{g_1}=j_{g_2}, g_2=g_1+1}\biggr.\right)^2=
$$

\vspace{1mm}
$$
=\sum\limits_{\stackrel{j_{1},\ldots,j_q,\ldots,j_{k}=0}{{}_{q\ne g_1, g_2}}}^p
\left(\sum_{j_{g_1}=p+1}^{\infty}
C_{j_k\ldots j_1}\biggl|_{j_{g_1}=j_{g_2}, g_2=g_1+1}\biggr.\right)^2=
$$

\vspace{1mm}
$$
= 
\sum\limits_{\stackrel{j_{1},\ldots,j_q,\ldots,j_{k}=0}{{}_{q\ne g_1, g_2}}}^p
\left(\sum\limits_{d=1}^4 
C_{j_k \ldots j_q \ldots j_1}^{**(d)}\biggl|_{q\ne g_1, g_2}\biggr.\right)^2
\le
\sum\limits_{\stackrel{j_{1},\ldots,j_q,\ldots,j_{k}=0}{{}_{q\ne g_1, g_2}}}^{\infty}
\left(\sum\limits_{d=1}^4
C_{j_k \ldots j_q \ldots j_1}^{**(d)}\biggl|_{q\ne g_1, g_2}\biggr.\right)^2
=
$$

\vspace{3mm}
$$
=\sum\limits_{\stackrel{j_{1},\ldots,j_q,\ldots,j_{k}=0}{{}_{q\ne g_1, g_2}}}^{\infty}
\left(~
\int\limits_{[t, T]^{k-2}} \sum\limits_{d=1}^4
H_p^{(d)}(t_1,\ldots,t_{g_1-1},t_{g_1+2},\ldots,t_k)\times\right.
$$

\vspace{1mm}
$$
\left.\times 
\prod_{\stackrel{q=1}{{}_{q\ne g_1, g_2}}}^{k}
\psi_q(t_q)\phi_{j_q}(t_q)
dt_1\ldots dt_{g_1-1}dt_{g_1+2}\ldots dt_k\right)^2=
$$

\vspace{1mm}
$$
=\int\limits_{[t, T]^{k-2}} 
\hspace{-1.5mm}\left(~\hspace{-1.5mm}
\sum\limits_{d=1}^4
H_p^{(d)}(t_1,\ldots,t_{g_1-1},t_{g_1+2},\ldots,t_k)\prod_{\stackrel{q=1}{{}_{q\ne g_1, g_2}}}^{k}
\psi_q(t_q)\right)^2
dt_1\ldots dt_{g_1-1}dt_{g_1+2}\ldots \ldots dt_k\le
$$

\vspace{1mm}
$$
\le 4\sum\limits_{d=1}^4 \int\limits_{[t, T]^{k-2}} 
\left(
H_p^{(d)}(t_1,\ldots,t_{g_1-1},t_{g_1+2},\ldots,t_k)\prod_{\stackrel{q=1}{{}_{q\ne g_1, g_2}}}^{k}
\psi_q(t_q)\right)^2
dt_1\ldots dt_{g_1-1}dt_{g_1+2}\ldots dt_k\le
$$

\vspace{1mm}
\begin{equation}
\label{after12000}
\le \frac{K}{p^{2-\varepsilon}}\ \to\  0
\end{equation}

\vspace{4mm}
\noindent
if $p\to\infty,$ where $\varepsilon$ is an arbitrary small positive real number
for the polynomial case and $\varepsilon=0$ for the 
trigonometric case, constant $K$ does not depend on $p.$
The cases (\ref{afterx300}), (\ref{afterx301}) are considered
analogously.

From (\ref{after12001}), (\ref{after12000}) and their
analogues for the cases (\ref{afterx200})--(\ref{afterx202}),
(\ref{afterx300}), (\ref{afterx301})
we obtain 

\vspace{-1mm}
\begin{equation}
\label{after17000}
\sum\limits_{j_{q_1},\ldots,j_{q_{k-2}}=0}^p
\left(\sum_{j_{g_1}=p+1}^{\infty}
C_{j_k\ldots j_1}\biggl|_{j_{g_1}=j_{g_2}}\biggr.\right)^2
\le \frac{K}{p^{2-\varepsilon}},
\end{equation}

\vspace{4mm}
\noindent
where constant $K$ is independent of $p.$ 
Thus the equality (\ref{after14000}) is proved.

Let us prove the equality (\ref{after14001}).
Consider the following cases

\vspace{4mm}
\centerline{1.\ $g_2>g_1+1$,\ $g_4=g_3+1$,\ \ \ \ 2.\ $g_2=g_1+1$,\ $g_4>g_3+1$,}

\vspace{2mm}
\centerline{3.\ $g_2>g_1+1$,\ $g_4>g_3+1$,\ \ \ \ 4.\ $g_2=g_1+1$,\ $g_4=g_3+1$.}

\vspace{4mm}

The proof for Cases 1--3 will be similar. Consider, for example, Case 2.
Using (\ref{after79}), we obtain

$$
\sum\limits_{j_{q_1}=0}^p
\left(\sum_{j_{g_1}=p+1}^{\infty}\sum_{j_{g_3}=p+1}^{\infty}
C_{j_5\ldots j_1}\biggl|_{j_{g_1}=j_{g_2},j_{g_3}=j_{g_4}, g_4>g_3+1, g_2=g_1+1}\biggr.\right)^2=
$$

\vspace{1mm}
$$
=
\sum\limits_{j_{q_1}=0}^p
\left(\sum_{j_{g_1}=p+1}^{\infty}\sum_{j_{g_3}=0}^{p}
C_{j_5\ldots j_1}\biggl|_{j_{g_1}=j_{g_2},j_{g_3}=j_{g_4}, g_4>g_3+1, g_2=g_1+1}\biggr.\right)^2=
$$

\vspace{1mm}
\begin{equation}
\label{afterpp1}
=\sum\limits_{j_{q_1}=0}^p
\left(\sum_{j_{g_3}=0}^{p}\sum_{j_{g_1}=p+1}^{\infty}
C_{j_5\ldots j_1}\biggl|_{j_{g_1}=j_{g_2},j_{g_3}=j_{g_4}, g_4>g_3+1, g_2=g_1+1}\biggr.\right)^2\le
\end{equation}

\vspace{1mm}
$$
\le (p+1)\sum\limits_{j_{q_1}=0}^p\sum_{j_{g_3}=0}^{p}
\left(\sum_{j_{g_1}=p+1}^{\infty}
C_{j_5\ldots j_1}\biggl|_{j_{g_1}=j_{g_2},j_{g_3}=j_{g_4}, g_4>g_3+1, g_2=g_1+1}\biggr.\right)^2=
$$

\vspace{1mm}
$$
=(p+1)\sum\limits_{j_{q_1}=0}^p\sum_{j_{g_3}, j_{g_4}=0}^{p}
\left(\sum_{j_{g_1}=p+1}^{\infty}
C_{j_5\ldots j_1}\biggl|_{j_{g_1}=j_{g_2}, g_4>g_3+1, g_2=g_1+1}\biggr.
\right)^2\Biggl|_{j_{g_3}=j_{g_4}}\le
$$

\vspace{1mm}
\begin{equation}
\label{after15000}
\le(p+1)\sum\limits_{j_{q_1}=0}^p\sum_{j_{g_3}, j_{g_4}=0}^{p}
\left(\sum_{j_{g_1}=p+1}^{\infty}
C_{j_5\ldots j_1}\biggl|_{j_{g_1}=j_{g_2}, g_4>g_3+1, g_2=g_1+1}\biggr.
\right)^2.
\end{equation}

\vspace{5mm}

It is easy to see that the expression 
(\ref{after15000}) (without the multiplier $p+1$) 
is a particular case $(g_4>g_3+1, g_2=g_1+1)$ of the left-hand side
of (\ref{after17000}).
Combining (\ref{after17000}) and (\ref{after15000}), we have

$$
\sum\limits_{j_{q_1}=0}^p
\left(\sum_{j_{g_1}=p+1}^{\infty}\sum_{j_{g_3}=p+1}^{\infty}
C_{j_5\ldots j_1}\biggl|_{j_{g_1}=j_{g_2},j_{g_3}=j_{g_4}, g_4>g_3+1, g_2=g_1+1}\biggr.\right)^2\le
$$

\vspace{1mm}
\begin{equation}
\label{after15000a}
\le \frac{(p+1)K}{p^{2-\varepsilon}}\le
\frac{K_1}{p^{1-\varepsilon}}\ \to\  0
\end{equation}

\vspace{5mm}
\noindent
if $p\to\infty,$ where constant $K_1$ does not depend on $p.$

Consider Case 4 ($g_2=g_1+1$,\ $g_4=g_3+1$). We have (see (\ref{after500}))

$$
\sum\limits_{j_{q_1}=0}^p
\left(\sum_{j_{g_1}=p+1}^{\infty}\sum_{j_{g_3}=p+1}^{\infty}
C_{j_5\ldots j_1}\biggl|_{j_{g_1}=j_{g_2},j_{g_3}=j_{g_4}}\biggr.\right)^2=
$$

\vspace{1mm}
$$
=\sum\limits_{j_{q_1}=0}^p
\left(\sum_{j_{g_1}=p+1}^{\infty}\left(\sum_{j_{g_3}=0}^{\infty}-\sum_{j_{g_3}=0}^{p}\right)
C_{j_5\ldots j_1}\biggl|_{j_{g_1}=j_{g_2},j_{g_3}=j_{g_4}}\biggr.\right)^2=
$$

\vspace{1mm}
$$
=\sum\limits_{j_{q_1}=0}^p
\left(\frac{1}{2}
\sum_{j_{g_1}=p+1}^{\infty}
C_{j_5\ldots j_1}\biggl|_{j_{g_1}=j_{g_2},
(j_{g_3} j_{g_3})\curvearrowright (\cdot)}\biggr.
-\sum_{j_{g_3}=0}^{p}\sum_{j_{g_1}=p+1}^{\infty}
C_{j_5\ldots j_1}\biggl|_{j_{g_1}=j_{g_2},j_{g_3}=j_{g_4}}\biggr.\right)^2\le
$$

\vspace{1mm}
\begin{equation}
\label{after19000}
\le
\frac{1}{2}\sum\limits_{j_{q_1}=0}^p
\left(
\sum_{j_{g_1}=p+1}^{\infty}
C_{j_5\ldots j_1}\biggl|_{j_{g_1}=j_{g_2},
(j_{g_3} j_{g_3})\curvearrowright (\cdot)}\biggr.\right)^2+
\end{equation}

\vspace{1mm}
\begin{equation}
\label{after19001}
+2\sum\limits_{j_{q_1}=0}^p
\left(\sum_{j_{g_3}=0}^{p}\sum_{j_{g_1}=p+1}^{\infty}
C_{j_5\ldots j_1}\biggl|_{j_{g_1}=j_{g_2},j_{g_3}=j_{g_4}}\biggr.\right)^2.
\end{equation}

\vspace{5mm}

An expression similar to (\ref{after19001}) was estimated 
(see (\ref{afterpp1})--(\ref{after15000a})). Let us estimate 
(\ref{after19000}). We have

$$
\sum\limits_{j_{q_1}=0}^p
\left(
\sum_{j_{g_1}=p+1}^{\infty}
C_{j_5\ldots j_1}\biggl|_{j_{g_1}=j_{g_2},
(j_{g_3} j_{g_3})\curvearrowright (\cdot)}\biggr.\right)^2=
$$

\vspace{1mm}
$$
=(T-t)\sum\limits_{j_{q_1}=0}^p
\left(
\sum_{j_{g_1}=p+1}^{\infty}
C_{j_5\ldots j_1}\biggl|_{j_{g_1}=j_{g_2},
(j_{g_3} j_{g_3})\curvearrowright 0}\biggr.\right)^2\le
$$

\vspace{1mm}
\begin{equation}
\label{after19005}
\le (T-t)\sum\limits_{j_{q_1}=0}^p
\sum_{j_{g_3}=0}^{p}\left(
\sum_{j_{g_1}=p+1}^{\infty}
C_{j_5\ldots j_1}\biggl|_{j_{g_1}=j_{g_2},
(j_{g_3} j_{g_3})\curvearrowright j_{g_3}}\biggr.\right)^2,
\end{equation}

\vspace{4mm}
\noindent
where the notations are the same as in the proof of Theorem~13.

The expression (\ref{after19005}) 
without the multiplier $T-t$ is an expression of type 
(\ref{after2501})--(\ref{after2506})
before passing to the limit $\lim\limits_{p\to\infty}$ 
(the only difference is the replacement of one of the weight functions 
$\psi_1(\tau),\ldots, \psi_4(\tau)$ in 
(\ref{after2501})--(\ref{after2506}) by the product 
$\psi_{l+1}(\tau)\psi_l(\tau)$ $(l=1,\ldots,4).$
Therefore, for Case 4 ($g_2=g_1+1$,\ $g_4=g_3+1$), we obtain the estimate

$$
\sum\limits_{j_{q_1}=0}^p
\left(\sum_{j_{g_1}=p+1}^{\infty}\sum_{j_{g_3}=p+1}^{\infty}
C_{j_5\ldots j_1}\biggl|_{j_{g_1}=j_{g_2},j_{g_3}=j_{g_4}, g_4=g_3+1, g_2=g_1+1}\biggr.\right)^2\le
$$

\vspace{1mm}
\begin{equation}
\label{after20000}
\le
\frac{K}{p^{1-\varepsilon}},
\end{equation}

\vspace{4mm}
\noindent
where constant $K$ is independent of $p.$

The estimates (\ref{after15000a}), (\ref{after20000}) 
prove (\ref{after14001}). 
Let us prove (\ref{afterafter001}). By analogy with (\ref{after19005}) we have

\vspace{-1mm}
$$
\sum\limits_{j_{q_1}=0}^p
\left(\sum_{j_{g_3}=p+1}^{\infty}
C_{j_5\ldots j_1}\biggl|_{(j_{g_2}j_{g_1})\curvearrowright (\cdot),
j_{g_1}=j_{g_2},
j_{g_3}=j_{g_4}, g_2=g_1+1}\biggr.\right)^2=
$$

\vspace{1mm}
$$
=\sum\limits_{j_{q_1}=0}^p
\left(\sum_{j_{g_3}=p+1}^{\infty}
C_{j_5\ldots j_1}\biggl|_{(j_{g_1}j_{g_1})\curvearrowright (\cdot),
j_{g_3}=j_{g_4}, g_2=g_1+1}\biggr.\right)^2=
$$

\vspace{1mm}
$$
=(T-t)\sum\limits_{j_{q_1}=0}^p
\left(\sum_{j_{g_3}=p+1}^{\infty}
C_{j_5\ldots j_1}\biggl|_{(j_{g_1}j_{g_1})\curvearrowright 0,
j_{g_3}=j_{g_4}, g_2=g_1+1}\biggr.\right)^2\le
$$

\vspace{1mm}
\begin{equation}
\label{afterafter120}
\le
(T-t)\sum\limits_{j_{q_1}=0}^p \sum\limits_{j_{g_1}=0}^p
\left(\sum_{j_{g_3}=p+1}^{\infty}
C_{j_5\ldots j_1}\biggl|_{(j_{g_1}j_{g_1})\curvearrowright j_{g_1},
j_{g_3}=j_{g_4}, g_2=g_1+1}\biggr.\right)^2.
\end{equation}

\vspace{5mm}

Thus, we obtain the estimate (see (\ref{after19005}) and the proof of Theorem~17)

\vspace{-1mm}
$$
\sum\limits_{j_{q_1}=0}^p
\left(\sum_{j_{g_3}=p+1}^{\infty}
C_{j_5\ldots j_1}\biggl|_{(j_{g_2}j_{g_1})\curvearrowright (\cdot),
j_{g_1}=j_{g_2},
j_{g_3}=j_{g_4}, g_2=g_1+1}\biggr.\right)^2
\le
$$

\vspace{1mm}
\begin{equation}
\label{afterafter121}
\le \frac{K}{p^{2-\varepsilon}},
\end{equation}

\vspace{4mm}
\noindent
where 
$\varepsilon$ is an arbitrary small positive real number
for the polynomial case and $\varepsilon=0$ for the 
trigonometric case,
constant $K$ does not depend on $p.$

The estimate (\ref{afterafter121}) proves (\ref{afterafter001}).
Theorem~18 is proved.

\vspace{5mm}

\section{Estimates for the Mean-Square Approximation Error of Expansions of Iterated
Stra\-to\-no\-vich Stochastic Integrals of Multiplicity $k$
in Theorems 13, 15}

\vspace{5mm}

In this section, we estimate the mean-square approximation error
for iterated Stratonovich stochastic integrals of multiplicity
$k$ ($k\in \mathbb{N}$) in Theorems~13, 15.

\vspace{2mm}                           
                                                             
{\bf Theorem~19}\ \cite{2018a}, \cite{arxiv-5}, \cite{arxiv-9}, \cite{new-art-1-xxy}.\  
{\it Suppose that every $\psi_l(\tau)$ $(l=1,\ldots,k)$
is a continuously differentiable nonrandom func\-ti\-on
at the interval $[t, T].$ Furthermore$,$ let
$\{\phi_j(x)\}_{j=0}^{\infty}$ is a complete orthonormal system of 
Legendre polynomials or trigonometric functions in the space $L_2([t, T]).$
Then the following estimates 
$$
{\sf M}\left\{\left(
J^{*}[\psi^{(k)}]_{T,t}^{(i_1\ldots i_k)}-
\sum\limits_{j_1,\ldots,j_k=0}^{p}
C_{j_k \ldots j_1}\prod\limits_{l=1}^k \zeta_{j_l}^{(i_l)}
\right)^2\right\}\le
$$

\vspace{-0.5mm}
\begin{equation}
\label{after3407}
\le K_1 \left(\frac{1}{p} +\sum_{r=1}^{\left[k/2\right]}
\sum_{\stackrel{(\{\{g_1, g_2\}, \ldots, 
\{g_{2r-1}, g_{2r}\}\}, \{q_1, \ldots, q_{k-2r}\})}
{{}_{\{g_1, g_2, \ldots, 
g_{2r-1}, g_{2r}, q_1, \ldots, q_{k-2r}\}=\{1, 2, \ldots, k\}}}}
{\sf M}\left\{\left(R_{T,t}^{(p)r, g_1,g_2,\ldots,g_{2r-1}, g_{2r}}\right)^2\right\}\right),
\end{equation}

\vspace{2mm}
$$
{\sf M}\left\{\left(
J^{*}[\psi^{(k)}]_{s,t}^{(i_1\ldots i_k)}-
\sum\limits_{j_1,\ldots,j_k=0}^{p}
C_{j_k \ldots j_1}(s)\prod\limits_{l=1}^k \zeta_{j_l}^{(i_l)}
\right)^2\right\}\le 
$$

\vspace{-0.5mm}
\begin{equation}
\label{after3408}
\le K_2(s) \left(\frac{1}{p} +\sum_{r=1}^{\left[k/2\right]}
\sum_{\stackrel{(\{\{g_1, g_2\}, \ldots, 
\{g_{2r-1}, g_{2r}\}\}, \{q_1, \ldots, q_{k-2r}\})}
{{}_{\{g_1, g_2, \ldots, 
g_{2r-1}, g_{2r}, q_1, \ldots, q_{k-2r}\}=\{1, 2, \ldots, k\}}}}
{\sf M}\left\{\left(R_{s,t}^{(p)r, g_1,g_2,\ldots,g_{2r-1}, g_{2r}}\right)^2\right\}\right)
\end{equation}

\vspace{5mm}
\noindent
hold, where $s\in(t,T]$ $(s$ is fixed$),$ $i_1,\ldots,i_k=1,\ldots,m,$

\vspace{-1mm}
$$
R_{s,t}^{(p)r, g_1,g_2,\ldots,g_{2r-1}, g_{2r}}=
R_{T,t}^{(p)r, g_1,g_2,\ldots,g_{2r-1}, g_{2r}}\biggl|_{T=s}\biggr.,
$$

\vspace{2mm}
\noindent
$R_{T,t}^{(p)r, g_1,g_2,\ldots,g_{2r-1}, g_{2r}}$ is defined by {\rm (\ref{afterr1}),}
$J^{*}[\psi^{(k)}]_{T,t}^{(i_1\ldots i_k)}$ and $J^{*}[\psi^{(k)}]_{s,t}^{(i_1\ldots i_k)}$
are iterated Stratonovich stochastic integrals {\rm (\ref{afterstr})} and {\rm (\ref{afterstr1}),}
$C_{j_k \ldots j_1}$ and $C_{j_k \ldots j_1}(s)$
are Fourier coefficients {\rm (\ref{after3000})} and {\rm (\ref{after1300}),} 
constants $K_1,$ $K_2(s)$ are independent of $p;$
another notations are the same as in Theorems {\rm 3,} 
{\rm 13,} {\rm 15.}
}

\vspace{2mm}

{\bf Proof.}\ As follows from Sect.~7 and 8, Conditions 1 and 2 of Theorems
13, 15 are satisfied under the conditions of Theorem~19.
Then from the proof of Theorem~13 it follows that
the expression (\ref{afteru11}) before passing to limit 
\hbox{\vtop{\offinterlineskip\halign{
\hfil#\hfil\cr
{\rm l.i.m.}\cr
$\stackrel{}{{}_{p\to \infty}}$\cr
}} } has the form

\vspace{-1mm}
$$
\sum_{j_1,\ldots,j_k=0}^{p}
C_{j_k\ldots j_1}
\prod\limits_{l=1}^k \zeta_{j_l}^{(i_l)}=
J[\psi^{(k)}]_{T,t}^{(i_1\ldots i_k)p}+
$$

\vspace{1.5mm}
$$
+
\sum_{r=1}^{\left[k/2\right]}\Biggl(\frac{1}{2^r}
\sum_{(s_r,\ldots,s_1)\in {\rm A}_{k,r}}
I[\psi^{(k)}]_{T,t}^{(i_1\ldots i_{s_1-1}i_{s_1+2} \ldots i_{s_r-1}i_{s_r+2}
\ldots i_k)p}
+\Biggr.
$$

\vspace{-1.5mm}
\begin{equation}
\label{after3400}
\Biggl.
+\sum_{\stackrel{(\{\{g_1, g_2\}, \ldots, 
\{g_{2r-1}, g_{2r}\}\}, \{q_1, \ldots, q_{k-2r}\})}
{{}_{\{g_1, g_2, \ldots, 
g_{2r-1}, g_{2r}, q_1, \ldots, q_{k-2r}\}=\{1, 2, \ldots, k\}}}}
R_{T,t}^{(p)r, g_1,g_2,\ldots,g_{2r-1}, g_{2r}}
\Biggr),
\end{equation}

\vspace{5mm}
\noindent                    
where 
$J[\psi^{(k)}]_{T,t}^{(i_1\ldots i_k)p}$
is the approximation for the iterated Ito
stochastic integral (\ref{ito}), which is obtained using Theorem~4, i.e.

$$
J[\psi^{(k)}]_{T,t}^{(i_1\ldots i_k)p}=
\sum\limits_{j_1,\ldots,j_k=0}^{p}
C_{j_k\ldots j_1}\Biggl(
\prod_{l=1}^k\zeta_{j_l}^{(i_l)}+\sum\limits_{r=1}^{[k/2]}
(-1)^r \times
\Biggr.
$$

\vspace{2mm}
\begin{equation}
\label{kkohh}
\times
\sum_{\stackrel{(\{\{g_1, g_2\}, \ldots, 
\{g_{2r-1}, g_{2r}\}\}, \{q_1, \ldots, q_{k-2r}\})}
{{}_{\{g_1, g_2, \ldots, 
g_{2r-1}, g_{2r}, q_1, \ldots, q_{k-2r}\}=\{1, 2, \ldots, k\}}}}
\prod\limits_{s=1}^r
{\bf 1}_{\{i_{g_{{}_{2s-1}}}=~i_{g_{{}_{2s}}}\ne 0\}}
\Biggl.{\bf 1}_{\{j_{g_{{}_{2s-1}}}=~j_{g_{{}_{2s}}}\}}
\prod_{l=1}^{k-2r}\zeta_{j_{q_l}}^{(i_{q_l})}\Biggr),
\end{equation}

\vspace{4mm}
\noindent
$I[\psi^{(k)}]_{T,t}^{(i_1\ldots i_{s_1-1}i_{s_1+2} \ldots i_{s_r-1}i_{s_r+2} \ldots i_k)p}$
is the approximation obtained using (\ref{kkohh}) for the 
iterated Ito
stochastic integral 
$J[\psi^{(k)}]_{T,t}^{s_r, \ldots, s_1}$
(see (\ref{afterito1})).

Using (\ref{after3400}) and Theorem~5, we have

\vspace{1mm}
$$
\sum_{j_1,\ldots,j_k=0}^{p}
C_{j_k\ldots j_1}
\prod\limits_{l=1}^k \zeta_{j_l}^{(i_l)}=
J[\psi^{(k)}]_{T,t}^{(i_1\ldots i_k)}+
\sum_{r=1}^{\left[k/2\right]}\frac{1}{2^r}
\sum_{(s_r,\ldots,s_1)\in {\rm A}_{k,r}}
I[\psi^{(k)}]_{T,t}^{(i_1\ldots i_{s_1-1}i_{s_1+2} \ldots i_{s_r-1}i_{s_r+2}\ldots i_k)}
+
$$

\vspace{3mm}
$$
+\biggl(J[\psi^{(k)}]_{T,t}^{(i_1\ldots i_k)p}-J[\psi^{(k)}]_{T,t}^{(i_1\ldots i_k)}\biggr)+
$$

\vspace{1mm}
$$
+
\sum_{r=1}^{\left[k/2\right]}\sum_{(s_r,\ldots,s_1)\in {\rm A}_{k,r}}
\frac{1}{2^r}
\Biggl(
I[\psi^{(k)}]_{T,t}^{(i_1\ldots i_{s_1-1}i_{s_1+2} \ldots i_{s_r-1}i_{s_r+2}\ldots i_k)p}-
I[\psi^{(k)}]_{T,t}^{(i_1\ldots i_{s_1-1}i_{s_1+2} \ldots i_{s_r-1}i_{s_r+2}\ldots i_k)}
\Biggr)+
$$

\vspace{3mm}
$$
+\sum_{r=1}^{\left[k/2\right]}\sum_{\stackrel{(\{\{g_1, g_2\}, \ldots, 
\{g_{2r-1}, g_{2r}\}\}, \{q_1, \ldots, q_{k-2r}\})}
{{}_{\{g_1, g_2, \ldots, 
g_{2r-1}, g_{2r}, q_1, \ldots, q_{k-2r}\}=\{1, 2, \ldots, k\}}}}
R_{T,t}^{(p)r, g_1,g_2,\ldots,g_{2r-1}, g_{2r}}
=
$$

\vspace{6mm}
$$
=J^{*}[\psi^{(k)}]_{T,t}^{(i_1\ldots i_k)}+
\biggl(J[\psi^{(k)}]_{T,t}^{(i_1\ldots i_k)p}-J[\psi^{(k)}]_{T,t}^{(i_1\ldots i_k)}\biggr)+
$$

\vspace{3mm}
$$
+
\sum_{r=1}^{\left[k/2\right]}\sum_{(s_r,\ldots,s_1)\in {\rm A}_{k,r}}
\frac{1}{2^r}
\Biggl(
I[\psi^{(k)}]_{T,t}^{(i_1\ldots i_{s_1-1}i_{s_1+2} \ldots i_{s_r-1}i_{s_r+2}\ldots i_k)p}-
I[\psi^{(k)}]_{T,t}^{(i_1\ldots i_{s_1-1}i_{s_1+2} \ldots i_{s_r-1}i_{s_r+2}\ldots i_k)}
\Biggr)+
$$

\vspace{3mm}
\begin{equation}
\label{after3401}
+\sum_{r=1}^{\left[k/2\right]}\sum_{\stackrel{(\{\{g_1, g_2\}, \ldots, 
\{g_{2r-1}, g_{2r}\}\}, \{q_1, \ldots, q_{k-2r}\})}
{{}_{\{g_1, g_2, \ldots, 
g_{2r-1}, g_{2r}, q_1, \ldots, q_{k-2r}\}=\{1, 2, \ldots, k\}}}}
R_{T,t}^{(p)r, g_1,g_2,\ldots,g_{2r-1}, g_{2r}}
\end{equation}

\vspace{6mm}
\noindent
w.~p.~1, where we denote $J[\psi^{(k)}]_{T,t}^{s_r, \ldots, s_1}$ as 
$I[\psi^{(k)}]_{T,t}^{(i_1\ldots i_{s_1-1}i_{s_1+2} \ldots i_{s_r-1}i_{s_r+2}\ldots i_k)}$.

In \cite{2018a} (Sect.~1.7.2, Remark~1.7) it is shown that under the conditions of Theorem~19
the following estimate

\vspace{-2mm}
\begin{equation}
\label{zsel1}
{\sf M}\left\{\left(
J[\psi^{(k)}]_{T,t}^{(i_1\ldots i_k)}-J[\psi^{(k)}]_{T,t}^{(i_1\ldots i_k)p}
\right)^2\right\}\le \frac{C}{p},
\end{equation}

\vspace{4mm}
\noindent
holds, where $J[\psi^{(k)}]_{T,t}^{(i_1\ldots i_k)}$ is defined by (\ref{ito}),
$J[\psi^{(k)}]_{T,t}^{(i_1\ldots i_k)p}$ has the form (\ref{kkohh}),
$i_1,\ldots,i_k=0,1,\ldots,m,$
constant $C$ depends only on $k$ and $T-t.$

Applying (\ref{zsel1}), we obtain the following estimates

\vspace{-1mm}
\begin{equation}
\label{after3404}
{\sf M}\left\{\biggl(
J[\psi^{(k)}]_{T,t}^{(i_1\ldots i_k)p}-J[\psi^{(k)}]_{T,t}^{(i_1\ldots i_k)}\biggr)^2\right\}
\le \frac{C}{p},
\end{equation}

\vspace{2mm}
$$
{\sf M}\left\{\Biggl(
I[\psi^{(k)}]_{T,t}^{(i_1\ldots i_{s_1-1}i_{s_1+2} \ldots i_{s_r-1}i_{s_r+2}\ldots i_k)p}
-I[\psi^{(k)}]_{T,t}^{(i_1\ldots i_{s_1-1}i_{s_1+2} \ldots i_{s_r-1}i_{s_r+2}\ldots i_k)}
\Biggr)^2\right\}\le 
$$

\begin{equation}
\label{after3405}
\le\frac{C}{p},
\end{equation}

\vspace{3mm}
\noindent
where constant $C$ does not depend on $p.$

From (\ref{after3401})--(\ref{after3405}) and
the elementary inequality

$$
\left(a_1+a_2+\ldots+a_n\right)^2 \le
n\left(a_1^2+a_2^2+\ldots+a_n^2\right),\ \ \ n\in \mathbb{N}
$$

\vspace{3mm}
\noindent
we obtain (\ref{after3407}).

The estimate 
(\ref{after3408}) is obtained similarly to the estimate 
(\ref{after3407})
using Theorems~1.11, 1.24 in \cite{2018a}, Theorem~15 and 
the estimate \cite{2018a} (Sect.~1.8.1, Remark~1.12)

\vspace{1mm}
$$
{\sf M}\left\{\left(
J[\psi^{(k)}]_{s,t}^{(i_1\ldots i_k)}-J[\psi^{(k)}]_{s,t}^{(i_1\ldots i_k)p}
\right)^2\right\}
\le \frac{C}{p},
$$

\vspace{4mm}
\noindent
where 

\vspace{-5mm}
$$
J[\psi^{(k)}]_{s,t}^{(i_1\ldots i_k)}=
\int\limits_t^s\psi_k(t_k) \ldots \int\limits_t^{t_{2}}
\psi_1(t_1) d{\bf f}_{t_1}^{(i_1)}\ldots
d{\bf f}_{t_k}^{(i_k)},
$$

\vspace{4mm}
$$
J[\psi^{(k)}]_{s,t}^{(i_1\ldots i_k)p}=
\sum\limits_{j_1,\ldots,j_k=0}^{p}
C_{j_k\ldots j_1}(s)\Biggl(
\prod_{l=1}^k\zeta_{j_l}^{(i_l)}+\sum\limits_{r=1}^{[k/2]}
(-1)^r \times
\Biggr.
$$

\vspace{4mm}
$$
\times
\sum_{\stackrel{(\{\{g_1, g_2\}, \ldots, 
\{g_{2r-1}, g_{2r}\}\}, \{q_1, \ldots, q_{k-2r}\})}
{{}_{\{g_1, g_2, \ldots, 
g_{2r-1}, g_{2r}, q_1, \ldots, q_{k-2r}\}=\{1, 2, \ldots, k\}}}}
\prod\limits_{s=1}^r
{\bf 1}_{\{i_{g_{{}_{2s-1}}}=~i_{g_{{}_{2s}}}\ne 0\}}
\Biggl.{\bf 1}_{\{j_{g_{{}_{2s-1}}}=~j_{g_{{}_{2s}}}\}}
\prod_{l=1}^{k-2r}\zeta_{j_{q_l}}^{(i_{q_l})}\Biggr),
$$

\vspace{6mm}

\noindent
where $s\in(t,T]$ $(s$ is fixed$),$\ $C_{j_k \ldots j_1}(s)$
is the Fourier coefficient (\ref{after1300}),
$i_1,\ldots,i_k=0,1,\ldots,m,$ constant $C$ depends only on $k$ and $s-t;$
another notations are the same as in Theorem 4, 15.
Theorem~19 is proved.

\vspace{5mm}

\section{Rate of the Mean-Square Convergence of Expansions of Iterated
Stra\-to\-no\-vich Stochastic Integrals of Multiplicities 3--5
in Theorems 16--18}

\vspace{5mm}

In this section, we consider the rate of convergence of 
approximations of iterated Stratonovich stochastic integrals in Theorems~16--18.
It is easy to see that in Theorems~16--18
the second term 
in parentheses
on the right-hand side of (\ref{after3407}) is estimated.
Combining these results with Theorem~19, we obtain the following theorems.

\vspace{2mm}

{\bf Theorem 20}\ \cite{2018a}, \cite{arxiv-5}, \cite{arxiv-9}, \cite{new-art-1-xxy}.\  
{\it Suppose 
that $\{\phi_j(x)\}_{j=0}^{\infty}$ is a complete orthonormal system of 
Legendre polynomials or trigonometric func\-ti\-ons in the space $L_2([t, T]).$
Furthermore, let $\psi_1(\tau), \psi_2(\tau),$ $\psi_3(\tau)$ are continuously dif\-ferentiable 
nonrandom functions on $[t, T].$ 
Then, for the 
iterated Stra\-to\-no\-vich stochastic integral of third multiplicity

$$
J^{*}[\psi^{(3)}]_{T,t}={\int\limits_t^{*}}^T\psi_3(t_3)
{\int\limits_t^{*}}^{t_3}\psi_2(t_2)
{\int\limits_t^{*}}^{t_2}\psi_1(t_1)
d{\bf f}_{t_1}^{(i_1)}
d{\bf f}_{t_2}^{(i_2)}d{\bf f}_{t_3}^{(i_3)}
$$

\vspace{3mm}
\noindent
the following 
estimate

\vspace{-1mm}
$$
{\sf M}\left\{\left(
J^{*}[\psi^{(3)}]_{T,t}-
\sum\limits_{j_1, j_2, j_3=0}^{p}
C_{j_3 j_2 j_1}\zeta_{j_1}^{(i_1)}\zeta_{j_2}^{(i_2)}\zeta_{j_3}^{(i_3)}\right)^2\right\}
\le \frac{C}{p}
$$

\vspace{3mm}
\noindent
is fulfilled, where $i_1, i_2, i_3=1,\ldots,m,$ 
constant $C$ is independent of $p,$

\vspace{-1mm}
$$
C_{j_3 j_2 j_1}=\int\limits_t^T\psi_3(t_3)\phi_{j_3}(t_3)
\int\limits_t^{t_3}\psi_2(t_2)\phi_{j_2}(t_2)
\int\limits_t^{t_2}\psi_1(t_1)\phi_{j_1}(t_1)dt_1dt_2dt_3
$$

\vspace{3mm}
\noindent
and
$$
\zeta_{j}^{(i)}=
\int\limits_t^T \phi_{j}(s) d{\bf f}_s^{(i)}
$$ 

\vspace{2mm}
\noindent
are independent standard Gaussian random variables for various 
$i$ or $j$.}

\vspace{2mm}

{\bf Theorem 21}\ \cite{2018a}, \cite{arxiv-5}, \cite{arxiv-9}, \cite{new-art-1-xxy}.\  
{\it Let
$\{\phi_j(x)\}_{j=0}^{\infty}$ be a complete orthonormal system of 
Legendre polynomials or trigonometric func\-ti\-ons in the space $L_2([t, T]).$
Furthermore, let $\psi_1(\tau), \ldots,$ $\psi_4(\tau)$ be continuously dif\-ferentiable 
nonrandom functions on $[t, T].$ 
Then, for the 
iterated Stra\-to\-no\-vich stochastic integral of fourth multiplicity

$$
J^{*}[\psi^{(4)}]_{T,t}={\int\limits_t^{*}}^T\psi_4(t_4)
{\int\limits_t^{*}}^{t_4}\psi_3(t_3)
{\int\limits_t^{*}}^{t_3}\psi_2(t_2)
{\int\limits_t^{*}}^{t_2}\psi_1(t_1)
d{\bf f}_{t_1}^{(i_1)}
d{\bf f}_{t_2}^{(i_2)}d{\bf f}_{t_3}^{(i_3)}d{\bf f}_{t_4}^{(i_4)}
$$

\vspace{3mm}
\noindent
the following 
estimate

\vspace{-1mm}
$$
{\sf M}\left\{\left(
J^{*}[\psi^{(4)}]_{T,t}-
\sum\limits_{j_1, j_2, j_3, j_4=0}^{p}
C_{j_4 j_3 j_2 j_1}\zeta_{j_1}^{(i_1)}\zeta_{j_2}^{(i_2)}\zeta_{j_3}^{(i_3)}
\zeta_{j_4}^{(i_4)}
\right)^2\right\}
\le \frac{C}{p^{1-\varepsilon}}
$$

\vspace{3mm}
\noindent
holds, where $i_1, i_2, i_3, i_4=1,\ldots,m,$ 
constant $C$ does not depend on $p,$
$\varepsilon$ is an arbitrary
small positive real number 
for the case of complete orthonormal system of 
Legendre polynomials in the space $L_2([t, T])$
and $\varepsilon=0$ for the case of
complete orthonormal system of 
trigonometric functions in the space $L_2([t, T]),$

\vspace{-1mm}
$$
C_{j_4 j_3 j_2 j_1}=
\int\limits_t^T\psi_4(t_4)\phi_{j_4}(t_4)
\int\limits_t^{t_4}\psi_3(t_3)\phi_{j_3}(t_3)
\int\limits_t^{t_3}\psi_2(t_2)\phi_{j_2}(t_2)
\int\limits_t^{t_2}\psi_1(t_1)\phi_{j_1}(t_1)dt_1\times
$$

$$
\times dt_2dt_3dt_4;
$$

\vspace{4mm}
\noindent
another notations are the same as in Theorem~{\rm 20}.}

\vspace{2mm}

{\bf Theorem 22}\ \cite{2018a}, \cite{arxiv-5}, \cite{arxiv-9}, \cite{new-art-1-xxy}.\  
{\it Assume 
that $\{\phi_j(x)\}_{j=0}^{\infty}$ is a complete orthonormal system of 
Legendre polynomials or trigonometric func\-ti\-ons in the space $L_2([t, T])$
and $\psi_1(\tau), \ldots,$ $\psi_5(\tau)$ are continuously dif\-ferentiable 
nonrandom functions on $[t, T].$ 
Then, for the 
iterated Stra\-to\-no\-vich stochastic integral of fifth multiplicity

$$
J^{*}[\psi^{(5)}]_{T,t}={\int\limits_t^{*}}^T\psi_5(t_5)
\ldots
{\int\limits_t^{*}}^{t_2}\psi_1(t_1)
d{\bf f}_{t_1}^{(i_1)}
\ldots d{\bf f}_{t_5}^{(i_5)}
$$

\vspace{3mm}
\noindent
the following 
estimate

\vspace{-1mm}
$$
{\sf M}\left\{\left(
J^{*}[\psi^{(5)}]_{T,t}-
\sum\limits_{j_1, \ldots, j_5=0}^{p}
C_{j_5 \ldots j_1}\zeta_{j_1}^{(i_1)}\ldots
\zeta_{j_5}^{(i_5)}
\right)^2\right\}
\le \frac{C}{p^{1-\varepsilon}}
$$

\vspace{3mm}
\noindent
is valid, where $i_1, \ldots, i_5=1,\ldots,m,$ 
constant $C$ is independent of $p,$
$\varepsilon$ is an arbitrary
small positive real number 
for the case of complete orthonormal system of 
Legendre polynomials in the space $L_2([t, T])$
and $\varepsilon=0$ for the case of
complete orthonormal system of 
trigonometric functions in the space $L_2([t, T]),$

\vspace{-1mm}
$$
C_{j_5 \ldots j_1}=
\int\limits_t^T\psi_5(t_5)\phi_{j_5}(t_5)\ldots
\int\limits_t^{t_2}\psi_1(t_1)\phi_{j_1}(t_1)dt_1\ldots dt_5;
$$

\vspace{3mm}
\noindent
another notations are the same as in Theorem~{\rm 20, 21}.}

\vspace{5mm}

\section{Expansion of Iterated Stratonovich Stochastic Integrals
of Multiplicity 6. The Case $p_1=\ldots =p_6\to \infty$ and 
$\psi_1(\tau),$ $\ldots,$ $\psi_6(\tau)\equiv 1$ 
(The Cases of Legendre 
Polynomials and Trigonometric Functions)}

\vspace{5mm}   
     
{\bf Theorem 23}\ \cite{2018a}, \cite{arxiv-5}, \cite{arxiv-9}, \cite{new-art-1xxys}.\
{\it Suppose that 
$\{\phi_j(x)\}_{j=0}^{\infty}$ is a complete orthonormal system of 
Legendre polynomials or trigonometric func\-ti\-ons in the space $L_2([t, T]).$
Then, for the 
iterated Stratonovich stochastic integral of sixth multiplicity

\begin{equation}
\label{after10001qu1}
J_{T,t}^{*(i_1\ldots i_6)}={\int\limits_t^{*}}^T
\ldots
{\int\limits_t^{*}}^{t_2}
d{\bf w}_{t_1}^{(i_1)}
\ldots d{\bf w}_{t_6}^{(i_6)}
\end{equation}

\vspace{3mm}
\noindent
the following 
expansion 

\vspace{-1mm}
$$
J_{T,t}^{*(i_1\ldots i_6)}
=\hbox{\vtop{\offinterlineskip\halign{
\hfil#\hfil\cr
{\rm l.i.m.}\cr
$\stackrel{}{{}_{p\to \infty}}$\cr
}} }
\sum\limits_{j_1, \ldots, j_6=0}^{p}
C_{j_6 \ldots j_1}\zeta_{j_1}^{(i_1)}\ldots
\zeta_{j_6}^{(i_6)}
$$

\vspace{4mm}
\noindent
that converges in the mean-square sense is valid, where
$i_1, \ldots, i_6=0, 1,\ldots,m,$

$$
C_{j_6 \ldots j_1}=
\int\limits_t^T\phi_{j_6}(t_6)\ldots
\int\limits_t^{t_2}\phi_{j_1}(t_1)dt_1\ldots dt_6
$$

\vspace{3mm}
\noindent
and
$$
\zeta_{j}^{(i)}=
\int\limits_t^T \phi_{j}(s) d{\bf w}_s^{(i)}
$$ 

\vspace{2mm}
\noindent
are independent standard Gaussian random variables for various 
$i$ or $j$
{\rm (}in the case when $i\ne 0${\rm ),}
${\bf w}_{\tau}^{(i)}={\bf f}_{\tau}^{(i)}$ for
$i=1,\ldots,m$ and 
${\bf w}_{\tau}^{(0)}=\tau.$}

\vspace{2mm}

{\bf Proof.}\ As noted in Sect.~7, Conditions 1 and 2
of Theorem~{\rm 13} are satisfied for complete
orthonormal systems of Legendre polynomials 
and trigonometric functions in the space
$L_2([t, T]).$ Let us verify Condition 3 of Theorem~13 for
the iterated Stratonovich stochastic integral (\ref{after10001qu1}). 
Thus, we have to check the following conditions

\vspace{-1mm}
\begin{equation}
\label{after14000qu1}
\lim\limits_{p\to\infty}
\sum\limits_{j_{q_1},j_{q_2},j_{q_3},j_{q_4}=0}^p
\left(\sum_{j_{g_1}=p+1}^{\infty}
C_{j_6\ldots j_1}\biggl|_{j_{g_1}=j_{g_2}}\biggr.\right)^2=0,
\end{equation}

\vspace{2mm}
\begin{equation}
\label{after14001qu2}
\lim\limits_{p\to\infty}
\sum\limits_{j_{q_1},j_{q_2}=0}^p
\left(\sum_{j_{g_1}=p+1}^{\infty}\sum_{j_{g_3}=p+1}^{\infty}
C_{j_6\ldots j_1}\biggl|_{j_{g_1}=j_{g_2},j_{g_3}=j_{g_4}}\biggr.\right)^2=0,
\end{equation}

\vspace{2mm}
\begin{equation}
\label{afterafter001qu3}
\lim\limits_{p\to\infty}
\sum\limits_{j_{q_1},j_{q_2}=0}^p
\left(\sum_{j_{g_1}=p+1}^{\infty}
C_{j_6\ldots j_1}\biggl|_{(j_{g_4}j_{g_3})\curvearrowright (\cdot),
j_{g_1}=j_{g_2},
j_{g_3}=j_{g_4}, g_4=g_3+1}\biggr.\right)^2=0,
\end{equation}

\vspace{2mm}
\begin{equation}
\label{after14001qu10}
\lim\limits_{p\to\infty}
\left(\sum_{j_{g_1}=p+1}^{\infty}\sum_{j_{g_3}=p+1}^{\infty}
\sum_{j_{g_5}=p+1}^{\infty}
C_{j_6\ldots j_1}\biggl|_{j_{g_1}=j_{g_2},j_{g_3}=j_{g_4},j_{g_5}=j_{g_6}}\biggr.\right)^2=0,
\end{equation}

\vspace{2mm}
\begin{equation}
\label{afterafter001qu33}
\lim\limits_{p\to\infty}
\left(\sum_{j_{g_1}=p+1}^{\infty}\sum_{j_{g_3}=p+1}^{\infty}
C_{j_6\ldots j_1}\biggl|_{(j_{g_6}j_{g_5})\curvearrowright (\cdot),
j_{g_1}=j_{g_2},
j_{g_3}=j_{g_4}, j_{g_5}=j_{g_6}, g_6=g_5+1}\biggr.\right)^2=0,
\end{equation}

\vspace{2mm}
\begin{equation}
\label{afterafter001qu36}
\lim\limits_{p\to\infty}
\left(\sum_{j_{g_1}=p+1}^{\infty}
C_{j_6\ldots j_1}\biggl|_{(j_{g_4}j_{g_3})\curvearrowright (\cdot)
(j_{g_6}j_{g_5})\curvearrowright (\cdot),
j_{g_1}=j_{g_2},
j_{g_3}=j_{g_4}, j_{g_5}=j_{g_6}, g_4=g_3+1, g_6=g_5+1}\biggr.\right)^2=0,
\end{equation}

\vspace{6mm}
\noindent
where the expressions

\vspace{-2mm}
$$
\left(\{g_1,g_2\},\{g_3,g_4\}, \{g_5,g_6\}\}\right),\ \ \
\left(\{g_1,g_2\},\{g_3,g_4\}, \{q_1,q_2\}\}\right),\ \ \
\left(\{g_1,g_2\}, \{q_1, q_2,q_3, q_4\}\right)
$$

\vspace{3mm}
\noindent
are partitions of the set $\{1,2,\ldots,6\}$ that is
$\{g_1,g_2,g_3,$ $g_4,g_5,g_6\}=\{g_1,g_2,g_3,g_4,q_1,q_2\}=
\{g_1,g_2,q_1,$ $q_2,q_3,q_4\}=
\{1,2,$ $\ldots,6\};$
braces mean an unordered 
set, and pa\-ren\-the\-ses mean an ordered set.

The equalities (\ref{after14000qu1}),
(\ref{afterafter001qu3}) were proved earlier (see the proof of equalities 
(\ref{after17000}), (\ref{after19005})).
The relation (\ref{afterafter001qu36}) follows 
from
the estimate (\ref{5tafter}) for the polynomial case and its analogue
for the
trigonometric case.
It is easy to see that the 
equalities (\ref{after14001qu2}) and 
(\ref{afterafter001qu33}) are proved in complete analogy with the proof of 
(\ref{after14001}), (\ref{after19005}).

Thus, we have to prove the relation (\ref{after14001qu10}).
The equality (\ref{after14001qu10}) is equivalent to the following equalities
\begin{equation}
\label{sixsix8}
\lim\limits_{p\to\infty}
\sum_{j_1=p+1}^{\infty}\sum_{j_2=p+1}^{\infty}
\sum_{j_3=p+1}^{\infty}
C_{j_3 j_2 j_1 j_3 j_2 j_1}=0,
\end{equation}

\begin{equation}
\label{sixsix9}
\lim\limits_{p\to\infty}
\sum_{j_1=p+1}^{\infty}\sum_{j_2=p+1}^{\infty}
\sum_{j_3=p+1}^{\infty}
C_{j_1 j_3 j_2 j_3 j_2 j_1}=0,
\end{equation}

\begin{equation}
\label{sixsix10}
\lim\limits_{p\to\infty}
\sum_{j_1=p+1}^{\infty}\sum_{j_2=p+1}^{\infty}
\sum_{j_3=p+1}^{\infty}
C_{j_3 j_2 j_3 j_1 j_2 j_1}=0,
\end{equation}

\begin{equation}
\label{sixsix4}
\lim\limits_{p\to\infty}
\sum_{j_1=p+1}^{\infty}\sum_{j_2=p+1}^{\infty}
\sum_{j_3=p+1}^{\infty}
C_{j_1 j_2 j_3 j_3 j_2 j_1}=0,
\end{equation}

\begin{equation}
\label{sixsix14}
\lim\limits_{p\to\infty}
\sum_{j_1=p+1}^{\infty}\sum_{j_2=p+1}^{\infty}
\sum_{j_3=p+1}^{\infty}
C_{j_1 j_2 j_2 j_3 j_3 j_1}=0,
\end{equation}

\begin{equation}
\label{sixsix3}
\lim\limits_{p\to\infty}
\sum_{j_1=p+1}^{\infty}\sum_{j_2=p+1}^{\infty}
\sum_{j_3=p+1}^{\infty}
C_{j_3 j_3 j_2 j_2 j_1 j_1}=0,
\end{equation}

\begin{equation}
\label{sixsix7}
\lim\limits_{p\to\infty}
\sum_{j_1=p+1}^{\infty}\sum_{j_2=p+1}^{\infty}
\sum_{j_3=p+1}^{\infty}
C_{j_2 j_3 j_3 j_2 j_1 j_1}=0,
\end{equation}

\begin{equation}
\label{sixsix6}
\lim\limits_{p\to\infty}
\sum_{j_1=p+1}^{\infty}\sum_{j_2=p+1}^{\infty}
\sum_{j_3=p+1}^{\infty}
C_{j_3 j_2 j_3 j_2 j_1 j_1}=0,
\end{equation}

\begin{equation}
\label{sixsix1}
\lim\limits_{p\to\infty}
\sum_{j_1=p+1}^{\infty}\sum_{j_2=p+1}^{\infty}
\sum_{j_3=p+1}^{\infty}
C_{j_3 j_3 j_2 j_1 j_2 j_1}=0,
\end{equation}

\begin{equation}
\label{sixsix2}
\lim\limits_{p\to\infty}
\sum_{j_1=p+1}^{\infty}\sum_{j_2=p+1}^{\infty}
\sum_{j_3=p+1}^{\infty}
C_{j_3 j_3 j_1 j_2 j_2 j_1}=0,
\end{equation}

\begin{equation}
\label{sixsix5}
\lim\limits_{p\to\infty}
\sum_{j_1=p+1}^{\infty}\sum_{j_2=p+1}^{\infty}
\sum_{j_3=p+1}^{\infty}
C_{j_2 j_1 j_3 j_3 j_2 j_1}=0,
\end{equation}

\begin{equation}
\label{sixsix12}
\lim\limits_{p\to\infty}
\sum_{j_1=p+1}^{\infty}\sum_{j_2=p+1}^{\infty}
\sum_{j_3=p+1}^{\infty}
C_{j_3 j_1 j_2 j_3 j_2 j_1}=0,
\end{equation}

\begin{equation}
\label{sixsix11}
\lim\limits_{p\to\infty}
\sum_{j_1=p+1}^{\infty}\sum_{j_2=p+1}^{\infty}
\sum_{j_3=p+1}^{\infty}
C_{j_2 j_3 j_1 j_3 j_2 j_1}=0,
\end{equation}

\begin{equation}
\label{sixsix13}
\lim\limits_{p\to\infty}
\sum_{j_1=p+1}^{\infty}\sum_{j_2=p+1}^{\infty}
\sum_{j_3=p+1}^{\infty}
C_{j_3 j_1 j_3 j_2 j_2 j_1}=0,
\end{equation}

\begin{equation}
\label{sixsix15}
\lim\limits_{p\to\infty}
\sum_{j_1=p+1}^{\infty}\sum_{j_2=p+1}^{\infty}
\sum_{j_3=p+1}^{\infty}
C_{j_2 j_3 j_3 j_1 j_2 j_1}=0.
\end{equation}

\vspace{5mm}

Consider in detail the case of Legendre polynomials 
(the case of trigonometric functions is considered in complete analogy).

First, we prove the following equality for the Fourier coefficients for the case 
$\psi_1(\tau),\ldots,\psi_6(\tau)\equiv 1$

$$
C_{j_6 j_5 j_4 j_3 j_2 j_1}+C_{j_1 j_2 j_3 j_4 j_5 j_6}=
C_{j_6}C_{j_5 j_4 j_3 j_2 j_1}-C_{j_5 j_6}C_{j_4 j_3 j_2 j_1}+
$$

\begin{equation}
\label{sixsix40}
+C_{j_4 j_5 j_6}C_{j_3 j_2 j_1}-C_{j_3 j_4 j_5 j_6}C_{j_2 j_1}+
C_{j_2 j_3 j_4 j_5 j_6}C_{j_1}.
\end{equation}

\vspace{5mm}

Using the integration order replacement, we have

$$
C_{j_6 j_5 j_4 j_3 j_2 j_1}=
$$

\vspace{2mm}
$$
=\int\limits_t^T\phi_{j_6}(t_6)\int\limits_t^{t_6}\phi_{j_5}(t_5)
\ldots
\int\limits_t^{t_2}\phi_{j_1}(t_1)dt_1 \ldots dt_5 dt_6=
$$

\vspace{2mm}
$$
=\int\limits_t^T\phi_{j_6}(t_6)\int\limits_t^{T}\phi_{j_5}(t_5)
\int\limits_t^{t_5}\phi_{j_4}(t_4)
\ldots \int\limits_t^{t_2}\phi_{j_1}(t_1)dt_1 \ldots  dt_4 dt_5 dt_6-
$$

\vspace{2mm}
$$
-\int\limits_t^T\phi_{j_6}(t_6)\int\limits_{t_6}^T\phi_{j_5}(t_5)
\int\limits_t^{t_5}\phi_{j_4}(t_4)
\ldots 
\int\limits_t^{t_2}\phi_{j_1}(t_1)dt_1 \ldots dt_4 dt_5 dt_6=
$$

\vspace{2mm}
$$
=C_{j_6}C_{j_5 j_4 j_3 j_2 j_1}-
$$

\vspace{2mm}
$$
-\int\limits_t^T\phi_{j_6}(t_6)\int\limits_{t_6}^T\phi_{j_5}(t_5)
\int\limits_t^{T}\phi_{j_4}(t_4)
\int\limits_t^{t_4}\phi_{j_3}(t_3)
\ldots 
\int\limits_t^{t_2}\phi_{j_1}(t_1)dt_1 \ldots dt_3 dt_4 dt_5 dt_6+
$$

\vspace{2mm}
$$
+\int\limits_t^T\phi_{j_6}(t_6)\int\limits_{t_6}^T\phi_{j_5}(t_5)
\int\limits_{t_5}^{T}\phi_{j_4}(t_4)\int\limits_{t}^{t_4}\phi_{j_3}(t_3)
\ldots 
\int\limits_t^{t_2}\phi_{j_1}(t_1)dt_1 \ldots dt_3 dt_4 dt_5 dt_6=
$$

\vspace{2mm}
$$
=C_{j_6}C_{j_5 j_4 j_3 j_2 j_1}-
$$

\vspace{2mm}
$$
-\int\limits_t^T\phi_{j_6}(t_6)\int\limits_{t_6}^T\phi_{j_5}(t_5)dt_5 dt_6\
C_{j_4 j_3 j_2 j_1}+
$$

\vspace{2mm}
$$
+\int\limits_t^T\phi_{j_6}(t_6)\int\limits_{t_6}^T\phi_{j_5}(t_5)
\int\limits_{t_5}^{T}\phi_{j_4}(t_4)\int\limits_{t}^{t_4}\phi_{j_3}(t_3)
\ldots 
\int\limits_t^{t_2}\phi_{j_1}(t_1)dt_1 \ldots dt_3 dt_4 dt_5 dt_6=
$$

\vspace{2mm}
$$
=C_{j_6}C_{j_5 j_4 j_3 j_2 j_1}-
C_{j_5 j_6}C_{j_4 j_3 j_2 j_1}+
$$

\vspace{2mm}
$$
+\int\limits_t^T\phi_{j_6}(t_6)\int\limits_{t_6}^T\phi_{j_5}(t_5)
\int\limits_{t_5}^{T}\phi_{j_4}(t_4)\int\limits_{t}^{t_4}\phi_{j_3}(t_3)
\ldots 
\int\limits_t^{t_2}\phi_{j_1}(t_1)dt_1 \ldots dt_3 dt_4 dt_5 dt_6=
$$

$$
\ldots
$$

$$
=
C_{j_6}C_{j_5 j_4 j_3 j_2 j_1}-C_{j_5 j_6}C_{j_4 j_3 j_2 j_1}+
C_{j_4 j_5 j_6}C_{j_3 j_2 j_1}-C_{j_3 j_4 j_5 j_6}C_{j_2 j_1}+
C_{j_2 j_3 j_4 j_5 j_6}C_{j_1}-
$$

\vspace{2mm}
$$
-\int\limits_t^T\phi_{j_6}(t_6)\int\limits_{t_6}^T\phi_{j_5}(t_5)
\ldots
\int\limits_{t_2}^T\phi_{j_1}(t_1)dt_1 \ldots dt_5 dt_6=
$$

\vspace{2mm}
$$
=
C_{j_6}C_{j_5 j_4 j_3 j_2 j_1}-C_{j_5 j_6}C_{j_4 j_3 j_2 j_1}+
C_{j_4 j_5 j_6}C_{j_3 j_2 j_1}-
$$

\begin{equation}
\label{sixsix41}
-C_{j_3 j_4 j_5 j_6}C_{j_2 j_1}+
C_{j_2 j_3 j_4 j_5 j_6}C_{j_1}-C_{j_1 j_2 j_3 j_4 j_5 j_6}.
\end{equation}

\vspace{5mm}

The equality (\ref{sixsix41}) completes the proof of the relation 
(\ref{sixsix40}).

Let us consider 
(\ref{sixsix8}). From (\ref{after80xx}) we obtain

\begin{equation}
\label{sixsix42}
\sum_{j_1=p+1}^{\infty}\sum_{j_2=p+1}^{\infty}
\sum_{j_3=p+1}^{\infty}
C_{j_3 j_2 j_1 j_3 j_2 j_1}=
-\sum_{j_1=0}^{p}\sum_{j_2=0}^{p}
\sum_{j_3=0}^{p}
C_{j_3 j_2 j_1 j_3 j_2 j_1}.
\end{equation}

\vspace{5mm}

Applying (\ref{sixsix40}), we get

$$
\sum_{j_1,j_2,j_3=0}^{p}
C_{j_3 j_2 j_1 j_3 j_2 j_1}+
\sum_{j_1,j_2,j_3=0}^{p}
C_{j_1 j_2 j_3 j_1 j_2 j_3}=
2\sum_{j_1,j_2,j_3=0}^{p}
C_{j_3 j_2 j_1 j_3 j_2 j_1}=
$$

\vspace{2mm}
$$
=
\sum_{j_1,j_2,j_3=0}^{p}\biggl(
C_{j_3}C_{j_2 j_1 j_3 j_2 j_1}-C_{j_2 j_3}C_{j_1 j_3 j_2 j_1}+
C_{j_1 j_2 j_3}C_{j_3 j_2 j_1}-\biggr.
$$

\begin{equation}
\label{sixsix43}
\biggl.-C_{j_3 j_1 j_2 j_3}C_{j_2 j_1}+
C_{j_2 j_3 j_1 j_2 j_3}C_{j_1}\biggr).
\end{equation}

\vspace{5mm}

The complete orthonormal system of Legendre polynomials in the 
space $L_2([t,T])$ looks as follows

\vspace{-1mm}
\begin{equation}
\label{4009d}
\phi_j(x)=\sqrt{\frac{2j+1}{T-t}}P_j\left(\left(
x-\frac{T+t}{2}\right)\frac{2}{T-t}\right),\ \ \ j=0, 1, 2,\ldots,
\end{equation}

\vspace{4mm}
\noindent
where 
$$
P_j(x)=\frac{1}{2^j j!} \frac{d^j}{dx^j}\left(x^2-1\right)^j
$$

\vspace{4mm}
\noindent
is the Legendre polynomial.

Note that
$$
C_{j_2j_1}=\int\limits_t^T\phi_{j_2}(\tau)\int\limits_t^{\tau}\phi_{j_1}(\theta)d\theta d\tau=
$$

\vspace{3mm}
\begin{equation}
\label{sixsix50}
=
\frac{T-t}{2}\left\{
\begin{matrix}
1/\sqrt{(2j_1+1)(2j_1+3)} &\hbox{if}\ j_2=j_1+1,\ j_1=0,1,2,\ldots \cr\cr
-1/\sqrt{4j_1^2-1} &\hbox{if}\ j_2=j_1-1,\ j_1=1,2,\ldots \cr\cr
1 &\hbox{if}\ j_1=j_2=0 \cr\cr
0 &\hbox{otherwise}
\end{matrix},\right.
\end{equation}

\vspace{5mm}
\begin{equation}
\label{zero1s}
C_{j_1}=\int\limits_t^T\phi_{j_1}(\tau)d\tau=
\left\{
\begin{matrix}
\sqrt{T-t} &\hbox{if}\ j_1=0\cr\cr
0 &\hbox{if}\ j_1\ne 0
\end{matrix}.\right.
\end{equation}

\vspace{5mm}

Moreover, the generalized Parseval equality gives

$$
\lim\limits_{p\to\infty}\sum_{j_1,j_2,j_3=0}^{p}
C_{j_1 j_2 j_3}C_{j_3 j_2 j_1}=
$$

\vspace{2mm}
$$
=\lim\limits_{p\to\infty}\sum_{j_1,j_2,j_3=0}^{p}
\int\limits_t^T\phi_{j_1}(t_3)\int\limits_t^{t_3}\phi_{j_2}(t_2)
\int\limits_t^{t_2}\phi_{j_3}(t_1)dt_1dt_2dt_3\times
$$

\vspace{2mm}
$$
\times
\int\limits_t^T\phi_{j_3}(t_3)\int\limits_t^{t_3}\phi_{j_2}(t_2)
\int\limits_t^{t_2}\phi_{j_1}(t_1)dt_1dt_2dt_3=
$$

\vspace{2mm}
$$
=\lim\limits_{p\to\infty}\sum_{j_1,j_2,j_3=0}^{p}
\int\limits_t^T\phi_{j_3}(t_3)\int\limits_{t_3}^T\phi_{j_2}(t_2)
\int\limits_{t_2}^T\phi_{j_1}(t_1)dt_1dt_2dt_3\times
$$

\vspace{2mm}
$$
\times
\int\limits_t^T\phi_{j_3}(t_3)\int\limits_t^{t_3}\phi_{j_2}(t_2)
\int\limits_t^{t_2}\phi_{j_1}(t_1)dt_1dt_2dt_3=
$$

\vspace{2mm}
$$
=\lim\limits_{p\to\infty}\sum_{j_1,j_2,j_3=0}^{p}~
\int\limits_{[t,T]^3}{\bf 1}_{\{t_3<t_2<t_1\}}
\prod\limits_{l=1}^3
\phi_{j_l}(t_l)dt_1dt_2dt_3\times
$$

\vspace{2mm}
$$
\times
\int\limits_{[t,T]^3}{\bf 1}_{\{t_1<t_2<t_3\}}
\prod\limits_{l=1}^3
\phi_{j_l}(t_l)dt_1dt_2dt_3=
$$

\vspace{2mm}
\begin{equation}
\label{pars100}
=\int\limits_{[t,T]^3}{\bf 1}_{\{t_3<t_2<t_1\}}{\bf 1}_{\{t_1<t_2<t_3\}}
dt_1dt_2dt_3=0.
\end{equation}

\vspace{5mm}

Using the above arguments and also (\ref{after80xx}), (\ref{sixsix42}), and (\ref{sixsix43}), we get

$$
-\lim\limits_{p\to\infty}\sum_{j_1=p+1}^{\infty}\sum_{j_2=p+1}^{\infty}
\sum_{j_3=p+1}^{\infty}
C_{j_3 j_2 j_1 j_3 j_2 j_1}=
\lim\limits_{p\to\infty}\sum_{j_1,j_2,j_3=0}^{p}
C_{j_3 j_2 j_1 j_3 j_2 j_1}=
$$

\vspace{2mm}
$$
=
\frac{1}{2}\lim\limits_{p\to\infty}\sum_{j_1,j_2,j_3=0}^{p}\biggl(
C_{j_3}C_{j_2 j_1 j_3 j_2 j_1}-C_{j_2 j_3}C_{j_1 j_3 j_2 j_1}-\biggr.
$$

\vspace{2mm}
$$
\biggl.-C_{j_3 j_1 j_2 j_3}C_{j_2 j_1}+
C_{j_2 j_3 j_1 j_2 j_3}C_{j_1}\biggr)=
$$

\vspace{2mm}
$$
=
\lim\limits_{p\to\infty}\sum_{j_1,j_2,j_3=0}^{p}\biggl(
C_{j_3}C_{j_2 j_1 j_3 j_2 j_1}-C_{j_3 j_1 j_2 j_3}C_{j_2 j_1}\biggr)=
$$

\vspace{2mm}
$$
=
\sqrt{T-t}\lim\limits_{p\to\infty}\sum_{j_1,j_2=0}^{p}
C_{j_2 j_1 0 j_2 j_1}
-
\lim\limits_{p\to\infty}\sum_{j_1,j_2,j_3=0}^{p}
C_{j_3 j_1 j_2 j_3}C_{j_2 j_1}
=
$$

\vspace{2mm}
\begin{equation}
\label{sixsix71}
=\sqrt{T-t}\lim\limits_{p\to\infty}\sum_{j_1,j_2=0}^{p}
C_{j_2 j_1 0 j_2 j_1}+
\lim\limits_{p\to\infty}\sum_{j_1,j_2=0}^p\sum_{j_3=p+1}^{\infty}
C_{j_3 j_1 j_2 j_3}C_{j_2 j_1}.
\end{equation}

\vspace{5mm}

By analogy with the proof of (\ref{after2508}) (see the proof of Theorem~17)
we obtain

\vspace{-1mm}
\begin{equation}
\label{sixsix72}
\lim\limits_{p\to\infty}\sum_{j_1,j_2=0}^{p}
C_{j_2 j_1 0 j_2 j_1}=
\lim\limits_{p\to\infty}\sum_{j_1=p+1}^{\infty}\sum_{j_2=p+1}^{\infty} 
C_{j_2 j_1 0 j_2 j_1}=0,
\end{equation}

\vspace{4mm}
\noindent
where we used the following representation

$$
C_{j_2 j_1 0 j_2 j_1}=
$$

\vspace{2mm}
$$
=\frac{1}{\sqrt{T-t}}
\int\limits_t^T 
\phi_{j_2}(t_5)
\int\limits_t^{t_5} 
\phi_{j_1}(t_4)
\int\limits_t^{t_4} 
\int\limits_t^{t_3} 
\phi_{j_2}(t_2)
\int\limits_t^{t_2} 
\phi_{j_1}(t_1)
dt_1 dt_2 dt_3 dt_4 dt_5=
$$

\vspace{2mm}
$$
=\frac{1}{\sqrt{T-t}}
\int\limits_t^T 
\phi_{j_2}(t_5)
\int\limits_t^{t_5} 
\phi_{j_1}(t_4)
\int\limits_t^{t_4} 
\phi_{j_2}(t_2)
\int\limits_t^{t_2} 
\phi_{j_1}(t_1)
dt_1 
\int\limits_{t_2}^{t_4} 
dt_3 dt_2 dt_4 dt_5=
$$

\vspace{2mm}
$$
=\frac{1}{\sqrt{T-t}}
\int\limits_t^T 
\phi_{j_2}(t_5)
\int\limits_t^{t_5} 
\phi_{j_1}(t_4)(t_4-t)
\int\limits_t^{t_4} 
\phi_{j_2}(t_2)
\int\limits_t^{t_2} 
\phi_{j_1}(t_1)
dt_1 
dt_2 dt_4 dt_5+
$$

\vspace{2mm}
$$
+\frac{1}{\sqrt{T-t}}
\int\limits_t^T 
\phi_{j_2}(t_5)
\int\limits_t^{t_5} 
\phi_{j_1}(t_4)
\int\limits_t^{t_4} 
\phi_{j_2}(t_2)(t-t_2)
\int\limits_t^{t_2} 
\phi_{j_1}(t_1)
dt_1 
dt_2 dt_4 dt_5\stackrel{\sf def}{=}
$$

\vspace{2mm}
$$
\stackrel{\sf def}{=}
\bar C_{j_2 j_1 j_2 j_1}+\tilde C_{j_2 j_1 j_2 j_1}.
$$

\vspace{5mm}

Further, we have (see (\ref{sixsix50}))

$$
\lim\limits_{p\to\infty}\sum_{j_1,j_2=0}^p\sum_{j_3=p+1}^{\infty}
C_{j_3 j_1 j_2 j_3}C_{j_2 j_1}=
\lim\limits_{p\to\infty}\sum_{j_3=p+1}^{\infty}
\biggl(C_{00}C_{j_3 00 j_3}+\biggr.
$$

\vspace{2mm}
\begin{equation}
\label{sixsix70}
\biggl.+\sum\limits_{j_1=1}^p C_{j_1-1,j_1}C_{j_3j_1,j_1-1,j_3}+
\sum\limits_{j_1=1}^{p-1} C_{j_1+1,j_1}C_{j_3j_1,j_1+1,j_3}+
C_{1,0}C_{j_301j_3}
\biggr).
\end{equation}

\vspace{5mm}

Observe that
\begin{equation}
\label{sixsix60}
|C_{j_1-1,j_1}|+|C_{j_1+1,j_1}|\le \frac{K}{j_1}\ \ \ (j_1=1,\ldots,p),
\end{equation}

\begin{equation}
\label{sixsix61}
|C_{j_3 00 j_3}|+|C_{j_3j_1,j_1-1,j_3}|+|C_{j_3j_1,j_1+1,j_3}|+|C_{j_301j_3}|\le
\frac{K_1}{j_3^2}\ \ \ (j_3\ge p+1),
\end{equation}

\vspace{3mm}
\noindent
where constants $K, K_1$ do not depend on $j_1, j_3.$

The estimate (\ref{sixsix60}) follow from (\ref{sixsix50}).
At the same time, the estimate (\ref{sixsix61}) can be obtained using the following reasoning.
First note that the integration order replacement gives

$$
C_{j_3 j_1 j_2 j_3}=\int\limits_t^T\phi_{j_3}(t_4)\int\limits_{t}^{t_4}\phi_{j_1}(t_3)
\int\limits_t^{t_3}\phi_{j_2}(t_2)
\int\limits_t^{t_2}\phi_{j_3}(t_1)
dt_1 dt_2 dt_3 dt_4=
$$

\vspace{2mm}
\begin{equation}
\label{sixsix62}
=\int\limits_t^T\phi_{j_1}(t_3)\int\limits_{t}^{t_3}\phi_{j_2}(t_2)
\left(\int\limits_t^{t_2}\phi_{j_3}(t_1)dt_1\right) dt_2
\left(\int\limits_{t_3}^T\phi_{j_3}(t_4)
dt_4\right) dt_3.
\end{equation}

\vspace{5mm}

Consider the well-known estimate for Legendre
polynomials

\vspace{-1mm}
\begin{equation}
\label{otit987}
|P_j(y)| <\frac{K}{\sqrt{j+1}(1-y^2)^{1/4}},\ \ \ 
y\in (-1, 1),\ \ \ j\in \mathbb{N},
\end{equation}

\vspace{4mm}
\noindent
where constant $K$ does not depend on $y$ and $j$.

The estimate (\ref{otit987}) can be rewritten for the 
function $\phi_j(x)$ (see (\ref{4009d})) in 
the following form

\begin{equation}
\label{ogo24}
|\phi_j(x)|< \sqrt{\frac{2j+1}{j+1}}\frac{K}{\sqrt{T-t}}
\frac{1}
{\left(1-z^2(x)\right)^{1/4}}
<\frac{K_1}{\sqrt{T-t}}
\frac{1}
{\left(1-z^2(x)\right)^{1/4}},
\end{equation}

\vspace{4mm}
\noindent
where
$K_1=K\sqrt{2},$\  $x\in (t, T),$\ $j\in\mathbb{N},$

\vspace{-1mm}
$$
z(x)=\left(x-\frac{T+t}{2}\right)\frac{2}{T-t}.
$$

\vspace{4mm}

Note analogues of the estimates (\ref{101xx}), (\ref{101xxx})

\vspace{-1mm}
\begin{equation}
\label{101xxqq}
\left|
\int\limits_t^x \phi_{j_1}(s)ds
\right| <
\frac{C}{j_1 (1-(z(x))^2)^{1/4}},\ \ \
\left|
\int\limits_x^T \phi_{j_1}(s)ds
\right| <
\frac{C}{j_1 (1-(z(x))^2)^{1/4}},\ \ \
x\in (t, T),
\end{equation}

\vspace{4mm}
\noindent
where $j_1>0,$ constant $C$ does not depend on $j_1.$

Applying the estimates (\ref{ogo24}) and (\ref{101xxqq}) to (\ref{sixsix62}) 
gives the estimate (\ref{sixsix61}).
Using (\ref{sixsix70}), (\ref{sixsix60}), and (\ref{sixsix61}),
we obtain

$$
\left\vert
\sum_{j_1,j_2=0}^p\sum_{j_3=p+1}^{\infty}
C_{j_3 j_1 j_2 j_3}C_{j_2 j_1}\right\vert\le
K\sum\limits_{j_3=p+1}^{\infty}\frac{1}{j_3^2}
\left(1+\sum\limits_{j_1=1}^p \frac{1}{j_1}\right)\le
$$

\vspace{2mm}
\begin{equation}
\label{sixsix74}
\le K\int\limits_p^{\infty} \frac{dx}{x^2}
\left(2+\int\limits_1^p \frac{dx}{x}\right)=
\frac{K(2+ln p)}{p}\to 0
\end{equation}

\vspace{4mm}
\noindent
if $p\to\infty$, where constant $K$ is independent of $p.$
Thus, the equality (\ref{sixsix8}) is proved (see (\ref{sixsix71}), (\ref{sixsix72}),
(\ref{sixsix74})).

The relation (\ref{sixsix9}) is proved in complete analogy 
with the proof of equality (\ref{sixsix8}).
For (\ref{sixsix9}) we have (see (\ref{sixsix40}))

$$
\lim\limits_{p\to\infty}
\left(\sum_{j_1,j_2,j_3=0}^{p}
C_{j_1 j_3 j_2 j_3 j_2 j_1}+
\sum_{j_1,j_2,j_3=0}^{p}
C_{j_1 j_2 j_3 j_2 j_3 j_1}\right)=
2\lim\limits_{p\to\infty}\sum_{j_1,j_2,j_3=0}^{p}
C_{j_1 j_3 j_2 j_3 j_2 j_1}=
$$

\vspace{2mm}
$$
=
\lim\limits_{p\to\infty}\sum_{j_1,j_2,j_3=0}^{p}\biggl(
C_{j_1}C_{j_3 j_2 j_3 j_2 j_1}-C_{j_3 j_1}C_{j_2 j_3 j_2 j_1}+
C_{j_2 j_3 j_1}C_{j_3 j_2 j_1}-\biggr.
$$

\vspace{2mm}
$$
\biggl.-C_{j_3 j_2 j_3 j_1}C_{j_2 j_1}+
C_{j_2 j_3 j_2 j_3 j_1}C_{j_1}\biggr)=
$$

\vspace{2mm}
$$
=2\lim\limits_{p\to\infty}\left(
\sqrt{T-t}\sum_{j_2,j_3=0}^{p}C_{j_3 j_2 j_3 j_2 0}-
\sum_{j_1,j_2,j_3=0}^{p}
C_{j_2 j_1}C_{j_3 j_2 j_3 j_1}\right)=
$$

\vspace{2mm}
$$
=-2\lim\limits_{p\to\infty}
\sum_{j_1,j_2,j_3=0}^{p}
C_{j_2 j_1}C_{j_3 j_2 j_3 j_1}.
$$

\vspace{5mm}

To estimate the Fourier coefficient $C_{j_3 j_2 j_3 j_1}$, 
we use the following 
(see the proof of (\ref{sixsix8}) for more details)

$$
C_{j_3 j_2 j_3 j_1}=\int\limits_t^T\phi_{j_3}(t_4)\int\limits_{t}^{t_4}\phi_{j_2}(t_3)
\int\limits_t^{t_3}\phi_{j_3}(t_2)
\int\limits_t^{t_2}\phi_{j_1}(t_1)
dt_1 dt_2 dt_3 dt_4=
$$

\vspace{2mm}
$$
=\int\limits_t^T\phi_{j_3}(t_4)\int\limits_{t}^{t_4}\phi_{j_2}(t_3)
\int\limits_t^{t_3}\phi_{j_1}(t_1)
\int\limits_{t_1}^{t_3}\phi_{j_3}(t_2)
dt_2 dt_1 dt_3 dt_4=
$$

\vspace{2mm}
$$
=\int\limits_t^T\phi_{j_3}(t_4)\int\limits_{t}^{t_4}\phi_{j_2}(t_3)
\left(\int\limits_{t}^{t_3}\phi_{j_3}(t_2)
dt_2\right)
\int\limits_t^{t_3}\phi_{j_1}(t_1)
dt_1 dt_3 dt_4-
$$

\vspace{2mm}
$$
-\int\limits_t^T\phi_{j_3}(t_4)\int\limits_{t}^{t_4}\phi_{j_2}(t_3)
\int\limits_t^{t_3}\phi_{j_1}(t_1)
\left(\int\limits_{t}^{t_1}\phi_{j_3}(t_2)
dt_2\right) dt_1 dt_3 dt_4=
$$

\vspace{2mm}
$$
=\int\limits_t^T
\phi_{j_2}(t_3)
\left(\int\limits_{t}^{t_3}\phi_{j_3}(t_2)
dt_2\right)
\int\limits_t^{t_3}\phi_{j_1}(t_1)
dt_1
\left(\int\limits_{t_3}^{T}
\phi_{j_3}(t_4)
dt_4\right) dt_3-
$$

\vspace{2mm}
$$
-\int\limits_t^T
\phi_{j_2}(t_3)
\int\limits_t^{t_3}\phi_{j_1}(t_1)
\left(\int\limits_{t}^{t_1}\phi_{j_3}(t_2)
dt_2\right) dt_1
\left(\int\limits_{t_3}^{T}
\phi_{j_3}(t_4) dt_4\right) dt_3.
$$

\vspace{5mm}

Let us prove (\ref{sixsix10}). 
From (\ref{after80xx}) we obtain

\begin{equation}
\label{sixsix80}
\sum_{j_1=p+1}^{\infty}\sum_{j_2=p+1}^{\infty}
\sum_{j_3=p+1}^{\infty}
C_{j_3 j_2 j_3 j_1 j_2 j_1}=
-\sum_{j_1=0}^{p}\sum_{j_2=0}^{p}
\sum_{j_3=0}^{p}
C_{j_3 j_2 j_3 j_1 j_2 j_1}.
\end{equation}

\vspace{5mm}

Applying (\ref{sixsix40}) and (\ref{sixsix80}), we get (we replaced $j_3$ by $j_4$)

$$
\sum_{j_1,j_2,j_4=0}^{p}
C_{j_4 j_2 j_4 j_1 j_2 j_1}+
\sum_{j_1,j_2,j_4=0}^{p}
C_{j_1 j_2 j_1 j_4 j_2 j_4}=
2\sum_{j_1,j_2,j_4=0}^{p}
C_{j_4 j_2 j_4 j_1 j_2 j_1}=
$$

\vspace{2mm}
$$
=
\sum_{j_1,j_2,j_4=0}^{p}\biggl(
C_{j_4}C_{j_2 j_4 j_1 j_2 j_1}-C_{j_2 j_4}C_{j_4 j_1 j_2 j_1}+
C_{j_4 j_2 j_4}C_{j_1 j_2 j_1}-\biggr.
$$

\vspace{2mm}
$$
\biggl.-C_{j_1 j_4 j_2 j_4}C_{j_2 j_1}+
C_{j_2 j_1 j_4 j_2 j_4}C_{j_1}\biggr)=
$$

\vspace{2mm}
$$
=
2\sum_{j_1,j_2,j_4=0}^{p}\biggl(C_{j_2 j_1 j_4 j_2 j_4}C_{j_1}-C_{j_1 j_4 j_2 j_4}C_{j_2 j_1}
\biggr)+
$$

\vspace{2mm}
\begin{equation}
\label{sixsix83}
+
\sum_{j_1,j_2,j_4=0}^{p}
C_{j_4 j_2 j_4}C_{j_1 j_2 j_1}.
\end{equation}

\vspace{5mm}

Further, we have (see (\ref{after80xx}))

$$
\lim\limits_{p\to\infty}\sum_{j_1,j_2,j_4=0}^{p}
C_{j_4 j_2 j_4}C_{j_1 j_2 j_1}=
\lim\limits_{p\to\infty}\sum_{j_2=0}^{p}
\left(\sum_{j_1=0}^{p}C_{j_1 j_2 j_1}\right)^2=
$$

\vspace{2mm}
\begin{equation}
\label{sixsix84}
=\lim\limits_{p\to\infty}\sum_{j_2=0}^{p}
\left(\sum_{j_1=p+1}^{\infty}C_{j_1 j_2 j_1}\right)^2=0,
\end{equation}

\vspace{5mm}
\noindent
where we applied the equality (\ref{after1602}).

Furthermore, by analogy with the proof of (\ref{sixsix8}), we have

\begin{equation}
\label{sixsix85}
\lim\limits_{p\to\infty}
\sum_{j_1,j_2,j_4=0}^{p}\biggl(C_{j_2 j_1 j_4 j_2 j_4}C_{j_1}-C_{j_1 j_4 j_2 j_4}C_{j_2 j_1}
\biggr)=0.
\end{equation}

\vspace{4mm}

To estimate the Fourier coefficient $C_{j_1 j_4 j_2 j_4}$ in (\ref{sixsix85}),
we use the following 
(see the proof of (\ref{sixsix8}) for more details)

$$
C_{j_1 j_4 j_2 j_4}=\int\limits_t^T\phi_{j_1}(t_4)\int\limits_{t}^{t_4}\phi_{j_4}(t_3)
\int\limits_t^{t_3}\phi_{j_2}(t_2)
\left(\int\limits_t^{t_2}\phi_{j_4}(t_1)
dt_1\right) dt_2 dt_3 dt_4=
$$

\vspace{2mm}
$$
=\int\limits_t^T\phi_{j_1}(t_4)\int\limits_{t}^{t_4}
\phi_{j_2}(t_2)
\left(\int\limits_t^{t_2}\phi_{j_4}(t_1)
dt_1\right)
\int\limits_{t_2}^{t_4}
\phi_{j_4}(t_3)
dt_3 dt_2 dt_4=
$$

\vspace{2mm}
$$
=\int\limits_t^T\phi_{j_1}(t_4)
\left(\int\limits_{t}^{t_4}
\phi_{j_4}(t_3)
dt_3\right)
\int\limits_{t}^{t_4}
\phi_{j_2}(t_2)
\left(\int\limits_t^{t_2}\phi_{j_4}(t_1)
dt_1\right)
dt_2 dt_4-
$$

\vspace{2mm}
$$
-\int\limits_t^T\phi_{j_1}(t_4)\int\limits_{t}^{t_4}
\phi_{j_2}(t_2)\left(
\int\limits_{t}^{t_2}
\phi_{j_4}(t_3)
dt_3\right)
\left(\int\limits_t^{t_2}\phi_{j_4}(t_1)
dt_1\right)
dt_2 dt_4.
$$

\vspace{5mm}

The relations (\ref{sixsix80})--(\ref{sixsix85}) 
complete the proof of equality (\ref{sixsix10}).

Let us prove (\ref{sixsix4}).
Using (\ref{after80xx}), we get 

\begin{equation}
\label{sixsix90}
\sum_{j_1=p+1}^{\infty}\sum_{j_2=p+1}^{\infty}
\sum_{j_3=p+1}^{\infty}
C_{j_1 j_2 j_3 j_3 j_2 j_1}=
\sum_{j_1=0}^{p}\sum_{j_2=0}^{p}
\sum_{j_3=p+1}^{\infty}
C_{j_1 j_2 j_3 j_3 j_2 j_1}.
\end{equation}

\vspace{5mm}

Applying (\ref{sixsix40}) and (\ref{sixsix90}), we obtain

$$
2\sum_{j_1,j_2=0}^{p}\sum_{j_3=p+1}^{\infty}
C_{j_1 j_2 j_3 j_3 j_2 j_1}=
$$

\vspace{2mm}
$$
=
\sum_{j_1,j_2=0}^{p}\sum_{j_3=p+1}^{\infty}\biggl(
C_{j_1}C_{j_2 j_3 j_3 j_2 j_1}-C_{j_2 j_1}C_{j_3 j_3 j_2 j_1}+
\left(C_{j_3 j_2 j_1}\right)^2-\biggr.
$$

\vspace{2mm}
$$
\biggl.-C_{j_3 j_3 j_2 j_1}C_{j_2 j_1}+
C_{j_2 j_3 j_3 j_2 j_1}C_{j_1}\biggr)=
$$

\vspace{2mm}
$$
=
2\sum_{j_1,j_2=0}^{p}\sum_{j_3=p+1}^{\infty}
\biggl(C_{j_1}C_{j_2 j_3 j_3 j_2 j_1}-C_{j_2 j_1}C_{j_3 j_3 j_2 j_1}
\biggr)+
$$

\vspace{2mm}
\begin{equation}
\label{sixsix91}
+
\sum_{j_1,j_2=0}^{p}\sum_{j_3=p+1}^{\infty}
\left(C_{j_3 j_2 j_1}\right)^2.
\end{equation}

\vspace{5mm}

In \cite{2018a} (Sect.~1.7.2) the following estimate

\vspace{1mm}
$$
\sum_{j_1=0}^{\infty}\ldots
\sum_{j_{s-1}=0}^{\infty}
\sum_{j_s=p+1}^{\infty}\sum_{j_{s+1}=0}^{\infty}\ldots \sum_{j_k=0}^{\infty}
C_{j_k\ldots j_1}^2\le
$$

\vspace{1mm}
\begin{equation}
\label{fffoh}
\le L_k
\sum_{j_s=p+1}^{\infty}\frac{1}{j_s^2} \le L_k\int\limits_p^{\infty}
\frac{dx}{x^2}=
\frac{L_k}{p}
\end{equation}

\vspace{3mm}
\noindent
is proved for the polynomial and trigonometric cases, 
where $s=1,\ldots,k,$ constant $L_k$ depends on $k$ and $T-t.$

Using the estimate (\ref{fffoh}),  we get

\begin{equation}
\label{sixsix92}
\lim\limits_{p\to\infty}\sum_{j_1,j_2=0}^{p}\sum_{j_3=p+1}^{\infty}
\left(C_{j_3 j_2 j_1}\right)^2=0.
\end{equation}

\vspace{4mm}

By analogy with the proof of (\ref{sixsix8}), we have

\begin{equation}
\label{sixsix93}
\lim\limits_{p\to\infty}
\sum_{j_1,j_2=0}^{p}\sum_{j_3=p+1}^{\infty}
\biggl(C_{j_1}C_{j_2 j_3 j_3 j_2 j_1}-C_{j_2 j_1}C_{j_3 j_3 j_2 j_1}
\biggr)=0,
\end{equation}

\vspace{4mm}
\noindent
where we applied
the equality (\ref{after2509}).
To estimate the Fourier coefficient $C_{j_3 j_3 j_2 j_1}$ in (\ref{sixsix93}),
we used the following 
(see the proof of (\ref{sixsix8}) for more details)

$$
C_{j_3 j_3 j_2 j_1}=\int\limits_t^T\phi_{j_3}(t_4)\int\limits_{t}^{t_4}\phi_{j_3}(t_3)
\int\limits_t^{t_3}\phi_{j_2}(t_2)
\int\limits_t^{t_2}\phi_{j_1}(t_1)
dt_1dt_2 dt_3 dt_4=
$$

\vspace{2mm}
$$
=\int\limits_t^T\phi_{j_1}(t_1)\int\limits_{t_1}^{T}
\phi_{j_2}(t_2)
\int\limits_{t_2}^{T}
\phi_{j_3}(t_3)
\int\limits_{t_3}^{T}\phi_{j_3}(t_4)
dt_4dt_3 dt_2 dt_1=
$$

\vspace{2mm}
\begin{equation}
\label{sept10}
=\frac{1}{2}\int\limits_t^T\phi_{j_1}(t_1)\int\limits_{t_1}^{T}
\phi_{j_2}(t_2)
\left(\int\limits_{t_2}^{T}
\phi_{j_3}(t_3)
dt_3\right)^2 dt_2 dt_1.
\end{equation}

\vspace{5mm}

Combining the equalities (\ref{sixsix90})--(\ref{sixsix93}), 
we obtain (\ref{sixsix4}).

Let us prove (\ref{sixsix14}) (we replace $j_2$ by $j_4$ and $j_3$ by $j_2$
in (\ref{sixsix14})).
As noted in Sect.~7, the sequential order of the series
$$
\sum_{j_1=p+1}^{\infty}\sum_{j_2=p+1}^{\infty}
\sum_{j_4=p+1}^{\infty}
$$

\vspace{3mm}
\noindent
is not important. This follows directly from the formulas 
(\ref{after500}) and (\ref{after80xx}).

Applying the mentioned property and (\ref{after80xx}), we get

\begin{equation}
\label{sixsix94}
\sum_{j_1=p+1}^{\infty}\sum_{j_2=p+1}^{\infty}
\sum_{j_4=p+1}^{\infty}
C_{j_1 j_4 j_4 j_2 j_2 j_1}=
-\sum_{j_1=0}^{p}\sum_{j_2=p+1}^{\infty}
\sum_{j_4=p+1}^{\infty}
C_{j_1 j_4 j_4 j_2 j_2 j_1}.
\end{equation}

\vspace{4mm}

Observe that (see the above reasoning)

\vspace{-1mm}
\begin{equation}
\label{sixsix95}
\sum_{j_2=p+1}^{\infty}
\sum_{j_4=p+1}^{\infty}
C_{j_1 j_4 j_4 j_2 j_2 j_1}=
\sum_{j_4=p+1}^{\infty}
\sum_{j_2=p+1}^{\infty}
C_{j_1 j_4 j_4 j_2 j_2 j_1}.
\end{equation}

\vspace{4mm}

Using (\ref{sixsix40}) and (\ref{sixsix95}), we obtain

$$
\sum_{j_1=0}^{p}\sum_{j_2=p+1}^{\infty}
\sum_{j_4=p+1}^{\infty}
\biggl(C_{j_1 j_4 j_4 j_2 j_2 j_1}+
C_{j_1 j_2 j_2 j_4 j_4 j_1}\biggr)=
2\sum_{j_1=0}^{p}\sum_{j_2=p+1}^{\infty}
\sum_{j_4=p+1}^{\infty}
C_{j_1 j_4 j_4 j_2 j_2 j_1}=
$$

\vspace{2mm}
$$
=
\sum_{j_1=0}^{p}\sum_{j_2=p+1}^{\infty}
\sum_{j_4=p+1}^{\infty}
\biggl(
C_{j_1}C_{j_4 j_4 j_2 j_2 j_1}-C_{j_4 j_1}C_{j_4 j_2 j_2 j_1}+
C_{j_4 j_4 j_1}C_{j_2 j_2 j_1}-\biggr.
$$

\vspace{2mm}
$$
\biggl.-C_{j_2 j_4 j_4 j_1}C_{j_2 j_1}+
C_{j_2 j_2 j_4 j_4 j_1}C_{j_1}\biggr)=
$$

\vspace{2mm}
$$
=
\sum_{j_1=0}^{p}\sum_{j_2=p+1}^{\infty}
\sum_{j_4=p+1}^{\infty}
\biggl(
C_{j_1}C_{j_4 j_4 j_2 j_2 j_1}-C_{j_4 j_1}C_{j_4 j_2 j_2 j_1}
-C_{j_2 j_4 j_4 j_1}C_{j_2 j_1}+
C_{j_2 j_2 j_4 j_4 j_1}C_{j_1}\biggr)+
$$

\vspace{2mm}
\begin{equation}
\label{sixsix96}
+\sum_{j_1=0}^{p}\left(\sum_{j_2=p+1}^{\infty}
C_{j_2 j_2 j_1}\right)^2.
\end{equation}

\vspace{5mm}

The equality 
\begin{equation}
\label{sixsix97}
\lim\limits_{p\to\infty}\sum_{j_1=0}^{p}\left(\sum_{j_2=p+1}^{\infty}
C_{j_2 j_2 j_1}\right)^2=0
\end{equation}

\vspace{4mm}
\noindent
follows from 
the relation (\ref{after1601}).

By analogy with the proof of equality (\ref{sixsix8}) we obtain

$$
\lim\limits_{p\to\infty}
\sum_{j_1=0}^{p}\sum_{j_2=p+1}^{\infty}
\sum_{j_4=p+1}^{\infty}
\biggl(
C_{j_1}C_{j_4 j_4 j_2 j_2 j_1}-C_{j_4 j_1}C_{j_4 j_2 j_2 j_1}
-\biggr.
$$

\vspace{2mm}
\begin{equation}
\label{sixsix98}
\biggl.-
C_{j_2 j_4 j_4 j_1}C_{j_2 j_1}+
C_{j_2 j_2 j_4 j_4 j_1}C_{j_1}\biggr)=0,
\end{equation}

\vspace{3mm}
\noindent
where we applied
the equality (\ref{after2507}).
To estimate the Fourier coefficient $C_{j_2 j_4 j_4 j_1}$ in (\ref{sixsix98}),
we used the following 
(see the proof of (\ref{sixsix8}) for more details)

$$
C_{j_2 j_4 j_4 j_1}=\int\limits_t^T\phi_{j_2}(t_4)\int\limits_{t}^{t_4}\phi_{j_4}(t_3)
\int\limits_t^{t_3}\phi_{j_4}(t_2)
\int\limits_t^{t_2}\phi_{j_1}(t_1)
dt_1dt_2 dt_3 dt_4=
$$

\vspace{2mm}
$$
=\int\limits_t^T\phi_{j_2}(t_4)\int\limits_{t}^{t_4}
\phi_{j_1}(t_1)
\int\limits_{t_1}^{t_4}
\phi_{j_4}(t_2)
\int\limits_{t_2}^{t_4}\phi_{j_4}(t_3)
dt_3dt_2 dt_1 dt_4=
$$

\vspace{2mm}
$$
=\frac{1}{2}\int\limits_t^T\phi_{j_2}(t_4)\int\limits_{t}^{t_4}
\phi_{j_1}(t_1)
\left(\int\limits_{t_1}^{t_4}
\phi_{j_4}(t_2)dt_2\right)^2 dt_1 dt_4=
$$

\vspace{2mm}
$$
=\frac{1}{2}\int\limits_t^T\phi_{j_2}(t_4)
\left(\int\limits_{t}^{t_4}
\phi_{j_4}(t_2)dt_2\right)^2\int\limits_{t}^{t_4}
\phi_{j_1}(t_1)
dt_1 dt_4+
$$

\vspace{2mm}
$$
+\frac{1}{2}\int\limits_t^T\phi_{j_2}(t_4)\int\limits_{t}^{t_4}
\phi_{j_1}(t_1)
\left(\int\limits_{t}^{t_1}
\phi_{j_4}(t_2)dt_2\right)^2 dt_1 dt_4-
$$

\vspace{2mm}
$$
-\int\limits_t^T\phi_{j_2}(t_4)
\left(\int\limits_{t}^{t_4}
\phi_{j_4}(t_2)dt_2\right)\int\limits_{t}^{t_4}
\phi_{j_1}(t_1)
\left(\int\limits_{t}^{t_1}
\phi_{j_4}(t_2)dt_2\right)
dt_1 dt_4.
$$

\vspace{5mm}

The relation (\ref{sixsix14}) follows from (\ref{sixsix94}),
(\ref{sixsix96})--(\ref{sixsix98}).

Consider (\ref{sixsix3}). Using the integration 
order replacement, we obtain

$$
C_{j_3j_3j_2j_2j_1j_1}=
$$

\vspace{1mm}
$$
=
\frac{1}{2}\int\limits_t^T \phi_{j_3}(t_6)
\int\limits_t^{t_6} \phi_{j_3}(t_5)
\int\limits_t^{t_5} \phi_{j_2}(t_4)
\int\limits_t^{t_4} \phi_{j_2}(t_3)
\left(\int\limits_t^{t_3} \phi_{j_1}(t_1)dt_1\right)^2 dt_3 dt_4 dt_5 dt_6=
$$

\vspace{2mm}
$$
=
\frac{1}{2}\int\limits_t^T 
\phi_{j_2}(t_3)
\left(\int\limits_t^{t_3} \phi_{j_1}(t_1)dt_1\right)^2
\int\limits_{t_3}^T 
\phi_{j_2}(t_4) 
\int\limits_{t_4}^T 
\phi_{j_3}(t_5) 
\int\limits_{t_5}^T 
\phi_{j_3}(t_6) 
dt_6 dt_5 dt_4 dt_3=
$$

\vspace{2mm}
\begin{equation}
\label{sixsix100}
=
\frac{1}{4}\int\limits_t^T 
\phi_{j_2}(t_3)
\left(\int\limits_t^{t_3} \phi_{j_1}(t_1)dt_1\right)^2
\int\limits_{t_3}^T 
\phi_{j_2}(t_4) 
\left(\int\limits_{t_4}^T 
\phi_{j_3}(t_5) 
dt_5\right)^2 dt_4 dt_3.
\end{equation}

\vspace{5mm}

Applying the estimates (\ref{101xxqq}) to 
(\ref{sixsix100}) 
gives the following estimate 

\begin{equation}
\label{sixsix101}
|C_{j_3j_3j_2j_2j_1j_1}|\le \frac{K}{j_1^2 j_3^2}\ \ \ (j_1, j_3>0,\ j_2\ge 0),
\end{equation}

\vspace{4mm}
\noindent
where constant $K$ does not depend on $j_1, j_2, j_3$.

Further, we get (see (\ref{after500}))

$$
\sum_{j_1=p+1}^{\infty}\sum_{j_2=p+1}^{\infty}
\sum_{j_3=p+1}^{\infty}
C_{j_3 j_3 j_2 j_2 j_1 j_1}=
\sum_{j_1=p+1}^{\infty}\sum_{j_3=p+1}^{\infty}
\sum_{j_2=p+1}^{\infty}
C_{j_3 j_3 j_2 j_2 j_1 j_1}=
$$

\vspace{2mm}
\begin{equation}
\label{sixsix102}
=\frac{1}{2}
\sum_{j_1=p+1}^{\infty}\sum_{j_3=p+1}^{\infty}
C_{j_3j_3j_2j_2j_1j_1}\biggl|_{(j_2 j_2)\curvearrowright (\cdot)}\biggr. -
\sum_{j_2=0}^{p}\sum_{j_1=p+1}^{\infty}\sum_{j_3=p+1}^{\infty}
C_{j_3 j_3 j_2 j_2 j_1 j_1},
\end{equation}

\vspace{5mm}
\noindent
where
$$
C_{j_3j_3j_2j_2j_1j_1}\biggl|_{(j_2 j_2)\curvearrowright (\cdot)}\biggr.=
$$

\vspace{2mm}
$$
=
\int\limits_t^T \phi_{j_3}(t_6)
\int\limits_t^{t_6} \phi_{j_3}(t_5)
\int\limits_t^{t_5} 
\int\limits_t^{t_4} 
\phi_{j_1}(t_2)
\int\limits_t^{t_2} 
\phi_{j_1}(t_1)
dt_1 dt_2 dt_4 dt_5 dt_6=
$$

\vspace{2mm}
$$
=
\int\limits_t^T \phi_{j_3}(t_6)
\int\limits_t^{t_6} \phi_{j_3}(t_5)
\int\limits_t^{t_5} 
\phi_{j_1}(t_2)
\int\limits_t^{t_2} 
\phi_{j_1}(t_1)
dt_1
\int\limits_{t_2}^{t_5} 
dt_4 dt_2 dt_5 dt_6=
$$

\vspace{2mm}
$$
=
\int\limits_t^T \phi_{j_3}(t_6)
\int\limits_t^{t_6} \phi_{j_3}(t_5)(t_5-t)
\int\limits_t^{t_5} 
\phi_{j_1}(t_2)
\int\limits_t^{t_2} 
\phi_{j_1}(t_1)
dt_1dt_2 dt_5 dt_6+
$$

\vspace{2mm}
$$
+\int\limits_t^T \phi_{j_3}(t_6)
\int\limits_t^{t_6} \phi_{j_3}(t_5)
\int\limits_t^{t_5} 
\phi_{j_1}(t_2)(t-t_2)
\int\limits_t^{t_2} 
\phi_{j_1}(t_1)
dt_1dt_2 dt_5 dt_6\stackrel{\sf def}{=}
$$

\vspace{2mm}
\begin{equation}
\label{sixsix103}
\stackrel{\sf def}{=}
C'_{j_3 j_3 j_1 j_1}+C''_{j_3 j_3 j_1 j_1}.
\end{equation}

\vspace{5mm}

Let us substitute (\ref{sixsix103}) into (\ref{sixsix102})

$$
\sum_{j_1=p+1}^{\infty}\sum_{j_2=p+1}^{\infty}
\sum_{j_3=p+1}^{\infty}
C_{j_3 j_3 j_2 j_2 j_1 j_1}=
\frac{1}{2}
\sum_{j_1=p+1}^{\infty}\sum_{j_3=p+1}^{\infty}
C'_{j_3 j_3 j_1 j_1}+
$$

\vspace{2mm}
\begin{equation}
\label{sixsix104}
+
\frac{1}{2}
\sum_{j_1=p+1}^{\infty}\sum_{j_3=p+1}^{\infty}
C''_{j_3 j_3 j_1 j_1}
-
\sum_{j_2=0}^{p}\sum_{j_1=p+1}^{\infty}\sum_{j_3=p+1}^{\infty}
C_{j_3 j_3 j_2 j_2 j_1 j_1}.
\end{equation}

\vspace{5mm}

The relation (\ref{after2507})
implies that

\begin{equation}
\label{sixsix105}
\lim\limits_{p\to\infty}\sum_{j_1=p+1}^{\infty}\sum_{j_3=p+1}^{\infty}
C'_{j_3 j_3 j_1 j_1}=0,\ \ \ 
\lim\limits_{p\to\infty}\sum_{j_1=p+1}^{\infty}\sum_{j_3=p+1}^{\infty}
C''_{j_3 j_3 j_1 j_1}=0.
\end{equation}

\vspace{5mm}

From the estimate (\ref{sixsix101}) we get

$$
\left\vert\sum_{j_2=0}^{p}\sum_{j_1=p+1}^{\infty}\sum_{j_3=p+1}^{\infty}
C_{j_3 j_3 j_2 j_2 j_1 j_1}\right\vert\le
K(p+1) \sum_{j_1=p+1}^{\infty}\frac{1}{j_1^2}
\sum_{j_3=p+1}^{\infty}\frac{1}{j_3^2}\le 
$$

\vspace{2mm}
\begin{equation}
\label{sixsix106}
\le
K(p+1) \left(\int\limits_p^{\infty} \frac{dx}{x^2}\right)^2
\le \frac{K(p+1)}{p^2}\ \to\ 0
\end{equation}

\vspace{4mm}
\noindent
if $p\to\infty$, where constant $K$ is independent of $p$.

The relations (\ref{sixsix104})--(\ref{sixsix106})
complete the proof of (\ref{sixsix3}).

Let us prove (\ref{sixsix7}). Using the integration
order replacement, we get

$$
C_{j_2 j_3 j_3 j_2 j_1 j_1}=
$$

\vspace{1mm}
$$
=
\frac{1}{2}\int\limits_t^T 
\phi_{j_2}(t_6)
\int\limits_t^{t_6} 
\phi_{j_3}(t_5)
\int\limits_{t}^{t_5} 
\phi_{j_3}(t_4) 
\int\limits_t^{t_4}
\phi_{j_2}(t_3) 
\left(\int\limits_t^{t_3}
\phi_{j_1}(t_1)dt_1\right)^2 
dt_3 dt_4 dt_5 dt_6=
$$

\vspace{2mm}
$$
=
\frac{1}{2}\int\limits_t^T 
\phi_{j_2}(t_3) 
\left(\int\limits_t^{t_3}
\phi_{j_1}(t_1)dt_1\right)^2 
\int\limits_{t_3}^T
\phi_{j_3}(t_4) 
\int\limits_{t_4}^T
\phi_{j_3}(t_5)
\int\limits_{t_5}^T
\phi_{j_2}(t_6)dt_6
dt_5 dt_4 dt_3=
$$

\vspace{2mm}
$$
=
\frac{1}{2}\int\limits_t^T 
\phi_{j_2}(t_3) 
\left(\int\limits_t^{t_3}
\phi_{j_1}(t_1)dt_1\right)^2 
\int\limits_{t_3}^T
\phi_{j_3}(t_5)
\int\limits_{t_5}^T
\phi_{j_2}(t_6)dt_6
\int\limits_{t_3}^{t_5}
\phi_{j_3}(t_4) 
dt_4 dt_5 dt_3=
$$

\vspace{2mm}
$$
=
\frac{1}{2}\int\limits_t^T 
\phi_{j_2}(t_3) 
\left(\int\limits_t^{t_3}
\phi_{j_1}(t_1)dt_1\right)^2 
\int\limits_{t_3}^T
\phi_{j_3}(t_5)
\left(\int\limits_{t_5}^T
\phi_{j_2}(t_6)dt_6\right)
\left(\int\limits_{t}^{t_5}
\phi_{j_3}(t_4)dt_4\right)dt_5 dt_3-
$$

\vspace{2mm}
\begin{equation}
\label{sixsix107}
-
\frac{1}{2}\int\limits_t^T 
\phi_{j_2}(t_3) 
\left(\int\limits_t^{t_3}
\phi_{j_1}(t_1)dt_1\right)^2 
\left(\int\limits_{t}^{t_3}
\phi_{j_3}(t_4)dt_4\right)\int\limits_{t_3}^T
\phi_{j_3}(t_5)
\left(\int\limits_{t_5}^T
\phi_{j_2}(t_6)dt_6\right)
dt_5 dt_3.
\end{equation}

\vspace{5mm}

Applying (\ref{after80xx}) and (\ref{after500}), we obtain

$$
-\sum_{j_1=p+1}^{\infty}\sum_{j_2=p+1}^{\infty}
\sum_{j_3=p+1}^{\infty}
C_{j_2 j_3 j_3 j_2 j_1 j_1}=
-\sum_{j_1=p+1}^{\infty}\sum_{j_3=p+1}^{\infty}\sum_{j_2=p+1}^{\infty}
C_{j_2 j_3 j_3 j_2 j_1 j_1}=
$$

\vspace{2mm}
$$
=
\sum_{j_2=0}^{p}\sum_{j_1=p+1}^{\infty}\sum_{j_3=p+1}^{\infty}
C_{j_2 j_3 j_3 j_2 j_1 j_1}=
$$

\vspace{2mm}
$$
=\frac{1}{2}
\sum_{j_2=0}^{p}\sum_{j_1=p+1}^{\infty}
C_{j_2 j_3 j_3 j_2 j_1 j_1}\biggl|_{(j_3 j_3)\curvearrowright (\cdot)}\biggr. -
\sum_{j_2=0}^{p}\sum_{j_3=0}^{p}\sum_{j_1=p+1}^{\infty}
C_{j_2 j_3 j_3 j_2 j_1 j_1}=
$$

\vspace{2mm}
$$
=\frac{1}{2}
\sum_{j_2=0}^{p}\sum_{j_1=p+1}^{\infty}
C_{j_2 j_3 j_3 j_2 j_1 j_1}\biggl|_{(j_3 j_3)\curvearrowright (\cdot)}\biggr. -
\sum_{j_1=p+1}^{\infty}
C_{0000 j_1 j_1}-
$$

\vspace{2mm}
$$
-\sum_{j_3=1}^{p}\sum_{j_1=p+1}^{\infty}
C_{0 j_3 j_3 0 j_1 j_1}-
\sum_{j_2=1}^{p}\sum_{j_1=p+1}^{\infty}
C_{j_2 00 j_2 j_1 j_1}-
$$

\vspace{2mm}
\begin{equation}
\label{sixsix108}
-
\sum_{j_2=1}^{p}\sum_{j_3=1}^{p}\sum_{j_1=p+1}^{\infty}
C_{j_2 j_3 j_3 j_2 j_1 j_1}.
\end{equation}

\vspace{5mm}

The equality
\begin{equation}
\label{sixsix109a}
\lim\limits_{p\to\infty}
\frac{1}{2}
\sum_{j_2=0}^{p}\sum_{j_1=p+1}^{\infty}
C_{j_2 j_3 j_3 j_2 j_1 j_1}\biggl|_{(j_3 j_3)\curvearrowright (\cdot)}\biggr.=0
\end{equation}

\vspace{3mm}
\noindent
follows from the inequality similar to (\ref{after9043})
(see the proof of Theorem~17),
where we used the following representation

$$
C_{j_2j_3j_3j_2j_1j_1}\biggl|_{(j_3 j_3)\curvearrowright (\cdot)}\biggr.=
$$

\vspace{2mm}
$$
=
\int\limits_t^T \phi_{j_2}(t_6)
\int\limits_t^{t_6} 
\int\limits_t^{t_4} 
\phi_{j_2}(t_3)
\int\limits_t^{t_3} 
\phi_{j_1}(t_2)
\int\limits_t^{t_2} 
\phi_{j_1}(t_1)
dt_1 dt_2 dt_3 dt_4 dt_6=
$$

\vspace{2mm}
$$
=
\int\limits_t^T \phi_{j_2}(t_6)
\int\limits_t^{t_6} \phi_{j_2}(t_3)
\int\limits_t^{t_3} 
\phi_{j_1}(t_2)
\int\limits_t^{t_2} 
\phi_{j_1}(t_1)
dt_1 dt_2
\int\limits_{t_3}^{t_6} 
dt_4 dt_3 dt_6=
$$

\vspace{2mm}
$$
+\int\limits_t^T \phi_{j_2}(t_6)(t_6-t)
\int\limits_t^{t_6} \phi_{j_2}(t_3)
\int\limits_t^{t_3} 
\phi_{j_1}(t_2)
\int\limits_t^{t_2} 
\phi_{j_1}(t_1)
dt_1dt_2 dt_3 dt_6+
$$

\vspace{2mm}
$$
+\int\limits_t^T \phi_{j_2}(t_6)
\int\limits_t^{t_6} \phi_{j_2}(t_3)(t-t_3)
\int\limits_t^{t_3} 
\phi_{j_1}(t_2)
\int\limits_t^{t_2} 
\phi_{j_1}(t_1)
dt_1dt_2 dt_3 dt_6\stackrel{\sf def}{=}
$$

\vspace{2mm}
\begin{equation}
\label{sept12}
\stackrel{\sf def}{=}
C^{*}_{j_2 j_2 j_1 j_1}+C^{**}_{j_2 j_2 j_1 j_1}.
\end{equation}

\vspace{5mm}

Applying the estimates (\ref{101xxqq})
and (\ref{after1940}) ($\varepsilon=1/2$) to 
(\ref{sixsix107}) 
gives the following estimates 

\begin{equation}
\label{sixsix109}
|C_{j_2 j_3 j_3 j_2 j_1 j_1}|\le \frac{K}{j_1^2 j_2 j_3^{3/4}}\ \ \ (j_1, j_2, j_3>0),
\end{equation}

\begin{equation}
\label{sixsix110}
|C_{j_2 00 j_2 j_1 j_1}|\le \frac{K}{j_1^2 j_2}\ \ \ (j_1, j_2> 0),
\end{equation}

\begin{equation}
\label{sixsix111}
|C_{0 j_3 j_3 0 j_1 j_1}|\le \frac{K}{j_1^2 j_3}\ \ \ (j_1, j_3>0),
\end{equation}

\begin{equation}
\label{sixsix112}
|C_{0000 j_1 j_1}|\le \frac{K}{j_1^2}\ \ \ (j_1> 0).
\end{equation}

\vspace{5mm}

Using the estimate (\ref{sixsix109}), we have

$$
\left\vert
\sum_{j_2=1}^{p}\sum_{j_3=1}^{p}\sum_{j_1=p+1}^{\infty}
C_{j_2 j_3 j_3 j_2 j_1 j_1}
\right\vert \le 
K\sum_{j_1=p+1}^{\infty}\frac{1}{j_1^2} \sum_{j_2=1}^{p} \frac{1}{j_2}
\sum_{j_3=1}^{p}\frac{1}{j_3^{3/4}}\le
$$

\vspace{2mm}
\begin{equation}
\label{sixsix113}
\le K \int\limits_p^{\infty} \frac{dx}{x^2}\left(1+\int\limits_1^p \frac{dx}{x}\right)
\left(1+\int\limits_1^p \frac{dx}{x^{3/4}}\right)\le
K_1 \frac{1+ ln p}{p^{3/4}}\ \to 0
\end{equation}

\vspace{4mm}
\noindent
if $p\to\infty,$ where constants $K, K_1$ do not depend on $p.$

Similarly we get (see (\ref{sixsix110})--(\ref{sixsix112}))

\begin{equation}
\label{sixsix114}
\left\vert\sum_{j_1=p+1}^{\infty}
C_{0000 j_1 j_1}\right\vert+
\left\vert\sum_{j_3=1}^{p}\sum_{j_1=p+1}^{\infty}
C_{0 j_3 j_3 0 j_1 j_1}\right\vert+
\left\vert\sum_{j_2=1}^{p}\sum_{j_1=p+1}^{\infty}
C_{j_2 00 j_2 j_1 j_1}\right\vert\ \to 0
\end{equation}

\vspace{4mm}
\noindent
if $p\to\infty.$

The relations 
(\ref{sixsix108}), (\ref{sixsix109a}), (\ref{sixsix113}), 
(\ref{sixsix114}) prove (\ref{sixsix7}).

Consider (\ref{sixsix6}). Using the integration
order replacement, we get

$$
C_{j_3 j_2 j_3 j_2 j_1 j_1}=
$$

\vspace{1mm}
$$
=
\frac{1}{2}\int\limits_t^T 
\phi_{j_3}(t_6)
\int\limits_t^{t_6} 
\phi_{j_2}(t_5)
\int\limits_{t}^{t_5} 
\phi_{j_3}(t_4) 
\int\limits_t^{t_4}
\phi_{j_2}(t_3) 
\left(\int\limits_t^{t_3}
\phi_{j_1}(t_1)dt_1\right)^2 
dt_3 dt_4 dt_5 dt_6=
$$

\vspace{2mm}
$$
=
\frac{1}{2}\int\limits_t^T 
\phi_{j_2}(t_3) 
\left(\int\limits_t^{t_3}
\phi_{j_1}(t_1)dt_1\right)^2 
\int\limits_{t_3}^T
\phi_{j_3}(t_4) 
\int\limits_{t_4}^T
\phi_{j_2}(t_5)
\int\limits_{t_5}^T
\phi_{j_3}(t_6)dt_6
dt_5 dt_4 dt_3=
$$

\vspace{2mm}
$$
=
\frac{1}{2}\int\limits_t^T 
\phi_{j_2}(t_3) 
\left(\int\limits_t^{t_3}
\phi_{j_1}(t_1)dt_1\right)^2 
\int\limits_{t_3}^T
\phi_{j_2}(t_5)
\int\limits_{t_5}^T
\phi_{j_3}(t_6)dt_6
\int\limits_{t_3}^{t_5}
\phi_{j_3}(t_4) 
dt_4 dt_5 dt_3=
$$

\vspace{2mm}
$$
=
\frac{1}{2}\int\limits_t^T 
\phi_{j_2}(t_3) 
\left(\int\limits_t^{t_3}
\phi_{j_1}(t_1)dt_1\right)^2 
\int\limits_{t_3}^T
\phi_{j_2}(t_5)
\left(\int\limits_{t}^{t_5}
\phi_{j_3}(t_4) 
dt_4\right)
\left(\int\limits_{t_5}^T
\phi_{j_3}(t_6)dt_6\right)
dt_5 dt_3-
$$

\vspace{2mm}
\begin{equation}
\label{sixsix115}
-
\frac{1}{2}\int\limits_t^T 
\phi_{j_2}(t_3) 
\left(\int\limits_t^{t_3}
\phi_{j_1}(t_1)dt_1\right)^2 
\left(\int\limits_{t}^{t_3}
\phi_{j_3}(t_4) 
dt_4 \right)
\int\limits_{t_3}^T
\phi_{j_2}(t_5)
\left(\int\limits_{t_5}^T
\phi_{j_3}(t_6)dt_6\right)dt_5 dt_3.
\end{equation}

\vspace{5mm}

Applying (\ref{after80xx}), we obtain

$$
\sum_{j_1=p+1}^{\infty}\sum_{j_2=p+1}^{\infty}
\sum_{j_3=p+1}^{\infty}
C_{j_3 j_2 j_3 j_2 j_1 j_1}=
\sum_{j_1=p+1}^{\infty}\sum_{j_3=p+1}^{\infty}
\sum_{j_2=p+1}^{\infty}
C_{j_3 j_2 j_3 j_2 j_1 j_1}=
$$

\vspace{2mm}
\begin{equation}
\label{sixsix116}
=-
\sum_{j_2=0}^{p}
\sum_{j_1=p+1}^{\infty}\sum_{j_3=p+1}^{\infty}
C_{j_3 j_2 j_3 j_2 j_1 j_1}.
\end{equation}

\vspace{4mm}

Further proof of the equality (\ref{sixsix6})
is based on the relations (\ref{sixsix115}), (\ref{sixsix116}) and 
is similar to the proof of the formula (\ref{sixsix7}).

Let us prove (\ref{sixsix1}). Applying the integration 
order replacement, we obtain

$$
C_{j_3 j_3 j_2 j_1 j_2 j_1}=
$$

\vspace{1mm}
$$
=\hspace{-1mm}
\int\limits_t^T 
\phi_{j_3}(t_6)
\int\limits_t^{t_6} 
\phi_{j_3}(t_5)
\int\limits_{t}^{t_5} 
\phi_{j_2}(t_4) 
\int\limits_t^{t_4}
\phi_{j_1}(t_3) 
\int\limits_t^{t_3}
\phi_{j_2}(t_2)
\int\limits_t^{t_2}
\phi_{j_1}(t_1)
dt_1 dt_2 dt_3 dt_4 dt_5 dt_6=
$$

\vspace{2mm}
$$
=\hspace{-1mm}
\int\limits_t^T 
\phi_{j_1}(t_1)
\int\limits_{t_1}^T
\phi_{j_2}(t_2)
\int\limits_{t_2}^T
\phi_{j_1}(t_3) 
\int\limits_{t_3}^T
\phi_{j_2}(t_4) 
\int\limits_{t_4}^T
\phi_{j_3}(t_5)
\int\limits_{t_5}^T
\phi_{j_3}(t_6)
dt_6 dt_5 dt_4 dt_3 dt_2 dt_1=
$$

\vspace{2mm}
$$
=\frac{1}{2}
\int\limits_t^T 
\phi_{j_1}(t_1)
\int\limits_{t_1}^T
\phi_{j_2}(t_2)
\int\limits_{t_2}^T
\phi_{j_1}(t_3) 
\int\limits_{t_3}^T
\phi_{j_2}(t_4) 
\left(\int\limits_{t_4}^T
\phi_{j_3}(t_5)
dt_5\right)^2 dt_4 dt_3 dt_2 dt_1=
$$

\vspace{2mm}
$$
=\frac{1}{2}
\int\limits_t^T 
\phi_{j_2}(t_4) 
\left(\int\limits_{t_4}^T
\phi_{j_3}(t_5)
dt_5\right)^2
\int\limits_t^{t_4}
\phi_{j_1}(t_3) 
\int\limits_t^{t_3}
\phi_{j_2}(t_2)
\int\limits_t^{t_2}
\phi_{j_1}(t_1)
dt_1 dt_2 dt_3 dt_4=
$$

\vspace{2mm}
$$
=\frac{1}{2}
\int\limits_t^T 
\phi_{j_2}(t_4) 
\left(\int\limits_{t_4}^T
\phi_{j_3}(t_5)
dt_5\right)^2
\int\limits_t^{t_4}
\phi_{j_2}(t_2)
\int\limits_t^{t_2}
\phi_{j_1}(t_1)
dt_1
\int\limits_{t_2}^{t_4}
\phi_{j_1}(t_3) 
dt_3 dt_2 dt_4=
$$

\vspace{2mm}
$$
=\frac{1}{2}
\int\limits_t^T 
\phi_{j_2}(t_4) 
\left(\int\limits_{t_4}^T
\phi_{j_3}(t_5)
dt_5\right)^2
\left(\int\limits_{t}^{t_4}
\phi_{j_1}(t_3) 
dt_3\right)
\int\limits_t^{t_4}
\phi_{j_2}(t_2)
\left(\int\limits_t^{t_2}
\phi_{j_1}(t_1)
dt_1\right)
dt_2 dt_4-
$$

\vspace{2mm}
\begin{equation}
\label{sixsix120}
-\frac{1}{2}
\int\limits_t^T 
\phi_{j_2}(t_4) 
\left(\int\limits_{t_4}^T
\phi_{j_3}(t_5)
dt_5\right)^2
\int\limits_t^{t_4}
\phi_{j_2}(t_2)
\left(\int\limits_t^{t_2}
\phi_{j_1}(t_1)
dt_1\right)^2 dt_2 dt_4.
\end{equation}

\vspace{5mm}

Using (\ref{after80xx}), we get

$$
\sum_{j_1=p+1}^{\infty}\sum_{j_2=p+1}^{\infty}
\sum_{j_3=p+1}^{\infty}
C_{j_3 j_3 j_2 j_1 j_2 j_1}=
\sum_{j_1=p+1}^{\infty}\sum_{j_3=p+1}^{\infty}
\sum_{j_2=p+1}^{\infty}
C_{j_3 j_3 j_2 j_1 j_2 j_1}=
$$

\vspace{2mm}
\begin{equation}
\label{sixsix121}
=-
\sum_{j_2=0}^{p}
\sum_{j_1=p+1}^{\infty}\sum_{j_3=p+1}^{\infty}
C_{j_3 j_3 j_2 j_1 j_2 j_1}.
\end{equation}

\vspace{4mm}

Further proof of the equality (\ref{sixsix1})
is based on the relations (\ref{sixsix120}), (\ref{sixsix121}) and 
is similar to the proof of the relations (\ref{sixsix7}),
(\ref{sixsix6}).

Consider (\ref{sixsix2}). 
Using the integration 
order replacement, we have

$$
C_{j_3 j_3 j_1 j_2 j_2 j_1}=
$$

\vspace{1mm}
$$
=\hspace{-1mm}
\int\limits_t^T 
\phi_{j_3}(t_6)
\int\limits_t^{t_6} 
\phi_{j_3}(t_5)
\int\limits_{t}^{t_5} 
\phi_{j_1}(t_4) 
\int\limits_t^{t_4}
\phi_{j_2}(t_3) 
\int\limits_t^{t_3}
\phi_{j_2}(t_2)
\int\limits_t^{t_2}
\phi_{j_1}(t_1)
dt_1 dt_2 dt_3 dt_4 dt_5 dt_6=
$$

\vspace{2mm}
$$
=\hspace{-1mm}
\int\limits_t^T 
\phi_{j_1}(t_1)
\int\limits_{t_1}^T
\phi_{j_2}(t_2)
\int\limits_{t_2}^T
\phi_{j_2}(t_3) 
\int\limits_{t_3}^T
\phi_{j_1}(t_4) 
\int\limits_{t_4}^T
\phi_{j_3}(t_5)
\int\limits_{t_5}^T
\phi_{j_3}(t_6)
dt_6 dt_5 dt_4 dt_3 dt_2 dt_1=
$$

\vspace{2mm}
$$
=\frac{1}{2}
\int\limits_t^T 
\phi_{j_1}(t_1)
\int\limits_{t_1}^T
\phi_{j_2}(t_2)
\int\limits_{t_2}^T
\phi_{j_2}(t_3) 
\int\limits_{t_3}^T
\phi_{j_1}(t_4) 
\left(\int\limits_{t_4}^T
\phi_{j_3}(t_5)
dt_5\right)^2 dt_4 dt_3 dt_2 dt_1=
$$

\vspace{2mm}
$$
=\frac{1}{2}
\int\limits_t^T 
\phi_{j_1}(t_4) 
\left(\int\limits_{t_4}^T
\phi_{j_3}(t_5)
dt_5\right)^2
\int\limits_t^{t_4}
\phi_{j_2}(t_3) 
\int\limits_t^{t_3}
\phi_{j_2}(t_2)
\int\limits_t^{t_2}
\phi_{j_1}(t_1)
dt_1 dt_2 dt_3 dt_4=
$$

\vspace{2mm}
$$
=\frac{1}{2}
\int\limits_t^T 
\phi_{j_1}(t_4) 
\left(\int\limits_{t_4}^T
\phi_{j_3}(t_5)
dt_5\right)^2
\int\limits_t^{t_4}
\phi_{j_2}(t_2)
\int\limits_t^{t_2}
\phi_{j_1}(t_1)
dt_1
\int\limits_{t_2}^{t_4}
\phi_{j_2}(t_3) 
dt_3 dt_2 dt_4=
$$

\vspace{2mm}
$$
=\frac{1}{2}
\int\limits_t^T 
\phi_{j_1}(t_4) 
\left(\int\limits_{t_4}^T
\phi_{j_3}(t_5)
dt_5\right)^2
\left(\int\limits_{t}^{t_4}
\phi_{j_2}(t_3) 
dt_3\right)
\int\limits_t^{t_4}
\phi_{j_2}(t_2)
\left(\int\limits_t^{t_2}
\phi_{j_1}(t_1)
dt_1\right)
dt_2 dt_4-
$$

\vspace{2mm}
\begin{equation}
\label{sixsix123}
-\frac{1}{2}
\int\limits_t^T 
\phi_{j_1}(t_4) 
\left(\int\limits_{t_4}^T
\phi_{j_3}(t_5)
dt_5\right)^2
\int\limits_t^{t_4}
\phi_{j_2}(t_2)
\left(\int\limits_t^{t_2}
\phi_{j_1}(t_1)
dt_1\right)
\left(\int\limits_t^{t_2}
\phi_{j_2}(t_3)
dt_3\right)
dt_2 dt_4.
\end{equation}

\vspace{5mm}

Applying (\ref{after80xx}) and (\ref{after500}), we obtain

$$
-\sum_{j_1=p+1}^{\infty}\sum_{j_2=p+1}^{\infty}
\sum_{j_3=p+1}^{\infty}
C_{j_3 j_3 j_1 j_2 j_2 j_1}=
-\sum_{j_2=p+1}^{\infty}\sum_{j_3=p+1}^{\infty}\sum_{j_1=p+1}^{\infty}
C_{j_2 j_3 j_1 j_2 j_2 j_1}=
$$

\vspace{2mm}
$$
=
\sum_{j_1=0}^{p}\sum_{j_2=p+1}^{\infty}\sum_{j_3=p+1}^{\infty}
C_{j_2 j_3 j_1 j_2 j_2 j_1}=
\sum_{j_1=0}^{p}\sum_{j_3=p+1}^{\infty}\sum_{j_2=p+1}^{\infty}
C_{j_2 j_3 j_1 j_2 j_2 j_1}=
$$

\vspace{2mm}
\begin{equation}
\label{sixsix124}
=\frac{1}{2}
\sum_{j_1=0}^{p}\sum_{j_3=p+1}^{\infty}
C_{j_3 j_3 j_1 j_2 j_2 j_1}\biggl|_{(j_2 j_2)\curvearrowright (\cdot)}\biggr. -
\sum_{j_1=0}^{p}\sum_{j_2=0}^{p}\sum_{j_3=p+1}^{\infty}
C_{j_3 j_3 j_1 j_2 j_2 j_1}.
\end{equation}

\vspace{5mm}

The equality
\begin{equation}
\label{sixsix125}
\lim\limits_{p\to\infty}\frac{1}{2}
\sum_{j_1=0}^{p}\sum_{j_3=p+1}^{\infty}
C_{j_3 j_3 j_1 j_2 j_2 j_1}\biggl|_{(j_2 j_2)\curvearrowright (\cdot)}\biggr. =0
\end{equation}

\vspace{4mm}
\noindent       
follows from the inequality (\ref{after9043}),
where we proceed similarly to the proof of equality (\ref{sixsix109a})
(see (\ref{sept12})).

The relation
\begin{equation}
\label{sixsix126}
\lim\limits_{p\to\infty}
\sum_{j_1=0}^{p}\sum_{j_2=0}^{p}\sum_{j_3=p+1}^{\infty}
C_{j_3 j_3 j_1 j_2 j_2 j_1}=0
\end{equation}

\vspace{4mm}
\noindent
is proved on the basis of (\ref{sixsix123}) and similarly with the proof 
of (\ref{sixsix7}).
The equalities (\ref{sixsix124})--(\ref{sixsix126}) prove
(\ref{sixsix2}). 

Let us prove (\ref{sixsix5}). Using (\ref{after80xx}) and (\ref{after500}), we get

$$
\sum_{j_1=p+1}^{\infty}\sum_{j_2=p+1}^{\infty}
\sum_{j_3=p+1}^{\infty}
C_{j_2 j_1 j_3 j_3 j_2 j_1}=
\sum_{j_3=p+1}^{\infty}\sum_{j_1, j_2= 0}^{p}
C_{j_2 j_1 j_3 j_3 j_2 j_1}=
$$

\vspace{2mm}
\begin{equation}
\label{sixsix127}
=\frac{1}{2}
\sum_{j_1,j_2=0}^{p}
C_{j_2 j_1 j_3 j_3 j_2 j_1}\biggl|_{(j_3 j_3)\curvearrowright (\cdot)}\biggr.
-
\sum_{j_1,j_2,j_3=0}^{p}
C_{j_2 j_1 j_3 j_3 j_2 j_1}.
\end{equation}

\vspace{5mm}

Using the equality (\ref{after2508}) we have

\begin{equation}
\label{sixsix128}
\lim\limits_{p\to\infty}
\frac{1}{2}
\sum_{j_1,j_2=0}^{p}
C_{j_2 j_1 j_3 j_3 j_2 j_1}\biggl|_{(j_3 j_3)\curvearrowright (\cdot)}\biggr.
=0,
\end{equation}

\vspace{4mm}
\noindent
where we proceed similarly to the proof of equality (\ref{sixsix109a})
(see (\ref{sept12})).

Further, we will prove the following relation

\begin{equation}
\label{sixsix129}
\lim\limits_{p\to\infty}
\sum_{j_1,j_2,j_3=0}^{p}
C_{j_2 j_1 j_3 j_3 j_2 j_1}=0
\end{equation}

\vspace{4mm}
\noindent
using the equality (\ref{sixsix40}). From (\ref{sixsix40}) we have

$$
\sum_{j_1,j_2,j_3=0}^{p}
C_{j_2 j_1 j_3 j_3 j_2 j_1}
=\frac{1}{2}\sum_{j_1,j_2,j_3=0}^{p}
\biggl(C_{j_2 j_1 j_3 j_3 j_2 j_1}+
C_{j_1 j_2 j_3 j_3 j_1 j_2}\biggr)=
$$

\vspace{2mm}
$$
=
\frac{1}{2}\sum_{j_1,j_2,j_3=0}^{p}\biggl(
C_{j_2}C_{j_1 j_3 j_3 j_2 j_1}-C_{j_1 j_2}C_{j_3 j_3 j_2 j_1}
+C_{j_3 j_1 j_2}C_{j_3 j_2 j_1}-\biggr.
$$

\vspace{2mm}
$$
\biggl.-C_{j_3 j_3 j_1 j_2}C_{j_2 j_1}+
C_{j_2 j_3 j_3 j_1 j_2}C_{j_1}\biggr)=
$$

\vspace{2mm}
$$
=
\sum_{j_1,j_2,j_3=0}^{p}\biggl(
C_{j_2 j_3 j_3 j_1 j_2}C_{j_1}-C_{j_3 j_3 j_1 j_2}C_{j_2 j_1}\biggr)+
$$

\vspace{2mm}
\begin{equation}
\label{sixsix130}
+
\frac{1}{2}\sum_{j_1,j_2,j_3=0}^{p}
C_{j_3 j_1 j_2}C_{j_3 j_2 j_1}.
\end{equation}

\vspace{5mm}

The generalized Parseval equality gives (by analogy with (\ref{pars100}))

\begin{equation}
\label{sixsix131}
\lim\limits_{p\to\infty}\frac{1}{2}\sum_{j_1,j_2,j_3=0}^{p}
C_{j_3 j_1 j_2}C_{j_3 j_2 j_1}=0.
\end{equation}

\vspace{4mm}

Let us prove the following equality

\vspace{-1mm}
\begin{equation}
\label{sixsix132}
\lim\limits_{p\to\infty}\sum_{j_1,j_2,j_3=0}^{p}\biggl(
C_{j_2 j_3 j_3 j_1 j_2}C_{j_1}-C_{j_3 j_3 j_1 j_2}C_{j_2 j_1}\biggr)=0.
\end{equation}

\vspace{4mm}

The relation
\begin{equation}
\label{sixsix132s}
\lim\limits_{p\to\infty}\sum_{j_1,j_2,j_3=0}^{p}
C_{j_2 j_3 j_3 j_1 j_2}C_{j_1}=0
\end{equation}

\vspace{4mm}
\noindent
is proved by the same methods as 
in the proof of equality (\ref{sixsix8}) and also using 
Theorem~17 and (\ref{after500}).

Further, we have (see (\ref{after500}))

\begin{equation}
\label{sept2}
\sum_{j_3=0}^{p}
C_{j_3 j_3 j_1 j_2}=
\frac{1}{2}
C_{j_3 j_3 j_1 j_2}\biggl|_{(j_3 j_3)\curvearrowright (\cdot)}\biggr.-
\sum_{j_3=p+1}^{\infty}
C_{j_3 j_3 j_1 j_2}.
\end{equation}

\vspace{4mm}

Moreover,
$$
C_{j_3 j_3 j_1 j_2}\biggl|_{(j_3 j_3)\curvearrowright (\cdot)}\biggr.=
\int\limits_t^T \int\limits_t^{t_3}\phi_{j_1}(t_2)
\int\limits_t^{t_2}\phi_{j_2}(t_1)dt_1 dt_2 dt_3=
$$

\vspace{2mm}
$$
=
\int\limits_t^T\phi_{j_1}(t_2)
\int\limits_t^{t_2}\phi_{j_2}(t_1)dt_1 \int\limits_{t_2}^T
dt_3 dt_2=
\int\limits_t^T(T-t_2)\phi_{j_1}(t_2)
\int\limits_t^{t_2}\phi_{j_2}(t_1)dt_1 
dt_2=
$$

\vspace{2mm}
$$
=\int\limits_t^T\phi_{j_2}(t_1)
\int\limits_{t_1}^T(T-t_2)\phi_{j_1}(t_2)dt_2 dt_1= 
\int\limits_t^T\phi_{j_2}(t_2)
\int\limits_{t_2}^T(T-t_1)\phi_{j_1}(t_1)dt_1 dt_2= 
$$

\vspace{2mm}
\begin{equation}
\label{sept3}
=\int\limits_{[t,T]^2}
(T-t_1){\bf 1}_{\{t_2<t_1\}}\phi_{j_1}(t_1)\phi_{j_2}(t_2)dt_1 dt_2
\stackrel{\sf def}{=}
\tilde C_{j_2 j_1}.
\end{equation}

\vspace{5mm}

Using (\ref{sept2}), (\ref{sept3}), and the generalized Parseval equality, we obtain 

$$
\lim\limits_{p\to\infty}\sum_{j_1,j_2,j_3=0}^{p}
C_{j_3 j_3 j_1 j_2}C_{j_2 j_1}=\frac{1}{2}\lim\limits_{p\to\infty}\sum_{j_1,j_2=0}^{p}
C_{j_2 j_1}\tilde C_{j_2 j_1}-
$$

\vspace{2mm}
\begin{equation}
\label{sept7}
-\lim\limits_{p\to\infty}\sum_{j_1,j_2=0}^{p}\sum_{j_3=p+1}^{\infty}
C_{j_3 j_3 j_1 j_2}C_{j_2 j_1}=
-\lim\limits_{p\to\infty}\sum_{j_1,j_2=0}^{p}\sum_{j_3=p+1}^{\infty}
C_{j_3 j_3 j_1 j_2}C_{j_2 j_1}.
\end{equation}

\vspace{5mm}

We have (see (\ref{sept10}))
\begin{equation}
\label{sept5}
C_{j_3 j_3 j_1 j_2}=
\frac{1}{2}\int\limits_t^T \phi_{j_2}(t_1)
\int\limits_{t_1}^T \phi_{j_1}(t_2)
\left(\int\limits_{t_2}^T \phi_{j_3}(t_3)
dt_3\right)^2 dt_2 dt_1.
\end{equation}

\vspace{4mm}
                        
By analogy with (\ref{sixsix74}) and also using (\ref{sept5}), we get

\begin{equation}
\label{sept6}
\lim\limits_{p\to\infty}\sum_{j_1,j_2=0}^{p}\sum_{j_3=p+1}^{\infty}
C_{j_3 j_3 j_1 j_2}C_{j_2 j_1}=0.
\end{equation}

\vspace{4mm}

Combining (\ref{sept7}) and (\ref{sept6}), we obtain

\vspace{-1mm}
\begin{equation}
\label{sept9}
\lim\limits_{p\to\infty}\sum_{j_1,j_2,j_3=0}^{p}
C_{j_3 j_3 j_1 j_2}C_{j_2 j_1}=0.
\end{equation}

\vspace{4mm}

The relation (\ref{sixsix132}) follows from (\ref{sixsix132s}) and (\ref{sept9}).
From (\ref{sixsix130})--(\ref{sixsix132}) we get (\ref{sixsix129}).
The equalities (\ref{sixsix127})--(\ref{sixsix129})
complete the proof of (\ref{sixsix5}).

For the proof of (\ref{sixsix12})--(\ref{sixsix15})
we will use a new idea.  
More precisely, we will consider the sums of 
expressions (\ref{sixsix12})--(\ref{sixsix15}) with the expressions 
already studied throughout this proof.

Let us begin from (\ref{sixsix12}). 
Applying the integration order replacement, we obtain

$$
C_{j_3 j_1 j_2 j_3 j_2 j_1}+C_{j_3 j_1 j_2 j_3 j_1 j_2}=
$$

\vspace{1mm}
$$
=
\int\limits_t^T 
\phi_{j_3}(t_6)
\int\limits_t^{t_6} 
\phi_{j_1}(t_5)
\int\limits_{t}^{t_5} 
\phi_{j_2}(t_4) 
\int\limits_t^{t_4}
\phi_{j_3}(t_3) 
\left(\int\limits_t^{t_3}
\phi_{j_2}(t_2)dt_2\right)
\left(\int\limits_t^{t_3}
\phi_{j_1}(t_1)
dt_1\right)
dt_3 dt_4 dt_5 dt_6=
$$

\vspace{2mm}
$$
=
\int\limits_t^T 
\phi_{j_3}(t_6)
\int\limits_t^{t_6} 
\phi_{j_1}(t_5)
\int\limits_{t}^{t_5} 
\phi_{j_3}(t_3) 
\left(\int\limits_t^{t_3}
\phi_{j_2}(t_2)dt_2\right)
\left(\int\limits_t^{t_3}
\phi_{j_1}(t_1)
dt_1\right)
\int\limits_{t_3}^{t_5}
\phi_{j_2}(t_4) 
dt_4
dt_3 dt_5 dt_6=
$$

\vspace{2mm}
$$
=
\int\limits_t^T 
\phi_{j_3}(t_6)
\int\limits_t^{t_6} 
\phi_{j_1}(t_5)
\left(\int\limits_{t}^{t_5}
\phi_{j_2}(t_4) 
dt_4
\right)
\int\limits_{t}^{t_5} 
\phi_{j_3}(t_3) 
\left(\int\limits_t^{t_3}
\phi_{j_2}(t_2)dt_2\right)
\left(\int\limits_t^{t_3}
\phi_{j_1}(t_1)
dt_1\right) 
dt_3 dt_5 dt_6-
$$

\vspace{2mm}
$$
-
\int\limits_t^T 
\phi_{j_3}(t_6)
\int\limits_t^{t_6} 
\phi_{j_1}(t_5)
\int\limits_{t}^{t_5} 
\phi_{j_3}(t_3) 
\left(\int\limits_t^{t_3}
\phi_{j_2}(t_2)dt_2\right)^2
\left(\int\limits_t^{t_3}
\phi_{j_1}(t_1)
dt_1\right)
dt_3 dt_5 dt_6=
$$

\vspace{2mm}
$$
=
\int\limits_t^T 
\phi_{j_1}(t_5)
\left(\int\limits_{t}^{t_5}
\phi_{j_2}(t_4) 
dt_4
\right)
\int\limits_{t}^{t_5} 
\phi_{j_3}(t_3) 
\left(\int\limits_t^{t_3}
\phi_{j_2}(t_2)dt_2\right)
\left(\int\limits_t^{t_3}
\phi_{j_1}(t_1)
dt_1\right) 
dt_3
\left(\int\limits_{t_5}^T
\phi_{j_3}(t_6)dt_6\right) dt_5-
$$

\vspace{2mm}
\begin{equation}
\label{sixsix133}
-
\int\limits_t^T 
\phi_{j_1}(t_5)
\int\limits_{t}^{t_5} 
\phi_{j_3}(t_3) 
\left(\int\limits_t^{t_3}
\phi_{j_2}(t_2)dt_2\right)^2
\left(\int\limits_t^{t_3}
\phi_{j_1}(t_1)
dt_1\right)
dt_3
\left(
\int\limits_{t_5}^T
\phi_{j_3}(t_6)
dt_6\right)dt_5.
\end{equation}

\vspace{5mm}

Using (\ref{after80xx}), we get

$$
\sum_{j_1=p+1}^{\infty}\sum_{j_2=p+1}^{\infty}
\sum_{j_3=p+1}^{\infty}
\biggl(
C_{j_3 j_1 j_2 j_3 j_2 j_1}+C_{j_3 j_1 j_2 j_3 j_1 j_2}\biggr)=
$$

\vspace{2mm}
\begin{equation}
\label{sixsix134}
=
\sum_{j_1=0}^{p}\sum_{j_3=0}^{p}
\sum_{j_2=p+1}^{\infty}
\biggl(
C_{j_3 j_1 j_2 j_3 j_2 j_1}+C_{j_3 j_1 j_2 j_3 j_1 j_2}\biggr).
\end{equation}

\vspace{5mm}

Further, by analogy with the proof of equality (\ref{sixsix7}) 
and using (\ref{sixsix133}), we obtain

\begin{equation}
\label{sixsix135}
\lim\limits_{p\to\infty}\sum_{j_1=0}^{p}\sum_{j_3=0}^{p}
\sum_{j_2=p+1}^{\infty}
\biggl(
C_{j_3 j_1 j_2 j_3 j_2 j_1}+C_{j_3 j_1 j_2 j_3 j_1 j_2}\biggr)=0.
\end{equation}

\vspace{4mm}

From (\ref{sixsix134}) and (\ref{sixsix135}) we get

\begin{equation}
\label{sixsix136}
\lim\limits_{p\to\infty}\sum_{j_1=p+1}^{\infty}\sum_{j_2=p+1}^{\infty}
\sum_{j_3=p+1}^{\infty}
\biggl(
C_{j_3 j_1 j_2 j_3 j_2 j_1}+C_{j_3 j_1 j_2 j_3 j_1 j_2}\biggr)=0.
\end{equation}

\vspace{4mm}

Moreover (see (\ref{sixsix8})),
\begin{equation}
\label{sixsix137}
\lim\limits_{p\to\infty}\sum_{j_1=p+1}^{\infty}\sum_{j_2=p+1}^{\infty}
\sum_{j_3=p+1}^{\infty}
C_{j_3 j_1 j_2 j_3 j_1 j_2}=0.
\end{equation}

\vspace{4mm}

Combining (\ref{sixsix136}) and (\ref{sixsix137}), we have

$$
\lim\limits_{p\to\infty}\sum_{j_1=p+1}^{\infty}\sum_{j_2=p+1}^{\infty}
\sum_{j_3=p+1}^{\infty}
C_{j_3 j_1 j_2 j_3 j_2 j_1}=0.
$$

\vspace{4mm}
\noindent
The equality (\ref{sixsix12}) is proved.

Consider (\ref{sixsix11}).
Using the integration order replacement, we have

$$
C_{j_2 j_3 j_1 j_3 j_2 j_1}+C_{j_2 j_3 j_1 j_3 j_1 j_2}=
$$

\vspace{1mm}
$$
=
\int\limits_t^T 
\phi_{j_2}(t_6)
\int\limits_t^{t_6} 
\phi_{j_3}(t_5)
\int\limits_{t}^{t_5} 
\phi_{j_1}(t_4) 
\int\limits_t^{t_4}
\phi_{j_3}(t_3) 
\left(\int\limits_t^{t_3}
\phi_{j_2}(t_2)dt_2\right)
\left(\int\limits_t^{t_3}
\phi_{j_1}(t_1)
dt_1\right)
dt_3 dt_4 dt_5 dt_6=
$$

\vspace{2mm}
$$
=
\int\limits_t^T 
\phi_{j_2}(t_6)
\int\limits_t^{t_6} 
\phi_{j_3}(t_5)
\int\limits_{t}^{t_5} 
\phi_{j_3}(t_3) 
\left(\int\limits_t^{t_3}
\phi_{j_2}(t_2)dt_2\right)
\left(\int\limits_t^{t_3}
\phi_{j_1}(t_1)
dt_1\right)
\int\limits_{t_3}^{t_5}
\phi_{j_1}(t_4) 
dt_4
dt_3 dt_5 dt_6=
$$

$$
=
\int\limits_t^T 
\phi_{j_2}(t_6)
\int\limits_t^{t_6} 
\phi_{j_3}(t_5)
\left(\int\limits_{t}^{t_5}
\phi_{j_1}(t_4) 
dt_4
\right)
\int\limits_{t}^{t_5} 
\phi_{j_3}(t_3) 
\left(\int\limits_t^{t_3}
\phi_{j_2}(t_2)dt_2\right)
\left(\int\limits_t^{t_3}
\phi_{j_1}(t_1)
dt_1\right) 
dt_3 dt_5 dt_6-
$$

\vspace{2mm}
$$
-
\int\limits_t^T 
\phi_{j_2}(t_6)
\int\limits_t^{t_6} 
\phi_{j_3}(t_5)
\int\limits_{t}^{t_5} 
\phi_{j_3}(t_3) 
\left(\int\limits_t^{t_3}
\phi_{j_2}(t_2)dt_2\right)
\left(\int\limits_t^{t_3}
\phi_{j_1}(t_1)
dt_1\right)^2
dt_3 dt_5 dt_6=
$$

\vspace{2mm}
$$
=
\int\limits_t^T 
\phi_{j_3}(t_5)
\left(\int\limits_{t}^{t_5}
\phi_{j_1}(t_4) 
dt_4
\right)
\int\limits_{t}^{t_5} 
\phi_{j_3}(t_3) 
\left(\int\limits_t^{t_3}
\phi_{j_2}(t_2)dt_2\right)
\left(\int\limits_t^{t_3}
\phi_{j_1}(t_1)
dt_1\right) 
dt_3
\left(\int\limits_{t_5}^T
\phi_{j_2}(t_6)dt_6\right) dt_5-
$$

\vspace{2mm}
\begin{equation}
\label{sixsix150}
-
\int\limits_t^T 
\phi_{j_3}(t_5)
\int\limits_{t}^{t_5} 
\phi_{j_3}(t_3) 
\left(\int\limits_t^{t_3}
\phi_{j_2}(t_2)dt_2\right)
\left(\int\limits_t^{t_3}
\phi_{j_1}(t_1)
dt_1\right)^2
dt_3\left(
\int\limits_{t_5}^T
\phi_{j_2}(t_6)
dt_6\right)dt_5.
\end{equation}

\vspace{5mm}

Using (\ref{after80xx}), we obtain

$$
-\sum_{j_1=p+1}^{\infty}\sum_{j_2=p+1}^{\infty}
\sum_{j_3=p+1}^{\infty}
\biggl(
C_{j_2 j_3 j_1 j_3 j_2 j_1}+C_{j_2 j_3 j_1 j_3 j_1 j_2}\biggr)=
$$

\vspace{2mm}
\begin{equation}
\label{sixsix151}
=
\sum_{j_3=0}^{p}
\sum_{j_1=p+1}^{\infty}\sum_{j_2=p+1}^{\infty}
\biggl(
C_{j_2 j_3 j_1 j_3 j_2 j_1}+C_{j_2 j_3 j_1 j_3 j_1 j_2}\biggr).
\end{equation}

\vspace{5mm}

By analogy with the proof of (\ref{sixsix7}) 
and applying (\ref{sixsix150}), we get

\begin{equation}
\label{sixsix152}
\lim\limits_{p\to\infty}\sum_{j_3=0}^{p}
\sum_{j_1=p+1}^{\infty}\sum_{j_2=p+1}^{\infty}
\biggl(
C_{j_2 j_3 j_1 j_3 j_2 j_1}+C_{j_2 j_3 j_1 j_3 j_1 j_2}\biggr)=0.
\end{equation}

\vspace{4mm}

From (\ref{sixsix151}) and (\ref{sixsix152}) we have

\begin{equation}
\label{sixsix153}
\lim\limits_{p\to\infty}
\sum_{j_1=p+1}^{\infty}\sum_{j_2=p+1}^{\infty}
\sum_{j_3=p+1}^{\infty}
\biggl(
C_{j_2 j_3 j_1 j_3 j_2 j_1}+C_{j_2 j_3 j_1 j_3 j_1 j_2}\biggr)=0.
\end{equation}

\vspace{4mm}

Moreover (see (\ref{sixsix9})),
\begin{equation}
\label{sixsix154}
\lim\limits_{p\to\infty}\sum_{j_1=p+1}^{\infty}\sum_{j_2=p+1}^{\infty}
\sum_{j_3=p+1}^{\infty}
C_{j_2 j_3 j_1 j_3 j_1 j_2}=0.
\end{equation}

\vspace{4mm}

Combining (\ref{sixsix153}) and (\ref{sixsix154}), we finally obtain

$$
\lim\limits_{p\to\infty}\sum_{j_1=p+1}^{\infty}\sum_{j_2=p+1}^{\infty}
\sum_{j_3=p+1}^{\infty}
C_{j_2 j_3 j_1 j_3 j_2 j_1}=0.
$$

\vspace{4mm}
\noindent
The equality (\ref{sixsix11}) is proved.

Now consider (\ref{sixsix13}).
Using the integration order replacement, we obtain

$$
C_{j_3 j_1 j_3 j_2 j_2 j_1}+C_{j_3 j_1 j_3 j_2 j_1 j_2}=
$$

\vspace{1mm}
$$
=
\int\limits_t^T 
\phi_{j_3}(t_6)
\int\limits_t^{t_6} 
\phi_{j_1}(t_5)
\int\limits_{t}^{t_5} 
\phi_{j_3}(t_4) 
\int\limits_t^{t_4}
\phi_{j_2}(t_3) 
\left(\int\limits_t^{t_3}
\phi_{j_2}(t_2)dt_2\right)
\left(\int\limits_t^{t_3}
\phi_{j_1}(t_1)
dt_1\right)
dt_3 dt_4 dt_5 dt_6=
$$

\vspace{2mm}
$$
=
\int\limits_t^T 
\phi_{j_3}(t_6)
\int\limits_t^{t_6} 
\phi_{j_1}(t_5)
\int\limits_{t}^{t_5} 
\phi_{j_2}(t_3) 
\left(\int\limits_t^{t_3}
\phi_{j_2}(t_2)dt_2\right)
\left(\int\limits_t^{t_3}
\phi_{j_1}(t_1)
dt_1\right)
\int\limits_{t_3}^{t_5}
\phi_{j_3}(t_4) 
dt_4
dt_3 dt_5 dt_6=
$$

\vspace{2mm}
$$
=
\int\limits_t^T 
\phi_{j_3}(t_6)
\int\limits_t^{t_6} 
\phi_{j_1}(t_5)
\left(\int\limits_{t}^{t_5}
\phi_{j_3}(t_4) 
dt_4
\right)
\int\limits_{t}^{t_5} 
\phi_{j_2}(t_3) 
\left(\int\limits_t^{t_3}
\phi_{j_2}(t_2)dt_2\right)
\left(\int\limits_t^{t_3}
\phi_{j_1}(t_1)
dt_1\right) 
dt_3 dt_5 dt_6-
$$

\vspace{2mm}
$$
-
\int\limits_t^T 
\phi_{j_3}(t_6)
\int\limits_t^{t_6} 
\phi_{j_1}(t_5)
\int\limits_{t}^{t_5} 
\phi_{j_2}(t_3) 
\left(\int\limits_t^{t_3}
\phi_{j_2}(t_2)dt_2\right)
\left(\int\limits_t^{t_3}
\phi_{j_1}(t_1)
dt_1\right)
\left(\int\limits_t^{t_3}
\phi_{j_3}(t_4)
dt_4\right)
dt_3 dt_5 dt_6=
$$

\vspace{2mm}
$$
=
\int\limits_t^T 
\phi_{j_1}(t_5)
\left(\int\limits_{t}^{t_5}
\phi_{j_3}(t_4) 
dt_4
\right)
\int\limits_{t}^{t_5} 
\phi_{j_2}(t_3) 
\left(\int\limits_t^{t_3}
\phi_{j_2}(t_2)dt_2\right)
\left(\int\limits_t^{t_3}
\phi_{j_1}(t_1)
dt_1\right) 
dt_3
\left(\int\limits_{t_5}^T
\phi_{j_3}(t_6)dt_6\right) dt_5-
$$

\vspace{2mm}
\begin{equation}
\label{sixsix160}
-
\int\limits_t^T 
\phi_{j_1}(t_5)
\int\limits_{t}^{t_5} 
\phi_{j_2}(t_3) 
\left(\int\limits_t^{t_3}
\phi_{j_2}(t_2)dt_2\right)
\left(\int\limits_t^{t_3}
\phi_{j_1}(t_1)
dt_1\right)
\left(\int\limits_t^{t_3}
\phi_{j_3}(t_4)
dt_4\right)dt_3
\left(
\int\limits_{t_5}^T
\phi_{j_3}(t_6)
dt_6\right)dt_5.
\end{equation}

\vspace{5mm}

Applying (\ref{after80xx}) and (\ref{after500}), we obtain

$$
\sum_{j_1=p+1}^{\infty}\sum_{j_2=p+1}^{\infty}
\sum_{j_3=p+1}^{\infty}
\biggl(C_{j_3 j_1 j_3 j_2 j_2 j_1}+C_{j_3 j_1 j_3 j_2 j_1 j_2}\biggr)=
$$

\vspace{2mm}
$$
=-\sum_{j_1=0}^{p}\sum_{j_3=p+1}^{\infty}\sum_{j_2=p+1}^{\infty}
\biggl(C_{j_3 j_1 j_3 j_2 j_2 j_1}+C_{j_3 j_1 j_3 j_2 j_1 j_2}\biggr)=
$$

\vspace{2mm}
$$
=
\sum_{j_1=0}^{p}\sum_{j_2=0}^{p}\sum_{j_3=p+1}^{\infty}
\biggl(C_{j_3 j_1 j_3 j_2 j_2 j_1}+C_{j_3 j_1 j_3 j_2 j_1 j_2}\biggr)-
$$

\vspace{2mm}
\begin{equation}
\label{sixsix170}
-\frac{1}{2}
\sum_{j_1=0}^{p}\sum_{j_3=p+1}^{\infty}
C_{j_3 j_1 j_3 j_2 j_2 j_1}\biggl|_{(j_2 j_2)\curvearrowright (\cdot)}\biggr..
\end{equation}

\vspace{5mm}

The equality
\begin{equation}
\label{sixsix171}
\lim\limits_{p\to\infty}\frac{1}{2}
\sum_{j_1=0}^{p}\sum_{j_3=p+1}^{\infty}
C_{j_3 j_1 j_3 j_2 j_2 j_1}\biggl|_{(j_2 j_2)\curvearrowright (\cdot)}\biggr. =0
\end{equation}

\vspace{4mm}
\noindent
follows from the 
equality (\ref{after2508}),
where we proceed similarly to the proof of equality (\ref{sixsix109a})
(see (\ref{sept12})).

By analogy with the proof of (\ref{sixsix7}) 
and applying (\ref{sixsix160}), we get

\begin{equation}
\label{sixsix171s}
\lim\limits_{p\to\infty}\sum_{j_1=0}^{p}\sum_{j_2=0}^{p}\sum_{j_3=p+1}^{\infty}
\biggl(C_{j_3 j_1 j_3 j_2 j_2 j_1}+C_{j_3 j_1 j_3 j_2 j_1 j_2}\biggr)=0.
\end{equation}

\vspace{4mm}

From (\ref{sixsix170})--(\ref{sixsix171s}) we have

\begin{equation}
\label{sixsix300}
\lim\limits_{p\to\infty}
\sum_{j_1=p+1}^{\infty}\sum_{j_2=p+1}^{\infty}
\sum_{j_3=p+1}^{\infty}
\biggl(C_{j_3 j_1 j_3 j_2 j_2 j_1}+C_{j_3 j_1 j_3 j_2 j_1 j_2}\biggr)=0.
\end{equation}

\vspace{4mm}

Moreover (see (\ref{sixsix10})),
\begin{equation}
\label{sixsix301}
\lim\limits_{p\to\infty}
\sum_{j_1=p+1}^{\infty}\sum_{j_2=p+1}^{\infty}
\sum_{j_3=p+1}^{\infty}
C_{j_3 j_1 j_3 j_2 j_1 j_2}=0.
\end{equation}

\vspace{4mm}

Combining (\ref{sixsix300}) and (\ref{sixsix301}), we finally obtain

$$
\lim\limits_{p\to\infty}
\sum_{j_1=p+1}^{\infty}\sum_{j_2=p+1}^{\infty}
\sum_{j_3=p+1}^{\infty}
C_{j_3 j_1 j_3 j_2 j_2 j_1}=0.
$$

\vspace{4mm}
\noindent
The equality (\ref{sixsix13}) is proved.

Finally consider (\ref{sixsix15}).
Using the integration order replacement, we have

$$
C_{j_2 j_3 j_3 j_1 j_2 j_1}+C_{j_2 j_3 j_3 j_1 j_1 j_2}=
$$

\vspace{1mm}
$$
=
\int\limits_t^T 
\phi_{j_2}(t_6)
\int\limits_t^{t_6} 
\phi_{j_3}(t_5)
\int\limits_{t}^{t_5} 
\phi_{j_3}(t_4) 
\int\limits_t^{t_4}
\phi_{j_1}(t_3) 
\left(\int\limits_t^{t_3}
\phi_{j_2}(t_2)dt_2\right)
\left(\int\limits_t^{t_3}
\phi_{j_1}(t_1)
dt_1\right)
dt_3 dt_4 dt_5 dt_6=
$$

\vspace{2mm}
$$
=
\int\limits_t^T 
\phi_{j_2}(t_6)
\int\limits_t^{t_6} 
\phi_{j_3}(t_5)
\int\limits_{t}^{t_5} 
\phi_{j_1}(t_3) 
\left(\int\limits_t^{t_3}
\phi_{j_2}(t_2)dt_2\right)
\left(\int\limits_t^{t_3}
\phi_{j_1}(t_1)
dt_1\right)
\int\limits_{t_3}^{t_5}
\phi_{j_3}(t_4) 
dt_4
dt_3 dt_5 dt_6=
$$

\vspace{2mm}
$$
=
\int\limits_t^T 
\phi_{j_2}(t_6)
\int\limits_t^{t_6} 
\phi_{j_3}(t_5)
\left(\int\limits_{t}^{t_5}
\phi_{j_3}(t_4) 
dt_4
\right)
\int\limits_{t}^{t_5} 
\phi_{j_1}(t_3) 
\left(\int\limits_t^{t_3}
\phi_{j_2}(t_2)dt_2\right)
\left(\int\limits_t^{t_3}
\phi_{j_1}(t_1)
dt_1\right) 
dt_3 dt_5 dt_6-
$$

\vspace{2mm}
$$
-
\int\limits_t^T 
\phi_{j_2}(t_6)
\int\limits_t^{t_6} 
\phi_{j_3}(t_5)
\int\limits_{t}^{t_5} 
\phi_{j_1}(t_3) 
\left(\int\limits_t^{t_3}
\phi_{j_2}(t_2)dt_2\right)
\left(\int\limits_t^{t_3}
\phi_{j_1}(t_1)
dt_1\right)
\left(\int\limits_t^{t_3}
\phi_{j_3}(t_4)
dt_4\right)
dt_3 dt_5 dt_6=
$$

\vspace{2mm}
$$
=
\int\limits_t^T 
\phi_{j_3}(t_5)
\left(\int\limits_{t}^{t_5}
\phi_{j_3}(t_4) 
dt_4
\right)
\int\limits_{t}^{t_5} 
\phi_{j_1}(t_3) 
\left(\int\limits_t^{t_3}
\phi_{j_2}(t_2)dt_2\right)
\left(\int\limits_t^{t_3}
\phi_{j_1}(t_1)
dt_1\right) 
dt_3
\left(\int\limits_{t_5}^T
\phi_{j_2}(t_6)dt_6\right) dt_5-
$$

\vspace{2mm}
\begin{equation}
\label{sixsix190}
-
\int\limits_t^T 
\phi_{j_3}(t_5)
\int\limits_{t}^{t_5} 
\phi_{j_1}(t_3) 
\left(\int\limits_t^{t_3}
\phi_{j_2}(t_2)dt_2\right)
\left(\int\limits_t^{t_3}
\phi_{j_1}(t_1)
dt_1\right)
\left(\int\limits_t^{t_3}
\phi_{j_3}(t_4)
dt_4\right)dt_3
\left(
\int\limits_{t_5}^T
\phi_{j_2}(t_6)
dt_6\right)dt_5.
\end{equation}

\vspace{5mm}

Using (\ref{after80xx}) and (\ref{after500}), we get

$$
\sum_{j_1=p+1}^{\infty}\sum_{j_2=p+1}^{\infty}
\sum_{j_3=p+1}^{\infty}
\biggl(C_{j_2 j_3 j_3 j_1 j_2 j_1}+
C_{j_2 j_3 j_3 j_1 j_1 j_2}\biggr)=
$$

\vspace{2mm}
$$
=\frac{1}{2}
\sum_{j_1=p+1}^{\infty}\sum_{j_2=p+1}^{\infty}
\left(C_{j_2 j_3 j_3 j_1 j_2 j_1}\biggl|_{(j_3 j_3)\curvearrowright (\cdot)}\biggr. +
C_{j_2 j_3 j_3 j_1 j_1 j_2}\biggl|_{(j_3 j_3)\curvearrowright (\cdot)}\biggr.\right)-
$$

\vspace{2mm}
$$
-\sum_{j_3=0}^{p}\sum_{j_1=p+1}^{\infty}
\sum_{j_2=p+1}^{\infty}
\biggl(C_{j_2 j_3 j_3 j_1 j_2 j_1}+
C_{j_2 j_3 j_3 j_1 j_1 j_2}\biggr)=
$$

\vspace{2mm}
$$
=\frac{1}{2}
\sum_{j_1=p+1}^{\infty}\sum_{j_2=p+1}^{\infty}
\left(C_{j_2 j_3 j_3 j_1 j_2 j_1}\biggl|_{(j_3 j_3)\curvearrowright (\cdot)}\biggr. +
C_{j_2 j_3 j_3 j_1 j_1 j_2}\biggl|_{(j_3 j_3)\curvearrowright (\cdot)}\biggr.\right)+
$$

\vspace{2mm}
$$
+\sum_{j_1=0}^{p}\sum_{j_3=0}^{p}
\sum_{j_2=p+1}^{\infty}
\biggl(C_{j_2 j_3 j_3 j_1 j_2 j_1}+
C_{j_2 j_3 j_3 j_1 j_1 j_2}\biggr)-
$$

\vspace{2mm}
\begin{equation}
\label{sixsix191}
-\frac{1}{2}
\sum_{j_3=0}^{p}\sum_{j_2=p+1}^{\infty}
C_{j_2 j_3 j_3 j_1 j_1 j_2}\biggl|_{(j_1 j_1)\curvearrowright (\cdot)}\biggr..
\end{equation}

\vspace{5mm}

The equalities

\vspace{-1mm}
\begin{equation}
\label{sixsix192}
\lim\limits_{p\to\infty}\frac{1}{2}
\sum_{j_1=p+1}^{\infty}\sum_{j_2=p+1}^{\infty}
\left(C_{j_2 j_3 j_3 j_1 j_2 j_1}\biggl|_{(j_3 j_3)\curvearrowright (\cdot)}\biggr. +
C_{j_2 j_3 j_3 j_1 j_1 j_2}\biggl|_{(j_3 j_3)\curvearrowright (\cdot)}\biggr.\right)=0,
\end{equation}

\vspace{2mm}
$$
\lim\limits_{p\to\infty}\frac{1}{2}
\sum_{j_3=0}^{p}\sum_{j_2=p+1}^{\infty}
C_{j_2 j_3 j_3 j_1 j_1 j_2}\biggl|_{(j_1 j_1)\curvearrowright (\cdot)}\biggr.
=
$$

\vspace{2mm}
$$
=\lim\limits_{p\to\infty}\frac{1}{4}
\sum_{j_2=p+1}^{\infty}
C_{j_2 j_3 j_3 j_1 j_1 j_2}\biggl|_{(j_1 j_1)\curvearrowright (\cdot)
(j_3 j_3)\curvearrowright (\cdot)}\biggr.-
$$

\vspace{2mm}
\begin{equation}
\label{sixsix193}
-
\lim\limits_{p\to\infty}\frac{1}{2}
\sum_{j_3=p+1}^{\infty}\sum_{j_2=p+1}^{\infty}
C_{j_2 j_3 j_3 j_1 j_1 j_2}\biggl|_{(j_1 j_1)\curvearrowright (\cdot)}\biggr.
=0
\end{equation}

\vspace{5mm}
\noindent
follows from the 
equalities (\ref{after2508}), (\ref{after2509}),
where we used the same technique as in (\ref{sept12}).
When proving (\ref{sixsix193}), we also applied (\ref{after500})
and (\ref{5tafter}).

By analogy with the proof of (\ref{sixsix7}) and applying (\ref{sixsix190}), we obtain

\begin{equation}
\label{sixsix194}
\lim\limits_{p\to\infty}
\sum_{j_1=0}^{p}\sum_{j_3=0}^{p}
\sum_{j_2=p+1}^{\infty}
\biggl(C_{j_2 j_3 j_3 j_1 j_2 j_1}+
C_{j_2 j_3 j_3 j_1 j_1 j_2}\biggr)=0.
\end{equation}

\vspace{4mm}

From (\ref{sixsix191})--(\ref{sixsix194}) we have

\begin{equation}
\label{sixsix194a}
\lim\limits_{p\to\infty}
\sum_{j_1=p+1}^{\infty}\sum_{j_2=p+1}^{\infty}
\sum_{j_3=p+1}^{\infty}
\biggl(C_{j_2 j_3 j_3 j_1 j_2 j_1}+
C_{j_2 j_3 j_3 j_1 j_1 j_2}\biggr)=0.
\end{equation}

\vspace{4mm}

Furthermore (see (\ref{sixsix14})),

\vspace{-1mm}
\begin{equation}
\label{sixsix195}
\lim\limits_{p\to\infty}\sum_{j_1=p+1}^{\infty}\sum_{j_2=p+1}^{\infty}
\sum_{j_3=p+1}^{\infty}
C_{j_2 j_3 j_3 j_1 j_1 j_2}=0.
\end{equation}

\vspace{4mm}

Combining (\ref{sixsix194a}) and (\ref{sixsix195}), we finally obtain

$$
\lim\limits_{p\to\infty}\sum_{j_1=p+1}^{\infty}\sum_{j_2=p+1}^{\infty}
\sum_{j_3=p+1}^{\infty}
C_{j_2 j_3 j_3 j_1 j_2 j_1}=0.
$$

\vspace{4mm}
\noindent
The equality (\ref{sixsix15}) is proved. 
Theorem 23 is proved.

\vspace{5mm}

\section{Generalization of Theorem~16.
The Case $p_1,$ $p_2,$ $p_3\to \infty$ and Continuously
Differetiable Weight Functions (The Cases of Legendre Polynomials and 
Trigonometric Functions). Proof of Hypothesis~3 for the Case $k=3$}

\vspace{5mm}

This section is devoted to the following theorem.

\vspace{2mm}

{\bf Theorem~24}\ \cite{2018a}, \cite{arxiv-5}, \cite{arxiv-9}.\  {\it Suppose that 
$\{\phi_j(x)\}_{j=0}^{\infty}$ is a complete orthonormal system of 
Legendre polynomials or trigonometric func\-ti\-ons in the space $L_2([t, T]).$
Furthermore, let $\psi_1(\tau), \psi_2(\tau), \psi_3(\tau)$ are continuously dif\-ferentiable 
nonrandom functions on $[t, T]$.
Then, for the 
iterated Stra\-to\-no\-vich stochastic integral of third multiplicity

$$
J^{*}[\psi^{(3)}]_{T,t}^{(i_1 i_2 i_3)}={\int\limits_t^{*}}^T\psi_3(t_3)
{\int\limits_t^{*}}^{t_3}\psi_2(t_2)
{\int\limits_t^{*}}^{t_2}\psi_1(t_1)
d{\bf w}_{t_1}^{(i_1)}
d{\bf w}_{t_2}^{(i_2)}d{\bf w}_{t_3}^{(i_3)}
$$

\vspace{2mm}
\noindent
the following 
expansion 

\vspace{-3mm}
\begin{equation}
\label{20231}
J^{*}[\psi^{(3)}]_{T,t}^{(i_1 i_2 i_3)}
=\hbox{\vtop{\offinterlineskip\halign{
\hfil#\hfil\cr
{\rm l.i.m.}\cr
$\stackrel{}{{}_{p_1,p_2,p_3\to \infty}}$\cr
}} }
\sum\limits_{j_1=0}^{p_1}\sum\limits_{j_2=0}^{p_2}\sum\limits_{j_3=0}^{p_3}
C_{j_3 j_2 j_1}\zeta_{j_1}^{(i_1)}\zeta_{j_2}^{(i_2)}\zeta_{j_3}^{(i_3)}
\end{equation}

\vspace{3mm}
\noindent
that converges in the mean-square sense is valid, where
$i_1, i_2, i_3=0, 1,\ldots,m,$

$$
C_{j_3 j_2 j_1}=\int\limits_t^T\psi_3(t_3)\phi_{j_3}(t_3)
\int\limits_t^{t_3}\psi_2(t_2)\phi_{j_2}(t_2)
\int\limits_t^{t_2}\psi_1(t_1)\phi_{j_1}(t_1)dt_1dt_2dt_3
$$

\vspace{2mm}
\noindent
and
$$
\zeta_{j}^{(i)}=
\int\limits_t^T \phi_{j}(s) d{\bf w}_s^{(i)}
$$ 

\vspace{2mm}
\noindent
are independent standard Gaussian random variables for various 
$i$ or $j$
{\rm (}in the case when $i\ne 0${\rm ),}
${\bf w}_{\tau}^{(i)}={\bf f}_{\tau}^{(i)}$ for
$i=1,\ldots,m$ and 
${\bf w}_{\tau}^{(0)}=\tau.$}

\vspace{2mm}

{\bf Proof.} Let us consider the case of
Legendre polynomials (the trigonometric case is simpler
and can be considered similarly). Applying (\ref{after33}), we obtain

\vspace{0.5mm}
$$
\sum_{j_1=0}^{p_1}\sum_{j_2=0}^{p_2}\sum_{j_3=0}^{p_3}
C_{j_3j_2j_1}
\zeta_{j_1}^{(i_1)}\zeta_{j_2}^{(i_2)}\zeta_{j_3}^{(i_3)}=
J'[K_{p_1p_2p_3}]_{T,t}^{(i_1 i_2 i_3)}+
$$

\vspace{2mm}
$$
+{\bf 1}_{\{i_1=i_2\ne 0\}}
\sum_{j_3=0}^{p_3}\sum_{j_1=0}^{\min\{p_1,p_2\}}
C_{j_3j_1j_1}J'[\phi_{j_3}]^{(i_3)}_{T,t}+
$$

\vspace{2mm}
$$
+{\bf 1}_{\{i_2=i_3\ne 0\}}
\sum_{j_1=0}^{p_1}\sum_{j_3=0}^{\min\{p_2,p_3\}}
C_{j_3j_3j_1}J'[\phi_{j_1}]^{(i_1)}_{T,t}+
$$

\vspace{2mm}
\begin{equation}
\label{20232}
+{\bf 1}_{\{i_1=i_3\ne 0\}}
\sum_{j_2=0}^{p_2}\sum_{j_1=0}^{\min\{p_1,p_3\}}
C_{j_1j_2j_1}J'[\phi_{j_2}]^{(i_2)}_{T,t}
\end{equation}

\vspace{4mm}
\noindent
w.~p.~1, where notations are the same as in (\ref{after33}).

Using Theorem~5 (see (\ref{30.4}) for the case $k=3$), Theorem~3 (see (\ref{tyyyarg}))
as well as (\ref{after609}) (see the derivation
of (\ref{after609})) and (\ref{after500}), we get

$$
J^{*}[\psi^{(3)}]_{T,t}^{(i_1 i_2 i_3)}=J[\psi^{(3)}]_{T,t}^{(i_1 i_2 i_3)}
+ \frac{1}{2}{\bf 1}_{\{i_1=i_2\ne 0\}}
\int\limits_t^T\psi_3(t_3)
\int\limits_t^{t_3}\psi_2(t_2)\psi_1(t_2)dt_2
d{\bf w}_{t_3}^{(i_3)}+
$$

\vspace{2mm}
$$
+ \frac{1}{2}{\bf 1}_{\{i_2=i_3\ne 0\}}
\int\limits_t^T\psi_3(t_3)\psi_2(t_3)
\int\limits_t^{t_3}\psi_1(t_1)d{\bf w}_{t_1}^{(i_1)}
dt_3=
$$

\vspace{4mm}
$$
=J[\psi^{(3)}]_{T,t}^{(i_1 i_2 i_3)}
+ \frac{1}{2}J[\psi^{(3)}]_{T,t}^1+\frac{1}{2}J[\psi^{(3)}]_{T,t}^2
=
$$

\vspace{4mm}
$$
=
\hbox{\vtop{\offinterlineskip\halign{
\hfil#\hfil\cr
{\rm l.i.m.}\cr
$\stackrel{}{{}_{p_1,p_2,p_3\to \infty}}$\cr
}} }J'[K_{p_1p_2p_3}]_{T,t}^{(i_1 i_2 i_3)}+
$$

\vspace{3mm}
$$
+{\bf 1}_{\{i_1=i_2\ne 0\}}
\hbox{\vtop{\offinterlineskip\halign{
\hfil#\hfil\cr
{\rm l.i.m.}\cr
$\stackrel{}{{}_{p_3\to \infty}}$\cr
}} }\frac{1}{2}\sum_{j_3=0}^{p_3} 
C_{j_3 j_2 j_1}\biggl|_{(j_{2} j_{1})\curvearrowright (\cdot),
j_{1}=j_{2}}\biggr.
J'[\phi_{j_3}]^{(i_3)}_{T,t}+
$$

\vspace{4mm}
$$
+{\bf 1}_{\{i_2=i_3\ne 0\}}
\hbox{\vtop{\offinterlineskip\halign{
\hfil#\hfil\cr
{\rm l.i.m.}\cr
$\stackrel{}{{}_{p_1\to \infty}}$\cr
}} }
\frac{1}{2}\sum_{j_1=0}^{p_1} 
C_{j_3 j_2 j_1}\biggl|_{(j_{3} j_{2})\curvearrowright (\cdot),
j_{2}=j_{3}}\biggr.
J'[\phi_{j_1}]^{(i_1)}_{T,t}=
$$

\vspace{6mm}
$$
=
\hbox{\vtop{\offinterlineskip\halign{
\hfil#\hfil\cr
{\rm l.i.m.}\cr
$\stackrel{}{{}_{p_1,p_2,p_3\to \infty}}$\cr
}} }J'[K_{p_1p_2p_3}]_{T,t}^{(i_1 i_2 i_3)}+
$$

\vspace{4mm}
$$
+{\bf 1}_{\{i_1=i_2\ne 0\}}
\hbox{\vtop{\offinterlineskip\halign{
\hfil#\hfil\cr
{\rm l.i.m.}\cr
$\stackrel{}{{}_{p_3\to \infty}}$\cr
}} }\sum_{j_3=0}^{p_3}\sum_{j_1=0}^{\infty}  
C_{j_3 j_1 j_1}  
J'[\phi_{j_3}]^{(i_3)}_{T,t}+
$$

\vspace{4mm}
\begin{equation}
\label{20233}
+{\bf 1}_{\{i_2=i_3\ne 0\}}
\hbox{\vtop{\offinterlineskip\halign{
\hfil#\hfil\cr
{\rm l.i.m.}\cr
$\stackrel{}{{}_{p_1\to \infty}}$\cr
}} }
\sum_{j_1=0}^{p_1} \sum_{j_3=0}^{\infty}
C_{j_3 j_3 j_1}
J'[\phi_{j_1}]^{(i_1)}_{T,t}
\end{equation}

\vspace{4mm}
\noindent
w.~p.~1.

Using (\ref{20232}), (\ref{20233}) and the elementary inequality 

$$
(a+b+c+d)^2\le 4\left(a^2+b^2+c^2+d^2\right),
$$

\vspace{3mm}
\noindent
we obtain

$$
{\sf M}\left\{\left(J^{*}[\psi^{(3)}]_{T,t}^{(i_1 i_2 i_3)}-
\sum\limits_{j_1=0}^{p_1}\sum\limits_{j_2=0}^{p_2}\sum\limits_{j_3=0}^{p_3}
C_{j_3 j_2 j_1}\zeta_{j_1}^{(i_1)}\zeta_{j_2}^{(i_2)}\zeta_{j_3}^{(i_3)}\right)^2\right\}\le
$$

\vspace{3mm}
$$
\le 4{\sf M}\left\{\biggl(J[\psi^{(3)}]_{T,t}^{(i_1 i_2 i_3)}-
J'[K_{p_1p_2p_3}]_{T,t}^{(i_1 i_2 i_3)}
\biggr)^2\right\}+
$$

\vspace{3mm}
$$
+4\cdot {\bf 1}_{\{i_1=i_2\ne 0\}}\times
$$

$$
\times
{\sf M}\left\{\left(
\hbox{\vtop{\offinterlineskip\halign{
\hfil#\hfil\cr
{\rm l.i.m.}\cr
$\stackrel{}{{}_{p_3\to \infty}}$\cr
}} }\sum_{j_3=0}^{p_3}\sum_{j_1=0}^{\infty}  
C_{j_3 j_1 j_1}  
J'[\phi_{j_3}]^{(i_3)}_{T,t}-
\sum_{j_3=0}^{p_3}\sum_{j_1=0}^{\min\{p_1,p_2\}}
C_{j_3j_1j_1}J'[\phi_{j_3}]^{(i_3)}_{T,t}\right)^2\right\}+
$$

\vspace{3mm}
$$
+4\cdot {\bf 1}_{\{i_2=i_3\ne 0\}}\times
$$

$$
\times
{\sf M}\left\{\left(
\hbox{\vtop{\offinterlineskip\halign{
\hfil#\hfil\cr
{\rm l.i.m.}\cr
$\stackrel{}{{}_{p_1\to \infty}}$\cr
}} }\sum_{j_1=0}^{p_1}\sum_{j_3=0}^{\infty}  
C_{j_3 j_3 j_1}  
J'[\phi_{j_1}]^{(i_1)}_{T,t}-
\sum_{j_1=0}^{p_1}\sum_{j_3=0}^{\min\{p_2,p_3\}}
C_{j_3j_3j_1}J'[\phi_{j_1}]^{(i_1)}_{T,t}\right)^2\right\}+
$$

\vspace{3mm}
$$
+4\cdot {\bf 1}_{\{i_1=i_3\ne 0\}}
{\sf M}\left\{\left(
\sum_{j_2=0}^{p_2}\sum_{j_1=0}^{\min\{p_1,p_3\}}
C_{j_1j_2j_1}J'[\phi_{j_2}]^{(i_2)}_{T,t}\right)^2\right\}=
$$

\vspace{3mm}
\begin{equation}
\label{20234}
= 4 A_{p_1 p_2 p_3} + 4\cdot {\bf 1}_{\{i_1=i_2\ne 0\}}B_{p_1 p_2 p_3}+
4\cdot {\bf 1}_{\{i_2=i_3\ne 0\}}C_{p_1 p_2 p_3}+
4\cdot {\bf 1}_{\{i_1=i_3\ne 0\}}D_{p_1 p_2 p_3}.
\end{equation}

\vspace{6mm}

Theorem~3 gives (see (\ref{tyyyarg}))

\vspace{-1mm}
\begin{equation}
\label{20235}
\lim\limits_{p_1,p_2,p_3\to\infty}
A_{p_1 p_2 p_3}=0.
\end{equation}

\vspace{5mm}

Further, in complete analogy with (\ref{after5001}) and using (\ref{after80xx}), we obtain

\vspace{1mm}
$$
D_{p_1 p_2 p_3}=
\sum_{j_2=0}^{p_2}\left(\sum_{j_1=0}^{\min\{p_1,p_3\}}
C_{j_1 j_2 j_1}\right)^2=
\sum_{j_2=0}^{p_2}\left(\sum_{j_1=\min\{p_1,p_3\}+1}^{\infty}
C_{j_1 j_2 j_1}\right)^2\le
$$

\vspace{3mm}
\begin{equation}
\label{20236}
\le
\sum_{j_2=0}^{\infty}\left(\sum_{j_1=\min\{p_1,p_3\}+1}^{\infty}
C_{j_1 j_2 j_1}\right)^2\le
\frac{K}{\left(\min\{p_1,p_3\}\right)^{2-\varepsilon}}\ \to\  0
\end{equation}

\vspace{5mm}
\noindent
if $p_1,p_2,p_3\to\infty,$ where 
$\varepsilon$ is an arbitrary
small positive real number,
constant $K$ is independent of $p$.

We have

$$
B_{p_1 p_2 p_3}=
{\sf M}\left\{\Biggl(\Biggl(
\hbox{\vtop{\offinterlineskip\halign{
\hfil#\hfil\cr
{\rm l.i.m.}\cr
$\stackrel{}{{}_{p_3\to \infty}}$\cr
}} }\sum_{j_3=0}^{p_3}\sum_{j_1=0}^{\infty}  
C_{j_3 j_1 j_1}  
J'[\phi_{j_3}]^{(i_3)}_{T,t}-
\sum_{j_3=0}^{p_3}\sum_{j_1=0}^{\infty}  
C_{j_3 j_1 j_1}  
J'[\phi_{j_3}]^{(i_3)}_{T,t}\Biggr)+\Biggr.\right.
$$

\vspace{4mm}
$$
\left.\Biggl.+\Biggl(\sum_{j_3=0}^{p_3}\sum_{j_1=0}^{\infty}  
C_{j_3 j_1 j_1}  
J'[\phi_{j_3}]^{(i_3)}_{T,t}-
\sum_{j_3=0}^{p_3}\sum_{j_1=0}^{\min\{p_1,p_2\}}
C_{j_3j_1j_1}J'[\phi_{j_3}]^{(i_3)}_{T,t}\Biggr)\Biggr)^2\right\}\le
$$

\vspace{4mm}
\begin{equation}
\label{20237}
\le
2 E_{p_3}+ 2 F_{p_1 p_2 p_3},
\end{equation}

\vspace{5mm}
\noindent
where

$$
E_{p_3}=
{\sf M}\left\{\Biggl(
\hbox{\vtop{\offinterlineskip\halign{
\hfil#\hfil\cr
{\rm l.i.m.}\cr
$\stackrel{}{{}_{p_3\to \infty}}$\cr
}} }\sum_{j_3=0}^{p_3}\sum_{j_1=0}^{\infty}  
C_{j_3 j_1 j_1}  
J'[\phi_{j_3}]^{(i_3)}_{T,t}-
\sum_{j_3=0}^{p_3}\sum_{j_1=0}^{\infty}  
C_{j_3 j_1 j_1}  
J'[\phi_{j_3}]^{(i_3)}_{T,t}\Biggr)^2\right\},
$$

\vspace{4mm}
$$
F_{p_1 p_2 p_3}=
{\sf M}\left\{
\Biggl(\sum_{j_3=0}^{p_3}\sum_{j_1=0}^{\infty}  
C_{j_3 j_1 j_1}  
J'[\phi_{j_3}]^{(i_3)}_{T,t}-
\sum_{j_3=0}^{p_3}\sum_{j_1=0}^{\min\{p_1,p_2\}}
C_{j_3j_1j_1}J'[\phi_{j_3}]^{(i_3)}_{T,t}\Biggr)^2\right\}=
$$

\vspace{4mm}
$$
={\sf M}\left\{
\Biggl(
\sum_{j_3=0}^{p_3}\sum_{j_1=\min\{p_1,p_2\}+1}^{\infty}
C_{j_3j_1j_1}J'[\phi_{j_3}]^{(i_3)}_{T,t}\Biggr)^2\right\}=
$$

\vspace{4mm}
\begin{equation}
\label{20238}
=\sum_{j_3=0}^{p_3}\Biggl(\sum_{j_1=\min\{p_1,p_2\}+1}^{\infty}
C_{j_3j_1j_1}\Biggr)^2.
\end{equation}

\vspace{5mm}

By analogy with (\ref{after1804}) we get

\vspace{1mm}
$$
\sum_{j_3=0}^{p_3}\Biggl(\sum_{j_1=\min\{p_1,p_2\}+1}^{\infty}
C_{j_3j_1j_1}\Biggr)^2\le
\sum_{j_3=0}^{\infty}\Biggl(\sum_{j_1=\min\{p_1,p_2\}+1}^{\infty}
C_{j_3j_1j_1}\Biggr)^2\le
$$

\vspace{3mm}
\begin{equation}
\label{20239}
\le \frac{K}{\left(\min\{p_1,p_2\}\right)^2}\ \to\  0
\end{equation}

\vspace{5mm}
\noindent
if $p_1,p_2,p_3\to\infty,$ where constant $K$ does not depend on $p.$

Moreover,

\vspace{-2mm}
\begin{equation}
\label{202310}
\lim\limits_{p_3\to\infty}E_{p_3}=
\lim\limits_{p_1,p_2,p_3\to\infty}E_{p_3}=0.
\end{equation}

\vspace{5mm}

Combining (\ref{20237})--(\ref{202310}), we obtain

\vspace{-1mm}
\begin{equation}
\label{202311}
\lim\limits_{p_1,p_2,p_3\to\infty}B_{p_1 p_2 p_3}=0.
\end{equation}

\vspace{5mm}

Consider $C_{p_1 p_2 p_3}$. We have

\vspace{2mm}
$$
C_{p_1 p_2 p_3}=
{\sf M}\left\{\Biggl(\Biggl(
\hbox{\vtop{\offinterlineskip\halign{
\hfil#\hfil\cr
{\rm l.i.m.}\cr
$\stackrel{}{{}_{p_1\to \infty}}$\cr
}} }\sum_{j_1=0}^{p_1}\sum_{j_3=0}^{\infty}  
C_{j_3 j_3 j_1}  
J'[\phi_{j_1}]^{(i_1)}_{T,t}-
\sum_{j_1=0}^{p_1}\sum_{j_3=0}^{\infty}  
C_{j_3 j_3 j_1}  
J'[\phi_{j_1}]^{(i_1)}_{T,t}\Biggr)+\Biggr.\right.
$$

\vspace{4mm}
$$
\left.\Biggl.+\Biggl(\sum_{j_1=0}^{p_1}\sum_{j_3=0}^{\infty}  
C_{j_3 j_3 j_1}  
J'[\phi_{j_1}]^{(i_1)}_{T,t}-
\sum_{j_1=0}^{p_1}\sum_{j_3=0}^{\min\{p_2,p_3\}}
C_{j_3j_3j_1}J'[\phi_{j_1}]^{(i_1)}_{T,t}\Biggr)\Biggr)^2\right\}\le
$$

\vspace{3mm}
\begin{equation}
\label{202350}
\le
2 G_{p_1}+ 2 H_{p_1 p_2 p_3},
\end{equation}

\vspace{4mm}
\noindent
where

$$
G_{p_1}=
{\sf M}\left\{\Biggl(
\hbox{\vtop{\offinterlineskip\halign{
\hfil#\hfil\cr
{\rm l.i.m.}\cr
$\stackrel{}{{}_{p_1\to \infty}}$\cr
}} }\sum_{j_1=0}^{p_1}\sum_{j_3=0}^{\infty}  
C_{j_3 j_3 j_1}  
J'[\phi_{j_1}]^{(i_1)}_{T,t}-
\sum_{j_1=0}^{p_1}\sum_{j_3=0}^{\infty}  
C_{j_3 j_3 j_1}  
J'[\phi_{j_1}]^{(i_1)}_{T,t}\Biggr)^2\right\},
$$

\vspace{4mm}
$$
H_{p_1 p_2 p_3}=
{\sf M}\left\{
\Biggl(\sum_{j_1=0}^{p_1}\sum_{j_3=0}^{\infty}  
C_{j_3 j_3 j_1}  
J'[\phi_{j_1}]^{(i_1)}_{T,t}-
\sum_{j_1=0}^{p_1}\sum_{j_3=0}^{\min\{p_2,p_3\}}
C_{j_3j_3j_1}J'[\phi_{j_1}]^{(i_1)}_{T,t}\Biggr)^2\right\}=
$$

\vspace{4mm}
$$
={\sf M}\left\{
\Biggl(
\sum_{j_1=0}^{p_1}\sum_{j_3=\min\{p_2,p_3\}+1}^{\infty}
C_{j_3j_3j_1}J'[\phi_{j_1}]^{(i_1)}_{T,t}\Biggr)^2\right\}=
$$

\vspace{4mm}
\begin{equation}
\label{202351}
=\sum_{j_1=0}^{p_1}\Biggl(\sum_{j_3=\min\{p_2,p_3\}+1}^{\infty}
C_{j_3j_3j_1}\Biggr)^2.
\end{equation}

\vspace{5mm}

By analogy with (\ref{after1903}) we get

\vspace{1mm}
$$
\sum_{j_1=0}^{p_1}\Biggl(\sum_{j_3=\min\{p_2,p_3\}+1}^{\infty}
C_{j_3j_3j_1}\Biggr)^2\le
\sum_{j_1=0}^{\infty}\Biggl(\sum_{j_3=\min\{p_2,p_3\}+1}^{\infty}
C_{j_3j_3j_1}\Biggr)^2\le
$$

\vspace{3mm}
\begin{equation}
\label{202352}
\le \frac{K}{\left(\min\{p_2,p_3\}\right)^2}\ \to\  0
\end{equation}

\vspace{5mm}
\noindent
if $p_1,p_2,p_3\to\infty,$ where constant $K$ does not depend on $p.$

Moreover,

\vspace{-2mm}
\begin{equation}
\label{202353}
\lim\limits_{p_1\to\infty}G_{p_1}=
\lim\limits_{p_1,p_2,p_3\to\infty}G_{p_1}=0.
\end{equation}

\vspace{5mm}

Combining (\ref{202350})--(\ref{202353}), we obtain

\vspace{-1mm}
\begin{equation}
\label{202354}
\lim\limits_{p_1,p_2,p_3\to\infty}C_{p_1 p_2 p_3}=0.
\end{equation}

\vspace{5mm}

The relations (\ref{20234})--(\ref{20236}),  (\ref{202311}),  
(\ref{202354}) complete the proof of Theorem~24.
Theorem~24 is proved.

\vspace{5mm}

\section{Hypotheses 1--3 from the Point
of View of the Wong--Zakai Approximation}

\vspace{5mm}

The iterated Ito stochastic integrals and solutions
of Ito SDEs are complex and important func\-ti\-o\-nals
from the independent components ${\bf f}_{\tau}^{(i)},$
$i=1, \ldots,m$ of the multidimensional
Wiener process ${\bf f}_{\tau}.$
Let ${\bf f}_{\tau}^{(i)p}$ $(p\in\mathbb{N})$ 
be some approximation of
${\bf f}_{\tau}^{(i)},$
$i=1,\ldots,m$.
Suppose that 
${\bf f}_{\tau}^{(i)p}$
converges to
${\bf f}_{\tau}^{(i)},$
$i=1,\ldots,m$ in some sense and has
differentiable sample trajectories.

A natural question arises: if we replace 
${\bf f}_{\tau}^{(i)}$
by ${\bf f}_{\tau}^{(i)p},$
$i=1,\ldots,m$ in the functionals
mentioned above, will the resulting
functionals converge to the original
functionals from the components 
${\bf f}_{\tau}^{(i)},$
$i=1,\ldots,m$ of the multidimentional
Wiener process ${\bf f}_{\tau}$?

The answere to this question is negative 
in the general case. However, 
in the pioneering works of Wong E. and Zakai M. \cite{W-Z-1},
\cite{W-Z-2} (also see \cite{W-Z-3}-\cite{Gyon2}),
it was shown that under the special conditions and 
for some types of approximations 
of the Wiener process the answere is affirmative
with one peculiarity: the convergence takes place 
to the iterated Stratonovich stochastic integrals
and solutions of Stratonovich SDEs and not to the iterated 
Ito stochastic integrals and solutions
of Ito SDEs.

The piecewise 
linear approximation 
as well as the regularization by convolution 
\cite{W-Z-1}-\cite{Gyon2}) relate the 
mentioned types of approximations
of the Wiener process. The above approximation 
of stochastic integrals and solutions of SDEs 
is often called the Wong--Zakai approximation.

Let ${\bf w}_{\tau},$ $\tau\in[0, T]$ is a random vector with 
an $m+1$ components: ${\bf w}_{\tau}^{(i)}={\bf f}_{\tau}^{(i)}$ 
for $i=1,\ldots,m$ and 
${\bf w}_{\tau}^{(0)}=\tau,$\ 
${\bf f}_{\tau}^{(i)}$ $(i=1,\ldots,m)$
are independent standard Wiener processes.

It is well known that the following representation 
takes place \cite{Lipt}, \cite{7e} 

\vspace{-1mm}
\begin{equation}
\label{um1x}
{\bf w}_{\tau}^{(i)}-{\bf w}_{t}^{(i)}=
\sum_{j=0}^{\infty}\int\limits_t^{\tau}
\phi_j(s)ds\ \zeta_j^{(i)},\ \ \ \zeta_j^{(i)}=
\int\limits_t^T \phi_j(s)d{\bf w}_s^{(i)},
\end{equation}

\vspace{3mm}
\noindent
where $\tau\in[t, T],$ $t\ge 0,$
$\{\phi_j(x)\}_{j=0}^{\infty}$ is an arbitrary complete 
orthonormal system of functions in the space $L_2([t, T]),$ and
$\zeta_j^{(i)}$ are independent standard Gaussian 
random variables for various $i$ or $j.$
Moreover, the series (\ref{um1x}) converges for any $\tau\in [t, T]$
in the mean-square sense.

Let ${\bf w}_{\tau}^{(i)p}-{\bf w}_{t}^{(i)p}$ be 
the mean-square approximation of the process
${\bf w}_{\tau}^{(i)}-{\bf w}_{t}^{(i)},$
which has the following form
\begin{equation}
\label{um1xx}
{\bf w}_{\tau}^{(i)p}-{\bf w}_{t}^{(i)p}=
\sum_{j=0}^{p}\int\limits_t^{\tau}
\phi_j(s)ds\ \zeta_j^{(i)}.
\end{equation}

\vspace{3mm}

From (\ref{um1xx}) we obtain
\begin{equation}
\label{um1xxx}
d{\bf w}_{\tau}^{(i)p}=
\sum_{j=0}^{p}
\phi_j(\tau)\zeta_j^{(i)} d\tau.
\end{equation}

\vspace{3mm}

Consider the following iterated Riemann--Stieltjes
integral

\vspace{-1mm}
\begin{equation}
\label{um1xxxx}
\int\limits_t^T
\psi_k(t_k)\ldots 
\int\limits_t^{t_2}\psi_1(t_1)
d{\bf w}_{t_1}^{(i_1)p_1}
\ldots d{\bf w}_{t_k}^{(i_k)p_k},
\end{equation}

\vspace{3mm}
\noindent
where $p_1,\ldots,p_k\in\mathbb{N},$ $i_1,\ldots,i_k=0,1,\ldots,m$ and

\vspace{-1mm}
\begin{equation}
\label{um1xxx1}
d{\bf w}_{\tau}^{(i)p}=
\left\{\begin{matrix}
d{\bf f}_{\tau}^{(i)p}\ &\hbox{\rm for}\ \ \ i=1,2,\ldots,m\cr\cr\cr
d\tau^p\ &\hbox{\rm for}\ \ \ i=0
\end{matrix}
\right.,\ \ \ \ \ p\in\mathbb{N},
\end{equation}

\vspace{3mm}
\noindent
where $d{\bf f}_{\tau}^{(i)p},$ $d\tau^p$ are defined by (\ref{um1xxx}).

Let us substitute (\ref{um1xxx}) into (\ref{um1xxxx})

\vspace{-1mm}
\begin{equation}
\label{um1xxxx1}
\int\limits_t^T
\psi_k(t_k)\ldots 
\int\limits_t^{t_2}\psi_1(t_1)
d{\bf w}_{t_1}^{(i_1)p_1}
\ldots d{\bf w}_{t_k}^{(i_k)p_k}
=
\sum\limits_{j_1=0}^{p_1}\ldots \sum\limits_{j_k=0}^{p_k}
C_{j_k \ldots j_1}\prod\limits_{l=1}^k \zeta_{j_l}^{(i_l)},
\end{equation}

\vspace{3mm}
\noindent
where 
$\zeta_j^{(i)}$ are independent standard Gaussian 
random variables for various $i$ or $j$ (if $i\ne 0$),

\vspace{-1mm}
$$
C_{j_k \ldots j_1}=\int\limits_t^T\psi_k(t_k)\phi_{j_k}(t_k)\ldots
\int\limits_t^{t_2}
\psi_1(t_1)\phi_{j_1}(t_1)
dt_1 \ldots dt_k
$$

\vspace{3mm}
\noindent
is the Fourier coefficient.

To best of our knowledge \cite{W-Z-1}, \cite{W-Z-2}, \cite{Watanabe}
the approximations of the Wiener process
in the Wong--Zakai approximation must satisfy fairly strong
restrictions
\cite{Watanabe}
(see Definition 7.1 on Pages 480--481).
At least the proof of an analogue of Theorem 7.2 
(see \cite{Watanabe}, Page 497)
for approximations of the Wiener 
process based on its series expansion (\ref{um1x})
should be carried out separately.
Thus, the mean-square convergence of the right-hand side
of (\ref{um1xxxx1}) to the iterated Stratonovich stochastic integral 
(\ref{str})
does not follow from the results of the papers
\cite{W-Z-1}, \cite{W-Z-2} (also see \cite{Watanabe},
Theorems 7.1, 7.2).

From the other hand, Theorems 3, 4, 6--12, 16--18, 23, 24 from this 
paper (see proofs of Theorems 3, 4, 6--12
in Chapters 1 and 2 of the monographs
\cite{2018a}-\cite{12aa-afterxxx}) 
can be considered as the proof of the
Wong--Zakai approximation for the iterated 
Stratonovich stochastic integrals (\ref{str}) of multiplicities 1 to 6
(or of multiplicity $k$ under the condition of convergence 
of trace series (Theorem~13))
based on the approximation (\ref{um1xx}) of the Wiener process
in the form of its series expansion.
At that, the Riemann--Stieltjes integrals (\ref{um1xxxx}) converge
(according to Theorems 6--13, 16--18, 23, 24)
to the mentioned Stratonovich 
stochastic integrals (\ref{str}). Recall that
$\{\phi_j(x)\}_{j=0}^{\infty}$ (see (\ref{um1x}), (\ref{um1xx}), and
Theorems 6--12, 16--18, 23, 24)
is a complete 
orthonormal system of Legendre polynomials or 
trigonometric functions 
in the space $L_2([t, T])$.

To illustrate the above reasoning, 
consider two examples for the case $k=2,$
$\psi_1(\tau),$ $\psi_2(\tau)\equiv 1,$ $i_1, i_2=1,\ldots,m.$

The first example relates to the piecewise linear approximation
of the multidimensional Wiener process (these approximations 
were considered in \cite{W-Z-1}, \cite{W-Z-2}, \cite{Watanabe}).

Let ${\bf b}_{\Delta}^{(i)}(t),$ $t\in[0, T]$ be the piecewise
linear approximation of the $i$th component ${\bf f}_t^{(i)}$
of the multidimensional standard Wiener process ${\bf f}_t,$
$t\in [0, T]$ with independent components
${\bf f}_t^{(i)},$ $i=1,\ldots,m,$ i.e.

\vspace{-1mm}
$$
{\bf b}_{\Delta}^{(i)}(t)={\bf f}_{k\Delta}^{(i)}+
\frac{t-k\Delta}{\Delta}\Delta{\bf f}_{k\Delta}^{(i)},
$$

\vspace{2mm}
\noindent
where 

\vspace{-2mm}
$$
\Delta{\bf f}_{k\Delta}^{(i)}={\bf f}_{(k+1)\Delta}^{(i)}-
{\bf f}_{k\Delta}^{(i)},\ \ \
t\in[k\Delta, (k+1)\Delta),\ \ \ k=0, 1, 2,\ldots, N-1.
$$

\vspace{4mm}

Consider the following iterated Riemann--Stieltjes
integral

\vspace{-1mm}
$$
\int\limits_0^T
\int\limits_0^{s}
d{\bf b}_{\Delta}^{(i_1)}(\tau)d{\bf b}_{\Delta}^{(i_2)}(s),\ \ \ 
i_1,i_2=1,\ldots,m.
$$

\vspace{3mm}

We can write w.~p.~1

\vspace{-1mm}
$$
\int\limits_0^T
\int\limits_0^{s}
d{\bf b}_{\Delta}^{(i_1)}(\tau)d{\bf b}_{\Delta}^{(i_2)}(s)=
\int\limits_0^T
\int\limits_0^{s}
\frac{\partial {\bf b}_{\Delta}^{(i_1)}}{\partial \tau}(\tau)d\tau
\frac{\partial {\bf b}_{\Delta}^{(i_2)}}{\partial s}(s)
ds =
$$

\vspace{2mm}
$$
=
\sum\limits_{l=0}^{N-1}\int\limits_{l\Delta}^{(l+1)\Delta}
\left(
\sum\limits_{q=0}^{l-1}\int\limits_{q\Delta}^{(q+1)\Delta}
\frac{\Delta{\bf f}_{q\Delta}^{(i_1)}}{\Delta}d\tau+
\int\limits_{l\Delta}^{s}
\frac{\Delta{\bf f}_{l\Delta}^{(i_1)}}{\Delta}d\tau\right)
\frac{\Delta{\bf f}_{l\Delta}^{(i_2)}}{\Delta}ds=
$$

\vspace{2mm}
$$
=\sum\limits_{l=0}^{N-1}\sum\limits_{q=0}^{l-1}
\Delta{\bf f}_{q\Delta}^{(i_1)}
\Delta{\bf f}_{l\Delta}^{(i_2)}+
\frac{1}{\Delta^2}\sum\limits_{l=0}^{N-1}
\Delta{\bf f}_{l\Delta}^{(i_1)}
\Delta{\bf f}_{l\Delta}^{(i_2)}
\int\limits_{l\Delta}^{(l+1)\Delta}
\int\limits_{l\Delta}^{s}d\tau ds=
$$

\vspace{2mm}
\begin{equation}
\label{oh-ty}
=\sum\limits_{l=0}^{N-1}\sum\limits_{q=0}^{l-1}
\Delta{\bf f}_{q\Delta}^{(i_1)}
\Delta{\bf f}_{l\Delta}^{(i_2)}+
\frac{1}{2}\sum\limits_{l=0}^{N-1}
\Delta{\bf f}_{l\Delta}^{(i_1)}
\Delta{\bf f}_{l\Delta}^{(i_2)}.
\end{equation}

\vspace{5mm}

Using (\ref{oh-ty}), it 
is not difficult to show that

\vspace{-1mm}
$$
\hbox{\vtop{\offinterlineskip\halign{
\hfil#\hfil\cr
{\rm l.i.m.}\cr
$\stackrel{}{{}_{N\to \infty}}$\cr
}} }
\int\limits_0^T
\int\limits_0^{s}
d{\bf b}_{\Delta}^{(i_1)}(\tau)d{\bf b}_{\Delta}^{(i_2)}(s)=
\int\limits_0^T
\int\limits_0^{s}
d{\bf f}_{\tau}^{(i_1)}d{\bf f}_{s}^{(i_2)}+
\frac{1}{2}{\bf 1}_{\{i_1=i_2\}}\int\limits_0^T ds=
$$

\vspace{2mm}
\begin{equation}
\label{uh-111}
=
{\int\limits_0^{*}}^T
{\int\limits_0^{*}}^s
d{\bf f}_{\tau}^{(i_1)}d{\bf f}_{s}^{(i_2)},
\end{equation}

\vspace{5mm}
\noindent
where $\Delta\to 0$ if $N\to\infty$ ($N\Delta=T$).

Obviously, (\ref{uh-111}) agrees with Theorem 7.1 (see \cite{Watanabe},
Page 486).

The next example relates to the approximation (\ref{um1xx})
of the Wiener process based on its series expansion (\ref{um1x}), where
$t=0$ and 
$\{\phi_j(x)\}_{j=0}^{\infty}$ 
is a complete 
orthonormal system of Legendre po\-ly\-no\-mi\-als or 
trigonometric functions 
in the space $L_2([0, T]).$

Consider the following iterated Riemann--Stieltjes
integral
       
\vspace{-1mm}
\begin{equation}
\label{abcd1}
\int\limits_0^T
\int\limits_0^{s}
d{\bf f}_{\tau}^{(i_1)p}d{\bf f}_{s}^{(i_2)p},\ \ \ 
i_1,i_2=1,2,\ldots,m,
\end{equation}

\vspace{3mm}
\noindent
where $d{\bf f}_{\tau}^{(i)p}$ is defined by the
relation
(\ref{um1xxx}).

Let us substitute (\ref{um1xxx}) into (\ref{abcd1}) 

\vspace{-1mm}
\begin{equation}
\label{set18}
\int\limits_0^T
\int\limits_0^{s}
d{\bf f}_{\tau}^{(i_1)p}d{\bf f}_{s}^{(i_2)p}=
\sum\limits_{j_1,j_2=0}^p
C_{j_2 j_1} \zeta_{j_1}^{(i_1)}\zeta_{j_2}^{(i_2)},
\end{equation}

\vspace{3mm}
\noindent
where 
$$
C_{j_2 j_1}=
\int\limits_0^T \phi_{j_2}(s)\int\limits_0^s
\phi_{j_1}(\tau)d\tau ds
$$

\vspace{3mm}
\noindent
is the Fourier coefficient; another notations 
are the same as in (\ref{um1xxxx1}).

As we noted above, approximations of the Wiener process that are
similar to (\ref{um1xx})
were not considered in \cite{W-Z-1}, \cite{W-Z-2}
(also see Theorems 7.1, 7.2 in \cite{Watanabe}).
Furthermore, transferring of the results of Theorems 7.1 and 7.2
\cite{Watanabe} to the case under consideration is
not obvious.

On the other hand, we can apply the theory built in Chapters 1 and 2
of the monographs \cite{2018a}-\cite{12aa-afterxxx}. More precisely, 
using 
Theorems 6 and 7, we obtain from (\ref{set18}) the desired result

\vspace{-1mm}
$$
\hbox{\vtop{\offinterlineskip\halign{
\hfil#\hfil\cr
{\rm l.i.m.}\cr
$\stackrel{}{{}_{p\to \infty}}$\cr
}} }
\int\limits_0^T
\int\limits_0^{s}
d{\bf f}_{\tau}^{(i_1)p}d{\bf f}_{s}^{(i_2)p}=
\hbox{\vtop{\offinterlineskip\halign{
\hfil#\hfil\cr
{\rm l.i.m.}\cr
$\stackrel{}{{}_{p\to \infty}}$\cr
}} }
\sum\limits_{j_1,j_2=0}^p
C_{j_2 j_1} \zeta_{j_1}^{(i_1)}\zeta_{j_2}^{(i_2)}=
$$

\vspace{1mm}
\begin{equation}
\label{umen-bl}
=
{\int\limits_0^{*}}^T
{\int\limits_0^{*}}^s
d{\bf f}_{\tau}^{(i_1)}d{\bf f}_{s}^{(i_2)}.
\end{equation}

\vspace{4mm}

From the other hand, by Theorems 3, 4
(see (\ref{a2})) for the case
$k=2$ we obtain from (\ref{set18}) the following relation

\vspace{-1mm}
$$
\hbox{\vtop{\offinterlineskip\halign{
\hfil#\hfil\cr
{\rm l.i.m.}\cr
$\stackrel{}{{}_{p\to \infty}}$\cr
}} }
\int\limits_0^T
\int\limits_0^{s}
d{\bf f}_{\tau}^{(i_1)p}d{\bf f}_{s}^{(i_2)p}=
\hbox{\vtop{\offinterlineskip\halign{
\hfil#\hfil\cr
{\rm l.i.m.}\cr
$\stackrel{}{{}_{p\to \infty}}$\cr
}} }
\sum\limits_{j_1,j_2=0}^p
C_{j_2 j_1} \zeta_{j_1}^{(i_1)}\zeta_{j_2}^{(i_2)}=
$$

\vspace{1mm}
$$
=
\hbox{\vtop{\offinterlineskip\halign{
\hfil#\hfil\cr
{\rm l.i.m.}\cr
$\stackrel{}{{}_{p\to \infty}}$\cr
}} }
\sum\limits_{j_1,j_2=0}^p
C_{j_2 j_1} \biggl(\zeta_{j_1}^{(i_1)}\zeta_{j_2}^{(i_2)}-
{\bf 1}_{\{i_1=i_2\}}{\bf 1}_{\{j_1=j_2\}}\biggr)+
{\bf 1}_{\{i_1=i_2\}}\sum\limits_{j_1=0}^{\infty}
C_{j_1 j_1}=
$$

\vspace{1mm}
\begin{equation}
\label{umen-blx}
=
\int\limits_0^T
\int\limits_0^{s}
d{\bf f}_{\tau}^{(i_1)}d{\bf f}_{s}^{(i_2)}+
{\bf 1}_{\{i_1=i_2\}}\sum\limits_{j_1=0}^{\infty}
C_{j_1 j_1}.
\end{equation}

\vspace{4mm}

Since
$$
\sum\limits_{j_1=0}^{\infty}
C_{j_1 j_1}=\frac{1}{2}\sum\limits_{j_1=0}^{\infty}
\left(\int\limits_0^T \phi_j(\tau)d\tau\right)^2
=
$$

\vspace{2mm}
$$
=\frac{1}{2}
\left(\int\limits_0^T \phi_0(\tau)d\tau\right)^2=\frac{1}{2}
\int\limits_0^T ds,
$$

\vspace{5mm}
\noindent
then from Theorem 5 ($k=2$) and (\ref{umen-blx}) we obtain (\ref{umen-bl}).

\vspace{5mm}

\section{Wong--Zakai Type Theorems for Iterated Stratonovich Stochastic 
Integrals. The Case of Approximation of the Multidimensional Wiener Process 
Based on its Series Expansion Using Legendre Polynomials and 
Trigonometric Functions}

\vspace{5mm}

As we mentioned above, there exists a lot of publications on the subject 
of Wong--Zakai appro\-xi\-ma\-ti\-on of stochastic integrals
and SDEs \cite{W-Z-1}-\cite{Gyon2}.
However, these works did not consider the approximation of 
iterated stochastic integrals and systems of SDEs
for the case of approximation of the 
multidimensional Wiener process based on its series expansions.
Usually, as an approximation of the Wiener process in the theorems 
of the Wong--Zakai type, the authors \cite{W-Z-1}-\cite{Gyon2} 
choose a piecewise linear approximation or an 
approximation based on the regularization by convolution.

The Wong--Zakai approximation is widely used to approximate 
stochastic integrals and SDEs. 
In particular, the Wong--Zakai approximation can be used to 
approximate the iterated Stratonovich stochastic integrals 
in the context of numerical integration of Ito SDEs
in the framework of the approach 
based on the Taylor--Stratonovich expansion \cite{KlPl2}-\cite{new-art-1-xxy},
\cite{Zapad-2}, \cite{Zapad-9}.
It should be noted that the authors of the works
\cite{KlPl2}
(Sect.~5.8, pp.~202--204), \cite{KPS} (pp.~82-84),
\cite{Zapad-2} (pp.~438-439),  
\cite{Zapad-9} (pp.~263-264) mention 
the Wong--Zakai approximation 
\cite{W-Z-1}-\cite{W-Z-3}, \cite{Watanabe} within the frames
of approximation of iterated Stratonovich stochastic integrals
based on the Karhunen--Loeve expansion of the Brownian bridge
process. However, in these works there is no rigorous proof 
of convergence for approximations of the mentioned stochastic integrals
of miltiplicity 3 and higher.

From the other hand, the theory constructed in Chapters 1 and 2
of the monographs
\cite{2018a}-\cite{12aa-afterxxx} 
can be considered as the proof of the Wong--Zakai
approximation for iterated Stratonovich stochastic integrals
of multiplicities 1 to 6 based on the Wiener process series expansion 
using Legendre polynomials and trigonometric functions.

The subject of this section is to reformulate the results of 
Chapter 2 of the monograph \cite{2018a} (also see \cite{2018aa}--\cite{12aa-afterxxx})
in the form of theorems on convergence of iterated
Riemann--Stiltjes integrals 
to iterated Stratonovich stochastic integrals.

Let us reformulate Theorems 1, 2, 7--13, 16--18, 23, 24 of this paper as
theorems on the convergence 
of iterated Riemann--Stiltjes integrals (\ref{um1xxxx})
to the iterated Stratonovich stochastic integrals (\ref{str}).

\vspace{2mm}

{\bf Theorem 25}.
{\it Suppose that the following conditions are fulfilled{\rm :}

{\rm 1}. Every $\psi_l(\tau)$ $(l=1, 2)$ is a continuously 
differentiable func\-ti\-on at the interval $[t, T]$.

{\rm 2}. $\{\phi_j(x)\}_{j=0}^{\infty}$ is a complete orthonormal system 
of Legendre polynomials or trigonometric functions
in the space $L_2([t, T]).$

Then, for the iterated 
Stratonovich stochastic integral of second multiplicity

\vspace{-1mm}
$$
J^{*}[\psi^{(2)}]_{T,t}={\int\limits_t^{*}}^T\psi_2(t_2)
{\int\limits_t^{*}}^{t_2}\psi_1(t_1)d{\bf w}_{t_1}^{(i_1)}
d{\bf w}_{t_2}^{(i_2)}\ \ \ (i_1, i_2=0, 1,\ldots,m)
$$

\vspace{2mm}
\noindent
the following equality
\begin{equation}
\label{hhh1}
J^{*}[\psi^{(2)}]_{T,t}=
\hbox{\vtop{\offinterlineskip\halign{
\hfil#\hfil\cr
{\rm l.i.m.}\cr
$\stackrel{}{{}_{p_1,p_2\to \infty}}$\cr
}} }\int\limits_t^T
\psi_2(t_2)\int\limits_t^{t_2}\psi_1(t_1)
d{\bf w}_{t_1}^{(i_1)p_1}
d{\bf w}_{t_2}^{(i_2)p_2}
\end{equation}

\vspace{2mm}
\noindent
is valid$,$ here and further 
{\rm l.i.m.} is a limit in the mean-square sense.}

\vspace{2mm}

{\bf Theorem 26}.
{\it Suppose that
$\{\phi_j(x)\}_{j=0}^{\infty}$ is a complete orthonormal
system of Legendre polynomials or trigonometric func\-ti\-ons
in the space $L_2([t, T])$.
Then, for the iterated Stratonovich stochastic integral of 
third multiplicity

\vspace{-1mm}
$$
{\int\limits_t^{*}}^T
{\int\limits_t^{*}}^{t_3}
{\int\limits_t^{*}}^{t_2}
d{\bf f}_{t_1}^{(i_1)}
d{\bf f}_{t_2}^{(i_2)}d{\bf f}_{t_3}^{(i_3)}\ \ \ (i_1, i_2, i_3=1,\ldots,m)
$$

\vspace{2mm}
\noindent
the following 
formula

\vspace{0.5mm}
$$
{\int\limits_t^{*}}^T
{\int\limits_t^{*}}^{t_3}
{\int\limits_t^{*}}^{t_2}
d{\bf f}_{t_1}^{(i_1)}
d{\bf f}_{t_2}^{(i_2)}d{\bf f}_{t_3}^{(i_3)}\ 
=
\hbox{\vtop{\offinterlineskip\halign{
\hfil#\hfil\cr
{\rm l.i.m.}\cr
$\stackrel{}{{}_{p_1,p_2,p_3\to \infty}}$\cr
}} }\int\limits_t^T
\int\limits_t^{t_3}
\int\limits_t^{t_2}
d{\bf f}_{t_1}^{(i_1)p_1}
d{\bf f}_{t_2}^{(i_2)p_2}
d{\bf f}_{t_3}^{(i_3)p_3}
$$

\vspace{2.5mm}
\noindent
is valid.}

\vspace{2mm}

{\bf Theorem 27}.\ {\it Suppose that
$\{\phi_j(x)\}_{j=0}^{\infty}$ is a complete orthonormal
system of Legendre poly\-no\-mi\-als
in the space $L_2([t, T])$.
Then, for the iterated Stratonovich stochastic integral of 
third multiplicity

\vspace{-1mm}
$$
I_{{l_1l_2l_3}_{T,t}}^{*(i_1i_2i_3)}={\int\limits_t^{*}}^T(t-t_3)^{l_3}
{\int\limits_t^{*}}^{t_3}(t-t_2)^{l_2}
{\int\limits_t^{*}}^{t_2}(t-t_1)^{l_1}
d{\bf f}_{t_1}^{(i_1)}
d{\bf f}_{t_2}^{(i_2)}d{\bf f}_{t_3}^{(i_3)}\ \ \ (i_1, i_2, i_3=1,\ldots,m)
$$

\vspace{3mm}
\noindent
the following 
equality

\vspace{-1mm}
$$
I_{{l_1l_2l_3}_{T,t}}^{*(i_1i_2i_3)}=
\hbox{\vtop{\offinterlineskip\halign{
\hfil#\hfil\cr
{\rm l.i.m.}\cr
$\stackrel{}{{}_{p_1,p_2,p_3\to \infty}}$\cr
}} }
\int\limits_t^T
(t-t_3)^{l_3}
\int\limits_t^{t_3}
(t-t_2)^{l_2}
\int\limits_t^{t_2}
(t-t_1)^{l_1}
d{\bf f}_{t_1}^{(i_1)p_1}
d{\bf f}_{t_2}^{(i_2)p_2}
d{\bf f}_{t_3}^{(i_3)p_3}
$$

\vspace{3mm}
\noindent
holds
for each of the following cases

\vspace{2mm}
\noindent
{\rm 1}.\ $i_1\ne i_2,\ i_2\ne i_3,\ i_1\ne i_3$\ and
$l_1, l_2, l_3=0, 1, 2,\ldots $\\
{\rm 2}.\ $i_1=i_2\ne i_3$ and $l_1=l_2\ne l_3$\ and
$l_1, l_2, l_3=0, 1, 2,\ldots $\\
{\rm 3}.\ $i_1\ne i_2=i_3$ and $l_1\ne l_2=l_3$\ and
$l_1, l_2, l_3=0, 1, 2,\ldots $\\
{\rm 4}.\ $i_1, i_2, i_3=1,\ldots,m;$ $l_1=l_2=l_3=l$\ and $l=0, 1, 
2,\ldots$\\
}

\vspace{2mm}

{\bf Theorem 28}.\
{\it Suppose that
$\{\phi_j(x)\}_{j=0}^{\infty}$ is a complete orthonormal
system of Legendre poly\-no\-mi\-als or trigonometric func\-ti\-ons
in the space $L_2([t, T])$
and $\psi_l(\tau)$ $(l=1, 2, 3)$ are continuously
differentiable functions at the interval $[t, T]$.
Then, for the iterated Stratonovich stochastic 
integral of third multiplicity

\vspace{-1mm}
$$
J^{*}[\psi^{(3)}]_{T,t}={\int\limits_t^{*}}^T\psi_3(t_3)
{\int\limits_t^{*}}^{t_3}\psi_2(t_2)
{\int\limits_t^{*}}^{t_2}\psi_1(t_1)
d{\bf f}_{t_1}^{(i_1)}
d{\bf f}_{t_2}^{(i_2)}d{\bf f}_{t_3}^{(i_3)}\ \ \ (i_1, i_2, i_3=1,\ldots,m)
$$

\vspace{3mm}
\noindent
the following 
formula

\vspace{-1mm}
$$
J^{*}[\psi^{(3)}]_{T,t}=
\hbox{\vtop{\offinterlineskip\halign{
\hfil#\hfil\cr
{\rm l.i.m.}\cr
$\stackrel{}{{}_{p\to \infty}}$\cr
}} }
\int\limits_t^T
\psi_3(t_3)
\int\limits_t^{t_3}
\psi_2(t_2)
\int\limits_t^{t_2}
\psi_1(t_1)
d{\bf f}_{t_1}^{(i_1)p}
d{\bf f}_{t_2}^{(i_2)p}
d{\bf f}_{t_3}^{(i_3)p}
$$

\vspace{3mm}
\noindent
holds
for each of the following cases{\rm :}

\vspace{2mm}
\noindent
{\rm 1}.\ $i_1\ne i_2,\ i_2\ne i_3,\ i_1\ne i_3,$\\
{\rm 2}.\ $i_1=i_2\ne i_3$ and
$\psi_1(\tau)\equiv\psi_2(\tau),$\\
{\rm 3}.\ $i_1\ne i_2=i_3$ and
$\psi_2(\tau)\equiv\psi_3(\tau),$\\
{\rm 4}.\ $i_1, i_2, i_3=1,\ldots,m$
and $\psi_1(\tau)\equiv\psi_2(\tau)\equiv\psi_3(\tau).$
}

\vspace{2mm}

{\bf Theorem 29}.\
{\it Suppose that
$\{\phi_j(x)\}_{j=0}^{\infty}$ is a complete orthonormal
system of Legendre polynomials or trigonomertic func\-ti\-ons
in the space $L_2([t, T])$. Furthermore, let
the function $\psi_2(\tau)$ is continuously
differentiable at the interval $[t, T]$ and
the functions $\psi_1(\tau),$ $\psi_3(\tau)$ are twice continuously
differentiable at the interval $[t, T]$.
Then, for the iterated Stra\-to\-no\-vich stochastic integral 
of third multiplicity

$$
J^{*}[\psi^{(3)}]_{T,t}={\int\limits_t^{*}}^T\psi_3(t_3)
{\int\limits_t^{*}}^{t_3}\psi_2(t_2)
{\int\limits_t^{*}}^{t_2}\psi_1(t_1)
d{\bf f}_{t_1}^{(i_1)}
d{\bf f}_{t_2}^{(i_2)}d{\bf f}_{t_3}^{(i_3)}
$$

\vspace{2mm}
\noindent
the following 
formula

\vspace{-1mm}
$$
J^{*}[\psi^{(3)}]_{T,t}=
\hbox{\vtop{\offinterlineskip\halign{
\hfil#\hfil\cr
{\rm l.i.m.}\cr
$\stackrel{}{{}_{p\to \infty}}$\cr
}} }
\int\limits_t^T
\psi_3(t_3)
\int\limits_t^{t_3}
\psi_2(t_2)
\int\limits_t^{t_2}
\psi_1(t_1)
d{\bf f}_{t_1}^{(i_1)p}
d{\bf f}_{t_2}^{(i_2)p}
d{\bf f}_{t_3}^{(i_3)p}
$$

\vspace{3mm}
\noindent
is valid$,$ where $i_1, i_2, i_3=1,\ldots,m$.}  

\vspace{2mm}

{\bf Theorem 30}.\ 
{\it Suppose that
$\{\phi_j(x)\}_{j=0}^{\infty}$ is a complete orthonormal
system of Legendre polynomials or trigonometric func\-ti\-ons
in the space $L_2([t, T])$.
Then, for the iterated Stra\-to\-no\-vich stochastic integral 
of fourth multiplicity

$$
I_{T,t}^{*(i_1 i_2 i_3 i_4)}=
{\int\limits_t^{*}}^T
{\int\limits_t^{*}}^{t_4}
{\int\limits_t^{*}}^{t_3}
{\int\limits_t^{*}}^{t_2}
d{\bf w}_{t_1}^{(i_1)}
d{\bf w}_{t_2}^{(i_2)}d{\bf w}_{t_3}^{(i_3)}d{\bf w}_{t_4}^{(i_4)}\ \ \
(i_1, i_2, i_3, i_4=0, 1,\ldots,m)
$$

\vspace{2mm}
\noindent
the following 
equality

\vspace{-1mm}
$$
I_{T,t}^{*(i_1 i_2 i_3 i_4)}=
\hbox{\vtop{\offinterlineskip\halign{
\hfil#\hfil\cr
{\rm l.i.m.}\cr
$\stackrel{}{{}_{p\to \infty}}$\cr
}} }
\int\limits_t^T
\int\limits_t^{t_4}
\int\limits_t^{t_3}
\int\limits_t^{t_2}
d{\bf w}_{t_1}^{(i_1)p}
d{\bf w}_{t_2}^{(i_2)p}
d{\bf w}_{t_3}^{(i_3)p}
d{\bf w}_{t_4}^{(i_4)p}
$$

\vspace{3mm}
\noindent
holds$,$ where 
${\bf w}_{\tau}^{(i)}={\bf f}_{\tau}^{(i)}$ $(i=1,\ldots,m)$ 
are independent 
standard Wiener processes and 
${\bf w}_{\tau}^{(0)}=\tau.$}       

\vspace{2mm}

{\bf Theorem 31}.\
{\it Suppose that $\psi_1(\tau), \ldots, \psi_k(\tau)$ are twice continuously
differentiable func\-ti\-ons at the interval
$[t, T]$ and
$\{\phi_j(x)\}_{j=0}^{\infty}$ is a complete
orthonormal system of
trigonometric functions in the space $L_2([t, T])$. 
Then, for the iterated Stratonovich stochastic integral

$$
J^{*}[\psi^{(k)}]_{T,t}=
{\int\limits_t^{*}}^T
\psi_k(t_k) \ldots 
{\int\limits_t^{*}}^{t_2}
\psi_1(t_1) d{\bf w}_{t_1}^{(i_1)}
\ldots
d{\bf w}_{t_k}^{(i_k)}
$$

\vspace{2mm}
\noindent
the following relation

\vspace{-1mm}
$$
\lim\limits_{p_1\to\infty}\varlimsup\limits_{p_2\to\infty}
\ldots\varlimsup\limits_{p_k\to\infty}
{\sf M}\left\{\Biggl(J^{*}[\psi^{(k)}]_{T,t}-
J^{*}[\psi^{(k)}]_{T,t}^{p_1,\ldots,p_k}
\Biggr)^{2n}\right\}=0
$$

\vspace{3mm}
\noindent
is valid, where $i_1,\ldots, i_k=0, 1, \ldots,m,$

$$
J^{*}[\psi^{(k)}]_{T,t}^{p_1,\ldots,p_k}=\int\limits_t^T
\psi_k(t_k)
\ldots
\int\limits_t^{t_2}
\psi_1(t_1)
d{\bf w}_{t_1}^{(i_1)p_1}
\ldots
d{\bf w}_{t_k}^{(i_k)p_k},
$$

\vspace{2mm}
\noindent
$n\in\mathbb{N},$ and $\varlimsup$ means $\limsup$.
}
                         
\vspace{2mm}

{\bf Theorem 32}.\
{\it Suppose that $\psi_1(\tau), \ldots, \psi_k(\tau)$ are twice continuously
differentiable func\-ti\-ons at the interval
$[t, T]$ and
$\{\phi_j(x)\}_{j=0}^{\infty}$ is a complete
orthonormal system of
trigonometric functions in the space $L_2([t, T])$. 
Then, for the iterated Stratonovich stochastic integral 

$$
J^{*}[\psi^{(k)}]_{T,t}=
{\int\limits_t^{*}}^T
\psi_k(t_k) \ldots 
{\int\limits_t^{*}}^{t_2}
\psi_1(t_1) d{\bf w}_{t_1}^{(i_1)}
\ldots
d{\bf w}_{t_k}^{(i_k)}
$$

\vspace{2mm}
\noindent
the following equality

\vspace{-1mm}
$$
\lim\limits_{p_k\to\infty}\varlimsup\limits_{p_{k-1}\to\infty}
\ldots
\varlimsup\limits_{p_1\to\infty}
{\sf M}\left\{\Biggl(J^{*}[\psi^{(k)}]_{T,t}-
J^{*}[\psi^{(k)}]_{T,t}^{p_1,\ldots,p_k}
\Biggr)^{2n}\right\}=0
$$

\vspace{3mm}
\noindent
is valid, where $i_1,\ldots, i_k=0, 1, \ldots,m,$

$$
J^{*}[\psi^{(k)}]_{T,t}^{p_1,\ldots,p_k}=\int\limits_t^T
\psi_k(t_k)
\ldots
\int\limits_t^{t_2}
\psi_1(t_1)
d{\bf w}_{t_1}^{(i_1)p_1}
\ldots
d{\bf w}_{t_k}^{(i_k)p_k},
$$

\vspace{2mm}
\noindent
$n\in\mathbb{N},$ and $\varlimsup$ means $\limsup$.
}

\vspace{2mm}

{\bf Theorem 33}.\
{\it Assume that
the continuously differentiable func\-ti\-ons 
$\psi_l(\tau)$ $(l=1,\ldots,k)$ and 
the complete orthonormal system $\{\phi_j(x)\}_{j=0}^{\infty}$
of continuous functions $(\phi_0(x)=1/\sqrt{T-t})$ 
in the space $L_2([t, T])$ are such that the following 
conditions are satisfied{\rm :}

\vspace{2mm}

{\rm 1.}\ The equality 

\vspace{-2mm}
\begin{equation}
\label{ssss111}
\frac{1}{2}
\int\limits_t^s \Phi_1(t_1)\Phi_2(t_1)dt_1
=\sum_{j_1=0}^{\infty}
\int\limits_t^s
\Phi_2(t_2)\phi_{j_1}(t_2)\int\limits_t^{t_2}
\Phi_1(t_1)\phi_{j_1}(t_1)dt_1 dt_2
\end{equation}

\vspace{3mm}
\noindent
holds for all $s\in (t, T],$ where the nonrandom functions 
$\Phi_1(\tau),$ $\Phi_2(\tau)$
are continuously differentiable on $[t, T]$
and the series on the right-hand side of {\rm (\ref{ssss111})}
converges absolutely.

{\rm 2.}\ The estimates

\vspace{-1mm}
$$
\left|\int\limits_t^s \phi_{j}(\tau)\Phi_1(\tau)d\tau\right|
\le \frac{\Psi_1(s)}{j^{1/2+\alpha}},\ \ \ 
\left|\int\limits_s^T \phi_{j}(\tau)\Phi_2(\tau)d\tau\right|\le
\frac{\Psi_1(s)}{j^{1/2+\alpha}},
$$

\vspace{2mm}
$$
\left|\sum_{j=p+1}^{\infty}\int\limits_t^s
\Phi_2(\tau)\phi_{j}(\tau)\int\limits_t^{\tau}
\Phi_1(\theta)\phi_{j}(\theta)d\theta d\tau\right|\le \frac{\Psi_2(s)}{p^{\beta}}
$$

\vspace{4mm}
\noindent
hold for all $s\in (t, T)$ and for some $\alpha, \beta >0,$ where 
$\Phi_1(\tau),$ $\Phi_2(\tau)$
are continuously differentiable nonrandom functions on $[t, T],$\ $j, p\in \mathbb{N},$
and

\vspace{-1mm}
$$
\int\limits_t^T \Psi_1^2(\tau) d\tau<\infty,\ \ \ 
\int\limits_t^T \left|\Psi_2(\tau)\right| d\tau<\infty.
$$

\vspace{2mm}

{\rm 3.}\ The condition

\vspace{1mm}
$$
\lim\limits_{p\to\infty}
\sum\limits_{\stackrel{j_1,\ldots,j_q,\ldots,j_k=0}{{}_{q\ne g_1, g_2, \ldots, g_{2r-1},
g_{2r}}}}^p
\left(S_{l_1}S_{l_2}\ldots S_{l_{d}}
\left\{\bar C^{(p)}_{j_k\ldots j_q \ldots j_1}\biggl|_{q\ne g_1,g_2,\ldots,g_{2r-1}, g_{2r}}
\right\}\right)^2=0
$$

\vspace{4mm}
\noindent
holds for all possible $g_1,g_2,\ldots,g_{2r-1},g_{2r}$ {\rm (}see {\rm (\ref{leto5007xxx}))}
and $l_1, l_2, \ldots, l_{d}$ such that
$l_1, l_2, \ldots, l_{d}\in \{1,2,$ $\ldots,$ $r\},$\
$l_1>l_2>\ldots >l_{d},$\ $d=0, 1, 2,\ldots, r-1,$\ 
where $r=1, 2,\ldots,[k/2]$ and

\vspace{1mm}
$$
S_{l_1}S_{l_2}\ldots S_{l_{d}}
\left\{\bar C^{(p)}_{j_k\ldots j_q \ldots j_1}\biggl|_{q\ne g_1,g_2,\ldots,g_{2r-1}, g_{2r}}
\right\}\stackrel{\sf def}{=}
\bar C^{(p)}_{j_k\ldots j_q \ldots j_1}\biggl|_{q\ne g_1,g_2,\ldots,g_{2r-1}, g_{2r}}
$$

\vspace{3mm}
\noindent
for $d=0.$

Then, for the iterated Stratonovich stochastic integral 
of arbitrary multiplicity $k$

$$
J^{*}[\psi^{(k)}]_{T,t}^{(i_1\ldots i_k)}=
{\int\limits_t^{*}}^T
\psi_k(t_k) \ldots 
{\int\limits_t^{*}}^{t_2}
\psi_1(t_1) d{\bf w}_{t_1}^{(i_1)}\ldots
d{\bf w}_{t_k}^{(i_k)}
$$

\vspace{2mm}
\noindent
the following 
formula

$$
J^{*}[\psi^{(k)}]_{T,t}^{(i_1\ldots i_k)}=
\hbox{\vtop{\offinterlineskip\halign{
\hfil#\hfil\cr
{\rm l.i.m.}\cr
$\stackrel{}{{}_{p\to \infty}}$\cr
}} }
\int\limits_t^{T}\psi_k(t_k) \ldots 
\int\limits_t^{t_{2}}
\psi_1(t_1) d{\bf w}_{t_1}^{(i_1)p}\ldots
d{\bf w}_{t_k}^{(i_k)p}
$$

\vspace{4mm}
\noindent
is valid$,$ where $i_1,\ldots, i_k=0, 1,\ldots,m$.}

\vspace{2mm}

{\bf Theorem 34.}\ 
{\it Suppose that 
$\{\phi_j(x)\}_{j=0}^{\infty}$ is a complete orthonormal system of 
Legendre polynomials or trigonometric func\-ti\-ons in the space $L_2([t, T]).$
Furthermore, let $\psi_1(\tau), \psi_2(\tau), \psi_3(\tau)$ are continuously dif\-ferentiable 
nonrandom functions on $[t, T].$ 
Then, for the 
iterated Stratonovich stochastic integral of third multiplicity

\vspace{1mm}
$$
J^{*}[\psi^{(3)}]_{T,t}={\int\limits_t^{*}}^T\psi_3(t_3)
{\int\limits_t^{*}}^{t_3}\psi_2(t_2)
{\int\limits_t^{*}}^{t_2}\psi_1(t_1)
d{\bf w}_{t_1}^{(i_1)}
d{\bf w}_{t_2}^{(i_2)}d{\bf w}_{t_3}^{(i_3)}
$$

\vspace{3mm}
\noindent
the following 
equality

$$
J^{*}[\psi^{(3)}]_{T,t}
=\hbox{\vtop{\offinterlineskip\halign{
\hfil#\hfil\cr
{\rm l.i.m.}\cr
$\stackrel{}{{}_{p\to \infty}}$\cr
}} }
\int\limits_t^T\psi_3(t_3)
\int\limits_t^{t_3}\psi_2(t_2)
\int\limits_t^{t_2}\psi_1(t_1)
d{\bf w}_{t_1}^{(i_1)p}
d{\bf w}_{t_2}^{(i_2)p}d{\bf w}_{t_3}^{(i_3)p}
$$

\vspace{4mm}
\noindent
holds$,$ where $i_1, i_3, i_3=0, 1,\ldots,m$.
}

\vspace{2mm}

{\bf Theorem 35.}\ 
{\it Suppose that 
$\{\phi_j(x)\}_{j=0}^{\infty}$ is a complete orthonormal system of 
Legendre polynomials or trigonometric func\-ti\-ons in the space $L_2([t, T]).$
Furthermore, let $\psi_1(\tau), \ldots, \psi_4(\tau)$ are continuously dif\-ferentiable 
nonrandom functions on $[t, T].$ 
Then, for the 
iterated Stratonovich stochastic integral of fourth multiplicity

$$
J^{*}[\psi^{(4)}]_{T,t}={\int\limits_t^{*}}^T\psi_4(t_4)
{\int\limits_t^{*}}^{t_4}\psi_3(t_3)
{\int\limits_t^{*}}^{t_3}\psi_2(t_2)
{\int\limits_t^{*}}^{t_2}\psi_1(t_1)
d{\bf w}_{t_1}^{(i_1)}
d{\bf w}_{t_2}^{(i_2)}d{\bf w}_{t_3}^{(i_3)}d{\bf w}_{t_4}^{(i_4)}
$$

\vspace{4mm}
\noindent
the following 
relation

$$
J^{*}[\psi^{(4)}]_{T,t}
=\hbox{\vtop{\offinterlineskip\halign{
\hfil#\hfil\cr
{\rm l.i.m.}\cr
$\stackrel{}{{}_{p\to \infty}}$\cr
}} }
\int\limits_t^T\psi_4(t_4)
\int\limits_t^{t_4}\psi_3(t_3)
\int\limits_t^{t_3}\psi_2(t_2)
\int\limits_t^{t_2}\psi_1(t_1)
d{\bf w}_{t_1}^{(i_1)p}
d{\bf w}_{t_2}^{(i_2)p}d{\bf w}_{t_3}^{(i_3)p}d{\bf w}_{t_4}^{(i_4)p}
$$

\vspace{4mm}
\noindent
is valid$,$ where $i_1, \ldots, i_4=0, 1,\ldots,m$.
}

\vspace{2mm}

{\bf Theorem 36.}\ 
{\it Suppose that 
$\{\phi_j(x)\}_{j=0}^{\infty}$ is a complete orthonormal system of 
Legendre polynomials or trigonometric func\-ti\-ons in the space $L_2([t, T]).$
Furthermore, let $\psi_1(\tau), \ldots, \psi_5(\tau)$ are continuously dif\-ferentiable 
nonrandom functions on $[t, T].$ 
Then, for the 
iterated Stratonovich stochastic integral of fifth multiplicity

\vspace{1mm}
$$
J^{*}[\psi^{(5)}]_{T,t}=
{\int\limits_t^{*}}^T\psi_5(t_5)
\ldots
{\int\limits_t^{*}}^{t_2}\psi_1(t_1)
d{\bf w}_{t_1}^{(i_1)}
\ldots d{\bf w}_{t_5}^{(i_5)}
$$

\vspace{3mm}
\noindent
the following 
equality

$$
J^{*}[\psi^{(5)}]_{T,t}=\hbox{\vtop{\offinterlineskip\halign{
\hfil#\hfil\cr
{\rm l.i.m.}\cr
$\stackrel{}{{}_{p\to \infty}}$\cr
}} }
\int\limits_t^T\psi_5(t_5)
\ldots
\int\limits_t^{t_2}\psi_1(t_1)
d{\bf w}_{t_1}^{(i_1)p}
\ldots
d{\bf w}_{t_5}^{(i_5)p}
$$

\vspace{4mm}
\noindent
holds$,$ where $i_1, \ldots, i_5=0, 1,\ldots,m$.
}

\vspace{2mm}

{\bf Theorem 37.}\ 
{\it Suppose that 
$\{\phi_j(x)\}_{j=0}^{\infty}$ is a complete orthonormal system of 
Legendre polynomials or trigonometric func\-ti\-ons in the space $L_2([t, T]).$
Then$,$ for the 
iterated Stratonovich stochastic integral of sixth multiplicity 

\vspace{1mm}
$$
J_{T,t}^{*(i_1\ldots i_6)}={\int\limits_t^{*}}^T
{\int\limits_t^{*}}^{t_6}
{\int\limits_t^{*}}^{t_5}
{\int\limits_t^{*}}^{t_4}
{\int\limits_t^{*}}^{t_3}
{\int\limits_t^{*}}^{t_2}
d{\bf w}_{t_1}^{(i_1)}
d{\bf w}_{t_2}^{(i_2)}
d{\bf w}_{t_3}^{(i_3)}
d{\bf w}_{t_4}^{(i_4)}
d{\bf w}_{t_5}^{(i_5)}
d{\bf w}_{t_6}^{(i_6)}
$$

\vspace{3mm}
\noindent
the following 
formula

$$
J_{T,t}^{*(i_1\ldots i_6)}
=\hbox{\vtop{\offinterlineskip\halign{
\hfil#\hfil\cr
{\rm l.i.m.}\cr
$\stackrel{}{{}_{p\to \infty}}$\cr
}} }
\int\limits_t^T
\int\limits_t^{t_6}
\int\limits_t^{t_5}
\int\limits_t^{t_4}
\int\limits_t^{t_3}
\int\limits_t^{t_2}
d{\bf w}_{t_1}^{(i_1)p}
d{\bf w}_{t_2}^{(i_2)p}d{\bf w}_{t_3}^{(i_3)p}d{\bf w}_{t_4}^{(i_4)p}
d{\bf w}_{t_5}^{(i_5)p}d{\bf w}_{t_6}^{(i_6)p}
$$

\vspace{4mm}
\noindent
is valid$,$ where $i_1, \ldots, i_6=0, 1,\ldots,m$.
}

\vspace{2mm}

{\bf Theorem~38.}\ {\it Suppose that 
$\{\phi_j(x)\}_{j=0}^{\infty}$ is a complete orthonormal system of 
Legendre polynomials or trigonometric func\-ti\-ons in the space $L_2([t, T]).$
Furthermore, let $\psi_1(s), \psi_2(s), \psi_3(s)$ are continuously dif\-ferentiable 
nonrandom functions on $[t, T]$.
Then, for the 
iterated Stratonovich stochastic integral of third multiplicity

$$
J^{*}[\psi^{(3)}]_{T,t}={\int\limits_t^{*}}^T\psi_3(t_3)
{\int\limits_t^{*}}^{t_3}\psi_2(t_2)
{\int\limits_t^{*}}^{t_2}\psi_1(t_1)
d{\bf w}_{t_1}^{(i_1)}
d{\bf w}_{t_2}^{(i_2)}d{\bf w}_{t_3}^{(i_3)}
$$

\vspace{3mm}
\noindent
the following 
formula

$$
J^{*}[\psi^{(3)}]_{T,t}
=\hbox{\vtop{\offinterlineskip\halign{
\hfil#\hfil\cr
{\rm l.i.m.}\cr
$\stackrel{}{{}_{p_1,p_2,p_3\to \infty}}$\cr
}} }
\int\limits_t^T\psi_3(t_3)
\int\limits_t^{t_3}\psi_2(t_2)
\int\limits_t^{t_2}\psi_1(t_1)
d{\bf w}_{t_1}^{(i_1)p_1}
d{\bf w}_{t_2}^{(i_2)p_2}d{\bf w}_{t_3}^{(i_3)p_3}
$$

\vspace{3mm}
\noindent
is valid$,$ where $i_1, i_2, i_3=0, 1,\ldots,m$.}

\vspace{2mm}

Let us reformulate Hypotheses 1--3 in terms of 
the convergence of iterated Riemann--Stiltjes
integrals to iterated Stratonovich stochastic integrals.

\vspace{2mm}

{\bf Hypothesis 4}.\ {\it Assume that
$\{\phi_j(x)\}_{j=0}^{\infty}$ is a complete orthonormal
system of Legendre polynomials or trigonometric functions
in the space $L_2([t, T])$.
Then, for the iterated Stratonovich stochastic 
integral of $k$th multiplicity

\vspace{1mm}
$$
{\int\limits_t^{*}}^T
\ldots
{\int\limits_t^{*}}^{t_2}
d{\bf w}_{t_1}^{(i_1)}
\ldots d{\bf w}_{t_k}^{(i_k)}\ \ \ 
(i_1,\ldots, i_k=0, 1,\ldots,m)
$$

\vspace{3mm}
\noindent
the following 
formula

\vspace{0.5mm}
$$
{\int\limits_t^{*}}^T
\ldots
{\int\limits_t^{*}}^{t_2}
d{\bf w}_{t_1}^{(i_1)}
\ldots d{\bf w}_{t_k}^{(i_k)}
=
\hbox{\vtop{\offinterlineskip\halign{
\hfil#\hfil\cr
{\rm l.i.m.}\cr
$\stackrel{}{{}_{p\to \infty}}$\cr
}} }
\int\limits_t^T
\ldots
\int\limits_t^{t_2}
d{\bf w}_{t_1}^{(i_1)p}
\ldots
d{\bf w}_{t_k}^{(i_k)p}
$$

\vspace{3mm}
\noindent
is valid, where  ${\rm l.i.m.}$ is a limit in the mean-square sense,
${\bf w}_{\tau}^{(i)}={\bf f}_{\tau}^{(i)}$ are independent 
standard Wiener processes
$(i=1,\ldots,m)$ and 
${\bf w}_{\tau}^{(0)}=\tau.$}

\vspace{2mm}

{\bf Hypothesis 5}.\ {\it Assume that
$\{\phi_j(x)\}_{j=0}^{\infty}$ is a complete orthonormal
system of Legendre polynomials or trigonometric functions
in the space $L_2([t, T])$. Moreover,
every $\psi_l(\tau)$ $(l=1,\ldots,k)$ is 
an enough smooth nonrandom function
on $[t,T].$
Then, for the iterated Stratonovich stochastic integral {\rm (\ref{str})}
of $k$th multiplicity

\vspace{1mm}
$$
J^{*}[\psi^{(k)}]_{T,t}=
{\int\limits_t^{*}}^T
\psi_k(t_k) \ldots 
{\int\limits_t^{*}}^{t_2}
\psi_1(t_1) d{\bf w}_{t_1}^{(i_1)}
\ldots
d{\bf w}_{t_k}^{(i_k)}\ \ \ (i_1,\ldots, i_k=0, 1,\ldots,m)
$$

\vspace{3mm}
\noindent
the following 
relation

$$
J^{*}[\psi^{(k)}]_{T,t}=
\hbox{\vtop{\offinterlineskip\halign{
\hfil#\hfil\cr
{\rm l.i.m.}\cr
$\stackrel{}{{}_{p\to \infty}}$\cr
}} }
\int\limits_t^T
\psi_k(t_k)
\ldots
\int\limits_t^{t_2}
\psi_1(t_1)
d{\bf w}_{t_1}^{(i_1)p}
\ldots
d{\bf w}_{t_k}^{(i_k)p}
$$

\vspace{3.5mm}
\noindent
holds, where ${\rm l.i.m.}$ is a limit in the mean-square sense,
${\bf w}_{\tau}^{(i)}={\bf f}_{\tau}^{(i)}$ are independent 
standard Wiener processes
$(i=1,2,\ldots,m)$ and 
${\bf w}_{\tau}^{(0)}=\tau.$}

\vspace{2mm}

{\bf Hypothesis 6}.\ {\it Assume that
$\{\phi_j(x)\}_{j=0}^{\infty}$ is a complete orthonormal
system of Legendre polynomials or trigonometric functions
in the space $L_2([t, T])$. Moreover,
every $\psi_l(\tau)$ $(l=1,\ldots,k)$ is 
an enough smooth nonrandom function
on $[t,T].$
Then, for the iterated Stratonovich stochastic integral {\rm (\ref{str})}
of $k$th multiplicity

\vspace{1mm}
$$
J^{*}[\psi^{(k)}]_{T,t}=
{\int\limits_t^{*}}^T
\psi_k(t_k) \ldots 
{\int\limits_t^{*}}^{t_2}
\psi_1(t_1) d{\bf w}_{t_1}^{(i_1)}
\ldots
d{\bf w}_{t_k}^{(i_k)}\ \ \ (i_1,\ldots, i_k=0, 1,\ldots,m)
$$

\vspace{3mm}
\noindent
the following 
equality

$$
J^{*}[\psi^{(k)}]_{T,t}=
\hbox{\vtop{\offinterlineskip\halign{
\hfil#\hfil\cr
{\rm l.i.m.}\cr
$\stackrel{}{{}_{p_1,\ldots,p_k\to \infty}}$\cr
}} }
\int\limits_t^T
\psi_k(t_k)
\ldots
\int\limits_t^{t_2}
\psi_1(t_1)
d{\bf w}_{t_1}^{(i_1)p_1}
\ldots
d{\bf w}_{t_k}^{(i_k)p_k}
$$

\vspace{3.5mm}
\noindent
holds, where ${\rm l.i.m.}$ is a limit in the mean-square sense,
${\bf w}_{\tau}^{(i)}={\bf f}_{\tau}^{(i)}$ are independent 
standard Wiener processes
$(i=1,2,\ldots,m)$ and 
${\bf w}_{\tau}^{(0)}=\tau.$}

\vspace{5mm}

\section{Modification of Condition 3 of Theorem 13 Using Parseval's Equality}

\vspace{5mm}

Let us make some remarks about the development of the 
approach based on Theorem~13 and describe the algorithm
of the verification of Condition~3 of Theorem~13.
First, consider the case $k=2n+1,$ $n=3, 4, \ldots$
($k$ is the multiplicity 
of the iterated Stratonovich stochastic integral (\ref{afterstr})).
Let Conditions 1 and 2 of Theorem~13 be satisfied.
Consider the equality (\ref{drdr1000}). The right-hand side of (\ref{drdr1000})
has the form

\vspace{1mm}
$$
\sum\limits_{j_{g_1}, j_{g_3},\ldots,j_{g_{2r-1}}=0}^p
C_{j_k\ldots j_1}\biggl|_{j_{g_1}=j_{g_2},\ldots, j_{g_{2r-1}}=j_{g_{2r}}}-
$$

\vspace{3mm}
$$
-\frac{1}{2^r} \prod\limits_{l=1}^r {\bf 1}_{\{g_{2l}=g_{2l-1}+1\}}
C_{j_k \ldots j_1}\biggl|_{(j_{g_2} j_{g_1})\curvearrowright (\cdot)
\ldots (j_{g_{2r}} j_{g_{2r-1}})\curvearrowright (\cdot),
j_{g_{{}_{1}}}=~j_{g_{{}_{2}}},\ldots, j_{g_{{}_{2r-1}}}=~j_{g_{{}_{2r}}}
}\biggr..
$$

\vspace{8mm}

Iterated application
of the formulas (\ref{after81}), (\ref{after82}), (\ref{after9031})
separately to the values

\vspace{1mm}
$$
\sum\limits_{j_{g_1}, j_{g_3},\ldots,j_{g_{2r-1}}=0}^p
C_{j_k\ldots j_1}\biggl|_{j_{g_1}=j_{g_2},\ldots, j_{g_{2r-1}}=j_{g_{2r}}}
$$

\vspace{1mm}
and 

$$
\frac{1}{2^r} \prod\limits_{l=1}^r {\bf 1}_{\{g_{2l}=g_{2l-1}+1\}}
C_{j_k \ldots j_1}\biggl|_{(j_{g_2} j_{g_1})\curvearrowright (\cdot)
\ldots (j_{g_{2r}} j_{g_{2r-1}})\curvearrowright (\cdot),
j_{g_{{}_{1}}}=~j_{g_{{}_{2}}},\ldots, j_{g_{{}_{2r-1}}}=~j_{g_{{}_{2r}}}
}\biggr.
$$

\vspace{8mm}
\noindent
($g_1, g_2,\ldots, g_{2r-1}, g_{2r}$ as in (\ref{leto5007xxx}), $r=1,2,\ldots,[k/2],$ $2r<k$)
gives the following representation (see (\ref{drdr1001}))

\vspace{1mm}
$$
\sum\limits_{\stackrel{j_1,\ldots,j_q,\ldots,j_k=0}{{}_{q\ne g_1, g_2, \ldots, g_{2r-1},
g_{2r}}}}^p
\Biggl(\sum\limits_{j_{g_1}, j_{g_3},\ldots,j_{g_{2r-1}}=0}^p
C_{j_k\ldots j_1}\biggl|_{j_{g_1}=j_{g_2},\ldots, j_{g_{2r-1}}=j_{g_{2r}}}-\Biggr.
$$

\vspace{3mm}
$$
\Biggl.-\frac{1}{2^r} \prod\limits_{l=1}^r {\bf 1}_{\{g_{2l}=g_{2l-1}+1\}}
C_{j_k \ldots j_1}\biggl|_{(j_{g_2} j_{g_1})\curvearrowright (\cdot)
\ldots (j_{g_{2r}} j_{g_{2r-1}})\curvearrowright (\cdot),
j_{g_{{}_{1}}}=~j_{g_{{}_{2}}},\ldots, j_{g_{{}_{2r-1}}}=~j_{g_{{}_{2r}}}
}\biggr.\Biggr)^2\le
$$

\vspace{8mm}
$$
\le \sum\limits_{\stackrel{j_1,\ldots,j_q,\ldots,j_k=0}{{}_{q\ne g_1, g_2, \ldots, g_{2r-1},
g_{2r}}}}^{\infty}
\Biggl(\sum\limits_{j_{g_1}, j_{g_3},\ldots,j_{g_{2r-1}}=0}^p
C_{j_k\ldots j_1}\biggl|_{j_{g_1}=j_{g_2},\ldots, j_{g_{2r-1}}=j_{g_{2r}}}-\Biggr.
$$

\vspace{3mm}
$$
\Biggl.-\frac{1}{2^r} \prod\limits_{l=1}^r {\bf 1}_{\{g_{2l}=g_{2l-1}+1\}}
C_{j_k \ldots j_1}\biggl|_{(j_{g_2} j_{g_1})\curvearrowright (\cdot)
\ldots (j_{g_{2r}} j_{g_{2r-1}})\curvearrowright (\cdot),
j_{g_{{}_{1}}}=~j_{g_{{}_{2}}},\ldots, j_{g_{{}_{2r-1}}}=~j_{g_{{}_{2r}}}
}\biggr.\Biggr)^2=
$$

\vspace{7mm}
$$
=\sum\limits_{\stackrel{j_1,\ldots,j_q,\ldots,j_k=0}{{}_{q\ne g_1, g_2, \ldots, g_{2r-1},
g_{2r}}}}^{\infty}
\left(~
\int\limits_{[t, T]^{k-2r}} 
R_p(t_1,\ldots,t_{g_1-1},t_{g_1+1},\ldots, t_{g_{2r}-1},t_{g_{2r}+1},\ldots,t_k)\times\right.
$$

\vspace{3mm}
\begin{equation}
\label{pars0}
\left.\times 
\prod_{\stackrel{q=1}{{}_{q\ne g_1, g_2, \ldots, g_{2r-1},
g_{2r}}}}^{k}
\psi_{q}(t_q)\phi_{j_q}(t_q)\ 
dt_1\ldots dt_{g_1-1}dt_{g_1+1}\ldots dt_{g_{2r}-1}dt_{g_{2r}+1}\ldots dt_k\right)^2,
\end{equation}

\vspace{4mm}
\noindent
where 

\vspace{-3mm}
$$
R_p(t_1,\ldots,t_{g_1-1},t_{g_1+1},\ldots, t_{g_{2r}-1},t_{g_{2r}+1},\ldots,t_k)=
$$

\vspace{2mm}
$$
=\sum\limits_{d=1}^{4^r}
\bar R_p^{(d)}(t_1,\ldots,t_{g_1-1},t_{g_1+1},\ldots, t_{g_{2r}-1},t_{g_{2r}+1},\ldots,t_k)-
$$

\vspace{2mm}
$$
-
\sum\limits_{d=1}^{2^r}
\tilde R_p^{(d)}(t_1,\ldots,t_{g_1-1},t_{g_1+1},\ldots, t_{g_{2r}-1},t_{g_{2r}+1},\ldots,t_k)
\in L_2([t, T]^{k-2r})
$$

\vspace{5mm}
\noindent
and 
$$
\int\limits_{[t, T]^{k-2r}} 
R_p(t_1,\ldots,t_{g_1-1},t_{g_1+1},\ldots, 
t_{g_{2r}-1},t_{g_{2r}+1},\ldots,t_k)
\times
$$

$$
\times 
\prod_{\stackrel{q=1}{{}_{q\ne g_1, g_2, \ldots, g_{2r-1},
g_{2r}}}}^{k}
\psi_{q}(t_q)\phi_{j_q}(t_q)\ 
dt_1\ldots dt_{g_1-1}dt_{g_1+1}\ldots dt_{g_{2r}-1}dt_{g_{2r}+1}\ldots dt_k
$$

\vspace{5mm}
\noindent
is the Fourier coefficient of 

\vspace{2mm}
$$
\hat R_p(t_1,\ldots,t_{g_1-1},t_{g_1+1},\ldots, t_{g_{2r}-1},t_{g_{2r}+1},\ldots,t_k)=
$$

\vspace{2mm}
$$
=
R_p(t_1,\ldots,t_{g_1-1},t_{g_1+1},\ldots, t_{g_{2r}-1},t_{g_{2r}+1},\ldots,t_k)
\prod_{\stackrel{q=1}{{}_{q\ne g_1, g_2, \ldots, g_{2r-1},
g_{2r}}}}^{k}
\psi_{q}(t_q).
$$

\vspace{4mm}

Also note that some of the functions

$$
\bar R_p^{(d)}(t_1,\ldots,t_{g_1-1},t_{g_1+1},\ldots, t_{g_{2r}-1},t_{g_{2r}+1},\ldots,t_k)
$$

\vspace{2mm}
\noindent
and 
$$
\tilde R_p^{(d)}(t_1,\ldots,t_{g_1-1},t_{g_1+1},\ldots, t_{g_{2r}-1},t_{g_{2r}+1},\ldots,t_k)
$$

\vspace{5mm}
\noindent
can be identically equal to zero.

Obviously, we could use another representation for the function

\vspace{-1mm}
\begin{equation}
\label{de200}
R_p(t_1,\ldots,t_{g_1-1},t_{g_1+1},\ldots, t_{g_{2r}-1},t_{g_{2r}+1},\ldots,t_k)
\end{equation}

\vspace{4mm}
\noindent
based on the left-hand side of the equality (\ref{drdr1000})
and (\ref{after81}), (\ref{after82}), (\ref{after9031}) (see Sect.~7, 10 for details).
In Sect.~10, we considered the function (\ref{de200}) in detail
for the case $k\ge 5,$ $r=1.$

Parseval's equality gives

\vspace{1mm}
$$
\sum\limits_{\stackrel{j_1,\ldots,j_q,\ldots,j_k=0}{{}_{q\ne g_1, g_2, \ldots, g_{2r-1},
g_{2r}}}}^{\infty}
\left(~
\int\limits_{[t, T]^{k-2r}} 
R_p(t_1,\ldots,t_{g_1-1},t_{g_1+1},\ldots, t_{g_{2r}-1},t_{g_{2r}+1},\ldots,t_k)\times\right.
$$

\vspace{3mm}
$$
\left.\times 
\prod_{\stackrel{q=1}{{}_{q\ne g_1, g_2, \ldots, g_{2r-1},
g_{2r}}}}^{k}
\psi_{q}(t_q)\phi_{j_q}(t_q)\ 
dt_1\ldots dt_{g_1-1}dt_{g_1+1}\ldots dt_{g_{2r}-1}dt_{g_{2r}+1}\ldots dt_k\right)^2=
$$

\vspace{7mm}
$$
=\int\limits_{[t, T]^{k-2r}} 
\left(
\hat R_p(t_1,\ldots,t_{g_1-1},t_{g_1+1},\ldots, t_{g_{2r}-1},t_{g_{2r}+1},\ldots,t_k)\right)^2
\times
$$

\vspace{3mm}
$$
\times
dt_1\ldots dt_{g_1-1}dt_{g_1+1}\ldots dt_{g_{2r}-1}dt_{g_{2r}+1}\ldots dt_k =
$$

\vspace{3mm}
\begin{equation}
\label{pars1}
=\bigl\Vert \hat R_p \bigr\Vert_{L_2([t, T]^{k-2r})}^2.
\end{equation}

\vspace{7mm}

Combining (\ref{pars0}) and (\ref{pars1}), we obtain

$$
\sum\limits_{\stackrel{j_1,\ldots,j_q,\ldots,j_k=0}{{}_{q\ne g_1, g_2, \ldots, g_{2r-1},
g_{2r}}}}^p
\Biggl(\sum\limits_{j_{g_1}, j_{g_3},\ldots,j_{g_{2r-1}}=0}^p
C_{j_k\ldots j_1}\biggl|_{j_{g_1}=j_{g_2},\ldots, j_{g_{2r-1}}=j_{g_{2r}}}-\Biggr.
$$

\vspace{3mm}
$$
\Biggl.-\frac{1}{2^r} \prod\limits_{l=1}^r {\bf 1}_{\{g_{2l}=g_{2l-1}+1\}}
C_{j_k \ldots j_1}\biggl|_{(j_{g_2} j_{g_1})\curvearrowright (\cdot)
\ldots (j_{g_{2r}} j_{g_{2r-1}})\curvearrowright (\cdot),
j_{g_{{}_{1}}}=~j_{g_{{}_{2}}},\ldots, j_{g_{{}_{2r-1}}}=~j_{g_{{}_{2r}}}
}\biggr.\Biggr)^2\le
$$

\vspace{4mm}
\begin{equation}
\label{pars2}
\le \bigl\Vert \hat R_p \bigr\Vert_{L_2([t, T]^{k-2r})}^2.
\end{equation}

\vspace{6mm}

Assume that we have succeeded in proving the following equality

\begin{equation}
\label{pars3s}
\lim\limits_{p\to\infty}
\bigl\Vert \hat R_p \bigr\Vert_{L_2([t, T]^{k-2r})}^2=0.
\end{equation}

\vspace{4mm}

Applying (\ref{pars2}) and (\ref{pars3s}), we get (compare with (\ref{drdr1001}))

$$
\lim\limits_{p\to\infty}
\sum\limits_{\stackrel{j_1,\ldots,j_q,\ldots,j_k=0}{{}_{q\ne g_1, g_2, \ldots, g_{2r-1},
g_{2r}}}}^p
\Biggl(\sum\limits_{j_{g_1}, j_{g_3},\ldots,j_{g_{2r-1}}=0}^p
C_{j_k\ldots j_1}\biggl|_{j_{g_1}=j_{g_2},\ldots, j_{g_{2r-1}}=j_{g_{2r}}}-\Biggr.
$$

\vspace{3mm}
\begin{equation}
\label{for900}
\Biggl.-\frac{1}{2^r} \prod\limits_{l=1}^r {\bf 1}_{\{g_{2l}=g_{2l-1}+1\}}
C_{j_k \ldots j_1}\biggl|_{(j_{g_2} j_{g_1})\curvearrowright (\cdot)
\ldots (j_{g_{2r}} j_{g_{2r-1}})\curvearrowright (\cdot),
j_{g_{{}_{1}}}=~j_{g_{{}_{2}}},\ldots, j_{g_{{}_{2r-1}}}=~j_{g_{{}_{2r}}}
}\biggr.\Biggr)^2=0.
\end{equation}

\vspace{5mm}

As noted in Sect.~7, Condition 3 of Theorem~13 can be replaced by a weaker condition 
(\ref{drdr1001}) (or (\ref{for900})).
Also Condition 3 of Theorem~13 can be replaced by (\ref{pars3s}).
From (\ref{for900}) we obviously obtain 

$$
\lim\limits_{p\to\infty} \sum\limits_{j_{g_1}, j_{g_3},\ldots,j_{g_{2r-1}}=0}^p
C_{j_k\ldots j_1}\biggl|_{j_{g_1}=j_{g_2},\ldots, j_{g_{2r-1}}=j_{g_{2r}}}=\Biggr.
$$

\vspace{3mm}
\begin{equation}
\label{pars900}
\Biggl.=\frac{1}{2^r} \prod\limits_{l=1}^r {\bf 1}_{\{g_{2l}=g_{2l-1}+1\}}
C_{j_k \ldots j_1}\biggl|_{(j_{g_2} j_{g_1})\curvearrowright (\cdot)
\ldots (j_{g_{2r}} j_{g_{2r-1}})\curvearrowright (\cdot),
j_{g_{{}_{1}}}=~j_{g_{{}_{2}}},\ldots, j_{g_{{}_{2r-1}}}=~j_{g_{{}_{2r}}}
}\biggr..
\end{equation}

\vspace{7mm}

According to (\ref{drdr1000}), the equality (\ref{pars900}) will be satisfied 
if

\begin{equation}
\label{sars10}
\lim\limits_{p\to\infty}
S_{l_1}S_{l_2}\ldots S_{l_{d}}
\left\{\bar C^{(p)}_{j_k\ldots j_q \ldots j_1}\biggl|_{q\ne g_1,g_2,\ldots,g_{2r-1}, g_{2r}}
\right\}=0,
\end{equation}

\vspace{5mm}
\noindent
where $g_1,g_2,\ldots,g_{2r-1},g_{2r}$ as in (\ref{leto5007xxx}),
$l_1, l_2, \ldots, l_{d}$ such that
$l_1, l_2, \ldots, l_{d}\in \{1,2,\ldots, r\},$\
$l_1>l_2>\ldots >l_{d},$\ $d=0, 1, 2,\ldots, r-1,$\ 
$r=1, 2,\ldots,[k/2],$

\vspace{3mm}
$$
S_{l_1}S_{l_2}\ldots S_{l_{d}}
\left\{\bar C^{(p)}_{j_k\ldots j_q \ldots j_1}\biggl|_{q\ne g_1,g_2,\ldots,g_{2r-1}, g_{2r}}
\right\}\stackrel{\sf def}{=}
\bar C^{(p)}_{j_k\ldots j_q \ldots j_1}\biggl|_{q\ne g_1,g_2,\ldots,g_{2r-1}, g_{2r}}
$$

\vspace{5mm}
\noindent
for $d=0,$ where 

$$
\bar C^{(p)}_{j_k\ldots j_q \ldots j_1}\biggl|_{q\ne g_1,g_2,\ldots,g_{2r-1}, g_{2r}},\ \ \
S_l \left\{\bar C^{(p)}_{j_k\ldots j_q \ldots j_1}\biggl|_{q\ne g_1,g_2,\ldots,g_{2r-1}, g_{2r}}
\right\}$$

\vspace{5mm}
\noindent
are defined by (\ref{nov100}), (\ref{nov101}), 
$l=1,2,\ldots,r$ (see Sect.~7 for details).

Let us make some remarks about the function
(\ref{de200})
for the case $k>5,$ $r=2.$
In this case, using the left-hand side of 
the equality (\ref{drdr1000})
and (\ref{after81}), (\ref{after82}), (\ref{after9031}), 
we represent the function (\ref{de200})
as the sum of several functions.
In particular, among these functions 
will be the following functions

$$
Q_p(t_1,\ldots,t_{s-1},t_{s+1},\ldots ,t_{l-1},t_{l+1},\ldots,
t_{q-1},t_{q+1},\ldots, t_{g-1},t_{g+1},\ldots,t_k)=
$$

\vspace{-1mm}
$$
=
{\bf 1}_{\{t_1< \ldots <t_{s-1}<t_{s+1}< \ldots <t_{l-1}<t_{l+1}< \ldots
<t_{q-1}<t_{q+1}< \ldots <t_{g-1}<t_{g+1}<\ldots<t_k\}}\times
$$

$$
\times
\sum_{j_l=p+1}^{\infty}~
\int\limits_{t}^{t_{s+1}} \psi_s(\tau) \phi_{j_{l}}(\tau)d\tau
\int\limits_{t}^{t_{l-1}} \psi_l(\tau) \phi_{j_{l}}(\tau)d\tau\times
$$

\vspace{1mm}
\begin{equation}
\label{func500}
\times
\sum_{j_q=p+1}^{\infty}~
\int\limits_{t}^{t_{q+1}} \psi_q(\tau) \phi_{j_{q}}(\tau)d\tau
\int\limits_{t}^{t_{g-1}} \psi_g(\tau) \phi_{j_{q}}(\tau)d\tau,
\end{equation}

\vspace{6mm}

$$
\bar Q_p(t_1,\ldots,t_{l-2},t_{l+3},\ldots,t_k)=
{\bf 1}_{\{t_1<\ldots<t_{l-2}<t_{l+3}<\ldots<t_k\}}\times
$$

\vspace{-1mm}
$$
\times
\sum_{j_l=p+1}^{\infty}~\left(\int\limits_{t}^{t_{l-2}} \psi_{l-1}(\theta)\phi_{j_{l}}(\theta)
\int\limits_{t}^{\theta} \psi_l(u)\phi_{j_{l}}(u)du d\theta\right)\times
$$

\vspace{1mm}
\begin{equation}
\label{func501}
\times
\sum_{j_q=p+1}^{\infty}~\left(\int\limits_{t}^{t_{l-2}} \psi_{l+1}(\theta)\phi_{j_{q}}(\theta)
\int\limits_{t}^{\theta} \psi_{l+2}(u)\phi_{j_{q}}(u)du d\theta\right),
\end{equation}

\vspace{7mm}

$$
\tilde Q_p(t_1,\ldots,t_{l-2},t_{l+3},\ldots,t_k)=
{\bf 1}_{\{t_1<\ldots<t_{l-2}<t_{l+3}<\ldots<t_k\}}\times
$$

\vspace{-1mm}
$$
\times
\sum_{j_l=p+1}^{\infty}\sum_{j_q=p+1}^{\infty}~
\int\limits_{t}^{t_{l+3}} \psi_{l+1}(\tau)\phi_{j_{q}}(\tau)
\left(\int\limits_{t}^{\tau} \psi_{l-1}(\theta)\phi_{j_{l}}(\theta)
\int\limits_{t}^{\theta} \psi_l(u)\phi_{j_{l}}(u)du d\theta\right)\times
$$

\vspace{1mm}
\begin{equation}
\label{func502}
\times
\int\limits_{t}^{\tau} \psi_{l+2}(u)\phi_{j_{q}}(u)du d\tau,
\end{equation}

\vspace{6mm}

$$
\hat Q_p(t_1,\ldots,t_{l-1},t_{l+2},\ldots, t_{q-1}, t_{q+2}, \ldots, t_k)=
$$

\vspace{-1mm}
$$
=
{\bf 1}_{\{t_1<\ldots<t_{l-1}<t_{l+2}<\ldots<t_{q-1}<t_{q+2}<\ldots <t_k\}}\times
$$

\vspace{1mm}
$$
\times
\sum_{j_{l}=p+1}^{\infty}
\sum_{j_{l+1}=p+1}^{\infty}~\left(\int\limits_{t}^{t_{l+2}} \psi_{l+1}(\theta)\phi_{j_{l+1}}(\theta)
\int\limits_{t}^{\theta} \psi_l(u)\phi_{j_{l}}(u)du d\theta\right)\times
$$

\vspace{1mm}
\begin{equation}
\label{func501qq}
\times
\left(\int\limits_{t}^{t_{q+2}} \psi_{q+1}(\theta)\phi_{j_{l+1}}(\theta)
\int\limits_{t}^{\theta} \psi_{q}(u)\phi_{j_{l}}(u)du d\theta\right).
\end{equation}

\vspace{4mm}

Note that the pairs $(g_1,g_2),$ $(g_3,g_4)$ for the functions (\ref{func501}) and (\ref{func502}) 
have the property: $g_2=g_1+1,$ $g_4=g_3+1,$ $g_3=g_2+1.$
At the same time, the pairs $(g_1,g_2),$ $(g_3,g_4)$ for the function (\ref{func500})
have the following property: $g_2>g_1+1,$ $g_4>g_3+1,$ $g_3\ge g_2+1.$
For the function (\ref{func501qq}), the pairs $(g_1,g_2),$ $(g_3,g_4)$
chosen as follows: $g_2>g_1+1,$ $g_4>g_3+1,$ $g_4=g_2+1,$ $g_3=g_1+1.$
Generally speaking, all possible pairs $(g_1,g_2),$ $(g_3,g_4)$ must be considered.
We consider the functions (\ref{func500})--(\ref{func501qq}) only as an example.

Suppose that $s+1=l-1,$ $l+1=q-1,$ $q+1=g-1$ in (\ref{func500}). Let us show that
(we consider the case of Legendre polynomials; the trigonometric case is simpler
and can be considered similarly)

\vspace{-1mm}
\begin{equation}
\label{func503}
\lim\limits_{p\to\infty}\bigl\Vert Q_p \bigr\Vert_{L_2([t, T]^{k-4})}^2=0,
\end{equation}

\begin{equation}
\label{func504}
\lim\limits_{p\to\infty}\bigl\Vert \bar Q_p \bigr\Vert_{L_2([t, T]^{k-4})}^2=0,
\end{equation}

\begin{equation}
\label{func505}
\lim\limits_{p\to\infty}\bigl\Vert \tilde Q_p \bigr\Vert_{L_2([t, T]^{k-4})}^2=0,
\end{equation}

\begin{equation}
\label{func505qq}
\lim\limits_{p\to\infty}\bigl\Vert \hat Q_p \bigr\Vert_{L_2([t, T]^{k-4})}^2=0.
\end{equation}

\vspace{4mm}

First consider the proof of (\ref{func503}). We have
($s+1=l-1,$ $l+1=q-1,$ $q+1=g-1$)

\vspace{2mm}
$$
\left(Q_p(t_1,\ldots,t_{l-3},t_{l-1},t_{l+1},t_{l+3},t_{l+5},\ldots,t_k)\right)^2=
$$

\vspace{1mm}
$$
=
{\bf 1}_{\{t_1< \ldots <t_{l-3}<t_{l-1}< t_{l+1}<t_{l+3}<t_{l+5}<\ldots<t_k\}}\times
$$

\vspace{1mm}
$$
\times
\left(\sum_{j_l=p+1}^{\infty}~
\int\limits_{t}^{t_{l-1}} \psi_{l-2}(\tau) \phi_{j_{l}}(\tau)d\tau
\int\limits_{t}^{t_{l-1}} \psi_l(\tau) \phi_{j_{l}}(\tau)d\tau\times\right.
$$

\vspace{1mm}
\begin{equation}
\label{func506}
\left.\times \sum_{j_q=p+1}^{\infty}~
\int\limits_{t}^{t_{l+3}} \psi_{l+2}(\tau) \phi_{j_{q}}(\tau)d\tau
\int\limits_{t}^{t_{l+3}} \psi_{l+4}(\tau) \phi_{j_{q}}(\tau)d\tau\right)^2.
\end{equation}

\vspace{6mm}

Using the estimate (\ref{after1940}), we obtain

\begin{equation}
\label{func507}
\left|
\int\limits_t^s\psi(\tau)\phi_{j}(\tau)d\tau
\right| <
\frac{K}{j^{1-\varepsilon/2} (1-z^2(s))^{1/4-\varepsilon/4}},
\end{equation}

\vspace{4mm}
\noindent
where $j\in\mathbb{N},$ $s\in (t, T),$
$z(s)$ is defined by (\ref{zz1}), $\varepsilon\in (0,1),$ constant $K$ does not depend on $j,$
$\{\phi_j(x)\}_{j=0}^{\infty}$ is a complete orthonormal system of 
Legendre polynomials in the space $L_2([t, T]),$
$\psi(\tau)$ is a continuously dif\-ferentiable 
nonrandom function on $[t, T].$ 

Applying (\ref{func507}) and 
(\ref{after1944}) (we take $\varepsilon$ instead of $\varepsilon/2$ in (\ref{after1944})), we get

$$
\left(\sum_{j_l=p+1}^{\infty}~
\int\limits_{t}^{t_{l-1}} \psi_{l-2}(\tau) \phi_{j_{l}}(\tau)d\tau
\int\limits_{t}^{t_{l-1}} \psi_l(\tau) \phi_{j_{l}}(\tau)d\tau\times\right.
$$

\vspace{1mm}
$$
\left.\times \sum_{j_q=p+1}^{\infty}~
\int\limits_{t}^{t_{l+3}} \psi_{l+2}(\tau) \phi_{j_{q}}(\tau)d\tau
\int\limits_{t}^{t_{l+3}} \psi_{l+4}(\tau) \phi_{j_{q}}(\tau)d\tau\right)^2\le
$$

\vspace{1mm}
\begin{equation}
\label{func508}
\le \frac{K_1}{p^{4(1-\varepsilon)}(1-z^2(t_{l-1}))^{1-\varepsilon}
(1-z^2(t_{l+3}))^{1-\varepsilon}},
\end{equation}

\vspace{4mm}
\noindent
where $t_{l-1}, t_{l+3}\in (t, T),$ constant $K_1$ is independent of $p.$
Combining (\ref{func506}) and (\ref{func508}), we have (\ref{func503}).

Let us prove (\ref{func504}).
The following equality is proved in Sect.~12 \cite{arxiv-5} (also see Sect.~2.9 \cite{2018a})
for the case of Legendre polynomials
($n>m;$\ $n, m\in \mathbb{N}$)

$$
\sum\limits_{j=m+1}^n
C_{j j}(s)=
\sum\limits_{j=m+1}^n
\int\limits_t^s \psi_2(\theta)\phi_{j}(\theta)
\int\limits_t^{\theta} \psi_1(\tau)\phi_{j}(\tau)d\tau d\theta=
$$

\vspace{2mm}
$$
=
\frac{T-t}{4}
\int\limits_{-1}^{z(s)}
\psi_1(u(x))\psi_2(u(x))
\left(P_{n+1}(x)P_{n}(x)
-P_{m+1}(x)P_{m}(x)\right)dx-
$$

\vspace{2mm}
$$
-\frac{(T-t)^2}{8}
\sum\limits_{j=m+1}^n
\frac{1}{2j+1}
\int\limits_{-1}^{z(s)}
\left(P_{j+1}(y)-P_{j-1}(y)\right)\psi_1'(u(y))\times
$$

\vspace{2mm}
$$
\times
\Biggl(\left(P_{j+1}(z(s))-P_{j-1}(z(s))\right)\psi_2(s)-
\left(P_{j+1}(y)-P_{j-1}(y)\right)\psi_2(u(y))-\Biggr.
$$

\vspace{2mm}
\begin{equation}
\label{agent14}
\Biggl.
-
\frac{T-t}{2}
\int\limits_{y}^{z(s)}
\left(P_{j+1}(x)-
P_{j-1}(x)\right)\psi_2'(u(x))dx\Biggr)dy,
\end{equation}

\vspace{5mm}
\noindent
where $s\in (t, T),$ 
$$
C_{j j}(s)=\int\limits_t^s\psi_2(\tau)\phi_{j}(\tau)
\int\limits_t^{\tau}\psi_1(\theta)\phi_{j}(\theta)d\theta d\tau,
$$

\vspace{4mm}
$$
u(y)=\frac{T-t}{2}y+\frac{T+t}{2},\ \ \
z(s)=\left(s-\frac{T+t}{2}\right)\frac{2}{T-t},
$$

\vspace{5mm}
\noindent
and $\psi_1'$, $\psi_2'$ are
derivatives of the functions $\psi_1(\tau)$, $\psi_2(\tau)$ with respect 
to the variable
$u(y)$.

Applying the estimate (\ref{after5000}) in (\ref{agent14}) and tak\-ing into account 
the boundedness of the functions $\psi_1(\tau),$ $\psi_2(\tau)$ 
and their derivatives, we obtain

$$
\left\vert\sum\limits_{j=m+1}^n
C_{j j}(s)\right\vert\le
C_1\left(\frac{1}{n^{1-\varepsilon}}+\frac{1}{m^{1-\varepsilon}}\right)
\int\limits_{-1}^{z(s)} \frac{dx}{\left(1-x^2\right)^{1/2-\varepsilon/2}}+
$$

\vspace{2mm}
$$
+C_2 \sum\limits_{j=m+1}^n \frac{1}{j^{2-\varepsilon}}\left(
\int\limits_{-1}^{z(s)}\frac{dy}{\left(1-y^2\right)^{1/2-\varepsilon/2}}
+\frac{1}{\left(1-z^2(s)\right)^{1/4-\varepsilon/4}}
\int\limits_{-1}^{z(s)}\frac{dy}{\left(1-y^2\right)^{1/4-\varepsilon/4}}+\right.
$$

\vspace{2mm}
\begin{equation}
\label{func800}
+
\left.\int\limits_{-1}^{z(s)}\frac{1}{\left(1-y^2\right)^{1/4-\varepsilon/4}}
\int\limits_{y}^{z(s)}\frac{dx}{\left(1-x^2\right)^{1/4-\varepsilon/4}}dy\right),
\end{equation}

\vspace{5mm}
\noindent
where 
$s\in (t, T),$ constants $C_1, C_2$ do not depend on $n$ and $m$.

From (\ref{func800}) we have

\begin{equation}
\label{func800x}
\left\vert\sum\limits_{j=m+1}^{\infty}
C_{j j}(s)\right\vert\le
\frac{K_1}{m^{1-\varepsilon}}
+K_2 \sum\limits_{j=m+1}^{\infty} \frac{1}{j^{2-\varepsilon}}\left(
1+\frac{1}{\left(1-z^2(s)\right)^{1/4-\varepsilon/4}}\right),
\end{equation}

\vspace{4mm}
\noindent
where $s\in (t, T),$ constants $K_1, K_2$ do not depend on $m$.

Applying (\ref{after1944}) (we take $\varepsilon$ instead of $\varepsilon/2$ 
in (\ref{after1944})) in (\ref{func800x}),  we get

\begin{equation}
\label{func801}
\left\vert\sum\limits_{j=m+1}^{\infty}
C_{j j}(s)\right\vert\le
\frac{K}{m^{1-\varepsilon}\left(1-z^2(s)\right)^{1/4-\varepsilon/4}},
\end{equation}

\vspace{4mm}
\noindent
where $s\in (t, T),$ constant $K$ is independent of $m$.

Using the estimate (\ref{func801}), we obtain (see (\ref{func501}))

\vspace{2mm}
$$
\left(\bar Q_p(t_1,\ldots,t_{l-2},t_{l+3},\ldots,t_k)\right)^2=
{\bf 1}_{\{t_1<\ldots<t_{l-2}<t_{l+3}<\ldots<t_k\}}\times
$$

\vspace{1mm}
$$
\times
\left(\sum_{j_l=p+1}^{\infty}~\left(\int\limits_{t}^{t_{l-2}} \psi_{l-1}(\theta)\phi_{j_{l}}(\theta)
\int\limits_{t}^{\theta} \psi_l(u)\phi_{j_{l}}(u)du d\theta\right)\times\right.
$$

\vspace{1mm}
$$
\left.
\times
\sum_{j_q=p+1}^{\infty}~\left(\int\limits_{t}^{t_{l-2}} \psi_{l+1}(\theta)\phi_{j_{q}}(\theta)
\int\limits_{t}^{\theta} \psi_{l+2}(u)\phi_{j_{q}}(u)du d\theta\right)\right)^2\le
$$

\vspace{2mm}

\begin{equation}
\label{func990}
\le \frac{K_1}{p^{4(1-\varepsilon)}(1-z^2(t_{l-2}))^{1-\varepsilon}},
\end{equation}

\vspace{5mm}
\noindent
where $t_{l-2}\in (t, T),$ constant $K_1$ is independent of $p.$
The inequality (\ref{func990}) completes the proof of (\ref{func504}).

Let us prove (\ref{func505}). 
Using (\ref{after1400}),
we obtain the following equality 
for the cases of Legendre polynomials and trigonometric functions 

\vspace{-1mm}
\begin{equation}
\label{agent0101}
\frac{1}{2}
\int\limits_t^s\psi_1(t_1)\psi_2(t_1)dt_1
-\sum_{j_1=0}^{p}
C_{j_1j_1}(s)=\sum_{j_1=p+1}^{\infty}
C_{j_1j_1}(s),
\end{equation}

\vspace{4mm}
\noindent
where $s\in (t, T)$ and 
$$
C_{j j}(s)=\int\limits_t^s\psi_2(\tau)\phi_{j}(\tau)
\int\limits_t^{\tau}\psi_1(\theta)\phi_{j}(\theta)d\theta d\tau.
$$

\vspace{4mm}

Applying (\ref{agent0101}) in (\ref{func502}), we get

\vspace{2mm}
$$
\left(\tilde Q_p(t_1,\ldots,t_{l-2},t_{l+3},\ldots,t_k)\right)^2\le
$$

\vspace{2mm}
$$
\le
\left(\sum_{j_l=p+1}^{\infty}\sum_{j_q=p+1}^{\infty}~
\int\limits_{t}^{t_{l+3}} \psi_{l+1}(\tau)\phi_{j_{q}}(\tau)
\left(\int\limits_{t}^{\tau} \psi_{l-1}(\theta)\phi_{j_{l}}(\theta)
\int\limits_{t}^{\theta} \psi_l(u)\phi_{j_{l}}(u)du d\theta\right)\times\right.
$$

\vspace{2mm}
$$
\left.\times
\int\limits_{t}^{\tau} \psi_{l+2}(u)\phi_{j_{q}}(u)du d\tau\right)^2=
$$

\vspace{2mm}
$$
=\left(\frac{1}{2}\sum_{j_l=p+1}^{\infty}~
\int\limits_{t}^{t_{l+3}} \psi_{l+1}(\tau)
\left(\int\limits_{t}^{\tau} \psi_{l-1}(\theta)\phi_{j_{l}}(\theta)
\int\limits_{t}^{\theta} \psi_l(u)\phi_{j_{l}}(u)du d\theta\right)
\psi_{l+2}(\tau)d\tau-\right.
$$

\vspace{2mm}
$$
-\sum_{j_q=0}^{p}~
\int\limits_{t}^{t_{l+3}}\psi_{l+1}(\tau)\phi_{j_{q}}(\tau)
\sum_{j_l=p+1}^{\infty}\left(\int\limits_{t}^{\tau} \psi_{l-1}(\theta)\phi_{j_{l}}(\theta)
\int\limits_{t}^{\theta}\psi_l(u)\phi_{j_{l}}(u)du d\theta\right)\times
$$

\vspace{2mm}
$$
\left.\times
\int\limits_{t}^{\tau} \psi_{l+2}(u)\phi_{j_{q}}(u)du d\tau\right)^2= 
$$

\vspace{2mm}
\begin{equation}
\label{func1000}
=(a-b)^2\le
2(|a|^2 + |b|^2).
\end{equation}

\vspace{5mm}

Further, we have

\begin{equation}
\label{func1001}
|a|\le 
\frac{1}{2}
\int\limits_{t}^{t_{l+3}} \left\vert\psi_{l+1}(\tau)\right\vert
\left\vert\sum_{j_l=p+1}^{\infty}~\int\limits_{t}^{\tau} \psi_{l-1}(\theta)\phi_{j_{l}}(\theta)
\int\limits_{t}^{\theta} \psi_l(u)\phi_{j_{l}}(u)du d\theta\right\vert
\left\vert\psi_{l+2}(\tau)\right\vert d\tau,
\end{equation}

\vspace{4mm}
$$
|b|\le
\sum_{j_q=0}^{p}~
\int\limits_{t}^{t_{l+3}}\left\vert \psi_{l+1}(\tau)\phi_{j_{q}}(\tau)\right\vert
\left\vert\sum_{j_l=p+1}^{\infty}\int\limits_{t}^{\tau} \psi_{l-1}(\theta)\phi_{j_{l}}(\theta)
\int\limits_{t}^{\theta}\psi_l(u)\phi_{j_{l}}(u)du d\theta\right\vert \times
$$

\vspace{2mm}
\begin{equation}
\label{func1002}
\times
\left\vert\int\limits_{t}^{\tau} \psi_{l+2}(u)\phi_{j_{q}}(u)du \right\vert d\tau.
\end{equation}

\vspace{5mm}

Combining (\ref{func801}) and (\ref{func1001}), we obtain

\begin{equation}
\label{func1002x}
|a|\le \frac{C}{p^{1-\varepsilon}},
\end{equation}

\vspace{3mm}
\noindent
where constant $C$ is independent of $p.$

Separating in (\ref{func1002}) the term with the number $j_q=0$ and then applying 
(\ref{ogo24}), (\ref{101xx}), (\ref{func801}), we obtain

$$
|b|\le
\frac{K}{p^{1-\varepsilon}}
\left(\int\limits_t^{t_{l+3}} \frac{d\tau}{\left(1-z^2(\tau)\right)^{1/2-\varepsilon/4}}+
\sum_{j_q=1}^{p}\frac{1}{j_q}
\int\limits_t^{t_{l+3}} \frac{d\tau}{\left(1-z^2(\tau)\right)^{3/4-\varepsilon/4}}\right)\le
$$

\vspace{4mm}
$$
\le \frac{K_1}{p^{1-\varepsilon}}\left(1+\sum_{j_q=1}^{p}\frac{1}{j_q}\right)\le
\frac{K_1}{p^{1-\varepsilon}}\left(2+\int\limits_1^p \frac{dx}{x}\right)=
$$

\vspace{4mm}
\begin{equation}
\label{func1003}
=\frac{K_1\left(2+ ln p\right)}{p^{1-\varepsilon}}\ \to\ 0
\end{equation}

\vspace{5mm}

\noindent
if $p\to\infty.$ The estimates (\ref{func1000}), (\ref{func1002x}), (\ref{func1003}) 
complete the proof of (\ref{func505}). 

Finally, consider the proof of (\ref{func505qq}).
Using the elementary inequality $|ab|\le (a^2+b^2)/2$ and Par\-se\-val's equality, we have 

\vspace{1mm}
$$
\left(\hat Q_p(t_1,\ldots,t_{l-1},t_{l+2},\ldots, t_{q-1}, t_{q+2}, \ldots, t_k)\right)^2\le
$$

\vspace{2mm}
$$
\le 
\left(\sum_{j_{l}=p+1}^{\infty}
\sum_{j_{l+1}=p+1}^{\infty}~\left\vert\int\limits_{t}^{t_{l+2}} \psi_{l+1}(\theta)\phi_{j_{l+1}}(\theta)
\int\limits_{t}^{\theta} \psi_l(u)\phi_{j_{l}}(u)du d\theta\right\vert\times\right.
$$

\vspace{2mm}
$$
\times
\left.\left\vert\int\limits_{t}^{t_{q+2}} \psi_{q+1}(\theta)\phi_{j_{l+1}}(\theta)
\int\limits_{t}^{\theta} \psi_{q}(u)\phi_{j_{l}}(u)du d\theta\right\vert\right)^2\le
$$

\vspace{2mm}

$$
\le 
\frac{1}{4}\left(\sum_{j_{l}=p+1}^{\infty}
\sum_{j_{l+1}=p+1}^{\infty}~\left(\int\limits_{t}^{t_{l+2}} \psi_{l+1}(\theta)\phi_{j_{l+1}}(\theta)
\int\limits_{t}^{\theta} \psi_l(u)\phi_{j_{l}}(u)du d\theta\right)^2+\right.
$$

\vspace{2mm}
$$
+
\left.
\sum_{j_{l}=p+1}^{\infty}
\sum_{j_{l+1}=p+1}^{\infty}~
\left(\int\limits_{t}^{t_{q+2}} \psi_{q+1}(\theta)\phi_{j_{l+1}}(\theta)
\int\limits_{t}^{\theta} \psi_{q}(u)\phi_{j_{l}}(u)du d\theta\right)^2\right)^2\le
$$

\vspace{2mm}
$$
\le 
\frac{1}{4}\left(\sum_{j_{l}=p+1}^{\infty}
\sum_{j_{l+1}=0}^{\infty}~\left(\int\limits_{t}^{t_{l+2}} \psi_{l+1}(\theta)\phi_{j_{l+1}}(\theta)
\int\limits_{t}^{\theta} \psi_l(u)\phi_{j_{l}}(u)du d\theta\right)^2+\right.
$$

\vspace{2mm}
$$
+
\left.
\sum_{j_{l}=p+1}^{\infty}
\sum_{j_{l+1}=0}^{\infty}~
\left(\int\limits_{t}^{t_{q+2}} \psi_{q+1}(\theta)\phi_{j_{l+1}}(\theta)
\int\limits_{t}^{\theta} \psi_{q}(u)\phi_{j_{l}}(u)du d\theta\right)^2\right)^2\le
$$

\vspace{2mm}
$$
\le 
\frac{1}{4}\left(\sum_{j_{l}=p+1}^{\infty}~\int\limits_{t}^{t_{l+2}} 
\psi_{l+1}^2(\theta)
\left(\int\limits_{t}^{\theta} \psi_l(u)\phi_{j_{l}}(u)du\right)^2 d\theta+\right.
$$

\vspace{2mm}
\begin{equation}
\label{func651}
+
\left.
\sum_{j_{l}=p+1}^{\infty}~\int\limits_{t}^{t_{q+2}} 
\psi_{q+1}^2(\theta)
\left(\int\limits_{t}^{\theta} \psi_{q}(u)\phi_{j_{l}}(u)du\right)^2 d\theta\right)^2.
\end{equation}

\vspace{5mm}

Note that
\begin{equation}
\label{obana}
\sum\limits_{j=p+1}^{\infty}\frac{1}{j^2}
\le \int\limits_{p}^{\infty}\frac{dx}{x^2}=\frac{1}{p}.
\end{equation}

\vspace{4.5mm}

From (\ref{func651}) and (\ref{obana}), (\ref{101xx}) we obtain

$$
\left(\hat Q_p(t_1,\ldots,t_{l-1},t_{l+2},\ldots, t_{q-1}, t_{q+2}, \ldots, t_k)\right)^2\le
$$
$$
\le
\frac{K}{p^2}\ \to\ 0
$$

\vspace{3mm}
\noindent
if $p\to\infty,$ where constant $K$ does not depend on $p.$
Thus the equalities 
(\ref{func503})--(\ref{func505qq}) are proved.

Recall that the function 
(\ref{de200}) (this function is defined using the left-hand side of the equality (\ref{drdr1000}))
for the case $k > 5,$ $r=2$
is represented
as the sum of several functions. Four of them, namely $Q_p,$ $\bar Q_p,$ $\tilde Q_p,$
$\hat Q_p$
(these functions correspond to the particular case of choosing the pairs $(g_1,g_2),$ $(g_3,g_4)$;
generally speaking, all possible pairs $(g_1,g_2),$ $(g_3,g_4)$ must be considered),
have been studied above. Absolutely similarly, we can consider
the remaining functions (for all possible pairs $(g_1,g_2),$ $(g_3,g_4)$)
whose sum is the function 
(\ref{de200}) 
for the case $k > 5,$ $r=2.$ As a result, we will have

$$
\lim\limits_{p\to\infty}
\bigl\Vert \hat R_p \bigr\Vert_{L_2([t, T]^{k-2r})}^2=0\ \ \ (k > 5,\  r=2).
$$

\vspace{4mm}

After that, we can go to the function 
(\ref{de200}) 
for the case $k > 5,$ $r=3,$ $2r<k$
(this function is defined using the left-hand side of the equality (\ref{drdr1000}))
and follow the same steps as above. This will lead us to the following 
equality

$$
\lim\limits_{p\to\infty}
\bigl\Vert \hat R_p \bigr\Vert_{L_2([t, T]^{k-2r})}^2=0\ \ \ (k > 5,\  r=3,\ 2r<k).
$$

\vspace{3mm}

Then we can move on to the next step and so on.
As a result, we get the equality (\ref{pars3s}) ($r=1,2,\ldots,[k/2]$). Thus 
the condition (\ref{drdr1001}) is satisfied for the case $k=2n+1,$ $n=3, 4, \ldots $
(recall that the condition (\ref{drdr1001}) is weaker than Condition~3 of Theorem~13
and the condition (\ref{drdr1001}) can be used in Theorem~13 instead of
Condition~3).

For the case $k=2n,$ $n=3, 4, \ldots$ we follow the above steps
for $r=1,2,\ldots,[k/2]-1$ $(2r\le k-2$).
For $2r=k$ we use the same technique as in the proof of the equalities
(\ref{after2508})--(\ref{after2507}). Recall that we used 
(\ref{after80xx}), (\ref{after500}) and 
Parseval's equality in the proof of (\ref{after2508})--(\ref{after2507}).

The obvious disadvantage of the proposed algorithm is the drastic 
increase of complexity of the proof when moving from $r=1$ to $r=2,$
$r=2$ to $r=3$ and so on.

The proofs of Theorems~17 and 18 contain a rather simple trick
of passing from $r=1$ to $r=2.$
Unfortunately, this procedure cannot be 
applied already at the transition from $r=2$ to $r=3.$

Note that the case $k=6,$ $r=3$ was successfully 
considered in Theorem~23 under the following 
simplifying assumption:
$\psi_1(\tau),\ldots,\psi_6(\tau)\equiv 1.$

Nevertheless, the results obtained in this paper
are quite sufficient for practical needs (see Chapters~4 and 5 \cite{2018a} for 
details).

\vspace{5mm}

\section{Generalization of Theorem~13 for 
Complete Or\-tho\-nor\-mal Sys\-tems of Functions in $L_2([t, T])$
and $\psi_1(\tau),$ $\ldots,$ $\psi_k(\tau)$ $\in $ $L_2([t, T])$
such that the Condition (\ref{drdr1001xyx}) is Satisfied}

\vspace{5mm}

First, note that (see the proof of Theorem~13 and (\ref{after906}))

\vspace{2mm}
$$
\hbox{\vtop{\offinterlineskip\halign{
\hfil#\hfil\cr
{\rm l.i.m.}\cr
$\stackrel{}{{}_{p\to \infty}}$\cr
}} }\sum_{j_1,\ldots,j_k=0}^{p}
C_{j_k\ldots j_1}
\prod\limits_{s=1}^r{\bf 1}_{\{j_{g_{{}_{2s-1}}}=~j_{g_{{}_{2s}}}\}}
\prod\limits_{s=1}^r{\bf 1}_{\{i_{g_{{}_{2s-1}}}=~i_{g_{{}_{2s}}}\ne 0\}}
J'[\phi_{j_{q_1}}\ldots \phi_{j_{q_{k-2r}}}]_{T,t}^{(i_{q_1}\ldots i_{q_{k-2r}})}=
$$

\vspace{5mm}
$$
=\hbox{\vtop{\offinterlineskip\halign{
\hfil#\hfil\cr
{\rm l.i.m.}\cr
$\stackrel{}{{}_{p\to \infty}}$\cr
}} }\hspace{-2.5mm}
\sum\limits_{\stackrel{j_1,\ldots,j_q,\ldots,j_k=0}{{}_{q\ne g_1, g_2,\ldots, g_{2r-1}, g_{2r}}}}^p
\sum\limits_{j_{g_1},j_{g_3},\ldots,j_{g_{2r-1}}=0}^p
\hspace{-2.5mm}
C_{j_k \ldots j_1}\biggl|_{j_{g_{{}_{1}}}=~j_{g_{{}_{2}}},\ldots, j_{g_{{}_{2r-1}}}=~j_{g_{{}_{2r}}}
}\biggr.
\times
$$

\vspace{5mm}
$$
\times
\prod\limits_{s=1}^r
{\bf 1}_{\{i_{g_{{}_{2s-1}}}=~i_{g_{{}_{2s}}}\ne 0\}}
J'[\phi_{j_{q_1}}\ldots \phi_{j_{q_{k-2r}}}]_{T,t}^{(i_{q_1}\ldots i_{q_{k-2r}})}=
$$

\vspace{5mm}
$$
=\hbox{\vtop{\offinterlineskip\halign{
\hfil#\hfil\cr
{\rm l.i.m.}\cr
$\stackrel{}{{}_{p\to \infty}}$\cr
}} }\hspace{-2.5mm}
\sum\limits_{\stackrel{j_1,\ldots,j_q,\ldots,j_k=0}{{}_{q\ne g_1, g_2,\ldots, g_{2r-1}, g_{2r}}}}^p
\Biggl(\sum\limits_{j_{g_1},j_{g_3},\ldots,j_{g_{2r-1}}=0}^p
\hspace{-2.5mm}
C_{j_k \ldots j_1}\biggl|_{j_{g_{{}_{1}}}=~j_{g_{{}_{2}}},\ldots, j_{g_{{}_{2r-1}}}=~j_{g_{{}_{2r}}}
}\biggr.-\Biggr.
$$

\vspace{5mm}
$$
\Biggl.-\frac{1}{2^r} \prod\limits_{l=1}^r {\bf 1}_{\{g_{2l}=g_{2l-1}+1\}}
C_{j_k \ldots j_1}\biggl|_{(j_{g_2} j_{g_1})\curvearrowright (\cdot)
\ldots (j_{g_{2r}} j_{g_{2r-1}})\curvearrowright (\cdot),
j_{g_{{}_{1}}}=~j_{g_{{}_{2}}},\ldots, j_{g_{{}_{2r-1}}}=~j_{g_{{}_{2r}}}}\Biggr)\times
$$

\vspace{5mm}
$$
\times
\prod\limits_{s=1}^r
{\bf 1}_{\{i_{g_{{}_{2s-1}}}=~i_{g_{{}_{2s}}}\ne 0\}}
J'[\phi_{j_{q_1}}\ldots \phi_{j_{q_{k-2r}}}]_{T,t}^{(i_{q_1}\ldots i_{q_{k-2r}})}+
$$

\vspace{5mm}
$$
+\hbox{\vtop{\offinterlineskip\halign{
\hfil#\hfil\cr
{\rm l.i.m.}\cr
$\stackrel{}{{}_{p\to \infty}}$\cr
}} }
\sum\limits_{\stackrel{j_1,\ldots,j_q,\ldots,j_k=0}{{}_{q\ne g_1, g_2,\ldots, g_{2r-1}, g_{2r}}}}^p
\frac{1}{2^r}
C_{j_k \ldots j_1}\biggl|_{(j_{g_2} j_{g_1})\curvearrowright (\cdot)
\ldots (j_{g_{2r}} j_{g_{2r-1}})\curvearrowright (\cdot),
j_{g_{{}_{1}}}=~j_{g_{{}_{2}}},\ldots, j_{g_{{}_{2r-1}}}=~j_{g_{{}_{2r}}}}\biggr.
\times 
$$

\vspace{5mm}
$$
\times
\prod\limits_{s=1}^r
{\bf 1}_{\{i_{g_{{}_{2s-1}}}=~i_{g_{{}_{2s}}}\ne 0\}}
\prod\limits_{s=1}^r
{\bf 1}_{\{g_{2s}=g_{2s-1}+1\}}
J'[\phi_{j_{q_1}}\ldots \phi_{j_{q_{k-2r}}}]_{T,t}^{(i_{q_1}\ldots i_{q_{k-2r}})}=
$$

\vspace{5mm}
$$
=\hbox{\vtop{\offinterlineskip\halign{
\hfil#\hfil\cr
{\rm l.i.m.}\cr
$\stackrel{}{{}_{p\to \infty}}$\cr
}} }\hspace{-2.5mm}
\sum\limits_{\stackrel{j_1,\ldots,j_q,\ldots,j_k=0}{{}_{q\ne g_1, g_2,\ldots, g_{2r-1}, g_{2r}}}}^p
\Biggl(\sum\limits_{j_{g_1},j_{g_3},\ldots,j_{g_{2r-1}}=0}^p
\hspace{-2.5mm}
C_{j_k \ldots j_1}\biggl|_{j_{g_{{}_{1}}}=~j_{g_{{}_{2}}},\ldots, j_{g_{{}_{2r-1}}}=~j_{g_{{}_{2r}}}
}\biggr.-\Biggr.
$$

\vspace{5mm}
$$
\Biggl.-\frac{1}{2^r} \prod\limits_{l=1}^r {\bf 1}_{\{g_{2l}=g_{2l-1}+1\}}
C_{j_k \ldots j_1}\biggl|_{(j_{g_2} j_{g_1})\curvearrowright (\cdot)
\ldots (j_{g_{2r}} j_{g_{2r-1}})\curvearrowright (\cdot),
j_{g_{{}_{1}}}=~j_{g_{{}_{2}}},\ldots, j_{g_{{}_{2r-1}}}=~j_{g_{{}_{2r}}}}\Biggr)\times
$$

\vspace{5mm}
$$
\times
\prod\limits_{s=1}^r
{\bf 1}_{\{i_{g_{{}_{2s-1}}}=~i_{g_{{}_{2s}}}\ne 0\}}
J'[\phi_{j_{q_1}}\ldots \phi_{j_{q_{k-2r}}}]_{T,t}^{(i_{q_1}\ldots i_{q_{k-2r}})}+
$$

\vspace{5mm}
\begin{equation}
\label{dydy11}
+
\frac{1}{2^r}\prod\limits_{s=1}^r
{\bf 1}_{\{g_{2s}=g_{2s-1}+1\}}
J[\psi^{(k)}]_{T,t}^{s_r, \ldots, s_1}\ \ \ \hbox{w.~p.~1.}
\end{equation}

\vspace{5mm}

Using (\ref{dydy11}) and the condition (\ref{drdr1001}), we 
obtain (\ref{after801}). This means that we get (\ref{afteru11}).
Thus the expansion (\ref{after1}) is proved.

Analyzing the proof of Theorems~13 and conditions of Theorem~5
as well as taking into account the above arguments,
it is easy to see that the following theorem is true.

\vspace{2mm}

{\bf Theorem~39}\ \cite{2018a}, \cite{arxiv-9}.\ {\it Assume that
the continuous functions 
$\psi_1(\tau),\ldots,\psi_k(\tau)$ at the interval $[t, T]$ and 
the complete orthonormal system $\{\phi_j(x)\}_{j=0}^{\infty}$
of func\-ti\-ons $(\phi_0(x)=1/\sqrt{T-t})$ 
in the space $L_2([t, T])$ are such that the following 
condition 

\vspace{1mm}
$$
\lim\limits_{p_1,\ldots,p_k\to\infty}~
\sum\limits_{j_1=0}^{p_1}\ldots \sum\limits_{j_q=0}^{p_q}\ldots \sum\limits_{j_k=0}^{p_k}~
\biggl|_{q\ne g_1, g_2, \ldots, g_{2r-1},g_{2r}}\times
$$

\vspace{4mm}
$$
\times
\Biggl(~\sum\limits_{j_{g_1}=0}^{\min\{p_{g_1}, p_{g_2}\}} \sum\limits_{j_{g_3}=0}^{\min\{p_{g_3}, p_{g_4}\}} \ldots \Biggr.
\sum\limits_{j_{g_{2r-1}}=0}^{\min\{p_{g_{2r-1}}, p_{g_{2r}}\}}
C_{j_k\ldots j_1}\biggl|_{j_{g_1}=j_{g_2},\ldots, j_{g_{2r-1}}=j_{g_{2r}}}-
$$

\vspace{2mm}
\begin{equation}
\label{drdr1001xyx}
\Biggl.-\frac{1}{2^r} \prod\limits_{l=1}^r {\bf 1}_{\{g_{2l}=g_{2l-1}+1\}}
C_{j_k \ldots j_1}\biggl|_{(j_{g_2} j_{g_1})\curvearrowright (\cdot)
\ldots (j_{g_{2r}} j_{g_{2r-1}})\curvearrowright (\cdot),
j_{g_{{}_{1}}}=~j_{g_{{}_{2}}},\ldots, j_{g_{{}_{2r-1}}}=~j_{g_{{}_{2r}}}
}\biggr.\Biggr)^2=0
\end{equation}

\vspace{5mm}
\noindent
is satisfied for all $r=1, 2,\ldots,[k/2]$.
Then, for the iterated Stratonovich stochastic integral 
of arbitrary multiplicity $k$

\vspace{1mm}
$$
J^{*}[\psi^{(k)}]_{T,t}^{(i_1\ldots i_k)}=
{\int\limits_t^{*}}^T
\psi_k(t_k) \ldots 
{\int\limits_t^{*}}^{t_{2}}
\psi_1(t_1) d{\bf w}_{t_1}^{(i_1)}\ldots
d{\bf w}_{t_k}^{(i_k)}
$$

\vspace{4mm}
\noindent
the following 
expansion

$$
J^{*}[\psi^{(k)}]_{T,t}^{(i_1\ldots i_k)}=
\hbox{\vtop{\offinterlineskip\halign{
\hfil#\hfil\cr
{\rm l.i.m.}\cr
$\stackrel{}{{}_{p_1,\ldots,p_k\to \infty}}$\cr
}} }
\sum\limits_{j_1=0}^{p_1}\ldots\sum\limits_{j_k=0}^{p_k}
C_{j_k \ldots j_1}\prod\limits_{l=1}^k \zeta_{j_l}^{(i_l)}
$$

\vspace{4mm}
\noindent
that converges in the mean-square sense is valid, where 

\vspace{0.5mm}
$$
C_{j_k \ldots j_1}=\int\limits_t^T\psi_k(t_k)\phi_{j_k}(t_k)\ldots
\int\limits_t^{t_2}
\psi_1(t_1)\phi_{j_1}(t_1)
dt_1\ldots dt_k
$$

\vspace{4mm}
\noindent
is the Fourier coefficient, 
${\rm l.i.m.}$ is a limit in the mean-square sense,
$i_1, \ldots, i_k=0, 1,\ldots,m,$

$$
\zeta_{j}^{(i)}=
\int\limits_t^T \phi_{j}(\tau) d{\bf w}_{\tau}^{(i)}
$$ 

\vspace{3mm}
\noindent
are independent standard Gaussian random variables for various 
$i$ or $j$ {\rm (}in the case when $i\ne 0${\rm )},
${\bf w}_{\tau}^{(i)}={\bf f}_{\tau}^{(i)}$ 
for $i=1,\ldots,m$ and 
${\bf w}_{\tau}^{(0)}=\tau.$}

\vspace{2mm}

Further in this section, we generalize Theorems~13, 39 to the case of
complete ortho\-nor\-mal systems of functions in the space $L_2([t, T])$
and $\psi_1(\tau),$ $\ldots,$ $\psi_k(\tau)$ $\in $ $L_2([t, T])$
such that the condition (\ref{drdr1001xyx}) is satisfied.

Let $(\Omega,{\rm F},{\sf P})$ be a complete probability
space and let $f(t,\omega)\stackrel{\sf def}{=}f_t:$ 
$[0, T]\times \Omega\rightarrow \mathbb{R}$
be the standard Wiener process
defined on the probability space $(\Omega,{\rm F},{\sf P}).$

Let us consider the family of $\sigma$-algebras
$\left\{{\rm F}_t,\ t\in[0,T]\right\}$ defined
on the probability space $(\Omega,{\rm F},{\sf P})$ and
connected
with the Wiener process $f_t$ in such a way that

\vspace{2mm}

1.\ ${\rm F}_s\subset {\rm F}_t\subset {\rm F}$\ for
$s<t.$

\vspace{2mm}

2.\ The Wiener process $f_t$ is ${\rm F}_t$-measurable for all
$t\in[0,T].$

\vspace{2mm}

3.\ The process $f_{t+\Delta}-f_{t}$ for all
$t\ge 0,$ $\Delta>0$ is independent with
the events of $\sigma$-algebra
${\rm F}_{t}.$

\vspace{3mm}

Let $\xi(\tau,\omega)\stackrel{\sf def}{=}
\xi_{\tau}:$ $[0, T]\times\Omega \to \mathbb{R}$ 
be some random process, which is measurable
with respect to the pair of variables
$(\tau,\omega)$ and satisfies to the following
condition

$$
\int\limits_t^T |\xi_{\tau}| d\tau <\infty\ \ \ \hbox{w.~p.~1}\ \ \ (t\ge 0).
$$

\vspace{3mm}

Let $\tau_j^{(N)},$ $j=0, 1, \ldots, N$ 
be a partition of the interval $[t, T],$ $t\ge 0$ such that

\begin{equation}
\label{dsds4}
t=\tau_0^{(N)}<\tau_1^{(N)}<\ldots <\tau_N^{(N)}=T,\ \ \ \
\max\limits_{0\le j\le N-1}\left|\tau_{j+1}^{(N)}-\tau_j^{(N)}\right|\to 0\ \
\hbox{if}\ \ N\to \infty.
\end{equation}

\vspace{3mm}
\noindent
Further, for simplicity, we write $\tau_j$ instead of 
$\tau_j^{(N)}.$

Consider the definition of the Stratonovich stochastic integral, 
which differs from the definition given in \cite{KlPl2}
(recall that we use definition \cite{KlPl2} above in this article).

The mean-square limit (if it exists)

\vspace{-2mm}
\begin{equation}
\label{dsds5}
\hbox{\vtop{\offinterlineskip\halign{
\hfil#\hfil\cr
{\rm l.i.m.}\cr
$\stackrel{}{{}_{N\to \infty}}$\cr
}} }\sum_{j=0}^{N-1}
\frac{1}{\tau_{j+1}-\tau_j}
\int\limits_{\tau_j}^{\tau_{j+1}}\xi_s ds\
\left(f_{\tau_{j+1}}-
f_{\tau_j}\right)
\stackrel{\sf def}{=}\int\limits_t^T \xi_{\tau} \circ df_\tau
\end{equation}

\vspace{1mm}
\noindent
is called \cite{bardina10} the Stratonovich stochastic integral 
of the process $\xi_{\tau}$, $\tau\in [t, T]$,
where $\tau_j,$ $j=0, 1, \ldots, N$
is a partition of the interval $[t, T]$ 
satisfying the condition (\ref{dsds4}).

We also denote by
$$
\int\limits_t^{\tau} \xi_{s} \circ df_s
$$

\vspace{0.5mm}
\noindent
the Stratonovich stochastic integral like (\ref{dsds5}) (if it exists) 
of $\xi_s {\bf 1}_{\{s\in [t, \tau]\}}$ for $\tau\in [t, T],$ $t\ge 0.$

It is known \cite{bardina10} (Lemma~A.2) that the following 
iterated Stratonovich stochastic integral 

\vspace{-2mm}
\begin{equation}
\label{dsds7}
J^{S}[\psi^{(k)}]_{\tau,t}^{(i_1\ldots i_k)}=
\int\limits_t^{\tau}
\psi_k(t_k)\ldots \int\limits_t^{t_2}
\psi_1(t_1) \circ d{\bf w}_{t_1}^{(i_1)}\ldots \circ d{\bf w}_{t_k}^{(i_k)}
\end{equation}

\vspace{2mm}
\noindent
exists for the case $i_1=\ldots=i_k\ne 0$, 
where $\tau\in [t, T],$ $\psi_1(\tau),\ldots,\psi_k(\tau)\in L_2([t, T]),$
$i_1,\ldots,i_k=0,1,\ldots,m,$\ 
${\bf w}_{\tau}^{(i)}=
{\bf f}_{\tau}^{(i)}$ for $i=1,\ldots,m$ and ${\bf w}_{\tau}^{(0)}=\tau$,
${\bf f}_{\tau}^{(i)}$ $(i=1,\ldots,m)$ are independent 
standard Wiener processes defined as above in this section.

In \cite{new-new-18} (2021) an analogue of Theorem~5 (1997) 
is proved for the case $i_1=\ldots=i_k\ne 0$ and 
$\psi_1(\tau),\ldots,\psi_k(\tau)\in L_2([t, T]).$

Let us denote

\vspace{-2mm}
\begin{equation}
\label{dsds9}
J[\psi^{(k)}]_{T,t}^{(i_1\ldots i_k)}+
\sum_{r=1}^{\left[k/2\right]}\frac{1}{2^r}
\sum_{(s_r,\ldots,s_1)\in {\rm A}_{k,r}}
J[\psi^{(k)}]_{T,t}^{s_r,\ldots,s_1}\stackrel{\sf def}{=}\bar J^{*}[\psi^{(k)}]_{T,t}^{(i_1\ldots i_k)},
\end{equation}

\vspace{2mm}
\noindent
where $\psi_1(\tau),\ldots,\psi_k(\tau)\in L_2([t, T]),$
$\psi_l(\tau)\psi_{l-1}(\tau)\in L_2([t, T])$ $(l=2, 3,\ldots, k),$
$J[\psi^{(k)}]_{T,t}^{(i_1\ldots i_k)}$ is the iterated Ito stochastic
integral (\ref{dsds12}),
$\sum\limits_{\emptyset}$ is supposed to be equal to zero; 
another notations are the same as in Theorem~{\rm 5.}

Further, by analogy with (\ref{after8}), (\ref{after7})
and using the version of (\ref{2023abc300}) for the case of an arbitrary 
complete orthonormal system $\{\phi_j(x)\}_{j=0}^{\infty}$
in $L_2([t, T])$ (see \cite{2018a} or \cite{12aa-afterxxx}, Sect.~1.11)
instead of (\ref{2023abc300}), we obtain 
the following generalization of (\ref{after8}) to the case 
of an arbitrary 
complete ortho\-nor\-mal system $\{\phi_j(x)\}_{j=0}^{\infty}$ in $L_2([t, T])$
and $\psi_1(\tau),$ $\ldots,$ $\psi_k(\tau)$ $\in $ $L_2([t, T])$

\vspace{1mm}
$$
\sum_{j_1=0}^{p_1}\ldots\sum_{j_k=0}^{p_k}
C_{j_k\ldots j_1}
\prod_{l=1}^k \zeta_{j_l}^{(i_l)}=
\sum_{j_1=0}^{p_1}\ldots\sum_{j_k=0}^{p_k}
C_{j_k\ldots j_1}
J'[\phi_{j_1}\ldots \phi_{j_k}]_{T,t}^{(i_1\ldots i_k)}+
$$

\vspace{2mm}
$$
+\sum_{j_1=0}^{p_1}\ldots\sum_{j_k=0}^{p_k}
C_{j_k\ldots j_1}
\sum\limits_{r=1}^{[k/2]}
\sum_{\stackrel{(\{\{g_1, g_2\}, \ldots, 
\{g_{2r-1}, g_{2r}\}\}, \{q_1, \ldots, q_{k-2r}\})}
{{}_{\{g_1, g_2, \ldots, 
g_{2r-1}, g_{2r}, q_1, \ldots, q_{k-2r}\}=\{1, 2, \ldots, k\}}}}
\prod\limits_{s=1}^r
{\bf 1}_{\{i_{g_{{}_{2s-1}}}=~i_{g_{{}_{2s}}}\ne 0\}}\times
$$

\vspace{2mm}
\begin{equation}
\label{after8xxds1}
\times{\bf 1}_{\{j_{g_{{}_{2s-1}}}=~j_{g_{{}_{2s}}}\}}
J'[\phi_{j_{q_1}}\ldots \phi_{j_{q_{k-2r}}}]_{T,t}^{(i_{q_1}\ldots i_{q_{k-2r}})}\ \ \ \hbox{w.~p.~1,}
\end{equation}

\vspace{4mm}
\noindent
where $J'[\phi_{j_1}\ldots \phi_{j_k}]_{T,t}^{(i_1\ldots i_k)},$
$J'[\phi_{j_{q_1}}\ldots \phi_{j_{q_{k-2r}}}]_{T,t}^{(i_{q_1}\ldots i_{q_{k-2r}})}$
are multiple Wiener sto\-chas\-tic integrals
de\-fi\-ned as in \cite{ito1951} (1951). Note that in
\cite{ito1951} the case of a scalar Wiener process has been considered.

It should be noted that 
Theorem~1.16 \cite{2018a} (Sect.~1.11) and Theorem~4 can be reformulated as follows
(also see \cite{arxiv-1}, Sect.~15)

\vspace{-1mm}
\begin{equation}
\label{dsds11}
J[\psi^{(k)}]_{T,t}^{(i_1\ldots i_k)}=
\hbox{\vtop{\offinterlineskip\halign{
\hfil#\hfil\cr
{\rm l.i.m.}\cr
$\stackrel{}{{}_{p_1,\ldots,p_k\to \infty}}$\cr
}} }\sum_{j_1=0}^{p_1}\ldots\sum_{j_k=0}^{p_k}
C_{j_k\ldots j_1}J'[\phi_{j_1}\ldots \phi_{j_k}]_{T,t}^{(i_1\ldots i_k)}\ \ \ \hbox{w.~p.~1},
\end{equation}

\vspace{3mm}
\noindent
where $\{\phi_j(x)\}_{j=0}^{\infty}$ is
an arbitrary 
complete ortho\-nor\-mal system in $L_2([t, T]),$
$\psi_1(\tau),$ $\ldots,$ $\psi_k(\tau)$ $\in $ $L_2([t, T]),$
$J'[\phi_{j_1}\ldots \phi_{j_k}]_{T,t}^{(i_1\ldots i_k)}$ is the 
multiple Wiener stochastic integral
defined as in \cite{ito1951} (1951) and
$J[\psi^{(k)}]_{T,t}^{(i_1\ldots i_k)}$ is the iterated Ito stochastic
integral

\begin{equation}
\label{dsds12}
J[\psi^{(k)}]_{T,t}^{(i_1\ldots i_k)}=
\int\limits_t^{T}\psi_k(t_k) \ldots 
\int\limits_t^{t_{2}}
\psi_1(t_1) d{\bf w}_{t_1}^{(i_1)}\ldots
d{\bf w}_{t_k}^{(i_k)};
\end{equation}

\vspace{3mm}
\noindent
another notations are the same as in Theorem~4.

Passing to the limit 
$\hbox{\vtop{\offinterlineskip\halign{
\hfil#\hfil\cr
{\rm l.i.m.}\cr
$\stackrel{}{{}_{p_1,\ldots,p_k\to \infty}}$\cr
}} }$ in (\ref{after8xxds1}) and using the equality (\ref{dsds11}), we get w.~p.~1 

\vspace{1mm}
$$
\hbox{\vtop{\offinterlineskip\halign{
\hfil#\hfil\cr
{\rm l.i.m.}\cr
$\stackrel{}{{}_{p_1,\ldots,p_k\to \infty}}$\cr
}} }\sum_{j_1=0}^{p_1}\ldots\sum_{j_k=0}^{p_k}
C_{j_k\ldots j_1}
\zeta_{j_1}^{(i_1)}\ldots \zeta_{j_k}^{(i_k)}
=J[\psi^{(k)}]_{T,t}^{(i_1\ldots i_k)}
+
$$

\vspace{2mm}
$$
+
\sum\limits_{r=1}^{[k/2]}
\sum_{\stackrel{(\{\{g_1, g_2\}, \ldots, 
\{g_{2r-1}, g_{2r}\}\}, \{q_1, \ldots, q_{k-2r}\})}
{{}_{\{g_1, g_2, \ldots, 
g_{2r-1}, g_{2r}, q_1, \ldots, q_{k-2r}\}=\{1, 2, \ldots, k\}}}}
\prod\limits_{s=1}^r
{\bf 1}_{\{i_{g_{{}_{2s-1}}}=~i_{g_{{}_{2s}}}\ne 0\}}\times
$$

\vspace{4mm}
\begin{equation}
\label{after501ds1}
\times \hbox{\vtop{\offinterlineskip\halign{
\hfil#\hfil\cr
{\rm l.i.m.}\cr
$\stackrel{}{{}_{p_1,\ldots,p_k\to \infty}}$\cr
}} }\sum_{j_1=0}^{p_1}\ldots\sum_{j_k=0}^{p_k}
C_{j_k\ldots j_1}
\prod\limits_{s=1}^r{\bf 1}_{\{j_{g_{{}_{2s-1}}}=~j_{g_{{}_{2s}}}\}}
J'[\phi_{j_{q_1}}\ldots \phi_{j_{q_{k-2r}}}]_{T,t}^{(i_{q_1}\ldots i_{q_{k-2r}})},
\end{equation}

\vspace{5mm}
\noindent
where 
$J'[\phi_{j_{q_1}}\ldots \phi_{j_{q_{k-2r}}}]_{T,t}^{(i_{q_1}\ldots i_{q_{k-2r}})}$ is the 
multiple Wiener stochastic integral
defined as in \cite{ito1951} (1951) and
$J[\psi^{(k)}]_{T,t}^{(i_1\ldots i_k)}$ is the iterated Ito stochastic
integral (\ref{dsds12}).

Suppose that $\{\phi_j(x)\}_{j=0}^{\infty}$ is an arbitrary
complete orthonormal system of functions in $L_2([t, T])$
and $\Phi_1(\tau), \Phi_2(\tau)\in L_2([t, T])$.
Then we have

$$
\sum_{j=0}^{\infty}\left|\int\limits_t^s \phi_j(\tau)\Phi_1(\tau)d\tau 
\int\limits_s^T \phi_j(\tau)\Phi_2(\tau)d\tau\right|\le 
$$

\vspace{2mm}
$$
\le \frac{1}{2}\sum_{j=0}^{\infty}
\left(\left(\int\limits_t^T {\bf 1}_{\{\tau<s\}}\phi_j(\tau)\Phi_1(\tau)d\tau\right)^2+ 
\left(\int\limits_t^T {\bf 1}_{\{\tau>s\}}\phi_j(\tau)\Phi_2(\tau)d\tau\right)^2\right)=
$$

\vspace{2mm}
\begin{equation}
\label{dsds14}
=\frac{1}{2}\left(\int\limits_t^s \Phi_1^2(\tau)d\tau+
\int\limits_s^T\Phi_2^2(\tau)d\tau\right)\le
\frac{1}{2}\left(\left\Vert\Phi_1\right\Vert_{L_2([t,T])}^2+
\left\Vert\Phi_2\right\Vert_{L_2([t,T])}^2\right)<\infty,
\end{equation}

\vspace{4mm}
\noindent
i.e. 
\begin{equation}
\label{dsds14fffff}
\left|\sum_{j=0}^{p}\int\limits_t^s \phi_j(\tau)\Phi_1(\tau)d\tau 
\int\limits_s^T \phi_j(\tau)\Phi_2(\tau)d\tau\right|\le C<\infty,
\end{equation}

\vspace{3mm}
\noindent
where $p\in\mathbb{N}.$

By interpreting the integrals in (\ref{after9})--(\ref{after400}) as 
Lebesgue integrals, using Fubini's theorem in (\ref{after9}) and
Lebesgue's 
Dominated Convergence Theorem in (\ref{after11}), we 
obtain (\ref{after80}) (see (\ref{dsds14fffff})) for
the case of an arbitrary complete 
orthonormal system of functions in the space $L_2([t, T])$
and $\psi_1(\tau),\ldots, \psi_k(\tau)\in L_2([t, T])$.

Using the equality (\ref{after1400}) 
for the case of an arbitrary complete 
orthonormal system of functions in the space $L_2([t, T])$
and $\Phi_1(\tau),\Phi_2(\tau)\in L_2([t, T])$ as well as
Fubini's Theorem when deriving (\ref{r12345x}), we obtain the generalization of
(\ref{after500}) for the case of an arbitrary complete 
orthonormal system of functions in the space $L_2([t, T])$
and $\psi_1(\tau),\ldots, \psi_k(\tau)\in L_2([t, T])$.

Repeating the steps of the proof of Theorem~13 below
the formula (\ref{after501}) using (\ref{dsds9}), (\ref{after501ds1})
or steps 
of the proof of Theorem~39 using (\ref{dsds9}), (\ref{after501ds1}),
we obtain
for complete 
orthonormal systems $\{\phi_j(x)\}_{j=0}^{\infty}$
$(\phi_0(x)=1/\sqrt{T-t})$ 
in the space $L_2([t, T])$
and $\psi_1(\tau),\ldots, \psi_k(\tau)\in L_2([t, T]),$
$\psi_l(\tau)\psi_{l-1}(\tau)\in L_2([t, T])$ $(l=2, 3,\ldots, k)$
(for which the condition (\ref{drdr1001xyx}) is satisfied) the following equality

$$
\hbox{\vtop{\offinterlineskip\halign{
\hfil#\hfil\cr
{\rm l.i.m.}\cr
$\stackrel{}{{}_{p_1,\ldots,p_k\to \infty}}$\cr
}} }
\sum_{j_1=0}^{p_1}\ldots\sum_{j_k=0}^{p_k}
C_{j_k \ldots j_1}\prod\limits_{l=1}^k \zeta_{j_l}^{(i_l)}
=
$$

\vspace{2mm}
\begin{equation}
\label{after333ds1}
=
J[\psi^{(k)}]_{T,t}^{(i_1\ldots i_k)}+
\sum_{r=1}^{\left[k/2\right]}\frac{1}{2^r}
\sum_{(s_r,\ldots,s_1)\in {\rm A}_{k,r}}
J[\psi^{(k)}]_{T,t}^{s_r,\ldots,s_1}=
\bar J^{*}[\psi^{(k)}]_{T,t}^{(i_1\ldots i_k)}
\end{equation}

\vspace{3mm}
\noindent
w.~p.~1, where notations in (\ref{after333ds1}) are the same as in Theorem~5
and $\bar J^{*}[\psi^{(k)}]_{T,t}^{(i_1\ldots i_k)}$ is defined by
(\ref{dsds9}).

Thus the following two theorems are proved.

\vspace{2mm}

{\bf Theorem~40}\ \cite{2018a}, \cite{arxiv-4}, \cite{arxiv-9}.\ {\it Assume that
the complete orthonormal system $\{\phi_j(x)\}_{j=0}^{\infty}$
$(\phi_0(x)=1/\sqrt{T-t})$ 
in the space $L_2([t, T])$ and 
$\psi_1(\tau),\ldots, \psi_k(\tau)\in L_2([t, T]),$
$\psi_l(\tau)\psi_{l-1}(\tau)\in L_2([t, T])$ $(l=2, 3,\ldots, k)$
are such that the folowing condition

\vspace{1mm}
$$
\lim\limits_{p_1,\ldots,p_k\to\infty}~
\sum\limits_{j_1=0}^{p_1}\ldots \sum\limits_{j_q=0}^{p_q}\ldots \sum\limits_{j_k=0}^{p_k}~
\biggl|_{q\ne g_1, g_2, \ldots, g_{2r-1},g_{2r}}\times
$$

\vspace{4mm}
$$
\times
\Biggl(~\sum\limits_{j_{g_1}=0}^{\min\{p_{g_1}, p_{g_2}\}} \sum\limits_{j_{g_3}=0}^{\min\{p_{g_3}, p_{g_4}\}}\ldots \Biggr.
\sum\limits_{j_{g_{2r-1}}=0}^{\min\{p_{g_{2r-1}}, p_{g_{2r}}\}}
C_{j_k\ldots j_1}\biggl|_{j_{g_1}=j_{g_2},\ldots, j_{g_{2r-1}}=j_{g_{2r}}}-
$$

\vspace{2mm}
\begin{equation}
\label{novorigin1}
\Biggl.-\frac{1}{2^r} \prod\limits_{l=1}^r {\bf 1}_{\{g_{2l}=g_{2l-1}+1\}}
C_{j_k \ldots j_1}\biggl|_{(j_{g_2} j_{g_1})\curvearrowright (\cdot)
\ldots (j_{g_{2r}} j_{g_{2r-1}})\curvearrowright (\cdot),
j_{g_{{}_{1}}}=~j_{g_{{}_{2}}},\ldots, j_{g_{{}_{2r-1}}}=~j_{g_{{}_{2r}}}
}\biggr.\Biggr)^2=0
\end{equation}

\vspace{5mm}
\noindent
is satisfied for all $r=1, 2,\ldots,[k/2]$.
Then, for the sum $\bar J^{*}[\psi^{(k)}]_{T,t}^{(i_1\ldots i_k)}$ of iterated Ito stochastic integrals 
defined by {\rm (\ref{dsds9})}
the following 
expansion 

\vspace{-1mm}
\begin{equation}
\label{january19c}
\bar J^{*}[\psi^{(k)}]_{T,t}^{(i_1\ldots i_k)}=
\hbox{\vtop{\offinterlineskip\halign{
\hfil#\hfil\cr
{\rm l.i.m.}\cr
$\stackrel{}{{}_{p_1,\ldots,p_k\to \infty}}$\cr
}} }
\sum_{j_1=0}^{p_1}\ldots\sum_{j_k=0}^{p_k}
C_{j_k \ldots j_1}\prod\limits_{l=1}^k \zeta_{j_l}^{(i_l)}
\end{equation}

\vspace{3mm}
\noindent
that converges in the mean-square sense is valid, where 

$$
C_{j_k \ldots j_1}=\int\limits_t^T\psi_k(t_k)\phi_{j_k}(t_k)\ldots
\int\limits_t^{t_2}
\psi_1(t_1)\phi_{j_1}(t_1)
dt_1\ldots dt_k
$$

\vspace{3mm}
\noindent
is the Fourier coefficient, 
${\rm l.i.m.}$ is a limit in the mean-square sense,
$i_1, \ldots, i_k=0, 1,\ldots,m,$

\vspace{-1mm}
$$
\zeta_{j}^{(i)}=
\int\limits_t^T \phi_{j}(\tau) d{\bf w}_{\tau}^{(i)}
$$ 

\vspace{2mm}
\noindent
are independent standard Gaussian random variables for various 
$i$ or $j$ {\rm (}in the case when $i\ne 0${\rm )},
${\bf w}_{\tau}^{(i)}={\bf f}_{\tau}^{(i)}$ 
for $i=1,\ldots,m$ and 
${\bf w}_{\tau}^{(0)}=\tau.$}

\vspace{2mm}

{\bf Theorem~41}\ \cite{2018a}, \cite{arxiv-9}.\ {\it Assume that
the complete orthonormal system $\{\phi_j(x)\}_{j=0}^{\infty}$
$(\phi_0(x)=1/\sqrt{T-t})$ 
in the space $L_2([t, T])$ and 
$\psi_1(\tau),\ldots, \psi_k(\tau)\in L_2([t, T]),$
$\psi_l(\tau)\psi_{l-1}(\tau)\in L_2([t, T])$ $(l=2, 3,\ldots, k)$
are such that 
the condition

\vspace{1mm}
$$
\lim\limits_{p\to\infty}
\sum\limits_{\stackrel{j_1,\ldots,j_q,\ldots,j_k=0}{{}_{q\ne g_1, g_2, \ldots, g_{2r-1},
g_{2r}}}}^p
\left(S_{l_1}S_{l_2}\ldots S_{l_{d}}
\left\{\bar C^{(p)}_{j_k\ldots j_q \ldots j_1}\biggl|_{q\ne g_1,g_2,\ldots,g_{2r-1}, g_{2r}}
\right\}\right)^2=0
$$

\vspace{3mm}
\noindent
holds for all possible $g_1,g_2,\ldots,g_{2r-1},g_{2r}$ {\rm (}see {\rm (\ref{leto5007xxx}))}
and $l_1, l_2, \ldots, l_{d}$ such that
$l_1, l_2, \ldots, l_{d}\in \{1,2,\ldots,$ $r\},$\
$l_1>l_2>\ldots >l_{d},$\ $d=0, 1, 2,\ldots, r-1,$\ 
where $r=1, 2,\ldots,[k/2]$ and

\vspace{1mm}
$$
S_{l_1}S_{l_2}\ldots S_{l_{d}}
\left\{\bar C^{(p)}_{j_k\ldots j_q \ldots j_1}\biggl|_{q\ne g_1,g_2,\ldots,g_{2r-1}, g_{2r}}
\right\}\stackrel{\sf def}{=}
\bar C^{(p)}_{j_k\ldots j_q \ldots j_1}\biggl|_{q\ne g_1,g_2,\ldots,g_{2r-1}, g_{2r}}
$$

\vspace{4mm}
\noindent
for $d=0.$

Then, for the sum $\bar J^{*}[\psi^{(k)}]_{T,t}^{(i_1\ldots i_k)}$ of iterated Ito stochastic integrals 
defined by {\rm (\ref{dsds9})}
the following 
expansion 
$$
\bar J^{*}[\psi^{(k)}]_{T,t}^{(i_1\ldots i_k)}=
\hbox{\vtop{\offinterlineskip\halign{
\hfil#\hfil\cr
{\rm l.i.m.}\cr
$\stackrel{}{{}_{p\to \infty}}$\cr
}} }
\sum_{j_1,\ldots,j_k=0}^{p}
C_{j_k \ldots j_1}\prod\limits_{l=1}^k \zeta_{j_l}^{(i_l)}
$$

\vspace{4mm}
\noindent
that converges in the mean-square sense is valid, where 

$$
C_{j_k \ldots j_1}=\int\limits_t^T\psi_k(t_k)\phi_{j_k}(t_k)\ldots
\int\limits_t^{t_2}
\psi_1(t_1)\phi_{j_1}(t_1)
dt_1\ldots dt_k
$$

\vspace{2mm}
\noindent
is the Fourier coefficient, 
${\rm l.i.m.}$ is a limit in the mean-square sense,
$i_1, \ldots, i_k=0, 1,\ldots,m,$

$$
\zeta_{j}^{(i)}=
\int\limits_t^T \phi_{j}(\tau) d{\bf w}_{\tau}^{(i)}
$$ 

\vspace{2mm}
\noindent
are independent standard Gaussian random variables for various 
$i$ or $j$ {\rm (}in the case when $i\ne 0${\rm )},
${\bf w}_{\tau}^{(i)}={\bf f}_{\tau}^{(i)}$ 
for $i=1,\ldots,m$ and 
${\bf w}_{\tau}^{(0)}=\tau.$}

\vspace{2mm}

Note that in Theorems~40, 41 (the case $k=2$)
the condition 
$\psi_1(\tau)\psi_{2}(\tau)\in L_2([t, T])$ 
can be omitted.

Using Theorem 5 together with Proposition 3.1 \cite{new-new-18} and the proof of
Lemma A.2 \cite{bardina10}, we can write 
$\bar J^{*}[\psi^{(k)}]_{T,t}^{(i_1\ldots i_k)}=J^{S}[\psi^{(k)}]_{T,t}^{(i_1\ldots i_k)}$\
w.~p.~1 and reformulate Theorems~40, 41 for 
$J^{S}[\psi^{(k)}]_{T,t}^{(i_1\ldots i_k)}$
($J^{S}[\psi^{(k)}]_{T,t}^{(i_1\ldots i_k)}$
is defined by (\ref{dsds7})).

Let us consider the special case $k=2$ of Theorem~40 in more detail.
In this case, the condition (\ref{novorigin1}) takes the following form
(compare with (\ref{5t}))

\vspace{-1mm}
\begin{equation}
\label{novorigin20}
\sum_{j_1=0}^{\infty}
C_{j_1j_1}=\frac{1}{2}
\int\limits_t^T\psi_1(t_1)\psi_2(t_1)dt_1.
\end{equation}

\vspace{3mm}

As follows from \cite{2018a} (Sect.~2.1.4), the equality (\ref{novorigin20})
is valid for the case 
of an arbitrary complete orthonormal 
system of functions in $L_2([t, T])$ and
$\psi_1(\tau), \psi_2(\tau)\in L_2([t, T]).$

From Proposition 3.1 \cite{new-new-18} for the case $k=2$ we obtain

\vspace{-1mm}
$$
\int\limits_t^{T}
\psi_2(t_2)\int\limits_t^{t_2}
\psi_1(t_1) \circ d{\bf w}_{t_1}^{(i)}\circ d{\bf w}_{t_2}^{(i)}
=\int\limits_t^{T}
\psi_2(t_2)\int\limits_t^{t_2}
\psi_1(t_1) d{\bf w}_{t_1}^{(i)} d{\bf w}_{t_2}^{(i)}+
$$

\begin{equation}
\label{novorigin3}
+
\frac{1}{2}
\int\limits_t^T\psi_1(t_1)\psi_2(t_1)dt_1
\end{equation}

\vspace{3mm}
\noindent
w.~p.~1, where $\psi_1(\tau), \psi_2(\tau)\in L_2([t, T]),$
$i=1,\ldots,m,$

\vspace{1mm}
$$
\int\limits_t^{T}
\psi_2(t_2)\int\limits_t^{t_2}
\psi_1(t_1) \circ d{\bf w}_{t_1}^{(i)}\circ d{\bf w}_{t_2}^{(i)}
$$

\vspace{3mm}
\noindent
is defined by (\ref{dsds5}), (\ref{dsds7}) and 

\vspace{-2mm}
$$
\int\limits_t^{T}
\psi_2(t_2)\int\limits_t^{t_2}
\psi_1(t_1) d{\bf w}_{t_1}^{(i)} d{\bf w}_{t_2}^{(i)}
$$

\vspace{3mm}
\noindent
is the iterated Ito stochastic integral of the form (\ref{ito})
($k=2$).

On the other hand, it is not difficult to show that

\vspace{-1mm}
\begin{equation}
\label{novorigin4}
\int\limits_t^{T}
\psi_2(t_2)\int\limits_t^{t_2}
\psi_1(t_1) \circ d{\bf w}_{t_1}^{(i)}\circ d{\bf w}_{t_2}^{(j)}
=\int\limits_t^{T}
\psi_2(t_2)\int\limits_t^{t_2}
\psi_1(t_1) d{\bf w}_{t_1}^{(i)} d{\bf w}_{t_2}^{(j)}
\end{equation}

\vspace{3mm}
\noindent
w.~p.~1, where $\psi_1(\tau), \psi_2(\tau)\in L_2([t, T]),$
$i\ne j$ 
$(i, j=1,\ldots,m),$ another notations are the same as in (\ref{novorigin3}).

Combining (\ref{novorigin3}) and (\ref{novorigin4}), we get
(see (\ref{dsds9}))

\vspace{-1mm}
$$
\int\limits_t^{T}
\psi_2(t_2)\int\limits_t^{t_2}
\psi_1(t_1) \circ d{\bf w}_{t_1}^{(i_1)}\circ d{\bf w}_{t_2}^{(i_2)}
=\int\limits_t^{T}
\psi_2(t_2)\int\limits_t^{t_2}
\psi_1(t_1) d{\bf w}_{t_1}^{(i_1)} d{\bf w}_{t_2}^{(i_2)}+
$$

\vspace{1mm}
\begin{equation}
\label{novorigin5}
+
\frac{1}{2}{\bf 1}_{\{i_1=i_2\}}
\int\limits_t^T\psi_1(t_1)\psi_2(t_1)dt_1
\stackrel{\sf def}{=}\bar J^{*}[\psi^{(2)}]_{T,t}^{(i_1 i_2)}
\end{equation}

\vspace{3mm}
\noindent
w.~p.~1, where $\psi_1(\tau),\psi_2(\tau)\in L_2([t, T]),$
$i_1, i_2=1,\ldots,m.$

It is easy to see that the condition 
$\phi_0(x)=1/\sqrt{T-t}$ 
can be omitted in Theorems~40, 41 for the case $k=2$
(see the proof of Theorem~13).

Summing up the above arguments, we obtain
the following generalization of Theorem~7 to the case 
of an arbitrary complete orthonormal 
system of functions in $L_2([t, T])$ and
$\psi_1(\tau), \psi_2(\tau)\in L_2([t, T]).$

\vspace{2mm}

{\bf Theorem~42}\ \cite{2018a}.\ {\it Suppose that 
$\{\phi_j(x)\}_{j=0}^{\infty}$ is an arbitrary complete orthonormal system of 
func\-ti\-ons in the space $L_2([t, T])$ and
$\psi_1(\tau), \psi_2(\tau)\in L_2([t, T])$.
Then$,$ 
for the iterated Stra\-to\-novich stochastic integral

\vspace{-1mm}
$$
J^{S}[\psi^{(2)}]_{T,t}^{(i_1 i_2)}=
\int\limits_t^{T}
\psi_2(t_2)\int\limits_t^{t_2}
\psi_1(t_1) \circ d{\bf f}_{t_1}^{(i_1)}\circ d{\bf f}_{t_2}^{(i_2)}\ \ \ (i_1, i_2=1,\ldots,m)
$$

\vspace{3mm}
\noindent
the following expansion  

\vspace{-1mm}
\begin{equation}
\label{novorigin10}
J^{S}[\psi^{(2)}]_{T,t}^{(i_1 i_2)}=\hbox{\vtop{\offinterlineskip\halign{
\hfil#\hfil\cr
{\rm l.i.m.}\cr
$\stackrel{}{{}_{p_1,p_2\to \infty}}$\cr
}} }\sum_{j_1=0}^{p_1}\sum_{j_2=0}^{p_2}
C_{j_2j_1}\zeta_{j_1}^{(i_1)}\zeta_{j_2}^{(i_2)}
\end{equation}

\vspace{4mm}
\noindent
that converges in the mean-square
sence is valid$,$ where the notations are the same as in Theorems {\rm 6, 7}
and $J^{S}[\psi^{(2)}]_{T,t}^{(i_1 i_2)}$ is defined by {\rm (\ref{dsds7})}. 
}

\vspace{2mm}

In this section, it is also appropriate to mention 
the so-called multiple Stratonovich stochastic integral
\cite{bardina10} (also see \cite{bugh1}).

The mean-square limit (if it exists)

\vspace{2mm}
$$
\hbox{\vtop{\offinterlineskip\halign{
\hfil#\hfil\cr
{\rm l.i.m.}\cr
$\stackrel{}{{}_{N\to \infty}}$\cr
}} }\sum_{l_1=0}^{N-1}\ldots \sum_{l_k=0}^{N-1}
\frac{1}{\Delta\tau_{l_1}\ldots \Delta\tau_{l_k}}
\int\limits_{[\tau_{l_1},\tau_{l_1+1}]\times \ldots \times [\tau_{l_k},\tau_{l_k+1}]}
K(t_1,\ldots,t_k)
dt_1\ldots dt_k\ \Delta{\bf w}_{\tau_{l_1}}^{(i_1)}\ldots
\Delta{\bf w}_{\tau_{l_k}}^{(i_k)}\stackrel{\sf def}{=}
$$

\vspace{4mm}
\begin{equation}
\label{january19}
\stackrel{\sf def}{=}\bar J^{S}[K]_{T,t}^{(i_1\ldots i_k)} 
\end{equation}  

\vspace{2mm}
\noindent
is called \cite{bardina10} the multiple Stratonovich stochastic integral 
of the function  $K(t_1,\ldots,t_k)\in L_2([t, T]^k)$, where
$\Delta {\bf w}_{\tau_j}^{(i)}={\bf w}_{\tau_{j+1}}^{(i)}-{\bf w}_{\tau_{j}}^{(i)}$
$(i=0, 1,\ldots,m),$
$\Delta\tau_j=\tau_{j+1}-\tau_{j},$
$\left\{\tau_j\right\}_{j=0}^N$ 
is a partition of the interval $[t, T]$ 
satisfying the condition (\ref{dsds4}),
$i_1,\ldots,i_k=0,1,\ldots,m,$\ 
${\bf w}_{\tau}^{(i)}=
{\bf f}_{\tau}^{(i)}$ for $i=1,\ldots,m$ and ${\bf w}_{\tau}^{(0)}=\tau$,
${\bf f}_{\tau}^{(i)}$ $(i=1,\ldots,m)$ are independent 
standard Wiener processes defined as above in this section.

Note that in \cite{bardina10} the case $i_1=\ldots=i_k\ne 0$
was considered.
We also denote by $\bar J^{S}[K]_{s,t}^{(i_1\ldots i_k)}$
the multiple Stratonovich stochastic integral 
(\ref{january19}) (if it exists) 
of $K(t_1,\ldots,t_k){\bf 1}_{\{(t_1,\ldots,t_k)\in [t, s]^k\}},$ 
where $K(t_1,\ldots,t_k)\in L_2([t, T]^k),$
$s\in [t, T],$ $t\ge 0.$

Let the function $K(t_1,\ldots,t_k)$ be chosen as follows

\vspace{1mm}
\begin{equation}
\label{january19a}
K(t_1,\ldots,t_k)=
\left\{\begin{matrix}
\psi_1(t_1)\ldots \psi_k(t_k),\ &t_1\le \ldots \le t_k\cr\cr\cr
0,\ &\hbox{\rm otherwise}
\end{matrix}
\right.,
\end{equation}

\vspace{5mm}
\noindent
where $\psi_1(\tau),\ldots, \psi_k(\tau)\in L_2([t, T]),$
$t_1,\ldots,t_k\in [t, T]$ $(k\ge 2)$ and 
$K(t_1)\equiv\psi_1(t_1)$ for $t_1\in[t, T].$

We will denote the 
multiple Stratonovich stochastic integral (\ref{january19})
of the function (\ref{january19a}) as follows
$\bar J^{S}[\psi^{(k)}]_{T,t}^{(i_1\ldots i_k)}$.
It is known \cite{bardina10} (Lemma~A.2) that the Stratonovich 
stochastic integrals
$J^{S}[\psi^{(k)}]_{T,t}^{(i_1\ldots i_k)}$ and
$\bar J^{S}[\psi^{(k)}]_{T,t}^{(i_1\ldots i_k)}$
exist for the case $i_1=\ldots=i_k\ne 0$.
Moreover, 

\vspace{2mm}
$$
J^{S}[\psi^{(k)}]_{T,t}^{(i_1\ldots i_k)}=
\bar J^{S}[\psi^{(k)}]_{T,t}^{(i_1\ldots i_k)}\ \ \ \hbox{w.~p.~1}
$$

\vspace{4mm}
\noindent
for this case \cite{bardina10} (Lemma~A.2).

Recall that an expansion similar to (\ref{after1})
was obtained in {\rm \cite{Rybakov3000}} for the multiple
Stratonovich stochastic integral
(\ref{january19}) under the condition of convergence of trace series.

Recently, 
another approach to the expansion of integral (\ref{january19})
has been proposed (assuming that the integral (\ref{january19}) exists),
where multiple Fourier--Walsh and Fourier--Haar series $(k\in \mathbb{N})$ have been applied
\cite{Rybakov3000xxx}.
The convergence was proved with respect to the special 
subsequence ($p_1=\ldots=p_k=p=2^m,$ $m\to\infty$ in a formula
similar to (\ref{january19c}) \cite{Rybakov3000xxx}).

\vspace{5mm}

\section{Expansion of Iterated Stratonovich Stochastic Integrals
of Multiplicity 3. The Case of an Arbitrary Complete Orthonormal System of 
Functions $(\phi_0(x)=1/\sqrt{T-t})$ 
in the Space $L_2([t,T])$ and $\psi_1(\tau), \psi_2(\tau), \psi_3(\tau)
\equiv 1$}

\vspace{5mm}

In this section, we will prove the following theorem.

\vspace{2mm}

{\bf Theorem~43}\ \cite{2018a}, \cite{arxiv-5}, \cite{arxiv-9}.\ {\it Suppose that
$\{\phi_j(x)\}_{j=0}^{\infty}$ $(\phi_0(x)=1/\sqrt{T-t})$ is an arbitrary complete orthonormal system of 
func\-ti\-ons in the space $L_2([t,T]).$
Then$,$ for the iterated Stra\-to\-no\-vich stochastic integral
of third multiplicity 

$$
{\int\limits_t^{*}}^T
{\int\limits_t^{*}}^{t_3}
{\int\limits_t^{*}}^{t_2}
d{\bf w}_{t_1}^{(i_1)}
d{\bf w}_{t_2}^{(i_2)}d{\bf w}_{t_3}^{(i_3)}\ \ \ (i_1,i_2,i_3=0,1,\ldots,m)
$$

\vspace{3mm}
\noindent
the following expansion 

\vspace{-1mm}
\begin{equation}
\label{2023novem1}
{\int\limits_t^{*}}^T
{\int\limits_t^{*}}^{t_3}
{\int\limits_t^{*}}^{t_2}
d{\bf w}_{t_1}^{(i_1)}
d{\bf w}_{t_2}^{(i_2)}d{\bf w}_{t_3}^{(i_3)}=
\hbox{\vtop{\offinterlineskip\halign{
\hfil#\hfil\cr
{\rm l.i.m.}\cr
$\stackrel{}{{}_{p\to \infty}}$\cr
}} }\sum_{j_1,j_2,j_3=0}^{p}
C_{j_3 j_2 j_1}\zeta_{j_1}^{(i_1)}\zeta_{j_2}^{(i_2)}\zeta_{j_3}^{(i_3)}
\end{equation}

\vspace{3mm}
\noindent
that converges in the mean-square sense is valid, where 

\vspace{-1mm}
$$
C_{j_3 j_2 j_1}=\int\limits_t^T
\phi_{j_3}(t_3)\int\limits_t^{t_3}
\phi_{j_2}(t_2)
\int\limits_t^{t_2}
\phi_{j_1}(t_1)dt_1dt_2dt_3
$$

\vspace{3mm}
and
$$
\zeta_{j}^{(i)}=
\int\limits_t^T \phi_{j}(\tau) d{\bf w}_{\tau}^{(i)}
$$ 

\vspace{3mm}
\noindent
are independent standard Gaussian random variables for various 
$i$ or $j$ {\rm (}in the case when $i\ne 0${\rm ),}
${\bf w}_{\tau}^{(i)}={\bf f}_{\tau}^{(i)}$ for
$i=1,\ldots,m$ and 
${\bf w}_{\tau}^{(0)}=\tau.$}

\vspace{2mm}

{\bf Proof.} First, note that under the conditions of Theorem~43
the equality

$$
\bar J^{*}[\psi^{(3)}]_{T,t}^{(i_1 i_2 i_3)}=
{\int\limits_t^{*}}^T
{\int\limits_t^{*}}^{t_3}
{\int\limits_t^{*}}^{t_2}
d{\bf w}_{t_1}^{(i_1)}
d{\bf w}_{t_2}^{(i_2)}d{\bf w}_{t_3}^{(i_3)}
$$

\vspace{3mm}
\noindent
is true w.~p.~1 (see Theorem~5), where $\bar J^{*}[\psi^{(3)}]_{T,t}^{(i_1 i_2 i_3)}$
is defined by (\ref{dsds9}).

According to Theorem~40, we come to the conclusion that 
Theorem~43 will be proved if we prove the following
equalities

\begin{equation}
\label{2023novem2}
\lim\limits_{p\to\infty}
\sum\limits_{j_3=0}^{p}
\left(~\sum\limits_{j_1=0}^{p} 
C_{j_3 j_2 j_1}\biggl|_{j_{1}=j_{2}}-
\frac{1}{2} 
C_{j_3 j_2 j_1}\biggl|_{(j_{1} j_{2})\curvearrowright (\cdot),j_{1}=j_{2}}
\biggr.\right)^2=0,
\end{equation}

\vspace{2mm}
\begin{equation}
\label{2023novem3}
\lim\limits_{p\to\infty}
\sum\limits_{j_1=0}^{p}
\left(~\sum\limits_{j_3=0}^{p} 
C_{j_3 j_2 j_1}\biggl|_{j_{2}=j_{3}}-
\frac{1}{2} 
C_{j_3 j_2 j_1}\biggl|_{(j_{2} j_{3})\curvearrowright (\cdot), j_{2}=j_{3}}
\biggr.\right)^2=0,
\end{equation}

\vspace{2mm}
\begin{equation}
\label{2023novem4}
\lim\limits_{p\to\infty}
\sum\limits_{j_2=0}^{p}
\left(~\sum\limits_{j_1=0}^{p} 
C_{j_3 j_2 j_1}\biggl|_{j_{1}=j_{3}}\right)^2=0.
\end{equation}

\vspace{5mm}

Note that using Theorem~41 (also see (\ref{after1400})), we can rewrite
the relations (\ref{2023novem2})--(\ref{2023novem4})) in the form
(compare with (\ref{after1600})--(\ref{after1602}))

$$
\lim\limits_{p\to\infty}
\sum\limits_{j_3=0}^{p}
\left(~\sum\limits_{j_1=p+1}^{\infty} 
C_{j_3 j_2 j_1}\biggl|_{j_{1}=j_{2}}\right)^2=0,\ \ \
\lim\limits_{p\to\infty}
\sum\limits_{j_1=0}^{p}
\left(~\sum\limits_{j_3=p+1}^{\infty} 
C_{j_3 j_2 j_1}\biggl|_{j_{2}=j_{3}}\right)^2=0,
$$

\vspace{2mm}
$$
\label{2023novem7}
\lim\limits_{p\to\infty}
\sum\limits_{j_2=0}^{p}
\left(~\sum\limits_{j_1=p+1}^{\infty} 
C_{j_3 j_2 j_1}\biggl|_{j_{1}=j_{3}}\right)^2=0.
$$

\vspace{5mm}

Let us prove (\ref{2023novem2}). Using Fubini's Theorem and Parseval's equality, we have

$$
\lim\limits_{p\to\infty}
\sum\limits_{j_3=0}^{p}
\left(~\sum\limits_{j_1=0}^{p} 
C_{j_3 j_2 j_1}\biggl|_{j_{1}=j_{2}}-
\frac{1}{2} 
C_{j_3 j_2 j_1}\biggl|_{(j_{1} j_{2})\curvearrowright (\cdot),j_{1}=j_{2}}
\biggr.\right)^2=
$$

\vspace{2mm}
$$
=\lim\limits_{p\to\infty}
\sum\limits_{j_3=0}^{p}
\left(
\frac{1}{2} 
C_{j_3 j_2 j_1}\biggl|_{(j_{1} j_{2})\curvearrowright (\cdot),j_{1}=j_{2}}
\biggr.-
\sum\limits_{j_1=0}^{p} C_{j_3 j_1 j_1}
\right)^2=
$$

\vspace{2mm}
$$
=\lim\limits_{p\to\infty}
\sum\limits_{j_3=0}^{p}
\left(\int\limits_t^T 
\phi_{j_3}(\tau)\left(\frac{1}{2}\int\limits_t^{\tau} ds
-
\sum\limits_{j_1=0}^p
\frac{1}{2}\left(\int\limits_t^{\tau}\phi_{j_1}(s)ds\right)^2
\right)d\tau\right)^2\le
$$

\vspace{2mm}
$$
\le\lim\limits_{p\to\infty}
\sum\limits_{j_3=0}^{\infty}
\left(\int\limits_t^T 
\phi_{j_3}(\tau)\left(\frac{1}{2}(\tau-t)
-
\sum\limits_{j_1=0}^p
\frac{1}{2}\left(\int\limits_t^{\tau}\phi_{j_1}(s)ds\right)^2
\right)d\tau\right)^2=
$$

\vspace{2mm}
\begin{equation}
\label{2023novem8}
=\lim\limits_{p\to\infty}
\int\limits_t^T 
\left(\frac{1}{2}(\tau-t)
-
\sum\limits_{j_1=0}^p
\frac{1}{2}\left(\int\limits_t^{\tau}\phi_{j_1}(s)ds\right)^2
\right)^2 d\tau.
\end{equation}

\vspace{5mm}

Applying the Parseval equality, we have

\vspace{-1mm}
$$
\sum\limits_{j_1=0}^{\infty}
\frac{1}{2}\left(\int\limits_t^{\tau}\phi_{j_1}(s)ds\right)^2
=\sum\limits_{j_1=0}^{\infty}
\frac{1}{2}\left(\int\limits_t^T {\bf 1}_{\{s<\tau\}}\phi_{j_1}(s)ds\right)^2=
$$

\vspace{2mm}
\begin{equation}
\label{2023novem10}
=\frac{1}{2}\int\limits_t^T \left({\bf 1}_{\{s<\tau\}}\right)^2 ds=
\frac{1}{2}(\tau-t).
\end{equation}

\vspace{2mm}

Moreover, 
\begin{equation}
\label{2023novem11}
\left\vert
\frac{1}{2}(\tau-t)
-
\sum\limits_{j_1=0}^p
\frac{1}{2}\left(\int\limits_t^{\tau}\phi_{j_1}(s)ds\right)^2\right\vert\le
\frac{1}{2}(\tau-t)\le \frac{1}{2}(T-t)<\infty.
\end{equation}

\vspace{5mm}

Using (\ref{2023novem10}), (\ref{2023novem11}) and
applying Lebesgue's 
Dominated Convergence Theorem in (\ref{2023novem8}), we obtain
the equality (\ref{2023novem2}).

Note that we could  use Dini's Theorem 
instead of Lebesgue's 
Dominated Convergence Theorem. Using
the continuity of the functions $u_p(\tau)$ (see below),
the nondecreasing property of
the functional sequence

$$
u_p(\tau)=
\sum\limits_{j_1=0}^{p}
\frac{1}{2}\left(\int\limits_t^{\tau}\phi_{j_1}(s)ds\right)^2,
$$

\vspace{3mm}
\noindent
and the continuity of the limit function
$u(\tau)=(\tau-t)/2$
according to Dini's 
Theorem,
we have the uniform convergence 
$u_p(\tau)$ to $u(\tau)$ at the interval $[t, T]$.
Then we can swap the limit and integral in (\ref{2023novem8})
and get (\ref{2023novem2}).

Let us prove (\ref{2023novem3}). Using Fubini's Theorem and Parseval's equality, we obtain

$$
\lim\limits_{p\to\infty}
\sum\limits_{j_1=0}^{p}
\left(~\sum\limits_{j_3=0}^{p} 
C_{j_3 j_2 j_1}\biggl|_{j_{2}=j_{3}}-
\frac{1}{2} 
C_{j_3 j_2 j_1}\biggl|_{(j_{2} j_{3})\curvearrowright (\cdot),j_{2}=j_{3}}
\biggr.\right)^2=
$$

\vspace{2mm}
$$
=\lim\limits_{p\to\infty}
\sum\limits_{j_1=0}^{p}
\left(
\frac{1}{2} 
C_{j_3 j_2 j_1}\biggl|_{(j_{2} j_{3})\curvearrowright (\cdot),j_{2}=j_{3}}
\biggr.-
\sum\limits_{j_3=0}^{p} C_{j_3 j_3 j_1}
\right)^2=
$$

\vspace{2mm}
$$
=\lim\limits_{p\to\infty}
\sum\limits_{j_1=0}^{p}
\left(\frac{1}{2}\int\limits_t^T 
\int\limits_t^{\tau} \phi_{j_1}(s)dsd\tau
-
\sum\limits_{j_3=0}^p
\int\limits_t^T \phi_{j_3}(\theta)\int\limits_t^{\theta} \phi_{j_3}(\tau)
\int\limits_t^{\tau} \phi_{j_1}(s)ds d\tau d\theta\right)^2=
$$

\vspace{2mm}
$$
=\lim\limits_{p\to\infty}
\sum\limits_{j_1=0}^{p}
\left(\frac{1}{2}\int\limits_t^T 
\phi_{j_1}(s)(T-s)ds
-
\sum\limits_{j_3=0}^p
\int\limits_t^T \phi_{j_1}(s)\int\limits_{s}^T \phi_{j_3}(\tau)
\int\limits_{\tau}^T \phi_{j_3}(\theta)d\theta d\tau ds\right)^2=
$$

\vspace{2mm}
$$
=\lim\limits_{p\to\infty}
\sum\limits_{j_1=0}^{p}
\left(\int\limits_t^T 
\phi_{j_1}(s)\left(\frac{1}{2}(T-s)
-
\sum\limits_{j_3=0}^p
\frac{1}{2}\left(\int\limits_{s}^T \phi_{j_3}(\tau)d\tau\right)^2\right) ds\right)^2\le
$$

\vspace{2mm}
$$
\le \lim\limits_{p\to\infty}
\sum\limits_{j_1=0}^{\infty}
\left(\int\limits_t^T 
\phi_{j_1}(s)\left(\frac{1}{2}(T-s)
-
\sum\limits_{j_3=0}^p
\frac{1}{2}\left(\int\limits_{s}^T \phi_{j_3}(\tau)d\tau\right)^2\right) ds\right)^2=
$$

\vspace{2mm}
\begin{equation}
\label{2023novem14}
=\lim\limits_{p\to\infty}
\int\limits_t^T 
\left(\frac{1}{2}(T-s)
-
\sum\limits_{j_3=0}^p
\frac{1}{2}\left(\int\limits_{s}^T \phi_{j_3}(\tau)d\tau\right)^2\right)^2 ds.
\end{equation}

\vspace{5mm}

Using the Parseval equality, we get

$$
\sum\limits_{j_3=0}^{\infty}
\frac{1}{2}\left(\int\limits_{s}^T \phi_{j_3}(\tau)d\tau\right)^2=
\sum\limits_{j_3=0}^{\infty}
\frac{1}{2}\left(\int\limits_t^T {\bf 1}_{\{s<\tau\}}\phi_{j_3}(\tau)d\tau\right)^2=
$$

\vspace{2mm}
\begin{equation}
\label{2023novem15}
=\frac{1}{2}\int\limits_t^T \left({\bf 1}_{\{s<\tau\}}\right)^2 d\tau=
\frac{1}{2}(T-s).
\end{equation}

\vspace{2mm}

Moreover, 
\begin{equation}
\label{2023novem16}
\left\vert
\frac{1}{2}(T-s)
-
\sum\limits_{j_3=0}^p
\frac{1}{2}\left(\int\limits_{s}^T \phi_{j_3}(\tau)d\tau\right)^2
\right\vert\le
\frac{1}{2}(T-s)\le \frac{1}{2}(T-t)<\infty.
\end{equation}

\vspace{5mm}

Combining (\ref{2023novem14})--(\ref{2023novem16}) and using
the same reasoning as in the proof of (\ref{2023novem2}), we obtain

$$
\lim\limits_{p\to\infty}
\int\limits_t^T 
\left(\frac{1}{2}(T-s)
-
\sum\limits_{j_3=0}^p
\frac{1}{2}\left(\int\limits_{s}^T \phi_{j_3}(\tau)d\tau\right)^2\right)^2 ds=0.
$$

\vspace{5mm}

The equality (\ref{2023novem3}) is proved.

Let us prove (\ref{2023novem4}). Applying Fubini's Theorem and Parseval's equality, we have

$$
\lim\limits_{p\to\infty}
\sum\limits_{j_2=0}^{p}
\left(~\sum\limits_{j_1=0}^{p} 
C_{j_1 j_2 j_1}\right)^2=
$$

\vspace{2mm}
$$
=\lim\limits_{p\to\infty}
\sum\limits_{j_2=0}^{p}
\left(~\sum\limits_{j_1=0}^{p} 
\int\limits_t^T \phi_{j_1}(\theta)\int\limits_t^{\theta}
\phi_{j_2}(\tau)\int\limits_t^{\tau} \phi_{j_1}(s)ds d\tau d\theta\right)^2=
$$

\vspace{2mm}
$$
=\lim\limits_{p\to\infty}
\sum\limits_{j_2=0}^{p}
\left(~\sum\limits_{j_1=0}^{p} 
\int\limits_t^T 
\phi_{j_2}(\tau)\int\limits_t^{\tau} \phi_{j_1}(s)ds
\int\limits_{\tau}^T
\phi_{j_1}(\theta)d\theta d\tau \right)^2\le
$$

\vspace{2mm}
$$
\le\lim\limits_{p\to\infty}
\sum\limits_{j_2=0}^{\infty}
\left(
\int\limits_t^T 
\phi_{j_2}(\tau)\sum\limits_{j_1=0}^{p}\int\limits_t^{\tau} \phi_{j_1}(s)ds
\int\limits_{\tau}^T
\phi_{j_1}(\theta)d\theta d\tau \right)^2=
$$

\vspace{2mm}
\begin{equation}
\label{2023novem18}
=\lim\limits_{p\to\infty}
\int\limits_t^T 
\left(\sum\limits_{j_1=0}^{p}\int\limits_t^{\tau} \phi_{j_1}(s)ds
\int\limits_{\tau}^T
\phi_{j_1}(\theta)d\theta\right)^2 d\tau.
\end{equation}

\vspace{5mm}

Applying (\ref{dsds14}), we obtain

$$
\left\vert
\sum\limits_{j_1=0}^{p}\int\limits_t^{\tau} \phi_{j_1}(s)ds
\int\limits_{\tau}^T
\phi_{j_1}(\theta)d\theta\right\vert\le
\sum\limits_{j_1=0}^{p}\left\vert
\int\limits_t^{\tau} \phi_{j_1}(s)ds
\int\limits_{\tau}^T
\phi_{j_1}(\theta)d\theta\right\vert\le
$$

\vspace{2mm}
\begin{equation}
\label{2023novem19}
\le 
\sum\limits_{j_1=0}^{\infty}\left\vert
\int\limits_t^{\tau} \phi_{j_1}(s)ds
\int\limits_{\tau}^T
\phi_{j_1}(\theta)d\theta\right\vert\le \frac{1}{2}(T-t)<\infty.
\end{equation}

\vspace{5mm}

Using the generalized Parseval equality, we get

$$
\lim\limits_{p\to\infty}
\sum\limits_{j_1=0}^{p}\int\limits_t^{\tau} \phi_{j_1}(s)ds
\int\limits_{\tau}^T
\phi_{j_1}(\theta)d\theta= 
\sum\limits_{j_1=0}^{\infty}\int\limits_t^T {\bf 1}_{\{s<\tau\}}\phi_{j_1}(s)ds
\int\limits_{t}^T
{\bf 1}_{\{s>\tau\}}
\phi_{j_1}(s)ds= 
$$

\vspace{2mm}
\begin{equation}
\label{2023novem20}
=
\int\limits_t^T {\bf 1}_{\{s<\tau\}}{\bf 1}_{\{s>\tau\}}ds=0.
\end{equation}

\vspace{5mm}

Taking into account (\ref{2023novem19}), (\ref{2023novem20}) and
applying Lebesgue's 
Dominated Convergence Theorem in (\ref{2023novem18}), we obtain
the equality (\ref{2023novem4}). Theorem~43 is proved.

\vspace{5mm}

\section{Expansion of Iterated Stratonovich Stochastic Integrals
of Multiplicity 4. The Case of an Arbitrary Complete Orthonormal System of 
Functions $(\phi_0(x)=1/\sqrt{T-t})$ 
in the Space $L_2([t,T])$ and $\psi_1(\tau),\ldots, \psi_4(\tau)
\equiv 1$}

\vspace{5mm}

In this section, we will prove the following theorem.

\vspace{2mm}                                       

{\bf Theorem~44}\ \cite{2018a}, \cite{arxiv-5}, \cite{arxiv-9}.\ {\it Suppose that
$\{\phi_j(x)\}_{j=0}^{\infty}$ $(\phi_0(x)=1/\sqrt{T-t})$ is an arbitrary complete orthonormal system of 
func\-ti\-ons in the space $L_2([t,T]).$
Then$,$ for the iterated Stra\-to\-no\-vich stochastic integral
of fourth multiplicity 

$$
J^{*}[\psi^{(4)}]_{T,t}=
{\int\limits_t^{*}}^T
{\int\limits_t^{*}}^{t_4}
{\int\limits_t^{*}}^{t_3}
{\int\limits_t^{*}}^{t_2}
d{\bf w}_{t_1}^{(i_1)}
d{\bf w}_{t_2}^{(i_2)}d{\bf w}_{t_3}^{(i_3)}d{\bf w}_{t_4}^{(i_4)}\ \ \ 
(i_1, i_2, i_3, i_4=0, 1,\ldots,m)
$$

\vspace{3mm}
\noindent
the following 
expansion 

\vspace{-1mm}
$$
J^{*}[\psi^{(4)}]_{T,t}=
\hbox{\vtop{\offinterlineskip\halign{
\hfil#\hfil\cr
{\rm l.i.m.}\cr
$\stackrel{}{{}_{p\to \infty}}$\cr
}} }
\sum\limits_{j_1, j_2, j_3, j_4=0}^{p}
C_{j_4 j_3 j_2 j_1}\zeta_{j_1}^{(i_1)}\zeta_{j_2}^{(i_2)}\zeta_{j_3}^{(i_3)}
\zeta_{j_4}^{(i_4)}
$$

\vspace{3mm}
\noindent
that converges in the mean-square sense is valid, where 

\vspace{-1mm}
$$
C_{j_4 j_3 j_2 j_1}=\int\limits_t^T
\phi_{j_4}(t_4)\int\limits_t^{t_4}
\phi_{j_3}(t_3)\int\limits_t^{t_3}
\phi_{j_2}(t_2)\int\limits_t^{t_2}
\phi_{j_1}(t_1)dt_1dt_2dt_3 dt_4
$$

\vspace{3mm}
\noindent
and
$$
\zeta_{j}^{(i)}=
\int\limits_t^T \phi_{j}(\tau) d{\bf w}_{\tau}^{(i)}
$$ 

\vspace{3mm}
\noindent
are independent standard Gaussian random variables for various 
$i$ or $j$ {\rm (}in the case when $i\ne 0${\rm ),}
${\bf w}_{\tau}^{(i)}={\bf f}_{\tau}^{(i)}$ for
$i=1,\ldots,m$ and 
${\bf w}_{\tau}^{(0)}=\tau.$}

\vspace{2mm}

{\bf Proof.} First, note that under the conditions of Theorem~44
the equality

$$
\bar J^{*}[\psi^{(4)}]_{T,t}^{(i_1 i_2 i_3 i_4)}=
{\int\limits_t^{*}}^T
{\int\limits_t^{*}}^{t_4}
{\int\limits_t^{*}}^{t_3}
{\int\limits_t^{*}}^{t_2}
d{\bf w}_{t_1}^{(i_1)}
d{\bf w}_{t_2}^{(i_2)}d{\bf w}_{t_3}^{(i_3)}d{\bf w}_{t_4}^{(i_4)}
$$

\vspace{3mm}
\noindent
is valid w.~p.~1 (see Theorem~5), where $\bar J^{*}[\psi^{(4)}]_{T,t}^{(i_1 i_2 i_3 i_4)}$
is defined by (\ref{dsds9}).

It is easy to see that Theorem~44 will be proved if we prove the following
equalities (see Theorem~40)

\begin{equation}
\label{2023novem200}
\lim\limits_{p\to\infty}
\sum\limits_{j_3,j_4=0}^{p}
\left(~\sum\limits_{j_1=0}^{p} 
C_{j_4 j_3 j_1 j_1}-
\frac{1}{2} 
C_{j_4 j_3 j_1 j_1}\biggl|_{(j_{1} j_{1})\curvearrowright (\cdot)}
\biggr.\right)^2=0,
\end{equation}

\vspace{2mm}
\begin{equation}
\label{2023novem201}
\lim\limits_{p\to\infty}
\sum\limits_{j_2,j_4=0}^{p}
\left(~\sum\limits_{j_1=0}^{p} 
C_{j_4 j_1 j_2 j_1}\right)^2=0,
\end{equation}

\vspace{2mm}
\begin{equation}
\label{2023novem202}
\lim\limits_{p\to\infty}
\sum\limits_{j_2, j_3=0}^{p}
\left(~\sum\limits_{j_1=0}^{p} 
C_{j_1 j_3 j_2 j_1}\right)^2=0,
\end{equation}

\vspace{2mm}
\begin{equation}
\label{2023novem203}
\lim\limits_{p\to\infty}
\sum\limits_{j_1,j_4=0}^{p}
\left(~\sum\limits_{j_2=0}^{p} 
C_{j_4 j_2 j_2 j_1}-
\frac{1}{2} 
C_{j_4 j_2 j_2 j_1}\biggl|_{(j_{2} j_{2})\curvearrowright (\cdot)}
\biggr.\right)^2=0,
\end{equation}

\vspace{2mm}
\begin{equation}
\label{2023novem204}
\lim\limits_{p\to\infty}
\sum\limits_{j_1, j_3=0}^{p}
\left(~\sum\limits_{j_2=0}^{p} 
C_{j_2 j_3 j_2 j_1}\right)^2=0,
\end{equation}

\vspace{2mm}
\begin{equation}
\label{2023novem205}
\lim\limits_{p\to\infty}
\sum\limits_{j_1,j_2=0}^{p}
\left(~\sum\limits_{j_3=0}^{p} 
C_{j_3 j_3 j_2 j_1}-
\frac{1}{2} 
C_{j_3 j_3 j_2 j_1}\biggl|_{(j_{3} j_{3})\curvearrowright (\cdot)}
\biggr.\right)^2=0,
\end{equation}

\vspace{2mm}
\begin{equation}
\label{2023novem206}
\lim\limits_{p\to\infty}
\sum\limits_{j_1, j_3=0}^{p}
C_{j_3 j_3 j_1 j_1}=\frac{1}{4} 
C_{j_3 j_3 j_1 j_1}\biggl|_{(j_{3} j_{3})\curvearrowright (\cdot)
(j_{1} j_{1})\curvearrowright (\cdot)}=\frac{1}{8}(T-t)^2,
\biggr.
\end{equation}

\vspace{2mm}
\begin{equation}
\label{2023novem207}
\lim\limits_{p\to\infty}
\sum\limits_{j_1, j_3=0}^{p}
C_{j_1 j_3 j_3 j_1}=0,
\biggr.
\end{equation}

\vspace{2mm}
\begin{equation}
\label{2023novem208}
\lim\limits_{p\to\infty}
\sum\limits_{j_1, j_2=0}^{p}
C_{j_2 j_1 j_2 j_1}=0.
\biggr.
\end{equation}

\vspace{5mm}

Let us prove the equalities (\ref{2023novem200})--(\ref{2023novem205}).
Using Fubini's Theorem and Parseval's equality, we obtain
the following relations for the prelimit
expressions on the left-hand sides of (\ref{2023novem200})--(\ref{2023novem205})

$$
\sum\limits_{j_3,j_4=0}^{p}
\left(~\sum\limits_{j_1=0}^{p} 
C_{j_4 j_3 j_1 j_1}-
\frac{1}{2} 
C_{j_4 j_3 j_1 j_1}\biggl|_{(j_{1} j_{1})\curvearrowright (\cdot)}
\biggr.\right)^2=
$$

\vspace{2mm}
$$
=\sum\limits_{j_3,j_4=0}^{p}\left(
\frac{1}{2}\int\limits_t^T \phi_{j_4}(t_4)
\int\limits_t^{t_4}\phi_{j_3}(t_3)(t_3-t)dt_3 dt_4-\right.
$$

\vspace{2mm}
$$
\left.-\sum\limits_{j_1=0}^p\int\limits_t^T \phi_{j_4}(t_4)
\int\limits_t^{t_4}\phi_{j_3}(t_3)
\int\limits_t^{t_3}\phi_{j_1}(t_2)
\int\limits_t^{t_2}\phi_{j_1}(t_1)dt_1 dt_2 dt_3 dt_4\right)^2=
$$

\vspace{2mm}
$$
=\sum\limits_{j_3,j_4=0}^{p}\left(
\int\limits_t^T \phi_{j_4}(t_4)
\int\limits_t^{t_4}\phi_{j_3}(t_3)\Biggl(\frac{1}{2}(t_3-t)-\Biggr.\right.
$$

\vspace{2mm}
$$
\left.\Biggl.-\sum\limits_{j_1=0}^p
\int\limits_t^{t_3}\phi_{j_1}(t_2)
\int\limits_t^{t_2}\phi_{j_1}(t_1)dt_1 dt_2\Biggr) dt_3 dt_4\right)^2=
$$

\vspace{2mm}
$$
=\sum\limits_{j_3,j_4=0}^{p}\left(
\int\limits_t^T \phi_{j_4}(t_4)
\int\limits_t^{t_4}\phi_{j_3}(t_3)\left(\frac{1}{2}(t_3-t)
-\sum\limits_{j_1=0}^p
\frac{1}{2}\left(\int\limits_t^{t_3}\phi_{j_1}(s)ds\right)^2
\right) dt_3 dt_4\right)^2\le
$$

\vspace{2mm}
$$
\le\sum\limits_{j_3,j_4=0}^{\infty}\left(
\int\limits_t^T \phi_{j_4}(t_4)
\int\limits_t^{t_4}\phi_{j_3}(t_3)\left(\frac{1}{2}(t_3-t)
-\sum\limits_{j_1=0}^p
\frac{1}{2}\left(\int\limits_t^{t_3}\phi_{j_1}(s)ds\right)^2
\right) dt_3 dt_4\right)^2=
$$

\vspace{2mm}
\begin{equation}
\label{2023novem209}
=
\int\limits_{[t, T]^2}{\bf 1}_{\{t_3<t_4\}}\left(\frac{1}{2}(t_3-t)
-\sum\limits_{j_1=0}^p
\frac{1}{2}\left(\int\limits_t^{t_3}\phi_{j_1}(s)ds\right)^2
\right)^2 dt_3 dt_4,
\end{equation}

\vspace{6mm}

$$
\sum\limits_{j_2,j_4=0}^{p}
\left(~\sum\limits_{j_1=0}^{p} 
C_{j_4 j_1 j_2 j_1}\right)^2=
$$

\vspace{2mm}
$$
=
\sum\limits_{j_2,j_4=0}^{p}\left(
\sum\limits_{j_1=0}^p\int\limits_t^T \phi_{j_4}(t_4)
\int\limits_t^{t_4}\phi_{j_1}(t_3)
\int\limits_t^{t_3}\phi_{j_2}(t_2)
\int\limits_t^{t_2}\phi_{j_1}(t_1)dt_1 dt_2 dt_3 dt_4\right)^2=
$$

\vspace{2mm}
$$
=
\sum\limits_{j_2,j_4=0}^{p}\left(
\sum\limits_{j_1=0}^p
\int\limits_t^T \phi_{j_4}(t_4)
\int\limits_t^{t_4}\phi_{j_2}(t_2)
\int\limits_t^{t_2}\phi_{j_1}(t_1)dt_1
\int\limits_{t_2}^{t_4}\phi_{j_1}(t_3)dt_3 dt_2 dt_4\right)^2=
$$

\vspace{2mm}
$$
=
\sum\limits_{j_2,j_4=0}^{p}\left(
\int\limits_t^T \phi_{j_4}(t_4)
\int\limits_t^{t_4}\phi_{j_2}(t_2)
\sum\limits_{j_1=0}^p
\int\limits_t^{t_2}\phi_{j_1}(t_1)dt_1
\int\limits_{t_2}^{t_4}\phi_{j_1}(t_3)dt_3 dt_2 dt_4\right)^2\le
$$

\vspace{2mm}
$$
\le
\sum\limits_{j_2,j_4=0}^{\infty}\left(
\int\limits_t^T \phi_{j_4}(t_4)
\int\limits_t^{t_4}\phi_{j_2}(t_2)
\sum\limits_{j_1=0}^p
\int\limits_t^{t_2}\phi_{j_1}(t_1)dt_1
\int\limits_{t_2}^{t_4}\phi_{j_1}(t_3)dt_3 dt_2 dt_4\right)^2=
$$

\vspace{2mm}
\begin{equation}
\label{2023novem210}
=
\int\limits_{[t,T]^2} {\bf 1}_{\{t_2<t_4\}}
\left(\sum\limits_{j_1=0}^p
\int\limits_t^{t_2}\phi_{j_1}(t_1)dt_1
\int\limits_{t_2}^{t_4}\phi_{j_1}(t_3)dt_3\right)^2 dt_2 dt_4,
\end{equation}

\vspace{6mm}

$$
\sum\limits_{j_2,j_3=0}^{p}
\left(~\sum\limits_{j_1=0}^{p} 
C_{j_1 j_3 j_2 j_1}\right)^2=
$$

\vspace{2mm}
$$
=
\sum\limits_{j_2,j_3=0}^{p}\left(
\sum\limits_{j_1=0}^p
\int\limits_t^T \phi_{j_1}(t_4)
\int\limits_t^{t_4}\phi_{j_3}(t_3)
\int\limits_t^{t_3}\phi_{j_2}(t_2)
\int\limits_t^{t_2}\phi_{j_1}(t_1)dt_1 dt_2 dt_3 dt_4\right)^2=
$$

\vspace{2mm}
$$
=
\sum\limits_{j_2,j_3=0}^{p}\left(
\sum\limits_{j_1=0}^p
\int\limits_t^T \phi_{j_3}(t_3)
\int\limits_t^{t_3}\phi_{j_2}(t_2)
\int\limits_t^{t_2}\phi_{j_1}(t_1)dt_1
\int\limits_{t_3}^{T}\phi_{j_1}(t_4)dt_4 dt_2 dt_3\right)^2=
$$

\vspace{2mm}
$$
=
\sum\limits_{j_2,j_3=0}^{p}\left(
\int\limits_t^T \phi_{j_3}(t_3)
\int\limits_t^{t_3}\phi_{j_2}(t_2)
\sum\limits_{j_1=0}^p\int\limits_t^{t_2}\phi_{j_1}(t_1)dt_1
\int\limits_{t_3}^{T}\phi_{j_1}(t_4)dt_4 dt_2 dt_3\right)^2\le
$$

\vspace{2mm}
$$
\le
\sum\limits_{j_2,j_3=0}^{\infty}\left(
\int\limits_t^T \phi_{j_3}(t_3)
\int\limits_t^{t_3}\phi_{j_2}(t_2)
\sum\limits_{j_1=0}^p\int\limits_t^{t_2}\phi_{j_1}(t_1)dt_1
\int\limits_{t_3}^{T}\phi_{j_1}(t_4)dt_4 dt_2 dt_3\right)^2=
$$

\vspace{2mm}
\begin{equation}
\label{2023novem211}
=
\int\limits_{[t, T]^2} {\bf 1}_{\{t_2<t_3\}}
\left(\sum\limits_{j_1=0}^p\int\limits_t^{t_2}\phi_{j_1}(t_1)dt_1
\int\limits_{t_3}^{T}\phi_{j_1}(t_4)dt_4\right)^2 dt_2 dt_3,
\end{equation}

\vspace{6mm}

$$
\sum\limits_{j_1,j_4=0}^{p}
\left(~\sum\limits_{j_2=0}^{p} 
C_{j_4 j_2 j_2 j_1}-
\frac{1}{2} 
C_{j_4 j_2 j_2 j_1}\biggl|_{(j_{2} j_{2})\curvearrowright (\cdot)}
\biggr.\right)^2=
$$

\vspace{2mm}
$$
=\sum\limits_{j_1,j_4=0}^{p}\left(
\frac{1}{2}\int\limits_t^T \phi_{j_4}(t_4)
\int\limits_t^{t_4} \int\limits_t^{t_2}\phi_{j_1}(t_1)dt_1 dt_2 dt_4-\right.
$$

\vspace{2mm}
$$
\left.-\sum\limits_{j_2=0}^p\int\limits_t^T \phi_{j_4}(t_4)
\int\limits_t^{t_4}\phi_{j_2}(t_3)
\int\limits_t^{t_3}\phi_{j_2}(t_2)
\int\limits_t^{t_2}\phi_{j_1}(t_1)dt_1 dt_2 dt_3 dt_4\right)^2=
$$

\vspace{2mm}
$$
=\sum\limits_{j_1,j_4=0}^{p}\left(
\frac{1}{2}\int\limits_t^T \phi_{j_4}(t_4)
\int\limits_t^{t_4}\phi_{j_1}(t_1)\int\limits_{t_1}^{t_4} dt_2 dt_1 dt_4-\right.
$$

\vspace{2mm}
$$
\left.-\sum\limits_{j_2=0}^p\int\limits_t^T \phi_{j_4}(t_4)
\int\limits_t^{t_4}\phi_{j_1}(t_1)
\int\limits_{t_1}^{t_4}\phi_{j_2}(t_2)
\int\limits_{t_2}^{t_4}\phi_{j_2}(t_3)dt_3 dt_2 dt_1 dt_4\right)^2=
$$

\vspace{2mm}
$$
=\sum\limits_{j_1,j_4=0}^{p}\left(
\int\limits_t^T \phi_{j_4}(t_4)
\int\limits_t^{t_4}\phi_{j_1}(t_1)\left(\frac{t_4-t_1}{2}
-\sum\limits_{j_2=0}^p
\frac{1}{2}\left(\int\limits_{t_1}^{t_4}\phi_{j_2}(s)ds\right)^2
\right) dt_1 dt_4\right)^2\le
$$

\vspace{2mm}
$$
\le\sum\limits_{j_1,j_4=0}^{\infty}\left(
\int\limits_t^T \phi_{j_4}(t_4)
\int\limits_t^{t_4}\phi_{j_1}(t_1)\left(\frac{t_4-t_1}{2}
-\sum\limits_{j_2=0}^p
\frac{1}{2}\left(\int\limits_{t_1}^{t_4}\phi_{j_2}(s)ds\right)^2
\right) dt_1 dt_4\right)^2=
$$

\vspace{2mm}
\begin{equation}
\label{2023novem212}
=
\int\limits_{[t,T]^2} {\bf 1}_{\{t_1<t_4\}}\left(\frac{1}{2}(t_4-t_1)
-\sum\limits_{j_2=0}^p
\frac{1}{2}\left(\int\limits_{t_1}^{t_4}\phi_{j_2}(s)ds\right)^2
\right)^2 dt_1 dt_4,
\end{equation}

\vspace{6mm}

$$
\sum\limits_{j_1,j_3=0}^{p}
\left(~\sum\limits_{j_2=0}^{p} 
C_{j_2 j_3 j_2 j_1}\right)^2=
$$

\vspace{2mm}
$$
=
\sum\limits_{j_1,j_3=0}^{p}\left(
\sum\limits_{j_2=0}^p
\int\limits_t^T \phi_{j_2}(t_4)
\int\limits_t^{t_4}\phi_{j_3}(t_3)
\int\limits_t^{t_3}\phi_{j_2}(t_2)
\int\limits_t^{t_2}\phi_{j_1}(t_1)dt_1 dt_2 dt_3 dt_4\right)^2=
$$

\vspace{2mm}
$$
=
\sum\limits_{j_1,j_3=0}^{p}\left(
\sum\limits_{j_2=0}^p
\int\limits_t^T \phi_{j_3}(t_3)
\int\limits_t^{t_3}\phi_{j_2}(t_2)
\int\limits_t^{t_2}\phi_{j_1}(t_1)dt_1 dt_2
\int\limits_{t_3}^{T}\phi_{j_2}(t_4)dt_4 dt_3 \right)^2=
$$

\vspace{2mm}
$$
=
\sum\limits_{j_1,j_3=0}^{p}\left(
\sum\limits_{j_2=0}^p
\int\limits_t^T \phi_{j_3}(t_3)
\int\limits_t^{t_3}\phi_{j_1}(t_1)
\int\limits_{t_1}^{t_3}\phi_{j_2}(t_2)dt_2
\int\limits_{t_3}^{T}\phi_{j_2}(t_4)dt_4 dt_1 dt_3 \right)^2=
$$

\vspace{2mm}
$$
=
\sum\limits_{j_1,j_3=0}^{p}\left(
\int\limits_t^T \phi_{j_3}(t_3)
\int\limits_t^{t_3}\phi_{j_1}(t_1)
\sum\limits_{j_2=0}^p
\int\limits_{t_1}^{t_3}\phi_{j_2}(t_2)dt_2
\int\limits_{t_3}^{T}\phi_{j_2}(t_4)dt_4 dt_1 dt_3 \right)^2\le
$$

\vspace{2mm}
$$
\le
\sum\limits_{j_1,j_3=0}^{\infty}\left(
\int\limits_t^T \phi_{j_3}(t_3)
\int\limits_t^{t_3}\phi_{j_1}(t_1)
\sum\limits_{j_2=0}^p
\int\limits_{t_1}^{t_3}\phi_{j_2}(t_2)dt_2
\int\limits_{t_3}^{T}\phi_{j_2}(t_4)dt_4 dt_1 dt_3 \right)^2=
$$

\vspace{2mm}
\begin{equation}
\label{2023novem213}
=
\int\limits_{[t,T]^2}{\bf 1}_{\{t_1<t_3\}}
\left(\sum\limits_{j_2=0}^p
\int\limits_{t_1}^{t_3}\phi_{j_2}(t_2)dt_2
\int\limits_{t_3}^{T}\phi_{j_2}(t_4)dt_4\right)^2 dt_1 dt_3,
\end{equation}

\vspace{6mm}

$$
\sum\limits_{j_1,j_2=0}^{p}
\left(~\sum\limits_{j_3=0}^{p} 
C_{j_3 j_3 j_2 j_1}-
\frac{1}{2} 
C_{j_3 j_3 j_2 j_1}\biggl|_{(j_{3} j_{3})\curvearrowright (\cdot)}
\biggr.\right)^2=
$$

\vspace{2mm}
$$
=\sum\limits_{j_1,j_2=0}^{p}\left(
\frac{1}{2}\int\limits_t^T
\int\limits_t^{t_3}\phi_{j_2}(t_2)\int\limits_t^{t_2}\phi_{j_1}(t_1)dt_1 dt_2 dt_3-\right.
$$

\vspace{2mm}
$$
\left.-\sum\limits_{j_3=0}^p\int\limits_t^T \phi_{j_3}(t_4)
\int\limits_t^{t_4}\phi_{j_3}(t_3)
\int\limits_t^{t_3}\phi_{j_2}(t_2)
\int\limits_t^{t_2}\phi_{j_1}(t_1)dt_1 dt_2 dt_3 dt_4\right)^2=
$$

\vspace{2mm}
$$
=\sum\limits_{j_1,j_2=0}^{p}\left(
\frac{1}{2}\int\limits_t^T\phi_{j_1}(t_1)
\int\limits_{t_1}^T\phi_{j_2}(t_2)\int\limits_{t_2}^T dt_3 dt_2 dt_1-\right.
$$

\vspace{2mm}
$$
\left.-\sum\limits_{j_3=0}^p\int\limits_t^T \phi_{j_1}(t_1)
\int\limits_{t_1}^T \phi_{j_2}(t_2)
\int\limits_{t_2}^T\phi_{j_3}(t_3)
\int\limits_{t_3}^T\phi_{j_3}(t_4)dt_4 dt_3 dt_2 dt_1\right)^2=
$$

\vspace{2mm}
$$
=\sum\limits_{j_1,j_2=0}^{p}\left(
\int\limits_t^T \phi_{j_1}(t_1)
\int\limits_{t_1}^T\phi_{j_2}(t_2)\left(\frac{T-t_2}{2}
-\sum\limits_{j_3=0}^p
\frac{1}{2}\left(\int\limits_{t_2}^T\phi_{j_3}(s)ds\right)^2
\right) dt_2 dt_1\right)^2\le
$$

\vspace{2mm}
$$
\le\sum\limits_{j_1,j_2=0}^{\infty}\left(
\int\limits_t^T \phi_{j_1}(t_1)
\int\limits_{t_1}^T\phi_{j_2}(t_2)\left(\frac{T-t_2}{2}
-\sum\limits_{j_3=0}^p
\frac{1}{2}\left(\int\limits_{t_2}^T\phi_{j_3}(s)ds\right)^2
\right) dt_2 dt_1\right)^2=
$$

\vspace{2mm}
\begin{equation}
\label{2023novem214}
=
\int\limits_{[t, T]^2}{\bf 1}_{\{t_1<t_2\}}\left(\frac{1}{2}(T-t_2)
-\sum\limits_{j_3=0}^p
\frac{1}{2}\left(\int\limits_{t_2}^T\phi_{j_3}(s)ds\right)^2
\right)^2 dt_2 dt_1.
\end{equation}

\vspace{5mm}

Using Parseval's equality, generalized Parseval's equality and 
Lebesgue's Dominated Convergence Theorem, 
as well as applying the same reasoning as in the proof of Theorem~43, 
we obtain that the right-hand sides of (\ref{2023novem209})--(\ref{2023novem214}) 
tend to zero when $p\to\infty.$
The equalities (\ref{2023novem200})--(\ref{2023novem205}) are proved.

Let us prove the equalities (\ref{2023novem206})--(\ref{2023novem208}).
We will use our idea from Sect.~13. More precisely, we consider
the following analogue of the equality (\ref{sixsix40})

\begin{equation}
\label{2023novem215}
C_{j_4 j_3 j_2 j_1}+C_{j_1 j_2 j_3 j_4}=
C_{j_4}C_{j_3 j_2 j_1}-C_{j_3 j_4}C_{j_2 j_1}+
C_{j_2 j_3 j_4}C_{j_1}.
\end{equation}

\vspace{5mm}

Using Fubini's Theorem, we have

$$
C_{j_4 j_3 j_2 j_1}=
$$

\vspace{2mm}
$$
=\int\limits_t^T\phi_{j_4}(t_4)\int\limits_t^{t_4}\phi_{j_3}(t_3)
\int\limits_t^{t_3}\phi_{j_2}(t_2)
\int\limits_t^{t_2}\phi_{j_1}(t_1)dt_1 dt_2 dt_3 dt_4=
$$

\vspace{2mm}
$$
=\int\limits_t^T\phi_{j_4}(t_4)\int\limits_t^{T}\phi_{j_3}(t_3)
\int\limits_t^{t_3}\phi_{j_2}(t_2)
\int\limits_t^{t_2}\phi_{j_1}(t_1)dt_1 dt_2 dt_3 dt_4-
$$

\vspace{2mm}
$$
-\int\limits_t^T\phi_{j_4}(t_4)\int\limits_{t_4}^T\phi_{j_3}(t_3)
\int\limits_t^{t_3}\phi_{j_2}(t_2)
\int\limits_t^{t_2}\phi_{j_1}(t_1)dt_1 dt_2 dt_3 dt_4=
$$

\vspace{2mm}
$$
=C_{j_4}C_{j_3 j_2 j_1}-
$$

\vspace{2mm}
$$
-\int\limits_t^T\phi_{j_4}(t_4)\int\limits_{t_4}^T\phi_{j_3}(t_3)
\int\limits_t^{T}\phi_{j_2}(t_2)
\int\limits_t^{t_2}\phi_{j_1}(t_1)dt_1 dt_2 dt_3 dt_4+
$$

\vspace{2mm}
$$
+\int\limits_t^T\phi_{j_4}(t_4)\int\limits_{t_4}^T\phi_{j_3}(t_3)
\int\limits_{t_3}^{T}\phi_{j_2}(t_2)
\int\limits_t^{t_2}\phi_{j_1}(t_1)dt_1 dt_2 dt_3 dt_4=
$$

\vspace{2mm}
$$
=C_{j_4}C_{j_3 j_2 j_1}-C_{j_3 j_4}C_{j_2 j_1}+
$$

\vspace{2mm}
$$
+\int\limits_t^T\phi_{j_4}(t_4)\int\limits_{t_4}^T\phi_{j_3}(t_3)
\int\limits_{t_3}^{T}\phi_{j_2}(t_2)\int\limits_{t}^{T}\phi_{j_1}(t_1)
dt_1 dt_2 dt_3 dt_4-
$$

\vspace{2mm}
$$
-\int\limits_t^T\phi_{j_4}(t_4)\int\limits_{t_4}^T\phi_{j_3}(t_3)
\int\limits_{t_3}^{T}\phi_{j_2}(t_2)\int\limits_{t_2}^{T}\phi_{j_1}(t_1)
dt_1 dt_2 dt_3 dt_4=
$$

\vspace{2mm}
\begin{equation}
\label{2023novem216}
=C_{j_4}C_{j_3 j_2 j_1}-C_{j_3 j_4}C_{j_2 j_1}+C_{j_2 j_3 j_4}C_{j_1}-
C_{j_1 j_2 j_3 j_4}.
\end{equation}

\vspace{5mm}

The equality (\ref{2023novem216}) completes the proof of the relation 
(\ref{2023novem215}).

Let us prove (\ref{2023novem206}). Substitute $j_4=j_3,$ $j_2=j_1$ into
(\ref{2023novem215})

\begin{equation}
\label{2023novem217}
C_{j_3 j_3 j_1 j_1}+C_{j_1 j_1 j_3 j_3}=
C_{j_3}C_{j_3 j_1 j_1}-C_{j_3 j_3}C_{j_1 j_1}+
C_{j_1 j_3 j_3}C_{j_1}.
\end{equation}

\vspace{4mm}

From (\ref{2023novem217}) we obtain

$$
\sum\limits_{j_1,j_3=0}^p \bigl(C_{j_3 j_3 j_1 j_1}+C_{j_1 j_1 j_3 j_3}\bigr)=
\sum\limits_{j_1,j_3=0}^p C_{j_3}C_{j_3 j_1 j_1}-\sum\limits_{j_1,j_3=0}^p C_{j_3 j_3}C_{j_1 j_1}+
$$

\vspace{2mm}
$$
+\sum\limits_{j_1,j_3=0}^p C_{j_1 j_3 j_3}C_{j_1}.
$$

\vspace{5mm}

Then
\begin{equation}
\label{2023novem218}
2\sum\limits_{j_1,j_3=0}^p C_{j_3 j_3 j_1 j_1}=
2\sum\limits_{j_1,j_3=0}^p C_{j_3}C_{j_3 j_1 j_1}-\left(\sum\limits_{j_1=0}^p C_{j_1 j_1}\right)^2.
\end{equation}

\vspace{5mm}

From (\ref{2023novem218}) we get

$$
\sum\limits_{j_1,j_3=0}^p C_{j_3 j_3 j_1 j_1}=
\sum\limits_{j_1,j_3=0}^p C_{j_3}C_{j_3 j_1 j_1}-\frac{1}{2}
\left(\sum\limits_{j_1=0}^p C_{j_1 j_1}\right)^2=
$$
\begin{equation}
\label{2023novem219}
=\sum\limits_{j_1,j_3=0}^p C_{j_3}C_{j_3 j_1 j_1}-\frac{1}{2}
\left(\sum\limits_{j_1=0}^p \frac{1}{2}\bigl(C_{j_1}\bigr)^2\right)^2=
\sum\limits_{j_1,j_3=0}^p C_{j_3}C_{j_3 j_1 j_1}-\frac{1}{8}
\left(\sum\limits_{j_1=0}^p \bigl(C_{j_1}\bigr)^2\right)^2.
\end{equation}

\vspace{5mm}

Recall that $\phi_0(\tau)=1/\sqrt{T-t}.$ Then 

\begin{equation}
\label{2023novem220}
C_j=\int\limits_t^T \phi_j(\tau)d\tau=
\left\{
\begin{matrix}
\sqrt{T-t} &\hbox{if}\ j=0\cr\cr
0 &\hbox{if}\ j\ne 0
\end{matrix}.\right.
\end{equation}

\vspace{5mm}

Combining (\ref{2023novem219}), (\ref{2023novem220}) and using Fubini's Theorem, we obtain

$$
\sum\limits_{j_1,j_3=0}^p C_{j_3 j_3 j_1 j_1}=
\sqrt{T-t}\sum\limits_{j_1=0}^p C_{0 j_1 j_1}-\frac{1}{8}(T-t)^2=
$$

\vspace{2mm}
$$
=\sum\limits_{j_1=0}^p \int\limits_t^T \int\limits_t^{t_3}
\phi_{j_1}(t_2)
\int\limits_t^{t_2}
\phi_{j_1}(t_1)dt_1 dt_2 dt_3-\frac{1}{8}(T-t)^2=
$$

\vspace{2mm}
$$
=\sum\limits_{j_1=0}^p \int\limits_t^T \phi_{j_1}(t_1) \int\limits_{t_1}^T
\phi_{j_1}(t_2)
\int\limits_{t_2}^T
dt_3 dt_2 dt_1-\frac{1}{8}(T-t)^2=
$$

\vspace{2mm}
$$
=\sum\limits_{j_1=0}^p \int\limits_t^T \phi_{j_1}(t_1) \int\limits_{t_1}^T
\phi_{j_1}(t_2)
(T-t_2)dt_2 dt_1-\frac{1}{8}(T-t)^2=
$$

\vspace{2mm}
\begin{equation}
\label{2023novem221}
=\sum\limits_{j_1=0}^p \int\limits_t^T 
\phi_{j_1}(t_2)
(T-t_2)\int\limits_{t}^{t_2}
\phi_{j_1}(t_1) dt_1 dt_2-\frac{1}{8}(T-t)^2.
\end{equation}

\vspace{5mm}

Finally applying (\ref{after1400}) and (\ref{2023novem221}), we have

$$
\lim\limits_{p\to\infty}\sum\limits_{j_1,j_3=0}^p C_{j_3 j_3 j_1 j_1}=
\frac{1}{2}\int\limits_t^T
(T-t_2)dt_2-\frac{1}{8}(T-t)^2=\frac{1}{8}(T-t)^2.
$$

\vspace{5mm}

The equality (\ref{2023novem206}) is proved.

Let us prove (\ref{2023novem207}). Substitute $j_4=j_1,$ $j_2=j_3$ into
(\ref{2023novem215})

\begin{equation}
\label{2023novem222}
C_{j_1 j_3 j_3 j_1}+C_{j_1 j_3 j_3 j_1}=
C_{j_1}C_{j_3 j_3 j_1}-C_{j_3 j_1}C_{j_3 j_1}+
C_{j_3 j_3 j_1}C_{j_1}.
\end{equation}

\vspace{5mm}

Using (\ref{2023novem222}), we get

\begin{equation}
\label{2023novem223}
2\sum\limits_{j_1,j_3=0}^p C_{j_1 j_3 j_3 j_1}=
2\sum\limits_{j_1,j_3=0}^p C_{j_1}C_{j_3 j_3 j_1}-\sum\limits_{j_1,j_3=0}^p \bigl(C_{j_3 j_1}\bigr)^2.
\end{equation}

\vspace{5mm}

Then applying (\ref{2023novem223}), (\ref{2023novem220}), Parseval's equality,
and (\ref{after1400}), we obtain

$$
\lim\limits_{p\to\infty}\sum\limits_{j_1,j_3=0}^p C_{j_1 j_3 j_3 j_1}=
\lim\limits_{p\to\infty}\sum\limits_{j_1,j_3=0}^p C_{j_1}C_{j_3 j_3 j_1}-
\frac{1}{2}\lim\limits_{p\to\infty}\sum\limits_{j_1,j_3=0}^p \bigl(C_{j_3 j_1}\bigr)^2=
$$

\vspace{2mm}
$$
=\sqrt{T-t}\sum\limits_{j_3=0}^{\infty} C_{j_3 j_3 0}-
\frac{1}{2}\sum\limits_{j_1,j_3=0}^{\infty}
\left(\int\limits_t^T \phi_{j_3}(t_2)\int\limits_t^{t_2}
\phi_{j_1}(t_1)dt_1 dt_2
\right)^2=
$$

\vspace{2mm}
$$
=\sum\limits_{j_3=0}^{\infty} \int\limits_t^T \phi_{j_3}(t_3)\int\limits_t^{t_3}
\phi_{j_3}(t_2)\int\limits_t^{t_2}dt_1 dt_2 dt_3-
$$

\vspace{2mm}
$$
-
\frac{1}{2}\sum\limits_{j_1,j_3=0}^{\infty}
\left(\int\limits_{[t,T]^2}{\bf 1}_{\{t_1<t_2\}} 
\phi_{j_1}(t_1) \phi_{j_3}(t_2) dt_1 dt_2
\right)^2=
$$

\vspace{2mm}
$$
=\sum\limits_{j_3=0}^{\infty} \int\limits_t^T \phi_{j_3}(t_3)\int\limits_t^{t_3}
\phi_{j_3}(t_2)(t_2-t) dt_2 dt_3-
\frac{1}{2}
\int\limits_{[t,T]^2}\left({\bf 1}_{\{t_1<t_2\}}\right)^2 dt_1 dt_2
=
$$

\vspace{2mm}
$$
=\frac{1}{2}\int\limits_t^T 
(t_2-t) dt_2-
\frac{1}{2}
\int\limits_t^T \int\limits_t^{t_2} dt_1 dt_2=0.
$$

\vspace{5mm}

The equality (\ref{2023novem207}) is proved.

Let us prove (\ref{2023novem208}). Substitute $j_3=j_1,$ $j_4=j_2$ into
(\ref{2023novem215})

\begin{equation}
\label{2023novem224}
C_{j_2 j_1 j_2 j_1}+C_{j_1 j_2 j_1 j_2}=
C_{j_2}C_{j_1 j_2 j_1}-C_{j_1 j_2}C_{j_2 j_1}+
C_{j_2 j_1 j_2}C_{j_1}.
\end{equation}

\vspace{5mm}

Then
$$
\sum\limits_{j_1,j_2=0}^p \bigl(C_{j_2 j_1 j_2 j_1}+C_{j_1 j_2 j_1 j_2}\bigr)=
\sum\limits_{j_1,j_2=0}^p \bigl(C_{j_2}C_{j_1 j_2 j_1}+C_{j_2 j_1 j_2}C_{j_1}\bigr)-
$$

\vspace{2mm}
\begin{equation}
\label{2023novem225}
-\sum\limits_{j_1,j_2=0}^p C_{j_1 j_2}C_{j_2 j_1}.
\end{equation}

\vspace{5mm}

From (\ref{2023novem225}) we have

$$
2\sum\limits_{j_1,j_2=0}^p C_{j_2 j_1 j_2 j_1}=
2\sum\limits_{j_1,j_2=0}^p C_{j_1}C_{j_2 j_1 j_2}-
$$

\vspace{2mm}
$$
-\sum\limits_{j_1,j_2=0}^p \frac{1}{2} \left( \bigl(C_{j_1 j_2} + C_{j_2 j_1}\bigr)^2 
-\bigl(C_{j_1 j_2}\bigr)^2 - \bigl(C_{j_2 j_1}\bigr)^2\right)=
$$

\vspace{2mm}
$$
=
2\sum\limits_{j_1,j_2=0}^p C_{j_1}C_{j_2 j_1 j_2}
-\frac{1}{2}\sum\limits_{j_1,j_2=0}^p \bigl(C_{j_1 j_2} + C_{j_2 j_1}\bigr)^2 
+
$$

\vspace{2mm}
\begin{equation}
\label{2023novem226}
+\sum\limits_{j_1,j_2=0}^p \bigl(C_{j_2 j_1}\bigr)^2.
\end{equation}

\vspace{5mm}

Using Fubini's Theorem, we obtain

\begin{equation}
\label{2023novem227}
C_{j_1 j_2} + C_{j_2 j_1}=C_{j_1}C_{j_2}.
\end{equation}

\vspace{5mm}

Applying (\ref{2023novem226}), (\ref{2023novem227}), (\ref{2023novem220}),
Fubini's Theorem, Parseval's equality,
and (\ref{after1400}), we get

$$
\lim\limits_{p\to \infty}\sum\limits_{j_1,j_2=0}^p C_{j_2 j_1 j_2 j_1}=
\lim\limits_{p\to \infty}\sum\limits_{j_1,j_2=0}^p C_{j_1}C_{j_2 j_1 j_2}
-\frac{1}{4}\lim\limits_{p\to \infty}\sum\limits_{j_1,j_2=0}^p \bigl(C_{j_1 j_2} + C_{j_2 j_1}\bigr)^2 
+
$$

\vspace{2mm}
$$
+\frac{1}{2}\lim\limits_{p\to \infty}\sum\limits_{j_1,j_2=0}^p \bigl(C_{j_2 j_1}\bigr)^2=
$$

\vspace{2mm}
$$
=\sqrt{T-t}\sum\limits_{j_2=0}^{\infty} C_{j_2 0 j_2}
-\frac{1}{4}\sum\limits_{j_1,j_2=0}^{\infty} \bigl(C_{j_1}C_{j_2}\bigr)^2 
+\frac{1}{2}\sum\limits_{j_1,j_2=0}^{\infty} \bigl(C_{j_2 j_1}\bigr)^2=
$$

\vspace{2mm}
$$
=\sum\limits_{j_2=0}^{\infty} \int\limits_t^T \phi_{j_2}(t_3)\int\limits_t^{t_3}\int\limits_t^{t_2}
\phi_{j_2}(t_1)dt_1 dt_2 dt_3
-\frac{1}{4}(T-t)^2 
+\frac{1}{2}\int\limits_{[t,T]^2} \bigl({\bf 1}_{\{t_1<t_2\}}\bigr)^2 dt_1 dt_2=
$$

\vspace{2mm}
$$
=\sum\limits_{j_2=0}^{\infty} \int\limits_t^T \phi_{j_2}(t_3)\int\limits_t^{t_3}
\phi_{j_2}(t_1)\int\limits_{t_1}^{t_3} dt_2 dt_1 dt_3=
$$

\vspace{2mm}
$$
=\sum\limits_{j_2=0}^{\infty} \int\limits_t^T \phi_{j_2}(t_3)(t_3-t)\int\limits_t^{t_3}
\phi_{j_2}(t_1)dt_1 dt_3+
\sum\limits_{j_2=0}^{\infty} \int\limits_t^T \phi_{j_2}(t_3)\int\limits_t^{t_3}
\phi_{j_2}(t_1)(t-t_1)dt_1 dt_3=
$$

\vspace{2mm}
$$
=\frac{1}{2}\int\limits_t^T (t_3-t)dt_3+
\frac{1}{2}\int\limits_t^T(t-t_3)dt_3=0.
$$

\vspace{5mm}

The equality (\ref{2023novem208}) is proved. The equalities (\ref{2023novem200})--(\ref{2023novem208})
are proved. Theorem~44 is proved.

\vspace{5mm}

\section{Generalization of Theorems~39--41, 43, 44 to the Case When the
Conditions $\phi_0(x)=1/\sqrt{T-t}$ and
$\psi_l(\tau)\psi_{l-1}(\tau)\in L_2([t, T])$ $(l=2, 3,\ldots, k)$
are Omitted}

\vspace{5mm}

In this section, we will show that the 
conditions $\phi_0(x)=1/\sqrt{T-t}$ and
$\psi_l(\tau)\psi_{l-1}(\tau)\in L_2([t, T])$ $(l=2, 3,\ldots, k)$
in Theorems~39--41, 43, 44
can be omitted.

\vspace{2mm}   

{\bf Theorem~45}\ \cite{2018a}, \cite{arxiv-5}.\  {\it Suppose that
$\{\phi_j(x)\}_{j=0}^{\infty}$ is an arbitrary complete orthonormal system of 
func\-ti\-ons in the space $L_2([t,T]).$
Then$,$ for the iterated Stra\-to\-no\-vich stochastic integral
of third multiplicity 

\vspace{-1mm}
$$
{\int\limits_t^{*}}^T
{\int\limits_t^{*}}^{t_3}
{\int\limits_t^{*}}^{t_2}
d{\bf w}_{t_1}^{(i_1)}
d{\bf w}_{t_2}^{(i_2)}d{\bf w}_{t_3}^{(i_3)}\ \ \ (i_1,i_2,i_3=0,1,\ldots,m)
$$

\vspace{3mm}
\noindent
the following expansion 

\vspace{-1mm}
\begin{equation}
\label{2023novem1xyz}
{\int\limits_t^{*}}^T
{\int\limits_t^{*}}^{t_3}
{\int\limits_t^{*}}^{t_2}
d{\bf w}_{t_1}^{(i_1)}
d{\bf w}_{t_2}^{(i_2)}d{\bf w}_{t_3}^{(i_3)}=
\hbox{\vtop{\offinterlineskip\halign{
\hfil#\hfil\cr
{\rm l.i.m.}\cr
$\stackrel{}{{}_{p\to \infty}}$\cr
}} }\sum_{j_1,j_2,j_3=0}^{p}
C_{j_3 j_2 j_1}\zeta_{j_1}^{(i_1)}\zeta_{j_2}^{(i_2)}\zeta_{j_3}^{(i_3)}
\end{equation}

\vspace{3mm}
\noindent
that converges in the mean-square sense is valid, where 

\vspace{-1mm}
$$
C_{j_3 j_2 j_1}=\int\limits_t^T
\phi_{j_3}(t_3)\int\limits_t^{t_3}
\phi_{j_2}(t_2)
\int\limits_t^{t_2}
\phi_{j_1}(t_1)dt_1dt_2dt_3
$$

\vspace{1mm}
\noindent
and
$$
\zeta_{j}^{(i)}=
\int\limits_t^T \phi_{j}(\tau) d{\bf w}_{\tau}^{(i)}
$$ 

\vspace{2mm}
\noindent
are independent standard Gaussian random variables for various 
$i$ or $j$ {\rm (}in the case when $i\ne 0${\rm ),}
${\bf w}_{\tau}^{(i)}={\bf f}_{\tau}^{(i)}$ for
$i=1,\ldots,m$ and 
${\bf w}_{\tau}^{(0)}=\tau.$}

\vspace{2mm}

{\bf Proof.} Analyzing the proof of Theorems~40 and 43 (also see the derivation of (\ref{after906})
and (\ref{dydy11})),
we notice that Theorem~45 will be proved if we prove that

\vspace{-1mm}
\begin{equation}
\label{febr1}
\int\limits_t^T \int\limits_t^{t_3}
dt_2 d{\bf w}_{t_3}^{(i_3)}=
\hbox{\vtop{\offinterlineskip\halign{
\hfil#\hfil\cr
{\rm l.i.m.}\cr
$\stackrel{}{{}_{p\to \infty}}$\cr
}} }\sum_{j_3=0}^{p}\int\limits_t^T \phi_{j_3}(t_3)\int\limits_t^{t_3}
dt_2 dt_3\ \zeta_{j_3}^{(i_3)},
\end{equation}

\vspace{2mm}
\begin{equation}
\label{febr2}
\int\limits_t^T \int\limits_t^{t_2}
d{\bf w}_{t_1}^{(i_1)}dt_2=
\hbox{\vtop{\offinterlineskip\halign{
\hfil#\hfil\cr
{\rm l.i.m.}\cr
$\stackrel{}{{}_{p\to \infty}}$\cr
}} }\sum_{j_1=0}^{p}\int\limits_t^T \int\limits_t^{t_2}\phi_{j_1}(t_1)dt_1 dt_2\
\zeta_{j_1}^{(i_1)}.
\end{equation}

\vspace{3mm}

The equality (\ref{febr1}) immediately follows from (\ref{dsds11}) for $k=1$.
Let us prove (\ref{febr2}).
Using the theorem on replacement of the integration order in iterated 
Ito stochastic integrals (see Theorems 3.1, 3.3 in \cite{2018a}) or 
the Ito formula, (\ref{dsds11}) for $k=1$,
and Fubini's Theorem, we obtain~w.~p.~1 

\vspace{-1mm}
$$
\int\limits_t^T \int\limits_t^{t_2}
d{\bf w}_{t_1}^{(i_1)}dt_2=\int\limits_t^T \int\limits_{t_1}^{T}
dt_2 d{\bf w}_{t_1}^{(i_1)}=
\hbox{\vtop{\offinterlineskip\halign{
\hfil#\hfil\cr
{\rm l.i.m.}\cr
$\stackrel{}{{}_{p\to \infty}}$\cr
}} }\sum_{j_1=0}^{p}\int\limits_t^T \phi_{j_1}(t_1)\int\limits_{t_1}^{T}dt_2 dt_1\
\zeta_{j_1}^{(i_1)}=
$$

\vspace{2mm}
$$
=
\hbox{\vtop{\offinterlineskip\halign{
\hfil#\hfil\cr
{\rm l.i.m.}\cr
$\stackrel{}{{}_{p\to \infty}}$\cr
}} }\sum_{j_1=0}^{p}\int\limits_t^T \int\limits_{t}^{t_2}\phi_{j_1}(t_1)dt_1 dt_2\
\zeta_{j_1}^{(i_1)}.
$$

\vspace{4mm}

The equality (\ref{febr2}) is proved. Theorem~45 is proved.

Let us develop this approach and prove the following
generalization of Theorem~44.

\vspace{2mm}

{\bf Theorem~46}\ \cite{2018a}, \cite{arxiv-5}.\ {\it Suppose that
$\{\phi_j(x)\}_{j=0}^{\infty}$ is an arbitrary complete orthonormal system of 
func\-ti\-ons in the space $L_2([t,T]).$
Then$,$ for the iterated Stra\-to\-no\-vich stochastic integral
of fourth multiplicity 

\vspace{-1mm}
$$
J^{*}[\psi^{(4)}]_{T,t}=
{\int\limits_t^{*}}^T
{\int\limits_t^{*}}^{t_4}
{\int\limits_t^{*}}^{t_3}
{\int\limits_t^{*}}^{t_2}
d{\bf w}_{t_1}^{(i_1)}
d{\bf w}_{t_2}^{(i_2)}d{\bf w}_{t_3}^{(i_3)}d{\bf w}_{t_4}^{(i_4)}\ \ \ 
(i_1, i_2, i_3, i_4=0, 1,\ldots,m)
$$

\vspace{3mm}
\noindent
the following 
expansion 

\vspace{-1mm}
$$
J^{*}[\psi^{(4)}]_{T,t}=
\hbox{\vtop{\offinterlineskip\halign{
\hfil#\hfil\cr
{\rm l.i.m.}\cr
$\stackrel{}{{}_{p\to \infty}}$\cr
}} }
\sum\limits_{j_1, j_2, j_3, j_4=0}^{p}
C_{j_4 j_3 j_2 j_1}\zeta_{j_1}^{(i_1)}\zeta_{j_2}^{(i_2)}\zeta_{j_3}^{(i_3)}
\zeta_{j_4}^{(i_4)}
$$

\vspace{3mm}
\noindent
that converges in the mean-square sense is valid, where 

\vspace{-1mm}
$$
C_{j_4 j_3 j_2 j_1}=\int\limits_t^T
\phi_{j_4}(t_4)\int\limits_t^{t_4}
\phi_{j_3}(t_3)\int\limits_t^{t_3}
\phi_{j_2}(t_2)\int\limits_t^{t_2}
\phi_{j_1}(t_1)dt_1dt_2dt_3 dt_4
$$

\vspace{1mm}
\noindent
and
$$
\zeta_{j}^{(i)}=
\int\limits_t^T \phi_{j}(\tau) d{\bf w}_{\tau}^{(i)}
$$ 

\vspace{2mm}
\noindent
are independent standard Gaussian random variables for various 
$i$ or $j$ {\rm (}in the case when $i\ne 0${\rm ),}
${\bf w}_{\tau}^{(i)}={\bf f}_{\tau}^{(i)}$ for
$i=1,\ldots,m$ and 
${\bf w}_{\tau}^{(0)}=\tau.$}

\vspace{2mm}
                                                           
{\bf Proof.}\ Considering the proof of Theorems~40 and 44 (also see the derivation of (\ref{after906})
and (\ref{dydy11})),
we conclude that Theorem~46 will be proved if we prove that

\vspace{-1mm}
\begin{equation}
\label{febr5}
\int\limits_t^T \int\limits_t^{t_3} \int\limits_t^{t_2}
dt_1 d{\bf w}_{t_2}^{(i_2)}d{\bf w}_{t_3}^{(i_3)}=
\hbox{\vtop{\offinterlineskip\halign{
\hfil#\hfil\cr
{\rm l.i.m.}\cr
$\stackrel{}{{}_{p\to \infty}}$\cr
}} }\sum_{j_2,j_3=0}^{p}\int\limits_t^T \phi_{j_3}(t_3)\int\limits_t^{t_3} 
\phi_{j_2}(t_2)\int\limits_t^{t_2}
dt_1 dt_2 dt_3
J'[\phi_{j_2}\phi_{j_3}]_{T,t}^{(i_2 i_3)},
\end{equation}

\vspace{2mm}
\begin{equation}
\label{febr6}
\int\limits_t^T \int\limits_t^{t_3} \int\limits_t^{t_2}
d{\bf w}_{t_1}^{(i_1)} dt_2 d{\bf w}_{t_3}^{(i_3)}=
\hbox{\vtop{\offinterlineskip\halign{
\hfil#\hfil\cr
{\rm l.i.m.}\cr
$\stackrel{}{{}_{p\to \infty}}$\cr
}} }\sum_{j_1,j_3=0}^{p}\int\limits_t^T \phi_{j_3}(t_3)\int\limits_t^{t_3} 
\int\limits_t^{t_2}
\phi_{j_1}(t_1)dt_1 
dt_2 dt_3
J'[\phi_{j_1}\phi_{j_3}]_{T,t}^{(i_1 i_3)},
\end{equation}

\vspace{2mm}
\begin{equation}
\label{febr7}
\int\limits_t^T \int\limits_t^{t_3} \int\limits_t^{t_2}
d{\bf w}_{t_1}^{(i_1)}d{\bf w}_{t_2}^{(i_2)}dt_3=
\hbox{\vtop{\offinterlineskip\halign{
\hfil#\hfil\cr
{\rm l.i.m.}\cr
$\stackrel{}{{}_{p\to \infty}}$\cr
}} }\sum_{j_1,j_2=0}^{p}\int\limits_t^T \int\limits_t^{t_3} 
\phi_{j_2}(t_2)\int\limits_t^{t_2}
\phi_{j_1}(t_1)dt_1
dt_2 dt_3
J'[\phi_{j_1}\phi_{j_2}]_{T,t}^{(i_1 i_2)},
\end{equation}

\vspace{2mm}
\begin{equation}
\label{march10}
\lim\limits_{p\to\infty}
\sum\limits_{j_1, j_3=0}^{p}
C_{j_3 j_3 j_1 j_1}=\frac{1}{4} 
C_{j_3 j_3 j_1 j_1}\biggl|_{(j_{3} j_{3})\curvearrowright (\cdot)
(j_{1} j_{1})\curvearrowright (\cdot)}=\frac{1}{8}(T-t)^2,
\biggr.
\end{equation}

\vspace{2mm}
\begin{equation}
\label{march11}
\lim\limits_{p\to\infty}
\sum\limits_{j_1, j_2=0}^{p}
C_{j_2 j_1 j_2 j_1}=0
\biggr.,
\end{equation}

\vspace{1mm}
\begin{equation}
\label{march12}
\lim\limits_{p\to\infty}
\sum\limits_{j_1, j_3=0}^{p}
C_{j_1 j_3 j_3 j_1}=0
\biggr.
\end{equation}

\vspace{4mm}
\noindent 
where we use the same notations as in (\ref{dsds11}).

Moreover, for $k=4, r=2, g_1=1, g_2=2, g_3=3, g_4=4$ we can write 
(see the derivation of (\ref{after906}))

$$
\hbox{\vtop{\offinterlineskip\halign{
\hfil#\hfil\cr
{\rm l.i.m.}\cr
$\stackrel{}{{}_{p\to \infty}}$\cr
}} }
\sum\limits_{\stackrel{j_1,\ldots,j_q,\ldots,j_k=0}{{}_{q\ne g_1, g_2,\ldots, g_{2r-1}, g_{2r}}}}^p
\frac{1}{2^r}
C_{j_k \ldots j_1}\biggl|_{(j_{g_2} j_{g_1})\curvearrowright (\cdot)
\ldots (j_{g_{2r}} j_{g_{2r-1}})\curvearrowright (\cdot),
j_{g_{{}_{1}}}=~j_{g_{{}_{2}}},\ldots, j_{g_{{}_{2r-1}}}=~j_{g_{{}_{2r}}}}\biggr.
\times 
$$

\vspace{3mm}
$$
\times
\prod\limits_{s=1}^r
{\bf 1}_{\{i_{g_{{}_{2s-1}}}=~i_{g_{{}_{2s}}}\ne 0\}}
J'[\phi_{j_{q_1}}\ldots \phi_{j_{q_{k-2r}}}]_{T,t}^{(i_{q_1}\ldots i_{q_{k-2r}})}=
$$

\vspace{2mm}

$$
=\frac{1}{4}{\bf 1}_{\{i_{1}=i_{2}\ne 0\}}{\bf 1}_{\{i_{3}=i_{4}\ne 0\}}
C_{j_3 j_3 j_1 j_1}\biggl|_{(j_{3} j_{3})\curvearrowright (\cdot)
(j_{1} j_{1})\curvearrowright (\cdot)}\biggr.=
$$

\vspace{2mm}
$$
=\frac{1}{4}{\bf 1}_{\{i_{1}=i_{2}\ne 0\}}{\bf 1}_{\{i_{3}=i_{4}\ne 0\}}
\int\limits_t^T\int\limits_t^{t_2} dt_1 dt_2=
{\bf 1}_{\{i_{1}=i_{2}\ne 0\}}{\bf 1}_{\{i_{3}=i_{4}\ne 0\}}
\frac{(T-t)^2}{8},
$$

\vspace{4mm}
\noindent
where
$J'[\phi_{j_{q_1}}\ldots \phi_{j_{q_{k-2r}}}]_{T,t}^{(i_{q_1}\ldots i_{q_{k-2r}})}
\stackrel{\sf def}{=}1$
for $k=2r$.

The equality (\ref{febr5}) immediately follows from (\ref{dsds11})
for $k=2$.
Let us prove (\ref{febr7}).
Using the theorem on replacement of the integration order in iterated 
Ito stochastic integrals (see Theorems 3.1, 3.3 in \cite{2018a}) or 
the Ito formula, (\ref{dsds11}) for $k=2$, 
and Fubini's Theorem, we get~w.~p.~1

\vspace{1mm}
$$
\int\limits_t^T \int\limits_t^{t_3} \int\limits_t^{t_2}
d{\bf w}_{t_1}^{(i_1)}d{\bf w}_{t_2}^{(i_2)}dt_3=
\int\limits_t^T (T-t_2)\int\limits_t^{t_2} 
d{\bf w}_{t_1}^{(i_1)}d{\bf w}_{t_2}^{(i_2)}=
$$

\vspace{3mm}
$$
=\hbox{\vtop{\offinterlineskip\halign{
\hfil#\hfil\cr
{\rm l.i.m.}\cr
$\stackrel{}{{}_{p\to \infty}}$\cr
}} }\sum_{j_1,j_2=0}^{p}\int\limits_t^T (T-t_2)
\phi_{j_2}(t_2)\int\limits_t^{t_2}
\phi_{j_1}(t_1)dt_1 dt_2
J'[\phi_{j_1}\phi_{j_2}]_{T,t}^{(i_1 i_2)}=
$$

\vspace{3mm}
$$
=\hbox{\vtop{\offinterlineskip\halign{
\hfil#\hfil\cr
{\rm l.i.m.}\cr
$\stackrel{}{{}_{p\to \infty}}$\cr
}} }\sum_{j_1,j_2=0}^{p}\int\limits_t^T 
\phi_{j_1}(t_1)\int\limits_{t_1}^T
\phi_{j_2}(t_2)(T-t_2)dt_2 dt_1
J'[\phi_{j_1}\phi_{j_2}]_{T,t}^{(i_1 i_2)}=
$$

\vspace{3mm}
$$
=\hbox{\vtop{\offinterlineskip\halign{
\hfil#\hfil\cr
{\rm l.i.m.}\cr
$\stackrel{}{{}_{p\to \infty}}$\cr
}} }\sum_{j_1,j_2=0}^{p}\int\limits_t^T 
\phi_{j_1}(t_1)\int\limits_{t_1}^T
\phi_{j_2}(t_2)\int\limits_{t_2}^T dt_3 dt_2 dt_1
J'[\phi_{j_1}\phi_{j_2}]_{T,t}^{(i_1 i_2)}=
$$

\vspace{3mm}
$$
=
\hbox{\vtop{\offinterlineskip\halign{
\hfil#\hfil\cr
{\rm l.i.m.}\cr
$\stackrel{}{{}_{p\to \infty}}$\cr
}} }\sum_{j_1,j_2=0}^{p}\int\limits_t^T \int\limits_t^{t_3} 
\phi_{j_2}(t_2)\int\limits_t^{t_2}
\phi_{j_1}(t_1)dt_1
dt_2 dt_3
J'[\phi_{j_1}\phi_{j_2}]_{T,t}^{(i_1 i_2)}.
$$

\vspace{5mm}

The equality (\ref{febr7}) is proved. To prove (\ref{febr6}) we will use the above arguments
((\ref{febr12}) (see below) also directly follows from the Ito formula)

$$
\int\limits_t^T \int\limits_t^{t_3} \int\limits_t^{t_2}
d{\bf w}_{t_1}^{(i_1)} dt_2 d{\bf w}_{t_3}^{(i_3)}=\hbox{[by Theorems~3.1, 3.3 in \cite{2018a}]}=
\int\limits_t^T \int\limits_t^{t_3} 
d{\bf w}_{t_1}^{(i_1)}\int\limits_{t_1}^{t_3} dt_2 d{\bf w}_{t_3}^{(i_3)}=
$$

\vspace{2mm}
$$
=\int\limits_t^T \int\limits_t^{t_3} (t_3-t_1)
d{\bf w}_{t_1}^{(i_1)}d{\bf w}_{t_3}^{(i_3)}=
$$

\vspace{2mm}
\begin{equation}
\label{febr12}
=\int\limits_t^T (t_3-t)\int\limits_t^{t_3} 
d{\bf w}_{t_1}^{(i_1)}d{\bf w}_{t_3}^{(i_3)}-
\int\limits_t^T \int\limits_t^{t_3} (t_1-t)
d{\bf w}_{t_1}^{(i_1)}d{\bf w}_{t_3}^{(i_3)}=
\end{equation}

\vspace{2mm}
$$
=\hbox{\vtop{\offinterlineskip\halign{
\hfil#\hfil\cr
{\rm l.i.m.}\cr
$\stackrel{}{{}_{p\to \infty}}$\cr
}} }\sum_{j_1,j_3=0}^{p}\int\limits_t^T (t_3-t)\phi_{j_3}(t_3)\int\limits_t^{t_3} 
\phi_{j_1}(t_1)dt_1 
dt_3
J'[\phi_{j_1}\phi_{j_3}]_{T,t}^{(i_1 i_3)}-
$$

\vspace{2mm}
$$
-\hbox{\vtop{\offinterlineskip\halign{
\hfil#\hfil\cr
{\rm l.i.m.}\cr
$\stackrel{}{{}_{p\to \infty}}$\cr
}} }\sum_{j_1,j_3=0}^{p}\int\limits_t^T \phi_{j_3}(t_3)\int\limits_t^{t_3} 
(t_1-t)\phi_{j_1}(t_1)dt_1 
dt_3
J'[\phi_{j_1}\phi_{j_3}]_{T,t}^{(i_1 i_3)}=
$$

\vspace{2mm}
$$
=\hbox{\vtop{\offinterlineskip\halign{
\hfil#\hfil\cr
{\rm l.i.m.}\cr
$\stackrel{}{{}_{p\to \infty}}$\cr
}} }\sum_{j_1,j_3=0}^{p}\left(\int\limits_t^T (t_3-t)\phi_{j_3}(t_3)\int\limits_t^{t_3} 
\phi_{j_1}(t_1)dt_1 
dt_3-\right.
$$

\vspace{2mm}
$$
\left.-\int\limits_t^T \phi_{j_3}(t_3)\int\limits_t^{t_3} 
(t_1-t)\phi_{j_1}(t_1)dt_1 
dt_3\right)
J'[\phi_{j_1}\phi_{j_3}]_{T,t}^{(i_1 i_3)}=
$$

\vspace{2mm}
$$
=\hbox{\vtop{\offinterlineskip\halign{
\hfil#\hfil\cr
{\rm l.i.m.}\cr
$\stackrel{}{{}_{p\to \infty}}$\cr
}} }\sum_{j_1,j_3=0}^{p}\int\limits_t^T \phi_{j_3}(t_3)\int\limits_t^{t_3} (t_3-t+t-t_1)
\phi_{j_1}(t_1)dt_1 dt_3
J'[\phi_{j_1}\phi_{j_3}]_{T,t}^{(i_1 i_3)}=
$$

\vspace{2mm}
$$
=\hbox{\vtop{\offinterlineskip\halign{
\hfil#\hfil\cr
{\rm l.i.m.}\cr
$\stackrel{}{{}_{p\to \infty}}$\cr
}} }\sum_{j_1,j_3=0}^{p}\int\limits_t^T \phi_{j_3}(t_3)\int\limits_t^{t_3} 
\phi_{j_1}(t_1)\int\limits_{t_1}^{t_3}dt_2
dt_1 dt_3
J'[\phi_{j_1}\phi_{j_3}]_{T,t}^{(i_1 i_3)}=
$$

\vspace{2mm}
$$
=
\hbox{\vtop{\offinterlineskip\halign{
\hfil#\hfil\cr
{\rm l.i.m.}\cr
$\stackrel{}{{}_{p\to \infty}}$\cr
}} }\sum_{j_1,j_3=0}^{p}\int\limits_t^T \phi_{j_3}(t_3)\int\limits_t^{t_3} 
\int\limits_t^{t_2}
\phi_{j_1}(t_1)dt_1 
dt_2 dt_3
J'[\phi_{j_1}\phi_{j_3}]_{T,t}^{(i_1 i_3)}.
$$

\vspace{5mm}

The equality (\ref{febr6}) is proved. 
Let us prove (\ref{march10})--(\ref{march12}).
Using (\ref{2023novem219}), we obtain

\vspace{-1mm}
\begin{equation}
\label{march13}
\sum\limits_{j_1,j_3=0}^p C_{j_3 j_3 j_1 j_1}=
\sum\limits_{j_1,j_3=0}^p C_{j_3}C_{j_3 j_1 j_1}-\frac{1}{8}
\left(\sum\limits_{j_1=0}^p \bigl(C_{j_1}\bigr)^2\right)^2.
\end{equation}

\vspace{4mm}

Applying Parseval's equality, we have

\vspace{-1mm}
\begin{equation}
\label{march14}
\lim\limits_{p\to\infty}
\sum\limits_{j_1=0}^p \bigl(C_{j_1}\bigr)^2=\int\limits_t^T 1^2 d\tau=T-t.
\end{equation}

\vspace{4mm}

Combining (\ref{march13}) and (\ref{march14}), we get

\vspace{-1mm}
\begin{equation}
\label{march15}
\lim\limits_{p\to\infty}\sum\limits_{j_1,j_3=0}^p C_{j_3 j_3 j_1 j_1}=
\lim\limits_{p\to\infty}\sum\limits_{j_1,j_3=0}^p C_{j_3}C_{j_3 j_1 j_1}-\frac{(T-t)^2}{8}.
\end{equation}

\vspace{4mm}

Further, we have
$$
\lim\limits_{p\to\infty}\sum\limits_{j_1,j_3=0}^p C_{j_3}C_{j_3 j_1 j_1}=
$$

\vspace{2mm}
\begin{equation}
\label{march16}
=
\frac{1}{2}\lim\limits_{p\to\infty}
\sum\limits_{j_3=0}^{p}C_{j_3}C_{j_3 j_1 j_1}\biggl|_{(j_{1} j_{1})\curvearrowright (\cdot)}-
\lim\limits_{p\to\infty}
\sum\limits_{j_3=0}^{p}C_{j_3}\left(\frac{1}{2}
C_{j_3 j_1 j_1}\biggl|_{(j_{1} j_{1})\curvearrowright (\cdot)}-
\sum\limits_{j_1=0}^p C_{j_3 j_1 j_1}\right).
\end{equation}

\vspace{5mm}

Applying the generalized Parseval equality, we obtain

\vspace{-1mm}
$$
\lim\limits_{p\to\infty}
\sum\limits_{j_3=0}^{p}C_{j_3}C_{j_3 j_1 j_1}\biggl|_{(j_{1} j_{1})\curvearrowright (\cdot)}=
\lim\limits_{p\to\infty}
\sum\limits_{j_3=0}^{p}\int\limits_t^T \phi_{j_3}(\tau)d\tau
\int\limits_t^T \phi_{j_3}(\tau)\int\limits_t^{\tau}d\theta d\tau
=
$$

\vspace{2mm}
\begin{equation}
\label{march17}
=
\int\limits_t^T 1\cdot \int\limits_t^{\tau}d\theta d\tau=
\frac{(T-t)^2}{2}.
\end{equation}

\vspace{4mm}

From (\ref{march16}) and (\ref{march17}) we have

\vspace{-1mm}
$$
\lim\limits_{p\to\infty}\sum\limits_{j_1,j_3=0}^p C_{j_3}C_{j_3 j_1 j_1}=
$$

\vspace{2mm}
\begin{equation}
\label{march18}
=
\frac{(T-t)^2}{4}-
\lim\limits_{p\to\infty}
\sum\limits_{j_3=0}^{p}C_{j_3}\left(\frac{1}{2}
C_{j_3 j_1 j_1}\biggl|_{(j_{1} j_{1})\curvearrowright (\cdot)}-
\sum\limits_{j_1=0}^p C_{j_3 j_1 j_1}\right).
\end{equation}

\vspace{5mm}

Combining (\ref{march15}) and (\ref{march18}), we obtain

\vspace{-1mm}
\begin{equation}
\label{march19}
\lim\limits_{p\to\infty}\sum\limits_{j_1,j_3=0}^p C_{j_3 j_3 j_1 j_1}=
\frac{(T-t)^2}{8}-
\lim\limits_{p\to\infty}
\sum\limits_{j_3=0}^{p}C_{j_3}\left(\frac{1}{2}
C_{j_3 j_1 j_1}\biggl|_{(j_{1} j_{1})\curvearrowright (\cdot)}-
\sum\limits_{j_1=0}^p C_{j_3 j_1 j_1}\right).
\end{equation}

\vspace{5mm}

Due to the inequality of Cauchy--Bunyakovsky
and (\ref{2023novem2}), (\ref{march14}), we get

\vspace{-1mm}
$$
\lim\limits_{p\to\infty}
\left(\sum\limits_{j_3=0}^{p}C_{j_3}\left(\frac{1}{2}
C_{j_3 j_1 j_1}\biggl|_{(j_{1} j_{1})\curvearrowright (\cdot)}-
\sum\limits_{j_1=0}^p C_{j_3 j_1 j_1}\right)\right)^2\le
$$

\vspace{2mm}
$$
\le \lim\limits_{p\to\infty}
\sum\limits_{j_3=0}^{p}\left(C_{j_3}\right)^2\
\sum\limits_{j_3=0}^{p}
\left(\frac{1}{2}
C_{j_3 j_1 j_1}\biggl|_{(j_{1} j_{1})\curvearrowright (\cdot)}-
\sum\limits_{j_1=0}^p C_{j_3 j_1 j_1}\right)^2\le
$$

\vspace{2mm}
$$
\le \lim\limits_{p\to\infty}
\sum\limits_{j_3=0}^{\infty}\left(C_{j_3}\right)^2\
\sum\limits_{j_3=0}^{p}
\left(\frac{1}{2}
C_{j_3 j_1 j_1}\biggl|_{(j_{1} j_{1})\curvearrowright (\cdot)}-
\sum\limits_{j_1=0}^p C_{j_3 j_1 j_1}\right)^2=
$$

\vspace{2mm}
\begin{equation}
\label{march20}
=(T-t)\lim\limits_{p\to\infty}
\sum\limits_{j_3=0}^{p}
\left(\frac{1}{2}
C_{j_3 j_1 j_1}\biggl|_{(j_{1} j_{1})\curvearrowright (\cdot)}-
\sum\limits_{j_1=0}^p C_{j_3 j_1 j_1}\right)^2=0.
\end{equation}
         
\vspace{5mm}

Taking into account (\ref{march19}) and (\ref{march20}), we obtain (\ref{march10}).
It is not difficult to see that by analogy with (\ref{march10}) we get

\vspace{-1mm}
\begin{equation}
\label{march20a}
\lim\limits_{p\to\infty}
\sum\limits_{j_1, j_3=0}^{p}
C_{j_3 j_3 j_1 j_1}(s)=\frac{1}{8}(s-t)^2,
\biggr.
\end{equation}

\vspace{3mm}
\noindent
where $s\in (t, T]$ and

\vspace{-2mm}
\begin{equation}
\label{march21a}
C_{j_4 j_3 j_2 j_1}(s)=\int\limits_t^s
\phi_{j_4}(t_4)\int\limits_t^{t_4}
\phi_{j_3}(t_3)\int\limits_t^{t_3}
\phi_{j_2}(t_2)\int\limits_t^{t_2}
\phi_{j_1}(t_1)dt_1dt_2dt_3 dt_4.
\end{equation}

\vspace{4mm}

Let us prove (\ref{march11}). Using (\ref{2023novem225}), we have

\vspace{-1mm}
\begin{equation}
\label{march21}
\sum\limits_{j_1,j_2=0}^p C_{j_2 j_1 j_2 j_1}=
\sum\limits_{j_1,j_2=0}^p C_{j_2}C_{j_1 j_2 j_1}
-\frac{1}{2}\sum\limits_{j_1,j_2=0}^p C_{j_1 j_2}C_{j_2 j_1}.
\end{equation}

\vspace{4mm}

Fubini's Theorem and the generalized Parseval equality give

\vspace{-1mm}
$$
\lim_{p\to\infty}\sum\limits_{j_1,j_2=0}^p C_{j_1 j_2}C_{j_2 j_1}=
$$

\vspace{2mm}
$$
=
\lim_{p\to\infty}
\sum\limits_{j_1,j_2=0}^p \int\limits_t^T \phi_{j_2}(t_2)\int\limits_{t_2}^{T}\phi_{j_1}(t_1)dt_1 dt_2
\int\limits_t^T \phi_{j_2}(t_2)\int\limits_t^{t_2}\phi_{j_1}(t_1)dt_1 dt_2=
$$

\vspace{2mm}
$$
=\lim_{p\to\infty}\sum\limits_{j_1,j_2=0}^p 
\int\limits_{[t,T]^2}{\bf 1}_{\{t_2<t_1\}}\phi_{j_1}(t_1)\phi_{j_2}(t_2)dt_1 dt_2
\int\limits_{[t,T]^2} {\bf 1}_{\{t_1<t_2\}} \phi_{j_1}(t_1)\phi_{j_2}(t_2)dt_1 dt_2=
$$

\vspace{2mm}
\begin{equation}
\label{march22}
=\int\limits_{[t,T]^2}{\bf 1}_{\{t_2<t_1\}}{\bf 1}_{\{t_1<t_2\}}dt_1 dt_2=0.
\end{equation}

\vspace{5mm}

The equalities (\ref{march21}) and (\ref{march22}) imply the relation

\vspace{-1mm}
\begin{equation}
\label{march23}
\lim_{p\to\infty}\sum\limits_{j_1,j_2=0}^p C_{j_2 j_1 j_2 j_1}=
\lim_{p\to\infty}\sum\limits_{j_1,j_2=0}^p C_{j_2}C_{j_1 j_2 j_1}.
\end{equation}

\vspace{4mm}

Further, we have (see the derivation of (\ref{march20}))

\vspace{-1mm}
$$
\lim_{p\to\infty}\left(\sum\limits_{j_2=0}^p C_{j_2} \sum\limits_{j_1=0}^p  C_{j_1 j_2 j_1}\right)^2\le
\lim_{p\to\infty}\sum\limits_{j_2=0}^p \left(C_{j_2}\right)^2 
\sum\limits_{j_2=0}^p\left(\sum\limits_{j_1=0}^p  C_{j_1 j_2 j_1}\right)^2\le
$$

\vspace{2mm}
\begin{equation}
\label{march24}
\le \lim_{p\to\infty}\sum\limits_{j_2=0}^{\infty} \left(C_{j_2}\right)^2 
\sum\limits_{j_2=0}^p\left(\sum\limits_{j_1=0}^p  C_{j_1 j_2 j_1}\right)^2=
(T-t)\lim_{p\to\infty}\sum\limits_{j_2=0}^p\left(\sum\limits_{j_1=0}^p  C_{j_1 j_2 j_1}\right)^2=0,
\end{equation}

\vspace{5mm}
\noindent
where (\ref{march24}) follows from (\ref{2023novem4}).

The relations (\ref{march23}) and (\ref{march24}) complete the proof of (\ref{march11}).
By analogy with the above reasoning, we obviously get
\begin{equation}
\label{march25}
\lim\limits_{p\to\infty}
\sum\limits_{j_1, j_2=0}^{p}
C_{j_2 j_1 j_2 j_1}(s)=0
\biggr.,
\end{equation}

\vspace{4mm}
\noindent
where $s\in (t, T]$ and $C_{j_2 j_1 j_2 j_1}(s)$ is defined by (\ref{march21a}).

Let us prove (\ref{march12}). Using (\ref{2023novem223}), we obtain

\vspace{-1mm}
\begin{equation}
\label{march26}
\sum\limits_{j_1,j_3=0}^p C_{j_1 j_3 j_3 j_1}=
\sum\limits_{j_1,j_3=0}^p C_{j_1}C_{j_3 j_3 j_1}-\frac{1}{2}
\sum\limits_{j_1,j_3=0}^p \bigl(C_{j_3 j_1}\bigr)^2.
\end{equation}

\vspace{4mm}

Parseval's equality gives

\vspace{-1mm}
$$
\lim\limits_{p\to\infty}\sum\limits_{j_1,j_3=0}^p \bigl(C_{j_3 j_1}\bigr)^2=
\lim\limits_{p\to\infty}\sum\limits_{j_1,j_3=0}^p
\left(~\int\limits_{[t,T]^2}{\bf 1}_{\{t_1<t_2\}}
\phi_{j_1}(t_1)\phi_{j_3}(t_2)dt_1 dt_2\right)^2=
$$

\vspace{2mm}
\begin{equation}
\label{march27}
=\int\limits_{[t,T]^2}\left({\bf 1}_{\{t_1<t_2\}}\right)^2 dt_1 dt_2=
\frac{(T-t)^2}{2}.
\end{equation}

\vspace{5mm}

Combining (\ref{march26}) and (\ref{march27}), we have

\vspace{-1mm}
\begin{equation}
\label{march28}
\lim\limits_{p\to\infty}\sum\limits_{j_1,j_3=0}^p C_{j_1 j_3 j_3 j_1}=
\lim\limits_{p\to\infty}\sum\limits_{j_1,j_3=0}^p C_{j_1}C_{j_3 j_3 j_1}-\frac{(T-t)^2}{4}.
\end{equation}

\vspace{4mm}

Further, we have
$$
\lim\limits_{p\to\infty}\sum\limits_{j_1,j_3=0}^p C_{j_1}C_{j_3 j_3 j_1}=
$$

\vspace{2mm}
\begin{equation}
\label{march29}
=
\frac{1}{2}\lim\limits_{p\to\infty}
\sum\limits_{j_1=0}^{p}C_{j_1}C_{j_3 j_3 j_1}\biggl|_{(j_{3} j_{3})\curvearrowright (\cdot)}-
\lim\limits_{p\to\infty}
\sum\limits_{j_1=0}^{p}C_{j_1}\left(\frac{1}{2}
C_{j_3 j_3 j_1}\biggl|_{(j_{3} j_{3})\curvearrowright (\cdot)}-
\sum\limits_{j_3=0}^p C_{j_3 j_3 j_1}\right).
\end{equation}

\vspace{5mm}

Applying Fubini's Theorem and the generalized Parseval equality, we obtain

$$
\lim\limits_{p\to\infty}
\sum\limits_{j_1=0}^{p}C_{j_1}C_{j_3 j_3 j_1}\biggl|_{(j_{3} j_{3})\curvearrowright (\cdot)}=
\lim\limits_{p\to\infty}
\sum\limits_{j_1=0}^{p}\int\limits_t^T \phi_{j_1}(\tau)d\tau
\int\limits_t^T \int\limits_t^{t_2}\phi_{j_1}(\tau)d\tau dt_2
=
$$

\vspace{2mm}
\begin{equation}
\label{march30}
=\lim\limits_{p\to\infty}
\sum\limits_{j_1=0}^{p}\int\limits_t^T \phi_{j_1}(\tau)d\tau
\int\limits_t^T \phi_{j_1}(\tau) \int\limits_{\tau}^{T} dt_2 d\tau
=
\int\limits_t^T 1\cdot \int\limits_{\tau}^T d\theta d\tau=
\frac{(T-t)^2}{2}.
\end{equation}

\vspace{5mm}

From (\ref{march29}) and (\ref{march30}) we have

\vspace{-1mm}
$$
\lim\limits_{p\to\infty}\sum\limits_{j_1,j_3=0}^p C_{j_1}C_{j_3 j_3 j_1}=
$$

\vspace{2mm}
\begin{equation}
\label{march31}
=
\frac{(T-t)^2}{4}-
\lim\limits_{p\to\infty}
\sum\limits_{j_1=0}^{p}C_{j_1}\left(\frac{1}{2}
C_{j_3 j_3 j_1}\biggl|_{(j_{3} j_{3})\curvearrowright (\cdot)}-
\sum\limits_{j_3=0}^p C_{j_3 j_3 j_1}\right).
\end{equation}

\vspace{5mm}

Combining (\ref{march28}) and (\ref{march31}), we obtain

\vspace{-1mm}
\begin{equation}
\label{march32}
\lim\limits_{p\to\infty}\sum\limits_{j_1,j_3=0}^p C_{j_1 j_3 j_3 j_1}=
-
\lim\limits_{p\to\infty}
\sum\limits_{j_1=0}^{p}C_{j_1}\left(\frac{1}{2}
C_{j_3 j_3 j_1}\biggl|_{(j_{3} j_{3})\curvearrowright (\cdot)}-
\sum\limits_{j_3=0}^p C_{j_3 j_3 j_1}\right).
\end{equation}

\vspace{5mm}

Due to the inequality of Cauchy--Bunyakovsky
and (\ref{2023novem3}), (\ref{march14}), we get

\vspace{-1mm}
$$
\lim\limits_{p\to\infty}
\left(\sum\limits_{j_1=0}^{p}C_{j_1}\left(\frac{1}{2}
C_{j_3 j_3 j_1}\biggl|_{(j_{3} j_{3})\curvearrowright (\cdot)}-
\sum\limits_{j_3=0}^p C_{j_3 j_3 j_1}\right)\right)^2\le
$$

\vspace{2mm}
$$
\le \lim\limits_{p\to\infty}
\sum\limits_{j_1=0}^{p}\left(C_{j_1}\right)^2\
\sum\limits_{j_1=0}^{p}
\left(\frac{1}{2}
C_{j_3 j_3 j_1}\biggl|_{(j_{3} j_{3})\curvearrowright (\cdot)}-
\sum\limits_{j_3=0}^p C_{j_3 j_3 j_1}\right)^2\le
$$

\vspace{2mm}
$$
\le \lim\limits_{p\to\infty}
\sum\limits_{j_1=0}^{\infty}\left(C_{j_1}\right)^2\
\sum\limits_{j_1=0}^{p}
\left(\frac{1}{2}
C_{j_3 j_3 j_1}\biggl|_{(j_{3} j_{3})\curvearrowright (\cdot)}-
\sum\limits_{j_3=0}^p C_{j_3 j_3 j_1}\right)^2=
$$

\vspace{2mm}
\begin{equation}
\label{march33}
=(T-t)\lim\limits_{p\to\infty}
\sum\limits_{j_1=0}^{p}
\left(\frac{1}{2}
C_{j_3 j_3 j_1}\biggl|_{(j_{3} j_{3})\curvearrowright (\cdot)}-
\sum\limits_{j_3=0}^p C_{j_3 j_3 j_1}\right)^2=0.
\end{equation}

\vspace{5mm}         

The relations (\ref{march32}) and (\ref{march33}) complete the proof of (\ref{march12}).
By analogy with the above reasoning, we obviously have
\begin{equation}
\label{march34}
\lim\limits_{p\to\infty}
\sum\limits_{j_1, j_3=0}^{p}
C_{j_1 j_3 j_3 j_1}(s)=0
\biggr.,
\end{equation}

\vspace{4mm}
\noindent
where $s\in (t, T]$ and $C_{j_1 j_3 j_3 j_1}(s)$ is defined by (\ref{march21a}).

The equalities (\ref{febr5})--(\ref{march12}) are proved. Theorem~46 is proved.

Note that the equalities (\ref{march25}) and (\ref{march34}) can be proved by another way.
Using Fubini's Theorem, we obtain
\begin{equation}
\label{march35}
C_{j_2 j_1 j_2 j_1}(s)=\frac{1}{2}\left(C_{j_2 j_1}(s)\right)^2-2C_{j_2 j_2 j_1 j_1}(s),
\end{equation}

\vspace{1mm}
\begin{equation}
\label{march36}
\sum\limits_{(j_1,j_2,j_3,j_4)}C_{j_4 j_3 j_2 j_1}(s)=C_{j_1}(s)C_{j_2}(s)C_{j_3}(s)C_{j_4}(s),
\end{equation}

\vspace{3mm}
\noindent
where $s\in (t, T],$
$$
\sum\limits_{(j_1,j_2,j_3,j_4)}
$$

\vspace{3mm}
\noindent
means the sum with respect to all
possible permutations
$(j_1,j_2,j_3,j_4)$ and

\vspace{-1mm}
$$
C_{j_k \ldots j_1}(s)=\int\limits_t^s
\phi_{j_k}(t_k)\ldots
\int\limits_t^{t_2}
\phi_{j_1}(t_1)dt_1\ldots dt_k \ \ \ (k=1,\ldots,4).
$$

\vspace{2mm}

Taking into account (\ref{march20a}), (\ref{march27})
(for $s$ instead of $T$), (\ref{march35}), we get

$$
\lim\limits_{p\to\infty}\sum\limits_{j_1,j_2=0}^{p}
C_{j_2 j_1 j_2 j_1}(s)=\frac{1}{2}
\lim\limits_{p\to\infty}\sum\limits_{j_1,j_2=0}^{p}\left(C_{j_2 j_1}(s)\right)^2-2
\lim\limits_{p\to\infty}\sum\limits_{j_1,j_2=0}^{p}C_{j_2 j_2 j_1 j_1}(s)=
$$

\vspace{2mm}
$$
=\frac{1}{2}\cdot \frac{(s-t)^2}{2}-2\cdot \frac{(s-t)^2}{8}=0.
$$

\vspace{4mm}

The equality (\ref{march25}) is proved.
Let us substitute $j_2=j_1$ and $j_4=j_3$ into (\ref{march36}). Then we obtain

$$
4\biggl(C_{j_3 j_3 j_1 j_1}(s)+C_{j_1 j_1 j_3 j_3}(s)+C_{j_3 j_1 j_1 j_3}(s)
+C_{j_1 j_3 j_3 j_1}(s)+\biggr.
$$

\vspace{-1mm}
\begin{equation}
\label{march40}
\biggl.+C_{j_3 j_1 j_3 j_1}(s)+C_{j_1 j_3 j_1 j_3}(s)\biggr)=
\left(C_{j_1}(s)\right)^2\left(C_{j_3}(s)\right)^2.
\end{equation}

\vspace{4mm}

The equality (\ref{march40}) implies that

\vspace{-1mm}
\begin{equation}
\label{march41}
8\sum\limits_{j_1,j_3=0}^{p}\biggl(C_{j_3 j_3 j_1 j_1}(s)+
C_{j_1 j_3 j_3 j_1}(s)+C_{j_3 j_1 j_3 j_1}(s)\biggr)=
\sum\limits_{j_1=0}^{p}\left(C_{j_1}(s)\right)^2\sum\limits_{j_3=0}^{p}\left(C_{j_3}(s)\right)^2.
\end{equation}

\vspace{4mm}

Passing to the limit $\lim\limits_{p\to\infty}$ in (\ref{march41})
and taking into account (\ref{march14}) (for $s$ instead of $T$), 
(\ref{march20a}), (\ref{march25}), we get
$$
8 \Biggl(\frac{(s-t)^2}{8} + \lim\limits_{p\to\infty}
\sum\limits_{j_1,j_3=0}^{p}C_{j_1 j_3 j_3 j_1}(s)+0\Biggr)=
(s-t)^2.
$$

\vspace{4mm}

The equality (\ref{march34}) is  proved.

Consider the following generalization of Theorem~40.

\vspace{2mm}

{\bf Theorem~47}\ \cite{2018a}, \cite{arxiv-5}.\ {\it Assume that
the complete orthonormal system $\{\phi_j(x)\}_{j=0}^{\infty}$
in the space $L_2([t, T])$ and
$\psi_1(\tau),\ldots, \psi_k(\tau)\in L_2([t, T])$
are such that 

\vspace{1mm}
$$
\lim\limits_{p_1,\ldots,p_k\to\infty}~
\sum\limits_{j_1=0}^{p_1}\ldots \sum\limits_{j_q=0}^{p_q}\ldots \sum\limits_{j_k=0}^{p_k}~
\biggl|_{q\ne g_1, g_2, \ldots, g_{2r-1},g_{2r}}\times
$$

\vspace{5mm}
$$
\times
\Biggl(~\sum\limits_{j_{g_1}=0}^{\min\{p_{g_1}, p_{g_2}\}} \sum\limits_{j_{g_3}=0}^{\min\{p_{g_3}, p_{g_4}\}}\ldots \Biggr.
\sum\limits_{j_{g_{2r-1}}=0}^{\min\{p_{g_{2r-1}}, p_{g_{2r}}\}}
C_{j_k\ldots j_1}\biggl|_{j_{g_1}=j_{g_2},\ldots, j_{g_{2r-1}}=j_{g_{2r}}}-
$$

\vspace{2mm}
\begin{equation}
\label{june100}
\Biggl.-\frac{1}{2^r} \prod\limits_{l=1}^r {\bf 1}_{\{g_{2l}=g_{2l-1}+1\}}
C_{j_k \ldots j_1}\biggl|_{(j_{g_2} j_{g_1})\curvearrowright (\cdot)
\ldots (j_{g_{2r}} j_{g_{2r-1}})\curvearrowright (\cdot),
j_{g_{{}_{1}}}=~j_{g_{{}_{2}}},\ldots, j_{g_{{}_{2r-1}}}=~j_{g_{{}_{2r}}}
}\biggr.\Biggr)^2=0
\end{equation}

\vspace{6mm}
\noindent
for all $r=1, 2,\ldots,[k/2].$ 
Then$,$ for the sum $\bar J^{*}[\psi^{(k)}]_{T,t}^{(i_1\ldots i_k)}$
of iterated Ito stochastic integrals 
defined by {\rm (\ref{dsds9})}
the following 
expansion 

$$
\bar J^{*}[\psi^{(k)}]_{T,t}^{(i_1\ldots i_k)}=
\hbox{\vtop{\offinterlineskip\halign{
\hfil#\hfil\cr
{\rm l.i.m.}\cr
$\stackrel{}{{}_{p_1,\ldots,p_k\to \infty}}$\cr
}} }
\sum_{j_1=0}^{p_1}\ldots\sum_{j_k=0}^{p_k}
C_{j_k \ldots j_1}\prod\limits_{l=1}^k \zeta_{j_l}^{(i_l)}
$$

\vspace{3mm}
\noindent
that converges in the mean-square sense is valid, where 

\begin{equation}
\label{july15030}
C_{j_k \ldots j_1}=\int\limits_t^T\psi_k(t_k)\phi_{j_k}(t_k)\ldots
\int\limits_t^{t_2}
\psi_1(t_1)\phi_{j_1}(t_1)
dt_1\ldots dt_k
\end{equation}

\vspace{3mm}
\noindent
is the Fourier coefficient, 
${\rm l.i.m.}$ is a limit in the mean-square sense,
$i_1, \ldots, i_k=0, 1,\ldots,m,$

\vspace{-1mm}
$$
\zeta_{j}^{(i)}=
\int\limits_t^T \phi_{j}(\tau) d{\bf w}_{\tau}^{(i)}
$$ 

\vspace{2mm}
\noindent
are independent standard Gaussian random variables for various 
$i$ or $j$ {\rm (}in the case when $i\ne 0${\rm )},
${\bf w}_{\tau}^{(0)}=\tau.$}

\vspace{2mm}

{\bf Proof.}\ To prove Theorem~47, we need to prove that under
the conditions of Theorem~47 the following equality 

\vspace{1mm}
$$
\hbox{\vtop{\offinterlineskip\halign{
\hfil#\hfil\cr
{\rm l.i.m.}\cr
$\stackrel{}{{}_{p\to \infty}}$\cr
}} }
\sum\limits_{\stackrel{j_1,\ldots,j_q,\ldots,j_k=0}{{}_{q\ne g_1, g_2,\ldots, g_{2r-1}, g_{2r}}}}^p
\frac{1}{2^r}
C_{j_k \ldots j_1}\biggl|_{(j_{g_2} j_{g_1})\curvearrowright (\cdot)
\ldots (j_{g_{2r}} j_{g_{2r-1}})\curvearrowright (\cdot),
j_{g_{{}_{1}}}=~j_{g_{{}_{2}}},\ldots, j_{g_{{}_{2r-1}}}=~j_{g_{{}_{2r}}}}\biggr.
\times 
$$

\vspace{4mm}
$$
\times
\prod\limits_{s=1}^r
{\bf 1}_{\{i_{g_{{}_{2s-1}}}=~i_{g_{{}_{2s}}}\ne 0\}}
J'[\phi_{j_{q_1}}\ldots \phi_{j_{q_{k-2r}}}]_{T,t}^{(i_{q_1}\ldots i_{q_{k-2r}})}=
$$

\vspace{2mm}
\begin{equation}
\label{febr14}
=\frac{1}{2^r}
J[\psi^{(k)}]_{T,t}^{s_r,\ldots,s_1}
\end{equation}

\vspace{2mm}
\noindent
holds w.~p.~1, where $g_{2}=g_{1}+1,\ldots, g_{2r}=g_{2r-1}+1,$
$g_{2i-1}\stackrel{\sf def}{=}s_i,$\ $i=1,2,\ldots,r,$\
$r=1,2,\ldots,\left[k/2\right],$ 
$(s_r,\ldots,s_1)\in {\rm A}_{k,r},$ $J[\psi^{(k)}]_{T,t}^{s_r,\ldots,s_1}$ is
defined by (\ref{30.1}) and ${\rm A}_{k,r}$ is defined by (\ref{30.5550001});
also we put $p_1=\ldots=p_k=p$ in (\ref{febr14}) to simplify the notation;
another notations in (\ref{febr14}) are the same as in Sect.~7.

Using the Ito formula, we obtain w.~p.~1

\vspace{1mm}
$$
\int\limits_t^T \psi_k(t_k)\ldots \int\limits_t^{t_{l+2}}
\psi_{l+1}(t_{l+1}) \int\limits_t^{t_{l+1}}\psi_l(t_{l-1})\psi_{l-1}(t_{l-1})
\int\limits_t^{t_{l-1}}\psi_{l-2}(t_{l-2})\ldots
$$

\vspace{2mm}
$$
\ldots \int\limits_t^{t_2} \psi_1(t_1)d{\bf w}_{t_1}^{(i_1)}\ldots
d{\bf w}_{t_{l-2}}^{(i_{l-2})}dt_{l-1} d{\bf w}_{t_{l+1}}^{(i_{l+1})}\ldots
d{\bf w}_{t_k}^{(i_k)}=
$$

\vspace{2mm}
$$
=\int\limits_t^T \psi_k(t_k)\ldots \int\limits_t^{t_{l+2}}
\psi_{l+1}(t_{l+1}) \left(\int\limits_t^{t_{l+1}}\psi_l(t_{l-1})\psi_{l-1}(t_{l-1})dt_{l-1}\right)
\int\limits_t^{t_{l+1}}\psi_{l-2}(t_{l-2})\ldots
$$

\vspace{2mm}
$$
\ldots \int\limits_t^{t_2} \psi_1(t_1)d{\bf w}_{t_1}^{(i_1)}\ldots
d{\bf w}_{t_{l-2}}^{(i_{l-2})}d{\bf w}_{t_{l+1}}^{(i_{l+1})}\ldots
d{\bf w}_{t_k}^{(i_k)}-
$$

\vspace{2mm}
$$
-\int\limits_t^T \psi_k(t_k)\ldots \int\limits_t^{t_{l+2}}
\psi_{l+1}(t_{l+1}) 
\int\limits_t^{t_{l+1}}\psi_{l-2}(t_{l-2})
\left(\int\limits_t^{t_{l-2}}\psi_l(t_{l-1})\psi_{l-1}(t_{l-1})dt_{l-1}\right)\times
$$

\vspace{2mm}
\begin{equation}
\label{febr15}
\times\int\limits_t^{t_{l-2}}\psi_{l-3}(t_{l-3})
\ldots
\int\limits_t^{t_2} \psi_1(t_1)d{\bf w}_{t_1}^{(i_1)}\ldots
d{\bf w}_{t_{l-3}}^{(i_{l-3})}d{\bf w}_{t_{l-2}}^{(i_{l-2})}d{\bf w}_{t_{l+1}}^{(i_{l+1})}\ldots
d{\bf w}_{t_k}^{(i_k)},
\end{equation}

\vspace{3mm}
\noindent
where $l\ge 3.$ Note that the formula (\ref{febr15})
will change in an obvious way for the case $t_{l+1}=T.$
We will also assume that the transformation (\ref{febr15})
is not carried out for $l=2$ since the integral

\vspace{-1mm}
$$
\int\limits_t^{t_{3}}\psi_2(t_{1})\psi_{1}(t_{1})dt_{1}
$$

\vspace{2mm}
\noindent
is an internal integral on the left-hand side of 
(\ref{febr15}) for this case.

It is important to note that the transformation (\ref{febr15})
fully complies with the classical rules for replacing
the order of integration (Fubini's Theorem) if 
we replace
all differentials of the form
$d{\bf w}^{(i_j)}_{t_j}$ with $dt_j$
in (\ref{febr15}).

Indeed, formally changing the order of integration 
on the left-hand side of (\ref{febr15}) according 
to the classical rules, we have

\vspace{1mm}
\begin{equation}
\label{febr20}
\int\limits_t^T \psi_k(t_k)\ldots \int\limits_t^{t_{l+2}}
\psi_{l+1}(t_{l+1}) \int\limits_t^{t_{l+1}}\psi_l(t_{l-1})\psi_{l-1}(t_{l-1})
\int\limits_t^{t_{l-1}}\psi_{l-2}(t_{l-2})\ldots
\end{equation}

\vspace{2mm}
$$
\ldots \int\limits_t^{t_2} \psi_1(t_1)d{\bf w}_{t_1}^{(i_1)}\ldots
d{\bf w}_{t_{l-2}}^{(i_{l-2})}dt_{l-1} d{\bf w}_{t_{l+1}}^{(i_{l+1})}\ldots
d{\bf w}_{t_k}^{(i_k)}=
$$

\vspace{2mm}
$$
=\int\limits_t^T \psi_k(t_k)\ldots \int\limits_t^{t_{l+2}}
\psi_{l+1}(t_{l+1})\left(\int\limits_t^{t_{l+1}}\psi_1(t_{1})d{\bf w}_{t_1}^{(i_1)}
\ldots 
\int\limits_{t_{l-3}}^{t_{l+1}}\psi_{l-2}(t_{l-2})d{\bf w}_{t_{l-2}}^{(i_{l-2})}\times\right.
$$

\vspace{2mm}
$$
\left.\times \int\limits_{t_{l-2}}^{t_{l+1}}
\psi_{l}(t_{l-1})\psi_{l-1}(t_{l-1})dt_{l-1}\right)
d{\bf w}_{t_{l+1}}^{(i_{l+1})}\ldots
d{\bf w}_{t_k}^{(i_k)}=
$$

\vspace{2mm}
$$
=\int\limits_t^T \psi_k(t_k)\ldots \int\limits_t^{t_{l+2}}
\psi_{l+1}(t_{l+1})\left(\int\limits_t^{t_{l+1}}\psi_1(t_{1})d{\bf w}_{t_1}^{(i_1)}
\ldots 
\int\limits_{t_{l-3}}^{t_{l+1}}\psi_{l-2}(t_{l-2})d{\bf w}_{t_{l-2}}^{(i_{l-2})}\times\right.
$$

\vspace{2mm}
$$
\left.\times \left(\int\limits_t^{t_{l+1}}-\int\limits_t^{t_{l-2}}~\right)
\psi_{l}(t_{l-1})\psi_{l-1}(t_{l-1})dt_{l-1}\right)
d{\bf w}_{t_{l+1}}^{(i_{l+1})}\ldots
d{\bf w}_{t_k}^{(i_k)}=
$$

\vspace{2mm}
$$
=\int\limits_t^T \psi_k(t_k)\ldots \int\limits_t^{t_{l+2}}
\psi_{l+1}(t_{l+1})
\left(\int\limits_t^{t_{l+1}}
\psi_{l}(t_{l-1})\psi_{l-1}(t_{l-1})dt_{l-1}\right)
\int\limits_t^{t_{l+1}}\psi_1(t_{1})d{\bf w}_{t_1}^{(i_1)}
\ldots 
$$

\vspace{2mm}
$$
\ldots\int\limits_{t_{l-3}}^{t_{l+1}}\psi_{l-2}(t_{l-2})d{\bf w}_{t_{l-2}}^{(i_{l-2})}
d{\bf w}_{t_{l+1}}^{(i_{l+1})}\ldots
d{\bf w}_{t_k}^{(i_k)}-
$$

\vspace{2mm}
$$
-\int\limits_t^T \psi_k(t_k)\ldots \int\limits_t^{t_{l+2}}
\psi_{l+1}(t_{l+1})\int\limits_t^{t_{l+1}}\psi_1(t_{1})d{\bf w}_{t_1}^{(i_1)}
\ldots 
\int\limits_{t_{l-3}}^{t_{l+1}}\psi_{l-2}(t_{l-2})\times
$$

\vspace{2mm}
$$
\times \left(\int\limits_t^{t_{l-2}}
\psi_{l}(t_{l-1})\psi_{l-1}(t_{l-1})dt_{l-1}\right)
d{\bf w}_{t_{l-2}}^{(i_{l-2})}
d{\bf w}_{t_{l+1}}^{(i_{l+1})}\ldots
d{\bf w}_{t_k}^{(i_k)}=
$$

\vspace{2mm}
$$
=\int\limits_t^T \psi_k(t_k)\ldots \int\limits_t^{t_{l+2}}
\psi_{l+1}(t_{l+1}) \left(\int\limits_t^{t_{l+1}}\psi_l(t_{l-1})\psi_{l-1}(t_{l-1})dt_{l-1}\right)
\int\limits_t^{t_{l+1}}\psi_{l-2}(t_{l-2})\ldots
$$

\vspace{2mm}
$$
\ldots \int\limits_t^{t_2} \psi_1(t_1)d{\bf w}_{t_1}^{(i_1)}\ldots
d{\bf w}_{t_{l-2}}^{(i_{l-2})}d{\bf w}_{t_{l+1}}^{(i_{l+1})}\ldots
d{\bf w}_{t_k}^{(i_k)}-
$$

\vspace{2mm}
$$
-\int\limits_t^T \psi_k(t_k)\ldots \int\limits_t^{t_{l+2}}
\psi_{l+1}(t_{l+1}) 
\int\limits_t^{t_{l+1}}\psi_{l-2}(t_{l-2})
\left(\int\limits_t^{t_{l-2}}\psi_l(t_{l-1})\psi_{l-1}(t_{l-1})dt_{l-1}\right)\times
$$

\vspace{2mm}
\begin{equation}
\label{febr17}
\times\int\limits_t^{t_{l-2}}\psi_{l-3}(t_{l-3})
\ldots
\int\limits_t^{t_2} \psi_1(t_1)d{\bf w}_{t_1}^{(i_1)}\ldots
d{\bf w}_{t_{l-3}}^{(i_{l-3})}d{\bf w}_{t_{l-2}}^{(i_{l-2})}d{\bf w}_{t_{l+1}}^{(i_{l+1})}\ldots
d{\bf w}_{t_k}^{(i_k)}.
\end{equation}

\vspace{3mm}

Comparing the right-hand sides of (\ref{febr15}) and (\ref{febr17})
we come to the conclusion that we got the same result.

The strict mathematical meaning of the transformations 
leading to (\ref{febr17}) is explained in Chapter~3 \cite{2018a}
at least for the case when $\psi_1(\tau),\ldots,\psi_k(\tau)$
are continuous functions on the interval $[t, T].$

Note that under the conditions of Theorem~47, the derivation of the 
formulas (\ref{febr15}) and (\ref{febr17}) will 
remain valid if in (\ref{febr15}) and (\ref{febr17}) we replace all 
differentials of the form $d{\bf w}^{(i_j)}_{t_j}$ with $dt_j$
(this follows from Fubini's Theorem).

Recall that
$$
J[\psi^{(k)}]_{T,t}^{s_r,\ldots,s_1} \stackrel{\rm def}{=}\ 
\prod_{q=1}^r {\bf 1}_{\{i_{s_q}=
i_{s_{q}+1}\ne 0\}}\ \times
$$

\vspace{2mm}
$$
\times
\int\limits_t^T\psi_k(t_k)\ldots \int\limits_t^{t_{s_r+3}}
\psi_{s_r+2}(t_{s_r+2})
\int\limits_t^{t_{s_r+2}}\psi_{s_r}(t_{s_r+1})
\psi_{s_r+1}(t_{s_r+1}) \times
$$

\vspace{2mm}
$$
\times
\int\limits_t^{t_{s_r+1}}\psi_{s_r-1}(t_{s_r-1})
\ldots
\int\limits_t^{t_{s_1+3}}\psi_{s_1+2}(t_{s_1+2})
\int\limits_t^{t_{s_1+2}}\psi_{s_1}(t_{s_1+1})
\psi_{s_1+1}(t_{s_1+1}) \times
$$

\vspace{2mm}
$$
\times
\int\limits_t^{t_{s_1+1}}\psi_{s_1-1}(t_{s_1-1})
\ldots \int\limits_t^{t_2}\psi_1(t_1)
d{\bf w}_{t_1}^{(i_1)}\ldots d{\bf w}_{t_{s_1-1}}^{(i_{s_1-1})}
dt_{s_1+1}d{\bf w}_{t_{s_1+2}}^{(i_{s_1+2})}\ldots
$$

\vspace{2mm}
$$
\ldots\
d{\bf w}_{t_{s_r-1}}^{(i_{s_r-1})}
dt_{s_r+1}d{\bf w}_{t_{s_r+2}}^{(i_{s_r+2})}\ldots d{\bf w}_{t_k}^{(i_k)},
$$

\vspace{5mm}
\noindent                                     
where
${\rm A}_{k,r}$ is defined by (\ref{30.5550001}):

\vspace{1mm}
$$
{\rm A}_{k,r}
=\left\{(s_r,\ldots,s_1):\
s_r>s_{r-1}+1,\ldots,s_2>s_1+1,\ s_r,\ldots,s_1=1,\ldots,k-1\right\}.
$$

\vspace{4mm}

Temporarily denote $J[\psi^{(k)}]_{T,t}^{s_r,\ldots,s_1}$ as 
$I[\psi^{(k)}]_{T,t}^{(i_1\ldots i_{s_1-1}i_{s_1+2} \ldots i_{s_r-1}i_{s_r+2}\ldots i_k)}$.
Let us carry out the trans\-for\-ma\-ti\-on 
(\ref{febr15}) for the iterated Ito stochastic integral 

\vspace{1mm}
$$
I[\psi^{(k)}]_{T,t}^{(i_1\ldots i_{s_1-1}i_{s_1+2} \ldots i_{s_r-1}i_{s_r+2}\ldots i_k)}
$$

\vspace{4mm}
\noindent
iteratively for $s_1,\ldots,s_r.$ 
After this, apply (\ref{dsds11}) to each of the obtained iterated Ito
stochastic integrals.
As a result, we obtain w.~p.~1

\vspace{1mm}
$$
I[\psi^{(k)}]_{T,t}^{(i_1\ldots i_{s_1-1}i_{s_1+2} \ldots i_{s_r-1}i_{s_r+2}\ldots i_k)}=
$$

\vspace{2mm}
$$
=
\prod_{q=1}^r {\bf 1}_{\{i_{s_q}=
i_{s_{q}+1}\ne 0\}}
\times
$$

\vspace{2mm}
$$
\times\sum\limits_{d=1}^{2^r} \left(
\hat I[\psi^{(k)}]_{T,t}^{d(i_1\ldots i_{s_1-1}i_{s_1+2} \ldots i_{s_r-1}i_{s_r+2}\ldots i_k)}
-
\bar I[\psi^{(k)}]_{T,t}^{d(i_1\ldots i_{s_1-1}i_{s_1+2} \ldots i_{s_r-1}i_{s_r+2}\ldots i_k)}
\right)=
$$

\vspace{2mm}
$$
=\prod_{q=1}^r {\bf 1}_{\{i_{s_q}=
i_{s_{q}+1}\ne 0\}}\times
$$

\vspace{3mm}
$$
\times\hbox{\vtop{\offinterlineskip\halign{
\hfil#\hfil\cr
{\rm l.i.m.}\cr
$\stackrel{}{{}_{p\to \infty}}$\cr
}} }
\sum\limits_{j_1,\ldots, j_{s_1-1},j_{s_1+2}, \ldots, j_{s_r-1},j_{s_r+2},\ldots, j_k=0}^p
\ \sum\limits_{d=1}^{2^r}\ \Biggl(
\hat C_{j_1\ldots j_{s_1-1}j_{s_1+2}\ldots j_{s_r-1}j_{s_r+2}\ldots j_k}^{(d)}-\Biggr.
$$

\vspace{4mm}
$$
\Biggl.-
\bar C_{j_1\ldots j_{s_1-1}j_{s_1+2}\ldots j_{s_r-1}j_{s_r+2}\ldots j_k}^{(d)}\Biggr)\times
$$

\vspace{3mm}
\begin{equation}
\label{febr25}
\times
J'[\phi_{j_1}\ldots \phi_{j_{s_1-1}}\phi_{j_{s_1+2}}\ldots \phi_{j_{s_r-1}}
\phi_{j_{s_r+2}}\ldots 
\phi_{j_k}]_{T,t}^{(i_1\ldots i_{s_1-1}i_{s_1+2} \ldots i_{s_r-1}i_{s_r+2}\ldots i_k)},
\end{equation}

\vspace{5mm}
\noindent
where some terms in the sum
$$
\sum\limits_{d=1}^{2^r}
$$

\vspace{2mm}
\noindent
can be identically equal to zero due to
the remark to (\ref{febr15}).

Taking into account that the iterated Ito stochastic integrals
$\hat I[\psi^{(k)}]_{T,t}^{d(i_1\ldots i_{s_1-1}i_{s_1+2} \ldots i_{s_r-1}i_{s_r+2}\ldots i_k)}$
and the Fourier coefficients
$\hat C_{j_1\ldots j_{s_1-1}j_{s_1+2}\ldots j_{s_r-1}j_{s_r+2}\ldots j_k}^{(d)}$
are 
formed on the basis of the same kernels (the same applies to 
the iterated Ito stochastic integrals
$\bar I[\psi^{(k)}]_{T,t}^{d(i_1\ldots i_{s_1-1}i_{s_1+2} \ldots i_{s_r-1}i_{s_r+2}\ldots i_k)}$
and 
the Fourier coefficients
$\bar C_{j_1\ldots j_{s_1-1}j_{s_1+2}\ldots j_{s_r-1}j_{s_r+2}\ldots j_k}^{(d)}$), as well 
as a remark about the relationship of the transformation (\ref{febr15}) 
based on the Ito formula                  
and on the basis 
of classical rules for replacing
the order of integration (see the derivation of (\ref{febr17})), we obtain using Fubini's theorem 
(applying the inverse transformation from (\ref{febr17}) to (\ref{febr20})
in which all differentials of the form
$d{\bf w}^{(i_j)}_{t_j}$ are replaced with $dt_j$)

\vspace{-2mm}
$$
\sum\limits_{d=1}^{2^r}\ \Biggl(
\hat C_{j_1\ldots j_{s_1-1}j_{s_1+2}\ldots j_{s_r-1}j_{s_r+2}\ldots j_k}^{(d)}-
\bar C_{j_1\ldots j_{s_1-1}j_{s_1+2}\ldots j_{s_r-1}j_{s_r+2}\ldots j_k}^{(d)}\Biggr)=
$$

\begin{equation}
\label{febr26}
=C_{j_k \ldots j_1}\biggl|_{(j_{g_2} j_{g_1})\curvearrowright (\cdot)
\ldots (j_{g_{2r}} j_{g_{2r-1}})\curvearrowright (\cdot),
j_{g_{{}_{1}}}=~j_{g_{{}_{2}}},\ldots, j_{g_{{}_{2r-1}}}=~j_{g_{{}_{2r}}}}\biggr.,
\end{equation}

\vspace{5mm}
\noindent
where $g_{2}=g_{1}+1,\ldots, g_{2r}=g_{2r-1}+1.$
Combining (\ref{febr25}) and (\ref{febr26}), we get w.~p.~1

\vspace{2.5mm}
$$
I[\psi^{(k)}]_{T,t}^{(i_1\ldots i_{s_1-1}i_{s_1+2} \ldots i_{s_r-1}i_{s_r+2}\ldots i_k)}=
$$

\vspace{2mm}
$$
=\hbox{\vtop{\offinterlineskip\halign{
\hfil#\hfil\cr
{\rm l.i.m.}\cr
$\stackrel{}{{}_{p\to \infty}}$\cr
}} }
\sum\limits_{\stackrel{j_1,\ldots,j_q,\ldots,j_k=0}{{}_{q\ne g_1, g_2,\ldots, g_{2r-1}, g_{2r}}}}^p
C_{j_k \ldots j_1}\biggl|_{(j_{g_2} j_{g_1})\curvearrowright (\cdot)
\ldots (j_{g_{2r}} j_{g_{2r-1}})\curvearrowright (\cdot),
j_{g_{{}_{1}}}=~j_{g_{{}_{2}}},\ldots, j_{g_{{}_{2r-1}}}=~j_{g_{{}_{2r}}}}\biggr.
\times 
$$

\vspace{4mm}
$$
\times
\prod\limits_{s=1}^r
{\bf 1}_{\{i_{g_{{}_{2s-1}}}=~i_{g_{{}_{2s}}}\ne 0\}}
J'[\phi_{j_{q_1}}\ldots \phi_{j_{q_{k-2r}}}]_{T,t}^{(i_{q_1}\ldots i_{q_{k-2r}})},
$$

\vspace{5mm}
\noindent
where we use the notations from Sect.~7. The equality (\ref{febr14}) is proved
for the case when $\{\phi_j(x)\}_{j=0}^{\infty}$ is an arbitrary
complete orthonormal system of functions
in the space $L_2([t, T]).$ 
Thus, the condition $\phi_0(x)=1/\sqrt{T-t}$ in Theorems~39--41 can be omitted.

Let us separately explain why the condition 
$\psi_l(\tau)\psi_{l-1}(\tau)\in L_2([t, T])$ 
$(l=2, 3,\ldots, k)$
in Theorems 40, 41 can also be omitted.

It is easy to see that the kernels $\hat K_d(t_1,\ldots, t_{k-2r})$ and
$\bar K_d(t_1,\ldots, t_{k-2r})$ of the iterated Ito
stochastic integrals
$\hat I[\psi^{(k)}]_{T,t}^{d(i_1\ldots i_{s_1-1}i_{s_1+2} \ldots i_{s_r-1}i_{s_r+2}\ldots i_k)}$
and 
$\bar I[\psi^{(k)}]_{T,t}^{d(i_1\ldots i_{s_1-1}i_{s_1+2} \ldots i_{s_r-1}i_{s_r+2}\ldots i_k)}$
have the same structure as (\ref{ppp}) but with new wight 
functions $\hat \psi_1(\tau),\ldots, \hat\psi_{k-2r}(\tau)$
and $\bar \psi_1(\tau),\ldots, \bar\psi_{k-2r}(\tau)$, some of which 
possibly coincide with $\psi_1(\tau),\ldots, \psi_k(\tau)\in L_2([t, T])$
(see (\ref{febr15})).
Moreover, the conditions $\psi_1(\tau),\ldots,\psi_{k}(\tau)\in L_2([t, T])$
and $\psi_l(\tau)\psi_{l-1}(\tau)\in L_1([t, T])$ 
$(l=2, 3,\ldots, k)$  guarantees that
$\hat K_d(t_1,\ldots,$ $t_{k-2r}),
\bar K_d(t_1,\ldots, t_{k-2r})\in L_2([t, T])$ (see 
(\ref{febr15})).
This means that the formula (\ref{febr25}) is true if
$\psi_1(\tau),\ldots\hspace{-0.3mm},$ $\psi_{k}(\tau)\in L_2([t, T])$
and $\psi_l(\tau)\psi_{l-1}(\tau)\in L_1([t, T])$ 
$(l=2, 3,\ldots, k).$
Furthermore, the formula (\ref{febr26}) holds under
the conditions 
$\psi_1(\tau),\ldots,\psi_{k}(\tau)\in L_2([t, T])$
and $\psi_l(\tau)\psi_{l-1}(\tau)\in L_1([t, T])$ 
$(l=2, 3,\ldots, k).$

Since the condition $\psi_1(\tau),\ldots,\psi_{k}(\tau)\in L_2([t, T])$
implies the condition $\psi_l(\tau)\psi_{l-1}(\tau)\in L_1([t, T])$ 
$(l=2, 3,\ldots, k),$ then the condition $\psi_l(\tau)\psi_{l-1}(\tau)\in L_1([t, T])$ 
$(l=2, 3,\ldots, k)$ can be omitted in the above reasoning.

Thus, the equalities (\ref{febr25}) and (\ref{febr26})
are satisfied under the condition 
$\psi_1(\tau),\ldots,\psi_{k}(\tau)\in L_2([t, T])$
and the condition $\psi_l(\tau)\psi_{l-1}(\tau)\in L_2([t, T])$ 
$(l=2, 3,\ldots, k)$ can be omitted in 
Theorems~40, 41.
Theorem~47 is proved.

\vspace{5mm}

\section{Expansion of Iterated Stratonovich Stochastic Integrals
of Multiplicity 5. The Case of an Ar\-bit\-ra\-ry Complete Orthonormal System of 
Functions in the Space $L_2([t,T])$ and $\psi_1(\tau),\ldots, \psi_5(\tau)
\equiv 1$}

\vspace{5mm}

{\bf Theorem~48}\ \cite{2018a}.\ {\it Suppose that
$\{\phi_j(x)\}_{j=0}^{\infty}$ is an arbitrary complete orthonormal system of 
func\-ti\-ons in the space $L_2([t,T]).$
Then$,$ for the iterated Stra\-to\-no\-vich stochastic integral
of fifth multiplicity 

$$
J^{*}[\psi^{(5)}]_{T,t}=
{\int\limits_t^{*}}^T
\ldots
{\int\limits_t^{*}}^{t_2}
d{\bf w}_{t_1}^{(i_1)}
\ldots d{\bf w}_{t_5}^{(i_5)}
$$

\vspace{2mm}
\noindent
the following 
expansion 

\vspace{-1mm}
$$
J^{*}[\psi^{(5)}]_{T,t}=
\hbox{\vtop{\offinterlineskip\halign{
\hfil#\hfil\cr
{\rm l.i.m.}\cr
$\stackrel{}{{}_{p\to \infty}}$\cr
}} }
\sum\limits_{j_1,\ldots, j_5=0}^{p}
C_{j_5 \ldots j_1}\zeta_{j_1}^{(i_1)}\ldots \zeta_{j_5}^{(i_5)}
$$

\vspace{3mm}
\noindent
that converges in the mean-square sense is valid, where 
$i_1,\ldots,i_5=0, 1,\ldots,m,$

$$
C_{j_5\ldots j_1}=\int\limits_t^T
\phi_{j_5}(t_5)
\ldots
\int\limits_t^{t_2}
\phi_{j_1}(t_1)dt_1\ldots dt_5
$$

\vspace{1mm}
\noindent
and
$$
\zeta_{j}^{(i)}=
\int\limits_t^T \phi_{j}(\tau) d{\bf w}_{\tau}^{(i)}
$$ 

\vspace{2mm}
\noindent
are independent standard Gaussian random variables for various 
$i$ or $j$ {\rm (}in the case when $i\ne 0${\rm ),}
${\bf w}_{\tau}^{(i)}={\bf f}_{\tau}^{(i)}$ for
$i=1,\ldots,m$ and 
${\bf w}_{\tau}^{(0)}=\tau.$}

\vspace{2mm}

{\bf Proof.}\ {\bf Step~1.}\ According to Theorem~47,
we conclude that Theorem~48 will be proved if we prove
the following equalities (see (\ref{june100}) for $k=5, r=1$ and $k=5, r=2$
($p_1=\ldots=p_5=p$))
under the conditions of Theorem~48   

\vspace{-1mm}
\begin{equation}
\label{april11}
\lim\limits_{p\to\infty}
\sum\limits_{j_3, j_4,j_5=0}^{p}
\left(\sum\limits_{j_1=0}^{p}
C_{j_5 j_4 j_3 j_1 j_1}-\frac{1}{2} 
C_{j_5 j_4 j_3 j_1 j_1}\biggl|_{(j_{1} j_{1})\curvearrowright (\cdot)}
\biggr.\right)^2=0,
\end{equation}

\begin{equation}
\label{april12}
\lim\limits_{p\to\infty}
\sum\limits_{j_2, j_4,j_5=0}^{p}
\left(\sum\limits_{j_1=0}^{p}
C_{j_5 j_4 j_1 j_2 j_1}\right)^2=0,
\end{equation}

\begin{equation}
\label{april13}
\lim\limits_{p\to\infty}
\sum\limits_{j_2, j_3,j_5=0}^{p}
\left(\sum\limits_{j_1=0}^{p}
C_{j_5 j_1 j_3 j_2 j_1}\right)^2=0,
\end{equation}

\begin{equation}
\label{april14}
\lim\limits_{p\to\infty}
\sum\limits_{j_2, j_3,j_4=0}^{p}
\left(\sum\limits_{j_1=0}^{p}
C_{j_1 j_4 j_3 j_2 j_1}\right)^2=0,
\end{equation}

\begin{equation}
\label{april15}
\lim\limits_{p\to\infty}
\sum\limits_{j_1, j_4,j_5=0}^{p}
\left(\sum\limits_{j_2=0}^{p}
C_{j_5 j_4 j_2 j_2 j_1}-\frac{1}{2} 
C_{j_5 j_4 j_2 j_2 j_1}\biggl|_{(j_{2} j_{2})\curvearrowright (\cdot)}
\biggr.\right)^2=0,
\end{equation}

\begin{equation}
\label{april16}
\lim\limits_{p\to\infty}
\sum\limits_{j_1, j_3,j_5=0}^{p}
\left(\sum\limits_{j_2=0}^{p}
C_{j_5 j_2 j_3 j_2 j_1}\right)^2=0,
\end{equation}

\begin{equation}
\label{april17}
\lim\limits_{p\to\infty}
\sum\limits_{j_1, j_3,j_4=0}^{p}
\left(\sum\limits_{j_2=0}^{p}
C_{j_2 j_4 j_3 j_2 j_1}\right)^2=0,
\end{equation}

\begin{equation}
\label{april18}
\lim\limits_{p\to\infty}
\sum\limits_{j_1, j_2,j_5=0}^{p}
\left(\sum\limits_{j_3=0}^{p}
C_{j_5 j_3 j_3 j_2 j_1}-\frac{1}{2} 
C_{j_5 j_3 j_3 j_2 j_1}\biggl|_{(j_{3} j_{3})\curvearrowright (\cdot)}
\biggr.\right)^2=0,
\end{equation}

\begin{equation}
\label{april19}
\lim\limits_{p\to\infty}
\sum\limits_{j_1, j_2,j_4=0}^{p}
\left(\sum\limits_{j_3=0}^{p}
C_{j_3 j_4 j_3 j_2 j_1}\right)^2=0,
\end{equation}

\begin{equation}
\label{april20}
\lim\limits_{p\to\infty}
\sum\limits_{j_1, j_2,j_3=0}^{p}
\left(\sum\limits_{j_4=0}^{p}
C_{j_4 j_4 j_3 j_2 j_1}-\frac{1}{2} 
C_{j_4 j_4 j_3 j_2 j_1}\biggl|_{(j_{4} j_{4})\curvearrowright (\cdot)}
\biggr.\right)^2=0,
\end{equation}

\begin{equation}
\label{april21}
\lim\limits_{p\to\infty}
\sum\limits_{j_5=0}^{p}
\left(\sum\limits_{j_1,j_3=0}^{p}
C_{j_5 j_3 j_3 j_1 j_1}-\frac{1}{4} 
C_{j_5 j_3 j_3 j_1 j_1}\biggl|_{(j_{1} j_{1})\curvearrowright (\cdot),
(j_{3} j_{3})\curvearrowright (\cdot)}
\biggr.\right)^2=0,
\end{equation}

\begin{equation}
\label{april22}
\lim\limits_{p\to\infty}
\sum\limits_{j_4=0}^{p}
\left(\sum\limits_{j_1,j_3=0}^{p}
C_{j_3 j_4 j_3 j_1 j_1}\right)^2=0,
\end{equation}

\begin{equation}
\label{april23}
\lim\limits_{p\to\infty}
\sum\limits_{j_3=0}^{p}
\left(\sum\limits_{j_1,j_4=0}^{p}
C_{j_4 j_4 j_3 j_1 j_1}-\frac{1}{4} 
C_{j_4 j_4 j_3 j_1 j_1}\biggl|_{(j_{1} j_{1})\curvearrowright (\cdot),
(j_{4} j_{4})\curvearrowright (\cdot)}
\biggr.\right)^2=0,
\end{equation}

\begin{equation}
\label{april24}
\lim\limits_{p\to\infty}
\sum\limits_{j_5=0}^{p}
\left(\sum\limits_{j_1,j_2=0}^{p}
C_{j_5 j_2 j_1 j_2 j_1}\right)^2=0,
\end{equation}

\begin{equation}
\label{april25}
\lim\limits_{p\to\infty}
\sum\limits_{j_4=0}^{p}
\left(\sum\limits_{j_1,j_2=0}^{p}
C_{j_2 j_4 j_1 j_2 j_1}\right)^2=0,
\end{equation}

\begin{equation}
\label{april26}
\lim\limits_{p\to\infty}
\sum\limits_{j_2=0}^{p}
\left(\sum\limits_{j_1,j_4=0}^{p}
C_{j_4 j_4 j_1 j_2 j_1}\right)^2=0,
\end{equation}

\begin{equation}
\label{april27}
\lim\limits_{p\to\infty}
\sum\limits_{j_5=0}^{p}
\left(\sum\limits_{j_1,j_2=0}^{p}
C_{j_5 j_1 j_2 j_2 j_1}\right)^2=0,
\end{equation}

\begin{equation}
\label{april28}
\lim\limits_{p\to\infty}
\sum\limits_{j_3=0}^{p}
\left(\sum\limits_{j_1,j_2=0}^{p}
C_{j_2 j_1 j_3 j_2 j_1}\right)^2=0,
\end{equation}

\begin{equation}
\label{april29}
\lim\limits_{p\to\infty}
\sum\limits_{j_2=0}^{p}
\left(\sum\limits_{j_1,j_3=0}^{p}
C_{j_3 j_1 j_3 j_2 j_1}\right)^2=0,
\end{equation}

\begin{equation}
\label{april30}
\lim\limits_{p\to\infty}
\sum\limits_{j_4=0}^{p}
\left(\sum\limits_{j_1,j_2=0}^{p}
C_{j_1 j_4 j_2 j_2 j_1}\right)^2=0,
\end{equation}

\begin{equation}
\label{april31}
\lim\limits_{p\to\infty}
\sum\limits_{j_3=0}^{p}
\left(\sum\limits_{j_1,j_2=0}^{p}
C_{j_1 j_2 j_3 j_2 j_1}\right)^2=0,
\end{equation}

\begin{equation}
\label{april32}
\lim\limits_{p\to\infty}
\sum\limits_{j_2=0}^{p}
\left(\sum\limits_{j_1,j_3=0}^{p}
C_{j_1 j_3 j_3 j_2 j_1}\right)^2=0,
\end{equation}

\begin{equation}
\label{april33}
\lim\limits_{p\to\infty}
\sum\limits_{j_1=0}^{p}
\left(\sum\limits_{j_2,j_4=0}^{p}
C_{j_4 j_4 j_2 j_2 j_1}-\frac{1}{4} 
C_{j_4 j_4 j_2 j_2 j_1}\biggl|_{(j_{2} j_{2})\curvearrowright (\cdot),
(j_{4} j_{4})\curvearrowright (\cdot)}
\biggr.\right)^2=0,
\end{equation}

\begin{equation}
\label{april34}
\lim\limits_{p\to\infty}
\sum\limits_{j_1=0}^{p}
\left(\sum\limits_{j_2,j_3=0}^{p}
C_{j_3 j_2 j_3 j_2 j_1}\right)^2=0,
\end{equation}

\begin{equation}
\label{april35}
\lim\limits_{p\to\infty}
\sum\limits_{j_1=0}^{p}
\left(\sum\limits_{j_2,j_3=0}^{p}
C_{j_2 j_3 j_3 j_2 j_1}\right)^2=0.
\end{equation}

\vspace{4mm}

{\bf Step~2.}\ Let us prove the equalities (\ref{april11})--(\ref{april20}).
Using Fubini's Theorem and Parseval's equality, we obtain
the following relations 
for the prelimit
expressions on the left-hand sides of (\ref{april11})--(\ref{april20})

\vspace{-1mm}
$$
\sum\limits_{j_3, j_4,j_5=0}^{p}
\left(\sum\limits_{j_1=0}^{p}
C_{j_5 j_4 j_3 j_1 j_1}-\frac{1}{2} 
C_{j_5 j_4 j_3 j_1 j_1}\biggl|_{(j_{1} j_{1})\curvearrowright (\cdot)}
\biggr.\right)^2=
$$

\vspace{1mm}
$$
=\sum\limits_{j_3, j_4,j_5=0}^{p}
\left(\int\limits_t^T\phi_{j_5}(t_5)\int\limits_t^{t_5}\phi_{j_4}(t_4)
\int\limits_t^{t_4}\phi_{j_3}(t_3)\left(\sum\limits_{j_1=0}^{p}
\frac{1}{2}\left(\int\limits_t^{t_3}\phi_{j_1}(\tau)d\tau\right)^2-\frac{t_3-t}{2}\right)
dt_3 dt_4 dt_5\right)^2\le
$$

\vspace{1mm}
$$
\le\sum\limits_{j_3, j_4,j_5=0}^{\infty}
\left(\int\limits_t^T\phi_{j_5}(t_5)\int\limits_t^{t_5}\phi_{j_4}(t_4)
\int\limits_t^{t_4}\phi_{j_3}(t_3)\left(\sum\limits_{j_1=0}^{p}
\frac{1}{2}\left(\int\limits_t^{t_3}\phi_{j_1}(\tau)d\tau\right)^2-\frac{t_3-t}{2}\right)
dt_3 dt_4 dt_5\right)^2=
$$

\vspace{1mm}
\begin{equation}
\label{april36}
=\int\limits_{[t,T]^3}\left({\bf 1}_{\{t_3<t_4<t_5\}}\right)^2
\left(\sum\limits_{j_1=0}^{p}
\frac{1}{2}\left(\int\limits_t^{t_3}\phi_{j_1}(\tau)d\tau\right)^2-\frac{t_3-t}{2}\right)^2
dt_3 dt_4 dt_5,
\end{equation}

\vspace{7mm}
$$
\sum\limits_{j_2, j_4,j_5=0}^{p}
\left(\sum\limits_{j_1=0}^{p}
C_{j_5 j_4 j_1 j_2 j_1}\right)^2=
$$

\vspace{1mm}
$$
=\sum\limits_{j_2, j_4,j_5=0}^{p}
\left(\int\limits_t^T\phi_{j_5}(t_5)\int\limits_t^{t_5}\phi_{j_4}(t_4)
\int\limits_t^{t_4}\phi_{j_2}(t_2)\sum\limits_{j_1=0}^{p}
\int\limits_t^{t_2}\phi_{j_1}(t_1)dt_1
\int\limits_{t_2}^{t_4}\phi_{j_1}(t_3)dt_3
dt_2 dt_4 dt_5\right)^2\le
$$

\vspace{1mm}
$$
\le\sum\limits_{j_2, j_4,j_5=0}^{\infty}
\left(\int\limits_t^T\phi_{j_5}(t_5)\int\limits_t^{t_5}\phi_{j_4}(t_4)
\int\limits_t^{t_4}\phi_{j_2}(t_2)\sum\limits_{j_1=0}^{p}
\int\limits_t^{t_2}\phi_{j_1}(t_1)dt_1
\int\limits_{t_2}^{t_4}\phi_{j_1}(t_3)dt_3dt_2 dt_4 dt_5\right)^2=
$$

\vspace{1mm}
\begin{equation}
\label{april37}
=\int\limits_{[t,T]^3}\left({\bf 1}_{\{t_2<t_4<t_5\}}\right)^2
\left(\sum\limits_{j_1=0}^{p}
\int\limits_t^{t_2}\phi_{j_1}(t_1)dt_1
\int\limits_{t_2}^{t_4}\phi_{j_1}(t_3)dt_3\right)^2
dt_2 dt_4 dt_5,
\end{equation}

\vspace{7mm}

$$
\sum\limits_{j_2, j_3,j_5=0}^{p}
\left(\sum\limits_{j_1=0}^{p}
C_{j_5 j_1 j_3 j_2 j_1}\right)^2=
$$

\vspace{1mm}
$$
=\sum\limits_{j_2, j_3,j_5=0}^{p}
\left(\int\limits_t^T\phi_{j_5}(t_5)\int\limits_t^{t_5}\phi_{j_3}(t_3)
\int\limits_t^{t_3}\phi_{j_2}(t_2)\sum\limits_{j_1=0}^{p}
\int\limits_t^{t_2}\phi_{j_1}(t_1)dt_1
\int\limits_{t_3}^{t_5}\phi_{j_1}(t_4)dt_4dt_2 dt_3 dt_5\right)^2\le
$$

\vspace{1mm}
$$
\le\sum\limits_{j_2, j_3,j_5=0}^{\infty}
\left(\int\limits_t^T\phi_{j_5}(t_5)\int\limits_t^{t_5}\phi_{j_3}(t_3)
\int\limits_t^{t_3}\phi_{j_2}(t_2)\sum\limits_{j_1=0}^{p}
\int\limits_t^{t_2}\phi_{j_1}(t_1)dt_1
\int\limits_{t_3}^{t_5}\phi_{j_1}(t_4)dt_4dt_2 dt_3 dt_5\right)^2=
$$

\vspace{1mm}
\begin{equation}
\label{april38}
=\int\limits_{[t,T]^3}\left({\bf 1}_{\{t_2<t_3<t_5\}}\right)^2
\left(\sum\limits_{j_1=0}^{p}
\int\limits_t^{t_2}\phi_{j_1}(t_1)dt_1
\int\limits_{t_3}^{t_5}\phi_{j_1}(t_4)dt_4\right)^2
dt_2 dt_3 dt_5,
\end{equation}

\vspace{7mm}

$$
\sum\limits_{j_2, j_3,j_4=0}^{p}
\left(\sum\limits_{j_1=0}^{p}
C_{j_1 j_4 j_3 j_2 j_1}\right)^2=
$$

\vspace{1mm}
$$
=\sum\limits_{j_2, j_3,j_4=0}^{p}
\left(\int\limits_t^T\phi_{j_4}(t_4)\int\limits_t^{t_4}\phi_{j_3}(t_3)
\int\limits_t^{t_3}\phi_{j_2}(t_2)\sum\limits_{j_1=0}^{p}
\int\limits_t^{t_2}\phi_{j_1}(t_1)dt_1
\int\limits_{t_4}^{T}\phi_{j_1}(t_5)dt_5 dt_2 dt_3 dt_4\right)^2\le
$$

\vspace{1mm}
$$
\le\sum\limits_{j_2, j_3,j_4=0}^{\infty}
\left(\int\limits_t^T\phi_{j_4}(t_4)\int\limits_t^{t_4}\phi_{j_3}(t_3)
\int\limits_t^{t_3}\phi_{j_2}(t_2)\sum\limits_{j_1=0}^{p}
\int\limits_t^{t_2}\phi_{j_1}(t_1)dt_1
\int\limits_{t_4}^{T}\phi_{j_1}(t_5)dt_5 dt_2 dt_3 dt_4\right)^2=
$$

\vspace{1mm}
\begin{equation}
\label{april39}
=\int\limits_{[t,T]^3}\left({\bf 1}_{\{t_2<t_3<t_4\}}\right)^2
\left(\sum\limits_{j_1=0}^{p}
\int\limits_t^{t_2}\phi_{j_1}(t_1)dt_1
\int\limits_{t_4}^{T}\phi_{j_1}(t_5)dt_5\right)^2
dt_2 dt_3 dt_4,
\end{equation}

\vspace{7mm}

$$
\sum\limits_{j_1, j_4,j_5=0}^{p}
\left(\sum\limits_{j_2=0}^{p}
C_{j_5 j_4 j_2 j_2 j_1}-\frac{1}{2} 
C_{j_5 j_4 j_2 j_2 j_1}\biggl|_{(j_{2} j_{2})\curvearrowright (\cdot)}
\biggr.\right)^2=
$$

\vspace{4mm}
$$
=\sum\limits_{j_1, j_4,j_5=0}^{p}
\left(\int\limits_t^T\phi_{j_5}(t_5)\int\limits_t^{t_5}\phi_{j_4}(t_4)
\int\limits_t^{t_4}\phi_{j_1}(t_1)
\sum\limits_{j_2=0}^p \int\limits_{t_1}^{t_4}\phi_{j_2}(t_2)
\int\limits_{t_2}^{t_4}\phi_{j_2}(t_3)dt_3 dt_2 dt_1 dt_4 dt_5-\right.
$$
\vspace{1mm}

$$
\left.
-\frac{1}{2}\int\limits_t^T\phi_{j_5}(t_5)\int\limits_t^{t_5}\phi_{j_4}(t_4)
\int\limits_t^{t_4}\int\limits_t^{t_2}\phi_{j_1}(t_1)dt_1 dt_2 dt_4 dt_5
\right)^2=
$$

\vspace{1mm}
$$
=\sum\limits_{j_1, j_4,j_5=0}^{p}
\left(\int\limits_t^T\phi_{j_5}(t_5)\int\limits_t^{t_5}\phi_{j_4}(t_4)
\int\limits_t^{t_4}\phi_{j_1}(t_1)
\left(\sum\limits_{j_2=0}^p \frac{1}{2}\left(\int\limits_{t_1}^{t_4}\phi_{j_2}(t_2)
dt_2\right)^2-\frac{t_4-t_1}{2}\right) dt_1 dt_4 dt_5
\right)^2\le
$$

\vspace{1mm}
$$
\le\sum\limits_{j_1, j_4,j_5=0}^{\infty}
\left(\int\limits_t^T\phi_{j_5}(t_5)\int\limits_t^{t_5}\phi_{j_4}(t_4)
\int\limits_t^{t_4}\phi_{j_1}(t_1)
\left(\sum\limits_{j_2=0}^p \frac{1}{2}\left(\int\limits_{t_1}^{t_4}\phi_{j_2}(t_2)
dt_2\right)^2-\frac{t_4-t_1}{2}\right) dt_1 dt_4 dt_5
\right)^2=
$$

\vspace{1mm}
\begin{equation}
\label{april40}
=\int\limits_{[t,T]^3}\left({\bf 1}_{\{t_1<t_4<t_5\}}\right)^2
\left(\sum\limits_{j_2=0}^p \frac{1}{2}\left(\int\limits_{t_1}^{t_4}\phi_{j_2}(t_2)
dt_2\right)^2-\frac{t_4-t_1}{2}\right)^2
dt_1 dt_4 dt_5,
\end{equation}

\vspace{7mm}

$$
\sum\limits_{j_1, j_3,j_5=0}^{p}
\left(\sum\limits_{j_2=0}^{p}
C_{j_5 j_2 j_3 j_2 j_1}\right)^2=
$$

\vspace{1mm}
$$
=\sum\limits_{j_1, j_3,j_5=0}^{p}
\left(\sum\limits_{j_2=0}^p \int\limits_t^T\phi_{j_5}(t_5)\int\limits_t^{t_5}\phi_{j_3}(t_3)
\int\limits_t^{t_3}\phi_{j_2}(t_2)
\int\limits_t^{t_2}\phi_{j_1}(t_1)dt_1 dt_2
\int\limits_{t_3}^{t_5}\phi_{j_2}(t_4)dt_4 dt_3 dt_5\right)^2=
$$

\vspace{1mm}
$$
=\sum\limits_{j_1, j_3,j_5=0}^{p}
\left(\int\limits_t^T\phi_{j_5}(t_5)\int\limits_t^{t_5}\phi_{j_3}(t_3)
\int\limits_t^{t_3}\phi_{j_1}(t_1)
\sum\limits_{j_2=0}^p \int\limits_{t_1}^{t_3}\phi_{j_2}(t_2)dt_2
\int\limits_{t_3}^{t_5}\phi_{j_2}(t_4)dt_4 dt_1 dt_3 dt_5\right)^2\le
$$

\vspace{1mm}
$$
\le\sum\limits_{j_1, j_3,j_5=0}^{\infty}
\left(\int\limits_t^T\phi_{j_5}(t_5)\int\limits_t^{t_5}\phi_{j_3}(t_3)
\int\limits_t^{t_3}\phi_{j_1}(t_1)
\sum\limits_{j_2=0}^p \int\limits_{t_1}^{t_3}\phi_{j_2}(t_2)dt_2
\int\limits_{t_3}^{t_5}\phi_{j_2}(t_4)dt_4 dt_1 dt_3 dt_5\right)^2=
$$

\vspace{1mm}
\begin{equation}
\label{april41}
=\int\limits_{[t,T]^3}\left({\bf 1}_{\{t_1<t_3<t_5\}}\right)^2
\left(\sum\limits_{j_2=0}^p \int\limits_{t_1}^{t_3}\phi_{j_2}(t_2)dt_2
\int\limits_{t_3}^{t_5}\phi_{j_2}(t_4)dt_4\right)^2
dt_1 dt_3 dt_5,
\end{equation}

\vspace{7mm}

$$
\sum\limits_{j_1, j_3,j_4=0}^{p}
\left(\sum\limits_{j_2=0}^{p}
C_{j_2 j_4 j_3 j_2 j_1}\right)^2=
$$

\vspace{1mm}
$$
=\sum\limits_{j_1, j_3,j_4=0}^{p}
\left(\sum\limits_{j_2=0}^p \int\limits_t^T\phi_{j_4}(t_4)\int\limits_t^{t_4}\phi_{j_3}(t_3)
\int\limits_t^{t_3}\phi_{j_2}(t_2)
\int\limits_t^{t_2}\phi_{j_1}(t_1)dt_1 dt_2 dt_3
\int\limits_{t_4}^{T}\phi_{j_2}(t_5)dt_5 dt_4\right)^2=
$$

\vspace{1mm}
$$
=\sum\limits_{j_1, j_3,j_4=0}^{p}
\left(\int\limits_t^T\phi_{j_4}(t_4)\int\limits_t^{t_4}\phi_{j_3}(t_3)
\int\limits_t^{t_3}\phi_{j_1}(t_1)
\sum\limits_{j_2=0}^p \int\limits_{t_1}^{t_3}\phi_{j_2}(t_2)dt_2
\int\limits_{t_4}^{T}\phi_{j_2}(t_5)dt_5 dt_1 dt_3 dt_4\right)^2\le
$$

\vspace{1mm}
$$
\le\sum\limits_{j_1, j_3,j_4=0}^{\infty}
\left(\int\limits_t^T\phi_{j_4}(t_4)\int\limits_t^{t_4}\phi_{j_3}(t_3)
\int\limits_t^{t_3}\phi_{j_1}(t_1)
\sum\limits_{j_2=0}^p \int\limits_{t_1}^{t_3}\phi_{j_2}(t_2)dt_2
\int\limits_{t_4}^{T}\phi_{j_2}(t_5)dt_5 dt_1 dt_3 dt_4\right)^2=
$$

\vspace{1mm}
\begin{equation}
\label{april42}
=\int\limits_{[t,T]^3}\left({\bf 1}_{\{t_1<t_3<t_4\}}\right)^2
\left(\sum\limits_{j_2=0}^p \int\limits_{t_1}^{t_3}\phi_{j_2}(t_2)dt_2
\int\limits_{t_4}^{T}\phi_{j_2}(t_5)dt_5\right)^2
dt_1 dt_3 dt_4,
\end{equation}

\vspace{7mm}

$$
\sum\limits_{j_1, j_2,j_5=0}^{p}
\left(\sum\limits_{j_3=0}^{p}
C_{j_5 j_3 j_3 j_2 j_1}-\frac{1}{2} 
C_{j_5 j_3 j_3 j_2 j_1}\biggl|_{(j_{3} j_{3})\curvearrowright (\cdot)}
\biggr.\right)^2=
$$

\vspace{4mm}
$$
=\sum\limits_{j_1, j_2,j_5=0}^{p}
\left(\sum\limits_{j_3=0}^p
\int\limits_t^T\phi_{j_5}(t_5)\int\limits_t^{t_5}\phi_{j_1}(t_1)
\int\limits_{t_1}^{t_5}\phi_{j_2}(t_2)
\int\limits_{t_2}^{t_5}\phi_{j_3}(t_3)
\int\limits_{t_3}^{t_5}\phi_{j_3}(t_4)dt_4 dt_3 dt_2 dt_1 dt_5-\right.
$$

\vspace{1mm}

$$
\left.
-\frac{1}{2}\int\limits_t^T\phi_{j_5}(t_5)\int\limits_t^{t_5}\int\limits_t^{t_3}
\phi_{j_2}(t_2)
\int\limits_t^{t_2}\phi_{j_1}(t_1)dt_1 dt_2 dt_3 dt_5
\right)^2=
$$

\vspace{1mm}
$$
=\sum\limits_{j_1, j_2,j_5=0}^{p}
\left(\sum\limits_{j_3=0}^p
\int\limits_t^T\phi_{j_5}(t_5)\int\limits_t^{t_5}\phi_{j_1}(t_1)
\int\limits_{t_1}^{t_5}\phi_{j_2}(t_2)
\int\limits_{t_2}^{t_5}\phi_{j_3}(t_3)
\int\limits_{t_3}^{t_5}\phi_{j_3}(t_4)dt_4 dt_3 dt_2 dt_1 dt_5-\right.
$$

\vspace{1mm}

$$
\left.
-\frac{1}{2}\int\limits_t^T\phi_{j_5}(t_5)
\int\limits_{t}^{t_5}\phi_{j_1}(t_1)
\int\limits_{t_1}^{t_5}\phi_{j_2}(t_2)\int\limits_{t_2}^{t_5}dt_3 dt_2 dt_1 dt_5
\right)^2=
$$

\vspace{1mm}
$$
=\sum\limits_{j_1, j_2,j_5=0}^{p}
\left(\int\limits_t^T\phi_{j_5}(t_5)\int\limits_t^{t_5}\phi_{j_1}(t_1)
\int\limits_{t_1}^{t_5}\phi_{j_2}(t_2)
\left(\sum\limits_{j_3=0}^p \frac{1}{2}\left(\int\limits_{t_2}^{t_5}\phi_{j_3}(t_3)
dt_3\right)^2-\frac{t_5-t_2}{2}\right)dt_2 dt_1 dt_5
\right)^2\le
$$

\vspace{1mm}
$$
\le\sum\limits_{j_1, j_2,j_5=0}^{\infty}
\left(\int\limits_t^T\phi_{j_5}(t_5)\int\limits_t^{t_5}\phi_{j_1}(t_1)
\int\limits_{t_1}^{t_5}\phi_{j_2}(t_2)
\left(\sum\limits_{j_3=0}^p \frac{1}{2}\left(\int\limits_{t_2}^{t_5}\phi_{j_3}(t_3)
dt_3\right)^2-\frac{t_5-t_2}{2}\right)dt_2 dt_1 dt_5
\right)^2=
$$

\vspace{1mm}
\begin{equation}
\label{april43}
=\int\limits_{[t,T]^3}\left({\bf 1}_{\{t_1<t_2<t_5\}}\right)^2
\left(\sum\limits_{j_3=0}^p \frac{1}{2}\left(\int\limits_{t_2}^{t_5}\phi_{j_3}(t_3)
dt_3\right)^2-\frac{t_5-t_2}{2}\right)^2
dt_2 dt_1 dt_5,
\end{equation}

\vspace{7mm}

$$
\sum\limits_{j_1, j_2,j_4=0}^{p}
\left(\sum\limits_{j_3=0}^{p}
C_{j_3 j_4 j_3 j_2 j_1}\right)^2=
$$

\vspace{1mm}
$$
=\sum\limits_{j_1, j_2,j_4=0}^{p}
\left(\sum\limits_{j_3=0}^p \int\limits_t^T\phi_{j_1}(t_1)\int\limits_{t_1}^{T}\phi_{j_2}(t_2)
\int\limits_{t_2}^{T}\phi_{j_3}(t_3)
\int\limits_{t_3}^{T}\phi_{j_4}(t_4)
\int\limits_{t_4}^{T}\phi_{j_3}(t_5)dt_5 dt_4 dt_3 dt_2 dt_1\right)^2=
$$

\vspace{1mm}
$$
=\sum\limits_{j_1, j_2,j_4=0}^{p}
\left(\int\limits_t^T\phi_{j_1}(t_1)\int\limits_{t_1}^{T}\phi_{j_2}(t_2)
\int\limits_{t_2}^{T}\phi_{j_4}(t_4)
\sum\limits_{j_3=0}^p\int\limits_{t_4}^{T}\phi_{j_3}(t_5)dt_5
\int\limits_{t_2}^{t_4}\phi_{j_3}(t_3)dt_3 dt_4 dt_2 dt_1\right)^2\le
$$

\vspace{1mm}
$$
\le\sum\limits_{j_1, j_2,j_4=0}^{\infty}
\left(\int\limits_t^T\phi_{j_1}(t_1)\int\limits_{t_1}^{T}\phi_{j_2}(t_2)
\int\limits_{t_2}^{T}\phi_{j_4}(t_4)
\sum\limits_{j_3=0}^p\int\limits_{t_4}^{T}\phi_{j_3}(t_5)dt_5
\int\limits_{t_2}^{t_4}\phi_{j_3}(t_3)dt_3 dt_4 dt_2 dt_1\right)^2=
$$

\vspace{1mm}
\begin{equation}
\label{april44}
=\int\limits_{[t,T]^3}\left({\bf 1}_{\{t_1<t_2<t_4\}}\right)^2
\left(
\sum\limits_{j_3=0}^p\int\limits_{t_4}^{T}\phi_{j_3}(t_5)dt_5
\int\limits_{t_2}^{t_4}\phi_{j_3}(t_3)dt_3
\right)^2
dt_4 dt_2 dt_1,
\end{equation}

\vspace{7mm}

$$
\sum\limits_{j_1, j_2,j_3=0}^{p}
\left(\sum\limits_{j_4=0}^{p}
C_{j_4 j_4 j_3 j_2 j_1}-\frac{1}{2} 
C_{j_4 j_4 j_3 j_2 j_1}\biggl|_{(j_{4} j_{4})\curvearrowright (\cdot)}
\biggr.\right)^2=
$$

\vspace{4mm}
$$
=\sum\limits_{j_1, j_2,j_3=0}^{p}
\left(
\int\limits_t^T\phi_{j_3}(t_3)\int\limits_t^{t_3}\phi_{j_2}(t_2)
\int\limits_{t}^{t_2}\phi_{j_1}(t_1) dt_1 dt_2
\sum\limits_{j_4=0}^p\int\limits_{t_3}^{T}\phi_{j_4}(t_4)
\int\limits_{t_4}^{T}\phi_{j_4}(t_5)dt_5 dt_4 dt_3-
\right.
$$

\vspace{1mm}

$$
\left.-
\frac{1}{2}\int\limits_t^T\int\limits_t^{t_4}\phi_{j_3}(t_3)\int\limits_t^{t_3}
\phi_{j_2}(t_2)
\int\limits_t^{t_2}\phi_{j_1}(t_1)dt_1 dt_2 dt_3 dt_4
\right)^2=
$$

\vspace{1mm}
$$
=\sum\limits_{j_1, j_2,j_3=0}^{p}
\left(
\int\limits_t^T\phi_{j_3}(t_3)\int\limits_t^{t_3}\phi_{j_2}(t_2)
\int\limits_{t}^{t_2}\phi_{j_1}(t_1) 
\left(\sum\limits_{j_4=0}^p\frac{1}{2}\left(\int\limits_{t_3}^{T}\phi_{j_4}(t_4)dt_4\right)^2
-\frac{T-t_3}{2}\right)dt_1 dt_2dt_3
\right)^2\le
$$

\vspace{1mm}
$$
\le\sum\limits_{j_1, j_2,j_3=0}^{\infty}
\left(
\int\limits_t^T\phi_{j_3}(t_3)\int\limits_t^{t_3}\phi_{j_2}(t_2)
\int\limits_{t}^{t_2}\phi_{j_1}(t_1) 
\left(\sum\limits_{j_4=0}^p\frac{1}{2}\left(\int\limits_{t_3}^{T}\phi_{j_4}(t_4)dt_4\right)^2
-\frac{T-t_3}{2}\right)dt_1 dt_2dt_3
\right)^2=
$$

\vspace{1mm}
\begin{equation}
\label{april45}
=\int\limits_{[t,T]^3}\left({\bf 1}_{\{t_1<t_2<t_3\}}\right)^2
\left(
\sum\limits_{j_4=0}^p\frac{1}{2}\left(\int\limits_{t_3}^{T}\phi_{j_4}(t_4)dt_4\right)^2
-\frac{T-t_3}{2}
\right)^2
dt_1 dt_2 dt_3.
\end{equation}

\vspace{5mm}

Further, applying the Parseval equality and the generalized Parseval equality
as well as using the 
Cauchy--Bunyakovsky inequality,
we have (see the proof of Theorem~43)

\vspace{-1mm}
\begin{equation}
\label{april46}
\sum\limits_{j=0}^{\infty}
\left(\int\limits_{t_1}^{t_2}\phi_{j}(s)ds\right)^2
=\int\limits_t^T \left({\bf 1}_{\{t_1<s<t_2\}}\right)^2 ds=
t_2-t_1,
\end{equation}

\vspace{3mm}
$$
\sum\limits_{j=0}^{\infty}\int\limits_{t_1}^{t_2}\phi_{j}(s)ds
\int\limits_{t_3}^{t_4}\phi_{j}(s)ds=
\sum\limits_{j=0}^{\infty}\int\limits_t^T {\bf 1}_{\{t_1<s<t_2\}}\phi_{j}(s)ds
\int\limits_t^T {\bf 1}_{\{t_3<s<t_4\}}\phi_{j}(s)ds=
$$

\begin{equation}
\label{april47}
=
\int\limits_t^T {\bf 1}_{\{t_1<s<t_2\}}{\bf 1}_{\{t_3<s<t_4\}}ds=0,
\end{equation}

\vspace{1mm}

\begin{equation}
\label{april48}
\left\vert
(t_2-t_1)
-\sum\limits_{j=0}^{p}
\left(\int\limits_{t_1}^{t_2}\phi_{j}(s)ds\right)^2
\right\vert\le
t_2-t_1\le T-t<\infty,
\end{equation}

\vspace{3mm}
$$
\left(\sum\limits_{j=0}^{p}\int\limits_{t_1}^{t_2}\phi_{j}(s)ds
\int\limits_{t_3}^{t_4}\phi_{j}(s)ds\right)^2\le
\sum\limits_{j=0}^{p}\left(\int\limits_{t_1}^{t_2}\phi_{j}(s)ds\right)^2
\sum\limits_{j=0}^{p}\left(\int\limits_{t_3}^{t_4}\phi_{j}(s)ds\right)^2\le
$$

\vspace{2mm}
\begin{equation}
\label{april49}
\le
(t_2-t_1)(t_4-t_3)\le (T-t)^2<\infty,
\end{equation}

\vspace{4mm}
\noindent
where $t\le t_1<t_2\le t_3<t_4\le T.$

Using Lebesgue's Dominated Convergence Theorem and (\ref{april46})--(\ref{april49}), 
we obtain that the right-hand sides of (\ref{april36})--(\ref{april45}) 
tend to zero when $p\to\infty.$
The equalities (\ref{april11})--(\ref{april20}) are proved.

\vspace{2mm}

{\bf Step~3.}\ Before proving the equalities (\ref{april21})--(\ref{april35}), we show that

\vspace{-1mm}
\begin{equation}
\label{april50}
\left\vert\sum\limits_{j_1, j_3=0}^p C_{j_3 j_3 j_1 j_1}(s,\tau)\right\vert\le K,
\end{equation}

\begin{equation}
\label{april51}
\left\vert\sum\limits_{j_1, j_3=0}^p C_{j_1 j_3 j_3 j_1}(s,\tau)\right\vert\le K,
\end{equation}

\begin{equation}
\label{april52}
\left\vert\sum\limits_{j_1, j_2=0}^p C_{j_2 j_1 j_2 j_1}(s,\tau)\right\vert\le K,
\end{equation}

\vspace{1mm}
\begin{equation}
\label{april53}
\sum\limits_{j_2=0}^p
\left(\sum\limits_{j_1=0}^p C_{j_1 j_2 j_1}(s,\tau)\right)^2\le
\int\limits_{\tau}^s \left(\sum\limits_{j_1=0}^p
\int\limits_{\tau}^{t_2}\phi_{j_1}(t_1)dt_1\int\limits_{t_2}^{s}\phi_{j_1}(t_3)dt_3\right)^2 dt_2,
\end{equation}

\vspace{3mm}
\noindent
where constant $K$ does not depend on $p, t_1, t_2;$
here and further in this proof

\vspace{-1mm}
$$
C_{j_k \ldots j_1}(s,\tau)=\int\limits_{\tau}^s
\phi_{j_k}(t_k)\ldots
\int\limits_{\tau}^{t_2}
\phi_{j_1}(t_1)dt_1\ldots dt_k \ \ \ (k=1,\ldots,4,\ t\le\tau<s\le T).
$$

\vspace{3mm}

Further, by $K, K_1, K_2$ we will denote contants
that can change from line to line.

By analogy with (\ref{march13}), (\ref{march21}), (\ref{march26}) 
and (\ref{march20a}), (\ref{march25}), (\ref{march34})
we get

\vspace{-1mm}
\begin{equation}
\label{april54}
\sum\limits_{j_1,j_3=0}^p C_{j_3 j_3 j_1 j_1}(s,\tau)=
\sum\limits_{j_1,j_3=0}^p C_{j_3}(s,\tau)C_{j_3 j_1 j_1}(s,\tau)-\frac{1}{8}
\left(\sum\limits_{j_1=0}^p \bigl(C_{j_1}(s,\tau)\bigr)^2\right)^2,
\end{equation}

\vspace{1mm}
\begin{equation}
\label{april55}
\sum\limits_{j_1,j_2=0}^p C_{j_2 j_1 j_2 j_1}(s,\tau)=
\sum\limits_{j_1,j_2=0}^p C_{j_2}(s,\tau)C_{j_1 j_2 j_1}(s,\tau)
-\frac{1}{2}\sum\limits_{j_1,j_2=0}^p C_{j_1 j_2}(s,\tau)C_{j_2 j_1}(s,\tau),
\end{equation}

\vspace{1mm}
\begin{equation}
\label{april56}
\sum\limits_{j_1,j_3=0}^p C_{j_1 j_3 j_3 j_1}(s,\tau)=
\sum\limits_{j_1,j_3=0}^p C_{j_1}(s,\tau)C_{j_3 j_3 j_1}(s,\tau)-\frac{1}{2}
\sum\limits_{j_1,j_3=0}^p \bigl(C_{j_3 j_1}(s,\tau)\bigr)^2,
\end{equation}

\vspace{1mm}
\begin{equation}
\label{april56x}
\lim\limits_{p\to\infty}
\sum\limits_{j_1, j_3=0}^{p}
C_{j_3 j_3 j_1 j_1}(s,\tau)=\frac{1}{8}(s-\tau)^2,
\biggr.
\end{equation}

\vspace{1mm}
\begin{equation}
\label{april56xx}
\lim\limits_{p\to\infty}
\sum\limits_{j_1, j_2=0}^{p}
C_{j_2 j_1 j_2 j_1}(s,\tau)=0
\biggr.,
\end{equation}

\vspace{1mm}
\begin{equation}
\label{april56xxx}
\lim\limits_{p\to\infty}
\sum\limits_{j_1, j_3=0}^{p}
C_{j_1 j_3 j_3 j_1}(s,\tau)=0
\biggr..
\end{equation}

\vspace{3mm}

Using (\ref{april54}), Parseval's equality,
Cauchy--Bunyakovsky's inequality, as well as Fubini's Theorem and the elementary inequality
$(a+b)^2\le 2 a^2 + 2 b^2,$
we obtain

\vspace{-2mm}
$$
\left(\sum\limits_{j_1,j_3=0}^p C_{j_3 j_3 j_1 j_1}(s,\tau)\right)^2\le
2\left(\sum\limits_{j_1,j_3=0}^p C_{j_3}(s,\tau)C_{j_3 j_1 j_1}(s,\tau)\right)^2+
2\cdot \frac{1}{64}
\left(\sum\limits_{j_1=0}^p \bigl(C_{j_1}(s,\tau)\bigr)^2\right)^4\le
$$

\vspace{2mm}
$$
\le 2\sum\limits_{j_3=0}^p \left(C_{j_3}(s,\tau)\right)^2
\sum\limits_{j_3=0}^p \left(\sum\limits_{j_1=0}^p
C_{j_3 j_1 j_1}(s,\tau)\right)^2+K_1\le
K_2
\sum\limits_{j_3=0}^{\infty} \left(\sum\limits_{j_1=0}^p
C_{j_3 j_1 j_1}(s,\tau)\right)^2+K_1=
$$

\vspace{2mm}
$$
= K_2
\sum\limits_{j_3=0}^{\infty} \left(\int\limits_{\tau}^s
\phi_{j_3}(t_3)\sum\limits_{j_1=0}^p
\int\limits_{\tau}^{t_3}\phi_{j_1}(t_2)
\int\limits_{\tau}^{t_2}\phi_{j_1}(t_1)dt_1 dt_2 dt_3
\right)^2+K_1=
$$

\vspace{2mm}
$$
= K_2
\int\limits_{\tau}^s \left(\frac{1}{2}\sum\limits_{j_1=0}^p
\left(\int\limits_{\tau}^{t_3}\phi_{j_1}(t_2)dt_2\right)^2\right)^2 dt_3
+K_1\le 
K_2
\int\limits_{\tau}^s \left(\frac{1}{2}\sum\limits_{j_1=0}^{\infty}
\left(\int\limits_{\tau}^{t_3}\phi_{j_1}(t_2)dt_2\right)^2\right)^2 dt_3
+K_1=
$$

\vspace{2mm}
$$
=
K_2
\int\limits_{\tau}^s \left(\frac{1}{2}(t_3-\tau)\right)^2 dt_3
+K_1\le K<\infty,
$$

\vspace{3mm}
\noindent
where constants $K, K_1, K_2$ do not depend on 
$p, s, \tau.$ The equality (\ref{april50}) is proved.

Let us prove (\ref{april51}). Using (\ref{april56}) and the
above reasoning, we get

\vspace{-2mm}
$$
\left(\sum\limits_{j_1,j_3=0}^p C_{j_1 j_3 j_3 j_1}(s,\tau)\right)^2\le
2\left(\sum\limits_{j_1,j_3=0}^p C_{j_1}(s,\tau)C_{j_3 j_3 j_1}(s,\tau)\right)^2+
2\cdot \frac{1}{4}\left(\sum\limits_{j_1,j_3=0}^p \bigl(C_{j_3 j_1}(s,\tau)\bigr)^2\right)^2\le
$$

\vspace{2mm}
$$
\le 2\sum\limits_{j_1=0}^p \left(C_{j_1}(s,\tau)\right)^2
\sum\limits_{j_1=0}^p \left(\sum\limits_{j_3=0}^p
C_{j_3 j_3 j_1}(s,\tau)\right)^2+K_1\le
K_2
\sum\limits_{j_1=0}^{\infty} \left(\sum\limits_{j_3=0}^p
C_{j_3 j_3 j_1}(s,\tau)\right)^2+K_1=
$$

\vspace{2mm}
$$
= K_2
\sum\limits_{j_1=0}^{\infty} \left(\int\limits_{\tau}^s
\phi_{j_1}(t_1)\sum\limits_{j_3=0}^p
\int\limits_{t_1}^{s}\phi_{j_3}(t_2)
\int\limits_{t_2}^{s}\phi_{j_3}(t_3)dt_3 dt_2 dt_1
\right)^2+K_1=
$$

\vspace{2mm}
$$
= K_2
\int\limits_{\tau}^s \left(\frac{1}{2}\sum\limits_{j_3=0}^p
\left(\int\limits_{t_1}^{s}\phi_{j_3}(t_2)dt_2\right)^2\right)^2 dt_1
+K_1\le 
K_2
\int\limits_{\tau}^s \left(\frac{1}{2}\sum\limits_{j_3=0}^{\infty}
\left(\int\limits_{t_1}^{s}\phi_{j_3}(t_2)dt_2\right)^2\right)^2 dt_1
+K_1=
$$

\vspace{2mm}
$$
=
K_2
\int\limits_{\tau}^s \left(\frac{1}{2}(s-t_1)\right)^2 dt_1
+K_1\le K<\infty,
$$

\vspace{3mm}
\noindent
where constants $K, K_1, K_2$ do not depend on 
$p, s, \tau.$ The equality (\ref{april51}) is proved.

Let us prove (\ref{april52}), (\ref{april53}). Applying (\ref{april55}), (\ref{april49}) and the
above reasoning, we have

\vspace{-2mm}
$$
\left(\sum\limits_{j_1,j_2=0}^p C_{j_2 j_1 j_2 j_1}(s,\tau)\right)^2 \hspace{-1.5mm}\le
2\left(\sum\limits_{j_1,j_2=0}^p C_{j_2}(s,\tau)C_{j_1 j_2 j_1}(s,\tau)\right)^2+
2\cdot \frac{1}{4}\left(\sum\limits_{j_1,j_2=0}^p C_{j_1 j_2}(s,\tau)C_{j_2 j_1}(s,\tau)\right)^2
\hspace{-1.5mm}\le
$$

\vspace{2mm}
$$
\le 2\sum\limits_{j_2=0}^p \left(C_{j_2}(s,\tau)\right)^2
\sum\limits_{j_2=0}^p \left(\sum\limits_{j_1=0}^p
C_{j_1 j_2 j_1}(s,\tau)\right)^2+
\frac{1}{2}\sum\limits_{j_1,j_2=0}^p \left(C_{j_1 j_2}(s,\tau)\right)^2
\sum\limits_{j_1,j_2=0}^p \left(C_{j_2 j_1}(s,\tau)\right)^2
\le
$$

\vspace{2mm}
\begin{equation}
\label{april60}
\le K_2
\sum\limits_{j_2=0}^{p} \left(\sum\limits_{j_1=0}^p
C_{j_1 j_2 j_1}(s,\tau)\right)^2+K_1\le
K_2
\sum\limits_{j_2=0}^{\infty} \left(\sum\limits_{j_1=0}^p
C_{j_1 j_2 j_1}(s,\tau)\right)^2+K_1=
\end{equation}

\vspace{2mm}
$$
= K_2
\sum\limits_{j_2=0}^{\infty} \left(\int\limits_{\tau}^s
\phi_{j_2}(t_2)\sum\limits_{j_1=0}^p
\int\limits_{\tau}^{t_2}\phi_{j_1}(t_1)dt_1
\int\limits_{t_2}^{s}\phi_{j_1}(t_3)dt_3 dt_2 
\right)^2+K_1=
$$

\vspace{2mm}
\begin{equation}
\label{april61}
= K_2
\int\limits_{\tau}^s \left(\sum\limits_{j_1=0}^p
\int\limits_{\tau}^{t_2}\phi_{j_1}(t_1)dt_1\int\limits_{t_2}^{s}\phi_{j_1}(t_3)dt_3\right)^2 dt_2
+K_1\le 
\end{equation}

\vspace{2mm}
$$
\le K_2
\int\limits_{\tau}^s \left((t_2-\tau)(s-t_2)\right)^2 dt_2
+K_1\le K<\infty,
$$

\vspace{3mm}
\noindent
where constants $K, K_1, K_2$ do not depend on 
$p, s, \tau.$ The equalities (\ref{april52}) and (\ref{april53}) (see (\ref{april60}), (\ref{april61}))
are proved. 

\vspace{2mm}

{\bf Step~4.}\ Let us start proving the equalities (\ref{april21})--(\ref{april35}).
Using Fubini's Theorem and Parseval's equality, we obtain
the following relations 
for the prelimit
expressions on the left-hand sides of (\ref{april21}), 
(\ref{april24}), (\ref{april27}),
(\ref{april33})--(\ref{april35})

\vspace{-1mm}
$$
\sum\limits_{j_5=0}^{p}
\left(\sum\limits_{j_1,j_3=0}^{p}
C_{j_5 j_3 j_3 j_1 j_1}-\frac{1}{4} 
C_{j_5 j_3 j_3 j_1 j_1}\biggl|_{(j_{1} j_{1})\curvearrowright (\cdot),
(j_{3} j_{3})\curvearrowright (\cdot)}
\biggr.\right)^2=
$$

\vspace{2mm}
$$
=\sum\limits_{j_5=0}^{p}
\left(\int\limits_t^T \phi_{j_5}(t_5) \left(\sum\limits_{j_1,j_3=0}^{p}
C_{j_3 j_3 j_1 j_1} (t_5,t)-\frac{1}{4}\int\limits_t^{t_5}(\tau-t)d\tau\right)dt_5\right)^2\le
$$

\vspace{2mm}
$$
\le\sum\limits_{j_5=0}^{\infty}
\left(\int\limits_t^T \phi_{j_5}(t_5) \left(\sum\limits_{j_1,j_3=0}^{p}
C_{j_3 j_3 j_1 j_1} (t_5,t)-\frac{1}{4}\int\limits_t^{t_5}(\tau-t)d\tau\right)dt_5\right)^2=
$$

\vspace{2mm}
\begin{equation}
\label{april62}
=\int\limits_t^T
\left(\sum\limits_{j_1,j_3=0}^{p}
C_{j_3 j_3 j_1 j_1} (t_5,t)-\frac{1}{8}(t_5-t)^2\right)^2 dt_5,
\end{equation}

\vspace{6mm}
$$
\sum\limits_{j_5=0}^{p}
\left(\sum\limits_{j_1,j_2=0}^{p}
C_{j_5 j_2 j_1 j_2 j_1}\right)^2
=\sum\limits_{j_5=0}^{p}
\left(\int\limits_t^T \phi_{j_5}(t_5) \sum\limits_{j_1,j_2=0}^{p}
C_{j_2 j_1 j_2 j_1} (t_5,t)dt_5\right)^2\le
$$

\vspace{2mm}
\begin{equation}
\label{april63}
\le\sum\limits_{j_5=0}^{\infty}
\left(\int\limits_t^T \phi_{j_5}(t_5) \sum\limits_{j_1,j_2=0}^{p}
C_{j_2 j_1 j_2 j_1} (t_5,t)dt_5\right)^2
=\int\limits_t^T
\left(\sum\limits_{j_1,j_2=0}^{p}
C_{j_2 j_1 j_2 j_1} (t_5,t)\right)^2 dt_5,
\end{equation}

\vspace{6mm}
$$
\sum\limits_{j_5=0}^{p}
\left(\sum\limits_{j_1,j_2=0}^{p}
C_{j_5 j_1 j_2 j_2 j_1}\right)^2=
\sum\limits_{j_5=0}^{p}
\left(\int\limits_t^T \phi_{j_5}(t_5) \sum\limits_{j_1,j_2=0}^{p}
C_{j_1 j_2 j_2 j_1} (t_5,t)dt_5\right)^2\le
$$

\vspace{2mm}
\begin{equation}
\label{april64}
\le\sum\limits_{j_5=0}^{\infty}
\left(\int\limits_t^T \phi_{j_5}(t_5) \sum\limits_{j_1,j_2=0}^{p}
C_{j_1 j_2 j_2 j_1} (t_5,t)dt_5\right)^2
=\int\limits_t^T
\left(\sum\limits_{j_1,j_2=0}^{p}
C_{j_1 j_2 j_2 j_1} (t_5,t)\right)^2 dt_5,
\end{equation}

\vspace{6mm}

$$
\sum\limits_{j_1=0}^{p}
\left(\sum\limits_{j_2,j_4=0}^{p}
C_{j_4 j_4 j_2 j_2 j_1}-\frac{1}{4} 
C_{j_4 j_4 j_2 j_2 j_1}\biggl|_{(j_{2} j_{2})\curvearrowright (\cdot),
(j_{4} j_{4})\curvearrowright (\cdot)}
\biggr.\right)^2=
$$

\vspace{2mm}
$$
=\sum\limits_{j_1=0}^{p}
\left(\int\limits_t^T  \phi_{j_1}(t_1)
\sum\limits_{j_2,j_4=0}^{p}
\int\limits_{t_1}^T  \phi_{j_2}(t_2)
\int\limits_{t_2}^T  \phi_{j_2}(t_3)
\int\limits_{t_3}^T  \phi_{j_4}(t_4)
\int\limits_{t_4}^T  \phi_{j_4}(t_5)
dt_5 dt_4 dt_3 dt_2 dt_1-\right.
$$

\vspace{2mm}
$$
\left.-\frac{1}{4} \int\limits_{t}^T\int\limits_{t}^{t_5} 
\int\limits_{t}^{t_3}\phi_{j_1}(t_1) dt_1 dt_3 dt_5\right)^2=
$$

\vspace{2mm}
$$
=\sum\limits_{j_1=0}^{p}
\left(\int\limits_t^T  \phi_{j_1}(t_1)
\left(\sum\limits_{j_2,j_4=0}^{p}
C_{j_4 j_4 j_2 j_2}(T,t_1)-
\frac{1}{4} \int\limits_{t_1}^T(T-t_3)dt_3
\right)dt_1\right)^2\le
$$

\vspace{2mm}
$$
\le\sum\limits_{j_1=0}^{\infty}
\left(\int\limits_t^T  \phi_{j_1}(t_1)
\left(\sum\limits_{j_2,j_4=0}^{p}
C_{j_4 j_4 j_2 j_2}(T,t_1)-
\frac{1}{8}(T-t_1)^2
\right)dt_1\right)^2=
$$

\vspace{2mm}
\begin{equation}
\label{april65}
=\int\limits_t^T
\left(\sum\limits_{j_2,j_4=0}^{p}
C_{j_4 j_4 j_2 j_2}(T,t_1)-
\frac{1}{8}(T-t_1)^2
\right)^2 dt_1,
\end{equation}

\vspace{6mm}
$$
\sum\limits_{j_1=0}^{p}
\left(\sum\limits_{j_2,j_3=0}^{p}
C_{j_3 j_2 j_3 j_2 j_1}\right)^2=
$$

\vspace{2mm}
$$
=\sum\limits_{j_1=0}^{p}
\left(\int\limits_t^T  \phi_{j_1}(t_1)
\sum\limits_{j_2,j_3=0}^{p}
\int\limits_{t_1}^T  \phi_{j_2}(t_2)
\int\limits_{t_2}^T  \phi_{j_3}(t_3)
\int\limits_{t_3}^T  \phi_{j_2}(t_4)
\int\limits_{t_4}^T  \phi_{j_3}(t_5)
dt_5 dt_4
dt_3 dt_2 dt_1\right)^2=
$$

\vspace{2mm}
$$
=\sum\limits_{j_1=0}^{p}
\left(\int\limits_t^T  \phi_{j_1}(t_1)
\sum\limits_{j_2,j_3=0}^{p}
C_{j_3 j_2 j_3 j_2}(T,t_1)dt_1\right)^2\le
$$

\vspace{2mm}
\begin{equation}
\label{april66}
\le\sum\limits_{j_1=0}^{\infty}
\left(\int\limits_t^T  \phi_{j_1}(t_1)
\sum\limits_{j_2,j_3=0}^{p}
C_{j_3 j_2 j_3 j_2}(T,t_1)dt_1\right)^2=
\int\limits_t^T
\left(\sum\limits_{j_2,j_3=0}^{p}
C_{j_3 j_2 j_3 j_2}(T,t_1)
\right)^2 dt_1,
\end{equation}

\vspace{6mm}

$$
\sum\limits_{j_1=0}^{p}
\left(\sum\limits_{j_2,j_3=0}^{p}
C_{j_2 j_3 j_3 j_2 j_1}\right)^2=
$$

\vspace{2mm}
$$
=\sum\limits_{j_1=0}^{p}
\left(\int\limits_t^T  \phi_{j_1}(t_1)
\sum\limits_{j_2,j_3=0}^{p}
\int\limits_{t_1}^T  \phi_{j_2}(t_2)
\int\limits_{t_2}^T  \phi_{j_3}(t_3)
\int\limits_{t_3}^T  \phi_{j_3}(t_4)
\int\limits_{t_4}^T  \phi_{j_2}(t_5)
dt_5 dt_4\
dt_3 dt_2 dt_1\right)^2=
$$

\vspace{2mm}
$$
=\sum\limits_{j_1=0}^{p}
\left(\int\limits_t^T  \phi_{j_1}(t_1)
\sum\limits_{j_2,j_3=0}^{p}
C_{j_2 j_3 j_3 j_2}(T,t_1)dt_1\right)^2\le
$$

\vspace{2mm}
\begin{equation}
\label{april67}
\le\sum\limits_{j_1=0}^{\infty}
\left(\int\limits_t^T  \phi_{j_1}(t_1)
\sum\limits_{j_2,j_3=0}^{p}
C_{j_2 j_3 j_3 j_2}(T,t_1)dt_1\right)^2=
\int\limits_t^T
\left(\sum\limits_{j_2,j_3=0}^{p}
C_{j_2 j_3 j_3 j_2}(T,t_1)
\right)^2 dt_1.
\end{equation}

\vspace{4mm}

Using Lebesgue's Dominated Convergence Theorem and (\ref{april50})--(\ref{april52}), 
(\ref{april56x})--(\ref{april56xxx}),
we obtain that the right-hand sides of (\ref{april62})--(\ref{april67}) 
tend to zero when $p\to\infty.$
The equalities (\ref{april21}), 
(\ref{april24}), (\ref{april27}),
(\ref{april33})--(\ref{april35}) are proved.

Further, let us prove the equalities (\ref{april23}), (\ref{april25}), 
(\ref{april28}), (\ref{april29}), (\ref{april31}).
Using Fubini's Theorem, Parseval's equality
and Cauchy--Bunyakovsky's inequality, we have
the following relations 
for the prelimit
expressions on the left-hand sides of 
(\ref{april23}), (\ref{april25}), 
(\ref{april28}), (\ref{april29}), (\ref{april31}) 

\vspace{-1mm}
$$
\sum\limits_{j_3=0}^{p}
\left(\sum\limits_{j_1,j_4=0}^{p}
C_{j_4 j_4 j_3 j_1 j_1}-\frac{1}{4} 
C_{j_4 j_4 j_3 j_1 j_1}\biggl|_{(j_{1} j_{1})\curvearrowright (\cdot),
(j_{4} j_{4})\curvearrowright (\cdot)}
\biggr.\right)^2=
$$

\vspace{2mm}
$$
=\sum\limits_{j_3=0}^{p}
\left(\int\limits_t^T  \phi_{j_3}(t_3)
\sum\limits_{j_1,j_4=0}^{p}
\int\limits_t^{t_3}  \phi_{j_1}(t_2)
\int\limits_t^{t_2}  \phi_{j_1}(t_1)dt_1 dt_2
\int\limits_{t_3}^{T}  \phi_{j_4}(t_4)
\int\limits_{t_4}^{T}  \phi_{j_4}(t_5)
dt_5 dt_4 dt_3-\right.
$$

\vspace{2mm}
$$
\left.-\frac{1}{4}
\int\limits_{t}^{T}\int\limits_{t}^{t_4}\phi_{j_3}(t_3)\int\limits_{t}^{t_3}dt_1
dt_3 dt_4\right)^2\le
$$

\vspace{2mm}
$$
\le\sum\limits_{j_3=0}^{\infty}
\left(\int\limits_t^T  \phi_{j_3}(t_3)
\left(\sum\limits_{j_1,j_4=0}^{p}
\frac{1}{4}\left(\int\limits_t^{t_3}  \phi_{j_1}(t_2)dt_2\right)^2
\left(\int\limits_{t_3}^{T}  \phi_{j_4}(t_4)dt_4\right)^2
-\frac{1}{4}(t_3-t)
\int\limits_{t_3}^{T}dt_4\right) dt_3\right)^2=
$$

\vspace{2mm}
\begin{equation}
\label{april100}
=\int\limits_t^T
\left(
\frac{1}{4}\sum\limits_{j_1=0}^{p}
\left(\int\limits_t^{t_3}  \phi_{j_1}(t_2)dt_2\right)^2
\sum\limits_{j_4=0}^{p}
\left(\int\limits_{t_3}^{T}  \phi_{j_4}(t_4)dt_4\right)^2
-\frac{1}{4}(t_3-t)(T-t_3)\right)^2 dt_3,
\end{equation}

\vspace{6mm}

$$
\sum\limits_{j_4=0}^{p}
\left(\sum\limits_{j_1,j_2=0}^{p}
C_{j_2 j_4 j_1 j_2 j_1}\right)^2=
$$

\vspace{2mm}
$$
=\sum\limits_{j_4=0}^{p}
\left(\int\limits_t^T \phi_{j_4}(t_4)
\sum\limits_{j_1,j_2=0}^{p}
\int\limits_t^{t_4} \phi_{j_1}(t_3)
\int\limits_t^{t_3}  \phi_{j_2}(t_2)
\int\limits_{t}^{t_2}  \phi_{j_1}(t_1)dt_1 dt_2 dt_3
\int\limits_{t_4}^{T}  \phi_{j_2}(t_5)
dt_5 dt_4\right)^2\le
$$

\vspace{2mm}
$$
\le\sum\limits_{j_4=0}^{\infty}
\left(\int\limits_t^T \phi_{j_4}(t_4)
\sum\limits_{j_1,j_2=0}^{p}C_{j_1 j_2 j_1}(t_4,t)
C_{j_2}(T,t_4) dt_4\right)^2=
$$

\vspace{2mm}
$$
=\int\limits_t^T
\left(
\sum\limits_{j_2=0}^{p} \sum\limits_{j_1=0}^{p}C_{j_1 j_2 j_1}(t_4,t)
C_{j_2}(T,t_4)\right)^2 dt_4\le
$$

\vspace{2mm}
$$
\le\int\limits_t^T
\sum\limits_{j_2=0}^{p}\left(C_{j_2}(T,t_4)\right)^2 
\sum\limits_{j_2=0}^{p}\left(\sum\limits_{j_1=0}^{p}C_{j_1 j_2 j_1}(t_4,t)
\right)^2 dt_4\le
$$

\vspace{2mm}
$$
\le\int\limits_t^T
\sum\limits_{j_2=0}^{\infty}\left(C_{j_2}(T,t_4)\right)^2 
\sum\limits_{j_2=0}^{p}\left(\sum\limits_{j_1=0}^{p}C_{j_1 j_2 j_1}(t_4,t)
\right)^2 dt_4\le
$$

\vspace{2mm}
\begin{equation}
\label{april69}
\le
K_1\int\limits_t^T
\sum\limits_{j_2=0}^{p}\left(\sum\limits_{j_1=0}^{p}C_{j_1 j_2 j_1}(t_4,t)
\right)^2 dt_4\le
\end{equation}

\vspace{1mm}
\begin{equation}
\label{april70}
\le
K_1\int\limits_t^T
\int\limits_{t}^{t_4} \left(\sum\limits_{j_1=0}^p
\int\limits_{t}^{t_2}\phi_{j_1}(t_1)dt_1\int\limits_{t_2}^{t_4}\phi_{j_1}(t_3)dt_3\right)^2 dt_2
dt_4=
\end{equation}

\vspace{1mm}
\begin{equation}
\label{april101}
=
K_1\int\limits_{[t,T]^2}{\bf 1}_{\{t_2<t_4\}}
\left(\sum\limits_{j_1=0}^p
\int\limits_{t}^{t_2}\phi_{j_1}(t_1)dt_1\int\limits_{t_2}^{t_4}\phi_{j_1}(t_3)dt_3\right)^2 dt_2
dt_4,
\end{equation}

\vspace{3mm}
\noindent
where constant $K_1$ does not depend on $p$
and the transition from (\ref{april69}) to (\ref{april70}) is based on 
(\ref{april53});

$$
\sum\limits_{j_3=0}^{p}
\left(\sum\limits_{j_1,j_2=0}^{p}
C_{j_2 j_1 j_3 j_2 j_1}\right)^2=
$$

\vspace{2mm}
$$
=\sum\limits_{j_3=0}^{p}
\left(\int\limits_t^T  \phi_{j_3}(t_3)
\sum\limits_{j_1,j_2=0}^{p}
\int\limits_t^{t_3}  \phi_{j_2}(t_2)
\int\limits_t^{t_2}  \phi_{j_1}(t_1)dt_1 dt_2
\int\limits_{t_3}^{T}  \phi_{j_1}(t_4)
\int\limits_{t_4}^{T}  \phi_{j_2}(t_5)dt_5 dt_4 dt_3
\right)^2
\le
$$

\vspace{2mm}
$$
\le\sum\limits_{j_3=0}^{\infty}
\left(\int\limits_t^T \phi_{j_3}(t_3)
\sum\limits_{j_1,j_2=0}^{p}
\int\limits_t^{t_3}  \phi_{j_2}(t_2)
\int\limits_t^{t_2}  \phi_{j_1}(t_1)dt_1 dt_2
\int\limits_{t_3}^{T} \phi_{j_1}(t_1)
\int\limits_{t_1}^{T}  \phi_{j_2}(t_2)dt_2 dt_1 dt_3
\right)^2
=
$$

\vspace{2mm}
$$
=\int\limits_t^T
\left(
\sum\limits_{j_1,j_2=0}^{p}
\int\limits_t^{t_3}\phi_{j_2}(t_2)
\int\limits_t^{t_2} \phi_{j_1}(t_1)dt_1 dt_2
\int\limits_{t_3}^{T} \phi_{j_1}(t_1)
\int\limits_{t_1}^{T}\phi_{j_2}(t_2)dt_2 dt_1\right)^2 dt_3=
$$

\begin{equation}
\label{april102}
=\int\limits_t^T
\left(
\sum\limits_{j_1,j_2=0}^{p}~
\int\limits_{[t,T]^2}{\bf 1}_{\{t_1<t_2<t_3\}}
\phi_{j_2}(t_2)\phi_{j_1}(t_1)dt_1 dt_2
\int\limits_{[t,T]^2}{\bf 1}_{\{t_2>t_1>t_3\}}
\phi_{j_2}(t_2)\phi_{j_1}(t_1)dt_1 dt_2
\right)^2 dt_3,
\end{equation}

\vspace{4mm}
\noindent
where, using the generalized Parseval equality and the 
Cauchy--Bunyakovsky inequality, we obtain

\vspace{1mm}
$$
\lim\limits_{p\to\infty}\sum\limits_{j_1,j_2=0}^{p}~
\int\limits_{[t,T]^2}{\bf 1}_{\{t_1<t_2<t_3\}}
\phi_{j_2}(t_2)\phi_{j_1}(t_1)dt_1 dt_2
\int\limits_{[t,T]^2}{\bf 1}_{\{t_2>t_1>t_3\}}
\phi_{j_2}(t_2)\phi_{j_1}(t_1)dt_1 dt_2=
$$

\vspace{2mm}
$$
=
\int\limits_{[t,T]^2}{\bf 1}_{\{t_1<t_2<t_3\}}{\bf 1}_{\{t_2>t_1>t_3\}}dt_1 dt_2=0,
$$

$$
\left(\sum\limits_{j_1,j_2=0}^{p}~
\int\limits_{[t,T]^2}{\bf 1}_{\{t_1<t_2<t_3\}}
\phi_{j_2}(t_2)\phi_{j_1}(t_1)dt_1 dt_2
\int\limits_{[t,T]^2}{\bf 1}_{\{t_2>t_1>t_3\}}
\phi_{j_2}(t_2)\phi_{j_1}(t_1)dt_1 dt_2\right)^2\le
$$
$$
\le
\sum\limits_{j_1,j_2=0}^{p}~\left(\int\limits_{[t,T]^2}{\bf 1}_{\{t_1<t_2<t_3\}}
\phi_{j_2}(t_2)\phi_{j_1}(t_1)dt_1 dt_2
\right)^2\times
$$
$$
\times
\sum\limits_{j_1,j_2=0}^{p}~\left(
\int\limits_{[t,T]^2}{\bf 1}_{\{t_2>t_1>t_3\}}
\phi_{j_2}(t_2)\phi_{j_1}(t_1)dt_1 dt_2\right)^2\le K_1<\infty,
$$

\vspace{6mm}
\noindent
where constant $K_1$ does not depend on $p;$

\vspace{2mm}
$$
\sum\limits_{j_2=0}^{p}
\left(\sum\limits_{j_1,j_3=0}^{p}
C_{j_3 j_1 j_3 j_2 j_1}\right)^2=
$$

\vspace{2mm}
$$
=\sum\limits_{j_2=0}^{p}
\left(\int\limits_t^T \phi_{j_2}(t_2)
\hspace{-1.2mm}\sum\limits_{j_1,j_3=0}^{p}
\int\limits_t^{t_2}  \phi_{j_1}(t_1)dt_1
\int\limits_{t_2}^T  \phi_{j_3}(t_3)
\int\limits_{t_3}^{T} \phi_{j_1}(t_4)
\int\limits_{t_4}^{T}  \phi_{j_3}(t_5)dt_5 dt_4 dt_3 dt_2
\right)^2
\le
$$

\vspace{2mm}
$$
\le\sum\limits_{j_2=0}^{\infty}
\left(\int\limits_t^T  \phi_{j_2}(t_2)
\hspace{-1.2mm}\sum\limits_{j_1,j_3=0}^{p}
\int\limits_t^{t_2} \phi_{j_1}(t_1)dt_1
\int\limits_{t_2}^T  \phi_{j_3}(t_3)
\int\limits_{t_3}^{T} \phi_{j_1}(t_4)
\int\limits_{t_4}^{T}  \phi_{j_3}(t_5)dt_5 dt_4 dt_3 dt_2
\right)^2
=
$$

\vspace{2mm}
$$
=
\int\limits_t^T 
\left(\sum\limits_{j_1,j_3=0}^{p}
\int\limits_t^{t_2}\phi_{j_1}(t_1)dt_1
\int\limits_{t_2}^T\phi_{j_3}(t_3)
\int\limits_{t_3}^{T}\phi_{j_1}(t_4)
\int\limits_{t_4}^{T}\phi_{j_3}(t_5)dt_5 dt_4 dt_3\right)^2 dt_2
=
$$

\vspace{2mm}
$$
=
\int\limits_t^T 
\left(\sum\limits_{j_1=0}^{p}
C_{j_1}(t_2,t)
\sum\limits_{j_3=0}^{p}
\int\limits_{t_2}^T\phi_{j_3}(t_5)
\int\limits_{t_2}^{t_5}\phi_{j_1}(t_4)
\int\limits_{t_2}^{t_4}\phi_{j_3}(t_3)dt_3 dt_4 dt_5\right)^2 dt_2
=
$$

\vspace{2mm}
$$
=
\int\limits_t^T 
\left(\sum\limits_{j_1=0}^{p}
C_{j_1}(t_2,t)
\sum\limits_{j_3=0}^{p}
C_{j_3 j_1 j_3}(T,t_2)
\right)^2 dt_2
\le 
$$

\vspace{2mm}
$$
\le
\int\limits_t^T 
\sum\limits_{j_1=0}^{p}
\left(C_{j_1}(t_2,t)\right)^2
\sum\limits_{j_1=0}^{p}\left(\sum\limits_{j_3=0}^{p}
C_{j_3 j_1 j_3}(T,t_2)
\right)^2 dt_2\le
$$

\vspace{2mm}
\begin{equation}
\label{april72}
\le
K_1\int\limits_t^T 
\sum\limits_{j_1=0}^{p}\left(\sum\limits_{j_3=0}^{p}
C_{j_3 j_1 j_3}(T,t_2)
\right)^2 dt_2\le
\end{equation}

\vspace{2mm}
\begin{equation}
\label{april73}
\le
K_1\int\limits_t^T
\int\limits_{t_2}^{T} \left(\sum\limits_{j_3=0}^p
\int\limits_{t_2}^{\theta}\phi_{j_3}(t_1)dt_1\int\limits_{\theta}^{T}\phi_{j_3}(t_3)dt_3\right)^2 d\theta
dt_2=
\end{equation}

\vspace{2mm}
\begin{equation}
\label{april103}
=
K_1\int\limits_{[t,T]^2}{\bf 1}_{\{t_2<\theta\}}
\left(\sum\limits_{j_3=0}^p
\int\limits_{t_2}^{\theta}\phi_{j_3}(t_1)dt_1\int\limits_{\theta}^{T}\phi_{j_3}(t_3)dt_3\right)^2 
d\theta dt_2,
\end{equation}

\vspace{4mm}
\noindent
where constant $K_1$ does not depend on $p$
and the transition from (\ref{april72}) to (\ref{april73}) is based on 
(\ref{april53});

\vspace{1mm}
$$
\lim\limits_{p\to\infty}
\sum\limits_{j_3=0}^{p}
\left(\sum\limits_{j_1,j_2=0}^{p}
C_{j_1 j_2 j_3 j_2 j_1}\right)^2=
$$

\vspace{2mm}
$$
=\sum\limits_{j_3=0}^{p}
\left(\int\limits_t^T  \phi_{j_3}(t_3)
\sum\limits_{j_1,j_2=0}^{p}
\int\limits_t^{t_3}  \phi_{j_2}(t_2)
\int\limits_t^{t_2}  \phi_{j_1}(t_1)dt_1 dt_2
\int\limits_{t_3}^{T}  \phi_{j_2}(t_4)
\int\limits_{t_4}^{T}  \phi_{j_1}(t_5)dt_5 dt_4 dt_3
\right)^2
\le
$$

\vspace{2mm}
$$
\le\sum\limits_{j_3=0}^{\infty}
\left(\int\limits_t^T \phi_{j_3}(t_3)
\sum\limits_{j_1,j_2=0}^{p}
\int\limits_t^{t_3}  \phi_{j_2}(t_2)
\int\limits_t^{t_2}  \phi_{j_1}(t_1)dt_1 dt_2
\int\limits_{t_3}^{T}  \phi_{j_2}(t_2)
\int\limits_{t_2}^{T}  \phi_{j_1}(t_1)dt_1 dt_2 dt_3
\right)^2
=
$$

\vspace{2mm}
$$
=\int\limits_t^T
\left(
\sum\limits_{j_1,j_2=0}^{p}
\int\limits_t^{t_3}\phi_{j_2}(t_2)
\int\limits_t^{t_2} \phi_{j_1}(t_1)dt_1 dt_2
\int\limits_{t_3}^{T} \phi_{j_2}(t_2)
\int\limits_{t_2}^{T}\phi_{j_1}(t_1)dt_1 dt_2\right)^2 dt_3=
$$

\vspace{2mm}
\begin{equation}
\label{april104}
=\int\limits_t^T
\left(
\sum\limits_{j_1,j_2=0}^{p}~
\int\limits_{[t,T]^2}{\bf 1}_{\{t_1<t_2<t_3\}}
\phi_{j_2}(t_2)\phi_{j_1}(t_1)dt_1 dt_2
\int\limits_{[t,T]^2}{\bf 1}_{\{t_1>t_2>t_3\}}
\phi_{j_2}(t_2)\phi_{j_1}(t_1)dt_1 dt_2
\right)^2 dt_3,
\end{equation}

\vspace{4mm}
\noindent
where, using the generalized Parseval equality and the 
Cauchy--Bunyakovsky inequality, we obtain

\vspace{1mm}

$$
\lim\limits_{p\to\infty}\sum\limits_{j_1,j_2=0}^{p}~
\int\limits_{[t,T]^2}{\bf 1}_{\{t_1<t_2<t_3\}}
\phi_{j_2}(t_2)\phi_{j_1}(t_1)dt_1 dt_2
\int\limits_{[t,T]^2}{\bf 1}_{\{t_1>t_2>t_3\}}
\phi_{j_2}(t_2)\phi_{j_1}(t_1)dt_1 dt_2=
$$

\vspace{1mm}
$$
=
\int\limits_{[t,T]^2}{\bf 1}_{\{t_1<t_2<t_3\}}{\bf 1}_{\{t_1>t_2>t_3\}}dt_1 dt_2=0,
$$

$$
\left(\sum\limits_{j_1,j_2=0}^{p}~
\int\limits_{[t,T]^2}{\bf 1}_{\{t_1<t_2<t_3\}}
\phi_{j_2}(t_2)\phi_{j_1}(t_1)dt_1 dt_2
\int\limits_{[t,T]^2}{\bf 1}_{\{t_1>t_2>t_3\}}
\phi_{j_2}(t_2)\phi_{j_1}(t_1)dt_1 dt_2\right)^2\le
$$
$$
\le
\sum\limits_{j_1,j_2=0}^{p}~\left(\int\limits_{[t,T]^2}{\bf 1}_{\{t_1<t_2<t_3\}}
\phi_{j_2}(t_2)\phi_{j_1}(t_1)dt_1 dt_2
\right)^2\times
$$
$$
\times
\sum\limits_{j_1,j_2=0}^{p}~\left(
\int\limits_{[t,T]^2}{\bf 1}_{\{t_1>t_2>t_3\}}
\phi_{j_2}(t_2)\phi_{j_1}(t_1)dt_1 dt_2\right)^2\le K_1<\infty,
$$

\vspace{4mm}
\noindent
where constant $K_1$ does not depend on $p;$

Using Lebesgue's Dominated Convergence Theorem,
we obtain that the right-hand sides of 
(\ref{april100}), (\ref{april101}), 
(\ref{april102}), (\ref{april103}), (\ref{april104})
tend to zero when $p\to\infty.$
The equalities (\ref{april23}), (\ref{april25}), 
(\ref{april28}), (\ref{april29}), (\ref{april31}) are proved.

\vspace{2mm}

{\bf Step~5.}\ Finally, let us prove the equalities 
(\ref{april22}), (\ref{april26}), 
(\ref{april30}), (\ref{april32}).
Using Parseval's equality,
Cauchy--Bunyakovsky's inequality, as well as Fubini's Theorem and the elementary inequality
$(a+b)^2\le 2 a^2 + 2 b^2,$
we obtain for the prelimit expression on the left-hand side of 
(\ref{april22})

\vspace{-1mm}
$$
\sum\limits_{j_4=0}^{p}
\left(\sum\limits_{j_1,j_3=0}^{p}
C_{j_3 j_4 j_3 j_1 j_1}\right)^2=
$$

\vspace{2mm}
$$
=\sum\limits_{j_4=0}^{p}
\left(\int\limits_t^T \phi_{j_4}(t_4)
\sum\limits_{j_1,j_3=0}^{p}
\int\limits_t^{t_4}  \phi_{j_3}(t_3)
\int\limits_t^{t_3}  \phi_{j_1}(t_2)
\int\limits_{t}^{t_2} \phi_{j_1}(t_1)dt_1 dt_2dt_3
\int\limits_{t_4}^{T}  \phi_{j_3}(t_5)dt_5 dt_4
\right)^2
\le
$$

\vspace{2mm}
$$
\le\sum\limits_{j_4=0}^{\infty}
\left(\int\limits_t^T  \phi_{j_4}(t_4)
\sum\limits_{j_1,j_3=0}^{p}
\int\limits_t^{t_4}  \phi_{j_3}(t_3)
\int\limits_t^{t_3}  \phi_{j_1}(t_2)
\int\limits_{t}^{t_2}  \phi_{j_1}(t_1)dt_1 dt_2dt_3
\int\limits_{t_4}^{T}  \phi_{j_3}(t_5)dt_5 dt_4
\right)^2
=
$$

\vspace{2mm}
$$
=
\int\limits_t^T  
\left(\sum\limits_{j_1,j_3=0}^{p}
\int\limits_t^{t_4}\phi_{j_3}(t_3)
\int\limits_t^{t_3}\phi_{j_1}(t_2)
\int\limits_{t}^{t_2}\phi_{j_1}(t_1)dt_1 dt_2dt_3
\int\limits_{t_4}^{T}\phi_{j_3}(t_5)dt_5 
\right)^2 dt_4=
$$

\vspace{2mm}
$$
=
\int\limits_t^T  
\left(\sum\limits_{j_3=0}^{p}
\int\limits_t^{t_4}\phi_{j_3}(t_3)
\left(\frac{1}{2}\sum\limits_{j_1=0}^{p}\left(
\int\limits_t^{t_3}\phi_{j_1}(t_2)dt_2\right)^2 \mp \frac{t_3-t}{2}\right)dt_3
\int\limits_{t_4}^{T}\phi_{j_3}(t_5)dt_5 
\right)^2 dt_4\le
$$

\vspace{2mm}
$$
\le 2
\int\limits_t^T 
\left(\sum\limits_{j_3=0}^{p}
\int\limits_t^{t_4}\phi_{j_3}(t_3)
\left(\frac{1}{2}\sum\limits_{j_1=0}^{p}\left(
\int\limits_t^{t_3}\phi_{j_1}(t_2)dt_2\right)^2 -\frac{t_3-t}{2}\right)dt_3
\int\limits_{t_4}^{T}\phi_{j_3}(t_5)dt_5 
\right)^2 dt_4+
$$

\vspace{2mm}
$$
+2
\int\limits_t^T  
\left(\sum\limits_{j_3=0}^{p}
\int\limits_t^{t_4}\phi_{j_3}(t_3)
\frac{t_3-t}{2}dt_3
\int\limits_{t_4}^{T}\phi_{j_3}(t_5)dt_5 
\right)^2 dt_4 \le
$$

\vspace{2mm}
$$
\le 2
\int\limits_t^T  
\sum\limits_{j_3=0}^{p} \left(C_{j_3}(T,t_4)\right)^2
\sum\limits_{j_3=0}^{p} \left(
\int\limits_t^{t_4}\phi_{j_3}(t_3)
\left(\frac{1}{2}\sum\limits_{j_1=0}^{p}\left(
\int\limits_t^{t_3}\phi_{j_1}(t_2)dt_2\right)^2 -\frac{t_3-t}{2}\right)dt_3
\right)^2 dt_4 
+\varepsilon_p \le
$$

\vspace{2mm}
$$
\le K_1
\int\limits_t^T  
\sum\limits_{j_3=0}^{p} \left(
\int\limits_t^{t_4}\phi_{j_3}(t_3)
\left(\frac{1}{2}\sum\limits_{j_1=0}^{p}\left(
\int\limits_t^{t_3}\phi_{j_1}(t_2)dt_2\right)^2 -\frac{t_3-t}{2}\right)dt_3
\right)^2 dt_4 + \varepsilon_p\le
$$

\vspace{2mm}
$$
\le K_1
\int\limits_t^T  
\sum\limits_{j_3=0}^{\infty} \left(
\int\limits_t^{t_4}\phi_{j_3}(t_3)
\left(\frac{1}{2}\sum\limits_{j_1=0}^{p}\left(
\int\limits_t^{t_3}\phi_{j_1}(t_2)dt_2\right)^2 -\frac{t_3-t}{2}\right)dt_3
\right)^2 dt_4 +\varepsilon_p=
$$

\vspace{2mm}
$$
= K_1
\int\limits_t^T  
\int\limits_t^{t_4}
\left(\frac{1}{2}\sum\limits_{j_1=0}^{p}\left(
\int\limits_t^{t_3}\phi_{j_1}(t_2)dt_2\right)^2 -\frac{t_3-t}{2}\right)^2 
dt_3dt_4 + \varepsilon_p=
$$

\begin{equation}
\label{april105}
= K_1
\int\limits_{[t,T]^2}
{\bf 1}_{\{t_3<t_4\}}
\left(\frac{1}{2}\sum\limits_{j_1=0}^{p}\left(
\int\limits_t^{t_3}\phi_{j_1}(t_2)dt_2\right)^2 -\frac{t_3-t}{2}\right)^2 
dt_3dt_4 + \varepsilon_p,
\end{equation}

\vspace{4mm}
\noindent
where constant $K_1$ does not depend on $p,$

\vspace{-1mm}
$$
\varepsilon_p=2\int\limits_t^T  
\left(\sum\limits_{j_3=0}^{p}
\int\limits_t^{t_4}\phi_{j_3}(t_3)
\frac{t_3-t}{2}dt_3
\int\limits_{t_4}^{T}\phi_{j_3}(t_5)dt_5 
\right)^2 dt_4.
$$

\vspace{4mm}

By analogy with (\ref{april47}), (\ref{april49}) we get

\vspace{-1mm}
\begin{equation}
\label{april106}
\left(\sum\limits_{j_3=0}^{p}
\int\limits_t^{t_4}\phi_{j_3}(t_3)
\frac{t_3-t}{2}dt_3
\int\limits_{t_4}^{T}\phi_{j_3}(t_5)dt_5 
\right)^2\le K_2<\infty,
\end{equation}

\vspace{1mm}
\begin{equation}
\label{april107}
\sum\limits_{j_3=0}^{\infty}
\int\limits_t^{t_4}\phi_{j_3}(t_3)
\frac{t_3-t}{2}dt_3
\int\limits_{t_4}^{T}\phi_{j_3}(t_5)dt_5=0,
\end{equation}

\vspace{4mm}
\noindent
where constant $K_2$ does not depend on $p.$

Using Lebesgue's Dominated Convergence Theorem and (\ref{april46}), (\ref{april48}), 
(\ref{april106}), (\ref{april107}), 
we obtain that the right-hand side of (\ref{april105})
tends to zero when $p\to\infty.$
The equality (\ref{april22}) is proved.

Let us prove the equality (\ref{april26}).
Using Parseval's equality,
Cauchy--Bunyakovsky's inequality, as well as Fubini's Theorem and the elementary inequality
$(a+b)^2\le 2 a^2 + 2 b^2,$
we obtain for the prelimit expression on the left-hand side of 
(\ref{april26})

\vspace{-1mm}
$$
\sum\limits_{j_2=0}^{p}
\left(\sum\limits_{j_1,j_4=0}^{p}
C_{j_4 j_4 j_1 j_2 j_1}\right)^2=
$$

\vspace{2mm}
$$
=\sum\limits_{j_2=0}^{p}
\left(\int\limits_t^T  \phi_{j_2}(t_2)
\sum\limits_{j_1,j_4=0}^{p}
\int\limits_t^{t_2}  \phi_{j_1}(t_1)dt_1
\int\limits_{t_2}^{T}  \phi_{j_1}(t_3)
\int\limits_{t_3}^{T}  \phi_{j_4}(t_4)
\int\limits_{t_4}^{T}  \phi_{j_4}(t_5)dt_5 dt_4 dt_3dt_2
\right)^2
\le
$$

\vspace{2mm}
$$
\le\sum\limits_{j_2=0}^{\infty}
\left(\int\limits_t^T  \phi_{j_2}(t_2)
\sum\limits_{j_1,j_4=0}^{p}
\int\limits_t^{t_2}  \phi_{j_1}(t_1)dt_1
\int\limits_{t_2}^{T}  \phi_{j_1}(t_3)
\int\limits_{t_3}^{T} \phi_{j_4}(t_4)
\int\limits_{t_4}^{T}  \phi_{j_4}(t_5)dt_5 dt_4 dt_3dt_2
\right)^2
=
$$

\vspace{2mm}
$$
=\int\limits_t^T \left(\sum\limits_{j_1,j_4=0}^{p}
\int\limits_t^{t_2} \phi_{j_1}(t_1)dt_1
\int\limits_{t_2}^{T} \phi_{j_1}(t_3)
\int\limits_{t_3}^{T} \phi_{j_4}(t_4)
\int\limits_{t_4}^{T}  \phi_{j_4}(t_5)dt_5 dt_4 dt_3\right)^2 dt_2
=
$$

\vspace{2mm}
$$
=
\int\limits_t^T  
\left(\sum\limits_{j_1=0}^{p}
\int\limits_t^{t_2}\phi_{j_1}(t_1)dt_1
\int\limits_{t_2}^{T}\phi_{j_1}(t_3)
\left(\frac{1}{2}\sum\limits_{j_4=0}^{p}\left(
\int\limits_{t_3}^T\phi_{j_4}(t_4)dt_4\right)^2 \mp \frac{T-t_3}{2}\right)
dt_3 
\right)^2 dt_2\le
$$

\vspace{2mm}
$$
\le 2
\int\limits_t^T  
\left(\sum\limits_{j_1=0}^{p}
\int\limits_t^{t_2}\phi_{j_1}(t_1)dt_1
\int\limits_{t_2}^{T}\phi_{j_1}(t_3)
\left(\frac{1}{2}\sum\limits_{j_4=0}^{p}\left(
\int\limits_{t_3}^T\phi_{j_4}(t_4)dt_4\right)^2 -\frac{T-t_3}{2}
\right)dt_3
\right)^2 dt_2+
$$

\vspace{2mm}
$$
+2
\int\limits_t^T  
\left(\sum\limits_{j_1=0}^{p}
\int\limits_t^{t_2}\phi_{j_1}(t_1)dt_1
\int\limits_{t_2}^{T}\phi_{j_1}(t_3)\frac{T-t_3}{2}dt_3
\right)^2 dt_2 \le
$$

\vspace{2mm}
$$
\le 2
\int\limits_t^T  
\sum\limits_{j_1=0}^{p}\left(C_{j_1}(t_2,t)\right)^2
\sum\limits_{j_1=0}^{p}\left(
\int\limits_{t_2}^{T}\phi_{j_1}(t_3)
\left(\frac{1}{2}\sum\limits_{j_4=0}^{p}\left(
\int\limits_{t_3}^T\phi_{j_4}(t_4)dt_4\right)^2 -\frac{T-t_3}{2}
\right)dt_3
\right)^2 dt_2+\mu_p\le
$$

\vspace{2mm}
$$
\le K_1
\int\limits_t^T  
\sum\limits_{j_1=0}^{p}\left(
\int\limits_{t_2}^{T}\phi_{j_1}(t_3)
\left(\frac{1}{2}\sum\limits_{j_4=0}^{p}\left(
\int\limits_{t_3}^T\phi_{j_4}(t_4)dt_4\right)^2 -\frac{T-t_3}{2}
\right)dt_3
\right)^2 dt_2+\mu_p\le
$$

\vspace{2mm}
$$
\le K_1
\int\limits_t^T  
\sum\limits_{j_1=0}^{\infty}\left(
\int\limits_{t_2}^{T}\phi_{j_1}(t_3)
\left(\frac{1}{2}\sum\limits_{j_4=0}^{p}\left(
\int\limits_{t_3}^T\phi_{j_4}(t_4)dt_4\right)^2 -\frac{T-t_3}{2}
\right)dt_3
\right)^2 dt_2+\mu_p=
$$

\vspace{2mm}
$$
= K_1
\int\limits_t^T  
\int\limits_{t_2}^{T}
\left(\frac{1}{2}\sum\limits_{j_4=0}^{p}\left(
\int\limits_{t_3}^T\hspace{-0.4mm}\phi_{j_4}(t_4)dt_4\right)^2 -\frac{T-t_3}{2}
\right)^2 dt_3 dt_2+\mu_p=
$$

\begin{equation}
\label{april108}
= K_1
\int\limits_{[t,T]^2}{\bf 1}_{\{t_2<t_3\}}
\left(\frac{1}{2}\sum\limits_{j_4=0}^{p}\left(
\int\limits_{t_3}^T\phi_{j_4}(t_4)dt_4\right)^2 -\frac{T-t_3}{2}
\right)^2 dt_3 dt_2+\mu_p,
\end{equation}

\vspace{4mm}
\noindent
where constant $K_1$ does not depend on $p,$

\vspace{-1mm}
$$
\mu_p=2
\int\limits_t^T  
\left(\sum\limits_{j_1=0}^{p}
\int\limits_t^{t_2}\phi_{j_1}(t_1)dt_1
\int\limits_{t_2}^{T}\phi_{j_1}(t_3)\frac{T-t_3}{2}dt_3
\right)^2 dt_2.
$$

\vspace{4mm}

By analogy with (\ref{april47}), (\ref{april49}) we get

\vspace{-1mm}
\begin{equation}
\label{april109}
\left(\sum\limits_{j_1=0}^{p}
\int\limits_t^{t_2}\phi_{j_1}(t_1)dt_1
\int\limits_{t_2}^{T}\phi_{j_1}(t_3)\frac{T-t_3}{2}dt_3
\right)^2\le K_2<\infty,
\end{equation}

\vspace{1mm}
\begin{equation}
\label{april110}
\sum\limits_{j_1=0}^{\infty}
\int\limits_t^{t_2}\phi_{j_1}(t_1)dt_1
\int\limits_{t_2}^{T}\phi_{j_1}(t_3)\frac{T-t_3}{2}dt_3=0,
\end{equation}

\vspace{4mm}
\noindent
where constant $K_2$ does not depend on $p.$

Using Lebesgue's Dominated Convergence Theorem and (\ref{april46}), (\ref{april48}), 
(\ref{april109}), (\ref{april110}), 
we obtain that the right-hand side of (\ref{april108})
tends to zero when $p\to\infty.$
The equality (\ref{april26}) is proved.

Let us prove the equality (\ref{april30}).
Using Parseval's equality,
Cauchy--Bunyakovsky's inequality, as well as Fubini's Theorem and the elementary inequality
$(a+b)^2\le 2 a^2 + 2 b^2,$
we obtain for the prelimit expression on the left-hand side of 
(\ref{april30})

\vspace{-1mm}
$$
\sum\limits_{j_4=0}^{p}
\left(\sum\limits_{j_1,j_2=0}^{p}
C_{j_1 j_4 j_2 j_2 j_1}\right)^2=
$$

\vspace{2mm}
$$
=\sum\limits_{j_4=0}^{p}
\left(\int\limits_t^T  \phi_{j_4}(t_4)
\sum\limits_{j_1,j_2=0}^{p}
\int\limits_t^{t_4}  \phi_{j_2}(t_3)
\int\limits_t^{t_3} \phi_{j_2}(t_2)
\int\limits_t^{t_2} \phi_{j_1}(t_1)dt_1 dt_2 dt_3
\int\limits_{t_4}^{T}  \phi_{j_1}(t_5)dt_5 dt_4 
\right)^2
\le
$$

\vspace{2mm}
$$
\le\sum\limits_{j_4=0}^{\infty}
\left(\int\limits_t^T  \phi_{j_4}(t_4)
\sum\limits_{j_1,j_2=0}^{p}
\int\limits_t^{t_4} \phi_{j_2}(t_3)
\int\limits_t^{t_3} \phi_{j_2}(t_2)
\int\limits_t^{t_2} \phi_{j_1}(t_1)dt_1 dt_2 dt_3
\int\limits_{t_4}^{T}  \phi_{j_1}(t_5)dt_5 dt_4 
\right)^2
=
$$

\vspace{2mm}
$$
=
\int\limits_t^T  
\left(\sum\limits_{j_1,j_2=0}^{p}
\int\limits_t^{t_4}\phi_{j_2}(t_3)
\int\limits_t^{t_3}\phi_{j_2}(t_2)
\int\limits_t^{t_2}\phi_{j_1}(t_1)dt_1 dt_2 dt_3
\int\limits_{t_4}^{T}\phi_{j_1}(t_5)dt_5 
\right)^2 dt_4 =
$$

\vspace{2mm}
$$
=
\int\limits_t^T  
\left(\sum\limits_{j_1,j_2=0}^{p}
\int\limits_t^{t_4}\phi_{j_1}(t_1)
\int\limits_{t_1}^{t_4}\phi_{j_2}(t_2)
\int\limits_{t_2}^{t_4}\phi_{j_2}(t_3)dt_3 dt_2 dt_1
\int\limits_{t_4}^{T}\phi_{j_1}(t_5)dt_5 
\right)^2 dt_4 =
$$

\vspace{2mm}
$$
=
\int\limits_t^T  
\left(\sum\limits_{j_1=0}^{p}
\int\limits_t^{t_4}\phi_{j_1}(t_1)
\left(\frac{1}{2}\sum\limits_{j_2=0}^{p}\left(
\int\limits_{t_1}^{t_4}\phi_{j_2}(t_2)dt_2\right)^2 \mp \frac{t_4-t_1}{2}\right)
dt_1\int\limits_{t_4}^T\phi_{j_1}(t_5)dt_5
\right)^2 dt_4\le
$$

\vspace{2mm}
$$
\le 2
\int\limits_t^T  
\left(\sum\limits_{j_1=0}^{p}
\int\limits_t^{t_4}\phi_{j_1}(t_1)
\left(\frac{1}{2}\sum\limits_{j_2=0}^{p}\left(
\int\limits_{t_1}^{t_4}
\phi_{j_2}(t_2)dt_2\right)^2 -\frac{t_4-t_1}{2}\right)
dt_1\int\limits_{t_4}^T\phi_{j_1}(t_5)dt_5
\right)^2 dt_4+
$$

\vspace{2mm}
$$
+2\int\limits_t^T  
\left(\sum\limits_{j_1=0}^{p}
\int\limits_t^{t_4}\phi_{j_1}(t_1)
\frac{t_4-t_1}{2}
dt_1\int\limits_{t_4}^T\phi_{j_1}(t_5)dt_5
\right)^2 dt_4\le
$$

\vspace{2mm}
$$
\le 2
\int\limits_t^T  
\sum\limits_{j_1=0}^{p}\left(C_{j_1}(T,t_4)\right)^2
\sum\limits_{j_1=0}^{p}\left(
\int\limits_t^{t_4}\phi_{j_1}(t_1)
\left(\frac{1}{2}\sum\limits_{j_2=0}^{p}\left(
\int\limits_{t_1}^{t_4}
\phi_{j_2}(t_2)dt_2\right)^2-\frac{t_4-t_1}{2}\right)
dt_1
\right)^2 dt_4+\rho_p\le
$$

\vspace{2mm}
$$
\le K_1
\int\limits_t^T  
\sum\limits_{j_1=0}^{p}\left(
\int\limits_t^{t_4}\phi_{j_1}(t_1)
\left(\frac{1}{2}\sum\limits_{j_2=0}^{p}\left(
\int\limits_{t_1}^{t_4}
\phi_{j_2}(t_2)dt_2\right)^2-\frac{t_4-t_1}{2}\right)
dt_1
\right)^2 dt_4+\rho_p\le
$$

\vspace{2mm}
$$
\le K_1
\int\limits_t^T  
\sum\limits_{j_1=0}^{\infty}\left(
\int\limits_t^{t_4}\phi_{j_1}(t_1)
\left(\frac{1}{2}\sum\limits_{j_2=0}^{p}\left(
\int\limits_{t_1}^{t_4}
\phi_{j_2}(t_2)dt_2\right)^2-\frac{t_4-t_1}{2}\right)
dt_1
\right)^2 dt_4+\rho_p=
$$

\vspace{2mm}
$$
=K_1
\int\limits_t^T  
\int\limits_t^{t_4}
\left(\frac{1}{2}\sum\limits_{j_2=0}^{p}\left(
\int\limits_{t_1}^{t_4}
\phi_{j_2}(t_2)dt_2\right)^2-\frac{t_4-t_1}{2}
\right)^2 dt_1 dt_4+\rho_p=
$$

\begin{equation}
\label{april111}
=K_1
\int\limits_{[t,T]^2}
{\bf 1}_{\{t_1<t_4\}}
\left(\frac{1}{2}\sum\limits_{j_2=0}^{p}\left(
\int\limits_{t_1}^{t_4}
\phi_{j_2}(t_2)dt_2\right)^2-\frac{t_4-t_1}{2}
\right)^2 dt_1 dt_4+\rho_p,
\end{equation}

\vspace{5mm}
\noindent
where constant $K_1$ does not depend on $p,$

\vspace{-1mm}
$$
\rho_p=
2\int\limits_t^T  
\left(\sum\limits_{j_1=0}^{p}
\int\limits_t^{t_4}\phi_{j_1}(t_1)
\frac{t_4-t_1}{2}
dt_1\int\limits_{t_4}^T\phi_{j_1}(t_5)dt_5
\right)^2 dt_4.
$$

\vspace{4mm}

By analogy with (\ref{april47}), (\ref{april49}) we get
$(t_4-t_1=(t_4-t)+(t-t_1))$

\vspace{-1mm}
\begin{equation}
\label{april112}
\left(\sum\limits_{j_1=0}^{p}
\int\limits_t^{t_4}\phi_{j_1}(t_1)
\frac{t_4-t_1}{2}
dt_1\int\limits_{t_4}^T\phi_{j_1}(t_5)dt_5
\right)^2\le K_2<\infty,
\end{equation}

\begin{equation}
\label{april113}
\sum\limits_{j_1=0}^{\infty}
\int\limits_t^{t_4}\phi_{j_1}(t_1)
\frac{t_4-t_1}{2}
dt_1\int\limits_{t_4}^T\phi_{j_1}(t_5)dt_5=0,
\end{equation}

\vspace{4mm}
\noindent
where constant $K_2$ does not depend on $p.$

Using Lebesgue's Dominated Convergence Theorem and (\ref{april46}), (\ref{april48}), 
(\ref{april112}), (\ref{april113}), 
we obtain that the right-hand side of (\ref{april111})
tends to zero when $p\to\infty.$
The equality (\ref{april30}) is proved.

Let us prove the equality (\ref{april32}).
Using Parseval's equality,
Cauchy--Bunyakovsky's inequality, as well as Fubini's Theorem and the elementary inequality
$(a+b)^2\le 2 a^2 + 2 b^2,$
we obtain for the prelimit expression on the left-hand side of 
(\ref{april32})

\vspace{-1mm}
$$
\sum\limits_{j_2=0}^{p}
\left(\sum\limits_{j_1,j_3=0}^{p}
C_{j_1 j_3 j_3 j_2 j_1}\right)^2=
$$

\vspace{2mm}
$$
=\sum\limits_{j_2=0}^{p}
\left(\int\limits_t^T  \phi_{j_2}(t_2)
\sum\limits_{j_1,j_3=0}^{p}
\int\limits_t^{t_2}  \phi_{j_1}(t_1)dt_1
\int\limits_{t_2}^T \phi_{j_3}(t_3)
\int\limits_{t_3}^T  \phi_{j_3}(t_4)
\int\limits_{t_4}^{T}  \phi_{j_1}(t_5)dt_5 dt_4 dt_3 dt_2
\right)^2
\le
$$

\vspace{2mm}
$$
\le\sum\limits_{j_2=0}^{\infty}
\left(\int\limits_t^T  \phi_{j_2}(t_2)
\sum\limits_{j_1,j_3=0}^{p}
\int\limits_t^{t_2}  \phi_{j_1}(t_1)dt_1
\int\limits_{t_2}^T \phi_{j_3}(t_3)
\int\limits_{t_3}^T  \phi_{j_3}(t_4)
\int\limits_{t_4}^{T}  \phi_{j_1}(t_5)dt_5 dt_4 dt_3 dt_2
\right)^2
=
$$

\vspace{2mm}
$$
=\int\limits_t^T  
\left(\sum\limits_{j_1,j_3=0}^{p}
\int\limits_t^{t_2}\phi_{j_1}(t_1)dt_1
\int\limits_{t_2}^T\phi_{j_3}(t_3)
\int\limits_{t_3}^T \phi_{j_3}(t_4)
\int\limits_{t_4}^{T} \phi_{j_1}(t_5)dt_5 dt_4 dt_3 
\right)^2 dt_2=
$$

\vspace{2mm}
$$
=\int\limits_t^T  
\left(\sum\limits_{j_1,j_3=0}^{p}
\int\limits_t^{t_2}\phi_{j_1}(t_1)dt_1
\int\limits_{t_2}^T\phi_{j_1}(t_5)
\int\limits_{t_2}^{t_5} \phi_{j_3}(t_4)
\int\limits_{t_2}^{t_4} \phi_{j_3}(t_3)dt_3 dt_4 dt_5 
\right)^2 dt_2=
$$

\vspace{2mm}
$$
=
\int\limits_t^T  
\left(\sum\limits_{j_1=0}^{p}
\int\limits_t^{t_2}\phi_{j_1}(t_1)dt_1\int\limits_{t_2}^T\phi_{j_1}(t_5)
\left(\frac{1}{2}\sum\limits_{j_3=0}^{p}\left(
\int\limits_{t_2}^{t_5}\phi_{j_3}(t_4)dt_4\right)^2 \mp \frac{t_5-t_2}{2}\right)
dt_5
\right)^2 dt_2\le
$$

\vspace{2mm}
$$
\le 2
\int\limits_t^T  
\left(\sum\limits_{j_1=0}^{p}
\int\limits_t^{t_2}\phi_{j_1}(t_1)dt_1\int\limits_{t_2}^T\phi_{j_1}(t_5)
\left(\frac{1}{2}\sum\limits_{j_3=0}^{p}\left(
\int\limits_{t_2}^{t_5}
\phi_{j_3}(t_4)dt_4\right)^2 - \frac{t_5-t_2}{2}\right)
dt_5
\right)^2 dt_2+
$$

\vspace{2mm}
$$
+2\int\limits_t^T  
\left(\sum\limits_{j_1=0}^{p}
\int\limits_t^{t_2}\phi_{j_1}(t_1)dt_1\int\limits_{t_2}^T\phi_{j_1}(t_5)
\frac{t_5-t_2}{2}
dt_5
\right)^2 dt_2\le
$$

\vspace{2mm}
$$
\le 
2\int\limits_t^T  
\sum\limits_{j_1=0}^{p}\left(C_{j_1}(t_2,t)\right)^2
\sum\limits_{j_1=0}^{p}\left(
\int\limits_{t_2}^T\phi_{j_1}(t_5)
\left(\frac{1}{2}\sum\limits_{j_3=0}^{p}\left(
\int\limits_{t_2}^{t_5}
\phi_{j_3}(t_4)dt_4\right)^2 - \frac{t_5-t_2}{2}\right)
dt_5
\right)^2 dt_2+ \chi_p\le
$$

\vspace{2mm}
$$
\le K_1
\int\limits_t^T  
\sum\limits_{j_1=0}^{p}\left(
\int\limits_{t_2}^T\phi_{j_1}(t_5)
\left(\frac{1}{2}\sum\limits_{j_3=0}^{p}\left(
\int\limits_{t_2}^{t_5}
\phi_{j_3}(t_4)dt_4\right)^2 - \frac{t_5-t_2}{2}\right)
dt_5
\right)^2 dt_2+ \chi_p\le
$$

\vspace{2mm}
$$
\le K_1
\int\limits_t^T  
\sum\limits_{j_1=0}^{\infty}\left(
\int\limits_{t_2}^T\phi_{j_1}(t_5)
\left(\frac{1}{2}\sum\limits_{j_3=0}^{p}\left(
\int\limits_{t_2}^{t_5}
\phi_{j_3}(t_4)dt_4\right)^2 - \frac{t_5-t_2}{2}\right)
dt_5
\right)^2 dt_2+ \chi_p=
$$

\vspace{2mm}
$$
= K_1
\int\limits_t^T  
\int\limits_{t_2}^T
\left(\frac{1}{2}\sum\limits_{j_3=0}^{p}\left(
\int\limits_{t_2}^{t_5}
\phi_{j_3}(t_4)dt_4\right)^2 - \frac{t_5-t_2}{2}
\right)^2 dt_5 dt_2+ \chi_p=
$$

\begin{equation}
\label{april114}
= K_1
\int\limits_{[t,T]^2}
{\bf 1}_{\{t_2<t_5\}}
\left(\frac{1}{2}\sum\limits_{j_3=0}^{p}\left(
\int\limits_{t_2}^{t_5}
\phi_{j_3}(t_4)dt_4\right)^2 - \frac{t_5-t_2}{2}
\right)^2 dt_5 dt_2+ \chi_p,
\end{equation}

\vspace{4mm}
\noindent
where constant $K_1$ does not depend on $p,$

\vspace{-1mm}
$$
\chi_p=
2\int\limits_t^T  
\left(\sum\limits_{j_1=0}^{p}
\int\limits_t^{t_2}\phi_{j_1}(t_1)dt_1\int\limits_{t_2}^T\phi_{j_1}(t_5)
\frac{t_5-t_2}{2}
dt_5
\right)^2 dt_2.
$$

\vspace{4mm}

By analogy with (\ref{april47}), (\ref{april49}) we get
$(t_5-t_2=(t_5-t)+(t-t_2))$

\vspace{-1mm}
\begin{equation}
\label{april115}
\left(\sum\limits_{j_1=0}^{p}
\int\limits_t^{t_2}\phi_{j_1}(t_1)dt_1\int\limits_{t_2}^T\phi_{j_1}(t_5)
\frac{t_5-t_2}{2}
dt_5
\right)^2\le K_2<\infty,
\end{equation}

\begin{equation}
\label{april116}
\sum\limits_{j_1=0}^{\infty}
\int\limits_t^{t_2}\phi_{j_1}(t_1)dt_1\int\limits_{t_2}^T\phi_{j_1}(t_5)
\frac{t_5-t_2}{2}
dt_5
=0,
\end{equation}

\vspace{4mm}
\noindent
where constant $K_2$ does not depend on $p.$

Using Lebesgue's Dominated Convergence Theorem and (\ref{april46}), (\ref{april48}), 
(\ref{april115}), (\ref{april116}), 
we obtain that the right-hand side of (\ref{april114})
tends to zero when $p\to\infty.$
The equality (\ref{april32}) is proved.
The equalities (\ref{april11})--(\ref{april35})
are proved.
Theorem~48 is proved.

\vspace{5mm}

\section{Expansion of Iterated Stratonovich Stochastic Integrals
of Multiplicity 3. The Case of an Ar\-bit\-ra\-ry Complete Orthonormal System of 
Functions in the Space $L_2([t,T])$ and 
Binomial Weight Functions}

\vspace{5mm}                                                                   

In this section, we will consider a generalization of Theorems~43, 45.
Namely, we will prove the following theorem.

\vspace{2mm}

{\bf Theorem~49}\ \cite{2018a}.\ {\it Suppose that
$\{\phi_j(x)\}_{j=0}^{\infty}$ is an arbitrary complete orthonormal system of 
func\-ti\-ons in the space $L_2([t,T]).$
Then$,$ for the iterated Stra\-to\-no\-vich stochastic integral
of third multiplicity 

\vspace{-1mm}
\begin{equation}
\label{may99}
I_{{l_1l_2l_3}_{T,t}}^{*(i_1i_2i_3)}={\int\limits_t^{*}}^T (t_3-t)^{l_3}
{\int\limits_t^{*}}^{t_3}(t_2-t)^{l_2}
{\int\limits_t^{*}}^{t_2}(t_1-t)^{l_1}
d{\bf w}_{t_1}^{(i_1)}
d{\bf w}_{t_2}^{(i_2)}d{\bf w}_{t_3}^{(i_3)}
\end{equation}

\vspace{3mm}
\noindent
the following expansion 
\begin{equation}
\label{may100}
I_{{l_1l_2l_3}_{T,t}}^{*(i_1i_2i_3)}=
\hbox{\vtop{\offinterlineskip\halign{
\hfil#\hfil\cr
{\rm l.i.m.}\cr
$\stackrel{}{{}_{p\to \infty}}$\cr
}} }\sum_{j_1,j_2,j_3=0}^{p}
C_{j_3 j_2 j_1}\zeta_{j_1}^{(i_1)}\zeta_{j_2}^{(i_2)}\zeta_{j_3}^{(i_3)}
\end{equation}

\vspace{3mm}
\noindent
that converges in the mean-square sense is valid, where 
$i_1,i_2,i_3=0,1,\ldots,m;$ $l_1,l_2,l_3=0,1,2,\ldots,$

\vspace{-1mm}
$$
C_{j_3 j_2 j_1}=\int\limits_t^T
(t_3-t)^{l_3}\phi_{j_3}(t_3)\int\limits_t^{t_3}
(t_2-t)^{l_2}
\phi_{j_2}(t_2)
\int\limits_t^{t_2}
(t_1-t)^{l_1}\phi_{j_1}(t_1)dt_1dt_2dt_3
$$
and
$$
\zeta_{j}^{(i)}=
\int\limits_t^T \phi_{j}(\tau) d{\bf w}_{\tau}^{(i)}
$$ 

\vspace{2mm}
\noindent
are independent standard Gaussian random variables for various 
$i$ or $j$ {\rm (}in the case when $i\ne 0${\rm ),}
${\bf w}_{\tau}^{(i)}={\bf f}_{\tau}^{(i)}$ for
$i=1,\ldots,m$ and 
${\bf w}_{\tau}^{(0)}=\tau.$}

\vspace{2mm}

Note that the iterated Stratonovich stochastic integrals (\ref{may99}) are important 
for applications (see Chapter~4 in \cite{2018a}).

{\bf Proof.}\ According to Theorems~47 and 5, we come to the conclusion that 
Theorem~49 will be proved if we prove the following
equalities

\vspace{-1mm}
\begin{equation}
\label{may101}
\lim\limits_{p\to\infty}
\sum\limits_{j_3=0}^{p}
\left(
\frac{1}{2} 
C_{j_3 j_1 j_1}\biggl|_{(j_{1} j_{1})\curvearrowright (\cdot)}
\biggr. - \sum\limits_{j_1=0}^{p} C_{j_3 j_1 j_1}
\right)^2=0,
\end{equation}

\vspace{1mm}
\begin{equation}
\label{may102}
\lim\limits_{p\to\infty}
\sum\limits_{j_1=0}^{p}
\left(
\frac{1}{2} 
C_{j_2 j_2 j_1}\biggl|_{(j_{2} j_{2})\curvearrowright (\cdot)}
\biggr. - \sum\limits_{j_2=0}^{p}  C_{j_2 j_2 j_1}
\right)^2=0,
\end{equation}

\vspace{1mm}
\begin{equation}
\label{may103}
\lim\limits_{p\to\infty}
\sum\limits_{j_2=0}^{p}
\left(~\sum\limits_{j_1=0}^{p} 
C_{j_1 j_2 j_1}\right)^2=0.
\end{equation}

\vspace{4mm}

First, we prove that

\vspace{-1mm}
\begin{equation}
\label{may103x}
\left\vert\sum\limits_{j=0}^{p}\int\limits_{t_1}^{t_2}(s-t)^{l}\phi_{j}(s)
\int\limits_{t_1}^{s}(\tau-t)^{m}\phi_{j}(\tau)d\tau ds\right\vert\le K<\infty,
\end{equation}

\vspace{3mm}
\noindent
where $l, m=0,1,2,\ldots,$ $t\le t_1<t_2\le T,$ constant $K$ does not depend on $p, t_1, t_2.$

Using Fubini's Theorem and Parseval's equality, we have for $m>l$ $(l, m=0,1,2,\ldots)$

\vspace{-1mm}
$$
\sum\limits_{j=0}^{p}\int\limits_{t}^{t_2}(s-t)^{l}\phi_{j}(s)
\int\limits_{t}^{s}(\tau-t)^{m}\phi_{j}(\tau)d\tau ds=
$$

\vspace{2mm}
$$
=\sum\limits_{j=0}^{p}\int\limits_{t}^{t_2}(s-t)^{l}\phi_{j}(s)
\int\limits_{t}^{s}(\tau-t)^{l} (\tau-t)^{m-l}\phi_{j}(\tau)d\tau ds=
$$

\vspace{2mm}
$$
=\sum\limits_{j=0}^{p}\int\limits_{t}^{t_2}(s-t)^{l}\phi_{j}(s)
\int\limits_{t}^{s}(\tau-t)^{l}\phi_{j}(\tau)\int\limits_t^{\tau} (\theta-t)^{m-l-1}(m-l)d\theta
d\tau ds=
$$

\vspace{2mm}
$$
=(m-l)\sum\limits_{j=0}^{p}\int\limits_{t}^{t_2}(\theta-t)^{m-l-1}
\int\limits_{\theta}^{t_2}(\tau-t)^{l}\phi_{j}(\tau)\int\limits_{\tau}^{t_2} 
(s-t)^{l}\phi_j(s)ds d\tau d\theta=
$$

\vspace{2mm}
$$
=(m-l)\int\limits_{t}^{t_2}(\theta-t)^{m-l-1}
\frac{1}{2}\sum\limits_{j=0}^{p}\left(\int\limits_{\theta}^{t_2}(\tau-t)^{l}\phi_{j}(\tau)d\tau\right)^2
d\theta\le
$$

\vspace{2mm}
$$
\le\frac{m-l}{2}\int\limits_{t}^{t_2}(\theta-t)^{m-l-1}
\sum\limits_{j=0}^{\infty}\left(\int\limits_{\theta}^{t_2}(\tau-t)^{l}\phi_{j}(\tau)d\tau\right)^2
d\theta=
$$

\vspace{2mm}
\begin{equation}
\label{may104}
=\frac{m-l}{2}\int\limits_{t}^{t_2}(\theta-t)^{m-l-1}
\int\limits_{\theta}^{t_2}(\tau-t)^{2l}d\tau
d\theta\le K_1 < \infty,
\end{equation}

\vspace{4mm}
\noindent
where constant $K_1$ does not depend on $p, t_2.$

For $l>m$ $(l, m=0,1,2,\ldots)$ we get

\vspace{-1mm}
$$
\sum\limits_{j=0}^{p}\int\limits_{t}^{t_2}(s-t)^{l}\phi_{j}(s)
\int\limits_{t}^{s}(\tau-t)^{m}\phi_{j}(\tau)d\tau ds=
$$

\vspace{2mm}
$$
=\sum\limits_{j=0}^{p}\int\limits_{t}^{t_2}(s-t)^{l}\phi_{j}(s)ds
\int\limits_{t}^{t_2}(\tau-t)^{m}\phi_{j}(\tau)d\tau-
$$

\vspace{2mm}
$$
-\sum\limits_{j=0}^{p}\int\limits_{t}^{t_2}(s-t)^{l}\phi_{j}(s)
\int\limits_{s}^{t_2}(\tau-t)^{m}\phi_{j}(\tau)d\tau ds=
$$

\vspace{2mm}
$$
=\sum\limits_{j=0}^{p}\int\limits_{t}^{t_2}(s-t)^{l}\phi_{j}(s)ds
\int\limits_{t}^{t_2}(\tau-t)^{m}\phi_{j}(\tau)d\tau-
$$

\vspace{2mm}
\begin{equation}
\label{may105}
-\sum\limits_{j=0}^{p}\int\limits_{t}^{t_2}(\tau-t)^{m}\phi_{j}(\tau)
\int\limits_{t}^{\tau}(s-t)^{l}\phi_{j}(s)ds d\tau.
\end{equation}

\vspace{4mm}

Applying 
Cauchy--Bunyakovsky's inequality  and Parseval's equality, we obtain

\vspace{-1mm}
$$
\left(\sum\limits_{j=0}^{p}\int\limits_{t}^{t_2}(s-t)^{l}\phi_{j}(s)ds
\int\limits_{t}^{t_2}(\tau-t)^{m}\phi_{j}(\tau)d\tau\right)^2\le
$$

\vspace{2mm}
$$
\le \sum\limits_{j=0}^{p}\left(\int\limits_{t}^{t_2}(s-t)^{l}\phi_{j}(s)ds\right)^2
\sum\limits_{j=0}^{p}\left(\int\limits_{t}^{t_2}(\tau-t)^{m}\phi_{j}(\tau)d\tau\right)^2\le
$$

\vspace{2mm}
$$
\le \sum\limits_{j=0}^{\infty}\left(\int\limits_{t}^{t_2}(s-t)^{l}\phi_{j}(s)ds\right)^2
\sum\limits_{j=0}^{\infty}\left(\int\limits_{t}^{t_2}(\tau-t)^{m}\phi_{j}(\tau)d\tau\right)^2=
$$

\vspace{2mm}
\begin{equation}
\label{may106}
=\int\limits_{t}^{t_2}(s-t)^{2l}ds
\int\limits_{t}^{t_2}(\tau-t)^{2m}d\tau\le K_2 < \infty,
\end{equation}

\vspace{4mm}
\noindent
where constant $K_2$ does not depend on $p, t_2.$

Using (\ref{may104})--(\ref{may106}), we obtain

\vspace{-1mm}
\begin{equation}
\label{may107}
\left\vert
\sum\limits_{j=0}^{p}\int\limits_{t}^{t_2}(s-t)^{l}\phi_{j}(s)
\int\limits_{t}^{s}(\tau-t)^{m}\phi_{j}(\tau)d\tau ds
\right\vert\le K_3<\infty,
\end{equation}

\vspace{4mm}
\noindent
where $l>m$ $(l, m=0,1,2,\ldots),$ constant $K_3$ does not depend on $p, t_2.$

For the case $l=m$ we get

\vspace{-1mm}
$$
\sum\limits_{j=0}^{p}\int\limits_{t}^{t_2}(s-t)^{l}\phi_{j}(s)
\int\limits_{t}^{s}(\tau-t)^{l}\phi_{j}(\tau)d\tau ds=
$$

\vspace{2mm}
$$
=\sum\limits_{j=0}^{p}\frac{1}{2}\left(\int\limits_{t}^{t_2}(s-t)^{l}\phi_{j}(s)ds\right)^2\le
\sum\limits_{j=0}^{\infty}\frac{1}{2}\left(\int\limits_{t}^{t_2}(s-t)^{l}\phi_{j}(s)ds\right)^2=
$$

\vspace{2mm}
\begin{equation}
\label{may108}
=\frac{1}{2}\int\limits_{t}^{t_2}(s-t)^{2l}ds\le K_4<\infty,
\end{equation}

\vspace{4mm}
\noindent
where constant $K_4$ does not depend on $p, t_2.$

Combining (\ref{may104}), (\ref{may107}),  (\ref{may108}), we have 

\vspace{-1mm}
\begin{equation}
\label{may109}
\left\vert
\sum\limits_{j=0}^{p}\int\limits_{t}^{t_2}(s-t)^{l}\phi_{j}(s)
\int\limits_{t}^{s}(\tau-t)^{m}\phi_{j}(\tau)d\tau ds
\right\vert\le K_5<\infty,
\end{equation}

\vspace{4mm}
\noindent
where $l, m=0,1,2,\ldots,$ constant $K_5$ does not depend on $p, t_2.$

Note that
$$
\sum\limits_{j=0}^{p}\int\limits_{t_1}^{t_2}(s-t)^{l}\phi_{j}(s)
\int\limits_{t_1}^{s}(\tau-t)^{m}\phi_{j}(\tau)d\tau ds=
$$

\vspace{2mm}
$$
=\sum\limits_{j=0}^{p}\int\limits_{t}^{t_2}(s-t)^{l}\phi_{j}(s)
\int\limits_{t}^{s}(\tau-t)^{m}\phi_{j}(\tau)d\tau ds-
$$

\vspace{2mm}
$$
-\sum\limits_{j=0}^{p}\int\limits_{t}^{t_1}(s-t)^{l}\phi_{j}(s)
\int\limits_{t}^{s}(\tau-t)^{m}\phi_{j}(\tau)d\tau ds-
$$

\vspace{2mm}
\begin{equation}
\label{may110}
-\sum\limits_{j=0}^{p}\int\limits_{t_1}^{t_2}(s-t)^{l}\phi_{j}(s)ds
\int\limits_{t}^{t_1}(\tau-t)^{m}\phi_{j}(\tau)d\tau,
\end{equation}

\vspace{4mm}
\noindent
where $l, m=0,1,2,\ldots$ and $t\le t_1<t_2\le T.$

By analogy with (\ref{may106}) we get

\vspace{-1mm}
\begin{equation}
\label{may111}
\left\vert
\sum\limits_{j=0}^{p}\int\limits_{t_1}^{t_2}(s-t)^{l}\phi_{j}(s)ds
\int\limits_{t}^{t_1}(\tau-t)^{m}\phi_{j}(\tau)d\tau\right\vert
\le K_6<\infty,
\end{equation}

\vspace{4mm}
\noindent
where $l, m=0,1,2,\ldots,$ constant $K_6$ does not depend on $p, t_2.$
Combining (\ref{may110}), (\ref{may109}), and (\ref{may111}), we obtain (\ref{may103x}).

Let us prove (\ref{may101}). Using Parseval's equality, we have

\vspace{-1mm}
$$
\lim\limits_{p\to\infty}\sum\limits_{j_3=0}^{p}
\left(
\frac{1}{2} 
C_{j_3 j_1 j_1}\biggl|_{(j_{1} j_{1})\curvearrowright (\cdot)}
\biggr. - \sum\limits_{j_1=0}^{p} C_{j_3 j_1 j_1}
\right)^2=
$$

\vspace{2mm}
$$
=
\lim\limits_{p\to\infty}\sum\limits_{j_3=0}^{p}
\left(\int\limits_t^T (\tau-t)^{l_3}
\phi_{j_3}(\tau)\left(\frac{1}{2}\int\limits_t^{\tau}(s-t)^{l_1+l_2} ds
-\right.\right.
$$

\vspace{2mm}
$$
\left.\left.-
\sum\limits_{j_1=0}^p
\int\limits_t^{\tau}(s-t)^{l_2}\phi_{j_1}(s)\int\limits_t^s (\theta-t)^{l_1}
\phi_{j_1}(\theta)d\theta ds\right)
d\tau\right)^2\le
$$

\vspace{2mm}
$$
\le
\lim\limits_{p\to\infty}\sum\limits_{j_3=0}^{\infty}
\left(\int\limits_t^T (\tau-t)^{l_3}
\phi_{j_3}(\tau)\left(\frac{1}{2}\int\limits_t^{\tau}(s-t)^{l_1+l_2} ds
-\right.\right.
$$

\vspace{2mm}
$$
\left.\left.-
\sum\limits_{j_1=0}^p
\int\limits_t^{\tau}(s-t)^{l_2}\phi_{j_1}(s)\int\limits_t^s (\theta-t)^{l_1}
\phi_{j_1}(\theta)d\theta ds\right)
d\tau\right)^2=
$$

\begin{equation}
\label{may112}
=\lim\limits_{p\to\infty}
\int\limits_t^T (\tau-t)^{2 l_3}
\left(\frac{1}{2}\int\limits_t^{\tau}(s-t)^{l_1+l_2} ds
-
\sum\limits_{j_1=0}^p
\int\limits_t^{\tau}(s-t)^{l_2}\phi_{j_1}(s)\int\limits_t^s (\theta-t)^{l_1}
\phi_{j_1}(\theta)d\theta ds\right)^2 d\tau.
\end{equation}

\vspace{4mm}

Using (\ref{after1400}), (\ref{may103x}) and
applying Lebesgue's 
Dominated Convergence Theorem in (\ref{may112}), we obtain
the equality (\ref{may101}).

Let us prove (\ref{may102}). Using Fubini's Theorem and Parseval's equality, we obtain

\vspace{-1mm}
$$
\lim\limits_{p\to\infty}
\sum\limits_{j_1=0}^{p}
\left(
\frac{1}{2} 
C_{j_2 j_2 j_1}\biggl|_{(j_{2} j_{2})\curvearrowright (\cdot)}
\biggr. - \sum\limits_{j_2=0}^{p}  C_{j_2 j_2 j_1}
\right)^2=
$$

\vspace{3mm}
$$
=
\lim\limits_{p\to\infty}\sum\limits_{j_1=0}^{p}
\left(\frac{1}{2}\int\limits_t^T (s-t)^{l_2+l_3}
\int\limits_t^{s}
(\theta-t)^{l_1}\phi_{j_1}(\theta)d\theta ds
-\right.
$$

\vspace{2mm}
$$
\left.-
\sum\limits_{j_2=0}^p
\int\limits_t^T (s-t)^{l_3}\phi_{j_2}(s)\int\limits_t^s (\tau-t)^{l_2}
\phi_{j_2}(\tau)\int\limits_t^{\tau} (\theta-t)^{l_1}
\phi_{j_1}(\theta) d\theta d\tau ds\right)^2=
$$

\vspace{2mm}
$$
=
\lim\limits_{p\to\infty}\sum\limits_{j_1=0}^{p}
\left(
\int\limits_t^T (\theta-t)^{l_1}\phi_{j_1}(\theta)
\left(\frac{1}{2}
\int\limits_{\theta}^{T}
(s-t)^{l_2+l_3}ds
-\right.\right.
$$

\vspace{2mm}
$$
\left.\left.-
\sum\limits_{j_2=0}^p
\int\limits_{\theta}^T
(\tau-t)^{l_2}
\phi_{j_2}(\tau)
\int\limits_{\tau}^T
(s-t)^{l_3}\phi_{j_2}(s)
ds d\tau \right)d\theta \right)^2\le
$$

\vspace{2mm}
$$
\le
\lim\limits_{p\to\infty}\sum\limits_{j_1=0}^{\infty}
\left(
\int\limits_t^T (\theta-t)^{l_1}\phi_{j_1}(\theta)
\left(\frac{1}{2}
\int\limits_{\theta}^{T}
(s-t)^{l_2+l_3}ds
-\right.\right.
$$

\vspace{2mm}
$$
\left.\left.-
\sum\limits_{j_2=0}^p
\int\limits_{\theta}^T
(\tau-t)^{l_2}
\phi_{j_2}(\tau)
\int\limits_{\tau}^T
(s-t)^{l_3}\phi_{j_2}(s)
ds d\tau \right)d\theta \right)^2=
$$

\vspace{2mm}
$$
=
\lim\limits_{p\to\infty}
\int\limits_t^T(\theta-t)^{2 l_1}
\left(\frac{1}{2}\int\limits_{\theta}^{T}
(s-t)^{l_2+l_3}ds
-\right.
$$

\vspace{2mm}
$$
\left.
\sum\limits_{j_2=0}^p
\int\limits_{\theta}^T
(\tau-t)^{l_2}
\phi_{j_2}(\tau)
\int\limits_{\tau}^T
(s-t)^{l_3}\phi_{j_2}(s)
ds d\tau\right)^2d\theta =
$$

\vspace{2mm}
\begin{equation}
\label{may113}
=
\lim\limits_{p\to\infty}
\int\limits_t^T(\theta-t)^{2 l_1}
\left(\frac{1}{2}\int\limits_{\theta}^{T}
(s-t)^{l_2+l_3}ds
-
\sum\limits_{j_2=0}^p
\int\limits_{\theta}^T
(s-t)^{l_3}\phi_{j_2}(s)
\int\limits_{\theta}^s
(\tau-t)^{l_2}
\phi_{j_2}(\tau)
d\tau ds\right)^2d\theta.
\end{equation}

\vspace{4mm}

Applying (\ref{after1400}), (\ref{may103x}) and
using Lebesgue's 
Dominated Convergence Theorem in (\ref{may113}), we get
the equality (\ref{may102}).

Let us prove (\ref{may103}). Applying Fubini's Theorem and Parseval's equality, we have

\vspace{1mm}
$$
\lim\limits_{p\to\infty}
\sum\limits_{j_2=0}^{p}
\left(~\sum\limits_{j_1=0}^{p} 
C_{j_1 j_2 j_1}\right)^2=
$$

\vspace{2mm}
$$
=\lim\limits_{p\to\infty}
\sum\limits_{j_2=0}^{p}
\left(~\sum\limits_{j_1=0}^{p} 
\int\limits_t^T (\theta-t)^{l_3} \phi_{j_1}(\theta)\int\limits_t^{\theta}
(\tau-t)^{l_2} \phi_{j_2}(\tau)\int\limits_t^{\tau} (s-t)^{l_1}\phi_{j_1}(s)ds d\tau d\theta\right)^2
=
$$

\vspace{2mm}
$$
=\lim\limits_{p\to\infty}
\sum\limits_{j_2=0}^{p}
\left(~\sum\limits_{j_1=0}^{p} 
\int\limits_t^T 
(\tau-t)^{l_2}\phi_{j_2}(\tau)\int\limits_t^{\tau} (s-t)^{l_1}\phi_{j_1}(s)ds
\int\limits_{\tau}^T
(\theta-t)^{l_3}\phi_{j_1}(\theta)d\theta d\tau \right)^2
\le
$$

\vspace{2mm}
$$
\le\lim\limits_{p\to\infty}
\sum\limits_{j_2=0}^{\infty}
\left(~\int\limits_t^T 
(\tau-t)^{l_2}\phi_{j_2}(\tau)\sum\limits_{j_1=0}^{p} 
\int\limits_t^{\tau} (s-t)^{l_1}\phi_{j_1}(s)ds
\int\limits_{\tau}^T
(\theta-t)^{l_3}\phi_{j_1}(\theta)d\theta d\tau \right)^2
\le
$$

\vspace{2mm}
\begin{equation}
\label{may114}
=\lim\limits_{p\to\infty}
\int\limits_t^T 
(\tau-t)^{2 l_2}\left(\sum\limits_{j_1=0}^{p}\int\limits_t^{\tau} (s-t)^{l_1}\phi_{j_1}(s)ds
\int\limits_{\tau}^T
(\theta-t)^{l_3}\phi_{j_1}(\theta)d\theta\right)^2 d\tau.
\end{equation}

\vspace{4mm}

Applying (\ref{dsds14fffff}), we obtain

\vspace{-1mm}
\begin{equation}
\label{may115}
\left\vert
\sum\limits_{j_1=0}^{p}\int\limits_t^{\tau} (s-t)^{l_1}\phi_{j_1}(s)ds
\int\limits_{\tau}^T
(\theta-t)^{l_3}\phi_{j_1}(\theta)d\theta\right\vert\le C<\infty,
\end{equation}

\vspace{5mm}
\noindent
where constant $C$ does not depend on $p, \tau.$

Using the generalized Parseval equality, we get

\vspace{-1mm}
$$
\sum\limits_{j_1=0}^{\infty}\int\limits_t^{\tau} (s-t)^{l_1}\phi_{j_1}(s)ds
\int\limits_{\tau}^T
(\theta-t)^{l_3}\phi_{j_1}(\theta)d\theta= 
$$

\begin{equation}
\label{may116}
=
\int\limits_t^T (s-t)^{l_1+l_3}{\bf 1}_{\{s<\tau\}}{\bf 1}_{\{s>\tau\}}ds=0.
\end{equation}

\vspace{3mm}

Taking into account (\ref{may115}), (\ref{may116}) and
applying Lebesgue's 
Dominated Convergence Theorem in (\ref{may114}), we obtain
the equality (\ref{may103}). Theorem~49 is proved.

\vspace{5mm}

\section{Expansion of Iterated Stratonovich Stochastic Integrals
of Multiplicity 3. The Case of an Ar\-bit\-ra\-ry Complete Orthonormal System of 
Functions in the Space $L_2([t,T])$ and 
$\psi_1(\tau), \psi_{2}(\tau), \psi_3(\tau)\in L_2([t, T])$}

\vspace{5mm}

In this section, we will prove the following two theorems. 

\vspace{2mm}

{\bf Theorem~50}\ \cite{2018a}.\  {\it Suppose that
$\{\phi_j(x)\}_{j=0}^{\infty}$ is an arbitra\-ry complete ortho\-nor\-mal system of 
func\-ti\-ons in the space $L_2([t,T])$ and $\psi_1(\tau),$ $\psi_2(\tau), \psi_3(\tau)
\in L_2([t,T])$ are such that 

\vspace{-1mm}
\begin{equation}
\label{novemberxxx1}
\Biggl|\sum\limits_{j_1=0}^{p}
\int\limits_t^{s}\psi_2(\tau)\phi_{j_1}(\tau)
\int\limits_t^{\tau}\psi_1(\theta)\phi_{j_1}(\theta)
d\theta d\tau\Biggr|^2\le K<\infty,
\end{equation}

\vspace{1mm}
\begin{equation}
\label{novemberxxx2}
\Biggl|\sum\limits_{j_3=0}^{p}
\int\limits_{s}^T \psi_2(\tau)\phi_{j_3}(\tau)
\int\limits_{\tau}^T \psi_3(\theta)\phi_{j_3}(\theta)d\theta d\tau\Biggr|^2
\le K<\infty
\end{equation}

\vspace{3mm}
\noindent
$\forall p\in \mathbb{N},$ where constant $K$ does not depend on $p$ and $s$ $(t\le s\le T).$
Then$,$ for the sum $\bar J^{*}[\psi^{(3)}]_{T,t}^{(i_1 i_2 i_3)}$
$(i_1,i_2,i_3=0,1,\ldots,m)$
of iterated Ito stochastic integrals 
defined by {\rm (\ref{dsds9}) $(k=3)$}
the following 
expansion 

$$
\bar J^{*}[\psi^{(3)}]_{T,t}^{(i_1 i_2 i_3)}=
\hbox{\vtop{\offinterlineskip\halign{
\hfil#\hfil\cr
{\rm l.i.m.}\cr
$\stackrel{}{{}_{p\to \infty}}$\cr
}} }\sum_{j_1,j_2,j_3=0}^{p}
C_{j_3 j_2 j_1}\zeta_{j_1}^{(i_1)}\zeta_{j_2}^{(i_2)}\zeta_{j_3}^{(i_3)}
$$

\vspace{4mm}
\noindent
that converges in the mean-square sense is valid, where 

$$
C_{j_3 j_2 j_1}=\int\limits_t^T \psi_3(t_3)
\phi_{j_3}(t_3)\int\limits_t^{t_3}\psi_2(t_2)
\phi_{j_2}(t_2)
\int\limits_t^{t_2}\psi_1(t_1)
\phi_{j_1}(t_1)dt_1dt_2dt_3
$$

\vspace{2mm}
\noindent
and
$$
\zeta_{j}^{(i)}=
\int\limits_t^T \phi_{j}(\tau) d{\bf w}_{\tau}^{(i)}
$$ 

\vspace{3mm}
\noindent
are independent standard Gaussian random variables for various 
$i$ or $j$ {\rm (}in the case when $i\ne 0${\rm ),}
${\bf w}_{\tau}^{(i)}={\bf f}_{\tau}^{(i)}$ for
$i=1,\ldots,m$ and 
${\bf w}_{\tau}^{(0)}=\tau.$}

\vspace{2mm}

{\bf Theorem~51}\ \cite{2018a}.\  {\it Suppose that
$\{\phi_j(x)\}_{j=0}^{\infty}$ is an arbitrary complete ortho\-nor\-mal system of 
func\-ti\-ons in the space $L_2([t,T])$ and $\psi_1(\tau),$ $\psi_2(\tau), \psi_3(\tau)$
are continuous func\-ti\-ons on $[t, T].$
Furthermore$,$ let the conditions {\rm (\ref{novemberxxx1}), (\ref{novemberxxx2})}
are satisfied.
Then$,$ for the iterated Stra\-to\-no\-vich stochastic integral
of third multiplicity

$$
{\int\limits_t^{*}}^T \psi_3(t_3)
{\int\limits_t^{*}}^{t_3}\psi_2(t_2)
{\int\limits_t^{*}}^{t_2}\psi_1(t_1)
d{\bf w}_{t_1}^{(i_1)}
d{\bf w}_{t_2}^{(i_2)}d{\bf w}_{t_3}^{(i_3)}\ \ \ (i_1,i_2,i_3=0,1,\ldots,m)
$$

\vspace{3mm}
\noindent
the following 
expansion 

\vspace{-1mm}
$$
{\int\limits_t^{*}}^T \hspace{-1.5mm}\psi_3(t_3)
{\int\limits_t^{*}}^{t_3}\hspace{-1.5mm}\psi_2(t_2)
{\int\limits_t^{*}}^{t_2}\hspace{-1.5mm}\psi_1(t_1)
d{\bf w}_{t_1}^{(i_1)}
d{\bf w}_{t_2}^{(i_2)}d{\bf w}_{t_3}^{(i_3)}=
\hbox{\vtop{\offinterlineskip\halign{
\hfil#\hfil\cr
{\rm l.i.m.}\cr
$\stackrel{}{{}_{p\to \infty}}$\cr
}} }\hspace{-1.5mm}\sum_{j_1,j_2,j_3=0}^{p}\hspace{-1.5mm}
C_{j_3 j_2 j_1}\zeta_{j_1}^{(i_1)}\zeta_{j_2}^{(i_2)}\zeta_{j_3}^{(i_3)}
$$

\vspace{3mm}
\noindent
that converges in the mean-square sense is valid, where 
notations are the same as in Theorem~{\rm 50.}}

\vspace{2mm}

Note that Theorem~51 is a simple consequence of Theorem~50 and Theorem~5 $(k=3).$
Let us prove Theorem~50.

{\bf Proof.}\ First, let us note
some facts that follow from Monotone Convergence Theorem
(\cite{Pugach}, Theorem~3.5.1) and Lebesgue's Dominated Convergence Theorem.
Suppose that $\left\{g_j(x)\right\}_{j=0}^{\infty}$ is an arbitrary
sequence of real-valued measurable functions such that

\vspace{-1mm}
\begin{equation}
\label{july1000}
\sum\limits_{j=0}^{\infty}\left|g_j(x)\right|\le K<\infty
\end{equation} 

\vspace{2mm}
\noindent
almost everywhere on $X$ (with respect to Lebesgue's measure),
where
constant $K$ does not depend on $x.$

It is easy to see that under the above conditions
the following equality 

\vspace{-1mm}
\begin{equation}
\label{july999}
\lim\limits_{p\to\infty}\int\limits_X h^2(x)\left(\sum\limits_{j=0}^{p}g_j(x)\right)^2 dx=
\int\limits_X h^2(x)\left(\sum\limits_{j=0}^{\infty} g_j(x)\right)^2 dx
\end{equation}

\vspace{3mm}
\noindent
is true, where $h(x)\in L_2(X)$ (further, we put $h(x)\equiv 1$ for simplicity).
Indeed, we have $g_j(x)=g_j^{+}(x)-g_j^{-}(x),$ 
$\left\vert g_j(x)\right\vert =g_j^{+}(x)+g_j^{-}(x),$ where 
$g_j^{+}(x)=\max\{g_j(x), 0\}\ge 0,$
$g_j^{-}(x)=-\min\{g_j(x), 0\}\ge 0.$ Moreover,

$$
\sum\limits_{j=0}^{\infty}g_j(x) =
\sum\limits_{j=0}^{\infty} g_j^{+}(x)-\sum\limits_{j=0}^{\infty} g_j^{-}(x),
$$

\vspace{-1mm}
\begin{equation}
\label{july1002}
\sum\limits_{j=0}^{\infty}\left\vert g_j(x)\right\vert =
\sum\limits_{j=0}^{\infty} g_j^{+}(x)+\sum\limits_{j=0}^{\infty} g_j^{-}(x).
\end{equation}

\vspace{3mm}

Uning (\ref{july1000}), we obtain that
the series (with non-negative terms) on the right-hand side of 
(\ref{july1002}) satisfy the condition (\ref{july1000}).
Further, using Monotone Convergence Theorem, we obtain

\vspace{-1mm}
$$
\lim\limits_{p\to\infty}\int\limits_X \left(\sum\limits_{j=0}^{p}g_j(x)\right)^2 dx=
\lim\limits_{p\to\infty}\int\limits_X \left(
\sum\limits_{j=0}^{p} g_j^{+}(x)-\sum\limits_{j=0}^{p} g_j^{-}(x)
\right)^2 dx=
$$

\vspace{1mm}
$$
=
\lim\limits_{p\to\infty} \int\limits_X  \left(
\sum\limits_{j=0}^{p} g_j^{+}(x)
\right)^2 dx-
\lim\limits_{p\to\infty} 2\int\limits_X 
\sum\limits_{j=0}^{p} g_j^{+}(x)\sum\limits_{j=0}^{p} g_j^{-}(x)
dx +
\lim\limits_{p\to\infty} \int\limits_X \left(
\sum\limits_{j=0}^{p} g_j^{-}(x)
\right)^2 dx=
$$

\vspace{1mm}
$$
=
\int\limits_X \lim\limits_{p\to\infty} \left(
\sum\limits_{j=0}^{p} g_j^{+}(x)
\right)^2 dx-
2\int\limits_X \lim\limits_{p\to\infty}
\sum\limits_{j=0}^{p} g_j^{+}(x)\sum\limits_{j=0}^{p} g_j^{-}(x)
dx +
\int\limits_X \lim\limits_{p\to\infty}\left(
\sum\limits_{j=0}^{p} g_j^{-}(x)
\right)^2 dx=
$$

\vspace{1mm}
\begin{equation}
\label{slovo4}
=
\int\limits_X \left(
\sum\limits_{j=0}^{\infty} g_j^{+}(x)
\right)^2 dx-
2\int\limits_X 
\sum\limits_{j=0}^{\infty} g_j^{+}(x)\sum\limits_{j=0}^{\infty} g_j^{-}(x)
dx +
\int\limits_X \left(
\sum\limits_{j=0}^{\infty} g_j^{-}(x)
\right)^2 dx=
\end{equation}

\vspace{1mm}
$$
=
\int\limits_X \left(
\sum\limits_{j=0}^{\infty} g_j^{+}(x)-
\sum\limits_{j=0}^{\infty} g_j^{-}(x)
\right)^2 dx=
\int\limits_X \left(\sum\limits_{j=0}^{\infty} g_j(x)\right)^2 dx.
$$

\vspace{4mm}

The equality (\ref{july999}) can be obtained under another
conditions. If we replace
the condition (\ref{july1000}) with
\begin{equation}
\label{july1000aaa1}
\left|\sum\limits_{j=0}^{p}g_j(x)\right| \le K<\infty\ \ \forall p\in \mathbb{N}\ \ \ 
\hbox{and}\ \ \ \lim\limits_{p\to\infty}\sum\limits_{j=0}^{p}g_j(x)\ \ \hbox{exists}
\end{equation} 

\vspace{3mm}
\noindent
almost everywhere on $X$ (with respect to Lebesgue's measure),
then by Lebesgue's Dominated Convergence Theorem
we obtain (\ref{july999}). Here
constant $K$ does not depend on $x$ and $p.$

According to Theorem~47, we come to the conclusion that 
Theorem~50 will be proved if we prove the following
equalities

\vspace{-1mm}
\begin{equation}
\label{july1003}
\lim\limits_{p\to\infty}
\sum\limits_{j_3=0}^{p}
\left(\frac{1}{2} 
C_{j_3 j_1 j_1}\biggl|_{(j_{1} j_{1})\curvearrowright (\cdot)}
\biggr.-
\sum\limits_{j_1=0}^{p}  C_{j_3 j_1 j_1}
\right)^2=0,
\end{equation}

\vspace{1mm}
\begin{equation}
\label{july1004}
\lim\limits_{p\to\infty}
\sum\limits_{j_1=0}^{p}
\left(\frac{1}{2} 
C_{j_3 j_3 j_1}\biggl|_{(j_{3} j_{3})\curvearrowright (\cdot)}
\biggr.- \sum\limits_{j_3=0}^{p} C_{j_3 j_3 j_1}
\right)^2=0,
\end{equation}

\vspace{1mm}
\begin{equation}
\label{july1005}
\lim\limits_{p\to\infty}
\sum\limits_{j_2=0}^{p}
\left(\sum\limits_{j_1=0}^{p} 
C_{j_1 j_2 j_1}\right)^2=0.
\end{equation}

\vspace{4mm}

Let us prove (\ref{july1003}). Using Parseval's equality, we have

$$
\lim\limits_{p\to\infty}
\sum\limits_{j_3=0}^{p}
\left(\frac{1}{2} 
C_{j_3 j_1 j_1}\biggl|_{(j_{1} j_{1})\curvearrowright (\cdot)}
\biggr.-
\sum\limits_{j_1=0}^{p}  C_{j_3 j_1 j_1}
\right)^2=
$$

\vspace{2mm}
$$
=\lim\limits_{p\to\infty}
\sum\limits_{j_3=0}^{p}
\left(\int\limits_t^T \psi_3(s)
\phi_{j_3}(s)\left(\frac{1}{2}\int\limits_t^{s} \psi_2(\tau)\psi_1(\tau)d\tau
-
\sum\limits_{j_1=0}^p
\int\limits_t^{s}\psi_2(\tau)\phi_{j_1}(\tau)
\int\limits_t^{\tau}\psi_1(\theta)\phi_{j_1}(\theta)
d\theta d\tau\right)ds\right)^2\le
$$

\vspace{2mm}
$$
\le\lim\limits_{p\to\infty}
\sum\limits_{j_3=0}^{\infty}
\left(\int\limits_t^T \psi_3(s)
\phi_{j_3}(s)\left(\frac{1}{2}\int\limits_t^{s} \psi_2(\tau)\psi_1(\tau)d\tau
-
\sum\limits_{j_1=0}^p
\int\limits_t^{s}\psi_2(\tau)\phi_{j_1}(\tau)
\int\limits_t^{\tau}\psi_1(\theta)\phi_{j_1}(\theta)
d\theta d\tau\right)ds\right)^2=
$$

\vspace{2mm}
\begin{equation}
\label{july1006}
=\lim\limits_{p\to\infty}
\int\limits_t^T \psi_3^2(s)
\left(\frac{1}{2}\int\limits_t^{s} \psi_2(\tau)\psi_1(\tau)d\tau
-
\sum\limits_{j_1=0}^p
\int\limits_t^{s}\psi_2(\tau)\phi_{j_1}(\tau)
\int\limits_t^{\tau}\psi_1(\theta)\phi_{j_1}(\theta)
d\theta d\tau \right)^2 ds=
\end{equation}

\vspace{2mm}
\begin{equation}
\label{july1007}
=
\int\limits_t^T 
\psi_3^2(s) \lim\limits_{p\to\infty}\left(\frac{1}{2}\int\limits_t^{s} \psi_2(\tau)\psi_1(\tau)d\tau
-
\sum\limits_{j_1=0}^p
\int\limits_t^{s}\psi_2(\tau)\phi_{j_1}(\tau)
\int\limits_t^{\tau}\psi_1(\theta)\phi_{j_1}(\theta)
d\theta d\tau \right)^2 ds=0,
\end{equation}

\vspace{4mm}
\noindent
where (\ref{july1007}) follows from (\ref{after1400}) and
the transition from (\ref{july1006}) to (\ref{july1007}) 
is based on (\ref{july999}), (\ref{july1000aaa1}) and
Lebesgue's Dominated Convergence Theorem (see (\ref{novemberxxx1})).
The equality (\ref{july1003}) is proved.

Let us prove (\ref{july1004}). Using Fubini's Theorem and Parseval's equality, we obtain

$$
\lim\limits_{p\to\infty}
\sum\limits_{j_1=0}^{p}
\left(\frac{1}{2} 
C_{j_3 j_3 j_1}\biggl|_{(j_{3} j_{3})\curvearrowright (\cdot)}
\biggr.- \sum\limits_{j_3=0}^{p} C_{j_3 j_3 j_1}
\right)^2=
$$

\vspace{2mm}
$$
=\lim\limits_{p\to\infty}
\sum\limits_{j_1=0}^{p}
\left(\frac{1}{2}\int\limits_t^T \psi_3(\tau)\psi_2(\tau)
\int\limits_t^{\tau} \psi_1(s)\phi_{j_1}(s)dsd\tau
-\right.
$$

\vspace{2mm}
$$
-\left.
\sum\limits_{j_3=0}^p
\int\limits_t^T \psi_3(\theta)\phi_{j_3}(\theta)\int\limits_t^{\theta} \psi_2(\tau)\phi_{j_3}(\tau)
\int\limits_t^{\tau} \psi_1(s)\phi_{j_1}(s)ds d\tau d\theta\right)^2=
$$

\vspace{2mm}
$$
=\lim\limits_{p\to\infty}
\sum\limits_{j_1=0}^{p}
\left(\frac{1}{2}\int\limits_t^T 
\psi_1(s)\phi_{j_1}(s)\int\limits_s^T \psi_3(\tau)\psi_2(\tau) d\tau ds
-\right.
$$

\vspace{2mm}
$$
\left.
-\sum\limits_{j_3=0}^p
\int\limits_t^T \psi_1(s)\phi_{j_1}(s)\int\limits_{s}^T \psi_2(\tau)\phi_{j_3}(\tau)
\int\limits_{\tau}^T \psi_3(\theta)\phi_{j_3}(\theta)d\theta d\tau ds\right)^2=
$$

\vspace{2mm}
$$
=\lim\limits_{p\to\infty}
\sum\limits_{j_1=0}^{p}
\left(\int\limits_t^T 
\psi_1(s)\phi_{j_1}(s)\left(\frac{1}{2}\int\limits_s^T \psi_3(\tau)\psi_2(\tau) d\tau 
-\sum\limits_{j_3=0}^p
\int\limits_{s}^T \psi_2(\tau)\phi_{j_3}(\tau)
\int\limits_{\tau}^T \psi_3(\theta)\phi_{j_3}(\theta)d\theta d\tau
\right)ds\right)^2\le
$$

\vspace{2mm}
$$
\le\lim\limits_{p\to\infty}
\sum\limits_{j_1=0}^{\infty}
\left(\int\limits_t^T 
\psi_1(s)\phi_{j_1}(s)\left(\frac{1}{2}\int\limits_s^T \psi_3(\tau)\psi_2(\tau) d\tau 
-\sum\limits_{j_3=0}^p
\int\limits_{s}^T \psi_2(\tau)\phi_{j_3}(\tau)
\int\limits_{\tau}^T \psi_3(\theta)\phi_{j_3}(\theta)d\theta d\tau
\right)ds\right)^2=
$$

\vspace{2mm}
\begin{equation}
\label{july1008}
=\lim\limits_{p\to\infty}
\int\limits_t^T 
\psi_1^2(s)\left(\frac{1}{2}\int\limits_s^T \psi_3(\tau)\psi_2(\tau) d\tau 
-\sum\limits_{j_3=0}^p
\int\limits_{s}^T \psi_2(\tau)\phi_{j_3}(\tau)
\int\limits_{\tau}^T \psi_3(\theta)\phi_{j_3}(\theta)d\theta d\tau
\right)^2 ds=
\end{equation}

\vspace{2mm}
\begin{equation}
\label{july1009}
=
\int\limits_t^T 
\psi_1^2(s) \lim\limits_{p\to\infty}\left(\frac{1}{2}\int\limits_s^T \psi_3(\tau)\psi_2(\tau) d\tau 
-\sum\limits_{j_3=0}^p
\int\limits_{s}^T \psi_2(\tau)\phi_{j_3}(\tau)
\int\limits_{\tau}^T \psi_3(\theta)\phi_{j_3}(\theta)d\theta d\tau
\right)^2 ds=0,
\end{equation}

\vspace{4mm}
\noindent
where (\ref{july1009}) follows from (\ref{after1400}) and
the transition from (\ref{july1008}) to (\ref{july1009}) 
is based on (\ref{july999}), (\ref{july1000aaa1}) and
Lebesgue's Dominated Convergence Theorem (see (\ref{novemberxxx2})).
The equality (\ref{july1004}) is proved.

Let us prove (\ref{july1005}).
Applying Fubini's Theorem and Parseval's equality, we have

$$
\lim\limits_{p\to\infty}
\sum\limits_{j_2=0}^{p}
\left(\sum\limits_{j_1=0}^{p} 
C_{j_1 j_2 j_1}\right)^2=
$$

\vspace{2mm}
$$
=\lim\limits_{p\to\infty}
\sum\limits_{j_2=0}^{p}
\left(\sum\limits_{j_1=0}^{p} 
\int\limits_t^T \psi_3(\theta)\phi_{j_1}(\theta)\int\limits_t^{\theta}
\psi_2(\tau)\phi_{j_2}(\tau)\int\limits_t^{\tau} \psi_1(s)\phi_{j_1}(s)ds d\tau d\theta\right)^2=
$$

\vspace{2mm}
$$
=\lim\limits_{p\to\infty}
\sum\limits_{j_2=0}^{p}
\left(\sum\limits_{j_1=0}^{p} 
\int\limits_t^T 
\psi_2(\tau)\phi_{j_2}(\tau)\int\limits_t^{\tau} \psi_1(s)\phi_{j_1}(s)ds
\int\limits_{\tau}^T
\psi_3(\theta)\phi_{j_1}(\theta)d\theta d\tau \right)^2\le
$$

\vspace{2mm}
$$
\le\lim\limits_{p\to\infty}
\sum\limits_{j_2=0}^{\infty}
\left(
\int\limits_t^T 
\psi_2(\tau)\phi_{j_2}(\tau) \sum\limits_{j_1=0}^{p} \int\limits_t^{\tau} \psi_1(s)\phi_{j_1}(s)ds
\int\limits_{\tau}^T
\psi_3(\theta)\phi_{j_1}(\theta)d\theta d\tau \right)^2=
$$

\vspace{2mm}
\begin{equation}
\label{july1010}
=\lim\limits_{p\to\infty}
\int\limits_t^T 
\psi_2^2(\tau)\left(\sum\limits_{j_1=0}^{p} \int\limits_t^{\tau} \psi_1(s)\phi_{j_1}(s)ds
\int\limits_{\tau}^T
\psi_3(\theta)\phi_{j_1}(\theta)d\theta \right)^2 d\tau =
\end{equation}

\vspace{2mm}
\begin{equation}
\label{july1011}
=
\int\limits_t^T 
\psi_2^2(\tau)\lim\limits_{p\to\infty}\left(\sum\limits_{j_1=0}^{p} 
\int\limits_t^{\tau} \psi_1(s)\phi_{j_1}(s)ds
\int\limits_{\tau}^T
\psi_3(\theta)\phi_{j_1}(\theta)d\theta \right)^2 d\tau = 0,
\end{equation}

\vspace{5mm}
\noindent
where (\ref{july1011}) follows from the equality

\vspace{-1mm}
\begin{equation}
\label{july1012}
\sum\limits_{j_1=0}^{\infty} 
\int\limits_t^{\tau} \psi_1(s)\phi_{j_1}(s)ds
\int\limits_{\tau}^T
\psi_3(\theta)\phi_{j_1}(\theta)d\theta=
\int\limits_t^T \psi_1(s){\bf 1}_{\{s<\tau\}}\psi_3(s){\bf 1}_{\{s>\tau\}}ds=0
\end{equation}

\vspace{4mm}
\noindent
(the relation (\ref{july1012}) follows from the generalized Parseval equality)
and
the transition from (\ref{july1010}) to (\ref{july1011}) 
is based on (\ref{july999}), (\ref{july1000aaa1}) and
Lebesgue's Dominated Convergence Theorem (see (\ref{dsds14fffff})).
The equality (\ref{july1005}) is proved.
Theorem~50 is proved.

\vspace{5mm}

\section{Expansion of Iterated Stratonovich Stochastic Integrals
of Multiplicities 4 and 5. The Case of an Ar\-bit\-ra\-ry Complete Orthonormal System of 
Functions in the Space $L_2([t,T])$ and 
$\psi_1(\tau),\ldots,\psi_5(\tau)\in L_2([t, T])$}

\vspace{5mm}

Let us develop the approach discussed in the previous section.
It is easy to see (according to Theorem~47) that analogues of Theorems~50 and 51 
for the cases $k=4$ and $k=5$ ($\psi_1(\tau), \ldots, \psi_5(\tau)
\in L_2([t,T])$) will be true if the relations
(\ref{2023novem200})--(\ref{2023novem205}), 
(\ref{april11})--(\ref{april35}) 
as well as the equalities

\vspace{-2mm}
\begin{equation}
\label{july13}
\lim\limits_{p\to\infty}
\sum\limits_{j_1, j_3=0}^{p}
C_{j_3 j_3 j_1 j_1}=\frac{1}{4}\int\limits_{t}^{T} 
\psi_4(t_3)\psi_3(t_3)\int\limits_{t}^{t_3}
\psi_2(t_1)\psi_1(t_1)dt_1 dt_3,
\biggr.
\end{equation}

\begin{equation}
\label{july14}
\lim\limits_{p\to\infty}
\sum\limits_{j_1, j_3=0}^{p}
C_{j_1 j_3 j_3 j_1}=0,
\biggr.
\end{equation}

\begin{equation}
\label{july15}
\lim\limits_{p\to\infty}
\sum\limits_{j_1, j_2=0}^{p}
C_{j_2 j_1 j_2 j_1}=0,
\biggr.
\end{equation}

\begin{equation}
\label{july16}
\lim\limits_{p\to\infty}
\sum\limits_{j_1, j_3=0}^{p}
C_{j_3 j_3 j_1 j_1}(s,\tau)
=\frac{1}{4}\int\limits_{\tau}^{s} \psi_4(t_3)\psi_3(t_3)\int\limits_{\tau}^{t_3}
\psi_2(t_1)\psi_1(t_1)dt_1 dt_3,
\biggr.
\end{equation}

\begin{equation}
\label{july17}
\lim\limits_{p\to\infty}
\sum\limits_{j_1, j_3=0}^{p}
C_{j_1 j_3 j_3 j_1}(s,\tau)=0,
\biggr.
\end{equation}

\begin{equation}
\label{july18}
\lim\limits_{p\to\infty}
\sum\limits_{j_1, j_2=0}^{p}
C_{j_2 j_1 j_2 j_1}(s,\tau)=0
\biggr.
\end{equation}

\vspace{3mm}
\noindent
are satisfied, provided that 
$\{\phi_j(x)\}_{j=0}^{\infty}$ is an arbitrary complete ortho\-nor\-mal system of 
functions in the space $L_2([t,T]),$ $\psi_1(\tau), \ldots, \psi_5(\tau)
\in L_2([t,T]),$ the series on the left-hand sides of (\ref{july13})--(\ref{july18}) 
converge absolutely, and

\vspace{-3mm}
$$
C_{j_4 \ldots j_1}=\int\limits_t^T
\psi_4(t_4)\phi_{j_4}(t_4)\ldots 
\int\limits_t^{t_2}
\psi_1(t_1)\phi_{j_1}(t_1)dt_1\ldots dt_4,
$$

\vspace{-1mm}
$$
C_{j_5 \ldots j_1}=\int\limits_t^T
\psi_5(t_5)\phi_{j_5}(t_5)\ldots 
\int\limits_t^{t_2}
\psi_1(t_1)\phi_{j_1}(t_1)dt_1\ldots dt_5,
$$

\vspace{-1mm}
$$
C_{j_4 \ldots j_1}(s,\tau)=\int\limits_{\tau}^s
\psi_4(t_4)\phi_{j_4}(t_4)\ldots 
\int\limits_{\tau}^{t_2}
\psi_1(t_1)\phi_{j_1}(t_1)dt_1\ldots dt_4
$$

\vspace{2mm}
\noindent
in (\ref{2023novem200})--(\ref{2023novem205}), 
(\ref{april11})--(\ref{april35}), (\ref{july13})--(\ref{july18}).

It is obvious that the equalities (\ref{july16})--(\ref{july18}) follow from the equalities 
(\ref{july13})--(\ref{july15})
if in (\ref{july13})--(\ref{july15}) we replace $\psi_4(t_4), 
\psi_3(t_3), \psi_2(t_2), \psi_1(t_1)$ with 
${\bf 1}_{\{\tau<t_4<s\}}\psi_4(t_4)$, 
${\bf 1}_{\{\tau<t_3\}}\psi_3(t_3)$, 
${\bf 1}_{\{\tau<t_2\}}\psi_2(t_2)$, 
${\bf 1}_{\{\tau<t_1\}}\psi_1(t_1),$
respectively.

Further, the proofs of Theorems~44 and 48 must be modified 
and carried out by analogy with the proof 
of Theorem~50, i.e. using the equality (\ref{july999}) and
Lebesgue's Dominated Convergence Theorem.
At that, the derivation of formulas similar to (\ref{2023novem209})--(\ref{2023novem214}),
(\ref{april36})--(\ref{april45}), (\ref{april62})--(\ref{april100}),
(\ref{april101}), (\ref{april102}), (\ref{april103}),
(\ref{april104}), (\ref{april105}), (\ref{april108}), (\ref{april111}), (\ref{april114})
is carried out completely similarly to (\ref{2023novem209})--(\ref{2023novem214}),
(\ref{april36})--(\ref{april45}), (\ref{april62})--(\ref{april100}),
(\ref{april101}), (\ref{april102}), (\ref{april103}),
(\ref{april104}), (\ref{april105}), (\ref{april108}), (\ref{april111}), (\ref{april114}), adjusted
for the fact that in (\ref{2023novem209})--(\ref{2023novem214}),
(\ref{april36})--(\ref{april45}), (\ref{april62})--(\ref{april100}),
(\ref{april101}), (\ref{april102}), (\ref{april103}),
(\ref{april104}), (\ref{april105}), (\ref{april108}), (\ref{april111}), (\ref{april114}) 
the functions $\psi_1(\tau),\ldots,\psi_5(\tau)\equiv 1$ are replaced by 
$\psi_1(\tau),\ldots,\psi_5(\tau)\in L_2([t, T])$.
Furthermore, the following conditions 

\vspace{-3mm}
\begin{equation}
\label{novemberxxx3}
\left|\sum\limits_{j=0}^{p}
\int\limits_{\tau}^{s}\psi_{m+1}(t_2)\phi_{j}(t_2)
\int\limits_{\tau}^{t_2}\psi_m(t_1)\phi_{j}(t_1)
dt_1 dt_2\right|^2\le K<\infty\ \ \ (m=1,2,3,4),
\end{equation}
\begin{equation}
\label{novemberxxx4}
\left|\sum\limits_{j_1,j_2=0}^{p}
C_{j_2 j_2 j_1 j_1}^{\psi_{m+3}\psi_{m+2}\psi_{m+1}\psi_m}(s,\tau)\right|^2\le K<\infty\ \ \ (m=1,2),
\end{equation}
\begin{equation}
\label{novemberxxx5}
\left|\sum\limits_{j_1,j_2=0}^{p}
C_{j_2 j_1 j_2 j_1}^{\psi_{m+3}\psi_{m+2}\psi_{m+1}\psi_m}(s,\tau)\right|^2\le K<\infty\ \ \ (m=1,2),
\end{equation}
\begin{equation}
\label{novemberxxx6}
\left|\sum\limits_{j_1,j_2=0}^{p}
C_{j_1 j_2 j_2 j_1}^{\psi_{m+3}\psi_{m+2}\psi_{m+1}\psi_m}(s,\tau)\right|^2\le K<\infty\ \ \ (m=1,2),
\end{equation}

\vspace{2mm}
\noindent
must be satisfied $\forall p\in \mathbb{N},$ 
where constant $K$ does not depend on $p, \tau, s$,

\vspace{-1mm}
$$
C_{j_4 j_3 j_2 j_1}^{\psi_{m+3}\psi_{m+2}\psi_{m+1}\psi_m}(s,\tau)=
\int\limits_{\tau}^{s}\psi_{m+3}(t_4)\phi_{j_4}(t_4)\times
$$

\vspace{-1mm}
$$
\times
\int\limits_{\tau}^{t_4}\psi_{m+2}(t_3)\phi_{j_3}(t_3)
\int\limits_{\tau}^{t_3}\psi_{m+1}(t_2)\phi_{j_2}(t_2)
\int\limits_{\tau}^{t_2}\psi_{m}(t_1)\phi_{j_1}(t_1)dt_1 dt_2 dt_3 dt_4,
$$

\vspace{3mm}
\noindent
where $m=1, 2$ and $t\le \tau<s\le T.$

The conditions (\ref{novemberxxx3})--(\ref{novemberxxx6}) are
required to perform the passage to the limit
using Lebesgue's Dominated Convergence Theorem
(see the proofs of Theorems~44, 48 for details).

The equality (\ref{july13}) is proved in \cite{rybakov5000}
for the case when $\{\phi_j(x)\}_{j=0}^{\infty}$ is an arbitrary complete ortho\-nor\-mal system of 
functions in the space $L_2([t,T])$ and $\psi_1(\tau),$ $\ldots, \psi_4(\tau)
\in L_2([t,T]).$ 
The equalities (\ref{july14}), (\ref{july15}) can also be obtained
\cite{rybakov7000xa}
using the approach from \cite{rybakov5000}. At that, the series on the left-hand
sides of (\ref{july13})--(\ref{july15}) converge absolutly.
We will return to these issues in Sect.~26. 
The part of Sect.~26 will be devoted to the method from \cite{rybakov5000}
based on trace class operators.
In Sect.~26, we will also prove the equalities 
(\ref{july13})--(\ref{july15}) using an approach based on the 
generalized Parseval equality and (\ref{after1400})
(the case when $\{\phi_j(x)\}_{j=0}^{\infty}$ is an arbitrary complete ortho\-nor\-mal system of 
functions in the space $L_2([t,T])$ and $\psi_1(\tau),$ $\ldots, \psi_4(\tau)
\in L_2([t,T])).$

Taking into account everything said above in this section
and the results of Sect.~26 (see below), we obtain the following four theorems.

\vspace{2mm}

{\bf Theorem~52}\ \cite{2018a}.\ {\it Suppose that
$\{\phi_j(x)\}_{j=0}^{\infty}$ is an arbitrary complete ortho\-nor\-mal system of 
func\-ti\-ons in the space $L_2([t,T])$ and $\psi_1(\tau),$ $\ldots, \psi_4(\tau)
\in L_2([t,T]).$ Furthermore, let the condition
{\rm (\ref{novemberxxx3})} $(m=1,2,3)$
is satisfied.
Then$,$ for the sum $\bar J^{*}[\psi^{(4)}]_{T,t}^{(i_1 \ldots i_4)}$
$(i_1,\ldots,i_4=0,1,\ldots,m)$
of iterated Ito stochastic integrals 
defined by {\rm (\ref{dsds9}) $(k=4)$}
the following 
expansion 

\vspace{-1mm}
$$
\bar J^{*}[\psi^{(4)}]_{T,t}^{(i_1 \ldots i_4)}=
\hbox{\vtop{\offinterlineskip\halign{
\hfil#\hfil\cr
{\rm l.i.m.}\cr
$\stackrel{}{{}_{p\to \infty}}$\cr
}} }\sum_{j_1,\ldots,j_4=0}^{p}
C_{j_4 \ldots j_1}\zeta_{j_1}^{(i_1)}\ldots \zeta_{j_4}^{(i_4)}
$$

\vspace{3mm}
\noindent
that converges in the mean-square sense is valid, where 

\vspace{-1mm}
$$
C_{j_4 \ldots j_1}=\int\limits_t^T \psi_4(t_4)\phi_{j_4}(t_4)
\ldots
\int\limits_t^{t_2}\psi_1(t_1)
\phi_{j_1}(t_1)dt_1\ldots dt_4
$$

\vspace{3mm}
and
$$
\zeta_{j}^{(i)}=
\int\limits_t^T \phi_{j}(\tau) d{\bf w}_{\tau}^{(i)}
$$ 

\vspace{2mm}
\noindent
are independent standard Gaussian random variables for various 
$i$ or $j$ {\rm (}in the case when $i\ne 0${\rm ),}
${\bf w}_{\tau}^{(i)}={\bf f}_{\tau}^{(i)}$ for
$i=1,\ldots,m$ and 
${\bf w}_{\tau}^{(0)}=\tau.$}

\vspace{2mm}

{\bf Theorem~53}\ \cite{2018a}.\ {\it Suppose that
$\{\phi_j(x)\}_{j=0}^{\infty}$ is an arbitrary complete ortho\-nor\-mal system of 
func\-ti\-ons in the space $L_2([t,T])$ and $\psi_1(\tau),$ $\ldots, \psi_4(\tau)$
are continuous functions on $[t, T].$
Furthermore, let the condition {\rm (\ref{novemberxxx3})} $(m=1,2,3)$
is satisfied.
Then$,$ for the iterated Stra\-to\-no\-vich stochastic integral
of fourth multiplicity

\vspace{-1mm}
$$
{\int\limits_t^{*}}^T \psi_4(t_4)
\ldots 
{\int\limits_t^{*}}^{t_2}\psi_1(t_1)
d{\bf w}_{t_1}^{(i_1)}
\ldots d{\bf w}_{t_4}^{(i_4)}\ \ \ (i_1,\ldots,i_4=0,1,\ldots,m)
$$

\vspace{3mm}
\noindent
the following 
expansion 

\vspace{-1mm}
$$
{\int\limits_t^{*}}^T \psi_4(t_4)
\ldots 
{\int\limits_t^{*}}^{t_2}\psi_1(t_1)
d{\bf w}_{t_1}^{(i_1)}
\ldots d{\bf w}_{t_4}^{(i_4)}=
\hbox{\vtop{\offinterlineskip\halign{
\hfil#\hfil\cr
{\rm l.i.m.}\cr
$\stackrel{}{{}_{p\to \infty}}$\cr
}} }\sum_{j_1,\ldots,j_4=0}^{p}
C_{j_4 \ldots j_1}\zeta_{j_1}^{(i_1)}\ldots \zeta_{j_4}^{(i_4)}
$$

\vspace{3mm}
\noindent
that converges in the mean-square sense is valid, where 
notations are the same as in Theorem~{\rm 52.}}

\vspace{2mm}

{\bf Theorem~54}\ \cite{2018a}.\ {\it Suppose that
$\{\phi_j(x)\}_{j=0}^{\infty}$ is an arbitrary complete ortho\-nor\-mal system of 
func\-ti\-ons in the space $L_2([t,T])$ and $\psi_1(\tau),$ $\ldots, \psi_5(\tau)
\in L_2([t,T]).$ Furthermore, let the conditions
{\rm (\ref{novemberxxx3})--(\ref{novemberxxx6})} are satisfied.
Then$,$ for the sum $\bar J^{*}[\psi^{(5)}]_{T,t}^{(i_1 \ldots i_5)}$
$(i_1,\ldots,i_5=0,1,\ldots,m)$
of iterated Ito stochastic integrals 
defined by {\rm (\ref{dsds9}) $(k=5)$}
the following 
expansion 

\vspace{-1mm}
$$
\bar J^{*}[\psi^{(5)}]_{T,t}^{(i_1 \ldots i_5)}=
\hbox{\vtop{\offinterlineskip\halign{
\hfil#\hfil\cr
{\rm l.i.m.}\cr
$\stackrel{}{{}_{p\to \infty}}$\cr
}} }\sum_{j_1,\ldots,j_5=0}^{p}
C_{j_5 \ldots j_1}\zeta_{j_1}^{(i_1)}\ldots \zeta_{j_5}^{(i_5)}
$$

\vspace{3mm}
\noindent
that converges in the mean-square sense is valid, where 

\vspace{-1mm}
$$
C_{j_5 \ldots j_1}=\int\limits_t^T \psi_5(t_5)\phi_{j_5}(t_5)
\ldots
\int\limits_t^{t_2}\psi_1(t_1)
\phi_{j_1}(t_1)dt_1\ldots dt_5
$$

\vspace{3mm}
and
$$
\zeta_{j}^{(i)}=
\int\limits_t^T \phi_{j}(\tau) d{\bf w}_{\tau}^{(i)}
$$ 

\vspace{2mm}
\noindent
are independent standard Gaussian random variables for various 
$i$ or $j$ {\rm (}in the case when $i\ne 0${\rm ),}
${\bf w}_{\tau}^{(i)}={\bf f}_{\tau}^{(i)}$ for
$i=1,\ldots,m$ and 
${\bf w}_{\tau}^{(0)}=\tau.$}

\vspace{2mm}

{\bf Theorem~55}\ \cite{2018a}.\ {\it Suppose that
$\{\phi_j(x)\}_{j=0}^{\infty}$ is an arbitrary complete ortho\-nor\-mal system of 
func\-ti\-ons in the space $L_2([t,T])$ and $\psi_1(\tau),$ $\ldots, \psi_5(\tau)$
are continuous functions on $[t, T].$
Furthermore, let the conditions
{\rm (\ref{novemberxxx3})--(\ref{novemberxxx6})} are satisfied.
Then$,$ for the iterated Stra\-to\-no\-vich stochastic integral
of fifth multiplicity 

\vspace{-1mm}
$$
{\int\limits_t^{*}}^T \psi_5(t_5)
\ldots 
{\int\limits_t^{*}}^{t_2}\psi_1(t_1)
d{\bf w}_{t_1}^{(i_1)}
\ldots d{\bf w}_{t_5}^{(i_5)}\ \ \ (i_1,\ldots,i_5=0,1,\ldots,m)
$$

\vspace{3mm}
\noindent
the following 
expansion 

\vspace{-1mm}
$$
{\int\limits_t^{*}}^T \psi_5(t_5)
\ldots 
{\int\limits_t^{*}}^{t_2}\psi_1(t_1)
d{\bf w}_{t_1}^{(i_1)}
\ldots d{\bf w}_{t_5}^{(i_5)}=
\hbox{\vtop{\offinterlineskip\halign{
\hfil#\hfil\cr
{\rm l.i.m.}\cr
$\stackrel{}{{}_{p\to \infty}}$\cr
}} }\sum_{j_1,\ldots,j_5=0}^{p}
C_{j_5 \ldots j_1}\zeta_{j_1}^{(i_1)}\ldots \zeta_{j_5}^{(i_5)}
$$

\vspace{3mm}
\noindent
that converges in the mean-square sense is valid, where 
notations are the same as in Theorem~{\rm 54.}}

\vspace{2mm}

Note that Theorems~53 and 55 are simple consequences of Theorems~52 and 54, respectively
(see Theorem~5 $(k=4,\ 5).$

\vspace{5mm}

\section{On the Calculation of Matrix Traces of Volterra--Type Integral Operators}

\vspace{5mm}

It is easy to see that the function (\ref{ppp}) for even $k=2r$ $(r\in\mathbb{N})$
forms a family of integral operators
$\mathbb{K}: L_2([t, T]^r) \rightarrow L_2([t, T]^r)$
(with the kernel (\ref{ppp}))
of the form

\vspace{1mm}
\begin{equation}
\label{july6999}
\left(\mathbb{K} f\right)(t_{g_1},\ldots,t_{g_r})=
\int\limits_{[t, T]^{r}}K(t_1,\ldots,t_k)f(t_{g_{r+1}}, \ldots, t_{g_k})
dt_{g_{r+1}} \ldots dt_{g_k},
\end{equation}

\vspace{4mm}
\noindent
where $\{g_1,\ldots,g_k\}=\{1,\ldots,k\},$
the kernel $K(t_1,\ldots,t_k)$ is defined by (\ref{ppp}), i.e. has
the form

\begin{equation}
\label{july7000}
K(t_1,\ldots,t_k)=
\begin{cases}
\psi_1(t_1)\ldots \psi_k(t_k)\ &\hbox{for}\ \ t_1<\ldots<t_k\\
~\\
0\ &\hbox{otherwise}
\end{cases},
\end{equation}

\vspace{4mm}
\noindent
where $\psi_1(\tau),\ldots,\psi_k(\tau)\in L_2([t,T])$,
$t_1,\ldots,t_k\in [t, T]$ $(k\ge 2)$ and 
$K(t_1)\equiv\psi_1(t_1)$ for $t_1\in[t, T].$

For example,
\begin{equation}
\label{july7013}
\left(\mathbb{K} f\right)(t_2)=
\int\limits_t^T K(t_1,t_2)f(t_1)dt_1
=\psi_2(t_2)\int\limits_t^{t_2}\psi_1(t_1)f(t_1)dt_1,
\end{equation}

\vspace{3mm}
$$
\left(\mathbb{K} f\right)(t_3,t_4)=
\int\limits_{[t, T]^{2}}K(t_1,\ldots,t_4)f(t_1, t_2)
dt_1 dt_2
=
$$

$$
=
{\bf 1}_{\{t_3<t_4\}}\psi_3(t_3)\psi_4(t_4)
\int\limits_t^{t_3} \psi_2(t_2)\int\limits_t^{t_2}\psi_1(t_1)
f(t_1, t_2)
dt_1 dt_2,
$$

\vspace{4mm}
$$
\left(\mathbb{K} f\right)(t_1,t_2)=
\int\limits_{[t, T]^{2}}K(t_1,\ldots,t_4)f(t_3, t_4)
dt_3 dt_4=
$$

$$
=
\psi_1(t_1)\psi_2(t_2){\bf 1}_{\{t_1<t_2\}}
\int\limits_{t_2}^{T} \psi_3(t_3)\int\limits_{t_3}^{T}\psi_4(t_4)
f(t_3, t_4)
dt_4 dt_3.
$$

\vspace{4mm}

The simplest representative of the family (\ref{july6999})
has the form
\begin{equation}
\label{july7003}
\left(\mathbb{V} f\right)(x)=\int\limits_0^x f(\tau)d\tau
\end{equation}

\vspace{2mm}
\noindent
and is called the Volterra integral operator, where $\mathbb{V}: L_2([0,1]) \rightarrow L_2([0,1]),$
$f(\tau)\in L_2([0,1]).$
The kernel of the Volterra integral operator has the following form

\vspace{-1mm}
$$
K(\tau,x)=
\left\{\begin{matrix}
1,\ &\tau < x\cr\cr
0,\ &\hbox{\rm otherwise}
\end{matrix}
\right.,\ \ \ \tau, x\in [0, 1].
$$

\vspace{3mm}

Suppose that $\mathbb{A}: H \rightarrow H$ is a linear bounded operator. 
Recall \cite{gohb} that $\mathbb{A}$ has a finite matrix trace
if for any orthonormal basis $\left\{\Psi_j(x)\right\}_{j=0}^{\infty}$
of the space $H$ the series

\vspace{-1mm}
\begin{equation}
\label{july7002}
\sum_{j=0}^{\infty} \left\langle
\mathbb{A}\Psi_j, \Psi_j\right\rangle_H
\end{equation}

\vspace{3mm}
\noindent
converges, where $\left\langle
\cdot , \cdot \right\rangle_H$ is a scalar probuct in $H$.

Note that the series (\ref{july7002}) converges absolutely
since its sum does not depend on the permutation of the terms
of the series (\ref{july7002})
(any permutation of basis functions $\Psi_j(x)$ forms a basis 
in $H)$ \cite{gohb}.

It is well known that the Volterra integral operator 
(\ref{july7003}) is not a trace class operator since
its singular values are equal to \cite{Brisl} 

\vspace{-1mm}
$$
s_j(\mathbb{V})=\frac{2}{\pi(2j+1)}.
$$

\vspace{3mm}

On the other hand, it is known \cite{Brisl} that for trace class
operators the equality of matrix and integral traces holds.
It turns out that for the Volterra integral operator 
(\ref{july7003}) (although it is not a trace class operator),
the equality of matrix and integral traces is also true
\cite{Brisl}.

Thus, one cannot count on the fact that operators of the more
general form (\ref{july6999}) (from the same 
family of operators as the Volterra integral operator (\ref{july7003})) 
are operators of the trace class.
Nevertheless, the proof of the equalities 
of matrix and integral traces 
for Volterra--type integral operators (\ref{july6999}) (which is obviously a problem) 
provides a way
to calculate the matrix traces of these operators.

Why do we talk so much in this section about matrix
traces of operators from the family (\ref{july6999})?
The point is that matrix traces of operators of the form
(\ref{july6999}) are of great importance for obtaining
of expansions
of iterated Stratonovich stochastic integrals.

Throughout this article, we have already 
considered the matrix traces mentioned above
(see the formulas (\ref{5t}), 
(\ref{sixsix8})--(\ref{sixsix15}),
(\ref{sixsix129}), 
(\ref{2023novem206})--(\ref{2023novem208}), 
(\ref{march10})--(\ref{march12}), 
(\ref{july13})--(\ref{july18})).

Let us consider some illustrative examples.
We have

\vspace{-1mm}
\begin{equation}
\label{july7015}
\sum_{j_1=0}^{\infty} \left\langle
\mathbb{K}\phi_{j_1}, \phi_{j_1}\right\rangle_{L_2([t,T])}=
\end{equation}

\begin{equation}
\label{july7020}
=
\sum_{j_1=0}^{\infty}
\int\limits_t^T
\psi_2(t_2)\phi_{j_1}(t_2)\int\limits_t^{t_2}
\psi_1(t_1)\phi_{j_1}(t_1)dt_1 dt_2=
\sum_{j_1=0}^{\infty}C_{j_1 j_1},
\end{equation}

\vspace{3mm}

\begin{equation}
\label{july7016}
\sum_{j_1,j_2=0}^{\infty} \left\langle
\mathbb{K}\Psi_{j_1 j_2}, \Psi_{j_1 j_2}\right\rangle_{L_2([t,T]^2)}=
\end{equation}

$$
=
\sum_{j_1,j_2=0}^{\infty}
\int\limits_t^T
\psi_4(t_4)\phi_{j_2}(t_4)\int\limits_t^{t_4}
\psi_3(t_3)\phi_{j_2}(t_3)\int\limits_t^{t_3}
\psi_2(t_2)\phi_{j_1}(t_2)
\int\limits_t^{t_2}
\psi_1(t_1)\phi_{j_1}(t_1)
dt_1 dt_2 dt_3 dt_4=
$$

\begin{equation}
\label{july7021}
=\sum_{j_1,j_2=0}^{\infty}
C_{j_2 j_2 j_1 j_1},
\end{equation}

\vspace{5mm}
\noindent
where $\left\{\Psi_{j_1 j_2}(x,y)\right\}_{j_1,j_2=0}^{\infty}=
\left\{\phi_{j_1}(x)\phi_{j_2}(y)\right\}_{j_1,j_2=0}^{\infty},$
$\left\{\phi_{j}(x)\right\}_{j=0}^{\infty}$
is an arbitrary complete orthonormal system of functions
in $L_2([t, T]),$ $\left(\mathbb{K}f\right)(t_2)$ in (\ref{july7015})
is defined by (\ref{july7013}),
and $\left(\mathbb{K}f\right)(t_2, t_3)$ in (\ref{july7016}) has the following form

\vspace{-1mm}
$$
\left(\mathbb{K} f\right)(t_2,t_3)=
\int\limits_{[t, T]^{2}}K(t_1,\ldots,t_4)f(t_1, t_4)
dt_1 dt_4=
$$

$$
=\psi_2(t_2)\psi_3(t_3){\bf 1}_{\{t_2<t_3\}}
\int\limits_t^{t_2} \psi_1(t_1)\int\limits_{t_3}^{T}\psi_4(t_4)
f(t_1, t_4)
dt_4 dt_1,
$$

\vspace{4mm}
\noindent
where $K(t_1,\ldots,t_4)$ is defined by (\ref{july7000}).

The expressions on the right-hand sides of (\ref{july7020}) and (\ref{july7021})
were considered earlier in this article under various assumptions
on $\left\{\phi_j(x)\right\}_{j=0}^{\infty}$ and $\psi_1(\tau),\ldots,\psi_4(\tau)$
(see the formulas (\ref{5t}), 
(\ref{2023novem206}), (\ref{march10}), 
(\ref{july13})).

Let us consider one of the possible ways to calculate
matix traces of Volterra-type integral operators 
(\ref{july6999}) based Fubini's Theorem, Parseval's equality
and generalized Parseval's equality.

Recall the equalities (\ref{sixsix40}) and (\ref{2023novem215})

\vspace{1mm}
$$
C_{j_6 j_5 j_4 j_3 j_2 j_1}+C_{j_1 j_2 j_3 j_4 j_5 j_6}=
C_{j_6}C_{j_5 j_4 j_3 j_2 j_1}-C_{j_5 j_6}C_{j_4 j_3 j_2 j_1}+
$$

\begin{equation}
\label{july7100}
+C_{j_4 j_5 j_6}C_{j_3 j_2 j_1}-C_{j_3 j_4 j_5 j_6}C_{j_2 j_1}+
C_{j_2 j_3 j_4 j_5 j_6}C_{j_1},
\end{equation}

\vspace{1mm}
\begin{equation}
\label{july7101}
C_{j_4 j_3 j_2 j_1}+C_{j_1 j_2 j_3 j_4}=
C_{j_4}C_{j_3 j_2 j_1}-C_{j_3 j_4}C_{j_2 j_1}+
C_{j_2 j_3 j_4}C_{j_1},
\end{equation}

\vspace{5mm}
\noindent
where $C_{j_k\ldots j_1}$ is defined by
the formula

$$
C_{j_k \ldots j_1}=\int\limits_t^T\psi_k(t_k)\phi_{j_k}(t_k)\ldots
\int\limits_t^{t_2}
\psi_1(t_1)\phi_{j_1}(t_1)
dt_1\ldots dt_k\ \ \ (k\in\mathbb{N})
$$
  
\vspace{3mm}
\noindent
for the case $\psi_1(\tau),\ldots,\psi_k(\tau)\equiv 1$.

It is easy to see (see the derivation of
(\ref{sixsix40}) and (\ref{2023novem215})) that analogues of the relations
(\ref{july7100}), (\ref{july7101}) (with appropriate changes) hold for 
$\psi_1(\tau),\ldots,\psi_6(\tau)\in L_2([t, T]).$

By analogy with (\ref{july7100}), (\ref{july7101})
(see the derivation of
(\ref{sixsix40}) and (\ref{2023novem215})) we obtain for $k=2r$ $(r=2,3,4,\ldots)$

\vspace{1mm}
$$
C_{j_k j_{k-1}\ldots j_1}^{\psi_k \psi_{k-1} \ldots \psi_1} +
C_{j_1 j_2\ldots j_k}^{\psi_1 \psi_{2} \ldots \psi_k}=
C_{j_k}^{\psi_k} \cdot  C_{j_{k-1} j_{k-2}\ldots j_1}^{\psi_{k-1} \psi_{k-2} \ldots \psi_1}
-C_{j_{k-1} j_k}^{\psi_{k-1} \psi_{k}} \cdot C_{j_{k-2} j_{k-3}\ldots j_1}
^{\psi_{k-2} \psi_{k-3} \ldots \psi_1}+
$$

\vspace{1mm}
\begin{equation}
\label{july7028}
+C_{j_{k-2} j_{k-1} j_k}^{\psi_{k-2} \psi_{k-1} \psi_k} \cdot
C_{j_{k-3} j_{k-4}\ldots j_1}^{\psi_{k-3} \psi_{k-4} \ldots \psi_1}
- \ldots - C_{j_3 j_4 \ldots j_k}^{\psi_3 \psi_{4} \ldots \psi_k} \cdot C_{j_2 j_1}^{\psi_2 \psi_1}+
C_{j_2 j_3 \ldots j_k}^{\psi_2 \psi_{3} \ldots \psi_k} \cdot C_{j_1}^{\psi_1},
\end{equation}

\vspace{3mm}
\noindent
where 
\begin{equation}
\label{july40000}
C_{j_k \ldots j_1}^{\psi_k\ldots \psi_1}=\int\limits_t^T\psi_k(t_k)\phi_{j_k}(t_k)\ldots
\int\limits_t^{t_2}
\psi_1(t_1)\phi_{j_1}(t_1)
dt_1\ldots dt_k\ \ \ (k\in\mathbb{N}).
\end{equation}

\vspace{4mm}

When proving Theorem~46, 
using (\ref{july7028}) (the case $k=4,$ $\psi_1(\tau),\ldots,\psi_4(\tau)\equiv 1$),
we obtained the following formulas
$$
\lim\limits_{p\to\infty}
\sum\limits_{j_1, j_3=0}^{p}
C_{j_3 j_3 j_1 j_1}=\frac{1}{8}(T-t)^2,
$$

$$
\lim\limits_{p\to\infty}
\sum\limits_{j_1, j_3=0}^{p}
C_{j_1 j_3 j_3 j_1}=0,
$$

\vspace{1mm}
$$
\lim\limits_{p\to\infty}
\sum\limits_{j_1, j_2=0}^{p}
C_{j_2 j_1 j_2 j_1}=0,
$$

\vspace{5mm}
\noindent
where
$\{\phi_j(x)\}_{j=0}^{\infty}$ is an arbitrary complete orthonormal system of 
functions in the space $L_2([t,T])$ and we use the notation
$C_{j_k \ldots j_1}$ instead of 
$C_{j_k \ldots j_1}^{\psi_k \ldots \psi_1}$ for the case $\psi_1(\tau),\ldots,\psi_k(\tau)\equiv 1.$

In principle, using (\ref{july7028}), we can 
calculate any matrix traces for which the following
symmetry condition 

\vspace{-3mm}
\begin{equation}
\label{july7030}
\psi_1(\tau)=\psi_k(\tau),\ \ \psi_2(\tau)=\psi_{k-1}(\tau),\ \ldots
,\ \ \psi_r(\tau)=\psi_{r+1}(\tau)\ \ \ (k=2r,\ r=2,3,4,\ldots)
\end{equation}

\vspace{4mm}
\noindent
is satisfied. Obviously, the case
$\psi_1(\tau),\ldots,\psi_k(\tau)\equiv 1$ 
is possible since it is a special case of (\ref{july7030}).
This case is important because it covers
the mean-square approximation of iterated Stratonovich
stochastic integrals from the classical Taylor--Stratonovich
expansions (see \cite{2018a}, Chapter~4).

Consider the case $k=4$ of (\ref{july7028})

\vspace{-1mm}
\begin{equation}
\label{july5000}
C_{j_4 j_3 j_2 j_1}^{\psi_4 \psi_3 \psi_2 \psi_1} + C_{j_1 j_2 j_3 j_4}^{\psi_1 \psi_2 \psi_3 \psi_4}=
C_{j_4}^{\psi_4}C_{j_3 j_2 j_1}^{\psi_3 \psi_2 \psi_1}-C_{j_3 j_4}^{\psi_3 \psi_4}
C_{j_2 j_1}^{\psi_2 \psi_1}+
C_{j_2 j_3 j_4}^{\psi_2 \psi_3 \psi_4}C_{j_1}^{\psi_1},
\end{equation}

\vspace{4mm}
\noindent
where $\psi_1(\tau), \ldots, \psi_4(\tau)\in L_2([t,T]).$

Substitute $j_4=j_3,$ $j_2=j_1$ into (\ref{july5000})

\vspace{-1mm}
\begin{equation}
\label{july5001}
C_{j_3 j_3 j_1 j_1}^{\psi_4 \psi_3 \psi_2 \psi_1} + C_{j_1 j_1 j_3 j_3}^{\psi_1 \psi_2 \psi_3 \psi_4}=
C_{j_3}^{\psi_4}C_{j_3 j_1 j_1}^{\psi_3 \psi_2 \psi_1}-C_{j_3 j_3}^{\psi_3 \psi_4}
C_{j_1 j_1}^{\psi_2 \psi_1}+
C_{j_1 j_3 j_3}^{\psi_2 \psi_3 \psi_4}C_{j_1}^{\psi_1},
\end{equation}

\vspace{4mm}

Applying (\ref{july5001}), we get

\vspace{-1mm}
$$
\lim\limits_{p\to\infty}\sum\limits_{j_1,j_3=0}^{p}
\left(C_{j_3 j_3 j_1 j_1}^{\psi_4 \psi_3 \psi_2 \psi_1} +
C_{j_1 j_1 j_3 j_3}^{\psi_1 \psi_2 \psi_3 \psi_4}\right)=
\lim\limits_{p\to\infty}\sum\limits_{j_1,j_3=0}^{p}
C_{j_3}^{\psi_4}C_{j_3 j_1 j_1}^{\psi_3 \psi_2 \psi_1}-
$$

\begin{equation}
\label{july5002}
-
\lim\limits_{p\to\infty}\sum\limits_{j_1,j_3=0}^{p}
C_{j_3 j_3}^{\psi_3 \psi_4}
C_{j_1 j_1}^{\psi_2 \psi_1}+
\lim\limits_{p\to\infty}\sum\limits_{j_1,j_3=0}^{p}
C_{j_1 j_3 j_3}^{\psi_2 \psi_3 \psi_4}C_{j_1}^{\psi_1}.
\end{equation}

\vspace{2mm}

From (\ref{after1400}) we have

$$
\lim\limits_{p\to\infty}\sum\limits_{j_3=0}^{p} C_{j_3 j_3}^{\psi_3 \psi_4}
\sum\limits_{j_1=0}^{p}C_{j_1 j_1}^{\psi_2 \psi_1}=
\lim\limits_{p\to\infty}\sum\limits_{j_3=0}^{p} C_{j_3 j_3}^{\psi_3 \psi_4}
\lim\limits_{p\to\infty}\sum\limits_{j_1=0}^{p} C_{j_1 j_1}^{\psi_2 \psi_1}=
$$

\begin{equation}
\label{july5003}
=
\frac{1}{4}\int\limits_t^T \psi_4(s) \psi_3(s)ds
\int\limits_t^T \psi_2(s) \psi_1(s)ds.
\end{equation}

\vspace{1mm}

Further, we obtain

\vspace{-1mm}
$$
\lim\limits_{p\to\infty}\sum\limits_{j_3=0}^{p} C_{j_3}^{\psi_4} 
\sum\limits_{j_1=0}^{p} C_{j_3 j_1 j_1}^{\psi_3 \psi_2 \psi_1}=
\frac{1}{2}\lim\limits_{p\to\infty}
\sum\limits_{j_3=0}^{p}C_{j_3}^{\psi_4} 
C_{j_3 j_1 j_1}^{\psi_3 \psi_2 \psi_1}\biggl|_{(j_{1} j_{1})\curvearrowright (\cdot)}-
$$

\begin{equation}
\label{july5004}
-
\lim\limits_{p\to\infty}
\sum\limits_{j_3=0}^{p}C_{j_3}^{\psi_4}
\left(\frac{1}{2}
C_{j_3 j_1 j_1}^{\psi_3 \psi_2 \psi_1}\biggl|_{(j_{1} j_{1})\curvearrowright (\cdot)}-
\sum\limits_{j_1=0}^{p} C_{j_3 j_1 j_1}^{\psi_3 \psi_2 \psi_1}\right).
\end{equation}

\vspace{4mm}

Applying the generalized Parseval equality, we have

$$
\lim\limits_{p\to\infty}
\sum\limits_{j_3=0}^{p}C_{j_3}^{\psi_4}
C_{j_3 j_1 j_1}^{\psi_3 \psi_2 \psi_1}\biggl|_{(j_{1} j_{1})\curvearrowright (\cdot)}=
\lim\limits_{p\to\infty}
\sum\limits_{j_3=0}^{p}\int\limits_t^T \psi_4(s)\phi_{j_3}(s)ds
\int\limits_t^T \psi_3(s)\phi_{j_3}(s)\int\limits_t^{s}
\psi_2(\tau)\psi_1(\tau)d\tau ds
=
$$

\begin{equation}
\label{july5005}
=
\int\limits_t^T \psi_4(s) \psi_3(s)
\int\limits_t^{s} \psi_2(\tau) \psi_1(\tau)d\tau ds.
\end{equation}

\vspace{4mm}

From (\ref{july5004}) and (\ref{july5005}) we obtain

\vspace{-1mm}
$$
\lim\limits_{p\to\infty}\sum\limits_{j_3=0}^{p} C_{j_3}^{\psi_4} 
\sum\limits_{j_1=0}^{p} C_{j_3 j_1 j_1}^{\psi_3 \psi_2 \psi_1}=
\frac{1}{2}\int\limits_t^T \psi_4(s) \psi_3(s)
\int\limits_t^{s} \psi_2(\tau) \psi_1(\tau)d\tau ds
-
$$

\begin{equation}
\label{july5006}
-\lim\limits_{p\to\infty}
\sum\limits_{j_3=0}^{p}C_{j_3}^{\psi_4}
\left(\frac{1}{2}
C_{j_3 j_1 j_1}^{\psi_3 \psi_2 \psi_1}\biggl|_{(j_{1} j_{1})\curvearrowright (\cdot)}-
\sum\limits_{j_1=0}^{p} C_{j_3 j_1 j_1}^{\psi_3 \psi_2 \psi_1}\right).
\end{equation}

\vspace{4mm}

Due to Cauchy--Bunyakovsky's inequality, Parseval's equality
and (\ref{july1003}), we get 

\vspace{-1mm}
$$
\lim\limits_{p\to\infty}
\left(\sum\limits_{j_3=0}^{p}C_{j_3}^{\psi_4}
\left(\frac{1}{2}
C_{j_3 j_1 j_1}^{\psi_3 \psi_2 \psi_1}\biggl|_{(j_{1} j_{1})\curvearrowright (\cdot)}-
\sum\limits_{j_1=0}^{p} C_{j_3 j_1 j_1}^{\psi_3 \psi_2 \psi_1}\right)\right)^2\le
$$

$$
\le \lim\limits_{p\to\infty}
\sum\limits_{j_3=0}^{p}\left(C_{j_3}^{\psi_4}\right)^2\
\sum\limits_{j_3=0}^{p}
\left(\frac{1}{2}
C_{j_3 j_1 j_1}^{\psi_3 \psi_2 \psi_1}\biggl|_{(j_{1} j_{1})\curvearrowright (\cdot)}-
\sum\limits_{j_1=0}^{p} C_{j_3 j_1 j_1}^{\psi_3 \psi_2 \psi_1}\right)^2\le
$$

$$
\le \lim\limits_{p\to\infty}
\sum\limits_{j_3=0}^{\infty}\left(C_{j_3}^{\psi_4}\right)^2\
\sum\limits_{j_3=0}^{p}
\left(\frac{1}{2}
C_{j_3 j_1 j_1}^{\psi_3 \psi_2 \psi_1}\biggl|_{(j_{1} j_{1})\curvearrowright (\cdot)}-
\sum\limits_{j_1=0}^{p} C_{j_3 j_1 j_1}^{\psi_3 \psi_2 \psi_1}\right)^2=
$$

\begin{equation}
\label{july5007}
=\int\limits_t^T \psi_4^2(s)ds\lim\limits_{p\to\infty}
\sum\limits_{j_3=0}^{p}
\left(\frac{1}{2}
C_{j_3 j_1 j_1}^{\psi_3 \psi_2 \psi_1}\biggl|_{(j_{1} j_{1})\curvearrowright (\cdot)}-
\sum\limits_{j_1=0}^{p} C_{j_3 j_1 j_1}^{\psi_3 \psi_2 \psi_1}\right)^2=0.
\end{equation}

\vspace{4mm}         

Combining (\ref{july5006}) and (\ref{july5007}), we obtain

\vspace{-1mm}
\begin{equation}
\label{july5008}
\lim\limits_{p\to\infty}\sum\limits_{j_3=0}^p C_{j_3}^{\psi_4} 
\sum\limits_{j_1=0}^p C_{j_3 j_1 j_1}^{\psi_3 \psi_2 \psi_1}=
\frac{1}{2}\int\limits_t^T \psi_4(s) \psi_3(s)
\int\limits_t^{s} \psi_2(\tau) \psi_1(\tau)d\tau ds.
\end{equation}

\vspace{4mm}

Absolutely similarly to (\ref{july5008}) we get

\vspace{-1mm}
\begin{equation}
\label{july5009}
\lim\limits_{p\to\infty}\sum\limits_{j_1=0}^p C_{j_1}^{\psi_1}
\sum\limits_{j_3=0}^p C_{j_1 j_3 j_3}^{\psi_2 \psi_3 \psi_4}=
\frac{1}{2}\int\limits_t^T \psi_2(s) \psi_1(s)
\int\limits_t^{s} \psi_3(\tau) \psi_4(\tau)d\tau ds.
\end{equation}

\vspace{4mm}

Combining (\ref{july5002}), (\ref{july5003}), (\ref{july5008}), (\ref{july5009}) and 
applying Fubini's Theorem, we have

\vspace{-1mm}
$$
\lim\limits_{p\to\infty}\sum\limits_{j_1,j_3=0}^p
\left(C_{j_3 j_3 j_1 j_1}^{\psi_4 \psi_3 \psi_2 \psi_1} +
C_{j_1 j_1 j_3 j_3}^{\psi_1 \psi_2 \psi_3 \psi_4}\right)=
\frac{1}{2}\int\limits_t^T \psi_4(s) \psi_3(s)
\int\limits_t^{s} \psi_2(\tau) \psi_1(\tau)d\tau ds+
$$

$$
+
\frac{1}{2}\int\limits_t^T \psi_2(s) \psi_1(s)
\int\limits_t^{s} \psi_3(\tau) \psi_4(\tau)d\tau ds
-\frac{1}{4}\int\limits_t^T \psi_4(s) \psi_3(s)ds
\int\limits_t^T \psi_2(s) \psi_1(s)ds=
$$

$$
=\frac{1}{4}\int\limits_t^T \psi_4(s) \psi_3(s)ds
\int\limits_t^T \psi_2(s) \psi_1(s)ds=
$$

\begin{equation}
\label{july5010}
=\frac{1}{4}\int\limits_t^T \psi_4(s) \psi_3(s)
\int\limits_t^{s} \psi_2(\tau) \psi_1(\tau)d\tau ds+
\frac{1}{4}\int\limits_t^T \psi_2(s) \psi_1(s)
\int\limits_t^{s} \psi_3(\tau) \psi_4(\tau)d\tau ds.
\end{equation}

\vspace{4mm}

Let us rewrite (\ref{july5010}) in the form

\vspace{-1mm}
$$
\lim\limits_{p\to\infty}\sum\limits_{j_1,j_3=0}^p
\left(C_{j_3 j_3 j_1 j_1}^{\psi_4 \psi_3 \psi_2 \psi_1} +
C_{j_3 j_3 j_1 j_1}^{\psi_1 \psi_2 \psi_3 \psi_4}\right)=
$$

\begin{equation}
\label{july5010x}
=\frac{1}{4}\int\limits_t^T \psi_4(s) \psi_3(s)
\int\limits_t^{s} \psi_2(\tau) \psi_1(\tau)d\tau ds+
\frac{1}{4}\int\limits_t^T \psi_2(s) \psi_1(s)
\int\limits_t^{s} \psi_3(\tau) \psi_4(\tau)d\tau ds.
\end{equation}

\vspace{3mm}

It is easy to see the left-hand side 
of (\ref{july5010x}) does not depend on 
the simultaneous rearrangement of $\psi_4$ with $\psi_1$
and $\psi_3$ with $\psi_2$.

Using the above arguments and using derivation method of (\ref{march11}) and (\ref{march12}), we get

\vspace{-1mm}
\begin{equation}
\label{july8000}
\lim\limits_{p\to\infty}\sum\limits_{j_1,j_3=0}^p
\left(C_{j_3 j_1 j_3 j_1}^{\psi_4 \psi_3 \psi_2 \psi_1} +
C_{j_3 j_1 j_3 j_1}^{\psi_1 \psi_2 \psi_3 \psi_4}\right)=0,
\end{equation}

\begin{equation}
\label{july8001}
\lim\limits_{p\to\infty}\sum\limits_{j_1,j_3=0}^p
\left(C_{j_1 j_3 j_3 j_1}^{\psi_4 \psi_3 \psi_2 \psi_1} +
C_{j_1 j_3 j_3 j_1}^{\psi_1 \psi_2 \psi_3 \psi_4}\right)=0.
\end{equation}

\vspace{4mm}

Using (\ref{july5010x})--(\ref{july8001}) under the conditions
$\psi_1(\tau)=\psi_4(\tau),$ $\psi_2(\tau)=\psi_3(\tau),$ we obtain

\vspace{-1mm}
$$
\lim\limits_{p\to\infty}\sum\limits_{j_1,j_3=0}^p
C_{j_3 j_3 j_1 j_1}^{\psi_1 \psi_2 \psi_2 \psi_1} 
=\frac{1}{4}\int\limits_t^T \psi_2(s) \psi_1(s)
\int\limits_t^{s} \psi_2(\tau) \psi_1(\tau)d\tau ds,
$$

$$
\lim\limits_{p\to\infty}\sum\limits_{j_1,j_3=0}^p
C_{j_3 j_1 j_3 j_1}^{\psi_1 \psi_2 \psi_2 \psi_1}=0,
$$

$$
\lim\limits_{p\to\infty}\sum\limits_{j_1,j_3=0}^p
C_{j_1 j_3 j_3 j_1}^{\psi_1 \psi_2 \psi_2 \psi_1}=0.
$$

\vspace{4mm}

An efficient method for calculating of matrix traces 
of Volterra--type integral operators of the form (\ref{july6999})
was proposed in \cite{rybakov5000}.
This method is based on Theorem~3.1 from \cite{Brisl}.
Theorem~3.1 \cite{Brisl} implies the following statement.

\vspace{2mm}

{\bf Theorem~A}\ (see \cite{Brisl} for details).\ {\it
Let $\mathbb{K}: L_2([t,T]^r)\rightarrow L_2([t,T]^r)$
$(r\in{\mathbb N})$ be a trace class operator with the kernel
$K\in L_2([t,T]^{2r}).$
Then
$\tilde K(t_1,\ldots,t_r,t_1,\ldots,t_r)$ exists
almost everywhere $[dt_1\ldots dt_r]$ and

\begin{equation}
\label{july15000}
tr\mathbb{K}=\int\limits_{[t,T]^r}\tilde K(t_1,\ldots,t_r,t_1,\ldots,t_r)
dt_1\ldots dt_r,
\end{equation}

\vspace{4mm}
\noindent
where 

\vspace{-1mm}
$$
\tilde F(x_1,\ldots,x_m)\stackrel{\sf def}{=}\lim\limits_{u\to 0}
A_u F(x_1,\ldots,x_m),
$$

\vspace{4mm}
$$
A_u F(x_1,\ldots,x_m)\stackrel{\sf def}{=}
\frac{1}{\left(2u\right)^{m}}
\int\limits_{[-u,u]^{m}}
F(x_1+\tau_1,\ldots,x_m+\tau_m)
d\tau_1\ldots d\tau_m\ \ \ (m\in\mathbb{N}).
$$ 
}

\vspace{3mm}

Let us consider the following statements.

\vspace{2mm}

{\bf Theorem~B} (\cite{goldberg}, P.~71).\ {\it 
Let $\mathbb{K}: L_2([t,T])\rightarrow L_2([t,T])$
be an integral operator defined by

$$
\left(\mathbb{K}f\right)(\tau)=\int\limits_t^T K(\tau,s)f(s)ds,
$$

\vspace{2.5mm}
\noindent
where the kernel $K(\tau,s)$ is continuous on $[t, T]\times[t, T]$
and satisfies the condition

\vspace{-1mm}
\begin{equation}
\label{july11002}
\left\vert K(\tau,s_2)-K(\tau,s_1)\right\vert \le C \left\vert s_2-s_1\right\vert^{\alpha},
\end{equation}

\vspace{3mm}
\noindent
where $0<\alpha\le 1$. If$,$ in addition$,$ $\mathbb{K}$ is a Hermitian
operator and $\alpha>1/2,$ then
$\mathbb{K}$ is a trace class operator.}

\vspace{2mm}

{\bf Theorem~C} (\cite{goldberg}, Theorem~5.6).\ {\it 
Let $\mathbb{K}: H\rightarrow H$ be a trace class operator.
Then

\vspace{-1mm}
\begin{equation}
\label{july11201}
tr\mathbb{A}=\sum_{j=0}^{\infty} \left\langle
\mathbb{A}\phi_j, \phi_j\right\rangle_H
\end{equation}

\vspace{3mm}
\noindent
for any orthonormal basis $\left\{\phi_j(x)\right\}_{j=0}^{\infty}$ of $H$.}

\vspace{2mm}

Consider an integral operator $\mathbb{K'}: L_2([t,T])\rightarrow L_2([t,T])$
defined by the equality

\vspace{-1mm}
$$
\left(\mathbb{K'}f\right)(\tau)=\int\limits_t^T K'(\tau,s)f(s)ds,
$$

\vspace{3mm}
\noindent
where the continuous kernel $K'(\tau,s)$ has the form 

\vspace{-1mm}
\begin{equation}
\label{ziko5001}
K'(t_1,t_2)=\left\{
\begin{matrix}
\psi_2(t_1)\psi_1(t_2),\ \ t_1\ge t_2\cr\cr\cr
\psi_1(t_1)\psi_2(t_2),\ \ t_1\le t_2
\end{matrix}
\right.\ \ \ (t_1,t_2\in[t,T])
\end{equation}

\vspace{4mm}
\noindent
and $\psi_1(\tau), \psi_2(\tau)$ are continuously differentiable functions on $[t, T].$
Recall that (see \cite{2018a}, Sect.~2.1.2)

\begin{equation}
\label{july11000}
\left|K'(t_2,s_2)-K'(t_1,s_1)\right|\le L
\left(|t_2-t_1|+|s_2-s_1|\right),
\end{equation}

\vspace{4mm}
\noindent
where $L<\infty$ and $(t_1,s_1)$, $(t_2,s_2)\in [t, T]^2.$
Let us substitute $t_1=t_2=\tau$ into (\ref{july11000})

\begin{equation}
\label{july11001}
\left|K'(\tau,s_2)-K'(\tau,s_1)\right|\le L
|s_2-s_1|.
\end{equation}

\vspace{4mm}

Thus, the condition (\ref{july11002}) is fulfilled $(\alpha=1$).
Further, using Fubini's Theorem, we have

\vspace{-1mm}
$$
\left\langle
\mathbb{K'}x, y\right\rangle_{L_2([t,T])}=
\int\limits_t^T\psi_2(t_2)y(t_2)\int\limits_t^{t_2}\psi_1(t_1)x(t_1)dt_1 dt_2+
\int\limits_t^T\psi_1(t_2)y(t_2)\int\limits_{t_2}^T\psi_2(t_1)x(t_1)dt_1 dt_2=
$$

\vspace{2mm}
\begin{equation}
\label{july11207}
=\int\limits_t^T\psi_1(t_1)x(t_1) \int\limits_{t_1}^T \psi_2(t_2)y(t_2) dt_2 dt_1+
\int\limits_t^T\psi_2(t_1)x(t_1) \int\limits_{t}^{t_2} \psi_1(t_2)y(t_2) dt_2 dt_1=
\left\langle\mathbb{K'}y, x\right\rangle_{L_2([t,T])}.
\end{equation}

\vspace{4mm}
\noindent
The conditions of Theorem~B are fulfilled. Then, $\mathbb{K'}$ is a trace class operator.

Let us prove the equality (\ref{july13})
using the method from \cite{rybakov5000} in our interpretation.
Consider two symmetric functions of the form (\ref{ziko5001})

\vspace{-1mm}
\begin{equation}
\label{july15004}
K'(t_1,t_2)=\psi_1(t_1)f_2(t_2){\bf 1}_{\{t_1\le t_2\}}+
\psi_1(t_2)f_2(t_1){\bf 1}_{\{t_1\ge t_2\}},
\end{equation}

\vspace{-2mm}
\begin{equation}
\label{july15005}
K''(t_3,t_4)=f_3(t_3)\psi_4(t_4){\bf 1}_{\{t_3\le t_4\}}+
f_3(t_4)\psi_4(t_3){\bf 1}_{\{t_3\ge t_4\}},
\end{equation}

\vspace{4mm}
\noindent
where we suppose that $\psi_1(\tau), \psi_4(\tau)$ are continuously differentiable
functions on $[t, T]$ (the case 
$\psi_1(\tau), \psi_4(\tau)\in L_2([t,T])$ will be considered further)
and $f_2(\tau), f_3(\tau)$  are polynomials of finite degrees.
As noted above, the 
kernels $K'(t_1,t_2)$ and $K''(t_3,t_4)$ (see (\ref{july15004}), (\ref{july15005}))
correspond to the 
trace class integral operators.

It is known \cite{Brisl} that the integral operator $\mathbb{A}$ is a trace class operator
if and only if the kernel $K(x,y)$ of $\mathbb{A}$ has the following
representation

\vspace{-1mm}
\begin{equation}
\label{july15007}
K(x,y)=\int\limits_{[t,T]^{2n}} K_1(x,\tau)K_2(\tau,y)d\tau
\end{equation}

\vspace{3mm}
\noindent
almost everywhere $[dxdy]$,
where $K_1(x,y), K_2(x,y)$ are kernels of Hilbert--Schmidt operators,
$x, y\in \mathbb{R}^n$ $(n\ge 1).$

Since $K'(t_1,t_2)$ and $K''(t_3,t_4)$ are kernels of
the trace class integral operators, then (see (\ref{july15007}))

\vspace{-1mm}
\begin{equation}
\label{july15008}
K'(t_1,t_2)=\int\limits_{[t,T]} K_1'(t_1,\tau)K_2'(\tau,t_2)d\tau,\ \ \
K''(t_1,t_2)=\int\limits_{[t,T]} K_1''(t_1,\tau)K_2''(\tau,t_2)d\tau
\end{equation}

\vspace{3mm}
\noindent
almost everywhere $[dt_1 dt_2],$ where $K_1', K_2', K_1'', K_2''\in L_2([t, T]^2)$.
Then, we have

$$
K'(t_1,t_2)K''(t_3,t_4)
=\int\limits_{[t,T]} K_1'(t_1,\tau_1)K_2'(\tau_1,t_2)d\tau_1
\int\limits_{[t,T]} K_1''(t_3,\tau_2)K_2''(\tau_2,t_4)d\tau_2=
$$

\begin{equation}
\label{july15010}
=\int\limits_{[t,T]^2} K_1'(t_1,\tau_1)K_1''(t_3,\tau_2)
K_2'(\tau_1,t_2)
K_2''(\tau_2,t_4) d\tau_1 d\tau_2.
\end{equation}

\vspace{3mm}

The equality (\ref{july15010}) can be written as follows

\vspace{-1mm}
$$
F(t_1,t_3,t_2,t_4)=\int\limits_{[t,T]^2} F_1(t_1,t_3,\tau_1,\tau_2)
F_2(\tau_1,\tau_2,t_2,t_4) d\tau_1 d\tau_2
$$

\vspace{3mm}
\noindent
almost everywhere $[dt_1 dt_2 dt_3 dt_4]$, where
$F(t_1,t_3,t_2,t_4)=K'(t_1,t_2)K''(t_3,t_4),$
$F_1(t_1,t_3, \tau_1,\tau_2)=K_1'(t_1,\tau_1)K_1''(t_3,\tau_2),$
and $F_2(\tau_1,\tau_2,t_2,t_4)$ $=K_2'(\tau_1,t_2)K_2''(\tau_2,t_4).$

As a result, the product $K'(t_1,t_2)K''(t_3,t_4)$
is also the kernel of the trace class operator (see (\ref{july15007})).
Let us denote it by $\mathbb{K'}.$

Suppose that $\left\{\phi_{j}(x)\right\}_{j=0}^{\infty}$
is an arbitrary complete orthonormal system of functions
in $L_2([t, T]).$ Then 
$\left\{\Psi_{j_1 j_2}(x,y)\right\}_{j_1,j_2=0}^{\infty}=
\left\{\phi_{j_1}(x)\phi_{j_2}(y)\right\}_{j_1,j_2=0}^{\infty}$
is an orthonormal basis in $L_2([t, T]^2).$

Consider matrix trace of $\mathbb{K'}.$ Using Fubini's Theorem, we obtain

\vspace{-1mm}
$$
\sum_{j_1,j_2=0}^{\infty} \left\langle
\Psi_{j_1 j_2}, \mathbb{K'}\Psi_{j_1 j_2}\right\rangle_{L_2([t,T]^2)}=
$$

\vspace{1mm}
$$
=\sum_{j_1,j_2=0}^{\infty}\int\limits_{[t,T]^2}
\phi_{j_2}(t_4)\phi_{j_1}(t_1)
\int\limits_{[t,T]^2}K'(t_1,t_2)K''(t_3,t_4)
\phi_{j_2}(t_3)\phi_{j_1}(t_2)dt_2 dt_3 dt_1 dt_4=
$$

\vspace{1mm}
$$
=
\sum_{j_1,j_2=0}^{\infty}\left(
\int\limits_t^T
\psi_4(t_4)\phi_{j_2}(t_4)\int\limits_t^{T}
\psi_1(t_1)\phi_{j_1}(t_1)\int\limits_t^{t_4}
f_3(t_3)\phi_{j_2}(t_3)
\int\limits_{t_1}^T
f_2(t_2)\phi_{j_1}(t_2)
dt_2 dt_3 dt_1 dt_4+\right.
$$

\vspace{1mm}
$$
+
\int\limits_t^T
f_3(t_4)\phi_{j_2}(t_4)\int\limits_t^{T}
\psi_1(t_1)\phi_{j_1}(t_1)\int\limits_{t_4}^T
\psi_4(t_3)\phi_{j_2}(t_3)
\int\limits_{t_1}^T
f_2(t_2)\phi_{j_1}(t_2)
dt_2 dt_3 dt_1 dt_4+
$$

\vspace{1mm}
$$
+
\int\limits_t^T
\psi_4(t_4)\phi_{j_2}(t_4)\int\limits_t^{T}
f_2(t_1)\phi_{j_1}(t_1)\int\limits_t^{t_4}
f_3(t_3)\phi_{j_2}(t_3)
\int\limits_t^{t_1}
\psi_1(t_2)\phi_{j_1}(t_2)
dt_2 dt_3 dt_1 dt_4+
$$

\vspace{1mm}
$$
\left.+
\int\limits_t^T
f_2(t_1)\phi_{j_1}(t_1)\int\limits_t^{T}
\psi_4(t_3)\phi_{j_2}(t_3)\int\limits_t^{t_3}
f_3(t_4)\phi_{j_2}(t_4)
\int\limits_t^{t_1}
\psi_1(t_2)\phi_{j_1}(t_2)
dt_2 dt_4 dt_3 dt_1\right)=
$$

\vspace{1mm}

$$
=
\sum_{j_1,j_2=0}^{\infty}\left(
\int\limits_t^T
\psi_4(t_4)\phi_{j_2}(t_4)\int\limits_t^{t_4}
f_3(t_3)\phi_{j_2}(t_3)\int\limits_t^{T}
f_2(t_2)\phi_{j_1}(t_2)
\int\limits_t^{t_2}
\psi_1(t_1)\phi_{j_1}(t_1)
dt_1 dt_2 dt_3 dt_4+\right.
$$

\vspace{1mm}
$$
+
\int\limits_t^T
\psi_4(t_3)\phi_{j_2}(t_3)\int\limits_t^{t_3}
f_3(t_4)\phi_{j_2}(t_4)\int\limits_t^{T}
f_2(t_2)\phi_{j_1}(t_2)
\int\limits_t^{t_2}
\psi_1(t_1)\phi_{j_1}(t_1)
dt_1 dt_2 dt_4 dt_3+
$$

\vspace{1mm}
$$
+
\int\limits_t^T
\psi_4(t_4)\phi_{j_2}(t_4)\int\limits_t^{t_4}
f_3(t_3)\phi_{j_2}(t_3)\int\limits_t^{T}
f_2(t_1)\phi_{j_1}(t_1)
\int\limits_t^{t_1}
\psi_1(t_2)\phi_{j_1}(t_2)
dt_2 dt_1 dt_3 dt_4+
$$

\vspace{1mm}
$$
\left. +
\int\limits_t^T
\psi_4(t_3)\phi_{j_2}(t_3)\int\limits_t^{t_3}
f_3(t_4)\phi_{j_2}(t_4)\int\limits_t^{T}
f_2(t_1)\phi_{j_1}(t_1)
\int\limits_t^{t_1}
\psi_1(t_2)\phi_{j_1}(t_2)
dt_2 dt_1 dt_4 dt_3\right)=
$$

\vspace{1mm}
\begin{equation}
\label{july15012}
=4
\sum_{j_1,j_2=0}^{\infty}
\int\limits_t^T
\psi_4(t_4)\phi_{j_2}(t_4)\int\limits_t^{t_4}
f_3(t_3)\phi_{j_2}(t_3)\int\limits_t^{T}
f_2(t_2)\phi_{j_1}(t_2)
\int\limits_t^{t_2}
\psi_1(t_1)\phi_{j_1}(t_1)
dt_1 dt_2 dt_3 dt_4.
\end{equation}

\vspace{4mm}

According to (\ref{july15012}), (\ref{july15000}), and Theorem~C, we get

$$
\sum_{j_1,j_2=0}^{\infty} \left\langle
\Psi_{j_1 j_2}, \mathbb{K'}\Psi_{j_1 j_2}\right\rangle_{L_2([t,T]^2)}=
$$

\vspace{1mm}
$$
=4\sum_{j_1,j_2=0}^{\infty}
\int\limits_t^T
\psi_4(t_4)\phi_{j_2}(t_4)\int\limits_t^{t_4}
f_3(t_3)\phi_{j_2}(t_3)\int\limits_t^{T}
f_2(t_2)\phi_{j_1}(t_2)
\int\limits_t^{t_2}
\psi_1(t_1)\phi_{j_1}(t_1)dt_1 dt_2 dt_3 dt_4=
$$

\vspace{1mm}
$$
=\int\limits_{[t,T]^2} 
\lim\limits_{u\to 0}
A_u K'(t_2,t_2)K''(t_4,t_4)dt_2 dt_4=
$$

\vspace{1mm}
$$
=
\int\limits_{[t,T]^2} 
\lim\limits_{u\to 0}
A_u K'(t_2,t_2) \lim\limits_{u\to 0}
A_u K''(t_4,t_4)dt_2 dt_4=\int\limits_{[t,T]^2} 
K'(t_2,t_2)K''(t_4,t_4)dt_2 dt_4=
$$

\vspace{1mm}
\begin{equation}
\label{july15014}
=\int\limits_{[t,T]^2} 
\psi_4(t_4)f_3(t_4)f_2(t_2)\psi_1(t_2)dt_2 dt_4.
\end{equation}

\vspace{4mm}

Recall that 
$f_2(\tau)$ and $f_3(\tau)$ are polynomials of finite degrees.
For example, $f_2(\tau)$ and $f_3(\tau)$ can be Legendre polynomials
that form a complete orthonormal system of functions in $L_2([t,T]).$

Denote
\begin{equation}
\label{july15015}
s_q(t_2,t_3)=\sum\limits_{l_1,l_2=0}^q
C_{l_2 l_1}\bar \phi_{l_1}(t_2)\bar \phi_{l_2}(t_3),
\end{equation}

\vspace{2mm}
\noindent 
where $\left\{\bar \phi_j(x)\right\}_{j=0}^{\infty}$ 
is a complete orthonormal system of Legendre polynomials in $L_2([t,T])$
and
$C_{l_2 l_1}$ are Fourier--Legendre coefficients for the function
$g(t_2,t_3)=\psi_2(t_2)\psi_3(t_3){\bf 1}_{\{t_2<t_3\}}$ 
($\psi_2(\tau), \psi_3(\tau)\in L_2([t,T])),$ i.e.
$$
C_{l_2 l_1}=\int\limits_t^T \psi_3(t_3)\bar \phi_{l_2}(t_3)
\int\limits_t^{t_3}\psi_2(t_2)\bar \phi_{l_1}(t_2)dt_2 dt_3.
$$

\vspace{2mm}

Further, we have

$$
\lim\limits_{q\to\infty}
\int\limits_{[t,T]^2}
\left(s_q(t_2,t_3)-g(t_2,t_3)\right)^2 dt_2 dt_3=0\ \ \ \hbox{or}\ \ \ 
\lim\limits_{q\to\infty}\left\Vert s_q - g\right\Vert^2_{L_2([t,T]^2)}=0.
$$

\vspace{3mm}

From (\ref{july15014}) we obtain (the sum on the right-hand side of (\ref{july15015}) is finite)

$$
\sum_{j_1,j_2=0}^{\infty} \left\langle
\Psi_{j_1 j_2}, \mathbb{K'}_q\Psi_{j_1 j_2}\right\rangle_{L_2([t,T]^2)}=
$$

$$
=4\sum_{j_1,j_2=0}^{\infty}
\int\limits_{[t,T]^4}{\bf 1}_{\{t_1<t_2\}}{\bf 1}_{\{t_3<t_4\}}
\psi_4(t_4)\phi_{j_2}(t_4)
s_q(t_2,t_3)\phi_{j_2}(t_3)
\phi_{j_1}(t_2)
\psi_1(t_1)\phi_{j_1}(t_1)
dt_1 dt_2 dt_3 dt_4
=
$$
\begin{equation}
\label{july15017}
=\int\limits_{[t,T]^2} 
\psi_4(t_4)s_q(t_2,t_4)\psi_1(t_2)dt_2 dt_4,
\end{equation}

\vspace{3mm}
\noindent
where the operator $\mathbb{K'}_q$ (more precisely, its kernel)
is obtained 
from the operator $\mathbb{K'}$ (more precisely, from its kernel) by replacing 
$f_2f_3$ with $s_q$.

Note that the equality (\ref{july15017}) remains true
when $s_q$ is a partial sum of the Fourier--Legendre series
of any function from $L_2([t,T]^2),$ i.e. the equality holds
on a dense subset in $L_2([t,T]^2).$

Trace class operators form a linear space. Therefore,
on the left-hand side of 
(\ref{july15017}) there is a matrix trace of a trace class
operator $\mathbb{K'}_q$. 
The mentioned matrix trace
is a linear bounded (and therefore continuous)
functional
in the space of trace class operators \cite{gohb},
\cite{goldberg}
(this functional can be extended to the space $L_2([t, T]^2)$ by continuity \cite{Pugach}).

From the other hand, the right-hand side of (\ref{july15017}) defines
(as a scalar product of $s_q(t_2,t_4)$ and $\psi_4(t_4)\psi_1(t_2)$
in $L_2([t, T]^2)$) a linear bounded (and therefore continuous)
functional in $L_2([t, T]^2),$
which is given by the function $\psi_4(t_4)\psi_1(t_2)$.
On the left-hand side of (\ref{july15017}) (by virtue of the equality (\ref{july15017}))
there is a linear continuous functional on a dense subset in 
$L_2([t,T]^2)$ (see above). This functional can be uniquely extended 
to a linear continuous functional in $L_2([t, T]^2)$
(see \cite{reed}, Theorem~I.7, P.~9).

Let us implement the passage to the limit $\lim\limits_{q\to\infty}$
in the equality (\ref{july15017}) (at that we suppose that $s_q$ is defined by (\ref{july15015}))

\vspace{-1mm}
$$
\sum_{j_1,j_2=0}^{\infty} \left\langle
\Psi_{j_1 j_2}, \mathbb{K''}\Psi_{j_1 j_2}\right\rangle_{L_2([t,T]^2)}=
$$

$$
=4\sum_{j_1,j_2=0}^{\infty}
\int\limits_{[t,T]^4}{\bf 1}_{\{t_1<t_2<t_3<t_4\}}
\psi_4(t_4)\psi_3(t_3)\psi_2(t_2)\psi_1(t_1)
\phi_{j_2}(t_4)
\phi_{j_2}(t_3)
\phi_{j_1}(t_2)
\phi_{j_1}(t_1)dt_1 dt_2 dt_3 dt_4=
$$
\begin{equation}
\label{july15019}
=\int\limits_t^T 
\psi_4(t_4) \psi_3(t_4) \int\limits_t^{t_4} \psi_2(t_2)\psi_1(t_2)dt_2 dt_4,
\end{equation}

\vspace{3mm}
\noindent
where the operator $\mathbb{K''}$ (more precisely, its kernel)
is obtained 
from the operator $\mathbb{K'}_q$ (more precisely, from its kernel) by replacing 
$s_q$ with $\lim\limits_{q\to\infty}s_q=g\in L_2([t,T]^2)$,
$\psi_2(\tau), \psi_3(\tau)\in L_2([t,T])$ and 
$\psi_1(\tau), \psi_4(\tau)$ are continuously differentiable
functions on $[t, T].$

Further, the formula (\ref{july15019}) will remain valid
if we choose

\vspace{1mm}
$$
\psi_1(\tau)=\bar \psi_1^{(p)}(\tau),\ \ \ \psi_4(\tau)=\bar \psi_4^{(p)}(\tau),
$$

\vspace{3mm}
\noindent
where
\begin{equation}
\label{july15018}
\bar \psi_1^{(p)}(\tau)=\sum\limits_{j=0}^p \bar \phi_j(\tau)\int\limits_t^T \bar \psi_1(s) 
\bar \phi_j(s)ds,\ \ \
\bar \psi_4^{(p)}(\tau)=\sum\limits_{j=0}^p \bar \phi_j(\tau)
\int\limits_t^T \bar \psi_4(s) \bar \phi_j(s)ds,
\end{equation}

\vspace{3mm}
\noindent
where $p\in \mathbb{N},$ 
$\bar \psi_1(\tau), \bar \psi_4(\tau)\in L_2([t,T]),$ and $\left\{\bar \phi_j(x)\right\}_{j=0}^{\infty}$ 
is a complete orthonormal system of Legendre polynomials in $L_2([t,T])$.

Substitute (\ref{july15018}) into (\ref{july15019})

\vspace{-1mm}
$$
\sum_{j_1,j_2=0}^{\infty} \left\langle
\Psi_{j_1 j_2}, \mathbb{K''}_p\Psi_{j_1 j_2}\right\rangle_{L_2([t,T]^2)}=
$$

\vspace{1mm}
$$
=4\hspace{-1.3mm}\sum_{j_1,j_2=0}^{\infty}
\int\limits_{[t,T]^4}\hspace{-1.3mm}{\bf 1}_{\{t_1<t_2<t_3<t_4\}}
\bar \psi_4^{(p)}(t_4)\psi_3(t_3)\psi_2(t_2)\bar\psi_1^{(p)}(t_1)
\phi_{j_2}(t_4)
\phi_{j_2}(t_3)
\phi_{j_1}(t_2)
\phi_{j_1}(t_1)dt_1 dt_2 dt_3 dt_4=
$$
\begin{equation}
\label{july15020}
=\int\limits_t^T 
\bar \psi_4^{(p)}(t_4) \psi_3(t_4) \int\limits_t^{t_4} \psi_2(t_2)\bar \psi_1^{(p)}(t_2)dt_2 dt_4.
\end{equation}

\vspace{3mm}
\noindent
where the operator $\mathbb{K''}_p$ (more precisely, its kernel)
is obtained 
from the operator $\mathbb{K''}$ (more precisely, from its kernel) by replacing 
$\psi_4$ and $\psi_1$ with $\bar \psi_4^{(p)}$ and $\bar \psi_1^{(p)},$
respectively.

Note that the equality (\ref{july15020}) will also remain true if
$\bar \psi_4^{(p)}\bar \psi_1^{(p)}$ is replaced by $s_p$
($s_p$ is the partial sum of the Fourier--Legendre series
of any function from $L_2([t, T]^2)$), i.e. 
the modified equality (\ref{july15020}) is true 
on a dense subset of $L_2([t,T]^2).$
Next, we can apply the reasoning below the formula 
(\ref{july15017}) and obtain the equality of two linear continuous functionals
in $L_2([t,T]^2).$
Let us implement the passage to the limit $\lim\limits_{p\to\infty}$
in the mentioned equality under the condition $s_p=\bar \psi_4^{(p)}\bar \psi_1^{(p)}$

\vspace{-1mm}
$$
4\sum_{j_1,j_2=0}^{\infty}
\int\limits_{[t,T]^4}{\bf 1}_{\{t_1<t_2<t_3<t_4\}}
\bar \psi_4(t_4)\psi_3(t_3)\psi_2(t_2)\bar\psi_1(t_1)
\phi_{j_2}(t_4)
\phi_{j_2}(t_3)
\phi_{j_1}(t_2)
\phi_{j_1}(t_1)dt_1 dt_2 dt_3 dt_4=
$$
\begin{equation}
\label{july15021}
=\int\limits_t^T 
\bar \psi_4(t_4) \psi_3(t_4) \int\limits_t^{t_4} \psi_2(t_2)\bar \psi_1(t_2)dt_2 dt_4,
\end{equation}

\vspace{3mm}
\noindent
where $\bar \psi_1(\tau), \psi_2(\tau), \psi_3(\tau), \bar \psi_4(\tau)\in L_2([t,T]).$

Rewrite the equality (\ref{july15021}) in the form

\vspace{-1mm}
$$
\lim\limits_{p\to\infty}\sum_{j_1,j_2=0}^{p}C_{j_2 j_2 j_1 j_1}=
$$

$$
=\sum_{j_1,j_2=0}^{\infty}
\int\limits_t^T
\psi_4(t_4)\phi_{j_2}(t_4) 
\int\limits_t^{t_4}
\psi_3(t_3)\phi_{j_2}(t_3)
\int\limits_t^{t_3}
\psi_2(t_2)\phi_{j_1}(t_2)
\int\limits_t^{t_2}
\psi_1(t_1)\phi_{j_1}(t_1)dt_1 dt_2 dt_3 dt_4=
$$
\begin{equation}
\label{july15022}
=\frac{1}{4}\int\limits_t^T 
\psi_4(t_4) \psi_3(t_4) \int\limits_t^{t_4} \psi_2(t_2)\psi_1(t_2)dt_2 dt_4,
\end{equation}

\vspace{4mm}
\noindent
where $\psi_1(\tau), \ldots, \psi_4(\tau)\in L_2([t,T]).$

Note that the series on the left-hand side of 
(\ref{july15022}) converges absolutly since
its sum does not depend 
on permutations of basis functions
(here the basis in $L_2([t,T]^2)$ is
$\left\{\phi_{j_1}(x)\phi_{j_2}(y)\right\}_{j_1,j_2=0}^{\infty}$).
The equality (\ref{july13}) is proved. 

In \cite{rybakov5000}, the equality (\ref{july15022})
is generalized as follows

\vspace{1mm}
$$
\lim\limits_{p\to\infty}\sum_{j_k,j_{k-2},\ldots, j_2=0}^{p}
C_{j_k j_k j_{k-2} j_{k-2} \ldots j_2 j_2}=
$$

\vspace{1mm}
\begin{equation}
\label{july15023}
=\frac{1}{2^r}\int\limits_t^T 
\psi_k(t_k) \psi_{k-1}(t_k) 
\int\limits_t^{t_k} \psi_{k-2}(t_{k-2})\psi_{k-3}(t_{k-2})
\ldots
\int\limits_t^{t_4} \psi_{2}(t_{2})\psi_{1}(t_{2})dt_2 \ldots 
dt_{k-2} dt_k,
\end{equation}

\vspace{3mm}
\noindent
where $k=2r$ $(r=2,3,\ldots),$
$\psi_1(\tau),\ldots, \psi_k(\tau)\in L_2([t,T]).$

The equalities (\ref{july14}), (\ref{july15}) can also be obtained
\cite{rybakov7000xa}            
using the approach from \cite{rybakov5000} and the series
on the left-hand sides of (\ref{july14}), (\ref{july15})
converge absolutely.

In the notations of Theorem~47, the equality
(\ref{july15023}) can be written in the form

\vspace{1mm}
$$
\lim\limits_{p\to\infty}\sum\limits_{j_{1}=0}^p \sum\limits_{j_{3}=0}^p\ldots 
\sum\limits_{j_{2r-1}=0}^p
C_{j_k\ldots j_1}\biggl|_{j_{1}=j_{2},\ldots, j_{2r-1}=j_{2r}}\biggr.=
$$

\vspace{1mm}
\begin{equation}
\label{july16000}
=\frac{1}{2^r} 
C_{j_k \ldots j_1}\biggl|_{(j_{2} j_{1})\curvearrowright (\cdot)
(j_{4} j_{3})\curvearrowright (\cdot)
\ldots (j_{2r} j_{2r-1})\curvearrowright (\cdot),
j_1=j_2,j_3=j_4,\ldots, j_{2r-1}=j_{2r}
}\biggr.,
\end{equation}

\vspace{3mm}
\noindent
where $k=2r$ $(r=2,3,\ldots)$ and $C_{j_k\ldots j_1}$ is defined by (\ref{july15030}).

In principle, using the method from \cite{rybakov5000}
the following equality can be obtained \cite{rybakov7000xa}

\vspace{1mm}
$$
\lim\limits_{p\to\infty}
\sum\limits_{j_{g_1}, j_{g_3},\ldots,j_{g_{2r-1}}=0}^p
C_{j_k\ldots j_1}\biggl|_{j_{g_1}=j_{g_2},\ldots, j_{g_{2r-1}}=j_{g_{2r}}}=
$$

\vspace{1mm}
$$
=\frac{1}{2^r} \prod\limits_{l=1}^r {\bf 1}_{\{g_{2l}=g_{2l-1}+1\}}
C_{j_k \ldots j_1}\biggl|_{(j_{g_2} j_{g_1})\curvearrowright (\cdot)
\ldots (j_{g_{2r}} j_{g_{2r-1}})\curvearrowright (\cdot),
j_{g_{{}_{1}}}=~j_{g_{{}_{2}}},\ldots, j_{g_{{}_{2r-1}}}=~j_{g_{{}_{2r}}}
}\biggr.
$$

\vspace{5mm}
\noindent
for all possible $g_1,g_2,\ldots,g_{2r-1},g_{2r}$ (see {\rm (\ref{leto5007xxx})),
where $k=2r$ $(r=2,3,\ldots),$ $C_{j_k\ldots j_1}$ is defined by (\ref{july15030}),
another notations are the same as in Theorem~47. We will prove this equality further.

Let us prove the equalities (\ref{july13})--(\ref{july15})
using a method based on generalized Parseval's equality and (\ref{after1400}).

Consider (\ref{july13}). Using (\ref{after1400}), we have

$$
\lim\limits_{p\to\infty}\sum_{j_1,j_2=0}^{p}
\int\limits_t^T
\psi_4(t_4)\phi_{j_2}(t_4)\int\limits_t^{t_4}
\psi_3(t_3)\phi_{j_2}(t_3)\int\limits_t^{T}
\psi_2(t_2)\phi_{j_1}(t_2)
\int\limits_{t}^{t_2}
\psi_1(t_1)\phi_{j_1}(t_1)
dt_1 dt_2 dt_3 dt_4=
$$

\vspace{1mm}
$$
=\lim\limits_{p\to\infty}\sum_{j_1,j_2=0}^{p}
\int\limits_t^T
\psi_4(t_4)\phi_{j_2}(t_4)\int\limits_t^{t_4}
\psi_3(t_3)\phi_{j_2}(t_3)dt_3 dt_4
\int\limits_t^{T}
\psi_2(t_2)\phi_{j_1}(t_2)
\int\limits_{t}^{t_2}
\psi_1(t_1)\phi_{j_1}(t_1)
dt_1 dt_2=
$$

\vspace{1mm}
$$
=\lim\limits_{p\to\infty}\sum_{j_2=0}^{p}
\int\limits_t^T
\psi_4(t_4)\phi_{j_2}(t_4)\int\limits_t^{t_4}
\psi_3(t_3)\phi_{j_2}(t_3)dt_3 dt_4
\lim\limits_{p\to\infty}\sum_{j_1=0}^{p}\int\limits_t^{T}
\psi_2(t_2)\phi_{j_1}(t_2)
\int\limits_{t}^{t_2}
\psi_1(t_1)\phi_{j_1}(t_1)
dt_1 dt_2=
$$

\vspace{1mm}
\begin{equation}
\label{july29000}
=\frac{1}{4}\int\limits_t^T
\psi_4(t_4)\psi_3(t_4)dt_4
\int\limits_t^T
\psi_2(t_2)\psi_1(t_2)dt_2=
\frac{1}{4}\int\limits_{[t,T]^2}
\psi_4(t_4)\psi_3(t_4)
\psi_2(t_2)\psi_1(t_2)dt_2 dt_4,
\end{equation}

\vspace{4mm}
\noindent
where $\psi_1(\tau),\ldots,\psi_4(\tau)\in L_2([t, T]).$ 

Suppose that $\psi_2(\tau)$ and $\psi_3(\tau)$ are polynomials of finite degrees.
For example, $\psi_2(\tau)$ and $\psi_3(\tau)$ can be Legendre polynomials
that form a complete orthonormal system of functions in $L_2([t,T]).$

Denote
\begin{equation}
\label{july30000}
s_q(t_2,t_3)=\sum\limits_{l_1,l_2=0}^q
C_{l_2 l_1}\bar \phi_{l_1}(t_2)\bar \phi_{l_2}(t_3),
\end{equation}

\vspace{3mm}
\noindent 
where $\left\{\bar \phi_j(x)\right\}_{j=0}^{\infty}$ 
is a complete orthonormal system of Legendre polynomials in $L_2([t,T])$
and
$C_{l_2 l_1}$ are Fourier--Legendre coefficients for the function
$g(t_2,t_3)=\bar \psi_2(t_2)\bar \psi_3(t_3){\bf 1}_{\{t_2<t_3\}}$ 
($\bar \psi_2(\tau), \bar \psi_3(\tau)\in L_2([t,T])),$ i.e.
$$
C_{l_2 l_1}=\int\limits_t^T \bar \psi_3(t_3)\bar \phi_{l_2}(t_3)
\int\limits_t^{t_3}\bar \psi_2(t_2)\bar \phi_{l_1}(t_2)dt_2 dt_3.
$$

\vspace{3mm}

Further, we have
$$
\lim\limits_{q\to\infty}\left\Vert s_q - g\right\Vert^2_{L_2([t,T]^2)}=0.
$$

\vspace{3mm}

From (\ref{july29000}) we obtain (the sum on the right-hand side of (\ref{july30000}) is finite)

$$
\sum_{j_1,j_2=0}^{\infty}
\int\limits_{[t,T]^4}{\bf 1}_{\{t_1<t_2\}}{\bf 1}_{\{t_3<t_4\}}
\psi_4(t_4)\phi_{j_2}(t_4)
s_q(t_2,t_3)\phi_{j_2}(t_3)
\phi_{j_1}(t_2)
\psi_1(t_1)\phi_{j_1}(t_1)dt_1 dt_2 dt_3 dt_4=
$$
\begin{equation}
\label{july30001}
=\frac{1}{4}\int\limits_{[t,T]^2} 
\psi_4(t_4)s_q(t_2,t_4)\psi_1(t_2)dt_2 dt_4.
\end{equation}

\vspace{4mm}

Note that the equality (\ref{july30001}) remains true
when $s_q$ is a partial sum of the Fourier--Legendre series
of any function from $L_2([t,T]^2),$ i.e. the equality holds
on a dense subset in $L_2([t,T]^2).$

The right-hand side of the equality (\ref{july30001}) defines
(as a scalar product of $s_q(t_2,t_4)$ and $\frac{1}{4}\psi_4(t_4)\psi_1(t_2)$
in $L_2([t, T]^2)$) a linear bounded (and therefore continuous)
functional in $L_2([t, T]^2),$
which is given by the function $\frac{1}{4}\psi_4(t_4)\psi_1(t_2)$.
On the left-hand side of (\ref{july30001}) (by virtue of the equality (\ref{july30001}))
there is a linear continuous functional on a dense subset in 
$L_2([t,T]^2).$ This functional can be uniquely extended 
to a linear continuous functional in $L_2([t, T]^2)$
(see \cite{reed}, Theorem~I.7, P.~9).

Let us implement the passage to the limit $\lim\limits_{q\to\infty}$
in (\ref{july30001}) (at that we suppose that $s_q$ is defined by (\ref{july30000}))

$$
\sum_{j_1,j_2=0}^{\infty}
\int\limits_{[t,T]^4}{\bf 1}_{\{t_1<t_2<t_3<t_4\}}
\psi_4(t_4)\bar \psi_3(t_3)\bar \psi_2(t_2)\psi_1(t_1)
\phi_{j_2}(t_4)
\phi_{j_2}(t_3)
\phi_{j_1}(t_2)
\phi_{j_1}(t_1)dt_1 dt_2 dt_3 dt_4=
$$
\begin{equation}
\label{july30002}
=\frac{1}{4}\int\limits_t^T 
\psi_4(t_4) \bar\psi_3(t_4) \int\limits_t^{t_4} \bar\psi_2(t_2)\psi_1(t_2)dt_2 dt_4,
\end{equation}

\vspace{4mm}
\noindent
where $\psi_1(\tau), \bar \psi_2(\tau), \bar \psi_3(\tau), \psi_4(\tau)\in L_2([t,T]).$

Rewrite the equality (\ref{july30002}) in the form

$$
\lim\limits_{p\to\infty}\sum_{j_1,j_2=0}^{p}C_{j_2 j_2 j_1 j_1}=
$$

\vspace{1mm}
$$
=\sum_{j_1,j_2=0}^{\infty}
\int\limits_t^T
\psi_4(t_4)\phi_{j_2}(t_4) 
\int\limits_t^{t_4}
\psi_3(t_3)\phi_{j_2}(t_3)
\int\limits_t^{t_3}
\psi_2(t_2)\phi_{j_1}(t_2)
\int\limits_t^{t_2}
\psi_1(t_1)\phi_{j_1}(t_1)dt_1 dt_2 dt_3 dt_4=
$$
\begin{equation}
\label{july30003}
=\frac{1}{4}\int\limits_t^T 
\psi_4(t_4) \psi_3(t_4) \int\limits_t^{t_4} \psi_2(t_2)\psi_1(t_2)dt_2 dt_4,
\end{equation}

\vspace{4mm}
\noindent
where $\psi_1(\tau), \ldots, \psi_4(\tau)\in L_2([t,T]).$

Note that the series on the left-hand side of 
(\ref{july30003}) converges absolutly since
its sum does not depend 
on permutations of basis functions
(here the basis in $L_2([t,T]^2)$ is
$\left\{\phi_{j_1}(x)\phi_{j_2}(y)\right\}_{j_1,j_2=0}^{\infty}$).
The equality (\ref{july13}) is proved.

Let us prove (\ref{july15}). Using the generalized Parseval equality, we obtain

$$
\lim\limits_{p\to\infty}\sum_{j_1,j_2=0}^{p}
\int\limits_t^T
\psi_4(t_4)\phi_{j_2}(t_4)\int\limits_t^{t_4}
\psi_3(t_3)\phi_{j_1}(t_3)\int\limits_t^{T}
\psi_2(t_2)\phi_{j_2}(t_2)
\int\limits_{t}^{t_2}
\psi_1(t_1)\phi_{j_1}(t_1)
dt_1 dt_2 dt_3 dt_4=
$$

\vspace{2mm}
$$
=\sum_{j_1,j_2=0}^{\infty}
\int\limits_t^T
\psi_4(t_4)\phi_{j_2}(t_4)\int\limits_t^{t_4}
\psi_3(t_3)\phi_{j_1}(t_3)dt_3 dt_4
\int\limits_t^{T}
\psi_2(t_2)\phi_{j_2}(t_2)
\int\limits_{t}^{t_2}
\psi_1(t_1)\phi_{j_1}(t_1)
dt_1 dt_2=
$$

\vspace{3mm}
$$
=\sum_{j_1,j_2=0}^{\infty}
\int\limits_{[t,T]^2}
\hspace{-1.2mm}{\bf 1}_{\{t_3<t_4\}}\psi_3(t_3)\psi_4(t_4)\phi_{j_1}(t_3)
\phi_{j_2}(t_4)dt_3 dt_4
\int\limits_{[t,T]^2}
\hspace{-1.2mm}{\bf 1}_{\{t_3<t_4\}}\psi_1(t_3)\psi_2(t_4)\phi_{j_1}(t_3)
\phi_{j_2}(t_4)dt_3 dt_4=
$$

\vspace{2mm}
\begin{equation}
\label{july30004}
=
\int\limits_{[t,T]^2}
{\bf 1}_{\{t_3<t_4\}}\psi_3(t_3)\psi_2(t_4)\psi_4(t_4)\psi_1(t_3)
dt_3 dt_4=
\int\limits_{[t,T]^2}
{\bf 1}_{\{t_3<t_2\}}\psi_3(t_3)\psi_2(t_2)\psi_4(t_2)\psi_1(t_3)
dt_3 dt_2,
\end{equation}

\vspace{4mm}
\noindent
where $\psi_1(\tau), \psi_2(\tau), \psi_3(\tau), \psi_4(\tau)\in L_2([t, T]).$

Suppose that $\psi_2(\tau)$ and $\psi_3(\tau)$ are Legendre polynomials of finite degrees.
Denote

\vspace{-1mm}
\begin{equation}
\label{july30005}
s_q(t_2,t_3)=\sum\limits_{l_1,l_2=0}^q
C_{l_2 l_1}\bar \phi_{l_1}(t_2)\bar \phi_{l_2}(t_3),
\end{equation}

\vspace{3mm}
\noindent 
where $\left\{\bar \phi_j(x)\right\}_{j=0}^{\infty}$ 
is a complete orthonormal system of Legendre polynomials in $L_2([t,T])$
and
$C_{l_2 l_1}$ are Fourier--Legendre coefficients for the function
$g(t_2,t_3)=\bar \psi_2(t_2)\bar \psi_3(t_3){\bf 1}_{\{t_2<t_3\}}$ 
($\bar \psi_2(\tau), \bar \psi_3(\tau)\in L_2([t,T])),$ i.e.
$$
C_{l_2 l_1}=\int\limits_t^T \bar \psi_3(t_3)\bar \phi_{l_2}(t_3)
\int\limits_t^{t_3}\bar \psi_2(t_2)\bar \phi_{l_1}(t_2)dt_2 dt_3.
$$

\vspace{4mm}

Moreover,
$$
\lim\limits_{q\to\infty}\left\Vert s_q - g\right\Vert^2_{L_2([t,T]^2)}=0.
$$

\vspace{5mm}

From (\ref{july30004}) we obtain (the sum on the right-hand side of (\ref{july30005}) is finite)

$$
\sum_{j_1,j_2=0}^{\infty}
\int\limits_{[t,T]^4}{\bf 1}_{\{t_1<t_2\}}{\bf 1}_{\{t_3<t_4\}}
\psi_4(t_4)
s_q(t_2,t_3)\psi_1(t_1)\phi_{j_2}(t_4)\phi_{j_1}(t_3)
\phi_{j_2}(t_2)
\phi_{j_1}(t_1)dt_1 dt_2 dt_3 dt_4=
$$

\begin{equation}
\label{july30006}
=
\int\limits_{[t,T]^2}
{\bf 1}_{\{t_3<t_2\}}s_q(t_2,t_3)\psi_1(t_3)\psi_4(t_2)
dt_3 dt_2.
\end{equation}

\vspace{4mm}

Note that the equality (\ref{july30006}) remains true
when $s_q$ is a partial sum of the Fourier--Legendre series
of any function from $L_2([t,T]^2),$ i.e. the equality holds
on a dense subset in $L_2([t,T]^2).$

The right-hand side of (\ref{july30006}) defines
(as a scalar product of $s_q(t_2,t_3)$ and ${\bf 1}_{\{t_3<t_2\}}\psi_1(t_3)\psi_4(t_2)$
in $L_2([t, T]^2)$) a linear bounded (and therefore continuous)
functional in $L_2([t, T]^2),$
which is given by the function ${\bf 1}_{\{t_3<t_2\}}\psi_1(t_3)\psi_4(t_2)$.
On the left-hand side of (\ref{july30006}) (by virtue of the equality (\ref{july30006}))
there is a linear continuous functional on a dense subset in 
$L_2([t,T]^2).$ This functional can be uniquely extended 
to a linear continuous functional in $L_2([t, T]^2)$
(see \cite{reed}, Theorem~I.7, P.~9).

Let us implement the passage to the limit $\lim\limits_{q\to\infty}$
in (\ref{july30006}) (at that we suppose that $s_q$ is defined by (\ref{july30005}))

\vspace{-1mm}
$$
\sum_{j_1,j_2=0}^{\infty}
\int\limits_{[t,T]^4}{\bf 1}_{\{t_1<t_2<t_3<t_4\}}
\psi_4(t_4)
\bar\psi_3(t_3)\bar \psi_2(t_2)\psi_1(t_1)\phi_{j_2}(t_4)\phi_{j_1}(t_3)
\phi_{j_2}(t_2)
\phi_{j_1}(t_1)dt_1 dt_2 dt_3 dt_4=
$$

\begin{equation}
\label{july30007}
=
\int\limits_{[t,T]^2}
{\bf 1}_{\{t_2>t_3\}}{\bf 1}_{\{t_2<t_3\}} \bar\psi_3(t_3)\bar \psi_2(t_2)\psi_1(t_3)\psi_4(t_2)
dt_3 dt_2=0.
\end{equation}

\vspace{4mm}

Rewrite the equality (\ref{july30007}) in the form

$$
\lim\limits_{p\to\infty}\sum_{j_1,j_2=0}^{p}C_{j_2 j_1 j_2 j_1}=
$$

\vspace{1mm}
\begin{equation}
\label{july30008}
=\sum_{j_1,j_2=0}^{\infty}
\int\limits_t^T
\psi_4(t_4)\phi_{j_2}(t_4) 
\int\limits_t^{t_4}
\psi_3(t_3)\phi_{j_1}(t_3)
\int\limits_t^{t_3}
\psi_2(t_2)\phi_{j_2}(t_2)
\int\limits_t^{t_2}
\psi_1(t_1)\phi_{j_1}(t_1)
dt_1 dt_2 dt_3 dt_4
=0,
\end{equation}

\vspace{4mm}
\noindent
where $\psi_1(\tau), \ldots, \psi_4(\tau)\in L_2([t,T]).$

Note that the series on the left-hand side of 
(\ref{july30008}) converges absolutly since
its sum does not depend 
on permutations of basis functions
(here the basis in $L_2([t,T]^2)$ is
$\left\{\phi_{j_1}(x)\phi_{j_2}(y)\right\}_{j_1,j_2=0}^{\infty}$).
The equality (\ref{july15}) is proved.

Let us prove (\ref{july14}). Using Fubini's Theorem and generalized Parseval's
equality, we get

$$
\lim\limits_{p\to\infty}\sum_{j_1,j_2=0}^{p}
\int\limits_t^T
\psi_4(t_4)\phi_{j_1}(t_4)\int\limits_t^{T}
\psi_3(t_3)\phi_{j_2}(t_3)\int\limits_t^{t_3}
\psi_2(t_2)\phi_{j_2}(t_2)
\int\limits_{t}^{t_2}
\psi_1(t_1)\phi_{j_1}(t_1)
dt_1 dt_2 dt_3 dt_4=
$$

\vspace{2mm}
$$
=
\lim\limits_{p\to\infty}\sum_{j_1,j_2=0}^{p}
C_{j_1}^{\psi_4}C_{j_2 j_2 j_1}^{\psi_3 \psi_2 \psi_1}
=
\frac{1}{2}\lim\limits_{p\to\infty}
\sum\limits_{j_1=0}^{p}C_{j_1}^{\psi_4} 
C_{j_2 j_2 j_1}^{\psi_3 \psi_2 \psi_1}\biggl|_{(j_{2} j_{2})\curvearrowright (\cdot)}-
$$

\vspace{2mm}
$$
-
\lim\limits_{p\to\infty}
\sum\limits_{j_1=0}^{p}C_{j_1}^{\psi_4}
\left(\frac{1}{2}
C_{j_2 j_2 j_1}^{\psi_3 \psi_2 \psi_1}\biggl|_{(j_{2} j_{2})\curvearrowright (\cdot)}-
\sum\limits_{j_2=0}^{p} C_{j_2 j_2 j_1}^{\psi_3 \psi_2 \psi_1}\right)=
$$

\vspace{2mm}
$$
=
\frac{1}{2}\lim\limits_{p\to\infty}
\sum\limits_{j_1=0}^{p}\int\limits_t^T \psi_4(s)\phi_{j_1}(s)ds
\int\limits_t^T \psi_3(\tau)\psi_2(\tau)\int\limits_t^{\tau}
\phi_{j_1}(s)\psi_1(s)ds d\tau-
$$

\vspace{2mm}
$$
-
\lim\limits_{p\to\infty}
\sum\limits_{j_1=0}^{p}C_{j_1}^{\psi_4}
\left(\frac{1}{2}
C_{j_2 j_2 j_1}^{\psi_3 \psi_2 \psi_1}\biggl|_{(j_{2} j_{2})\curvearrowright (\cdot)}-
\sum\limits_{j_2=0}^{p} C_{j_2 j_2 j_1}^{\psi_3 \psi_2 \psi_1}\right)=
$$

\vspace{2mm}
$$
=
\frac{1}{2}\lim\limits_{p\to\infty}
\sum\limits_{j_1=0}^{p}\int\limits_t^T \psi_4(s)\phi_{j_1}(s)ds
\int\limits_t^T \phi_{j_1}(s)\psi_1(s)
\int\limits_s^{T}\psi_3(\tau)\psi_2(\tau)
d\tau ds-
$$

\vspace{2mm}
$$
-
\lim\limits_{p\to\infty}
\sum\limits_{j_1=0}^{p}C_{j_1}^{\psi_4}
\left(\frac{1}{2}
C_{j_2 j_2 j_1}^{\psi_3 \psi_2 \psi_1}\biggl|_{(j_{2} j_{2})\curvearrowright (\cdot)}-
\sum\limits_{j_2=0}^{p} C_{j_2 j_2 j_1}^{\psi_3 \psi_2 \psi_1}\right)=
$$

\vspace{2mm}
$$
=
\frac{1}{2}\int\limits_t^T \psi_4(s) \psi_1(s)
\int\limits_s^{T} \psi_3(\tau) \psi_2(\tau)d\tau ds-
$$

\begin{equation}
\label{july30009}
-
\lim\limits_{p\to\infty}
\sum\limits_{j_1=0}^{p}C_{j_1}^{\psi_4}
\left(\frac{1}{2}
C_{j_2 j_2 j_1}^{\psi_3 \psi_2 \psi_1}\biggl|_{(j_{2} j_{2})\curvearrowright (\cdot)}-
\sum\limits_{j_2=0}^{p} C_{j_2 j_2 j_1}^{\psi_3 \psi_2 \psi_1}\right),
\end{equation}

\vspace{5mm}
\noindent
where $C_{j_1}^{\psi_4}$ and 
$C_{j_2 j_2 j_1}^{\psi_3 \psi_2 \psi_1}$ are defined by (\ref{july40000}).

Due to Cauchy--Bunyakovsky's inequality, Parseval's equality
and (\ref{july1004}), we get 

$$
\lim\limits_{p\to\infty}
\left(\sum\limits_{j_1=0}^{p}C_{j_1}^{\psi_4}
\left(\frac{1}{2}
C_{j_2 j_2 j_1}^{\psi_3 \psi_2 \psi_1}\biggl|_{(j_{2} j_{2})\curvearrowright (\cdot)}-
\sum\limits_{j_2=0}^{p} C_{j_2 j_2 j_1}^{\psi_3 \psi_2 \psi_1}\right)\right)^2\le
$$

$$
\le \lim\limits_{p\to\infty}
\sum\limits_{j_1=0}^{p}\left(C_{j_1}^{\psi_4}\right)^2\
\sum\limits_{j_1=0}^{p}
\left(\frac{1}{2}
C_{j_2 j_2 j_1}^{\psi_3 \psi_2 \psi_1}\biggl|_{(j_{2} j_{2})\curvearrowright (\cdot)}-
\sum\limits_{j_2=0}^{p} C_{j_2 j_2 j_1}^{\psi_3 \psi_2 \psi_1}\right)^2\le
$$

$$
\le \lim\limits_{p\to\infty}
\sum\limits_{j_1=0}^{\infty}\left(C_{j_1}^{\psi_4}\right)^2\
\sum\limits_{j_2=0}^{p}
\left(\frac{1}{2}
C_{j_2 j_2 j_1}^{\psi_3 \psi_2 \psi_1}\biggl|_{(j_{2} j_{2})\curvearrowright (\cdot)}-
\sum\limits_{j_2=0}^{p} C_{j_2 j_2 j_1}^{\psi_3 \psi_2 \psi_1}\right)^2=
$$

\begin{equation}
\label{july30010}
=\int\limits_t^T \psi_4^2(s)ds\lim\limits_{p\to\infty}
\sum\limits_{j_2=0}^{p}
\left(\frac{1}{2}
C_{j_2 j_2 j_1}^{\psi_3 \psi_2 \psi_1}\biggl|_{(j_{2} j_{2})\curvearrowright (\cdot)}-
\sum\limits_{j_2=0}^{p} C_{j_2 j_2 j_1}^{\psi_3 \psi_2 \psi_1}\right)^2=0.
\end{equation}

\vspace{4mm}         

Combining (\ref{july30009}) and (\ref{july30010}), we obtain

$$
\lim\limits_{p\to\infty}\sum_{j_1,j_2=0}^{p}
\int\limits_t^T
\psi_4(t_4)\phi_{j_1}(t_4)\int\limits_t^{T}
\psi_3(t_3)\phi_{j_2}(t_3)\int\limits_t^{t_3}
\psi_2(t_2)\phi_{j_2}(t_2)
\int\limits_{t}^{t_2}
\psi_1(t_1)\phi_{j_1}(t_1)dt_1 dt_2 dt_3 dt_4=
$$

\begin{equation}
\label{july30011}
=\frac{1}{2}\int\limits_t^T \psi_4(s) \psi_1(s)
\int\limits_s^{T} \psi_3(\tau) \psi_2(\tau)d\tau ds=
\frac{1}{2}\int\limits_{[t,T]^2} 
\psi_3(t_3)\psi_4(t_4){\bf 1}_{\{t_4<t_3\}}\psi_1(t_4)\psi_2(t_3)dt_4 dt_3,
\end{equation}

\vspace{4mm}
\noindent
where $\psi_1(\tau), \ldots, \psi_4(\tau)\in L_2([t,T]).$

Suppose that $\psi_3(\tau)$ and $\psi_4(\tau)$ are Legendre polynomials of finite degrees.
Denote

\vspace{-1mm}
\begin{equation}
\label{july30012}
s_q(t_3,t_4)=\sum\limits_{l_1,l_2=0}^q
C_{l_2 l_1}\bar \phi_{l_1}(t_3)\bar \phi_{l_2}(t_4),
\end{equation}

\vspace{3mm}
\noindent 
where $\left\{\bar \phi_j(x)\right\}_{j=0}^{\infty}$ 
is a complete orthonormal system of Legendre polynomials in $L_2([t,T])$
and
$C_{l_2 l_1}$ are Fourier--Legendre coefficients for the function
$g(t_3,t_4)=\bar \psi_3(t_3)\bar \psi_4(t_4){\bf 1}_{\{t_3<t_4\}}$ 
($\bar \psi_3(\tau), \bar \psi_4(\tau)\in L_2([t,T])),$ i.e.
$$
C_{l_2 l_1}=\int\limits_t^T \bar \psi_4(t_4)\bar \phi_{l_2}(t_4)
\int\limits_t^{t_4}\bar \psi_3(t_3)\bar \phi_{l_1}(t_3)dt_3 dt_4.
$$

\vspace{3mm}

Further, we have
$$
\lim\limits_{q\to\infty}\left\Vert s_q - g\right\Vert^2_{L_2([t,T]^2)}=0.
$$

\vspace{3mm}

From (\ref{july30011}) we obtain (the sum on the right-hand side of (\ref{july30012}) is finite)

$$
\sum_{j_1,j_2=0}^{\infty}
\int\limits_{[t,T]^4}{\bf 1}_{\{t_1<t_2<t_3\}}
\phi_{j_1}(t_4)\phi_{j_2}(t_3)
s_q(t_3,t_4)
\psi_2(t_2)\psi_1(t_1)\phi_{j_2}(t_2)
\phi_{j_1}(t_1)dt_1 dt_2 dt_3 dt_4=
$$

\begin{equation}
\label{july30013}
=\frac{1}{2}\int\limits_{[t,T]^2} 
s_q(t_3,t_4){\bf 1}_{\{t_4<t_3\}}\psi_1(t_4)\psi_2(t_3)dt_4 dt_3.
\end{equation}

\vspace{3mm}

Note that the equality (\ref{july30013}) remains true
when $s_q$ is a partial sum of the Fourier--Legendre series
of any function from $L_2([t,T]^2),$ i.e. the equality holds
on a dense subset in $L_2([t,T]^2).$

The right-hand side of (\ref{july30013}) defines
(as a scalar product of $s_q(t_3,t_4)$ and $\frac{1}{2}{\bf 1}_{\{t_4<t_3\}}\psi_1(t_4)\psi_2(t_3)$
in $L_2([t, T]^2)$) a linear bounded (and therefore continuous)
functional in $L_2([t, T]^2),$
which is given by the function $\frac{1}{2}{\bf 1}_{\{t_4<t_3\}}\psi_1(t_4)\psi_2(t_3)$.
On the left-hand side of (\ref{july30013}) (by virtue of the equality (\ref{july30013}))
there is a linear continuous functional on a dense subset in 
$L_2([t,T]^2).$ This functional can be uniquely extended 
to a linear continuous functional in $L_2([t, T]^2)$
(see \cite{reed}, Theorem~I.7, P.~9).

Let us implement the passage to the limit $\lim\limits_{q\to\infty}$
in (\ref{july30013}) (at that we suppose that $s_q$ is defined by (\ref{july30012}))

\vspace{-2mm}
$$
\sum_{j_1,j_2=0}^{\infty}
\int\limits_{[t,T]^4}{\bf 1}_{\{t_1<t_2<t_3<t_4\}}
\bar \psi_4(t_4)\phi_{j_1}(t_4)\bar \psi_3(t_3)\phi_{j_2}(t_3)
\psi_2(t_2)\phi_{j_2}(t_2)
\psi_{1}(t_1)\phi_{j_1}(t_1)dt_1 dt_2 dt_3 dt_4=
$$

\begin{equation}
\label{july30014}
=
\frac{1}{2}\int\limits_{[t,T]^2} 
\bar \psi_3(t_3)\bar \psi_4(t_4)
{\bf 1}_{\{t_3<t_4\}}{\bf 1}_{\{t_4<t_3\}}\psi_1(t_4)\psi_2(t_3)dt_4 dt_3=0.
\end{equation}

\vspace{4mm}

Rewrite the equality (\ref{july30014}) in the form

$$
\lim\limits_{p\to\infty}\sum_{j_1,j_2=0}^{p}C_{j_1 j_2 j_2 j_1}=
$$

\begin{equation}
\label{july30015}
=\sum_{j_1,j_2=0}^{\infty}
\int\limits_t^T
\psi_4(t_4)\phi_{j_1}(t_4) 
\int\limits_t^{t_4}
\psi_3(t_3)\phi_{j_2}(t_3)
\int\limits_t^{t_3}
\psi_2(t_2)\phi_{j_2}(t_2)
\int\limits_t^{t_2}
\psi_1(t_1)\phi_{j_1}(t_1)
dt_1 dt_2 dt_3 dt_4
=0,
\end{equation}

\vspace{4mm}
\noindent
where $\psi_1(\tau), \ldots, \psi_4(\tau)\in L_2([t,T]).$

Note that the series on the left-hand side of 
(\ref{july30015}) converges absolutly since
its sum does not depend 
on permutations of basis functions
(here the basis in $L_2([t,T]^2)$ is
$\left\{\phi_{j_1}(x)\phi_{j_2}(y)\right\}_{j_1,j_2=0}^{\infty}$).
The equality (\ref{july14}) is proved. The equalities
(\ref{july13})--(\ref{july15}) are proved.

By induction we prove the following equality
(i.e. by a different method com\-pa\-red 
with \cite{rybakov5000})

$$
\lim\limits_{p\to\infty}\sum_{j_{2r}, j_{2r-2}, \ldots, j_2=0}^{p}
C_{j_{2r}j_{2r} j_{2r-2}j_{2r-2} \ldots j_2 j_2}=
$$

\begin{equation}
\label{july30016}
=
\frac{1}{2^{r}}
\int\limits_t^T
\psi_{2r}(t_{2r})
\psi_{2r-1}(t_{2r})
\int\limits_t^{t_{2r}}
\psi_{2r-2}(t_{2r-2})
\psi_{2r-3}(t_{2r-2})\ldots
\int\limits_t^{t_{4}}
\psi_{2}(t_{2})
\psi_{1}(t_{2})dt_2\ldots dt_{2r-2}dt_{2r},
\end{equation}

\vspace{4mm}
\noindent
where $r\in\mathbb{N},$ $C_{j_{2r}j_{2r} j_{2r-2}j_{2r-2} \ldots j_2 j_2}$
is defined by 

$$
C_{j_k \ldots j_1}=\int\limits_t^T\psi_k(t_k)\phi_{j_k}(t_k)\ldots
\int\limits_t^{t_2}
\psi_1(t_1)\phi_{j_1}(t_1)
dt_1\ldots dt_k\ \ \ (k\in\mathbb{N}),
$$

\vspace{3mm}
\noindent
$\left\{\phi_j(x)\right\}_{j=0}^{\infty}$
is an arbitrary complete orthonormal system of 
functions in the space $L_2([t,T]),$ and
$\psi_1(\tau),\ldots ,\psi_{2r}(\tau)\in $ $L_2([t, T]).$

Note that the equality (\ref{july13}) is a particular case of 
(\ref{july30016}) for $r=2$ and the equality (\ref{after1400}) is a particular case of 
(\ref{july30016}) for $r=1$.
Thus, the equality
(\ref{july30016}) is true for $r=1, 2.$
Suppose that the equality (\ref{july30016}) is true for some 
$r>2.$ Then, using (\ref{after1400}), we get 

$$
\lim\limits_{p\to\infty}\sum_{j_{2r+2}, j_{2r}, \ldots, j_2=0}^{p}
\int\limits_t^T
\psi_{2r+2}(t_{2r+2})
\phi_{j_{2r+2}}(t_{2r+2})
\int\limits_t^{t_{2r+2}}
\psi_{2r+1}(t_{2r+1})
\phi_{j_{2r+2}}(t_{2r+1})\times
$$

\vspace{1mm}
$$
\times
\int\limits_t^T
\psi_{2r}(t_{2r})
\phi_{j_{2r}}(t_{2r})
\int\limits_t^{t_{2r}}
\psi_{2r-1}(t_{2r-1})
\phi_{j_{2r}}(t_{2r-1})\ldots
$$

\vspace{1mm}
$$
\ldots 
\int\limits_t^{t_3}
\psi_{2}(t_{2})
\phi_{j_{2}}(t_{2})
\int\limits_t^{t_{2}}
\psi_{1}(t_{1})
\phi_{j_{2}}(t_{1})
dt_1 dt_2\ldots dt_{2r-1}dt_{2r}dt_{2r+1}dt_{2r+2}=
$$

\vspace{1mm}
$$
=
\sum_{j_{2r+2}=0}^{\infty}
\int\limits_t^T
\psi_{2r+2}(t_{2r+2})
\phi_{j_{2r+2}}(t_{2r+2})
\int\limits_t^{t_{2r+2}}
\psi_{2r+1}(t_{2r+1})
\phi_{j_{2r+2}}(t_{2r+1})dt_{2r+1}dt_{2r+2}\times
$$

\vspace{1mm}
$$
\times\sum_{j_{2r}, j_{2r-2}, \ldots, j_2=0}^{\infty}
\int\limits_t^T
\psi_{2r}(t_{2r})
\phi_{j_{2r}}(t_{2r})
\int\limits_t^{t_{2r}}
\psi_{2r-1}(t_{2r-1})
\phi_{j_{2r}}(t_{2r-1})\times
$$

\vspace{1mm}
$$
\times
\int\limits_t^{t_{2r-1}}
\psi_{2r-2}(t_{2r-2})
\phi_{j_{2r-2}}(t_{2r-2})
\int\limits_t^{t_{2r-2}}
\psi_{2r-3}(t_{2r-3})
\phi_{j_{2r-2}}(t_{2r-3})\ldots
$$

\vspace{1mm}
$$
\ldots 
\int\limits_t^{t_3}
\psi_{2}(t_{2})
\phi_{j_{2}}(t_{2})
\int\limits_t^{t_{2}}
\psi_{1}(t_{1})
\phi_{j_{2}}(t_{1})
dt_1 dt_2\ldots dt_{2r-3}dt_{2r-2}dt_{2r-1}dt_{2r}=
$$

\vspace{1mm}
$$
=\frac{1}{2}
\int\limits_t^T
\psi_{2r+2}(t_{2r+2})
\psi_{2r+1}(t_{2r+2})
dt_{2r+2}\times
$$

\vspace{1mm}
\begin{equation}
\label{july30017}
\times
\frac{1}{2^{r}}
\int\limits_t^T
\psi_{2r}(t_{2r})
\psi_{2r-1}(t_{2r})
\int\limits_t^{t_{2r}}
\psi_{2r-2}(t_{2r-2})
\psi_{2r-3}(t_{2r-2})\ldots
\int\limits_t^{t_{4}}
\psi_{2}(t_{2})
\psi_{1}(t_{2})dt_2\ldots dt_{2r-2}dt_{2r}.
\end{equation}

\vspace{4mm}

Let us rewrite the equality (\ref{july30017}) in the form

$$
\lim\limits_{p\to\infty}\sum_{j_{2r+2}, j_{2r}, \ldots, j_2=0}^{p}
\int\limits_t^T
\psi_{2r+2}(t_{2r+2})
\phi_{j_{2r+2}}(t_{2r+2})
\int\limits_t^{t_{2r+2}}
\psi_{2r+1}(t_{2r+1})
\phi_{j_{2r+2}}(t_{2r+1})\times
$$

\vspace{1mm}
$$
\times
\int\limits_t^T
\psi_{2r}(t_{2r})
\phi_{j_{2r}}(t_{2r})
\int\limits_t^{t_{2r}}
\psi_{2r-1}(t_{2r-1})
\phi_{j_{2r}}(t_{2r-1})\ldots
$$

\vspace{1mm}
$$
\ldots 
\int\limits_t^{t_3}
\psi_{2}(t_{2})
\phi_{j_{2}}(t_{2})
\int\limits_t^{t_{2}}
\psi_{1}(t_{1})
\phi_{j_{2}}(t_{1})
dt_1 dt_2\ldots dt_{2r-1}dt_{2r}dt_{2r+1}dt_{2r+2}=
$$

\vspace{1mm}
$$
=\frac{1}{2^{r+1}}
\int\limits_t^T
\psi_{2r+2}(t_{2r+2})
\psi_{2r+1}(t_{2r+2})
\int\limits_t^T
\psi_{2r}(t_{2r})
\psi_{2r-1}(t_{2r})\times
$$

\vspace{1mm}
\begin{equation}
\label{july30018}
\times\int\limits_t^{t_{2r}}
\psi_{2r-2}(t_{2r-2})
\psi_{2r-3}(t_{2r-2})\ldots
\int\limits_t^{t_{4}}
\psi_{2}(t_{2})
\psi_{1}(t_{2})dt_2\ldots dt_{2r-2}dt_{2r}dt_{2r+2},
\end{equation}

\vspace{4mm}
\noindent
where $\psi_1(\tau),\ldots, \psi_{2r+2}(\tau)\in L_2([t, T]).$

Suppose that $\psi_{1}(\tau),\psi_{3}(\tau),\ldots, \psi_{2r-3}(\tau),
\psi_{2r}(\tau), \psi_{2r+1}(\tau)$ in (\ref{july30018}) are Legendre polynomials
of finite degrees.
Denote

\vspace{-2mm}
$$
h(t_2,t_4,\ldots,t_{2r-2},t_{2r-1},t_{2r+2})=
\psi_2(t_2)\psi_4(t_4)\ldots
\psi_{2r-2}(t_{2r-2})\psi_{2r-1}(t_{2r-1})\psi_{2r+2}(t_{2r+2}),
$$

\vspace{1mm}
\begin{equation}
\label{july50000}
g(t_1,t_3,\ldots,t_{2r-3},t_{2r},t_{2r+1})=
\bar \psi_{1}(t_1)\bar \psi_{3}(t_3)\ldots \bar\psi_{2r-3}(t_{2r-3})
\bar \psi_{2r}(t_{2r})\bar \psi_{2r+1}(t_{2r+1})
{\bf 1}_{\{t_{2r}<t_{2r+1}\}},
\end{equation}

\vspace{2mm}
$$
s_q(t_1,t_3,\ldots,t_{2r-3},t_{2r},t_{2r+1})=
$$
\begin{equation}
\label{july30019}
=
\sum\limits_{l_1, \ldots, l_{r+1}=0}^q
C_{l_{r+1}\ldots l_1}\bar \phi_{l_1}(t_{1})\bar \phi_{l_2}(t_{3})\ldots \bar \phi_{l_{r-1}}(t_{2r-3})
\bar \phi_{l_r}(t_{2r})\bar \phi_{l_{r+1}}(t_{2r+1}),
\end{equation}

\vspace{2.5mm}
\noindent
where $\left\{\bar \phi_j(x)\right\}_{j=0}^{\infty}$ 
is a complete orthonormal system of Legendre polynomials in $L_2([t,T]),$
$C_{l_{r+1}\ldots l_1}$ are Fourier--Legendre coefficients for the function
(\ref{july50000}), 
$\bar \psi_{1}(\tau),\bar \psi_{3}(\tau),\ldots ,\bar\psi_{2r-3}(\tau),
\bar \psi_{2r}(\tau),\bar \psi_{2r+1}(\tau)$ $\in L_2([t,T]).$ 
Then we have

\vspace{-2mm}
$$
\lim\limits_{q\to\infty}\left\Vert s_q - g\right\Vert^2_{L_2([t,T]^{r+1})}=0.
$$

\vspace{4mm}

From (\ref{july30018}) we obtain (the sum on the right-hand side of (\ref{july30019}) is finite)

$$
\lim\limits_{p\to\infty}\sum_{j_{2r+2}, j_{2r}, \ldots, j_2=0}^{p}
~\int\limits_{[t,T]^{2r+2}}
{\bf 1}_{\{t_1<t_2<\ldots <t_{2r}\}}{\bf 1}_{\{t_{2r+1}<t_{2r+2}\}}
s_q(t_1,t_3,\ldots,t_{2r-3},t_{2r},t_{2r+1})\times
$$

\vspace{2mm}
$$
\times
h(t_2,t_4,\ldots,t_{2r-2},t_{2r-1},t_{2r+2})
\times
$$

\vspace{2mm}
$$
\times
\prod\limits_{d=1}^{r+1}\phi_{j_{2d}}(t_{2d-1})\phi_{j_{2d}}(t_{2d})
dt_1 dt_2\ldots dt_{2r-1}dt_{2r}dt_{2r+1}dt_{2r+2}=
$$

\vspace{4mm}
$$
=\frac{1}{2^{r+1}}
\int\limits_{[t,T]^{r+1}}
{\bf 1}_{\{t_2<t_4<\ldots <t_{2r}\}}s_q(t_2,t_4,\ldots,t_{2r-2},t_{2r},t_{2r+2})\times
$$

\vspace{2mm}
\begin{equation}
\label{july30020}
\times
h(t_2,t_4,\ldots,t_{2r-2},t_{2r},t_{2r+2})dt_2 dt_4\ldots dt_{2r-2}dt_{2r}dt_{2r+2}.
\end{equation}

\vspace{7mm}

The right-hand side of the equality (\ref{july30020}) defines
(as a scalar product of

$$
s_q(t_2,t_4,\ldots,t_{2r-2},t_{2r},t_{2r+2})
$$

\vspace{1mm}
\noindent
and 
$$
\frac{1}{2^{r+1}}
{\bf 1}_{\{t_2<t_4<\ldots <t_{2r}\}}
h(t_2,t_4,\ldots,t_{2r-2},t_{2r},t_{2r+2})
$$

\vspace{4mm}
\noindent
in the space $L_2([t, T]^{r+1})$) a linear bounded (and therefore continuous)
functional in the space $L_2([t, T]^{r+1}).$
The mentioned functional is given by the function 

$$
\frac{1}{2^{r+1}}{\bf 1}_{\{t_2<t_4<\ldots <t_{2r}\}}h(t_2,t_4,\ldots,t_{2r-2},t_{2r},t_{2r+2}).
$$

\vspace{4mm}

Note that the equality (\ref{july30020}) will also remain true if
$s_q$ is replaced 
by $\bar s_q$ ($\bar s_q$ is the partial sum of the Fourier--Legendre series
of any function from $L_2([t, T]^{r+1})$), i.e. 
the modified equality (\ref{july30020}) is true 
on a dense subset in $L_2([t,T]^{r+1}).$
On the left-hand side of (\ref{july30020}) (by virtue of the equality (\ref{july30020}))
there is a linear continuous functional on a dense subset in 
$L_2([t,T]^{r+1}).$ This functional can be uniquely extended 
to a linear continuous functional in $L_2([t, T]^{r+1})$
(see \cite{reed}, Theorem~I.7, P.~9).
Thus, we have the equality of two 
linear continuous functionals in $L_2([t, T]^{r+1}).$
Let us implement the passage to the limit $\lim\limits_{q\to\infty}$
in the mentioned equality if instead of $\bar s_q$
we choose $s_q$ of the form (\ref{july30019}) (i.e. passage to the limit $\lim\limits_{q\to\infty}$
in (\ref{july30020}))

$$
\lim\limits_{p\to\infty}\sum_{j_{2r+2}, j_{2r}, \ldots, j_2=0}^{p}
~\int\limits_{[t,T]^{2r+2}}
{\bf 1}_{\{t_1<t_2<\ldots <t_{2r}\}}{\bf 1}_{\{t_{2r+1}<t_{2r+2}\}}
g(t_1,t_3,\ldots,t_{2r-3},t_{2r},t_{2r+1})\times
$$

\vspace{2mm}
$$
\times
h(t_2,t_4,\ldots,t_{2r-2},t_{2r-1},t_{2r+2})
\times
$$

\vspace{1mm}
$$
\times
\prod\limits_{d=1}^{r+1}\phi_{j_{2d}}(t_{2d-1})\phi_{j_{2d}}(t_{2d})
dt_1 dt_2\ldots dt_{2r-1}dt_{2r}dt_{2r+1}dt_{2r+2}=
$$

\vspace{4mm}
$$
=\frac{1}{2^{r+1}}
\int\limits_{[t,T]^{r+1}}
{\bf 1}_{\{t_2<t_4<\ldots <t_{2r}\}}g(t_2,t_4,\ldots,t_{2r-2},t_{2r},t_{2r+2})\times
$$

\vspace{1mm}
\begin{equation}
\label{july30022}
\times
h(t_2,t_4,\ldots,t_{2r-2},t_{2r},t_{2r+2})dt_2 dt_4\ldots dt_{2r-2}dt_{2r}dt_{2r+2},
\end{equation}

\vspace{5mm}
\noindent
where $\bar \psi_{1}(\tau),\bar \psi_{3}(\tau),\ldots ,\bar\psi_{2r-3}(\tau)
\bar \psi_{2r}(\tau), \bar \psi_{2r+1}(\tau)\in L_2([t,T]).$

It is easy to see that the equality (\ref{july30022}) (up to notations)
is the equality (\ref{july30016}) in which $r$ is replaced by $r+1.$
So, we proved the equality (\ref{july30016}) by induction.

Note that the series on the left-hand side of 
(\ref{july30016}) converges absolutly since
its sum does not depend 
on permutations of basis functions
(here the basis in $L_2([t,T]^{r})$ is
$\left\{\phi_{j_1}(x_1)\ldots \phi_{j_r}(x_r)\right\}_{j_1,\ldots,j_r=0}^{\infty}$).

Further, let us show that

\vspace{1mm}
$$
\lim\limits_{p\to\infty}
\sum\limits_{j_{g_1}, j_{g_3},\ldots,j_{g_{2r-1}}=0}^p
C_{j_k\ldots j_1}\biggl|_{j_{g_1}=j_{g_2},\ldots, j_{g_{2r-1}}=j_{g_{2r}}}=
$$

\vspace{2mm}
\begin{equation}
\label{july90000}
=\frac{1}{2^r} \prod\limits_{l=1}^r {\bf 1}_{\{g_{2l}=g_{2l-1}+1\}}
C_{j_k \ldots j_1}\biggl|_{(j_{g_2} j_{g_1})\curvearrowright (\cdot)
\ldots (j_{g_{2r}} j_{g_{2r-1}})\curvearrowright (\cdot),
j_{g_{{}_{1}}}=~j_{g_{{}_{2}}},\ldots, j_{g_{{}_{2r-1}}}=~j_{g_{{}_{2r}}}
}\biggr.
\end{equation}

\vspace{5mm}
\noindent
for all possible $g_1,g_2,\ldots,g_{2r-1},g_{2r}$ (see {\rm (\ref{leto5007xxx})),
where $k=2r$ $(r=2,3,\ldots),$ $C_{j_k\ldots j_1}$ is defined by (\ref{july15030}),
another notations are the same as in Theorem~47.

The case
$$
\prod\limits_{l=1}^r {\bf 1}_{\{g_{2l}=g_{2l-1}+1\}}=1
$$

\vspace{2mm}
\noindent
corresponds to (\ref{july30016}).

Thus, it remains to prove that

\vspace{-1mm}
\begin{equation}
\label{july80000}
\lim\limits_{p\to\infty}
\sum\limits_{j_{g_1}, j_{g_3},\ldots,j_{g_{2r-1}}=0}^p
C_{j_k\ldots j_1}\biggl|_{j_{g_1}=j_{g_2},\ldots, j_{g_{2r-1}}=j_{g_{2r}}}=0
\end{equation}

\vspace{4mm}
\noindent
for the case
$$
\prod\limits_{l=1}^r {\bf 1}_{\{g_{2l}=g_{2l-1}+1\}}=0.
$$

\vspace{3mm}

Below we consider two examples that clearly explain 
the algorithm for the proof of equality (\ref{july80000}).
After this we will formulate the algorithm.

First, let us prove that

$$
\lim\limits_{p\to\infty}\sum_{j_1,j_3,j_4=0}^{p}
C_{j_3 j_4 j_4 j_3 j_1 j_1}=
$$

\vspace{1mm}

$$
=\lim\limits_{p\to\infty}\sum_{j_1,j_3,j_4=0}^{p}
\int\limits_t^T
\psi_6(t_6)\phi_{j_3}(t_6)\int\limits_t^{t_6}
\psi_5(t_5)\phi_{j_4}(t_5)\int\limits_t^{t_5}
\psi_4(t_4)\phi_{j_4}(t_4)
\int\limits_{t}^{t_4}
\psi_3(t_3)\phi_{j_3}(t_3)\times\Biggr.
$$

\vspace{1mm}
\begin{equation}
\label{july80013}
\times
\int\limits_t^{t_3}\psi_2(t_2)\phi_{j_1}(t_2)\int\limits_t^{t_2}
\psi_1(t_1)\phi_{j_1}(t_1)dt_1 dt_2 dt_3 dt_4 dt_5 dt_6=0,
\end{equation}

\vspace{3mm}
\noindent
where $\left\{\phi_j(x)\right\}_{j=0}^{\infty}$
is an arbitrary complete orthonormal system of 
functions in the space $L_2([t,T])$ and
$\psi_1(\tau),\ldots ,\psi_{6}(\tau)\in L_2([t, T]).$
         
\vspace{2mm}

{\bf Step~1.}\ Using (\ref{july30016}) ($r=2$) and generalized Parseval's equality, we obtain

$$
\lim\limits_{p\to\infty}\sum_{j_1,j_3,j_4=0}^{p}
\int\limits_t^T
\psi_6(t_6)\phi_{j_3}(t_6)\int\limits_t^{T}
\psi_5(t_5)\phi_{j_4}(t_5)\int\limits_t^{t_5}
\psi_4(t_4)\phi_{j_4}(t_4)
\int\limits_{t}^{T}
\psi_3(t_3)\phi_{j_3}(t_3)\times\Biggr.
$$

\vspace{1mm}
\begin{equation}
\label{july80100}
\times
\int\limits_t^{T}\psi_2(t_2)\phi_{j_1}(t_2)\int\limits_t^{t_2}
\psi_1(t_1)\phi_{j_1}(t_1)dt_1 dt_2 dt_3 dt_4 dt_5 dt_6=
\end{equation}

\vspace{1mm}
$$
=\lim\limits_{p\to\infty}\sum_{j_3=0}^{p}
\int\limits_t^T
\psi_6(t_6)\phi_{j_3}(t_6)dt_6 \int\limits_{t}^{T}
\psi_3(t_3)\phi_{j_3}(t_3)dt_3\times
$$

\vspace{1mm}
$$
\times
\lim\limits_{p\to\infty}\sum_{j_4=0}^{p}
\int\limits_t^{T}
\psi_5(t_5)\phi_{j_4}(t_5)\int\limits_t^{t_5}
\psi_4(t_4)\phi_{j_4}(t_4)dt_4 dt_5\times
$$

\vspace{1mm}
$$
\times
\lim\limits_{p\to\infty}\sum_{j_1=0}^{p}
\int\limits_t^{T}\psi_2(t_2)\phi_{j_1}(t_2)\int\limits_t^{t_2}
\psi_1(t_1)\phi_{j_1}(t_1)dt_1 dt_2=
$$

\vspace{1mm}
\begin{equation}
\label{july80001}
=
\int\limits_t^T
\psi_6(t_6)\psi_3(t_6)dt_6
\cdot \frac{1}{2}
\int\limits_t^{T}
\psi_5(t_4)\psi_4(t_4)dt_4
\cdot \frac{1}{2}
\int\limits_t^{T}\psi_2(t_2)\psi_1(t_2)dt_2.
\end{equation}

\vspace{3mm}

Let us rewrite (\ref{july80001}) in the form

$$
\sum_{j_1,j_3,j_4=0}^{\infty}
\int\limits_t^T
\psi_6(t_6)\phi_{j_3}(t_6)\int\limits_t^{T}
\psi_5(t_5)\phi_{j_4}(t_5)\int\limits_t^{t_5}
\psi_4(t_4)\phi_{j_4}(t_4)
\int\limits_{t}^{T}
\psi_3(t_3)\phi_{j_3}(t_3)\times\Biggr.
$$

\vspace{1mm}
$$
\times
\int\limits_t^{T}\psi_2(t_2)\phi_{j_1}(t_2)\int\limits_t^{t_2}
\psi_1(t_1)\phi_{j_1}(t_1)dt_1 dt_2 dt_3 dt_4 dt_5 dt_6=
$$

\vspace{1mm}
\begin{equation}
\label{july80002}
=
\frac{1}{4}\int\limits_t^T
\psi_6(t_6)\psi_3(t_6)
\int\limits_t^{T}
\psi_5(t_4)\psi_4(t_4)
\int\limits_t^{T}\psi_2(t_2)\psi_1(t_2)dt_2 dt_4 dt_6.
\end{equation}

\vspace{3mm}

{\bf Step~2.}\ Suppose that $\psi_2(\tau),\psi_3(\tau),\psi_4(\tau)$ are Legendre
polynomials of finite degrees.
Denote

\vspace{-1mm}
\begin{equation}
\label{july80003}
s_q(t_2,t_3,t_4)=\sum\limits_{l_1,l_2.l_3=0}^q
C_{l_3 l_2 l_1}\bar \phi_{l_1}(t_2)\bar \phi_{l_2}(t_3) \bar \phi_{l_3}(t_4),
\end{equation}

\vspace{3mm}
\noindent 
where $\left\{\bar \phi_j(x)\right\}_{j=0}^{\infty}$ 
is a complete orthonormal system of Legendre polynomials in $L_2([t,T])$
and
$C_{l_3 l_2 l_1}$ are Fourier--Legendre coefficients for the function
$g(t_2,t_3,t_4)=\bar \psi_2(t_2)\bar \psi_3(t_3)\bar \psi_4(t_4){\bf 1}_{\{t_2<t_3\}}$ 
($\bar \psi_2(\tau), \bar \psi_3(\tau), \bar \psi_4(\tau)\in L_2([t,T])),$ i.e.
$\lim\limits_{q\to\infty}\left\Vert s_q - g\right\Vert^2_{L_2([t,T]^3)}=0.$

From (\ref{july80002}) we obtain (the sum on the right-hand side of (\ref{july80003}) is finite)

$$
\sum_{j_1,j_3,j_4=0}^{\infty}
\int\limits_{[t,T]^6}{\bf 1}_{\{t_1<t_2\}}{\bf 1}_{\{t_4<t_5\}}
s_q(t_2,t_3,t_4)\psi_6(t_6)\psi_5(t_5)\psi_1(t_1)
\phi_{j_3}(t_6)
\phi_{j_3}(t_3)
\phi_{j_4}(t_5)
\times\Biggr.
$$

\vspace{1mm}
$$
\times
\phi_{j_4}(t_4)\phi_{j_1}(t_2)
\phi_{j_1}(t_1)dt_1 dt_2 dt_3 dt_4 dt_5 dt_6=
$$

\vspace{1mm}
\begin{equation}
\label{july80004}
=
\frac{1}{4}\int\limits_{[t,T]^3}
s_q(t_2,t_6,t_4)\psi_6(t_6)
\psi_5(t_4)\psi_1(t_2)dt_2 dt_4 dt_6.
\end{equation}

\vspace{3mm}

Note that the equality (\ref{july80004}) remains true
when $s_q$ is a partial sum of the Fourier--Legendre series
of any function from $L_2([t,T]^3),$ i.e. the equality holds
on a dense subset in $L_2([t,T]^3).$

The right-hand side of (\ref{july80004}) defines
(as a scalar product of $s_q(t_2,t_6,t_4)$ and $\frac{1}{4}\psi_6(t_6)\psi_5(t_4)\psi_1(t_2)$
in $L_2([t, T]^3)$) a linear bounded (and therefore continuous)
functional in $L_2([t, T]^3),$
which is given by the function $\frac{1}{4}\psi_6(t_6)\psi_5(t_4)\psi_1(t_2)$.
On the left-hand side of (\ref{july80004}) (by virtue of the equality (\ref{july80004}))
there is a linear continuous functional on a dense subset in 
$L_2([t,T]^3).$ This functional can be uniquely extended 
to a linear continuous functional in $L_2([t, T]^3)$
(see \cite{reed}, Theorem~I.7, P.~9).

Let us implement the passage to the limit $\lim\limits_{q\to\infty}$
in (\ref{july80004}) (at that we suppose that $s_q$ is defined by (\ref{july80003}))

$$
\sum_{j_1,j_3,j_4=0}^{\infty}
\int\limits_{[t,T]^6}{\bf 1}_{\{t_1<t_2<t_3\}}{\bf 1}_{\{t_4<t_5\}}
\psi_6(t_6)\psi_5(t_5)\bar \psi_4(t_4)\bar \psi_3(t_3)\bar \psi_2(t_2) \psi_1(t_1)
\phi_{j_3}(t_6)
\phi_{j_3}(t_3)
\times\Biggr.
$$

\vspace{1mm}
$$
\times
\phi_{j_4}(t_5)\phi_{j_4}(t_4)\phi_{j_1}(t_2)
\phi_{j_1}(t_1)dt_1 dt_2 dt_3 dt_4 dt_5 dt_6=
$$

\vspace{1mm}
\begin{equation}
\label{july80005}
=
\frac{1}{4}\int\limits_{[t,T]^3}
{\bf 1}_{\{t_2<t_6\}}
\psi_6(t_6)\bar \psi_3(t_6)
\psi_5(t_4)\bar \psi_4(t_4)
\bar \psi_2(t_2)
\psi_1(t_2)dt_2 dt_4 dt_6.
\end{equation}

\vspace{3mm}

Rewrite the equality (\ref{july80005}) in the form

$$
\sum_{j_1,j_3,j_4=0}^{\infty}
\int\limits_{[t,T]^6}{\bf 1}_{\{t_1<t_2<t_3\}}{\bf 1}_{\{t_4<t_5\}}
\psi_6(t_6)\psi_5(t_5)\psi_4(t_4)\psi_3(t_3)\psi_2(t_2) \psi_1(t_1)
\phi_{j_3}(t_6)
\phi_{j_3}(t_3)
\times\Biggr.
$$

\vspace{1mm}
$$
\times
\phi_{j_4}(t_5)\phi_{j_4}(t_4)\phi_{j_1}(t_2)
\phi_{j_1}(t_1)dt_1 dt_2 dt_3 dt_4 dt_5 dt_6=
$$

\vspace{1mm}
\begin{equation}
\label{july80006}
=
\frac{1}{4}\int\limits_{[t,T]^3}
{\bf 1}_{\{t_2<t_6\}}
\psi_6(t_6)\psi_3(t_6)
\psi_5(t_4)\psi_4(t_4)
\psi_2(t_2)
\psi_1(t_2)dt_2 dt_4 dt_6,
\end{equation}

\vspace{3mm}
\noindent
where $\psi_1(\tau),\ldots,\psi_6(\tau)\in L_2([t,T]).$

\vspace{2mm}

{\bf Step~3.}\ Suppose that $\psi_3(\tau),\psi_4(\tau),\psi_1(\tau)$ are 
Legendre polynomials of finite degrees.
Denote

\vspace{-1mm}
\begin{equation}
\label{july80007}
s_q(t_3,t_4,t_1)=\sum\limits_{l_1,l_2.l_3=0}^q
C_{l_3 l_2 l_1}\bar \phi_{l_1}(t_3)\bar \phi_{l_2}(t_4) \bar \phi_{l_3}(t_1),
\end{equation}

\vspace{3mm}
\noindent 
where $\left\{\bar \phi_j(x)\right\}_{j=0}^{\infty}$ as in (\ref{july80003})
and
$C_{l_3 l_2 l_1}$ are Fourier--Legendre coefficients for the function
$g(t_3,t_4,t_1)=\bar \psi_3(t_3)\bar \psi_4(t_4)\bar \psi_1(t_1){\bf 1}_{\{t_3<t_4\}}$ 
($\bar \psi_3(\tau), \bar \psi_4(\tau), \bar \psi_1(\tau)\in L_2([t,T])),$ i.e.
$\lim\limits_{q\to\infty}\left\Vert s_q - g\right\Vert^2_{L_2([t,T]^3)}=0.$

From (\ref{july80006}) we obtain (the sum on the right-hand side of (\ref{july80007}) is finite)

$$
\sum_{j_1,j_3,j_4=0}^{\infty}
\int\limits_{[t,T]^6}{\bf 1}_{\{t_1<t_2<t_3\}}{\bf 1}_{\{t_4<t_5\}}
s_q(t_3,t_4,t_1)
\psi_6(t_6)\psi_5(t_5)\psi_2(t_2)
\phi_{j_3}(t_6)
\phi_{j_3}(t_3)
\times\Biggr.
$$

\vspace{1mm}
$$
\times
\phi_{j_4}(t_5)\phi_{j_4}(t_4)\phi_{j_1}(t_2)
\phi_{j_1}(t_1)dt_1 dt_2 dt_3 dt_4 dt_5 dt_6=
$$

\vspace{1mm}
\begin{equation}
\label{july80008}
=
\frac{1}{4}\int\limits_{[t,T]^3}
{\bf 1}_{\{t_2<t_6\}}
s_q(t_6,t_4,t_2)
\psi_6(t_6)
\psi_5(t_4)
\psi_2(t_2)
dt_2 dt_4 dt_6.
\end{equation}

\vspace{3mm}

Note that the equality (\ref{july80008}) remains true
when $s_q$ is a partial sum of the Fourier--Legendre series
of any function from $L_2([t,T]^3),$ i.e. the equality holds
on a dense subset in $L_2([t,T]^3).$

The right-hand side of (\ref{july80008}) defines
(as a scalar product of $s_q(t_6,t_4,t_2)$ and $\psi_6(t_6)\psi_5(t_4)\psi_2(t_2)$ $\times \frac{1}{4} {\bf 1}_{\{t_2<t_6\}}$
in $L_2([t, T]^3)$) a linear bounded (and therefore continuous)
functional in $L_2([t, T]^3),$
which is given by the function $\frac{1}{4}{\bf 1}_{\{t_2<t_6\}}\psi_6(t_6)\psi_5(t_4)\psi_2(t_2)$.
On the left-hand side of (\ref{july80008}) (by virtue of the equality (\ref{july80008}))
there is a linear continuous functional on a dense subset in 
$L_2([t,T]^3).$ This functional can be uniquely extended 
to a linear continuous functional in $L_2([t, T]^3)$
(see \cite{reed}, Theorem~I.7, P.~9).

Let us implement the passage to the limit $\lim\limits_{q\to\infty}$
in (\ref{july80008}) (at that we suppose that $s_q$ is defined by (\ref{july80007}))

$$
\sum_{j_1,j_3,j_4=0}^{\infty}
\int\limits_{[t,T]^6}{\bf 1}_{\{t_1<t_2<t_3<t_4<t_5\}}
\psi_6(t_6)\psi_5(t_5)\bar \psi_4(t_4)
\bar \psi_3(t_3)\psi_2(t_2)\bar \psi_1(t_1)
\phi_{j_3}(t_6)
\phi_{j_3}(t_3)
\times\Biggr.
$$

\vspace{1mm}
$$
\times
\phi_{j_4}(t_5)\phi_{j_4}(t_4)\phi_{j_1}(t_2)
\phi_{j_1}(t_1)dt_1 dt_2 dt_3 dt_4 dt_5 dt_6=
$$

\vspace{1mm}
\begin{equation}
\label{july80009}
=
\frac{1}{4}\int\limits_{[t,T]^3}
{\bf 1}_{\{t_2<t_6\}}{\bf 1}_{\{t_6<t_4\}}
\psi_6(t_6)\bar \psi_3(t_6)
\psi_5(t_4)\bar \psi_4(t_4)
\psi_2(t_2)\bar \psi_1(t_2)
dt_2 dt_4 dt_6.
\end{equation}

\vspace{3mm}

Rewrite (\ref{july80009}) in the form

$$
\sum_{j_1,j_3,j_4=0}^{\infty}
\int\limits_{[t,T]^6}{\bf 1}_{\{t_1<t_2<t_3<t_4<t_5\}}
\psi_6(t_6)\psi_5(t_5)\psi_4(t_4)
\psi_3(t_3)\psi_2(t_2)\psi_1(t_1)
\phi_{j_3}(t_6)
\phi_{j_3}(t_3)
\times\Biggr.
$$

\vspace{1mm}

$$
\times
\phi_{j_4}(t_5)\phi_{j_4}(t_4)\phi_{j_1}(t_2)
\phi_{j_1}(t_1)dt_1 dt_2 dt_3 dt_4 dt_5 dt_6=
$$

\vspace{1mm}
\begin{equation}
\label{july80009x}
=
\frac{1}{4}\int\limits_{[t,T]^3}
{\bf 1}_{\{t_2<t_6\}}{\bf 1}_{\{t_6<t_4\}}
\psi_6(t_6)\psi_3(t_6)
\psi_5(t_4)\psi_4(t_4)
\psi_2(t_2)\psi_1(t_2)
dt_2 dt_4 dt_6,
\end{equation}

\vspace{3mm}
\noindent
where $\psi_1(\tau),\ldots,\psi_6(\tau)\in L_2([t,T]).$

\vspace{2mm}

{\bf Step~4.}\ Suppose that $\psi_5(\tau),\psi_6(\tau),\psi_2(\tau)$ are 
Legendre polynomials of finite degrees.
Denote

\vspace{-1mm}
\begin{equation}
\label{july80010}
s_q(t_5,t_6,t_2)=\sum\limits_{l_1,l_2.l_3=0}^q
C_{l_3 l_2 l_1}\bar \phi_{l_1}(t_5)\bar \phi_{l_2}(t_6) \bar \phi_{l_3}(t_2),
\end{equation}

\vspace{3mm}
\noindent 
where $\left\{\bar \phi_j(x)\right\}_{j=0}^{\infty}$ as in (\ref{july80003})
and
$C_{l_3 l_2 l_1}$ are Fourier--Legendre coefficients for the function
$g(t_5,t_6,t_2)=\bar \psi_5(t_5)\bar \psi_6(t_6)\bar \psi_2(t_2){\bf 1}_{\{t_5<t_6\}}$ 
($\bar \psi_5(\tau), \bar \psi_6(\tau), \bar \psi_2(\tau)\in L_2([t,T])),$ i.e.
$\lim\limits_{q\to\infty}\left\Vert s_q - g\right\Vert^2_{L_2([t,T]^3)}=0.$

From (\ref{july80009x}) we obtain (the sum on the right-hand side of (\ref{july80010}) is finite)

$$
\sum_{j_1,j_3,j_4=0}^{\infty}
\int\limits_{[t,T]^6}{\bf 1}_{\{t_1<t_2<t_3<t_4<t_5\}}
s_q(t_5,t_6,t_2)
\psi_4(t_4)
\psi_3(t_3)\psi_1(t_1)
\phi_{j_3}(t_6)
\phi_{j_3}(t_3)
\times\Biggr.
$$

\vspace{1mm}
$$
\times
\phi_{j_4}(t_5)\phi_{j_4}(t_4)\phi_{j_1}(t_2)
\phi_{j_1}(t_1)dt_1 dt_2 dt_3 dt_4 dt_5 dt_6=
$$

\vspace{1mm}
\begin{equation}
\label{july80011}
=
\frac{1}{4}\int\limits_{[t,T]^3}
{\bf 1}_{\{t_2<t_6\}}{\bf 1}_{\{t_6<t_4\}}
s_q(t_4,t_6,t_2)
\psi_3(t_6)
\psi_4(t_4)
\psi_1(t_2)
dt_2 dt_4 dt_6.
\end{equation}

\vspace{3mm}

Note that the equality (\ref{july80011}) remains true
when $s_q$ is a partial sum of the Fourier--Legendre series
of any function from $L_2([t,T]^3),$ i.e. the equality holds
on a dense subset in $L_2([t,T]^3).$

The right-hand side of (\ref{july80011}) defines
(as a scalar product of $s_q(t_4,t_6,t_2)$ and $\psi_3(t_6)\psi_4(t_4)\psi_1(t_2)$ 
$\times \frac{1}{4} {\bf 1}_{\{t_2<t_6\}}{\bf 1}_{\{t_6<t_4\}}$
in $L_2([t, T]^3)$) a linear bounded (and therefore continuous)
functional in $L_2([t, T]^3),$
which is given by the function $\frac{1}{4} {\bf 1}_{\{t_2<t_6\}}{\bf 1}_{\{t_6<t_4\}}\psi_3(t_6)\psi_4(t_4)\psi_1(t_2)$.
On the left-hand side of (\ref{july80011}) (by virtue of the equality (\ref{july80011}))
there is a linear continuous functional on a dense subset in 
$L_2([t,T]^3).$ This functional can be uniquely extended 
to a linear continuous functional in $L_2([t, T]^3)$
(see \cite{reed}, Theorem~I.7, P.~9).

Let us implement the passage to the limit $\lim\limits_{q\to\infty}$
in (\ref{july80011}) (at that we suppose that $s_q$ is defined by (\ref{july80010}))

$$
\sum_{j_1,j_3,j_4=0}^{\infty}
\int\limits_{[t,T]^6}\hspace{-2mm}{\bf 1}_{\{t_1<t_2<t_3<t_4<t_5<t_6\}}
\bar \psi_6(t_6)\bar \psi_5(t_5)\psi_4(t_4)
\psi_3(t_3)\bar \psi_2(t_2)\psi_1(t_1)
\phi_{j_3}(t_6)
\phi_{j_3}(t_3)
\times\Biggr.
$$

\vspace{1mm}
$$
\times
\phi_{j_4}(t_5)\phi_{j_4}(t_4)\phi_{j_1}(t_2)
\phi_{j_1}(t_1)dt_1 dt_2 dt_3 dt_4 dt_5 dt_6=
$$

\vspace{1mm}
\begin{equation}
\label{july80012}
=
\frac{1}{4}\int\limits_{[t,T]^3}
{\bf 1}_{\{t_2<t_6\}}{\bf 1}_{\{t_6<t_4\}}{\bf 1}_{\{t_4<t_6\}}
\bar \psi_6(t_6)\psi_3(t_6)
\bar \psi_5(t_4)\psi_4(t_4)
\bar \psi_2(t_2)\psi_1(t_2)
dt_2 dt_4 dt_6=0.
\end{equation}

\vspace{3mm}

It is obvious that the equality (\ref{july80012}) (up to notations)
is (\ref{july80013}). The equality (\ref{july80013}) is proved.

As a second example, we will prove the equality (\ref{july15}).
In this case, we will use the same approach as in the proof
of equality (\ref{july80013}). Thus, we prove that

\vspace{-1mm}
\begin{equation}
\label{july80015}
\lim\limits_{p\to\infty}
\sum\limits_{j_1, j_2=0}^{p}
C_{j_2 j_1 j_2 j_1}=0.
\end{equation}

\vspace{4mm}

{\bf Step~1.}\ Using generalized Parseval's equality, we obtain

\begin{equation}
\label{july80101}
\lim\limits_{p\to\infty}\sum_{j_1,j_2=0}^{p}
\int\limits_t^T
\psi_4(t_4)\phi_{j_2}(t_4)\int\limits_t^{T}
\psi_3(t_3)\phi_{j_1}(t_3)\int\limits_t^{T}
\psi_2(t_2)\phi_{j_2}(t_2)
\int\limits_{t}^{T}
\psi_1(t_1)\phi_{j_1}(t_1)
dt_1 dt_2 dt_3 dt_4=
\end{equation}

\vspace{1mm}
$$
=\lim\limits_{p\to\infty}\sum_{j_2=0}^{p}
\int\limits_t^T
\psi_4(t_4)\phi_{j_2}(t_4)dt_4
\int\limits_t^{T}
\psi_2(t_2)\phi_{j_2}(t_2)dt_2\times
$$

\vspace{1mm}
$$
\times
\lim\limits_{p\to\infty}\sum_{j_1=0}^{p}
\int\limits_t^{T}
\psi_3(t_3)\phi_{j_1}(t_3)dt_3
\int\limits_{t}^{T}
\psi_1(t_1)\phi_{j_1}(t_1)dt_1=
$$

\vspace{2mm}
\begin{equation}
\label{july80016}
=
\int\limits_t^T
\psi_4(t_4)\psi_2(t_4)dt_4 
\int\limits_t^{T}
\psi_3(t_3)\psi_1(t_3)dt_3.
\end{equation}

\vspace{3mm}

Rewrite the equality (\ref{july80016}) in the form

$$
\sum_{j_1,j_2=0}^{\infty}
\int\limits_{[t,T]^4}
\psi_4(t_4)\psi_3(t_3)\psi_2(t_2)\psi_1(t_1)
\phi_{j_2}(t_4)
\phi_{j_1}(t_3)
\phi_{j_2}(t_2)
\phi_{j_1}(t_1)
dt_1 dt_2 dt_3 dt_4=
$$

\vspace{1mm}
\begin{equation}
\label{july80017}
=\int\limits_{[t,T]^2}
\psi_4(t_4)\psi_2(t_4)
\psi_3(t_2)\psi_1(t_2)dt_2 dt_4.
\end{equation}

\vspace{3mm}

{\bf Step~2.}\ Suppose that $\psi_1(\tau), \psi_2(\tau)$ are Legendre polynomials of finite degrees.
Denote

\vspace{-1mm}
$$
s_q(t_1,t_2)=\sum\limits_{l_1,l_2=0}^q
C_{l_2 l_1}\bar \phi_{l_1}(t_1)\bar \phi_{l_2}(t_2),
$$

\vspace{3mm}
\noindent 
where $\left\{\bar \phi_j(x)\right\}_{j=0}^{\infty}$ as in (\ref{july80003}),
$C_{l_2 l_1}$ are Fourier--Legendre coefficients for the function
$g(t_1,t_2)=\bar \psi_1(t_1)\bar \psi_2(t_2){\bf 1}_{\{t_1<t_2\}}$ 
($\bar \psi_1(\tau), \bar \psi_2(\tau)\in L_2([t,T])).$

From (\ref{july80017}) we obtain 

$$
\sum_{j_1,j_2=0}^{\infty}
\int\limits_{[t,T]^4}
s_q(t_1,t_2)\psi_4(t_4)\psi_3(t_3)
\phi_{j_2}(t_4)
\phi_{j_1}(t_3)
\phi_{j_2}(t_2)
\phi_{j_1}(t_1)
dt_1 dt_2 dt_3 dt_4=
$$

\vspace{1mm}
\begin{equation}
\label{july80019}
=\int\limits_{[t,T]^2}
s_q(t_2,t_4)\psi_4(t_4)
\psi_3(t_2)dt_2 dt_4.
\end{equation}

\vspace{3mm}

The left-hand and right-hand sides of (\ref{july80019}) define
linear continuous functionals in $L_2([t, T]^2)$ 
(see explanation earlier in this section).
Let us implement the passage to the limit $\lim\limits_{q\to\infty}$
in (\ref{july80019})

$$
\sum_{j_1,j_2=0}^{\infty}
\int\limits_{[t,T]^4}
{\bf 1}_{\{t_1<t_2\}}\psi_4(t_4)\psi_3(t_3)\bar \psi_2(t_2)\bar \psi_1(t_1)
\phi_{j_2}(t_4)
\phi_{j_1}(t_3)
\phi_{j_2}(t_2)
\phi_{j_1}(t_1)
dt_1 dt_2 dt_3 dt_4=
$$

\vspace{1mm}
\begin{equation}
\label{july80020}
=\int\limits_{[t,T]^2}
{\bf 1}_{\{t_2<t_4\}}\psi_4(t_4)\bar \psi_2(t_4)
\psi_3(t_2)\bar \psi_1(t_2)dt_2 dt_4.
\end{equation}

\vspace{3mm}

Rewrite the equality (\ref{july80020}) in the form

$$
\sum_{j_1,j_2=0}^{\infty}
\int\limits_{[t,T]^4}
{\bf 1}_{\{t_1<t_2\}}\psi_4(t_4)\psi_3(t_3)\psi_2(t_2)\psi_1(t_1)
\phi_{j_2}(t_4)
\phi_{j_1}(t_3)
\phi_{j_2}(t_2)
\phi_{j_1}(t_1)
dt_1 dt_2 dt_3 dt_4=
$$

\vspace{1mm}
\begin{equation}
\label{july80021}
=\int\limits_{[t,T]^2}
{\bf 1}_{\{t_2<t_4\}}\psi_4(t_4)\psi_2(t_4)
\psi_3(t_2)\psi_1(t_2)dt_2 dt_4,
\end{equation}

\vspace{3mm}
\noindent
where $\psi_1(\tau),\ldots,\psi_4(\tau)\in L_2([t, T]).$

\vspace{2mm}

{\bf Step~3.}\ Suppose that $\psi_2(\tau), \psi_3(\tau)$ are Legendre polynomials of finite degrees.
Denote

\vspace{-1mm}
$$
s_q(t_2,t_3)=\sum\limits_{l_1,l_2=0}^q
C_{l_2 l_1}\bar \phi_{l_1}(t_2)\bar \phi_{l_2}(t_3),
$$

\vspace{3mm}
\noindent 
where $\left\{\bar \phi_j(x)\right\}_{j=0}^{\infty}$ as in (\ref{july80003}),
$C_{l_2 l_1}$ are Fourier--Legendre coefficients for the function
$g(t_2,t_3)=\bar \psi_2(t_2)\bar \psi_3(t_3){\bf 1}_{\{t_2<t_3\}}$ 
($\bar \psi_2(\tau), \bar \psi_3(\tau)\in L_2([t,T])).$

From (\ref{july80021}) we obtain 

$$
\sum_{j_1,j_2=0}^{\infty}
\int\limits_{[t,T]^4}
{\bf 1}_{\{t_1<t_2\}}
s_q(t_2,t_3)
\psi_4(t_4)\psi_1(t_1)
\phi_{j_2}(t_4)
\phi_{j_1}(t_3)
\phi_{j_2}(t_2)
\phi_{j_1}(t_1)
dt_1 dt_2 dt_3 dt_4=
$$

\vspace{1mm}
\begin{equation}
\label{july80022}
=\int\limits_{[t,T]^2}
{\bf 1}_{\{t_2<t_4\}}s_q(t_4,t_2)\psi_4(t_4)
\psi_1(t_2)dt_2 dt_4.
\end{equation}

\vspace{3mm}

The left-hand and right-hand sides of (\ref{july80022}) define
linear continuous functionals in $L_2([t, T]^2)$.
Let us implement the passage to the limit $\lim\limits_{q\to\infty}$
in (\ref{july80022})

$$
\sum_{j_1,j_2=0}^{\infty}
\int\limits_{[t,T]^4}
{\bf 1}_{\{t_1<t_2<t_3\}}
\psi_4(t_4)\bar \psi_3(t_3)\bar \psi_2(t_2) \psi_1(t_1)
\phi_{j_2}(t_4)
\phi_{j_1}(t_3)
\phi_{j_2}(t_2)
\phi_{j_1}(t_1)
dt_1 dt_2 dt_3 dt_4=
$$

\vspace{1mm}
\begin{equation}
\label{july80023}
=\int\limits_{[t,T]^2}
{\bf 1}_{\{t_2<t_4\}}{\bf 1}_{\{t_4<t_2\}}\psi_4(t_4)\bar \psi_2(t_4)
\bar \psi_3(t_2)\psi_1(t_2)dt_2 dt_4=0.
\end{equation}

\vspace{3mm}

Rewrite the equality (\ref{july80023}) in the form

\begin{equation}
\label{july80024}
\sum_{j_1,j_2=0}^{\infty}
\int\limits_{[t,T]^4}
{\bf 1}_{\{t_1<t_2<t_3\}}
\psi_4(t_4)\psi_3(t_3)\psi_2(t_2) \psi_1(t_1)
\phi_{j_2}(t_4)
\phi_{j_1}(t_3)
\phi_{j_2}(t_2)
\phi_{j_1}(t_1)
dt_1 dt_2 dt_3 dt_4=0.
\end{equation}

\vspace{3mm}

{\bf Step~4.}\ Suppose that $\psi_3(\tau), \psi_4(\tau)$ are Legendre polynomials of finite degrees.
Denote

\vspace{-1mm}
$$
s_q(t_3,t_4)=\sum\limits_{l_1,l_2=0}^q
C_{l_2 l_1}\bar \phi_{l_1}(t_3)\bar \phi_{l_2}(t_4),
$$

\vspace{3mm}
\noindent 
where $\left\{\bar \phi_j(x)\right\}_{j=0}^{\infty}$ as in (\ref{july80003}),
$C_{l_2 l_1}$ are Fourier--Legendre coefficients for the function
$g(t_3,t_4)=\bar \psi_3(t_3)\bar \psi_4(t_4){\bf 1}_{\{t_3<t_4\}}$ 
($\bar \psi_3(\tau), \bar \psi_4(\tau)\in L_2([t,T])).$

From (\ref{july80024}) we obtain 

\begin{equation}
\label{july80025}
\sum_{j_1,j_2=0}^{\infty}
\int\limits_{[t,T]^4}
{\bf 1}_{\{t_1<t_2<t_3\}}
s_q(t_3,t_4)\psi_2(t_2) \psi_1(t_1)
\phi_{j_2}(t_4)
\phi_{j_1}(t_3)
\phi_{j_2}(t_2)
\phi_{j_1}(t_1)
dt_1 dt_2 dt_3 dt_4=0.
\end{equation}

\vspace{3mm}

The left-hand and right-hand sides of (\ref{july80025}) define
linear continuous functionals in $L_2([t, T]^2)$
(we interpret the right-hand side of (\ref{july80025})
as a zero functional in $L_2([t, T]^2)$).
Let us implement the passage to the limit $\lim\limits_{q\to\infty}$
in (\ref{july80025})

\begin{equation}
\label{july80026}
\sum_{j_1,j_2=0}^{\infty}
\int\limits_{[t,T]^4}
{\bf 1}_{\{t_1<t_2<t_3<t_4\}}
\bar \psi_4(t_4) \bar \psi_3(t_3)\psi_2(t_2) \psi_1(t_1)
\phi_{j_2}(t_4)
\phi_{j_1}(t_3)
\phi_{j_2}(t_2)
\phi_{j_1}(t_1)
dt_1 dt_2 dt_3 dt_4=0.
\end{equation}

\vspace{3mm}

It is easy to see that the equality (\ref{july80026}) (up to notations)
is the equality (\ref{july15}).
The equality (\ref{july15}) is proved.

Let us formulate the ideas used when considering
the two above examples in the form of an algorithm.

\vspace{2mm}

{\bf Step~1.}\ Suppose $k=2r$ $(r=2,3,4,\ldots)$, where $r$ is the number of pairs
$\{g_1, g_2\}, \ldots, 
\{g_{2r-1}, g_{2r}\}$ (see (\ref{leto5007xxx})). Let us select blocks
in the multi-index $j_k\ldots j_1$ that
correspond to the fulfillment of the condition
$$
\prod\limits_{l=1}^{r_d} {\bf 1}_{\{g_{2l}=g_{2l-1}+1\}}=1,
$$

\vspace{2mm}
\noindent
where $r_d$ is the number of pairs (see (\ref{leto5007xxx}))
in the block with number $d.$

\vspace{2mm}

{\bf Step~2.}\ Let us write the Volterra--type kernel (\ref{july7000}) in the form

\vspace{-1mm}
\begin{equation}
\label{july80027}
K(t_1,\ldots,t_k)=
\psi_1(t_1)\ldots \psi_k(t_k){\bf 1}_{\{t_1<t_{2}\}}{\bf 1}_{\{t_2<t_{3}\}}\ldots
{\bf 1}_{\{t_{k-1}<t_{k}\}},
\end{equation}

\vspace{3mm}
\noindent
where $\psi_1(\tau),\ldots,\psi_k(\tau)\in L_2([t,T])$,
$t_1,\ldots,t_k\in [t, T],$ $k\ge 4.$

Let us save multipliers of the form ${\bf 1}_{\{t_n<t_{n+1}\}}$
in the expression (\ref{july80027}) that correspond 
to the above blocks. At that, we remove the remaining 
multipliers of the form 
${\bf 1}_{\{t_n<t_{n+1}\}}$ from the expression (\ref{july80027}).
As a result, we get a modified kernel $\bar K(t_1,\ldots,t_k)$.
Let us write an analogue of the left-hand side
of equality (\ref{july80000}) for the modified kernel $\bar K(t_1,\ldots,t_k)$
(see (\ref{july80100}) and (\ref{july80101}) as examples).
For definiteness, let us denote this expression by
$({}^{-})$.

\vspace{2mm}

{\bf Step~3.}\ Using generalized Parseval's equality and (\ref{july30016}), we represent
the expression $({}^{-})$ as an integral over the hypercube $[t, T]^r$
(see the right-hand sides of (\ref{july80002}) and (\ref{july80017}) as examples).
For definiteness, let us denote the obtained equality by
$(\bar K)$ ((\ref{july80002}) and (\ref{july80017}) are examples of $(\bar K)$).

\vspace{2mm}

{\bf Step~4.}\ Further, transformations and passages to the limit
in the equality $(\bar K)$ are performed iteratively 
in such a way as to restore the removed multipliers ${\bf 1}_{\{t_n<t_{n+1}\}}$
on the left-hand side of $(\bar K)$
(for more details, see the proof of formulas (\ref{july80013}), (\ref{july80015})).
As a result, we obtain the equality (\ref{july80000}).
More precisely, we can move from right to left
along a multi-index corresponding to the left-hand side of $(\bar K)$.
Let us assume that at the $n$-th step we need to restore
the multiplier ${\bf 1}_{\{t_n<t_{n+1}\}}$. 
Then the function $g$ (see the proof of formulas (\ref{july80013}), (\ref{july80015})) 
will be the product
of ${\bf 1}_{\{t_n<t_{n+1}\}}\psi_n(t_n)\psi_{n+1}(t_{n+1})$
and $r-2$ weight functions that are chosen so that
on the right-hand side of the equality 
$(\bar K)$ there is a scalar product in $L_2([t, T]^r)$
involving $s_q$ ($s_q$ is an approximation of $g$).

\vspace{2mm}

Using the above algorithm, we prove the equality (\ref{july90000})
for the case $k=2r$ $(r=2,3,\ldots).$
The equality (\ref{july90000}) is proved.

Note that the series on the left-hand side of 
(\ref{july90000}) converges absolutly since
its sum does not depend 
on permutations of basis functions
(here the basis in $L_2([t,T]^{r})$ is
$\left\{\phi_{j_1}(x_1)\ldots \phi_{j_r}(x_r)\right\}_{j_1,\ldots,j_r=0}^{\infty}$).

\vspace{5mm}

\section{Revision of Hypotheses on Expansion of Iterated Stra\-to\-no\-vich Stochastic
Integrals of Multiplicity $k$ $(k\in\mathbb{N})$}

\vspace{5mm}

In Sect.~2, we formulated three 
hypotheses on expansion of iterated Stra\-to\-no\-vich stochastic 
integrals based on the results obtained by the author
in the 2010s. In light of recent results (Theorems~13--55),
a new vision of the above problem 
has appeared.
In particular, it became clear that it is possible 
to methodically obtain results related to the expansion
of iterated Stratonovich stochastic 
integrals for the case of an 
arbitrary 
complete ortho\-nor\-mal system of functions in the space $L_2([t, T])$
and $\psi_1(\tau),$ $\ldots,$ $\psi_k(\tau)$ $\in $ $L_2([t, T])$.

Definition of the Stratonovich stochastic 
integral from \cite{KlPl2} (also see \cite{2018a}, Sect.2.1), 
which we mainly use in this article, imposes
its own limitations. In particular, this definition assumes that
$\psi_1(\tau),$ $\ldots,$ $\psi_k(\tau)$ are continuous functions 
at the interval $[t, T]$.

Based on Theorems~42, 47, 48, 51--55, we formulate the following
hypothesis on expansion of the sum $\bar J^{*}[\psi^{(k)}]_{T,t}^{(i_1\ldots i_k)}$
of iterated Ito stochastic integrals (see (\ref{dsds9})).

\vspace{2mm}

{\bf Hypothesis~7.}\ {\it Suppose that $\{\phi_j(x)\}_{j=0}^{\infty}$
is an arbitrary complete orthonormal system of functions
in $L_2([t, T])$ and
$\psi_1(\tau),\ldots, \psi_k(\tau)\in L_2([t, T]).$
Then$,$ for the sum $\bar J^{*}[\psi^{(k)}]_{T,t}^{(i_1\ldots i_k)}$
of iterated Ito stochastic integrals 

\vspace{1mm}
$$
\bar J^{*}[\psi^{(k)}]_{T,t}^{(i_1\ldots i_k)}=J[\psi^{(k)}]_{T,t}^{(i_1\ldots i_k)}+
\sum_{r=1}^{\left[k/2\right]}\frac{1}{2^r}
\sum_{(s_r,\ldots,s_1)\in {\rm A}_{k,r}}
J[\psi^{(k)}]_{T,t}^{s_r,\ldots,s_1}
$$

\vspace{3mm}
\noindent
the following 
expansion 

\vspace{-1mm}

$$
\bar J^{*}[\psi^{(k)}]_{T,t}^{(i_1\ldots i_k)}=
\hbox{\vtop{\offinterlineskip\halign{
\hfil#\hfil\cr
{\rm l.i.m.}\cr
$\stackrel{}{{}_{p_1,\ldots,p_k\to \infty}}$\cr
}} }
\sum_{j_1=0}^{p_1}\ldots\sum_{j_k=0}^{p_k}
C_{j_k \ldots j_1}\prod\limits_{l=1}^k \zeta_{j_l}^{(i_l)}
$$

\vspace{4mm}
\noindent
that converges in the mean-square sense is valid, where 

\vspace{-1mm}
$$
C_{j_k \ldots j_1}=\int\limits_t^T\psi_k(t_k)\phi_{j_k}(t_k)\ldots
\int\limits_t^{t_2}
\psi_1(t_1)\phi_{j_1}(t_1)
dt_1\ldots dt_k
$$

\vspace{3mm}
\noindent
is the Fourier coefficient, 
${\rm l.i.m.}$ is a limit in the mean-square sense,
$i_1, \ldots, i_k=0, 1,\ldots,m,$

\vspace{-1mm}
$$
\zeta_{j}^{(i)}=
\int\limits_t^T \phi_{j}(\tau) d{\bf w}_{\tau}^{(i)}
$$ 

\vspace{2mm}
\noindent
are independent standard Gaussian random variables for various 
$i$ or $j$ {\rm (}in the case when $i\ne 0${\rm )},
${\bf w}_{\tau}^{(i)}={\bf f}_{\tau}^{(i)}$ 
for $i=1,\ldots,m$ and 
${\bf w}_{\tau}^{(0)}=\tau;$ another notations are the same as
in Theorem~{\rm 5}.}

\vspace{2mm}

Using Theorem~5, we obtain the following hypothesis.

\vspace{2mm}

{\bf Hypothesis~8.}\ {\it Suppose that $\{\phi_j(x)\}_{j=0}^{\infty}$
is an arbitrary complete orthonormal system of functions
in $L_2([t, T])$ and
$\psi_1(\tau),\ldots, \psi_k(\tau)$ are continuous functions
at the interval $[t, T].$
Then$,$ for the iterated Stratonovich sto\-chas\-tic integral 
of arbitrary multiplicity $k$

\vspace{1mm}
$$
J^{*}[\psi^{(k)}]_{T,t}^{(i_1\ldots i_k)}=
{\int\limits_t^{*}}^T
\psi_k(t_k) \ldots 
{\int\limits_t^{*}}^{t_{2}}
\psi_1(t_1) d{\bf w}_{t_1}^{(i_1)}\ldots
d{\bf w}_{t_k}^{(i_k)}
$$

\vspace{3mm}
\noindent
the following 
expansion 

\vspace{-1mm}
$$
J^{*}[\psi^{(k)}]_{T,t}^{(i_1\ldots i_k)}=
\hbox{\vtop{\offinterlineskip\halign{
\hfil#\hfil\cr
{\rm l.i.m.}\cr
$\stackrel{}{{}_{p_1,\ldots,p_k\to \infty}}$\cr
}} }
\sum\limits_{j_1=0}^{p_1}\ldots\sum\limits_{j_k=0}^{p_k}
C_{j_k \ldots j_1}\prod\limits_{l=1}^k \zeta_{j_l}^{(i_l)}
$$

\vspace{4mm}
\noindent
that converges in the mean-square sense is valid, where 

\vspace{-1mm}
$$
C_{j_k \ldots j_1}=\int\limits_t^T\psi_k(t_k)\phi_{j_k}(t_k)\ldots
\int\limits_t^{t_2}
\psi_1(t_1)\phi_{j_1}(t_1)
dt_1\ldots dt_k
$$

\vspace{3mm}
\noindent
is the Fourier coefficient, 
${\rm l.i.m.}$ is a limit in the mean-square sense,
$i_1, \ldots, i_k=0, 1,\ldots,m,$

\vspace{-1mm}
$$
\zeta_{j}^{(i)}=
\int\limits_t^T \phi_{j}(\tau) d{\bf w}_{\tau}^{(i)}
$$ 

\vspace{2mm}
\noindent
are independent standard Gaussian random variables for various 
$i$ or $j$ {\rm (}in the case when $i\ne 0${\rm )},
${\bf w}_{\tau}^{(i)}={\bf f}_{\tau}^{(i)}$ 
for $i=1,\ldots,m$ and 
${\bf w}_{\tau}^{(0)}=\tau.$}

\vspace{5mm}

\section{Proof of Hypotheses~7 and 8 Under the Condition (\ref{july90001})
for the Case $k\ge 2r,$ $p_1=\ldots=p_k=p$ and Under Some Additional Assumptions}

\vspace{5mm}

Suppose that the equality

\vspace{1mm}
$$
\lim\limits_{p\to\infty}
\sum\limits_{j_{g_1}, j_{g_3},\ldots,j_{g_{2r-1}}=0}^p
C_{j_k\ldots j_1}\biggl|_{j_{g_1}=j_{g_2},\ldots, j_{g_{2r-1}}=j_{g_{2r}}}=
$$

\vspace{2mm}
\begin{equation}
\label{july90001}
=\frac{1}{2^r} \prod\limits_{l=1}^r {\bf 1}_{\{g_{2l}=g_{2l-1}+1\}}
C_{j_k \ldots j_1}\biggl|_{(j_{g_2} j_{g_1})\curvearrowright (\cdot)
\ldots (j_{g_{2r}} j_{g_{2r-1}})\curvearrowright (\cdot),
j_{g_{{}_{1}}}=~j_{g_{{}_{2}}},\ldots, j_{g_{{}_{2r-1}}}=~j_{g_{{}_{2r}}}
}\biggr.
\end{equation}

\vspace{4mm}                                           
\noindent
is satisfied for all possible $g_1,g_2,\ldots,g_{2r-1},g_{2r}$ (see {\rm (\ref{leto5007xxx})),
where $k\ge 2r,$ $r=1,2,\ldots,[k/2],$ $C_{j_k\ldots j_1}$ is defined by (\ref{july15030}),
another notations are the same as in Theorem~47. Recall that the case
$k=2r$ is considered in Sect.~26.

Moreover, suppose that
the series 

\vspace{1mm}
$$
\lim\limits_{p\to\infty}\sum\limits_{j_{g_1}, j_{g_3},\ldots,j_{g_{2r-1}}=0}^p
C_{j_k\ldots j_1}\biggl|_{j_{g_1}=j_{g_2},\ldots, j_{g_{2r-1}}=j_{g_{2r}}}
$$

\vspace{4mm}
\noindent
converges absolutly for any fixed $j_1,\ldots,j_q,\ldots,j_k,$
where $q\ne g_1, g_2, \ldots, g_{2r-1},$ $g_{2r}$ and $k>2r.$
It should be noted that the above assumptions will be proved further (see Sect.~29).

Hypotheses~7 and 8 will be proved for the case $p_1=\ldots=p_k=p$
if we prove that (see 
Theorem~47 for the case $p_1=\ldots=p_k=p$)

\vspace{1mm}
$$
\lim\limits_{p\to\infty}
\sum\limits_{\stackrel{j_1,\ldots,j_q,\ldots,j_k=0}{{}_{q\ne g_1, g_2, \ldots, g_{2r-1},
g_{2r}}}}^p
\Biggl(\sum\limits_{j_{g_1}, j_{g_3},\ldots,j_{g_{2r-1}}=0}^p
C_{j_k\ldots j_1}\biggl|_{j_{g_1}=j_{g_2},\ldots, j_{g_{2r-1}}=j_{g_{2r}}}-\Biggr.
$$

\vspace{2mm}
\begin{equation}
\label{july90025}
\Biggl.-\frac{1}{2^r} \prod\limits_{l=1}^r {\bf 1}_{\{g_{2l}=g_{2l-1}+1\}}
C_{j_k \ldots j_1}\biggl|_{(j_{g_2} j_{g_1})\curvearrowright (\cdot)
\ldots (j_{g_{2r}} j_{g_{2r-1}})\curvearrowright (\cdot),
j_{g_{{}_{1}}}=~j_{g_{{}_{2}}},\ldots, j_{g_{{}_{2r-1}}}=~j_{g_{{}_{2r}}}
}\biggr.\Biggr)^2=0
\end{equation}

\vspace{5mm}
\noindent 
for all $r=1,2,\ldots,[k/2],$
where notations are the same as in 
(\ref{july90001}).

Further, we have

\vspace{1mm}
$$
\sum\limits_{\stackrel{j_1,\ldots,j_q,\ldots,j_k=0}{{}_{q\ne g_1, g_2, \ldots, g_{2r-1},
g_{2r}}}}^p
\Biggl(\sum\limits_{j_{g_1}, j_{g_3},\ldots,j_{g_{2r-1}}=0}^p
C_{j_k\ldots j_1}\biggl|_{j_{g_1}=j_{g_2},\ldots, j_{g_{2r-1}}=j_{g_{2r}}}-\Biggr.
$$

\vspace{2mm}
$$
\Biggl.-\frac{1}{2^r} \prod\limits_{l=1}^r {\bf 1}_{\{g_{2l}=g_{2l-1}+1\}}
C_{j_k \ldots j_1}\biggl|_{(j_{g_2} j_{g_1})\curvearrowright (\cdot)
\ldots (j_{g_{2r}} j_{g_{2r-1}})\curvearrowright (\cdot),
j_{g_{{}_{1}}}=~j_{g_{{}_{2}}},\ldots, j_{g_{{}_{2r-1}}}=~j_{g_{{}_{2r}}}
}\biggr.\Biggr)^2\le
$$

\vspace{5mm}
$$
\le \sum\limits_{\stackrel{j_1,\ldots,j_q,\ldots,j_k=0}{{}_{q\ne g_1, g_2, \ldots, g_{2r-1},
g_{2r}}}}^{\infty}
\Biggl(\sum\limits_{j_{g_1}, j_{g_3},\ldots,j_{g_{2r-1}}=0}^p
C_{j_k\ldots j_1}\biggl|_{j_{g_1}=j_{g_2},\ldots, j_{g_{2r-1}}=j_{g_{2r}}}-\Biggr.
$$

\vspace{2mm}
\begin{equation}
\label{july90009xx}
\Biggl.-\frac{1}{2^r} \prod\limits_{l=1}^r {\bf 1}_{\{g_{2l}=g_{2l-1}+1\}}
C_{j_k \ldots j_1}\biggl|_{(j_{g_2} j_{g_1})\curvearrowright (\cdot)
\ldots (j_{g_{2r}} j_{g_{2r-1}})\curvearrowright (\cdot),
j_{g_{{}_{1}}}=~j_{g_{{}_{2}}},\ldots, j_{g_{{}_{2r-1}}}=~j_{g_{{}_{2r}}}
}\biggr.\Biggr)^2,
\end{equation}

\vspace{2mm}
\noindent
where
\begin{equation}
\label{july90009}
\sum\limits_{\stackrel{j_1,\ldots,j_q,\ldots,j_k=0}{{}_{q\ne g_1, g_2, \ldots, g_{2r-1},
g_{2r}}}}^{\infty}\stackrel{\sf def}{=}
\lim\limits_{q\to\infty}
\sum\limits_{\stackrel{j_1,\ldots,j_q,\ldots,j_k=0}{{}_{q\ne g_1, g_2, \ldots, g_{2r-1},
g_{2r}}}}^{q}.
\end{equation}

\vspace{4mm}

Consider the following analogue
of Monotone Convergence Theorem for infinite series.

\vspace{2mm}

{\bf Proposition~1.}\ {\it Suppose that $x_{m,n}\ge 0$ for all $m,n\in\mathbb{N},$

\vspace{1mm}
$$
\lim\limits_{m\to\infty}x_{m,n}=y_n\ \ \ (\hbox{for any fixed}~n\in \mathbb{N}),
$$

\vspace{4mm}
\noindent
and $x_{m,n}\le x_{m+1,n}$ for all $m\in\mathbb{N}$ and for any fixed $n\in \mathbb{N}.$
Then

\begin{equation}
\label{july90002x}
\lim_{m\to\infty}\sum\limits_{n=1}^{\infty}
x_{m,n}=\sum\limits_{n=1}^{\infty}
\lim_{m\to\infty}x_{m,n}=
\sum\limits_{n=1}^{\infty}y_n.
\end{equation}
}

\vspace{3mm}

{\bf Proof.}\ Proposition~1 can be easily proved 
using the following version of Fatou's Lemma for infinite series
\begin{equation}
\label{july90002}
\sum\limits_{n=1}^{\infty}\liminf _{m\to\infty}
x_{m,n}\le \liminf _{m\to\infty}
\sum\limits_{n=1}^{\infty}x_{m,n},
\end{equation}

\vspace{4mm}
\noindent
where it is assumed that the conditions of Proposition~1
are fulfilled. Indeed, we have

\vspace{1mm}
$$
0\le x_{m,n}\le y_n.
$$

\vspace{4mm}

Then
$$
\sum\limits_{n=1}^{\infty}x_{m,n}\le 
\sum\limits_{n=1}^{\infty}y_n
$$

\vspace{3mm}
\noindent
and (see (\ref{july90002}))

\vspace{-1mm}
\begin{equation}
\label{july90003}
\limsup_{m\to\infty}\sum\limits_{n=1}^{\infty}x_{m,n}\le 
\sum\limits_{n=1}^{\infty}y_n=
\sum\limits_{n=1}^{\infty}\liminf_{m\to\infty}x_{m,n}
\le \liminf _{m\to\infty}
\sum\limits_{n=1}^{\infty}x_{m,n}.
\end{equation}

\vspace{4mm}

From (\ref{july90003}) we get

\vspace{1mm}
$$
\sum\limits_{n=1}^{\infty}y_n=\liminf_{m\to\infty}\sum\limits_{n=1}^{\infty}x_{m,n}=
\limsup_{m\to\infty}\sum\limits_{n=1}^{\infty}x_{m,n}=
\lim_{m\to\infty}\sum\limits_{n=1}^{\infty}x_{m,n},
$$

\vspace{4mm}
\noindent
i.e. the equality (\ref{july90002x}) is proved.

To prove (\ref{july90002}) we note that

\vspace{1mm}
$$
\inf_{j\ge m}x_{j,n}\le x_{k,n}\ \ \ (\forall k\ge m).
$$

\vspace{3mm}

Then
$$
\sum\limits_{n=1}^N \inf_{j\ge m}x_{j,n}\le \sum\limits_{n=1}^N x_{k,n}\ \ \ (\forall k\ge m)
$$

\vspace{1mm}
\noindent
and
\begin{equation}
\label{july90004}
\sum\limits_{n=1}^N \inf_{j\ge m}x_{j,n}\le\inf\limits_{k\ge m}
\sum\limits_{n=1}^{N} x_{k,n}\le
\inf\limits_{k\ge m}
\sum\limits_{n=1}^{\infty} x_{k,n}.
\end{equation}

\vspace{4mm}

Passing to the limit $\lim\limits_{m\to\infty}$ in (\ref{july90004}), we obtain

\vspace{-1mm}
\begin{equation}
\label{july90005}
\sum\limits_{n=1}^N \lim\limits_{m\to\infty} \inf_{j\ge m}x_{j,n}
\le
\lim\limits_{m\to\infty} \inf\limits_{k\ge m}
\sum\limits_{n=1}^{\infty} x_{k,n}.
\end{equation}

\vspace{4mm}

Passing to the limit $\lim\limits_{N\to\infty}$ in (\ref{july90005}), we get

\vspace{-1mm}
$$
\sum\limits_{n=1}^{\infty} \lim\limits_{m\to\infty} \inf_{j\ge m}x_{j,n}
\le
\lim\limits_{m\to\infty} \inf\limits_{k\ge m}
\sum\limits_{n=1}^{\infty} x_{k,n},
$$

\vspace{4mm}
\noindent
i.e. the equality (\ref{july90002}) is satisfied.
Proposition~1 is proved.

\vspace{2mm}
         
{\bf Proposition~2.}\ {\it Suppose that 

\vspace{-3mm}
\begin{equation}
\label{july90017}
\sum\limits_{j=1}^{\infty} g_{j,n}=0,
\end{equation}

\vspace{2mm}
\noindent
the series {\rm (\ref{july90017})} converges absolutely for any fixed $n\in\mathbb{N}$ 
and

\vspace{-1mm}
\begin{equation}
\label{july500001}
\sum\limits_{n=1}^{\infty}\left(\sum\limits_{j=1}^{\infty} \left\vert g_{j,n}\right\vert\right)^2
<\infty.
\end{equation}

\vspace{2mm}

Then
\begin{equation}
\label{july90015}
\lim_{m\to\infty}\sum\limits_{n=1}^{\infty}
\left(\sum\limits_{j=1}^m g_{j,n}\right)^2=
\sum\limits_{n=1}^{\infty}\lim_{m\to\infty}\left(\sum\limits_{j=1}^{m} g_{j,n}\right)^2=0.
\end{equation}
}

\vspace{4mm}

{\bf Proof.}\ We have 
$$
g_{j,n}=g_{j,n}^{+}-g_{j,n}^{-},\ \ \ 
\left\vert g_{j,n}\right\vert =g_{j,n}^{+}+g_{j,n}^{-},
$$ 

\vspace{3mm}
\noindent
where 
$$
g_{j,n}^{+}=\max\{g_{j,n}, 0\}=\frac{1}{2}\left(\left\vert g_{j,n}\right\vert + g_{j,n}\right)\ge 0,
$$

$$
g_{j,n}^{-}=-\min\{g_{j,n}, 0\}=\frac{1}{2}\left(\left\vert g_{j,n}\right\vert - g_{j,n}\right)
\ge 0.
$$ 

\vspace{3mm}

Moreover,
\begin{equation}
\label{july90016s}
\sum\limits_{j=1}^{\infty}g_{j,n}=
\sum\limits_{j=1}^{\infty} g_{j,n}^{+}-\sum\limits_{j=1}^{\infty} g_{j,n}^{-}=0,
\end{equation}

\begin{equation}
\label{july90016}
\sum\limits_{j=1}^{\infty}\left\vert g_{j,n}\right\vert =
\sum\limits_{j=1}^{\infty} g_{j,n}^{+}+\sum\limits_{j=1}^{\infty} g_{j,n}^{-}
=2\sum\limits_{j=1}^{\infty} g_{j,n}^{+}=
2\sum\limits_{j=1}^{\infty} g_{j,n}^{-}.
\end{equation}

\vspace{3mm}

Since the series (\ref{july90017}) converges absolutely, then by virtue 
of the equality (\ref{july90016}) the series (with non-negative terms) on the right-hand side of 
(\ref{july90016}) and on the right-hand side of (\ref{july90016s}) converge. 

Further, using Proposition~1 and (\ref{july90016s}), (\ref{july90016}), we obtain

$$
\lim_{m\to\infty}\sum\limits_{n=1}^{\infty}
\left(\sum\limits_{j=1}^m g_{j,n}\right)^2=
\lim_{m\to\infty}\sum\limits_{n=1}^{\infty}
\left(\sum\limits_{j=1}^m g_{j,n}^{+}- \sum\limits_{j=1}^m g_{j,n}^{-}\right)^2=
$$

\vspace{2mm}
$$
=\lim_{m\to\infty}\sum\limits_{n=1}^{\infty}
\left(\sum\limits_{j=1}^m g_{j,n}^{+}\right)^2-
\lim_{m\to\infty}\sum\limits_{n=1}^{\infty}
\left(2\sum\limits_{j=1}^m g_{j,n}^{+}\sum\limits_{j=1}^m g_{j,n}^{-}\right)
+
\lim_{m\to\infty}\sum\limits_{n=1}^{\infty}
\left(\sum\limits_{j=1}^m g_{j,n}^{-}\right)^2=
$$

\vspace{2mm}
$$
=\sum\limits_{n=1}^{\infty}\lim_{m\to\infty}
\left(\sum\limits_{j=1}^m g_{j,n}^{+}\right)^2-
\sum\limits_{n=1}^{\infty}\lim_{m\to\infty}
\left(2\sum\limits_{j=1}^m g_{j,n}^{+}\sum\limits_{j=1}^m g_{j,n}^{-}\right)+
$$

\vspace{2mm}
$$
+
\sum\limits_{n=1}^{\infty}\lim_{m\to\infty}
\left(\sum\limits_{j=1}^m g_{j,n}^{-}\right)^2=
$$

\vspace{2mm}

$$
=\sum\limits_{n=1}^{\infty}
\left(\sum\limits_{j=1}^{\infty} g_{j,n}^{+}\right)^2-
2\sum\limits_{n=1}^{\infty}
\left(\sum\limits_{j=1}^{\infty} g_{j,n}^{+}\sum\limits_{j=1}^{\infty} g_{j,n}^{-}\right)+
\sum\limits_{n=1}^{\infty}
\left(\sum\limits_{j=1}^{\infty} g_{j,n}^{-}\right)^2=
$$

\vspace{2mm}

$$
=\frac{1}{4}\sum\limits_{n=1}^{\infty}
\left(\sum\limits_{j=1}^{\infty} \left|g_{j,n}\right|\right)^2-
\frac{1}{2}\sum\limits_{n=1}^{\infty}
\left(\sum\limits_{j=1}^{\infty} \left|g_{j,n}\right|\right)^2
+\frac{1}{4}\sum\limits_{n=1}^{\infty}
\left(\sum\limits_{j=1}^{\infty} \left|g_{j,n}\right|\right)^2=0.
$$

\vspace{4mm}

Proposition~2 is proved.

It is easy to see that by analogy with the proof of Propositions~1 and 2
the following statements can be proved.

\vspace{2mm}

{\bf Proposition 3.}\ {\it Suppose that $h_{p,k_1,\ldots,k_d}\ge 0$ 
for all $p\in \mathbb{N}$ and for any fixed $k_1,\ldots,k_d\in \mathbb{N},$

$$
\lim\limits_{p\to\infty}h_{p,k_1,\ldots,k_d}=u_{k_1,\ldots,k_d}\ \ \ (\hbox{for any fixed}~
k_1,\ldots,k_d\in \mathbb{N}),
$$

\vspace{3mm}
\noindent
and $h_{p,k_1,\ldots,k_d}\le h_{p+1,k_1,\ldots,k_d}$ for all $p\in \mathbb{N}$
and for any fixed $k_1,\ldots,k_d\in \mathbb{N}.$ 
Then

\begin{equation}
\label{july90006}
\lim_{p\to\infty}\sum\limits_{k_1,\ldots,k_d=1}^{\infty}
h_{p,k_1,\ldots,k_d}=\sum\limits_{k_1,\ldots,k_d=1}^{\infty}
\lim_{p\to\infty}h_{p,k_1,\ldots,k_d}=
\sum\limits_{k_1,\ldots,k_d=1}^{\infty}
u_{k_1,\ldots,k_d},
\end{equation}

\vspace{4mm}
\noindent
where $h_{p,k_1,\ldots,k_d}, u_{k_1,\ldots,k_d}\in \mathbb{R},$ 
$d\in \mathbb{N},$ the series on the left-hand side of {\rm (\ref{july90006})}
is understood in the same sense as in {\rm (\ref{july90009})}.
}

\vspace{2mm}

{\bf Proposition 4.}\ {\it Suppose that 

\begin{equation}
\label{july700010}
\lim\limits_{p\to\infty}
\sum\limits_{j_1,\ldots,j_q=1}^{p} h_{j_1,\ldots,j_q,k_1,\ldots,k_d}
\stackrel{\sf def}{=}\sum\limits_{j_1,\ldots,j_q=1}^{\infty} h_{j_1,\ldots,j_q,k_1,\ldots,k_d}=0,
\end{equation}

\vspace{3mm}
\noindent
the series {\rm (\ref{july700010})}
converges absolutely for any fixed $k_1,\ldots,k_d\in\mathbb{N}$ 
and

$$
\sum\limits_{k_1,\ldots,k_d=1}^{\infty} 
\left(\sum\limits_{j_1,\ldots,j_q=1}^{\infty} \left|h_{j_1,\ldots,j_q,k_1,\ldots,k_d}\right|
\right)^2<\infty.
$$

\vspace{2mm}

Then
$$
\lim\limits_{p\to\infty}
\sum\limits_{k_1,\ldots,k_d=1}^{\infty} 
\left(\sum\limits_{j_1,\ldots,j_q=1}^{p} h_{j_1,\ldots,j_q,k_1,\ldots,k_d}\right)^2=
$$

\vspace{2mm}

$$
=\sum\limits_{k_1,\ldots,k_d=1}^{\infty} 
\lim\limits_{p\to\infty}
\left(\sum\limits_{j_1,\ldots,j_q=1}^{p} h_{j_1,\ldots,j_q,k_1,\ldots,k_d}\right)^2=
0,
$$

\vspace{4mm}
\noindent
where 
$$
\lim\limits_{n\to\infty}
\sum\limits_{k_1,\ldots,k_d=1}^{n}
\stackrel{\sf def}{=}\sum\limits_{k_1,\ldots,k_d=1}^{\infty},
$$

\vspace{4mm}
\noindent
$h_{j_1,\ldots,j_q,k_1,\ldots,k_d}\in \mathbb{R}$ and
$d, q\in \mathbb{N}.$
}

\vspace{2mm}

Obviously, Proposition~4 follows from Proposition~3
in the same way as Proposition~2 follows from Proposition~1.
Applying Proposition~4 to the right-hand side of (\ref{july90009xx})
(using (\ref{july90001}) and the absolute convergence of the series on the left-hand side of 
(\ref{july90001})),
we obtain (\ref{july90025}). 
At that, we used the conditions

\vspace{-1.5mm}
\begin{equation}
\label{slovo10}
\sum\limits_{\stackrel{j_1,\ldots,j_q,\ldots,j_k=0}{{}_{q\ne g_1, g_2, \ldots, g_{2r-1},
g_{2r}}}}^{\infty}
\left(\lim\limits_{p\to\infty}\sum\limits_{j_{g_1}, j_{g_3},\ldots,j_{g_{2r-1}}=0}^p
\left\vert C_{j_k\ldots j_1}\biggl|_{j_{g_1}=j_{g_2},\ldots, j_{g_{2r-1}}=j_{g_{2r}}}
\right\vert\right)^2<\infty,
\end{equation}

\vspace{2mm}
\begin{equation}
\label{july700011}
\sum\limits_{\stackrel{j_1,\ldots,j_q,\ldots,j_k=0}{{}_{q\ne g_1, g_2, \ldots, g_{2r-1},
g_{2r}}}}^{\infty}
\left(
C_{j_k \ldots j_1}\biggl|_{(j_{g_2} j_{g_1})\curvearrowright (\cdot)
\ldots (j_{g_{2r}} j_{g_{2r-1}})\curvearrowright (\cdot),
j_{g_{{}_{1}}}=~j_{g_{{}_{2}}},\ldots, j_{g_{{}_{2r-1}}}=~j_{g_{{}_{2r}}}
}\biggr.\right)^2<\infty.
\end{equation}

\vspace{5mm}

Note that (\ref{july700011}) follows from the Parseval equality
since the expression

$$
C_{j_k \ldots j_1}\biggl|_{(j_{g_2} j_{g_1})\curvearrowright (\cdot)
\ldots (j_{g_{2r}} j_{g_{2r-1}})\curvearrowright (\cdot),
j_{g_{{}_{1}}}=~j_{g_{{}_{2}}},\ldots, j_{g_{{}_{2r-1}}}=~j_{g_{{}_{2r}}}
}\biggr.\stackrel{\sf def}{=}H_{j_{q_1}\ldots j_{q_{k-2r}}}
$$

\vspace{3mm}
\noindent
is a finite linear combination 
of the Fourier coefficients
of $L_2([t,T]^{k-2r})$--func\-ti\-ons
after iteratively
applying
transformations (\ref{july100003}), (\ref{july100004}) (see Sect.~29) to
$H_{j_{q_1}\ldots j_{q_{k-2r}}}$
for integrations not involving the basis functions
$\phi_{j_{q_1}},\ldots, \phi_{j_{q_{k-2r}}}.$

Let us consider another sufficient condition under which the equality
(\ref{july90025}) is satisfied.
Suppose that 

\vspace{-1mm}
$$
\exists\ \ \ \lim\limits_{p,q\to\infty}
\sum\limits_{\stackrel{j_1,\ldots,j_m,\ldots,j_k=0}{{}_{m\ne g_1, g_2, \ldots, g_{2r-1},
g_{2r}}}}^q
\Biggl(\sum\limits_{j_{g_1}, j_{g_3},\ldots ,j_{g_{2r-1}}=0}^p
C_{j_k\ldots j_1}\biggl|_{j_{g_1}=j_{g_2},\ldots, j_{g_{2r-1}}=j_{g_{2r}}}-\Biggr.
$$

\vspace{2mm}
\begin{equation}
\label{july700001xyz}
\Biggl.-\frac{1}{2^r} \prod\limits_{l=1}^r {\bf 1}_{\{g_{2l}=g_{2l-1}+1\}}
C_{j_k \ldots j_1}\biggl|_{(j_{g_2} j_{g_1})\curvearrowright (\cdot)
\ldots (j_{g_{2r}} j_{g_{2r-1}})\curvearrowright (\cdot),
j_{g_{{}_{1}}}=~j_{g_{{}_{2}}},\ldots, j_{g_{{}_{2r-1}}}=~j_{g_{{}_{2r}}}
}\biggr.\Biggr)^2<\infty
\end{equation}

\vspace{5mm}
\noindent 
for all $r=1,2,\ldots,[k/2],$
where notations are the same as in 
(\ref{july90001}).
Then, by 
theorem on reducing 
of a limit to iterated one 
and (\ref{july90001})
we obtain

\vspace{-1mm}
$$
\lim\limits_{p,q\to\infty}
\sum\limits_{\stackrel{j_1,\ldots,j_m,\ldots,j_k=0}{{}_{m\ne g_1, g_2, \ldots, g_{2r-1},
g_{2r}}}}^q
\Biggl(\sum\limits_{j_{g_1}, j_{g_3},\ldots ,j_{g_{2r-1}}=0}^p
C_{j_k\ldots j_1}\biggl|_{j_{g_1}=j_{g_2},\ldots, j_{g_{2r-1}}=j_{g_{2r}}}-\Biggr.
$$

\vspace{2mm}
$$
\Biggl.-\frac{1}{2^r} \prod\limits_{l=1}^r {\bf 1}_{\{g_{2l}=g_{2l-1}+1\}}
C_{j_k \ldots j_1}\biggl|_{(j_{g_2} j_{g_1})\curvearrowright (\cdot)
\ldots (j_{g_{2r}} j_{g_{2r-1}})\curvearrowright (\cdot),
j_{g_{{}_{1}}}=~j_{g_{{}_{2}}},\ldots, j_{g_{{}_{2r-1}}}=~j_{g_{{}_{2r}}}
}\biggr.\Biggr)^2=
$$

\vspace{5mm}
$$
=\lim\limits_{q\to\infty}
\sum\limits_{\stackrel{j_1,\ldots,j_m,\ldots,j_k=0}{{}_{m\ne g_1, g_2, \ldots, g_{2r-1},
g_{2r}}}}^q
\lim\limits_{p\to\infty}\Biggl(\sum\limits_{j_{g_1}, j_{g_3},\ldots ,j_{g_{2r-1}}=0}^p
C_{j_k\ldots j_1}\biggl|_{j_{g_1}=j_{g_2},\ldots, j_{g_{2r-1}}=j_{g_{2r}}}-\Biggr.
$$

\vspace{2mm}
$$
\Biggl.-\frac{1}{2^r} \prod\limits_{l=1}^r {\bf 1}_{\{g_{2l}=g_{2l-1}+1\}}
C_{j_k \ldots j_1}\biggl|_{(j_{g_2} j_{g_1})\curvearrowright (\cdot)
\ldots (j_{g_{2r}} j_{g_{2r-1}})\curvearrowright (\cdot),
j_{g_{{}_{1}}}=~j_{g_{{}_{2}}},\ldots, j_{g_{{}_{2r-1}}}=~j_{g_{{}_{2r}}}
}\biggr.\Biggr)^2=0.
$$

\vspace{5mm}

Thus, we get
$$
\lim\limits_{p,q\to\infty}
\sum\limits_{\stackrel{j_1,\ldots,j_m,\ldots,j_k=0}{{}_{m\ne g_1, g_2, \ldots, g_{2r-1},
g_{2r}}}}^q
\Biggl(\sum\limits_{j_{g_1}, j_{g_3},\ldots ,j_{g_{2r-1}}=0}^p
C_{j_k\ldots j_1}\biggl|_{j_{g_1}=j_{g_2},\ldots, j_{g_{2r-1}}=j_{g_{2r}}}-\Biggr.
$$

\vspace{2mm}
\begin{equation}
\label{july700002xyz}
\Biggl.-\frac{1}{2^r} \prod\limits_{l=1}^r {\bf 1}_{\{g_{2l}=g_{2l-1}+1\}}
C_{j_k \ldots j_1}\biggl|_{(j_{g_2} j_{g_1})\curvearrowright (\cdot)
\ldots (j_{g_{2r}} j_{g_{2r-1}})\curvearrowright (\cdot),
j_{g_{{}_{1}}}=~j_{g_{{}_{2}}},\ldots, j_{g_{{}_{2r-1}}}=~j_{g_{{}_{2r}}}
}\biggr.\Biggr)^2=0.
\end{equation}

\vspace{5mm}
\noindent
Substituting $p=q$ in (\ref{july700002xyz}), we obtain 
(\ref{july90025}).

As a result, 
Hypotheses~7 and 8 are proved under 
the conditions formulated above in this section.

\vspace{5mm}

\section{Expansion of Iterated Stratonovich Stochastic Integrals
of Arbitrary Multiplicity $k$ $(k\in\mathbb{N})$. 
The Case of an Ar\-bit\-ra\-ry Complete Orthonormal System of 
Functions in $L_2([t,T]),$ $\psi_1(\tau),\ldots, \psi_k(\tau)
\in L_2([t,T]).$
Proof of Hypotheses~7, 8 for the Case $p_1=\ldots=p_k=p$
and Under the Condition (\ref{july700001xyz})}

\vspace{5mm}

This section is devoted to the following theorems.

\vspace{2mm}

{\bf Theorem~56}\ \cite{2018a}.\ {\it Suppose that the condition
{\rm (\ref{july700001xyz})}
is fulfilled$,$
$\{\phi_j(x)\}_{j=0}^{\infty}$
is an arbitrary complete orthonormal system of functions
in $L_2([t, T])$ and
$\psi_1(\tau),\ldots, \psi_k(\tau)\in L_2([t, T]).$
Then$,$ for the sum $\bar J^{*}[\psi^{(k)}]_{T,t}^{(i_1\ldots i_k)}$
of iterated Ito stochastic integrals 

\vspace{-1mm}
$$
\bar J^{*}[\psi^{(k)}]_{T,t}^{(i_1\ldots i_k)}=J[\psi^{(k)}]_{T,t}^{(i_1\ldots i_k)}+
\sum_{r=1}^{\left[k/2\right]}\frac{1}{2^r}
\sum_{(s_r,\ldots,s_1)\in {\rm A}_{k,r}}
J[\psi^{(k)}]_{T,t}^{s_r,\ldots,s_1}
$$

\vspace{3mm}
\noindent
the following 
expansion 
$$
\bar J^{*}[\psi^{(k)}]_{T,t}^{(i_1\ldots i_k)}=
\hbox{\vtop{\offinterlineskip\halign{
\hfil#\hfil\cr
{\rm l.i.m.}\cr
$\stackrel{}{{}_{p\to \infty}}$\cr
}} }
\sum_{j_1,\ldots,j_k=0}^{p}
C_{j_k \ldots j_1}\prod\limits_{l=1}^k \zeta_{j_l}^{(i_l)}
$$

\vspace{3mm}
\noindent
that converges in the mean-square sense is valid, where 

\vspace{-1mm}
\begin{equation}
\label{july99999}
C_{j_k \ldots j_1}=\int\limits_t^T\psi_k(t_k)\phi_{j_k}(t_k)\ldots
\int\limits_t^{t_2}
\psi_1(t_1)\phi_{j_1}(t_1)
dt_1\ldots dt_k
\end{equation}

\vspace{3mm}
\noindent
is the Fourier coefficient, 
${\rm l.i.m.}$ is a limit in the mean-square sense,
$i_1, \ldots, i_k=0, 1,\ldots,m,$

\vspace{-1mm}
$$
\zeta_{j}^{(i)}=
\int\limits_t^T \phi_{j}(\tau) d{\bf w}_{\tau}^{(i)}
$$ 

\vspace{2mm}
\noindent
are independent standard Gaussian random variables for various 
$i$ or $j$ {\rm (}in the case when $i\ne 0${\rm )},
${\bf w}_{\tau}^{(i)}={\bf f}_{\tau}^{(i)}$ 
for $i=1,\ldots,m$ and 
${\bf w}_{\tau}^{(0)}=\tau;$ another notations are the same as
in Theorem~{\rm 5}.}

\vspace{2mm}

Using Theorem~5, we obtain the following corollary of Theorem~56.

\vspace{2mm}
               
{\bf Theorem~57}\ \cite{2018a}.\ {\it Suppose that the condition
{\rm (\ref{july700001xyz})}
is fulfilled$,$
$\{\phi_j(x)\}_{j=0}^{\infty}$
is an arbitrary complete orthonormal system of functions
in $L_2([t, T])$ and
$\psi_1(\tau),\ldots, \psi_k(\tau)$ are continuous functions
at the interval $[t, T].$
Then$,$ for the iterated Stratonovich sto\-chas\-tic integral 
of arbitrary multiplicity $k$

\vspace{-1mm}
$$
J^{*}[\psi^{(k)}]_{T,t}^{(i_1\ldots i_k)}=
{\int\limits_t^{*}}^T
\psi_k(t_k) \ldots 
{\int\limits_t^{*}}^{t_{2}}
\psi_1(t_1) d{\bf w}_{t_1}^{(i_1)}\ldots
d{\bf w}_{t_k}^{(i_k)}
$$

\vspace{2.5mm}
\noindent
the following 
expansion 
\begin{equation}
\label{july300001}
J^{*}[\psi^{(k)}]_{T,t}^{(i_1\ldots i_k)}=
\hbox{\vtop{\offinterlineskip\halign{
\hfil#\hfil\cr
{\rm l.i.m.}\cr
$\stackrel{}{{}_{p\to \infty}}$\cr
}} }
\sum\limits_{j_1,\ldots,j_k=0}^{p}
C_{j_k \ldots j_1}\prod\limits_{l=1}^k \zeta_{j_l}^{(i_l)}
\end{equation}

\vspace{3mm}
\noindent
that converges in the mean-square sense is valid, where 

\vspace{-1mm}
$$
C_{j_k \ldots j_1}=\int\limits_t^T\psi_k(t_k)\phi_{j_k}(t_k)\ldots
\int\limits_t^{t_2}
\psi_1(t_1)\phi_{j_1}(t_1)
dt_1\ldots dt_k
$$

\vspace{3mm}
\noindent
is the Fourier coefficient, 
${\rm l.i.m.}$ is a limit in the mean-square sense,
$i_1, \ldots, i_k=0, 1,\ldots,m,$

\vspace{-1mm}
$$
\zeta_{j}^{(i)}=
\int\limits_t^T \phi_{j}(\tau) d{\bf w}_{\tau}^{(i)}
$$ 

\vspace{2mm}
\noindent
are independent standard Gaussian random variables for various 
$i$ or $j$ {\rm (}in the case when $i\ne 0${\rm )},
${\bf w}_{\tau}^{(i)}={\bf f}_{\tau}^{(i)}$ 
for $i=1,\ldots,m$ and 
${\bf w}_{\tau}^{(0)}=\tau.$}

\vspace{2mm}

{\bf Proof of Theorem~56.}\ According to the results
of Sect.~28, Theorem~56 will be proved if we prove that (see (\ref{july90001})):

\vspace{-1mm}
$$
\lim\limits_{p\to\infty}
\sum\limits_{j_{1},j_{3},\ldots,j_{2r-1}=0}^p
C_{j_k\ldots j_1}\biggl|_{j_{g_1}=j_{g_2},\ldots, j_{g_{2r-1}}=j_{g_{2r}}}=
$$

\vspace{3mm}
\begin{equation}
\label{july100000}
=\frac{1}{2^r} \prod\limits_{l=1}^r {\bf 1}_{\{g_{2l}=g_{2l-1}+1\}}
C_{j_k \ldots j_1}\biggl|_{(j_{g_2} j_{g_1})\curvearrowright (\cdot)
\ldots (j_{g_{2r}} j_{g_{2r-1}})\curvearrowright (\cdot),
j_{g_{{}_{1}}}=~j_{g_{{}_{2}}},\ldots, j_{g_{{}_{2r-1}}}=~j_{g_{{}_{2r}}}
}\biggr.
\end{equation}

\vspace{4mm}
\noindent
for all possible $g_1,g_2,\ldots,g_{2r-1},g_{2r}$ (see {\rm (\ref{leto5007})),
where $k\ge 2r,$ $r=1,2,\ldots,[k/2],$ $C_{j_k\ldots j_1}$ is defined by (\ref{july99999}),
another notations are the same as in Theorem~47.

Moreover (assuming that (\ref{july100000}) is proved), 
the series on the left-hand side of 
(\ref{july100000}) converges absolutly (the case $k=2r$)
and converges absolutly
for any fixed $j_1,\ldots,j_q,\ldots,j_k$
and $q\ne g_1, g_2, \ldots,$ $g_{2r-1},$ $g_{2r}$
(the case $k>2r$)
since
its sum does not depend 
on permutations of basis func\-ti\-ons
(here the basis in $L_2([t,T]^{r})$ is
$\left\{\phi_{j_1}(x_1)\ldots \phi_{j_r}(x_r)\right\}_{j_1,\ldots,j_r=0}^{\infty}$).
Recall that any permutation of basis functions in a Hilbert space forms a basis 
in this Hilbert space \cite{gohb}.

Also recall that the case
$k=2r$ of (\ref{july100000}) is considered in Sect.~26.
Consider the case $k>2r.$
Using Fubini's Theorem, we obtain

\vspace{-1mm}

$$
\int\limits_t^T h_{k}(t_k)\ldots \int\limits_t^{t_{l+2}} h_{l+1}(t_{l+1})
\int\limits_t^{t_{l+1}} h_{l}(t_{l})
\int\limits_t^{t_{l}} h_{l-1}(t_{l-1})\ldots
\int\limits_t^{t_2} h_{1}(t_1)
dt_1\ldots 
dt_{l-1}dt_{l}dt_{l+1}\ldots dt_k=
$$

\vspace{1mm}
$$
=\int\limits_t^T h_{k}(t_k)\ldots \int\limits_t^{t_{l+2}} h_{l+1}(t_{l+1})
\int\limits_t^{t_{l+1}} h_{1}(t_{1})
\int\limits_{t_1}^{t_{l+1}} h_{2}(t_{2})\ldots
\int\limits_{t_{l-2}}^{t_{l+1}} h_{l-1}(t_{l-1})
\int\limits_{t_{l-1}}^{t_{l+1}} h_{l}(t_{l})dt_l\times
$$

\vspace{1mm}
$$
\times dt_{l-1}\ldots dt_2dt_{1}dt_{l+1}\ldots dt_k=
$$

\vspace{1mm}
$$
=\int\limits_t^T h_{k}(t_k)\ldots \int\limits_t^{t_{l+2}} h_{l+1}(t_{l+1})
\left(\int\limits_{t}^{t_{l+1}} h_{l}(t_{l})dt_l\right)\int\limits_t^{t_{l+1}} h_{1}(t_{1})
\int\limits_{t_1}^{t_{l+1}} h_{2}(t_{2})\ldots
\int\limits_{t_{l-2}}^{t_{l+1}} h_{l-1}(t_{l-1})
\times
$$

\vspace{1mm}
$$
\times dt_{l-1}\ldots dt_2dt_{1}dt_{l+1}\ldots dt_k-
$$

\vspace{1mm}
$$
-\int\limits_t^T h_{k}(t_k)\ldots \int\limits_t^{t_{l+2}} h_{l+1}(t_{l+1})
\int\limits_t^{t_{l+1}} h_{1}(t_{1})
\int\limits_{t_1}^{t_{l+1}} h_{2}(t_{2})\ldots
\int\limits_{t_{l-2}}^{t_{l+1}} h_{l-1}(t_{l-1})
\left(\int\limits_{t}^{t_{l-1}} h_{l}(t_{l})dt_l\right)\times
$$

\vspace{1mm}
$$
\times dt_{l-1}\ldots dt_2dt_{1}dt_{l+1}\ldots dt_k=
$$

\vspace{1mm}
$$
=\int\limits_t^T h_{k}(t_k)\ldots \int\limits_t^{t_{l+2}} h_{l+1}(t_{l+1})
\left(\int\limits_t^{t_{l+1}} h_{l}(t_{l})dt_l\right)
\int\limits_t^{t_{l+1}} h_{l-1}(t_{l-1})\ldots
$$

\vspace{1mm}
$$
\ldots 
\int\limits_t^{t_2} h_{1}(t_1)
dt_1\ldots dt_{l-1}dt_{l+1}\ldots dt_k-
$$

\vspace{1mm}
$$
-\int\limits_t^T h_{k}(t_k)\ldots \int\limits_t^{t_{l+2}} h_{l+1}(t_{l+1})
\int\limits_t^{t_{l+1}} h_{l-1}(t_{l-1})\left(\int\limits_t^{t_{l-1}} h_{l}(t_{l})dt_l\right)
\int\limits_t^{t_{l-1}} h_{l-2}(t_{l-2})
\ldots
$$

\vspace{1mm}
\begin{equation}
\label{july100003}
\ldots
\int\limits_t^{t_2} h_{1}(t_1)
dt_1\ldots dt_{l-2}dt_{l-1}dt_{l+1}\ldots dt_k,
\end{equation}

\vspace{4mm}
\noindent
where $2<l<k-1$ and $h_1(\tau),\ldots,h_k(\tau)\in L_2([t, T]).$ 

By analogy with (\ref{july100003}) we have for $l=k$

$$
\int\limits_t^{T} h_{l}(t_{l})
\int\limits_t^{t_{l}} h_{l-1}(t_{l-1})\ldots
\int\limits_t^{t_2} h_{1}(t_1)
dt_1\ldots 
dt_{l-1}dt_{l}=
$$

\vspace{1mm}
$$
=\int\limits_{t}^{T} h_{1}(t_{1})
\int\limits_{t_1}^{T} h_{2}(t_{2})\ldots
\int\limits_{t_{l-2}}^{T} h_{l-1}(t_{l-1})\int\limits_{t_{l-1}}^{T} h_{l}(t_{l})
dt_ldt_{l-1}\ldots dt_2dt_{1}=
$$

\vspace{1mm}
$$
=\left(\int\limits_{t}^{T} h_{l}(t_{l})
dt_l\right)\int\limits_{t}^{T} h_{1}(t_{1})
\int\limits_{t_1}^{T} h_{2}(t_{2})\ldots
\int\limits_{t_{l-2}}^{T} h_{l-1}(t_{l-1})
dt_{l-1}\ldots dt_2dt_{1}-
$$

\vspace{1mm}
$$
-\int\limits_{t}^{T}h_{1}(t_{1})
\int\limits_{t_1}^{T} h_{2}(t_{2})\ldots
\int\limits_{t_{l-2}}^{T} h_{l-1}(t_{l-1})\left(
\int\limits_{t}^{t_{l-1}} h_{l}(t_{l})dt_l\right)
dt_{l-1}\ldots dt_2dt_{1}=
$$

\vspace{1mm}
$$
=\left(\int\limits_{t}^{T} h_{l}(t_{l})
dt_l\right)\int\limits_{t}^{T} h_{l-1}(t_{l-1})
\ldots
\int\limits_{t}^{t_2} h_{1}(t_{1})
dt_{1}\ldots dt_{l-1}-
$$

\vspace{1mm}
\begin{equation}
\label{july100004}
-\int\limits_{t}^{T}
h_{l-1}(t_{l-1})\left(
\int\limits_{t}^{t_{l-1}} h_{l}(t_{l})dt_l\right)
\int\limits_{t}^{t_{l-1}} h_{l-2}(t_{l-2})\ldots
\int\limits_{t}^{t_2} 
h_{1}(t_{1})
dt_{1}\ldots dt_{l-1}.
\end{equation}

\vspace{3mm}

We will assume that for $l=1$ the transformation 
(\ref{july100003}) is not carried out since

\vspace{-1mm}
$$
\int\limits_t^{t_2} h_{1}(t_1)
dt_1
$$

\vspace{3mm}
\noindent
is the innermost integral on the left-hand side of (\ref{july100003}).
The formulas (\ref{july100003}), (\ref{july100004})
will be used further.

Let us carry out the transformations 
(\ref{july100003}), (\ref{july100004}) 
for 

\vspace{-1mm}
$$
C_{j_k\ldots j_1}\biggl|_{j_{g_1}=j_{g_2},\ldots, j_{g_{2r-1}}=j_{g_{2r}}}
$$

\vspace{3mm}
\noindent
iteratively for $j_1,\ldots,j_q,\ldots,j_k$
$(q\ne g_1, g_2, \ldots, g_{2r-1},$ $g_{2r}).$
As a result, we obtain 

$$
C_{j_k\ldots j_1}\biggl|_{j_{g_1}=j_{g_2},\ldots, j_{g_{2r-1}}=j_{g_{2r}}}=
$$

\begin{equation}
\label{july100005}
=\sum\limits_{d=1}^{2^{k-2r}}(-1)^{d-1}
\left(\hat C^{(d)}_{j_k\ldots j_1}\biggl|_{j_{g_1}=j_{g_2},\ldots, j_{g_{2r-1}}=j_{g_{2r}}}-~
\bar C^{(d)}_{j_k\ldots j_1}\biggl|_{j_{g_1}=j_{g_2},\ldots, j_{g_{2r-1}}=j_{g_{2r}}}\right),
\end{equation}

\vspace{5mm}
\noindent
where some terms in the sum
$$
\sum\limits_{d=1}^{2^{k-2r}}
$$

\vspace{3mm}
\noindent
can be identically equal to zero due to
the remark to (\ref{july100003}), (\ref{july100004}).

Using
(\ref{july100005}), we obtain  

\vspace{-1mm}
$$
\lim\limits_{p\to\infty}
\sum\limits_{j_{g_1}, j_{g_3},\ldots,j_{g_{2r-1}}=0}^p
C_{j_k\ldots j_1}\biggl|_{j_{g_1}=j_{g_2},\ldots, j_{g_{2r-1}}=j_{g_{2r}}}=
$$

\vspace{4mm}
$$
=\lim\limits_{p\to\infty}
\sum\limits_{j_{g_1}, j_{g_3},\ldots,j_{g_{2r-1}}=0}^p
\sum\limits_{d=1}^{2^{k-2r}}(-1)^{d-1}
\left(\hat C^{(d)}_{j_k\ldots j_1}\biggl|_{j_{g_1}=j_{g_2},\ldots, j_{g_{2r-1}}=j_{g_{2r}}}-
\right.
$$

\vspace{2mm}
$$
\left.
-~\bar C^{(d)}_{j_k\ldots j_1}\biggl|_{j_{g_1}=j_{g_2},\ldots, j_{g_{2r-1}}=j_{g_{2r}}}\right)=
$$

\vspace{4mm}
$$
=\sum\limits_{d=1}^{2^{k-2r}}(-1)^{d-1} \lim\limits_{p\to\infty}
\sum\limits_{j_{g_1}, j_{g_3},\ldots,j_{g_{2r-1}}=0}^p
\left(\hat C^{(d)}_{j_k\ldots j_1}\biggl|_{j_{g_1}=j_{g_2},\ldots, j_{g_{2r-1}}=j_{g_{2r}}}-
\right.
$$

\vspace{2mm}
\begin{equation}
\label{july100010}
\left.
-~\bar C^{(d)}_{j_k\ldots j_1}\biggl|_{j_{g_1}=j_{g_2},\ldots, j_{g_{2r-1}}=j_{g_{2r}}}\right).
\end{equation}

\vspace{3mm}

Further, consider 3 possible cases.

\vspace{2mm}

{\bf Case~1.}\ The quantities

\vspace{-1mm}
\begin{equation}
\label{july100010x}
\hat C^{(d)}_{j_k\ldots j_1}\biggl|_{j_{g_1}=j_{g_2},\ldots, j_{g_{2r-1}}=j_{g_{2r}}},\ \ \ 
\bar C^{(d)}_{j_k\ldots j_1}\biggl|_{j_{g_1}=j_{g_2},\ldots, j_{g_{2r-1}}=j_{g_{2r}}}
\end{equation}

\vspace{3mm}
\noindent
are such that
\begin{equation}
\label{july100010y}
\prod\limits_{l=1}^r {\bf 1}_{\{g_{2l}=g_{2l-1}+1\}}=1
\end{equation}

\vspace{3mm}
\noindent
for $d=1,2,\ldots, 2^{k-2r}$ and 
\begin{equation}
\label{july100010z}
C_{j_k\ldots j_1}\biggl|_{j_{g_1}=j_{g_2},\ldots, j_{g_{2r-1}}=j_{g_{2r}}}
\end{equation}

\vspace{5mm}
\noindent
is such that the condition (\ref{july100010y}) is fulfilled for
(\ref{july100010z}).

\vspace{2mm}

{\bf Case~2.}\ The quantities (\ref{july100010x}) 
are such that the condition (\ref{july100010y})
is satisfied for $d=1,2,\ldots, 2^{k-2r}$ and 
(\ref{july100010z}) is such that the condition 

\vspace{-1mm}
\begin{equation}
\label{july100010v}
\prod\limits_{l=1}^r {\bf 1}_{\{g_{2l}=g_{2l-1}+1\}}=0
\end{equation}

\vspace{2mm}
\noindent
is fulfilled for (\ref{july100010z}).

\vspace{2mm}

{\bf Case~3.}\ The quantities (\ref{july100010x}) 
are such that the condition (\ref{july100010v}) 
is satisfied for $d=1,2,\ldots, 2^{k-2r}$ and 
(\ref{july100010z}) is such that the condition
(\ref{july100010v}) 
is fulfilled for (\ref{july100010z}).

\vspace{2mm}

For Case~1, applying 
(\ref{july100000}) for the case $k=2r$ and (\ref{july100010}),
we get
for any fixed $j_1,\ldots,j_q,\ldots,j_k$
$(q\ne g_1, g_2, \ldots, g_{2r-1},$ $g_{2r})$

\vspace{1mm}
$$
\lim\limits_{p\to\infty}
\sum\limits_{j_{g_1}, j_{g_3},\ldots,j_{g_{2r-1}}=0}^p
C_{j_k\ldots j_1}\biggl|_{j_{g_1}=j_{g_2},\ldots, j_{g_{2r-1}}=j_{g_{2r}}}=
$$

\vspace{4mm}
$$
=\sum\limits_{d=1}^{2^{k-2r}}(-1)^{d-1} \lim\limits_{p\to\infty}
\sum\limits_{j_{g_1}, j_{g_3},\ldots,j_{g_{2r-1}}=0}^p
\left(\hat C^{(d)}_{j_k\ldots j_1}\biggl|_{j_{g_1}=j_{g_2},\ldots, j_{g_{2r-1}}=j_{g_{2r}}}-
\right.
$$

\vspace{2mm}
$$
\left.
-~\bar C^{(d)}_{j_k\ldots j_1}\biggl|_{j_{g_1}=j_{g_2},\ldots, j_{g_{2r-1}}=j_{g_{2r}}}\right)=
$$

\vspace{4mm}
$$
=
\sum\limits_{d=1}^{2^{k-2r}}(-1)^{d-1}\frac{1}{2^r} 
\prod\limits_{l=1}^r {\bf 1}_{\{g_{2l}=g_{2l-1}+1\}}\times
$$

\vspace{1mm}
$$
\times
\left(\hat C^{(d)}_{j_k \ldots j_1}\biggl|_{(j_{g_2} j_{g_1})\curvearrowright (\cdot)
\ldots (j_{g_{2r}} j_{g_{2r-1}})\curvearrowright (\cdot),
j_{g_{{}_{1}}}=~j_{g_{{}_{2}}},\ldots, j_{g_{{}_{2r-1}}}=~j_{g_{{}_{2r}}}
}\biggr.-\right.
$$

\vspace{1mm}
\begin{equation}
\label{july100006abc}
\left.-~\bar C^{(d)}_{j_k \ldots j_1}\biggl|_{(j_{g_2} j_{g_1})\curvearrowright (\cdot)
\ldots (j_{g_{2r}} j_{g_{2r-1}})\curvearrowright (\cdot),
j_{g_{{}_{1}}}=~j_{g_{{}_{2}}},\ldots, j_{g_{{}_{2r-1}}}=~j_{g_{{}_{2r}}}
}\biggr.
\right)=
\end{equation}

\vspace{2mm}
$$
=
\sum\limits_{d=1}^{2^{k-2r}}(-1)^{d-1}\frac{1}{2^r} 
\left(\hat C^{(d)}_{j_k \ldots j_1}\biggl|_{(j_{g_2} j_{g_1})\curvearrowright (\cdot)
\ldots (j_{g_{2r}} j_{g_{2r-1}})\curvearrowright (\cdot),
j_{g_{{}_{1}}}=~j_{g_{{}_{2}}},\ldots, j_{g_{{}_{2r-1}}}=~j_{g_{{}_{2r}}}
}\biggr.-\right.
$$

\vspace{1mm}
\begin{equation}
\label{july100006}
\left.-~\bar C^{(d)}_{j_k \ldots j_1}\biggl|_{(j_{g_2} j_{g_1})\curvearrowright (\cdot)
\ldots (j_{g_{2r}} j_{g_{2r-1}})\curvearrowright (\cdot),
j_{g_{{}_{1}}}=~j_{g_{{}_{2}}},\ldots, j_{g_{{}_{2r-1}}}=~j_{g_{{}_{2r}}}
}\biggr.
\right),
\end{equation}

\vspace{5mm}
\noindent
where $g_1,g_2,\ldots,g_{2r-1},g_{2r}$ as in (\ref{leto5007xxx}),
$k>2r,$ $r=1,2,\ldots,[k/2].$

It is not difficult to see that 
the left-hand side of (\ref{july100010y}) is a constant for the
quantities (\ref{july100010x}) for all $d=1,2,\ldots, 2^{k-2r}$.

Using (\ref{july100003}), (\ref{july100004}), we obtain

$$
\sum\limits_{d=1}^{2^{k-2r}}(-1)^{d-1}\frac{1}{2^r} 
\left(\hat C^{(d)}_{j_k \ldots j_1}\biggl|_{(j_{g_2} j_{g_1})\curvearrowright (\cdot)
\ldots (j_{g_{2r}} j_{g_{2r-1}})\curvearrowright (\cdot),
j_{g_{{}_{1}}}=~j_{g_{{}_{2}}},\ldots, j_{g_{{}_{2r-1}}}=~j_{g_{{}_{2r}}}
}\biggr.-\right.
$$

\vspace{1mm}
$$
\left.-~\bar C^{(d)}_{j_k \ldots j_1}\biggl|_{(j_{g_2} j_{g_1})\curvearrowright (\cdot)
\ldots (j_{g_{2r}} j_{g_{2r-1}})\curvearrowright (\cdot),
j_{g_{{}_{1}}}=~j_{g_{{}_{2}}},\ldots, j_{g_{{}_{2r-1}}}=~j_{g_{{}_{2r}}}
}\biggr.
\right)=
$$

\vspace{2mm}
\begin{equation}
\label{july100007}
=\frac{1}{2^r} C_{j_k \ldots j_1}\biggl|_{(j_{g_2} j_{g_1})\curvearrowright (\cdot)
\ldots (j_{g_{2r}} j_{g_{2r-1}})\curvearrowright (\cdot),
j_{g_{{}_{1}}}=~j_{g_{{}_{2}}},\ldots, j_{g_{{}_{2r-1}}}=~j_{g_{{}_{2r}}}
}\biggr..
\end{equation}

\vspace{6mm}

Combining (\ref{july100006}) and (\ref{july100007}), we have
for any fixed $j_1,\ldots,j_q,\ldots,j_k$
$(q\ne g_1, g_2, \ldots, g_{2r-1},$ $g_{2r})$

$$
\lim\limits_{p\to\infty}
\sum\limits_{j_{g_1}, j_{g_3},\ldots,j_{g_{2r-1}}=0}^p
C_{j_k\ldots j_1}\biggl|_{j_{g_1}=j_{g_2},\ldots, j_{g_{2r-1}}=j_{g_{2r}}}=
$$

\vspace{2mm}
\begin{equation}
\label{july100008}
=\frac{1}{2^r}
C_{j_k \ldots j_1}\biggl|_{(j_{g_2} j_{g_1})\curvearrowright (\cdot)
\ldots (j_{g_{2r}} j_{g_{2r-1}})\curvearrowright (\cdot),
j_{g_{{}_{1}}}=~j_{g_{{}_{2}}},\ldots, j_{g_{{}_{2r-1}}}=~j_{g_{{}_{2r}}}
}\biggr.,
\end{equation}

\vspace{5mm}
\noindent
where $g_1,g_2,\ldots,g_{2r-1},g_{2r}$ as in (\ref{leto5007xxx}),
$k>2r,$ $r=1,2,\ldots,[k/2].$

From (\ref{july100000}) for the case $k=2r$ and (\ref{july100008}) ($k>2r$) we obtain 
(\ref{july100000}) for the case $k\ge 2r$. The equality
(\ref{july100000}) is proved for Case~1.

For Case~2, applying 
(\ref{july100000}) for the case $k=2r$ and (\ref{july100010}),
we get (\ref{july100006})
for any fixed $j_1,\ldots,j_q,\ldots,j_k$
$(q\ne g_1, g_2, \ldots, g_{2r-1},$ $g_{2r}).$
Further, note that

\vspace{-1mm}
$$
\hat C^{(d)}_{j_k \ldots j_1}\biggl|_{(j_{g_2} j_{g_1})\curvearrowright (\cdot)
\ldots (j_{g_{2r}} j_{g_{2r-1}})\curvearrowright (\cdot),
j_{g_{{}_{1}}}=~j_{g_{{}_{2}}},\ldots, j_{g_{{}_{2r-1}}}=~j_{g_{{}_{2r}}}
}\biggr.=
$$

\vspace{2mm}
\begin{equation}
\label{july100008xyz}
=~\bar C^{(d)}_{j_k \ldots j_1}\biggl|_{(j_{g_2} j_{g_1})\curvearrowright (\cdot)
\ldots (j_{g_{2r}} j_{g_{2r-1}})\curvearrowright (\cdot),
j_{g_{{}_{1}}}=~j_{g_{{}_{2}}},\ldots, j_{g_{{}_{2r-1}}}=~j_{g_{{}_{2r}}}
}\biggr.
\end{equation}

\vspace{5mm}
\noindent
for Case~2. Combining (\ref{july100006}) and (\ref{july100008xyz}), we
obtain (Case~2)
for any fixed $j_1,\ldots,j_q,\ldots,j_k$
$(q\ne g_1, g_2, \ldots, g_{2r-1},$ $g_{2r})$

\begin{equation}
\label{july100008xyz1}
\lim\limits_{p\to\infty}
\sum\limits_{j_{g_1}, j_{g_3},\ldots,j_{g_{2r-1}}=0}^p
C_{j_k\ldots j_1}\biggl|_{j_{g_1}=j_{g_2},\ldots, j_{g_{2r-1}}=j_{g_{2r}}}=0.
\end{equation}

\vspace{5mm}

From (\ref{july100000}) for the case $k=2r$ and (\ref{july100008xyz1}) $(k>2r)$ we obtain 
(\ref{july100008xyz1}) for the case $k\ge 2r$. The equality
(\ref{july100000}) is proved for Case~2.

For Case~3, applying 
(\ref{july100000}) for the case $k=2r$ and (\ref{july100010}),
we get (\ref{july100006abc})
for any fixed $j_1,\ldots,j_q,\ldots,j_k$
$(q\ne g_1, g_2, \ldots, g_{2r-1},$ $g_{2r}).$
Since 
\begin{equation}
\label{july400000}
\prod\limits_{l=1}^r {\bf 1}_{\{g_{2l}=g_{2l-1}+1\}}=0
\end{equation}

\vspace{3mm}
\noindent
for Case~3, then from (\ref{july100006abc})
we get (\ref{july100008xyz1}) for $k>2r$
(recall that the left-hand side of (\ref{july400000}) is a constant for the
quantities (\ref{july100010x}) for all $d=1,2,\ldots, 2^{k-2r}$).
From (\ref{july100000}) for $k=2r$ and (\ref{july100008xyz1}) for $k>2r$ (Case~3) we obtain 
(\ref{july100008xyz1}) for $k\ge 2r$ (Case~3). The equality
(\ref{july100000}) is proved for Case~3. 
Theorem~56 is proved. Theorem~57 is also proved.

In conclusion of this section, we will make a remark
about the condition (\ref{july90025}). It would seem
that according to (\ref{july90001}), we can write

$$
\lim\limits_{p\to\infty}
\sum\limits_{\stackrel{j_1,\ldots,j_q,\ldots,j_k=0}{{}_{q\ne g_1, g_2, \ldots, g_{2r-1},
g_{2r}}}}^p
\Biggl(\sum\limits_{j_{g_1}, j_{g_3},\ldots ,j_{g_{2r-1}}=0}^p
C_{j_k\ldots j_1}\biggl|_{j_{g_1}=j_{g_2},\ldots, j_{g_{2r-1}}=j_{g_{2r}}}-\Biggr.
$$

\vspace{2mm}
$$
\Biggl.-\frac{1}{2^r} \prod\limits_{l=1}^r {\bf 1}_{\{g_{2l}=g_{2l-1}+1\}}
C_{j_k \ldots j_1}\biggl|_{(j_{g_2} j_{g_1})\curvearrowright (\cdot)
\ldots (j_{g_{2r}} j_{g_{2r-1}})\curvearrowright (\cdot),
j_{g_{{}_{1}}}=~j_{g_{{}_{2}}},\ldots, j_{g_{{}_{2r-1}}}=~j_{g_{{}_{2r}}}
}\biggr.\Biggr)^2=
$$

\vspace{6mm}
$$
=\sum\limits_{\stackrel{j_1,\ldots,j_q,\ldots,j_k=0}{{}_{q\ne g_1, g_2, \ldots, g_{2r-1},
g_{2r}}}}^{\infty}
\Biggl(\sum\limits_{j_{g_1}, j_{g_3},\ldots ,j_{g_{2r-1}}=0}^{\infty}
C_{j_k\ldots j_1}\biggl|_{j_{g_1}=j_{g_2},\ldots, j_{g_{2r-1}}=j_{g_{2r}}}-\Biggr.
$$

\vspace{2mm}
$$
\Biggl.-\frac{1}{2^r} \prod\limits_{l=1}^r {\bf 1}_{\{g_{2l}=g_{2l-1}+1\}}
C_{j_k \ldots j_1}\biggl|_{(j_{g_2} j_{g_1})\curvearrowright (\cdot)
\ldots (j_{g_{2r}} j_{g_{2r-1}})\curvearrowright (\cdot),
j_{g_{{}_{1}}}=~j_{g_{{}_{2}}},\ldots, j_{g_{{}_{2r-1}}}=~j_{g_{{}_{2r}}}
}\biggr.\Biggr)^2=
$$

\vspace{4.5mm}
$$
=\sum\limits_{\stackrel{j_1,\ldots,j_q,\ldots,j_k=0}{{}_{q\ne g_1, g_2, \ldots, g_{2r-1},
g_{2r}}}}^{\infty}
\Bigl( 0 \Bigr)^2= 0
$$

\vspace{4mm}
\noindent 
for all $r=1,2,\ldots,[k/2]$ and 
for all possible $g_1,g_2,\ldots,g_{2r-1},g_{2r}$ (see {\rm (\ref{leto5007})),
where notations are the same as in 
(\ref{july90025}).

However, the above argument contains an error associated
with the replacement of the limit with an iterated one.
Let us consider this observation in more detail
using an example.

To begin, let us recall that the sum of an infinite number series
is defined as the limit of the partial sums of this series, i.e.

\vspace{-2mm}
$$
\lim\limits_{n\to\infty}
\sum\limits_{i=1}^{n}a_i
\stackrel{\sf def}{=}
\sum\limits_{i=1}^{\infty}a_i.
$$

\vspace{4mm}

Let $k=3,$ $r=1$ and $g_1=1, g_2=3.$ Further, we have

\vspace{-1mm}
$$
\lim\limits_{p\to\infty}\sum\limits_{j_2=0}^p
\Biggl(
\sum\limits_{j_1=0}^{p}C_{j_1 j_2 j_1}\Biggr)^2=
\lim\limits_{p\to\infty}\sum\limits_{j_2=0}^p
\sum\limits_{j_1=0}^{p}C_{j_1 j_2 j_1}\sum\limits_{j_3=0}^{p}C_{j_3 j_2 j_3}=
$$

\vspace{2mm}

\begin{equation}
\label{2025july2025}
=
\lim\limits_{p\to\infty}\sum\limits_{j_3, j_2, j_1=0}^p
C_{j_1 j_2 j_1}C_{j_3 j_2 j_3}
\stackrel{\sf def}{=}\sum\limits_{j_3, j_2, j_1=0}^{\infty}
C_{j_1 j_2 j_1}C_{j_3 j_2 j_3},
\end{equation}

\vspace{4mm}
$$
\lim\limits_{q\to\infty}\lim\limits_{p\to\infty}\sum\limits_{j_2=0}^q
\Biggl(
\sum\limits_{j_1=0}^{p}C_{j_1 j_2 j_1}\Biggr)^2=
\lim\limits_{q\to\infty}\sum\limits_{j_2=0}^q
\lim\limits_{p\to\infty}\Biggl(
\sum\limits_{j_1=0}^{p}C_{j_1 j_2 j_1}\Biggr)^2=
\lim\limits_{q\to\infty}\sum\limits_{j_2=0}^q
\Biggl(
\sum\limits_{j_1=0}^{\infty}C_{j_1 j_2 j_1}\Biggr)^2=
$$

\vspace{2mm}

\begin{equation}
\label{2025july2025a}
=
\sum\limits_{j_2=0}^{\infty}
\Biggl(
\sum\limits_{j_1=0}^{\infty}C_{j_1 j_2 j_1}\Biggr)^2=
\sum\limits_{j_2=0}^{\infty}
\sum\limits_{j_1=0}^{\infty}C_{j_1 j_2 j_1}
\sum\limits_{j_3=0}^{\infty}C_{j_3 j_2 j_3}=
\sum\limits_{j_2=0}^{\infty}
\sum\limits_{j_1=0}^{\infty}
\sum\limits_{j_3=0}^{\infty}C_{j_1 j_2 j_1}C_{j_3 j_2 j_3}.
\end{equation}

\vspace{5mm}

It is obvious that the right-hand sides of equalities
(\ref{2025july2025}) and (\ref{2025july2025a})
are generally not equal. The equality 
of the mentioned expressions requires separate proof.

In the next section, we will consider a fairly
efficient approach to proving the equality 
(\ref{july90025}).

\vspace{5mm}

\section{Expansion of Iterated Stratonovich Stochastic Integrals
of Arbitrary Multiplicity $k$ $(k\in\mathbb{N})$. 
The Case of an Ar\-bit\-ra\-ry Complete Orthonormal System of 
Functions in $L_2([t,T]),$ $\psi_1(\tau),\ldots, \psi_k(\tau)
\in L_2([t,T]).$
Proof of Hypotheses~7, 8 for the Case $p_1=\ldots=p_k=p$
Under the Condition (\ref{09091})}

\vspace{5mm}

We will start this section with an example. Let us assume that
$h_1(\tau),\ldots,$ $h_{12}(\tau)\in L_2([t, T])$ and consider the
following integral

$$
I\stackrel{\sf def}{=}\int\limits_t^T h_{12}(t_{12})
\int\limits_t^{t_{12}}h_{11}(t_{11})\ldots 
\int\limits_t^{t_{2}}h_{1}(t_{1})dt_1 \ldots dt_{11} dt_{12}.
$$

\vspace{3mm}

We want to transform the integral $I$ in such a way that

$$
I=\int\limits_t^T h_{10}(t_{10})
\int\limits_t^{t_{10}}h_{6}(t_{6})
\int\limits_t^{t_{6}}h_{4}(t_{4})
\int\limits_t^{t_{4}}h_{3}(t_{3})
\left(\ldots\right)
dt_3 dt_{4} dt_6 dt_{10},
$$

\vspace{3mm}
\noindent
where $\left(\ldots\right)$ is some expression.

Using Fubini's Theorem, we obtain

\vspace{-1mm}
$$
I=\int\limits_t^T h_{12}(t_{12})
\int\limits_t^{t_{12}}h_{11}(t_{11})
\int\limits_t^{t_{11}}\underline{h_{10}(t_{10})}
\int\limits_t^{t_{10}}h_{9}(t_{9})
\int\limits_t^{t_{9}}h_{8}(t_{8})
\int\limits_t^{t_{8}}h_{7}(t_{7})
\int\limits_t^{t_{7}}\underline{h_{6}(t_{6})}\times
$$

\vspace{2mm}
$$
\times
\int\limits_t^{t_{6}}h_{5}(t_{5})
\int\limits_t^{t_{5}}\underline{h_{4}(t_{4})}
\int\limits_t^{t_{4}}\underline{h_{3}(t_{3})}
\int\limits_t^{t_{3}}h_{2}(t_{2})
\int\limits_t^{t_{2}}h_{1}(t_{1})dt_1 dt_2 dt_3 dt_4
dt_5 dt_6 dt_7 dt_8 
dt_9 dt_{10} dt_{11} dt_{12}=
$$

\vspace{2mm}
$$
=\int\limits_t^T
\underline{h_{10}(t_{10})}
\int\limits_t^{t_{10}}h_{9}(t_{9})
\int\limits_t^{t_{9}}h_{8}(t_{8})
\int\limits_t^{t_{8}}h_{7}(t_{7})
\int\limits_t^{t_{7}}\underline{h_{6}(t_{6})}\int\limits_t^{t_{6}}h_{5}(t_{5})
\times
$$

\vspace{2mm}
$$
\times
\int\limits_t^{t_{5}}\underline{h_{4}(t_{4})}
\int\limits_t^{t_{4}}\underline{h_{3}(t_{3})}
\int\limits_t^{t_{3}}h_{2}(t_{2})
\int\limits_t^{t_{2}}h_{1}(t_{1})dt_1 dt_2 dt_3 dt_4
dt_5 dt_6 dt_7 dt_8 dt_9\times
$$

\vspace{2mm}
$$
\times
\left(\int\limits_{t_{10}}^T h_{11}(t_{11})
\int\limits_{t_{11}}^T h_{12}(t_{12})
dt_{12} dt_{11}\right) dt_{10}=
$$

\vspace{2mm}
$$
=\int\limits_t^T
\underline{h_{10}(t_{10})}
\int\limits_t^{t_{10}}\underline{h_{6}(t_{6})}\int\limits_t^{t_{6}}h_{5}(t_{5})
\int\limits_t^{t_{5}}\underline{h_{4}(t_{4})}
\int\limits_t^{t_{4}}\underline{h_{3}(t_{3})}
\int\limits_t^{t_{3}}h_{2}(t_{2})
\int\limits_t^{t_{2}}h_{1}(t_{1})\times
$$

\vspace{2mm}
$$
\times
dt_1 dt_2 dt_3 dt_4 dt_5 \left(\int\limits_{t_6}^{t_{10}}h_7(t_7)
\int\limits_{t_7}^{t_{10}}h_8(t_8)
\int\limits_{t_8}^{t_{10}}h_9(t_9)dt_9 dt_8 dt_7
\right) dt_6 \times
$$

\vspace{2mm}
$$
\times
\left(\int\limits_{t_{10}}^T h_{11}(t_{11})
\int\limits_{t_{11}}^T h_{12}(t_{12})
dt_{12} dt_{11}\right) dt_{10}=
$$

\vspace{2mm}
$$
=\int\limits_t^T
\underline{h_{10}(t_{10})}
\int\limits_t^{t_{10}}\underline{h_{6}(t_{6})}\int\limits_t^{t_{6}}
\underline{h_{4}(t_{4})}
\int\limits_t^{t_{4}}\underline{h_{3}(t_{3})}
\left(\int\limits_t^{t_{3}}h_{2}(t_{2})
\int\limits_t^{t_{2}}h_{1}(t_{1})
dt_1 dt_2\right) dt_3 \times
$$

\vspace{2mm}
$$
\times \left(\int\limits_{t_4}^{t_6} h_5(t_5)dt_5 \right)dt_4 
\left(\int\limits_{t_6}^{t_{10}}h_7(t_7)
\int\limits_{t_7}^{t_{10}}h_8(t_8)
\int\limits_{t_8}^{t_{10}}h_9(t_9)dt_9 dt_8 dt_7
\right) dt_6 \times
$$

\vspace{2mm}
$$
\times
\left(\int\limits_{t_{10}}^T h_{11}(t_{11})
\int\limits_{t_{11}}^T h_{12}(t_{12})
dt_{12} dt_{11}\right) dt_{10}=
$$

\vspace{2mm}
$$
=\int\limits_t^T
\underline{h_{10}(t_{10})}
\int\limits_t^{t_{10}}\underline{h_{6}(t_{6})}\int\limits_t^{t_{6}}
\underline{h_{4}(t_{4})}
\int\limits_t^{t_{4}}\underline{h_{3}(t_{3})}
\left(\int\limits_t^{t_{3}}h_{2}(t_{2})
\int\limits_t^{t_{2}}h_{1}(t_{1})
dt_1 dt_2\right) \left(\int\limits_{t_4}^{t_6} h_5(t_5)dt_5 \right) \times
$$

\vspace{2mm}
\begin{equation}
\label{copa1}
\times 
\left(\int\limits_{t_6}^{t_{10}}h_9(t_9)
\int\limits_{t_6}^{t_{9}}h_8(t_8)
\int\limits_{t_6}^{t_{8}}h_7(t_7)dt_7 dt_8 dt_9
\right)
\left(\int\limits_{t_{10}}^T h_{12}(t_{12})
\int\limits_{t_{10}}^{t_{12}} h_{11}(t_{11})
dt_{11} dt_{12}\right) dt_3 dt_4 dt_6 dt_{10}.
\end{equation}

\vspace{5mm}

Further, suppose that $h_l(\tau)=\psi_l(\tau)\phi_{j_l}(\tau)$ $(l=1,\ldots,12)$ in 
(\ref{copa1}) (here $\left\{\phi_j(x)\right\}_{j=0}^{\infty}$
is an ar\-bit\-ra\-ry complete orthonormal system of 
functions in the space $L_2([t,T])$ and
$\psi_1(\tau),\ldots ,\psi_{12}(\tau)\in $ $L_2([t, T])$).
Thus, we get

\vspace{-1mm}
$$
C_{j_{12} j_{11} \underline{j_{10}} j_9 j_8 j_7 \underline{j_6} j_5 \underline{j_4} 
\underline{j_3} j_2 j_1}=
\int\limits_t^T
\psi_{10}(t_{10})\phi_{j_{10}}(t_{10})
\int\limits_t^{t_{10}}
\psi_{6}(t_{6})\phi_{j_{6}}(t_{6})
\int\limits_t^{t_6}
\psi_{4}(t_{4})\phi_{j_{4}}(t_{4})\times
$$

$$
\times
\int\limits_t^{t_4}
\psi_{3}(t_{3})\phi_{j_{3}}(t_{3})
C_{j_{12}j_{11}}^{\psi_{12}\psi_{11}}(T,t_{10})
C_{j_{9}j_{8}j_7}^{\psi_{9}\psi_{8}\psi_7}(t_{10},t_6)
C_{j_{5}}^{\psi_{5}}(t_6,t_{4})
C_{j_{2}j_{1}}^{\psi_{2}\psi_{1}}(t_{3},t)\times
$$

\begin{equation}
\label{copa1a}
\times
dt_3 dt_4 dt_6 dt_{10},
\end{equation}

\vspace{2mm}
\noindent
where (here and further)
$$
C_{j_m \ldots j_l}^{\psi_m\ldots \psi_l}(s,\tau)=\int\limits_{\tau}^s
\psi_m(t_m)\phi_{j_m}(t_m)\ldots
\int\limits_{\tau}^{t_{l+1}}
\psi_l(t_l)\phi_{j_l}(t_l)dt_l\ldots dt_m,
$$

\vspace{3mm}
\noindent
where $t\le\tau<s\le T$, $m\ge l$ $(m, l\in {\mathbb N}).$

Suppose that 
$g_1,g_2,\ldots,g_{2r-1},g_{2r}$ as in (\ref{leto5007xxx}) 
and $k>2r,$ $r\ge 1$
(the case $k=2r$ see in Sect.~26).
Consider $d_1,e_1,$ $\ldots,d_f,e_f, f \in\mathbb{N}$
such that 

\vspace{-2mm}

$$
1\le d_1-e_1+1<\ldots < d_1-1 <  d_1 <\ldots < d_f-e_f+1 <\ldots < d_f-1 < d_f \le k,
$$

$$
e_1+e_2+\ldots+e_f=2r,
$$

$$
\left\{g_1,g_2,\ldots,g_{2r-1},g_{2r}\right\}=\left\{d_1-e_1+1,\ldots, d_1\right\}\cup\ldots
\cup  \left\{d_f-e_f+1,\ldots, d_f\right\},
$$

$$
\{1,\ldots, k\}~\backslash \left\{g_1,g_2,\ldots,g_{2r-1},g_{2r}\right\}=
\left\{q_1,\ldots,q_{k-2r}\right\}.
$$

\vspace{6mm}

{\it We will say that the condition $(A)$ is satisfied if~~$\forall$   
$\left\{g_{2l-1},g_{2l}\right\}$ $(l=1,\ldots,r)$ $\exists$ $h\in \left\{1,\ldots,f\right\}$
such that
\begin{equation}
\label{copa2}
\left\{g_{2l-1},g_{2l}\right\}\subset
\left\{d_h-e_h+1,\ldots, d_h\right\}.
\end{equation}

\vspace{4mm}
\noindent
Moreover$,$ $\forall$ $h\in \left\{1,\ldots,f\right\}$
$\exists$ $\left\{g_{2l-1},g_{2l}\right\}$ $(l=1,\ldots,r)$
such that  {\rm (\ref{copa2})} is fulfilled.}

\vspace{2mm}

If the condition $(A)$ is satisfied, then $e_1, \ldots, e_f$ are even and we can write

$$
\left\{d_1-e_1+1,\ldots, d_1\right\}=
\left\{g_1^{(1)},g_2^{(1)},\ldots,g_{2r_1-1}^{(1)},g_{2r_1}^{(1)}\right\},
$$

\vspace{-2mm}

$$
\ldots
$$

\vspace{-2mm}
$$
\left\{d_f-e_f+1,\ldots, d_f\right\}=
\left\{g_1^{(f)},g_2^{(f)},\ldots,g_{2r_f-1}^{(f)},g_{2r_f}^{(f)}\right\},
$$

\vspace{2mm}
$$
\bigl\{g_1,g_2,\ldots,g_{2r-1},g_{2r}\bigr\}=
$$

\vspace{-1mm}
$$
=
\left\{g_1^{(1)},g_2^{(1)},\ldots,g_{2r_1-1}^{(1)},g_{2r_1}^{(1)},\ldots, 
g_1^{(f)},g_2^{(f)},\ldots,g_{2r_f-1}^{(f)},g_{2r_f}^{(f)}\right\}.
$$

\vspace{4mm}

\noindent
If the condition $(A)$ is not fulfilled, then some of $e_1, \ldots, e_f$ can be uneven.

Using (\ref{july90000}) and a modification of the algorithm from Sect.~26 
(see below for details) it can be proved that

$$
\lim\limits_{p\to\infty}
\sum\limits_{j_{g_1}, j_{g_3},\ldots ,j_{g_{2r-1}}=0}^p
\left(
C_{j_{d_f}\ldots j_{d_f-e_f+1}}^{\psi_{d_f}\ldots \psi_{d_f-e_f+1}}(t_{d_f+1},t_{d_f-e_f})
\ldots \right.
$$

\vspace{2.5mm}
$$
\left.
\ldots C_{j_{d_1}\ldots j_{d_1-e_1+1}}^{\psi_{d_1}\ldots \psi_{d_1-e_1+1}}(t_{d_1+1},t_{d_1-e_1})
\right)\biggl|_{j_{g_1}=j_{g_2},\ldots, j_{g_{2r-1}}=j_{g_{2r}}}= 
$$

\vspace{4mm}
$$
=\prod\limits_{h=1}^f
\frac{1}{2^{r_h}}\prod\limits_{l=1}^{r_h}
{\bf 1}_{\{g_{2l}^{(h)}=g_{2l-1}^{(h)}+1\}}\times
$$

\vspace{2mm}
\begin{equation}
\label{copa3}
\times C_{j_{d_h}\ldots j_{d_h-e_h+1}}^{\psi_{d_h}\ldots \psi_{d_h-e_h+1}}(t_{d_h+1},t_{d_h-e_h})
\biggl|_{(j_{g_2^{(h)}} j_{g_1^{(h)}})\curvearrowright (\cdot)
\ldots (j_{g_{2r_h}^{(h)}} j_{g_{2r_h-1}^{(h)}})\curvearrowright (\cdot),
j_{g_{{}_{1}}^{(h)}}=j_{g_{{}_{2}}^{(h)}},\ldots, j_{g_{{}_{2r_h-1}}^{(h)}}=j_{g_{{}_{2r_h}}^{(h)}}
}\biggr.
\end{equation}

\vspace{4mm}
\noindent
if the condition $(A)$ is satisfied, and 

\vspace{2mm}
$$
\lim\limits_{p\to\infty}
\sum\limits_{j_{g_1}, j_{g_3},\ldots ,j_{g_{2r-1}}=0}^p
\left(
C_{j_{d_f}\ldots j_{d_f-e_f+1}}^{\psi_{d_f}\ldots \psi_{d_f-e_f+1}}(t_{d_f+1},t_{d_f-e_f})
\ldots \right.
$$

\vspace{2.5mm}
\begin{equation}
\label{copa4}
\left.
\ldots C_{j_{d_1}\ldots j_{d_1-e_1+1}}^{\psi_{d_1}\ldots \psi_{d_1-e_1+1}}(t_{d_1+1},t_{d_1-e_1})
\right)\biggl|_{j_{g_1}=j_{g_2},\ldots, j_{g_{2r-1}}=j_{g_{2r}}}=0
\end{equation}

\vspace{4mm}
\noindent
if the condition $(A)$ is not fulfilled, where $t_{k+1}\stackrel{\sf def}{=}T,$
$t_{0}\stackrel{\sf def}{=}t,$
$e_1+\ldots+e_f=2r$ in (\ref{copa3}), (\ref{copa4})
and
$e_h=2 r_h$ $(h=1,\ldots, f),$ $r_1+\ldots+r_f=r$ in (\ref{copa3}). 

Note that the series on the left-hand sides of 
(\ref{copa3}) and (\ref{copa4}) converge absolutly 
since
their sums do not depend 
on permutations of basis functions
(here the basis in $L_2([t,T]^{r})$ has the following form
$\left\{\phi_{j_1}(x_1)\ldots \phi_{j_r}(x_r)\right\}_{j_1,\ldots,j_r=0}^{\infty}$).
Recall that any permutation of basis functions in a Hilbert space forms a basis 
in this Hilbert space \cite{gohb}.

Let us prove the formulas (\ref{copa3}) and (\ref{copa4}).

\vspace{2mm}

1. Suppose that the condition $(A)$ is satisfied and 

\vspace{-1mm}
\begin{equation}
\label{copa5}
\prod\limits_{l=1}^{r_h}
{\bf 1}_{\{g_{2l}^{(h)}=g_{2l-1}^{(h)}+1\}}=1
\end{equation}

\vspace{3mm}
\noindent
for all $h=1,\ldots,f.$ In this case we can use the results from Sect.~26.
We have (see (\ref{july90000}))

$$
\lim\limits_{p\to\infty}
\sum\limits_{j_{g_1}, j_{g_3},\ldots ,j_{g_{2r-1}}=0}^p
\left(
C_{j_{d_f}\ldots j_{d_f-e_f+1}}^{\psi_{d_f}\ldots \psi_{d_f-e_f+1}}(t_{d_f+1},t_{d_f-e_f})
\ldots \right.
$$

\vspace{2.5mm}
$$
\left.
\ldots C_{j_{d_1}\ldots j_{d_1-e_1+1}}^{\psi_{d_1}\ldots \psi_{d_1-e_1+1}}(t_{d_1+1},t_{d_1-e_1})
\right)\biggl|_{j_{g_1}=j_{g_2},\ldots, j_{g_{2r-1}}=j_{g_{2r}}}= 
$$

\vspace{6mm}
$$
=\lim\limits_{p\to\infty}
\sum\limits_{j_{g_{{}_{1}}^{(1)}}, j_{g_{{}_{3}}^{(1)}},\ldots, j_{g_{{}_{2r_1-1}}^{(1)}}=0}^p
C_{j_{d_1}\ldots j_{d_1-e_1+1}}^{\psi_{d_1}\ldots \psi_{d_1-e_1+1}}(t_{d_1+1},t_{d_1-e_1})
\biggl|_{
j_{g_{{}_{1}}^{(1)}}=j_{g_{{}_{2}}^{(1)}},\ldots, j_{g_{{}_{2r_1-1}}^{(1)}}=j_{g_{{}_{2r_1}}^{(1)}}
}\times
$$
\vspace{-4mm}

$$
\ldots
$$
\vspace{-3.5mm}

$$
\times 
\lim\limits_{p\to\infty}
\sum\limits_{j_{g_{{}_{1}}^{(f)}}, j_{g_{{}_{3}}^{(f)}}, \ldots , j_{g_{{}_{2r_f-1}}^{(f)}}=0}^p
C_{j_{d_f}\ldots j_{d_f-e_f+1}}^{\psi_{d_f}\ldots \psi_{d_f-e_f+1}}(t_{d_f+1},t_{d_f-e_f})
\biggl|_{
j_{g_{{}_{1}}^{(f)}}=j_{g_{{}_{2}}^{(f)}},\ldots, j_{g_{{}_{2r_f-1}}^{(f)}}=j_{g_{{}_{2r_f}}^{(f)}}
}=
$$

\vspace{4mm}

$$
=\prod\limits_{h=1}^f
\frac{1}{2^{r_h}}\prod\limits_{l=1}^{r_h}
{\bf 1}_{\{g_{2l}^{(h)}=g_{2l-1}^{(h)}+1\}}\times
$$

\vspace{2mm}

$$
\times C_{j_{d_h}\ldots j_{d_h-e_h+1}}^{\psi_{d_h}\ldots \psi_{d_h-e_h+1}}(t_{d_h+1},t_{d_h-e_h})
\biggl|_{(j_{g_2^{(h)}} j_{g_1^{(h)}})\curvearrowright (\cdot)
\ldots (j_{g_{2r_h}^{(h)}} j_{g_{2r_h-1}^{(h)}})\curvearrowright (\cdot),
j_{g_{{}_{1}}^{(h)}}=j_{g_{{}_{2}}^{(h)}},\ldots, j_{g_{{}_{2r_h-1}}^{(h)}}=j_{g_{{}_{2r_h}}^{(h)}}
}\biggr..
$$

\vspace{6mm}

\noindent
Thus, we get the formula (\ref{copa3}).

\vspace{2mm}

2. Suppose that the condition $(A)$ is satisfied and for some $h=1,\ldots,f$

\vspace{-1mm}
\begin{equation}
\label{copa6}
\prod\limits_{l=1}^{r_h}
{\bf 1}_{\{g_{2l}^{(h)}=g_{2l-1}^{(h)}+1\}}=0.
\end{equation}

\vspace{3mm}

In this case, we act the same as in the previous case.
Applying (\ref{july90000}), we obtain

\vspace{1mm}
$$
\lim\limits_{p\to\infty}
\sum\limits_{j_{g_1}, j_{g_3}, \ldots,  j_{g_{2r-1}}=0}^p
\left(
C_{j_{d_f}\ldots j_{d_f-e_f+1}}^{\psi_{d_f}\ldots \psi_{d_f-e_f+1}}(t_{d_f+1},t_{d_f-e_f})
\ldots \right.
$$

\vspace{3mm}
$$
\left.
\ldots C_{j_{d_1}\ldots j_{d_1-e_1+1}}^{\psi_{d_1}\ldots \psi_{d_1-e_1+1}}(t_{d_1+1},t_{d_1-e_1})
\right)\biggl|_{j_{g_1}=j_{g_2},\ldots, j_{g_{2r-1}}=j_{g_{2r}}}= 
$$

\vspace{6mm}
$$
=\lim\limits_{p\to\infty}
\sum\limits_{j_{g_{{}_{1}}^{(1)}}, j_{g_{{}_{3}}^{(1)}}, \ldots ,
j_{g_{{}_{2r_1-1}}^{(1)}}=0}^p
C_{j_{d_1}\ldots j_{d_1-e_1+1}}^{\psi_{d_1}\ldots \psi_{d_1-e_1+1}}(t_{d_1+1},t_{d_1-e_1})
\biggl|_{
j_{g_{{}_{1}}^{(1)}}=j_{g_{{}_{2}}^{(1)}},\ldots, j_{g_{{}_{2r_1-1}}^{(1)}}=j_{g_{{}_{2r_1}}^{(1)}}
}\times
$$
\vspace{-4mm}

$$
\ldots
$$
\vspace{-3mm}

\begin{equation}
\label{2024october1}
\times 
\lim\limits_{p\to\infty}
\sum\limits_{j_{g_{{}_{1}}^{(f)}}, j_{g_{{}_{3}}^{(f)}},\ldots ,
j_{g_{{}_{2r_f-1}}^{(f)}}=0}^p
C_{j_{d_f}\ldots j_{d_f-e_f+1}}^{\psi_{d_f}\ldots \psi_{d_f-e_f+1}}(t_{d_f+1},t_{d_f-e_f})
\biggl|_{
j_{g_{{}_{1}}^{(f)}}=j_{g_{{}_{2}}^{(f)}},\ldots, j_{g_{{}_{2r_f-1}}^{(f)}}=j_{g_{{}_{2r_f}}^{(f)}}
}=0
\end{equation}

\vspace{5mm}
\noindent
(al least one of the multipliers is equal to zero on the right-hand side of (\ref{2024october1})).

The equality (\ref{copa3}) is proved in our case 
(the right-hand side of (\ref{copa3}) is equal to zero for the considered case (see (\ref{copa6}))).

\vspace{2mm}

3. Suppose that the condition $(A)$ is not satisfied. 
In this case, we act according to the algorithm
from Sect.~26.
More precisely, let us select blocks
in the multi-index $j_{d_h}\ldots j_{d_h-e_h+1}$ $(h=1,\ldots,f)$
that
correspond to the fulfillment of the condition 

\vspace{-1mm}
$$
\prod\limits_{l=1}^{r_{m,h}} {\bf 1}_{\{g_{2l}^{(h)}=g_{2l-1}^{(h)}+1\}}=1,
$$

\vspace{2mm}
\noindent
where $r_{m,h}$ is the number of pairs $\{g_{2l-1}^{(h)}, g_{2l}^{(h)}\}$ (from the set 
$\{g_1,g_2,\ldots,$ $g_{2r-1},g_{2r}\})$
in the block with number $m$ that corresponds to 
the multi-index $j_{d_h}\ldots j_{d_h-e_h+1}$.

Let us save multipliers of the form 
$$
{\bf 1}_{\{t_n<t_{n+1}\}}
$$

\vspace{2mm}
\noindent
in the  Volterra--type kernels corresponding to the Fourier
coefficients

\vspace{-1mm}
\begin{equation}
\label{copa7}
C_{j_{d_1}\ldots j_{d_1-e_1+1}}^{\psi_{d_1}\ldots \psi_{d_1-e_1+1}}(t_{d_1+1},t_{d_1-e_1}),
\ldots, 
C_{j_{d_f}\ldots j_{d_f-e_f+1}}^{\psi_{d_f}\ldots \psi_{d_f-e_f+1}}(t_{d_f+1},t_{d_f-e_f})
\end{equation}

\vspace{3mm}
\noindent
and corresponding to the above blocks.

At that, we remove the remaining 
multipliers of the form 

\vspace{-2mm}
$$
{\bf 1}_{\{t_n<t_{n+1}\}}
$$

\vspace{2mm}
\noindent
in the  Volterra--type kernels corresponding to the Fourier
coefficients (\ref{copa7}).

As a result, we get a modified left-hand side
of the equality (\ref{copa4}).
For definiteness, let us denote this expression by
$({}^{-})$.

Using generalized Parseval's equality
(Parseval's equality for two functions)
and (\ref{july30016}), we represent
the expression $({}^{-})$ as an integral over the hypercube $[t, T]^r$.

It is not difficult to see that the indicated integral over the hypercube $[t, T]^r$ is represented as a product
of integrals over hypercubes of smaller dimentions.
At that, at least one of these integrals 
is equal to zero
due to the generalized Parseval equality (Parseval's equality for two functions)
and the fulfillment of the condition

\vspace{-2mm}
$$
t\le t_{d_1-e_1}\le t_{d_1+1}\le \ldots \le  t_{d_f-e_f} \le  t_{d_f+1} \le T
$$

\vspace{2mm}
\noindent
(see the above example and (\ref{copa1}) and (\ref{copa1a})).
For definiteness, let us denote the equality of $({}^{-})$ to zero by $(\bar K)$.
We interpret the above zero as the zero functional in $L_2([t,T]^r).$
Further, transformations and passages to the limit
in the equality $(\bar K)$ are performed iteratively 
in such a way as to restore the removed multipliers ${\bf 1}_{\{t_n<t_{n+1}\}}$
on the left-hand side of $(\bar K)$
(for more details, see Sect.~26).
As a result, we obtain the equality (\ref{copa4}).
The equalities (\ref{copa3}) and (\ref{copa4}) are proved.

For definiteness, suppose that $q_1<\ldots <q_{k-2r}$
and $k>2r,$ $r\ge 1$
(the case $k=2r$ see in Sect.~26).
Using Fubini's Theorem (as in the above example (see (\ref{copa1})), we 
obtain

\vspace{1mm}
$$
\sum\limits_{j_{g_1},j_{g_3},\ldots, j_{g_{2r-1}}=0}^p
C_{j_k\ldots j_1}\biggl|_{j_{g_1}=j_{g_2},\ldots, j_{g_{2r-1}}=j_{g_{2r}}}=
$$

\vspace{3.5mm}
$$
=
\int\limits_t^T \psi_{q_{k-2r}}(t_{q_{k-2r}})
\phi_{j_{q_{k-2r}}}(t_{q_{k-2r}})\ldots
\int\limits_t^{t_{q_1+1}} \psi_{q_{1}}(t_{q_{1}})
\phi_{j_{q_{1}}}(t_{q_{1}})\times
$$

\vspace{3mm}
$$
\times 
\sum\limits_{j_{g_1},j_{g_3},\ldots, j_{g_{2r-1}}=0}^p\left(
C_{j_{d_f}\ldots j_{d_f-e_f+1}}^{\psi_{d_f}\ldots \psi_{d_f-e_f+1}}(t_{d_f+1},t_{d_f-e_f})
\ldots \right.
$$

\vspace{3.5mm}
$$
\left.
\ldots C_{j_{d_1}\ldots j_{d_1-e_1+1}}^{\psi_{d_1}\ldots \psi_{d_1-e_1+1}}(t_{d_1+1},t_{d_1-e_1})
\right)\biggl|_{j_{g_1}=j_{g_2},\ldots, j_{g_{2r-1}}=j_{g_{2r}}}\times
$$

\vspace{1mm}
\begin{equation}
\label{copa9}
\times dt_{q_1}\ldots dt_{q_{k-2r}},
\end{equation}

\vspace{4.5mm}
$$
\frac{1}{2^r} \prod\limits_{l=1}^r {\bf 1}_{\{g_{2l}=g_{2l-1}+1\}}
C_{j_k \ldots j_1}\biggl|_{(j_{g_2} j_{g_1})\curvearrowright (\cdot)
\ldots (j_{g_{2r}} j_{g_{2r-1}})\curvearrowright (\cdot),
j_{g_{{}_{1}}}=~j_{g_{{}_{2}}},\ldots, j_{g_{{}_{2r-1}}}=~j_{g_{{}_{2r}}}
}\biggr.=
$$

\vspace{3.5mm}
$$
=
\int\limits_t^T \psi_{q_{k-2r}}(t_{q_{k-2r}})
\phi_{j_{q_{k-2r}}}(t_{q_{k-2r}})\ldots
\int\limits_t^{t_{q_1+1}} \psi_{q_{1}}(t_{q_{1}})
\phi_{j_{q_{1}}}(t_{q_{1}})\times
$$

\vspace{2mm}
$$
\times {\bf 1}_{\{the\ condition\ (A)\ is\ satisfied\}}\prod\limits_{h=1}^f
\frac{1}{2^{r_h}}\prod\limits_{l=1}^{r_h}
{\bf 1}_{\{g_{2l}^{(h)}=g_{2l-1}^{(h)}+1\}}\times
$$

\vspace{2mm}
$$
\times C_{j_{d_h}\ldots j_{d_h-e_h+1}}^{\psi_{d_h}\ldots \psi_{d_h-e_h+1}}(t_{d_h+1},t_{d_h-e_h})
\biggl|_{(j_{g_2^{(h)}} j_{g_1^{(h)}})\curvearrowright (\cdot)
\ldots (j_{g_{2r_h}^{(h)}} j_{g_{2r_h-1}^{(h)}})\curvearrowright (\cdot),
j_{g_{{}_{1}}^{(h)}}=j_{g_{{}_{2}}^{(h)}},\ldots, j_{g_{{}_{2r_h-1}}^{(h)}}=j_{g_{{}_{2r_h}}^{(h)}}
}\biggr.\hspace{-0.7mm}\times
$$

\vspace{1mm}
\begin{equation}
\label{copa10}
\times
dt_{q_1}\ldots dt_{q_{k-2r}}.
\end{equation}

\vspace{6mm}

Suppose that

\vspace{-2.5mm}
$$
\Biggl|
\sum\limits_{j_{g_1}, j_{g_3},\ldots ,j_{g_{2r-1}}=0}^p
\left(
C_{j_{d_f}\ldots j_{d_f-e_f+1}}^{\psi_{d_f}\ldots \psi_{d_f-e_f+1}}(t_{d_f+1},t_{d_f-e_f})
\ldots\right.\Biggr.
$$

\vspace{2mm}
\begin{equation}
\label{09091}
\Biggl.\left. 
\ldots C_{j_{d_1}\ldots j_{d_1-e_1+1}}^{\psi_{d_1}\ldots \psi_{d_1-e_1+1}}(t_{d_1+1},t_{d_1-e_1})
\right)\biggl|_{j_{g_1}=j_{g_2},\ldots, j_{g_{2r-1}}=j_{g_{2r}}}\Biggr|\le K <\infty,
\end{equation}

\vspace{3.5mm}
\noindent
where constant $K$ does not depend on $p$ and 
$t_{d_1+1},t_{d_1-e_1},\ldots,t_{d_f+1},t_{d_f-e_f}$
(here $d_1-e_1\ge 1$ and $d_f+1\le k$).
In (\ref{09091}):\ 
$t_{k+1}\stackrel{\sf def}{=}T,$
$t_{0}\stackrel{\sf def}{=}t,$
$e_1+\ldots+e_f=2r;$ another notations as above in this section.

Applying (\ref{copa3}), (\ref{copa4}), (\ref{copa9}), (\ref{copa10}), we obtain 
($k>2r,$ $r\ge 1$)

\vspace{1mm}
$$
\lim\limits_{p\to\infty}
\sum\limits_{\stackrel{j_1,\ldots,j_q,\ldots,j_k=0}{{}_{q\ne g_1, g_2, \ldots, g_{2r-1},
g_{2r}}}}^p
\Biggl(\sum\limits_{j_{g_1},j_{g_3},\ldots, j_{g_{2r-1}}=0}^p
C_{j_k\ldots j_1}\biggl|_{j_{g_1}=j_{g_2},\ldots, j_{g_{2r-1}}=j_{g_{2r}}}-\Biggr.
$$

\vspace{3mm}
$$
\Biggl.-\frac{1}{2^r} \prod\limits_{l=1}^r {\bf 1}_{\{g_{2l}=g_{2l-1}+1\}}
C_{j_k \ldots j_1}\biggl|_{(j_{g_2} j_{g_1})\curvearrowright (\cdot)
\ldots (j_{g_{2r}} j_{g_{2r-1}})\curvearrowright (\cdot),
j_{g_{{}_{1}}}=~j_{g_{{}_{2}}},\ldots, j_{g_{{}_{2r-1}}}=~j_{g_{{}_{2r}}}
}\biggr.\Biggr)^2\le
$$

\vspace{6mm}
$$
\le \lim\limits_{p\to\infty}
\sum\limits_{\stackrel{j_1,\ldots,j_q,\ldots,j_k=0}{{}_{q\ne g_1, g_2, \ldots, g_{2r-1},
g_{2r}}}}^{\infty}
\Biggl(\sum\limits_{j_{g_1},j_{g_3},\ldots, j_{g_{2r-1}}=0}^p
C_{j_k\ldots j_1}\biggl|_{j_{g_1}=j_{g_2},\ldots, j_{g_{2r-1}}=j_{g_{2r}}}-\Biggr.
$$

\vspace{3mm}
$$
\Biggl.-\frac{1}{2^r} \prod\limits_{l=1}^r {\bf 1}_{\{g_{2l}=g_{2l-1}+1\}}
C_{j_k \ldots j_1}\biggl|_{(j_{g_2} j_{g_1})\curvearrowright (\cdot)
\ldots (j_{g_{2r}} j_{g_{2r-1}})\curvearrowright (\cdot),
j_{g_{{}_{1}}}=~j_{g_{{}_{2}}},\ldots, j_{g_{{}_{2r-1}}}=~j_{g_{{}_{2r}}}
}\biggr.\Biggr)^2=
$$

\vspace{8mm}
$$
=\lim\limits_{p\to\infty}
\sum\limits_{j_{q_1},\ldots,j_{q_{k-2r}}=0}^{\infty}
\Biggl(
\int\limits_t^T \psi_{q_{k-2r}}(t_{q_{k-2r}})
\phi_{j_{q_{k-2r}}}(t_{q_{k-2r}})\ldots
\int\limits_t^{t_{q_1+1}} \psi_{q_{1}}(t_{q_{1}})
\phi_{j_{q_{1}}}(t_{q_{1}})\times\Biggr.
$$

\vspace{3mm}
$$
\times 
\Biggl(\sum\limits_{j_{g_1},j_{g_3},\ldots, j_{g_{2r-1}}=0}^p\left(
C_{j_{d_f}\ldots j_{d_f-e_f+1}}^{\psi_{d_f}\ldots \psi_{d_f-e_f+1}}(t_{d_f+1},t_{d_f-e_f})
\ldots \right.\Biggr.
$$

\vspace{5mm}
$$
\left.
\ldots C_{j_{d_1}\ldots j_{d_1-e_1+1}}^{\psi_{d_1}\ldots \psi_{d_1-e_1+1}}(t_{d_1+1},t_{d_1-e_1})
\right)\biggl|_{j_{g_1}=j_{g_2},\ldots, j_{g_{2r-1}}=j_{g_{2r}}}-
$$

\vspace{5mm}
$$
-{\bf 1}_{\{the\ condition\ (A)\ is\ satisfied\}}\prod\limits_{h=1}^f
\frac{1}{2^{r_h}}\prod\limits_{l=1}^{r_h}
{\bf 1}_{\{g_{2l}^{(h)}=g_{2l-1}^{(h)}+1\}}\times
$$
$$
\Biggl.\times C_{j_{d_h}\ldots j_{d_h-e_h+1}}^{\psi_{d_h}\ldots \psi_{d_h-e_h+1}}(t_{d_h+1},t_{d_h-e_h})
\biggl|_{(j_{g_2^{(h)}} j_{g_1^{(h)}})\curvearrowright (\cdot)
\ldots (j_{g_{2r_h}^{(h)}} j_{g_{2r_h-1}^{(h)}})\curvearrowright (\cdot),
j_{g_{{}_{1}}^{(h)}}=j_{g_{{}_{2}}^{(h)}},\ldots, j_{g_{{}_{2r_h-1}}^{(h)}}=j_{g_{{}_{2r_h}}^{(h)}}
}\biggr.\hspace{-0.5mm}\Biggr)\hspace{-0.5mm}\times
$$

\vspace{3mm}
\begin{equation}
\label{2024dec1}
\Biggl.\times
dt_{q_1}\ldots dt_{q_{k-2r}}\Biggr)^2=
\end{equation}

\vspace{5mm}
$$
=\lim\limits_{p\to\infty}
\int\limits_t^T \psi_{q_{k-2r}}^2(t_{q_{k-2r}})
\ldots
\int\limits_t^{t_{q_1+1}} \psi_{q_{1}}^2(t_{q_{1}})\times
$$

\vspace{3mm}
$$
\times
\Biggl(\sum\limits_{j_{g_1},j_{g_3},\ldots, j_{g_{2r-1}}=0}^p\left(
C_{j_{d_f}\ldots j_{d_f-e_f+1}}^{\psi_{d_f}\ldots \psi_{d_f-e_f+1}}(t_{d_f+1},t_{d_f-e_f})
\ldots \right.\Biggr.
$$

\vspace{4mm}
$$
\left.
\ldots C_{j_{d_1}\ldots j_{d_1-e_1+1}}^{\psi_{d_1}\ldots \psi_{d_1-e_1+1}}(t_{d_1+1},t_{d_1-e_1})
\right)\biggl|_{j_{g_1}=j_{g_2},\ldots, j_{g_{2r-1}}=j_{g_{2r}}}-
$$

\vspace{5mm}
$$
-{\bf 1}_{\{the\ condition\ (A)\ is\ satisfied\}}\prod\limits_{h=1}^f
\frac{1}{2^{r_h}}\prod\limits_{l=1}^{r_h}
{\bf 1}_{\{g_{2l}^{(h)}=g_{2l-1}^{(h)}+1\}}\times
$$

\vspace{3mm}
$$
\Biggl.\times C_{j_{d_h}\ldots j_{d_h-e_h+1}}^{\psi_{d_h}\ldots \psi_{d_h-e_h+1}}(t_{d_h+1},t_{d_h-e_h})
\biggl|_{(j_{g_2^{(h)}} j_{g_1^{(h)}})\curvearrowright (\cdot)
\ldots (j_{g_{2r_h}^{(h)}} j_{g_{2r_h-1}^{(h)}})\curvearrowright (\cdot),
j_{g_{{}_{1}}^{(h)}}=j_{g_{{}_{2}}^{(h)}},\ldots, j_{g_{{}_{2r_h-1}}^{(h)}}=j_{g_{{}_{2r_h}}^{(h)}}
}\biggr.\hspace{-0.5mm}\Biggr)^2\hspace{-1.7mm}\times
$$

\begin{equation}
\label{copa14}
\Biggl.\times
dt_{q_1}\ldots dt_{q_{k-2r}}=
\end{equation}

\vspace{5mm}
$$
=
\int\limits_t^T \psi_{q_{k-2r}}^2(t_{q_{k-2r}})
\ldots
\int\limits_t^{t_{q_1+1}} \psi_{q_{1}}^2(t_{q_{1}})\times
$$

\vspace{3mm}
$$
\times
\lim\limits_{p\to\infty}\Biggl(\sum\limits_{j_{g_1},j_{g_3},\ldots, j_{g_{2r-1}}=0}^p\left(
C_{j_{d_f}\ldots j_{d_f-e_f+1}}^{\psi_{d_f}\ldots \psi_{d_f-e_f+1}}(t_{d_f+1},t_{d_f-e_f})
\ldots \right.\Biggr.
$$

\vspace{3mm}
$$
\left.
\ldots C_{j_{d_1}\ldots j_{d_1-e_1+1}}^{\psi_{d_1}\ldots \psi_{d_1-e_1+1}}(t_{d_1+1},t_{d_1-e_1})
\right)\biggl|_{j_{g_1}=j_{g_2},\ldots, j_{g_{2r-1}}=j_{g_{2r}}}-
$$

\vspace{3mm}
$$
-{\bf 1}_{\{the\ condition\ (A)\ is\ satisfied\}}\prod\limits_{h=1}^f
\frac{1}{2^{r_h}}\prod\limits_{l=1}^{r_h}
{\bf 1}_{\{g_{2l}^{(h)}=g_{2l-1}^{(h)}+1\}}\times
$$

\vspace{1mm}
$$
\Biggl.\times C_{j_{d_h}\ldots j_{d_h-e_h+1}}^{\psi_{d_h}\ldots \psi_{d_h-e_h+1}}(t_{d_h+1},t_{d_h-e_h})
\biggl|_{(j_{g_2^{(h)}} j_{g_1^{(h)}})\curvearrowright (\cdot)
\ldots (j_{g_{2r_h}^{(h)}} j_{g_{2r_h-1}^{(h)}})\curvearrowright (\cdot),
j_{g_{{}_{1}}^{(h)}}=j_{g_{{}_{2}}^{(h)}},\ldots, j_{g_{{}_{2r_h-1}}^{(h)}}=j_{g_{{}_{2r_h}}^{(h)}}
}\biggr.\hspace{-0.5mm}\Biggr)^2\hspace{-1.7mm}\times
$$

\begin{equation}
\label{copa15}
\Biggl.\times
dt_{q_1}\ldots dt_{q_{k-2r}}=0,
\end{equation}

\vspace{4mm}
\noindent
where 
the transition from (\ref{2024dec1}) to (\ref{copa14})
is based on the Parseval equality
and the transition from (\ref{copa14}) to (\ref{copa15}) 
is based on Lebesgue's Dominated Convergence Theorem (see
(\ref{july999}), (\ref{july1000aaa1}), (\ref{copa3}), (\ref{copa4}), (\ref{09091})) 
and also on
convergence to zero (almost everywhere on 
$X=\{(t_{q_1},\ldots ,t_{q_{k-2r}}):  t\le t_{q_1}\le \ldots \le t_{q_{k-2r}}\le T\}$
with respect to Lebesgue's measure)
of the integrand function in (\ref{copa14}).

Thus, the equality (\ref{june100}) and Hypotheses~7, 8
are proved for the case $p_1=\ldots=p_k=p$ 
under the condition (\ref{09091})
and we have the following theorem.

\vspace{2mm}

{\bf Theorem~58}\ \cite{2018a}.\ {\it Suppose that 
the condition {\rm (\ref{09091})} 
is fulfilled$,$
$\{\phi_j(x)\}_{j=0}^{\infty}$
is an arbitrary complete orthonormal system of func\-ti\-ons
in $L_2([t, T])$ and
$\psi_1(\tau),\ldots, \psi_k(\tau)\in L_2([t, T]).$
Then$,$ for the sum $\bar J^{*}[\psi^{(k)}]_{T,t}^{(i_1\ldots i_k)}$
of iterated Ito stochastic integrals 

$$
\bar J^{*}[\psi^{(k)}]_{T,t}^{(i_1\ldots i_k)}=J[\psi^{(k)}]_{T,t}^{(i_1\ldots i_k)}+
\sum_{r=1}^{\left[k/2\right]}\frac{1}{2^r}
\sum_{(s_r,\ldots,s_1)\in {\rm A}_{k,r}}
J[\psi^{(k)}]_{T,t}^{s_r,\ldots,s_1}
$$

\vspace{3mm}
\noindent
the following 
expansion 

\vspace{-3mm}
$$
\bar J^{*}[\psi^{(k)}]_{T,t}^{(i_1\ldots i_k)}=
\hbox{\vtop{\offinterlineskip\halign{
\hfil#\hfil\cr
{\rm l.i.m.}\cr
$\stackrel{}{{}_{p\to \infty}}$\cr
}} }
\sum_{j_1,\ldots,j_k=0}^{p}
C_{j_k \ldots j_1}\prod\limits_{l=1}^k \zeta_{j_l}^{(i_l)}
$$

\vspace{4mm}
\noindent
that converges in the mean-square sense is valid, where 

$$
C_{j_k \ldots j_1}=\int\limits_t^T\psi_k(t_k)\phi_{j_k}(t_k)\ldots
\int\limits_t^{t_2}
\psi_1(t_1)\phi_{j_1}(t_1)
dt_1\ldots dt_k
$$

\vspace{4mm}
\noindent
is the Fourier coefficient, 
${\rm l.i.m.}$ is a limit in the mean-square sense,
$i_1, \ldots, i_k=0, 1,\ldots,m,$

\vspace{-1mm}
$$
\zeta_{j}^{(i)}=
\int\limits_t^T \phi_{j}(\tau) d{\bf w}_{\tau}^{(i)}
$$ 

\vspace{3mm}
\noindent
are independent standard Gaussian random variables for various 
$i$ or $j$ {\rm (}in the case when $i\ne 0${\rm )},
${\bf w}_{\tau}^{(i)}={\bf f}_{\tau}^{(i)}$ 
for $i=1,\ldots,m$ and 
${\bf w}_{\tau}^{(0)}=\tau;$ another notations are the same as
in Theorem~{\rm 5}.}

\vspace{2mm}

Using Theorem~5, we obtain the following corollary of Theorem~58.

\vspace{2mm}

{\bf Theorem~59}\ \cite{2018a}.\ {\it Suppose that 
the condition {\rm (\ref{09091})} 
is fulfilled$,$
$\{\phi_j(x)\}_{j=0}^{\infty}$
is an arbitrary complete orthonormal system of func\-ti\-ons
in $L_2([t, T])$ and
$\psi_1(\tau),\ldots, \psi_k(\tau)$ are continuous functions
at the interval $[t, T].$
Then$,$ for the iterated Stratonovich sto\-chas\-tic integral 
of multiplicity $k$ $(k\in\mathbb{N})$

$$
J^{*}[\psi^{(k)}]_{T,t}^{(i_1\ldots i_k)}=
{\int\limits_t^{*}}^T
\psi_k(t_k) \ldots 
{\int\limits_t^{*}}^{t_{2}}
\psi_1(t_1) d{\bf w}_{t_1}^{(i_1)}\ldots
d{\bf w}_{t_k}^{(i_k)}
$$

\vspace{3mm}
\noindent
the following 
expansion 

\vspace{-2mm}
\begin{equation}
\label{march000195}
J^{*}[\psi^{(k)}]_{T,t}^{(i_1\ldots i_k)}=
\hbox{\vtop{\offinterlineskip\halign{
\hfil#\hfil\cr
{\rm l.i.m.}\cr
$\stackrel{}{{}_{p\to \infty}}$\cr
}} }
\sum\limits_{j_1,\ldots,j_k=0}^{p}
C_{j_k \ldots j_1}\prod\limits_{l=1}^k \zeta_{j_l}^{(i_l)}
\end{equation}

\vspace{4mm}
\noindent
that converges in the mean-square sense is valid, where 
notations are the same as in Theorem~{\rm 58}.}

\vspace{5mm}

\section{Expansion of Iterated Stratonovich Stochastic Integrals
of Multiplicity 6. The Case of an Ar\-bit\-ra\-ry Complete Orthonormal System of 
Functions in the Space $L_2([t,T])$ and $\psi_1(\tau),\ldots, \psi_6(\tau)
\equiv 1$}

\vspace{5mm}

This section is devoted to the following theorem.

\vspace{2mm}

{\bf Theorem~60}\ \cite{2018a}.\  {\it Suppose that
$\{\phi_j(x)\}_{j=0}^{\infty}$ is an arbitrary complete orthonormal system of 
func\-ti\-ons in the space $L_2([t,T]).$
Then$,$ for the iterated Stra\-to\-no\-vich stochastic integral
of sixth multiplicity

$$
J^{*}[\psi^{(6)}]_{T,t}=
{\int\limits_t^{*}}^T
\ldots
{\int\limits_t^{*}}^{t_2}
d{\bf w}_{t_1}^{(i_1)}
\ldots d{\bf w}_{t_6}^{(i_6)}
$$

\vspace{2mm}
\noindent
the following 
expansion 
$$
J^{*}[\psi^{(6)}]_{T,t}=
\hbox{\vtop{\offinterlineskip\halign{
\hfil#\hfil\cr
{\rm l.i.m.}\cr
$\stackrel{}{{}_{p\to \infty}}$\cr
}} }
\sum\limits_{j_1,\ldots, j_6=0}^{p}
C_{j_6 \ldots j_1}\zeta_{j_1}^{(i_1)}\ldots \zeta_{j_6}^{(i_6)}
$$

\vspace{3mm}
\noindent
that converges in the mean-square sense is valid, where 
$i_1,\ldots,i_6=0, 1,\ldots,m,$

\begin{equation}
\label{november1620241}
C_{j_6\ldots j_1}=\int\limits_t^T
\phi_{j_6}(t_6)
\ldots
\int\limits_t^{t_2}
\phi_{j_1}(t_1)dt_1\ldots dt_6
\end{equation}
and
$$
\zeta_{j}^{(i)}=
\int\limits_t^T \phi_{j}(\tau) d{\bf w}_{\tau}^{(i)}
$$ 

\vspace{2mm}
\noindent
are independent standard Gaussian random variables for various 
$i$ or $j$ {\rm (}in the case when $i\ne 0${\rm ),}
${\bf w}_{\tau}^{(i)}={\bf f}_{\tau}^{(i)}$ for
$i=1,\ldots,m$ and 
${\bf w}_{\tau}^{(0)}=\tau.$}

\vspace{2mm}

{\bf Proof.}\ Our proof will be based on Theorem~59
and verification of the equality (\ref{09091}) under the conditions
of Theorem~60 (the case $k=6>2r$, where $r=1, 2$). Recall that the case
$k=2r$ is considered in Sect.~26 (see (\ref{july90000})).
Under the conditions of Theorem~60, this means that
$k=6=2r$, where $r=3$.

Let throughout this proof 

\vspace{-1mm}
$$
C_{j_k \ldots j_1}(s,\tau)=\int\limits_{\tau}^s
\phi_{j_k}(t_k)\ldots
\int\limits_{\tau}^{t_2}
\phi_{j_1}(t_1)dt_1\ldots dt_k \ \ \ (k=1,\ldots,4,\ t\le\tau<s\le T),
$$

\vspace{3mm}
\noindent
and $C_{j_6\ldots j_1}$ is defined by (\ref{november1620241}).

Using Fubini's Theorem and the technique that leads to the formulas (\ref{copa1}),
(\ref{copa1a}),
we obtain (note that we find all possible
combinations of pairs using the equality (\ref{after36})):

\vspace{3mm}

1. $r=1$ (15 combinations)

$$
C_{j_1 j_5 j_4 j_3 j_2 j_1}=
\int\limits_t^T
\phi_{j_5}(t_5)
\int\limits_t^{t_5}
\phi_{j_4}(t_4)
\int\limits_t^{t_4}
\phi_{j_3}(t_3)
\int\limits_t^{t_3}
\phi_{j_2}(t_2)
C_{j_1}(t_2,t) C_{j_1}(T,t_5)dt_2 dt_3 dt_4 dt_5,
$$

\vspace{2mm}
$$
C_{j_2 j_5 j_4 j_3 j_2 j_1}=
\int\limits_t^T
\phi_{j_5}(t_5)
\int\limits_t^{t_5}
\phi_{j_4}(t_4)
\int\limits_t^{t_4}
\phi_{j_3}(t_3)
\int\limits_t^{t_3}
\phi_{j_1}(t_1)
C_{j_2}(t_3,t_1) C_{j_2}(T,t_5)dt_1 dt_3 dt_4 dt_5,
$$

\vspace{2mm}
$$
C_{j_3 j_5 j_4 j_3 j_2 j_1}=
\int\limits_t^T
\phi_{j_5}(t_5)
\int\limits_t^{t_5}
\phi_{j_4}(t_4)
\int\limits_t^{t_4}
\phi_{j_2}(t_2)
\int\limits_t^{t_2}
\phi_{j_1}(t_1)
C_{j_3}(t_4,t_2) C_{j_3}(T,t_5)dt_1 dt_2 dt_4 dt_5,
$$

\vspace{2mm}
$$
C_{j_4 j_5 j_4 j_3 j_2 j_1}=
\int\limits_t^T
\phi_{j_5}(t_5)
\int\limits_t^{t_5}
\phi_{j_3}(t_3)
\int\limits_t^{t_3}
\phi_{j_2}(t_2)
\int\limits_t^{t_2}
\phi_{j_1}(t_1)
C_{j_4}(t_5,t_3) C_{j_4}(T,t_5)dt_1 dt_2 dt_3 dt_5,
$$

\vspace{2mm}
$$
C_{j_5 j_5 j_4 j_3 j_2 j_1}=
\int\limits_t^T
\phi_{j_4}(t_4)
\int\limits_t^{t_4}
\phi_{j_3}(t_3)
\int\limits_t^{t_3}
\phi_{j_2}(t_2)
\int\limits_t^{t_2}
\phi_{j_1}(t_1)
C_{j_5 j_5}(T,t_4)dt_1 dt_2 dt_3 dt_4,
$$

\vspace{2mm}
$$
C_{j_6 j_5 j_4 j_3 j_1 j_1}=
\int\limits_t^T
\phi_{j_6}(t_6)
\int\limits_t^{t_6}
\phi_{j_5}(t_5)
\int\limits_t^{t_5}
\phi_{j_4}(t_4)
\int\limits_t^{t_4}
\phi_{j_3}(t_3)
C_{j_1 j_1}(t_3,t)dt_3 dt_4 dt_5 dt_6,
$$

\vspace{2mm}
$$
C_{j_6 j_5 j_4 j_1 j_2 j_1}=
\int\limits_t^T
\phi_{j_6}(t_6)
\int\limits_t^{t_6}
\phi_{j_5}(t_5)
\int\limits_t^{t_5}
\phi_{j_4}(t_4)
\int\limits_t^{t_4}
\phi_{j_2}(t_2)
C_{j_1}(t_2,t) C_{j_1}(t_4,t_2)dt_2 dt_4 dt_5 dt_6,
$$

\vspace{2mm}
$$
C_{j_6 j_5 j_1 j_3 j_2 j_1}=
\int\limits_t^T
\phi_{j_6}(t_6)
\int\limits_t^{t_6}
\phi_{j_5}(t_5)
\int\limits_t^{t_5}
\phi_{j_3}(t_3)
\int\limits_t^{t_3}
\phi_{j_2}(t_2)
C_{j_1}(t_2,t) C_{j_1}(t_5,t_3)dt_2 dt_3 dt_5 dt_6,
$$

\vspace{2mm}
$$
C_{j_6 j_1 j_4 j_3 j_2 j_1}=
\int\limits_t^T
\phi_{j_6}(t_6)
\int\limits_t^{t_6}
\phi_{j_4}(t_4)
\int\limits_t^{t_4}
\phi_{j_3}(t_3)
\int\limits_t^{t_3}
\phi_{j_2}(t_2)
C_{j_1}(t_2,t) C_{j_1}(t_6,t_4)dt_2 dt_3 dt_4 dt_6,
$$

\vspace{2mm}
$$
C_{j_6 j_5 j_4 j_2 j_2 j_1}=
\int\limits_t^T
\phi_{j_6}(t_6)
\int\limits_t^{t_6}
\phi_{j_5}(t_5)
\int\limits_t^{t_5}
\phi_{j_4}(t_4)
\int\limits_t^{t_4}
\phi_{j_1}(t_1)
C_{j_2 j_2}(t_4,t_1)dt_1 dt_4 dt_5 dt_6,
$$

\vspace{2mm}
$$
C_{j_6 j_5 j_2 j_3 j_2 j_1}=
\int\limits_t^T
\phi_{j_6}(t_6)
\int\limits_t^{t_6}
\phi_{j_5}(t_5)
\int\limits_t^{t_5}
\phi_{j_3}(t_3)
\int\limits_t^{t_3}
\phi_{j_1}(t_1)
C_{j_2}(t_3,t_1) C_{j_2}(t_5,t_3)dt_1 dt_3 dt_5 dt_6,
$$

\vspace{2mm}
$$
C_{j_6 j_2 j_4 j_3 j_2 j_1}=
\int\limits_t^T
\phi_{j_6}(t_6)
\int\limits_t^{t_6}
\phi_{j_4}(t_4)
\int\limits_t^{t_4}
\phi_{j_3}(t_3)
\int\limits_t^{t_3}
\phi_{j_1}(t_1)
C_{j_2}(t_3,t_1) C_{j_2}(t_6,t_4)dt_1 dt_3 dt_4 dt_6,
$$

\vspace{2mm}
$$
C_{j_6 j_5 j_3 j_3 j_2 j_1}=
\int\limits_t^T
\phi_{j_6}(t_6)
\int\limits_t^{t_6}
\phi_{j_5}(t_5)
\int\limits_t^{t_5}
\phi_{j_2}(t_2)
\int\limits_t^{t_2}
\phi_{j_1}(t_1)
C_{j_3 j_3}(t_5,t_2)dt_1 dt_2 dt_5 dt_6,
$$

\vspace{2mm}
$$
C_{j_6 j_3 j_4 j_3 j_2 j_1}=
\int\limits_t^T
\phi_{j_6}(t_6)
\int\limits_t^{t_6}
\phi_{j_4}(t_4)
\int\limits_t^{t_4}
\phi_{j_2}(t_2)
\int\limits_t^{t_2}
\phi_{j_1}(t_1)
C_{j_3}(t_4,t_2) C_{j_3}(t_6,t_4)dt_1 dt_2 dt_4 dt_6,
$$

\vspace{2mm}
$$
C_{j_6 j_4 j_4 j_3 j_2 j_1}=
\int\limits_t^T
\phi_{j_6}(t_6)
\int\limits_t^{t_6}
\phi_{j_3}(t_3)
\int\limits_t^{t_3}
\phi_{j_2}(t_2)
\int\limits_t^{t_2}
\phi_{j_1}(t_1)
C_{j_4 j_4}(t_6,t_3)dt_1 dt_2 dt_3 dt_6,
$$

\vspace{4mm}

2. $r=2$ (45 combinations)

\vspace{-1mm}
$$
C_{j_6 j_5 j_3 j_3 j_1 j_1}=
\int\limits_t^T
\phi_{j_6}(t_6)
\int\limits_t^{t_6}
\phi_{j_5}(t_5)
C_{j_3 j_3 j_1 j_1}(t_5,t)dt_5 dt_6,
$$

\vspace{2mm}
$$
C_{j_6 j_3 j_4 j_3 j_1 j_1}=
\int\limits_t^T
\phi_{j_6}(t_6)
\int\limits_t^{t_6}
\phi_{j_4}(t_4)
C_{j_3 j_1 j_1}(t_4,t)C_{j_3}(t_6,t_4)dt_4 dt_6,
$$

\vspace{2mm}
$$
C_{j_6 j_4 j_4 j_3 j_1 j_1}=
\int\limits_t^T
\phi_{j_6}(t_6)
\int\limits_t^{t_6}
\phi_{j_3}(t_3)
C_{j_1 j_1}(t_3,t)C_{j_4 j_4}(t_6,t_3)dt_3 dt_6,
$$

\vspace{2mm}
$$
C_{j_6 j_5 j_2 j_1 j_2 j_1}=
\int\limits_t^T
\phi_{j_6}(t_6)
\int\limits_t^{t_6}
\phi_{j_5}(t_5)
C_{j_2 j_1 j_2 j_1}(t_5,t)dt_5 dt_6,
$$

\vspace{2mm}
$$
C_{j_6 j_2 j_4 j_1 j_2 j_1}=
\int\limits_t^T
\phi_{j_6}(t_6)
\int\limits_t^{t_6}
\phi_{j_4}(t_4)
C_{j_1 j_2 j_1}(t_4,t)C_{j_2}(t_6,t_4)dt_4 dt_6,
$$

\vspace{2mm}
$$
C_{j_6 j_4 j_4 j_1 j_2 j_1}=
\int\limits_t^T
\phi_{j_6}(t_6)
\int\limits_t^{t_6}
\phi_{j_2}(t_2)
C_{j_1}(t_2,t)C_{j_4 j_4 j_1}(t_6,t_2)dt_2 dt_6,
$$

\vspace{2mm}
$$
C_{j_6 j_5 j_1 j_2 j_2 j_1}=
\int\limits_t^T
\phi_{j_6}(t_6)
\int\limits_t^{t_6}
\phi_{j_5}(t_5)
C_{j_1 j_2 j_2 j_1}(t_5,t)dt_5 dt_6,
$$

\vspace{2mm}
$$
C_{j_6 j_2 j_1 j_3 j_2 j_1}=
\int\limits_t^T
\phi_{j_6}(t_6)
\int\limits_t^{t_6}
\phi_{j_3}(t_3)
C_{j_2 j_1}(t_3,t)C_{j_2 j_1}(t_6,t_3)dt_3 dt_6,
$$

\vspace{2mm}
$$
C_{j_6 j_3 j_1 j_3 j_2 j_1}=
\int\limits_t^T
\phi_{j_6}(t_6)
\int\limits_t^{t_6}
\phi_{j_2}(t_2)
C_{j_1}(t_2,t)C_{j_3 j_1 j_3}(t_6,t_2)dt_2 dt_6,
$$

\vspace{2mm}
$$
C_{j_6 j_1 j_4 j_2 j_2 j_1}=
\int\limits_t^T
\phi_{j_6}(t_6)
\int\limits_t^{t_6}
\phi_{j_4}(t_4)
C_{j_2 j_2 j_1}(t_4,t)C_{j_1}(t_6,t_4)dt_4 dt_6,
$$

\vspace{2mm}
$$
C_{j_6 j_1 j_2 j_3 j_2 j_1}=
\int\limits_t^T
\phi_{j_6}(t_6)
\int\limits_t^{t_6}
\phi_{j_3}(t_3)
C_{j_2 j_1}(t_3,t)C_{j_1 j_2}(t_6,t_3)dt_3 dt_6,
$$

\vspace{2mm}
$$
C_{j_6 j_1 j_3 j_3 j_2 j_1}=
\int\limits_t^T
\phi_{j_6}(t_6)
\int\limits_t^{t_6}
\phi_{j_2}(t_2)
C_{j_1}(t_2,t)C_{j_1 j_3 j_3}(t_6,t_2)dt_2 dt_6,
$$

\vspace{2mm}
$$
C_{j_6 j_4 j_4 j_2 j_2 j_1}=
\int\limits_t^T
\phi_{j_6}(t_6)
\int\limits_t^{t_6}
\phi_{j_1}(t_1)
C_{j_4 j_4 j_2 j_2}(t_6,t_1)dt_1 dt_6,
$$

\vspace{2mm}
$$
C_{j_6 j_3 j_2 j_3 j_2 j_1}=
\int\limits_t^T
\phi_{j_6}(t_6)
\int\limits_t^{t_6}
\phi_{j_1}(t_1)
C_{j_3 j_2 j_3 j_2}(t_6,t_1)dt_1 dt_6,
$$

\vspace{2mm}
$$
C_{j_6 j_2 j_3 j_3 j_2 j_1}=
\int\limits_t^T
\phi_{j_6}(t_6)
\int\limits_t^{t_6}
\phi_{j_1}(t_1)
C_{j_2 j_3 j_3 j_2}(t_6,t_1)dt_1 dt_6,
$$

\vspace{2mm}
$$
C_{j_1 j_5 j_3 j_3 j_2 j_1}=
\int\limits_t^T
\phi_{j_5}(t_5)
\int\limits_t^{t_5}
\phi_{j_2}(t_2)
C_{j_1}(t_2,t)C_{j_3 j_3}(t_5,t_2)C_{j_1}(T,t_5)dt_2 dt_5,
$$

\vspace{2mm}
$$
C_{j_1 j_3 j_4 j_3 j_2 j_1}=
\int\limits_t^T
\phi_{j_4}(t_4)
\int\limits_t^{t_4}
\phi_{j_2}(t_2)
C_{j_1}(t_2,t)C_{j_3}(t_4,t_2)C_{j_1 j_3}(T,t_4)dt_2 dt_4,
$$

\vspace{2mm}
$$
C_{j_1 j_2 j_4 j_3 j_2 j_1}=
\int\limits_t^T
\phi_{j_4}(t_4)
\int\limits_t^{t_4}
\phi_{j_3}(t_3)
C_{j_2 j_1}(t_3,t)C_{j_1 j_2}(T,t_4)dt_3 dt_4,
$$

\vspace{2mm}
$$
C_{j_1 j_5 j_2 j_3 j_2 j_1}=
\int\limits_t^T
\phi_{j_5}(t_5)
\int\limits_t^{t_5}
\phi_{j_3}(t_3)
C_{j_2 j_1}(t_3,t)C_{j_2}(t_5,t_3)C_{j_1}(T,t_5)dt_3 dt_5,
$$

\vspace{2mm}
$$
C_{j_1 j_4 j_4 j_3 j_2 j_1}=
\int\limits_t^T
\phi_{j_3}(t_3)
\int\limits_t^{t_3}
\phi_{j_2}(t_2)
C_{j_1}(t_2,t)C_{j_1 j_4 j_4}(T,t_3)dt_2 dt_3,
$$

\vspace{2mm}
$$
C_{j_1 j_5 j_4 j_2 j_2 j_1}=
\int\limits_t^T
\phi_{j_5}(t_5)
\int\limits_t^{t_5}
\phi_{j_4}(t_4)
C_{j_2 j_2 j_1}(t_4,t)C_{j_1}(T,t_5)dt_4 dt_5,
$$

\vspace{2mm}
$$
C_{j_2 j_3 j_4 j_3 j_2 j_1}=
\int\limits_t^T
\phi_{j_4}(t_4)
\int\limits_t^{t_4}
\phi_{j_1}(t_1)
C_{j_2 j_3}(t_4,t_1)C_{j_2 j_3}(T,t_4)dt_1 dt_4,
$$

\vspace{2mm}
$$
C_{j_2 j_4 j_4 j_3 j_2 j_1}=
\int\limits_t^T
\phi_{j_3}(t_3)
\int\limits_t^{t_3}
\phi_{j_1}(t_1)
C_{j_2}(t_3,t_1)C_{j_2 j_4 j_4}(T,t_3)dt_1 dt_3,
$$

\vspace{2mm}
$$
C_{j_2 j_5 j_3 j_3 j_2 j_1}=
\int\limits_t^T
\phi_{j_5}(t_5)
\int\limits_t^{t_5}
\phi_{j_1}(t_1)
C_{j_2 j_3 j_3}(t_5,t_1)C_{j_2}(T,t_5)dt_1 dt_5,
$$

\vspace{2mm}
$$
C_{j_2 j_1 j_4 j_3 j_2 j_1}=
\int\limits_t^T
\phi_{j_4}(t_4)
\int\limits_t^{t_4}
\phi_{j_3}(t_3)
C_{j_2 j_1}(t_3,t)C_{j_2 j_1}(T,t_4)dt_3 dt_4,
$$

\vspace{2mm}
$$
C_{j_2 j_5 j_1 j_3 j_2 j_1}=
\int\limits_t^T
\phi_{j_5}(t_5)
\int\limits_t^{t_5}
\phi_{j_3}(t_3)
C_{j_2 j_1}(t_3,t)C_{j_1}(t_5,t_3)C_{j_2}(T,t_5)dt_3 dt_5,
$$

\vspace{2mm}
$$
C_{j_2 j_5 j_4 j_1 j_2 j_1}=
\int\limits_t^T
\phi_{j_5}(t_5)
\int\limits_t^{t_5}
\phi_{j_4}(t_4)
C_{j_1 j_2 j_1}(t_4,t)C_{j_2}(T,t_5)dt_4 dt_5,
$$

\vspace{2mm}
$$
C_{j_3 j_2 j_4 j_3 j_2 j_1}=
\int\limits_t^T
\phi_{j_4}(t_4)
\int\limits_t^{t_4}
\phi_{j_1}(t_1)
C_{j_3 j_2}(t_4,t_1)C_{j_3 j_2}(T,t_4)dt_1 dt_4,
$$

\vspace{2mm}
$$
C_{j_3 j_4 j_4 j_3 j_2 j_1}=
\int\limits_t^T
\phi_{j_2}(t_2)
\int\limits_t^{t_2}
\phi_{j_1}(t_1)
C_{j_3 j_4 j_4 j_3}(T,t_2)dt_1 dt_2,
$$

\vspace{2mm}
$$
C_{j_3 j_5 j_2 j_3 j_2 j_1}=
\int\limits_t^T
\phi_{j_5}(t_5)
\int\limits_t^{t_5}
\phi_{j_1}(t_1)
C_{j_2 j_3 j_2}(t_5,t_1)C_{j_3}(T,t_5)dt_1 dt_5,
$$

\vspace{2mm}
$$
C_{j_3 j_1 j_4 j_3 j_2 j_1}=
\int\limits_t^T
\phi_{j_4}(t_4)
\int\limits_t^{t_4}
\phi_{j_2}(t_2)
C_{j_1}(t_2,t)C_{j_3}(t_4,t_2)C_{j_3 j_1}(T,t_4)dt_2 dt_4,
$$

\vspace{2mm}
$$
C_{j_3 j_5 j_1 j_3 j_2 j_1}=
\int\limits_t^T
\phi_{j_5}(t_5)
\int\limits_t^{t_5}
\phi_{j_2}(t_2)
C_{j_1}(t_2,t)C_{j_1 j_3}(t_5,t_2)C_{j_3}(T,t_5)dt_2 dt_5,
$$

\vspace{2mm}
$$
C_{j_3 j_5 j_4 j_3 j_1 j_1}=
\int\limits_t^T
\phi_{j_5}(t_5)
\int\limits_t^{t_5}
\phi_{j_4}(t_4)
C_{j_3 j_1 j_1}(t_4,t)C_{j_3}(T,t_5)dt_4 dt_5,
$$

\vspace{2mm}
$$
C_{j_4 j_3 j_4 j_3 j_2 j_1}=
\int\limits_t^T
\phi_{j_2}(t_2)
\int\limits_t^{t_2}
\phi_{j_1}(t_1)
C_{j_4 j_3 j_4 j_3}(T,t_2)dt_1 dt_2,
$$

\vspace{2mm}
$$
C_{j_4 j_2 j_4 j_3 j_2 j_1}=
\int\limits_t^T
\phi_{j_3}(t_3)
\int\limits_t^{t_3}
\phi_{j_1}(t_1)
C_{j_2}(t_3,t_1)C_{j_4 j_2 j_4}(T,t_3)dt_1 dt_3,
$$

\vspace{2mm}
$$
C_{j_4 j_5 j_4 j_2 j_2 j_1}=
\int\limits_t^T
\phi_{j_5}(t_5)
\int\limits_t^{t_5}
\phi_{j_1}(t_1)
C_{j_4 j_2 j_2}(t_5,t_1)C_{j_4}(T,t_5)dt_1 dt_5,
$$

\vspace{2mm}
$$
C_{j_4 j_1 j_4 j_3 j_2 j_1}=
\int\limits_t^T
\phi_{j_3}(t_3)
\int\limits_t^{t_3}
\phi_{j_2}(t_2)
C_{j_1}(t_2,t)C_{j_4 j_1 j_4}(T,t_3)dt_2 dt_3,
$$

\vspace{2mm}
$$
C_{j_4 j_5 j_4 j_1 j_2 j_1}=
\int\limits_t^T
\phi_{j_5}(t_5)
\int\limits_t^{t_5}
\phi_{j_2}(t_2)
C_{j_1}(t_2,t)C_{j_4 j_1}(t_5,t_2)C_{j_4}(T,t_5)dt_2 dt_5,
$$

\vspace{2mm}
$$
C_{j_4 j_5 j_4 j_3 j_1 j_1}=
\int\limits_t^T
\phi_{j_5}(t_5)
\int\limits_t^{t_5}
\phi_{j_3}(t_3)
C_{j_1 j_1}(t_3,t)C_{j_4}(t_5,t_3)C_{j_4}(T,t_5)dt_3 dt_5,
$$

\vspace{2mm}
$$
C_{j_5 j_5 j_3 j_3 j_2 j_1}=
\int\limits_t^T
\phi_{j_2}(t_2)
\int\limits_t^{t_2}
\phi_{j_1}(t_1)
C_{j_5 j_5 j_3 j_3}(T,t_2)dt_1 dt_2,
$$

\vspace{2mm}
$$
C_{j_5 j_5 j_2 j_3 j_2 j_1}=
\int\limits_t^T
\phi_{j_3}(t_3)
\int\limits_t^{t_3}
\phi_{j_1}(t_1)
C_{j_2}(t_3,t_1)C_{j_5 j_5 j_2}(T,t_3)dt_1 dt_3,
$$

\vspace{2mm}
$$
C_{j_5 j_5 j_4 j_2 j_2 j_1}=
\int\limits_t^T
\phi_{j_4}(t_4)
\int\limits_t^{t_4}
\phi_{j_1}(t_1)
C_{j_2 j_2}(t_4,t_1)C_{j_5 j_5}(T,t_4)dt_1 dt_4,
$$

\vspace{2mm}
$$
C_{j_5 j_5 j_1 j_3 j_2 j_1}=
\int\limits_t^T
\phi_{j_3}(t_3)
\int\limits_t^{t_3}
\phi_{j_2}(t_2)
C_{j_1}(t_2,t)C_{j_5 j_5 j_1}(T,t_3)dt_2 dt_3,
$$

\vspace{2mm}
$$
C_{j_5 j_5 j_4 j_1 j_2 j_1}=
\int\limits_t^T
\phi_{j_4}(t_4)
\int\limits_t^{t_4}
\phi_{j_2}(t_2)
C_{j_1}(t_2,t)C_{j_1}(t_4,t_2)C_{j_5 j_5}(T,t_4)dt_2 dt_4,
$$

\vspace{2mm}
$$
C_{j_5 j_5 j_4 j_3 j_1 j_1}=
\int\limits_t^T
\phi_{j_4}(t_4)
\int\limits_t^{t_4}
\phi_{j_3}(t_3)
C_{j_1 j_1}(t_3,t)C_{j_5 j_5}(T,t_4)dt_3 dt_4.
$$

\vspace{4mm}

It is not difficult to see (based on the above equalities)
that the condition (\ref{09091}) will be satisfied under the conditions
of Theorem~60 if

\vspace{-2mm}
\begin{equation}
\label{2024december12}
\left\vert \sum\limits_{j_1=0}^p 
C_{j_1 j_1}(s,\tau)\right\vert\le K,
\end{equation}
\begin{equation}
\label{2024december13}
\left\vert \sum\limits_{j_1=0}^p 
C_{j_1}(s,\tau)C_{j_1}(\theta,u)\right\vert\le K,
\end{equation}
\begin{equation}
\label{2024december1}
\left\vert \sum\limits_{j_1, j_2=0}^p C_{j_2 j_2 j_1 j_1}(s,\tau)\right\vert\le K,
\end{equation}
\begin{equation}
\label{2024december2}
\left\vert \sum\limits_{j_1, j_2=0}^p C_{j_2 j_1 j_2 j_1}(s,\tau)\right\vert\le K,
\end{equation}
\begin{equation}
\label{2024december3}
\left\vert \sum\limits_{j_1, j_2=0}^p C_{j_1 j_2 j_2 j_1}(s,\tau)\right\vert\le K,
\end{equation}
\begin{equation}
\label{2024december4}
\left\vert \sum\limits_{j_1, j_2=0}^p C_{j_2 j_1 j_1}(s,\tau)C_{j_2}(\theta,u)\right\vert\le K,
\end{equation}
\begin{equation}
\label{2024december5}
\left\vert \sum\limits_{j_1, j_2=0}^p C_{j_1 j_2 j_1}(s,\tau)C_{j_2}(\theta,u)\right\vert\le K,
\end{equation}
\begin{equation}
\label{2024december6}
\left\vert \sum\limits_{j_1, j_2=0}^p C_{j_2 j_2 j_1}(s,\tau)C_{j_1}(\theta,u)\right\vert\le K,
\end{equation}
\begin{equation}
\label{2024december7}
\left\vert \sum\limits_{j_1, j_2=0}^p C_{j_1 j_1}(s,\tau)C_{j_2 j_2}(\theta,u)\right\vert\le K,
\end{equation}
\begin{equation}
\label{2024december8}
\left\vert \sum\limits_{j_1, j_2=0}^p C_{j_2 j_1}(s,\tau)C_{j_2 j_1}(\theta,u)\right\vert\le K,
\end{equation}
\begin{equation}
\label{2024december9}
\left\vert \sum\limits_{j_1, j_2=0}^p C_{j_2 j_1}(s,\tau)C_{j_1 j_2}(\theta,u)\right\vert\le K,
\end{equation}
\begin{equation}
\label{2024december10}
\left\vert \sum\limits_{j_1, j_2=0}^p C_{j_1}(s,\tau)C_{j_1}(\rho,v)
C_{j_2 j_2}(\theta,u)\right\vert\le K,
\end{equation}
\begin{equation}
\label{2024december11}
\left\vert \sum\limits_{j_1, j_2=0}^p C_{j_1}(s,\tau)C_{j_2}(\rho,v)
C_{j_1 j_2}(\theta,u)\right\vert\le K,
\end{equation}

\vspace{2mm}
\noindent
where $p\in\mathbb{N},$ $t\le \tau < s \le T,$ $t\le u<\theta \le T,$
$t\le v<\rho \le T,$ constant $K$ does not depend on 
$p, s, \tau, u, \theta, v, \rho$ (but only on $t, T$) and may differ from line to line.

The equalities (\ref{2024december1})--(\ref{2024december3}) 
have been proved earlier (see (\ref{april50})--(\ref{april52})).

Using Fubini's Theorem and Parseval's equality, we get

\vspace{-1mm}
$$
\left\vert \sum\limits_{j_1=0}^p 
C_{j_1 j_1}(s,\tau)\right\vert=
\frac{1}{2}\sum\limits_{j_1=0}^p C_{j_1}^2(s,\tau)\le
\frac{1}{2}\sum\limits_{j_1=0}^{\infty}
C_{j_1}^2(s,\tau)=
\frac{1}{2}(s-\tau) \le
\frac{1}{2}(T-t) \le K.
$$

\vspace{3mm}
\noindent
The equality (\ref{2024december12}) is proved.
Moreover, (\ref{2024december7}) follows from 
(\ref{2024december12}).

Using the inequality of Cauchy--Bunyakovsky and 
Parseval's equality, we obtain

\vspace{-1mm}
$$
\left(\sum\limits_{j_1=0}^p 
C_{j_1}(s,\tau)C_{j_1}(\theta,u)\right)^2\le
\sum\limits_{j_1=0}^p 
C_{j_1}^2(s,\tau)
\sum\limits_{j_1=0}^p 
C_{j_1}^2(\theta,u)\le
$$

\vspace{1mm}
$$
\le \sum\limits_{j_1=0}^{\infty} 
C_{j_1}^2(s,\tau)
\sum\limits_{j_1=0}^{\infty} 
C_{j_1}^2(\theta,u)=
(s-\tau)(\theta-u)\le 
(T-t)^2 \le K^2,
$$

\vspace{3mm}
$$
\left(\sum\limits_{j_1, j_2=0}^p C_{j_2 j_1}(s,\tau)C_{j_2 j_1}(\theta,u)\right)^2
\le
\sum\limits_{j_1,j_2=0}^p 
C_{j_2 j_1}^2(s,\tau)
\sum\limits_{j_1,j_2=0}^p 
C_{j_2 j_1}^2(\theta,u)\le
$$

\vspace{1mm}
$$
\le
\sum\limits_{j_1,j_2=0}^{\infty}
C_{j_2 j_1}^2(s,\tau)
\sum\limits_{j_1,j_2=0}^{\infty} 
C_{j_2 j_1}^2(\theta,u)=
\int\limits_{\tau}^s \int\limits_{\tau}^v dx dv
\int\limits_{u}^{\theta} \int\limits_{u}^v dx dv\le 
\frac{1}{4}(T-t)^4\le K^2.
$$

\vspace{4mm}

Thus, the inequalities (\ref{2024december13}), (\ref{2024december8}) are proved.
The inequalities (\ref{2024december9}), (\ref{2024december11}) are 
proved similarly to (\ref{2024december8}).
Moreover, (\ref{2024december10}) follows from
(\ref{2024december12}), (\ref{2024december13}).

Further, let us prove the equalities (\ref{2024december4})--(\ref{2024december6}).
Applying the 
Cau\-chy--Bunyakovsky inequality as well as 
Parseval's equality and (\ref{2024december12}), we have

\vspace{-1mm}
$$
\left(\sum\limits_{j_1, j_2=0}^p C_{j_2 j_1 j_1}(s,\tau)C_{j_2}(\theta,u)\right)^2
\le 
\sum\limits_{j_2=0}^p \left(\sum\limits_{j_1=0}^p C_{j_2 j_1 j_1}(s,\tau)\right)^2
\sum\limits_{j_2=0}^p C_{j_2}^2(\theta,u)\le
$$

\vspace{1mm}
$$
\le \sum\limits_{j_2=0}^{\infty} \left(\sum\limits_{j_1=0}^p C_{j_2 j_1 j_1}(s,\tau)\right)^2
\sum\limits_{j_2=0}^{\infty} C_{j_2}^2(\theta,u)=
\sum\limits_{j_2=0}^{\infty} \left(
\int\limits_{\tau}^s \phi_{j_2}(v) \sum\limits_{j_1=0}^p  C_{j_1 j_1}(v,\tau) dv
\right)^2 \cdot (\theta-u)=
$$

\vspace{1mm}
$$
=
(\theta-u)\int\limits_{\tau}^s \left(\sum\limits_{j_1=0}^p  C_{j_1 j_1}(v,\tau)\right)^2 dv
\le 
K^2 (\theta-u)(s-\tau)\le K^2(T-t)^2=K_1.
$$

\vspace{4mm}
\noindent
The equality (\ref{2024december4}) is proved.

Using the 
Cau\-chy--Bunyakovsky inequality as well as 
Fubini's Theorem, Parseval's equality and (\ref{2024december13}), we have

\vspace{-1mm}
$$
\left(\sum\limits_{j_1, j_2=0}^p C_{j_1 j_2 j_1}(s,\tau)C_{j_2}(\theta,u)\right)^2
\le 
\sum\limits_{j_2=0}^p \left(\sum\limits_{j_1=0}^p C_{j_1 j_2 j_1}(s,\tau)\right)^2
\sum\limits_{j_2=0}^p C_{j_2}^2(\theta,u)\le
$$

\vspace{1mm}
$$
\le \sum\limits_{j_2=0}^{\infty} \left(\sum\limits_{j_1=0}^p 
\int\limits_{\tau}^s \phi_{j_1}(z)\int\limits_{\tau}^{z} \phi_{j_2}(y)
\int\limits_{\tau}^{y} \phi_{j_1}(x)dx dy dz\right)^2
\sum\limits_{j_2=0}^{\infty} C_{j_2}^2(\theta,u)=
$$

\vspace{1mm}
$$
=\sum\limits_{j_2=0}^{\infty} \left(\sum\limits_{j_1=0}^p 
\int\limits_{\tau}^s 
\phi_{j_2}(y)
\int\limits_{\tau}^{y} \phi_{j_1}(x)dx
\int\limits_{y}^{s}
\phi_{j_1}(z) dz dy\right)^2
\cdot (\theta-u)=
$$

\vspace{1mm}
$$
=(\theta-u)\sum\limits_{j_2=0}^{\infty} \left(
\int\limits_{\tau}^s 
\phi_{j_2}(y)
\sum\limits_{j_1=0}^p C_{j_1}(y,\tau) C_{j_1}(s,y) dy\right)^2=
$$

\vspace{1mm}
$$
=(\theta-u)
\int\limits_{\tau}^s 
\left(\sum\limits_{j_1=0}^p C_{j_1}(y,\tau) C_{j_1}(s,y)\right)^2 dy
\le 
$$

\vspace{2mm}
$$
\le K^2 (\theta-u)(s-\tau)\le K^2(T-t)^2=K_1.
$$

\vspace{4mm}
\noindent
The equality (\ref{2024december5}) is proved.

Using the 
Cau\-chy--Bunyakovsky inequality as well as 
Fubini's Theorem, Parseval's equality and (\ref{2024december12}), we have

\vspace{-1mm}
$$
\left(\sum\limits_{j_1, j_2=0}^p C_{j_2 j_2 j_1}(s,\tau)C_{j_1}(\theta,u)\right)^2\le
\sum\limits_{j_1=0}^p \left(\sum\limits_{j_2=0}^p C_{j_2 j_2 j_1}(s,\tau)\right)^2
\sum\limits_{j_1=0}^p C_{j_1}^2(\theta,u)\le
$$

\vspace{1mm}
$$
\le \sum\limits_{j_1=0}^{\infty} \left(\sum\limits_{j_2=0}^p 
\int\limits_{\tau}^s \phi_{j_2}(z)\int\limits_{\tau}^{z} \phi_{j_2}(y)
\int\limits_{\tau}^{y} \phi_{j_1}(x)dx dy dz\right)^2
\sum\limits_{j_1=0}^{\infty} C_{j_1}^2(\theta,u)=
$$

\vspace{1mm}
$$
=\sum\limits_{j_1=0}^{\infty} \left(\sum\limits_{j_2=0}^p 
\int\limits_{\tau}^s \phi_{j_1}(x)\int\limits_{x}^{s} \phi_{j_2}(y)
\int\limits_{y}^{s} \phi_{j_2}(z)dz dy dx\right)^2
\cdot (\theta-u)=
$$

\vspace{1mm}
$$
=(\theta-u)\sum\limits_{j_1=0}^{\infty} \left(\sum\limits_{j_2=0}^p 
\int\limits_{\tau}^s \phi_{j_1}(x)\int\limits_{x}^{s} \phi_{j_2}(z)
\int\limits_{x}^{z} \phi_{j_2}(y)dy dz dx\right)^2=
$$

\vspace{1mm}
$$
=(\theta-u)\sum\limits_{j_1=0}^{\infty} \left(
\int\limits_{\tau}^s \phi_{j_1}(x)  
\sum\limits_{j_2=0}^p C_{j_2 j_2}(s,x)
dx\right)^2=
$$

\vspace{1mm}
$$
=(\theta-u)
\int\limits_{\tau}^s \left(
\sum\limits_{j_2=0}^p C_{j_2 j_2}(s,x)\right)^2
dx\le
$$

\vspace{2mm}
$$
\le K^2 (\theta-u)(s-\tau)\le K^2(T-t)^2=K_1.
$$

\vspace{4mm}
\noindent
The equality (\ref{2024december6}) is proved.
The equalities (\ref{2024december12})--(\ref{2024december11})
are proved.

Thus, the condition (\ref{09091}) of Theorem~59 is satisfied 
under the conditions of Theorem~60. The assertion
of Theorem~60 now follows from Theorem~59. 
Theorem~60 is proved.

\vspace{5mm}

\section{Expansion of Iterated Stratonovich Stochastic Integrals
of Multiplicity 4. The Case of an Ar\-bit\-ra\-ry Complete Orthonormal System of 
Functions in the Space $L_2([t,T])$ and Binomial Weight Functions}

\vspace{5mm}

Let us prove the following theorem.

\vspace{2mm}   

{\bf Theorem~61}\ \cite{2018a}.\ {\it Suppose that
$\{\phi_j(x)\}_{j=0}^{\infty}$ is an arbitrary complete orthonormal system of 
func\-ti\-ons in the space $L_2([t,T]).$
Then$,$ for the iterated Stra\-to\-no\-vich stochastic integral
of fourth multiplicity 

\vspace{-2mm}
$$
I_{{l_1l_2l_3 l_4}_{T,t}}^{*(i_1i_2i_3 i_4)}=
{\int\limits_t^{*}}^T (t_4-t)^{l_4}{\int\limits_t^{*}}^{t_4} (t_3-t)^{l_3}
{\int\limits_t^{*}}^{t_3}(t_2-t)^{l_2}
{\int\limits_t^{*}}^{t_2}(t_1-t)^{l_1}
d{\bf w}_{t_1}^{(i_1)}
d{\bf w}_{t_2}^{(i_2)}d{\bf w}_{t_3}^{(i_3)}d{\bf w}_{t_4}^{(i_4)}
$$

\vspace{3mm}
\noindent
the following expansion 

\vspace{-1mm}
$$
I_{{l_1l_2l_3 l_4}_{T,t}}^{*(i_1i_2i_3i_4)}=
\hbox{\vtop{\offinterlineskip\halign{
\hfil#\hfil\cr
{\rm l.i.m.}\cr
$\stackrel{}{{}_{p\to \infty}}$\cr
}} }\sum_{j_1,j_2,j_3,j_4=0}^{p}
C_{j_4 j_3 j_2 j_1}\zeta_{j_1}^{(i_1)}\zeta_{j_2}^{(i_2)}\zeta_{j_3}^{(i_3)}\zeta_{j_4}^{(i_4)}
$$

\vspace{3mm}
\noindent
that converges in the mean-square sense is valid, where 
$i_1,i_2,i_3,i_4=0,1,\ldots,m;$ $l_1,l_2,l_3,l_4=0,1,2,\ldots,$
\begin{equation}
\label{2024decem1}
C_{j_4 j_3 j_2 j_1}=\int\limits_t^T
(t_4-t)^{l_4}\phi_{j_4}(t_4)\int\limits_t^{t_4}
(t_3-t)^{l_3}\phi_{j_3}(t_3)\int\limits_t^{t_3}
(t_2-t)^{l_2}
\phi_{j_2}(t_2)
\int\limits_t^{t_2}
(t_1-t)^{l_1}\phi_{j_1}(t_1)
dt_1dt_2dt_3dt_4
\end{equation}

\vspace{1mm}
\noindent
and
$$
\zeta_{j}^{(i)}=
\int\limits_t^T \phi_{j}(\tau) d{\bf w}_{\tau}^{(i)}
$$ 

\vspace{2mm}
\noindent
are independent standard Gaussian random variables for various 
$i$ or $j$ {\rm (}in the case when $i\ne 0${\rm ),}
${\bf w}_{\tau}^{(i)}={\bf f}_{\tau}^{(i)}$ for
$i=1,\ldots,m$ and 
${\bf w}_{\tau}^{(0)}=\tau.$}

\vspace{2mm}

{\bf Proof.}\ The following proof will be based on Theorem~59 
and verification of the equality (\ref{09091}) under the conditions
of Theorem~61 (the case $k=4>2r$, where $r=1$). Note that the case
$k=2r$ is proved in Sect.~26 (see (\ref{july90000})).
Under the conditions of Theorem~61, the equality $k=2r$ means that
$k=4$ and $r=2$.

Let throughout this proof 

\vspace{-1mm}
$$
C_{j_1 j_1}^{\psi_{i+1}\psi_i}(s,\tau)=\int\limits_{\tau}^s 
\psi_{i+1}(y)\phi_{j_1}(y)
\int\limits_{\tau}^{y}
\psi_i(x)\phi_{j_1}(x)dx dy,
$$

\vspace{1mm}
$$
C_{j_1}^{\psi_q}(s,\tau)=\int\limits_{\tau}^s 
\psi_{q}(x)\phi_{j_1}(x)dx,
$$

\vspace{3mm}
\noindent
where $i=1,2,3,$ $t\le\tau<s\le T,$ $\psi_q(x)=(x-t)^{l_q},$ $l_q=0,1,2,\ldots,$ 
$q=1,\ldots,4,$ $x\in [t, T],$
and $C_{j_4 j_3 j_2 j_1}$ is defined by (\ref{2024decem1}).

Using Fubini's Theorem and the technique that leads to the formulas (\ref{copa1}),
(\ref{copa1a}),
we obtain (note that we find all possible
combinations of pairs using the equality (\ref{after34})):

$$
C_{j_4 j_3 j_1 j_1}=
\int\limits_t^T
\psi_4(t_4)\phi_{j_4}(t_4)\int\limits_t^{t_4}
\psi_3(t_3)\phi_{j_3}(t_3)C_{j_1 j_1}^{\psi_2 \psi_1}(t_3,t)dt_3 dt_4,
$$

\vspace{2mm}
$$
C_{j_4 j_1 j_2 j_1}=
\int\limits_t^T
\psi_4(t_4)\phi_{j_4}(t_4)\int\limits_t^{t_4}
\psi_2(t_2)\phi_{j_2}(t_2)C_{j_1}^{\psi_1}(t_2,t)C_{j_1}^{\psi_3}(t_4,t_2)dt_2 dt_4,
$$

\vspace{2mm}
$$
C_{j_1 j_3 j_2 j_1}=
\int\limits_t^T
\psi_3(t_3)\phi_{j_3}(t_3)\int\limits_t^{t_3}
\psi_2(t_2)\phi_{j_2}(t_2)C_{j_1}^{\psi_1}(t_2,t)C_{j_1}^{\psi_4}(T,t_3)dt_2 dt_3,
$$

\vspace{2mm}
$$
C_{j_4 j_2 j_2 j_1}=
\int\limits_t^T
\psi_4(t_4)\phi_{j_4}(t_4)\int\limits_t^{t_4}
\psi_1(t_1)\phi_{j_1}(t_1)C_{j_2 j_2}^{\psi_3 \psi_2}(t_4,t_1)dt_1 dt_4,
$$

\vspace{2mm}
$$
C_{j_2 j_3 j_2 j_1}=
\int\limits_t^T
\psi_3(t_3)\phi_{j_3}(t_3)\int\limits_t^{t_3}
\psi_1(t_1)\phi_{j_1}(t_1)C_{j_2}^{\psi_2}(t_3,t_1) C_{j_2}^{\psi_4}(T,t_3)dt_1 dt_3,
$$

\vspace{2mm}
$$
C_{j_3 j_3 j_1 j_1}=
\int\limits_t^T
\psi_2(t_2)\phi_{j_2}(t_2)\int\limits_t^{t_2}
\psi_1(t_1)\phi_{j_1}(t_1)C_{j_3 j_3}^{\psi_4 \psi_3}(T,t_2)dt_1 dt_2.
$$

\vspace{4mm}

It is easy to see (based on the above equalities)
that the condition (\ref{09091}) will be satisfied under the conditions
of Theorem~61 if

\vspace{-1mm}
\begin{equation}
\label{2024decem12}
\left\vert \sum\limits_{j_1=0}^p 
C_{j_1 j_1}^{\psi_{i+1} \psi_i}(s,\tau)\right\vert\le K,
\end{equation}
\begin{equation}
\label{2024decem13}
\left\vert \sum\limits_{j_1=0}^p 
C_{j_1}^{\psi_k}(s,\tau)C_{j_1}^{\psi_q}(\theta,u)\right\vert\le K,
\end{equation}

\vspace{2mm}
\noindent
where $p\in\mathbb{N},$ $i=1,2,3,$ $k,q=1,\ldots,4,$ $t\le \tau < s \le T,$ $t\le u<\theta \le T,$
constant $K$ does not depend on 
$p, s, \tau, u, \theta$ (but only on $t, T$). 

The equality (\ref{2024decem12}) has been proved earlier (see (\ref{may103x})).
Obviously,  the relation (\ref{2024decem13}) is proved in complete
analogy with (\ref{may106}).

Thus, the condition (\ref{09091}) of Theorem~59 is fulfilled
under the conditions of Theorem~61. Then
Theorem~61 follows from Theorem~59. 
Theorem~61 is proved.

\vspace{5mm}

\section{Another Proof of Theorem~48 Based on Theorem~59}

\vspace{5mm}

The following proof will be based on Theorem~59 
and verification of the equality (\ref{09091}) under the conditions
of Theorem~48 (the case $k=5>2r$, where $r=1$ or $r=2$). 

Further, suppose that

\vspace{-2mm}
$$
C_{j_k \ldots j_1}(s,\tau)=\int\limits_{\tau}^s
\phi_{j_k}(t_k)\ldots
\int\limits_{\tau}^{t_2}
\phi_{j_1}(t_1)dt_1\ldots dt_k, 
$$

\vspace{3mm}
\noindent
where $k=1,\ldots,4,$ $t\le\tau<s\le T$, and 

\vspace{-1mm}
$$
C_{j_5\ldots j_1}=\int\limits_t^T
\phi_{j_5}(t_5)
\ldots
\int\limits_t^{t_2}
\phi_{j_1}(t_1)dt_1\ldots dt_5.
$$

\vspace{3mm}

Applying the technique that leads to (\ref{copa1}), we obtain 
(note that we find all possible
combinations of pairs using the equality (\ref{after35}))

$$
C_{j_5 j_4 j_3 j_1 j_1}
=\int\limits_t^T\phi_{j_5}(t_5)\int\limits_t^{t_5}\phi_{j_4}(t_4)
\int\limits_t^{t_4}\phi_{j_3}(t_3)C_{j_1 j_1}(t_3,t)
dt_3 dt_4 dt_5,
$$

\vspace{2mm}
$$
C_{j_5 j_4 j_1 j_2 j_1}=
\int\limits_t^T\phi_{j_5}(t_5)\int\limits_t^{t_5}\phi_{j_4}(t_4)
\int\limits_t^{t_4}\phi_{j_2}(t_2)
C_{j_1}(t_2,t)C_{j_1}(t_4,t_2)
dt_2 dt_4 dt_5,
$$

\vspace{2mm}
$$
C_{j_5 j_1 j_3 j_2 j_1}=
\int\limits_t^T\phi_{j_5}(t_5)\int\limits_t^{t_5}\phi_{j_3}(t_3)
\int\limits_t^{t_3}\phi_{j_2}(t_2)
C_{j_1}(t_2,t)C_{j_1}(t_5,t_3)
dt_2 dt_3 dt_5,
$$

\vspace{2mm}
$$
C_{j_1 j_4 j_3 j_2 j_1}=
\int\limits_t^T\phi_{j_4}(t_4)\int\limits_t^{t_4}\phi_{j_3}(t_3)
\int\limits_t^{t_3}\phi_{j_2}(t_2)
C_{j_1}(t_2,t)C_{j_1}(T,t_4)
dt_2 dt_3 dt_4,
$$

\vspace{2mm}
$$
C_{j_5 j_4 j_2 j_2 j_1}=
\int\limits_t^T\phi_{j_5}(t_5)\int\limits_t^{t_5}\phi_{j_4}(t_4)
\int\limits_t^{t_4}\phi_{j_1}(t_1)
C_{j_2 j_2}(t_4,t_1)
dt_1 dt_4 dt_5,
$$

\vspace{2mm}
$$
C_{j_5 j_2 j_3 j_2 j_1}=
\int\limits_t^T\phi_{j_5}(t_5)\int\limits_t^{t_5}\phi_{j_3}(t_3)
\int\limits_t^{t_3}\phi_{j_1}(t_1)
C_{j_2}(t_3,t_1)C_{j_2}(t_5,t_3)
dt_1 dt_3 dt_5,
$$

\vspace{2mm}
$$
C_{j_2 j_4 j_3 j_2 j_1}=
\int\limits_t^T\phi_{j_4}(t_4)\int\limits_t^{t_4}\phi_{j_3}(t_3)
\int\limits_t^{t_3}\phi_{j_1}(t_1)
C_{j_2}(t_3,t_1)C_{j_2}(T,t_4)
dt_1 dt_3 dt_4,
$$

\vspace{2mm}
$$
C_{j_5 j_3 j_3 j_2 j_1}
\int\limits_t^T\phi_{j_5}(t_5)\int\limits_t^{t_5}\phi_{j_2}(t_2)
\int\limits_{t}^{t_2}\phi_{j_1}(t_1)
C_{j_3 j_3}(t_5,t_2)
dt_1 dt_2 dt_5,
$$

\vspace{2mm}
$$
C_{j_3 j_4 j_3 j_2 j_1}=
\int\limits_t^T\phi_{j_4}(t_4)\int\limits_{t}^{t_4}\phi_{j_2}(t_2)
\int\limits_{t}^{t_2}\phi_{j_1}(t_1)
C_{j_3}(t_4,t_2)C_{j_3}(T,t_4)
dt_1 dt_2 dt_4,
$$

\vspace{2mm}
$$
C_{j_4 j_4 j_3 j_2 j_1}=
\int\limits_t^T\phi_{j_3}(t_3)\int\limits_t^{t_3}\phi_{j_2}(t_2)
\int\limits_{t}^{t_2}\phi_{j_1}(t_1) 
C_{j_4 j_4}(T,t_3)dt_1 dt_2dt_3,
$$

\vspace{2mm}
$$
C_{j_5 j_3 j_3 j_1 j_1}=
\int\limits_t^T \phi_{j_5}(t_5)
C_{j_3 j_3 j_1 j_1}(t_5,t)dt_5,
$$

\vspace{2mm}
$$
C_{j_5 j_2 j_1 j_2 j_1}=
\int\limits_t^T \phi_{j_5}(t_5) 
C_{j_2 j_1 j_2 j_1} (t_5,t)dt_5,
$$

\vspace{2mm}
$$
C_{j_5 j_1 j_2 j_2 j_1}=
\int\limits_t^T \phi_{j_5}(t_5)
C_{j_1 j_2 j_2 j_1} (t_5,t)dt_5,
$$

\vspace{2mm}
$$
C_{j_4 j_4 j_2 j_2 j_1}=
\int\limits_t^T  \phi_{j_1}(t_1)
C_{j_4 j_4 j_2 j_2}(T,t_1)dt_1,
$$

\vspace{2mm}
$$
C_{j_3 j_2 j_3 j_2 j_1}=
\int\limits_t^T  \phi_{j_1}(t_1)
C_{j_3 j_2 j_3 j_2}(T,t_1)dt_1,
$$

\vspace{2mm}
$$
C_{j_2 j_3 j_3 j_2 j_1}=
\int\limits_t^T  \phi_{j_1}(t_1)
C_{j_2 j_3 j_3 j_2}(T,t_1)dt_1,
$$

\vspace{2mm}
$$
C_{j_4 j_4 j_3 j_1 j_1}
=\int\limits_t^T  \phi_{j_3}(t_3)
C_{j_1 j_1}(t_3,t)C_{j_4 j_4}(T,t_3)
dt_3,
$$

\vspace{2mm}
$$
C_{j_2 j_4 j_1 j_2 j_1}=
\int\limits_t^T \phi_{j_4}(t_4)
C_{j_1 j_2 j_1}(t_4,t)
C_{j_2}(T,t_4) dt_4,
$$

\vspace{2mm}
$$
C_{j_2 j_1 j_3 j_2 j_1}=
\int\limits_t^T
\phi_{j_3}(t_3)C_{j_2 j_1}(t_3,t)
C_{j_2 j_1}(T,t_3)
dt_3,
$$

\vspace{2mm}
$$
C_{j_3 j_1 j_3 j_2 j_1}=\int\limits_t^T
\phi_{j_2}(t_2)C_{j_1}(t_2,t)
C_{j_3 j_1 j_3}(T,t_2)
dt_2,
$$

\vspace{2mm}
$$
C_{j_1 j_2 j_3 j_2 j_1}=
\int\limits_t^T
\phi_{j_3}(t_3)
C_{j_2 j_1}(t_3,t)C_{j_1 j_2}(T,t_3)dt_3,
$$

\vspace{2mm}
$$
C_{j_3 j_4 j_3 j_1 j_1}=
\int\limits_t^T  
\phi_{j_4}(t_4)C_{j_3 j_1 j_1}(t_4,t)C_{j_3}(T,t_4)dt_4,
$$

\vspace{2mm}
$$
C_{j_4 j_4 j_1 j_2 j_1}=
\int\limits_t^T  \phi_{j_2}(t_2)C_{j_1}(t_2,t)
C_{j_4 j_4 j_1}(T,t_2)dt_2,
$$

\vspace{2mm}
$$
C_{j_1 j_4 j_2 j_2 j_1}=
\int\limits_t^T \phi_{j_4}(t_4)
C_{j_2 j_2 j_1}(t_4,t)C_{j_1}(T,t_4)dt_4,
$$

\vspace{2mm}
$$
C_{j_1 j_3 j_3 j_2 j_1}=
\int\limits_t^T  
\phi_{j_2}(t_2)
C_{j_1}(t_2,t) C_{j_1 j_3 j_3}(T,t_2)
dt_2.
$$

\vspace{4mm}

It is easy to see (based on the above relations)
that (\ref{09091}) will be satisfied (under the conditions
of Theorem~48) if 
(\ref{2024december12})--(\ref{2024december9})
are fulfilled.
The equalities 
(\ref{2024december12})--(\ref{2024december9})
are proved in Sect.~31.
The assertion
of Theorem~48 now follows from Theorem~59. 
Theorem~48 is proved.

Recall that for the case $k=6$, together with
(\ref{2024december12})--(\ref{2024december9}), the conditions 
(\ref{2024december10}), (\ref{2024december11}) and the equality 
(\ref{july90000}) ($k=2r,$ $k=6,$ $r=3$)  
must be satisfied (see the proof of Theorem~60).

\vspace{5mm}

\section{Partial Proof of the Condition (\ref{09091})}

\vspace{5mm}

In this section, we will prove (\ref{09091})
for the case when the condition $(A)$
and the relation (\ref{copa5}) are satisfied (see Sect.~30).

Suppose that $\{\phi_j(x)\}_{j=0}^{\infty}$
is an arbitrary complete orthonormal system of functions
in $L_2([t, T])$ and $\psi_1(\tau),\ldots,\psi_k(\tau)\equiv 1.$

It is easy to see that (\ref{09091}) will be proved
for the above case if we prove that

\vspace{-1mm}
\begin{equation}
\label{febr2025}
\left|\sum_{j_r,j_{r-2},\ldots, j_2=0}^{p}
C_{j_r j_r j_{r-2} j_{r-2} \ldots j_2 j_2}(s,\tau)\right|\le K <\infty,
\end{equation}

\vspace{3mm}
\noindent
where $p\in\mathbb{N},$ $r=2, 4, 6,\ldots,$ constant $K$ does not depend on $p, s, \tau$
(but only on $t, T$),

\begin{equation}
\label{utoch1}
C_{j_k \ldots j_1}(s,\tau)=\int\limits_{\tau}^s
\phi_{j_k}(t_k)\ldots
\int\limits_{\tau}^{t_2}
\phi_{j_1}(t_1)dt_1\ldots dt_k,
\end{equation}

\vspace{3mm}
\noindent
where $k\in \mathbb{N},$ $t\le\tau<s\le T$.

By analogy with (\ref{july7028}) we obtain

$$
C_{j_r j_r j_{r-2}j_{r-2}\ldots j_2 j_2}(s,\tau) +
C_{j_2 j_2\ldots j_{r-2}j_{r-2} j_r j_r}(s,\tau)=
$$

\vspace{3mm}
$$
=
C_{j_r}(s,\tau) \cdot  C_{j_{r} j_{r-2} j_{r-2} \ldots j_4 j_4 j_2 j_2}(s,\tau)
-C_{j_{r} j_r}(s,\tau) \cdot C_{j_{r-2} j_{r-2}\ldots j_4 j_4 j_2 j_2}(s,\tau)+
$$

\vspace{3mm}
$$
+C_{j_{r-2} j_{r} j_r}(s,\tau) \cdot
C_{j_{r-2} j_{r-4}j_{r-4}\ldots j_4 j_4 j_2 j_2}(s,\tau)
- \ldots 
$$

\vspace{3mm}
\begin{equation}
\label{febr2025a}
- C_{j_4 j_4 \ldots j_{r-2}j_{r-2} j_{r}j_{r}}(s,\tau) \cdot C_{j_2 j_2}(s,\tau)+
C_{j_2 j_4 j_4 \ldots j_{r-2}j_{r-2} j_r j_r}(s,\tau) \cdot C_{j_2}(s,\tau).
\end{equation}

\vspace{6mm}

Applying (\ref{febr2025a}), we get

$$
2 \sum_{j_r,j_{r-2},\ldots, j_4, j_2=0}^{p} C_{j_r j_r j_{r-2}j_{r-2}\ldots j_4 j_4 j_2 j_2}(s,\tau)=
$$

\vspace{3mm}
$$
=
\sum_{j_r=0}^p  C_{j_r}(s,\tau)\sum_{j_{r-2},\ldots, j_4, j_2=0}^{p}
 C_{j_{r} j_{r-2} j_{r-2} \ldots j_4 j_4 j_2 j_2}(s,\tau)
-
$$

\vspace{3mm}
$$
-\sum_{j_r=0}^{p} C_{j_{r} j_r}(s,\tau) \sum_{j_{r-2},\ldots, j_4, j_2=0}^{p}
C_{j_{r-2} j_{r-2}\ldots j_4 j_4 j_2 j_2}(s,\tau)+
$$

\vspace{3mm}
$$
+  \sum_{j_{r-2}=0}^{p} \sum_{j_{r}=0}^{p} C_{j_{r-2} j_{r} j_r}(s,\tau) 
\sum_{j_{r-4},\ldots, j_4, j_2=0}^{p} C_{j_{r-2} j_{r-4}j_{r-4}\ldots j_4 j_4 j_2 j_2}(s,\tau)
- \ldots 
$$

\vspace{3mm}
$$
- \sum_{j_r,j_{r-2},\ldots, j_4=0}^{p}
C_{j_4 j_4 \ldots j_{r-2}j_{r-2} j_{r}j_{r}}(s,\tau) \sum_{j_{2}=0}^{p}C_{j_2 j_2}(s,\tau)+
$$

\vspace{3mm}
\begin{equation}
\label{febr2025b}
+
\sum_{j_{2}=0}^{p} \sum_{j_r,j_{r-2},\ldots, j_4=0}^{p} 
C_{j_2 j_4 j_4\ldots j_{r-2}j_{r-2} j_r j_r}(s,\tau)\cdot  C_{j_2}(s,\tau).
\end{equation}

\vspace{5mm}

Let us prove (\ref{febr2025}) by induction.
The equality (\ref{febr2025}) is proved for $r=2, 4$ (see 
(\ref{april48}), (\ref{april50}) and the relation $C_{j_1 j_1}(s,\tau)=
\frac{1}{2}\left(C_{j_1}(s,\tau)\right)^2$
for the case under consideration).
Suppose that

\begin{equation}
\label{febr2025b1}
\left|\sum_{j_6, j_4, j_2=0}^{p}
C_{j_6 j_6 j_{4} j_{4} j_2 j_2}(s,\tau)\right|\le K <\infty,
\end{equation}

\vspace{2mm}
\begin{equation}
\label{febr2025b2}
\left|\sum_{j_8, j_6, j_4, j_2=0}^{p}
C_{j_8 j_8 j_6 j_6 j_{4} j_{4} j_2 j_2}(s,\tau)\right|\le K <\infty,
\end{equation}

\vspace{-1mm}
$$
\ldots
$$

\vspace{-5mm}
\begin{equation}
\label{febr2025b3}
\left|\sum_{j_{r-2},j_{r-4},\ldots, j_2=0}^{p}
C_{j_{r-2} j_{r-2} j_{r-4} j_{r-4} \ldots j_2 j_2}(s,\tau)\right|\le K <\infty
\end{equation}

\vspace{4mm}
\noindent
and prove (\ref{febr2025}).

Using the induction hypothesis (see (\ref{febr2025b1})--(\ref{febr2025b3})), we obtain

\begin{equation}
\label{febr2025b4}
\left|\sum_{j_r=0}^{p} C_{j_{r} j_r}(s,\tau) \sum_{j_{r-2},\ldots, j_4, j_2=0}^{p}
C_{j_{r-2} j_{r-2}\ldots j_4 j_4 j_2 j_2}(s,\tau)\right|\le K^2<\infty,
\end{equation}

\vspace{2mm}
\begin{equation}
\label{febr2025b5}
\left|
\sum_{j_r, j_{r-2}=0}^{p} C_{j_{r-2} j_{r-2} j_{r} j_{r}}(s,\tau) \sum_{j_{r-4},\ldots, j_4, j_2=0}^{p}
C_{j_{r-4} j_{r-4}\ldots j_4 j_4 j_2 j_2}(s,\tau)\right|\le K^2<\infty,
\end{equation}

\vspace{2mm}
$$
\ldots
$$

\vspace{-2mm}
\begin{equation}
\label{febr2025b6}
\left|\sum_{j_r,j_{r-2},\ldots, j_4=0}^{p}
C_{j_4 j_4 \ldots j_{r-2}j_{r-2} j_{r}j_{r}}(s,\tau) \sum_{j_{2}=0}^{p}C_{j_2 j_2}(s,\tau)
\right|\le K^2<\infty.
\end{equation}

\vspace{5mm}

Applying the inequality 
of Cauchy--Bunyakovsky, Parseval's equality and the induction hypothesis, we obtain

$$
\left(
\sum_{j_r=0}^p  C_{j_r}(s,\tau)\sum_{j_{r-2},\ldots, j_4, j_2=0}^{p}
C_{j_{r} j_{r-2} j_{r-2} \ldots j_4 j_4 j_2 j_2}(s,\tau)\right)^2\le
$$

\vspace{3mm}
$$
\le \sum_{j_r=0}^p  \left(C_{j_r}(s,\tau)\right)^2
\sum_{j_r=0}^p \left(
\sum_{j_{r-2},\ldots, j_4, j_2=0}^{p}
C_{j_{r} j_{r-2} j_{r-2} \ldots j_4 j_4 j_2 j_2}(s,\tau)\right)^2\le
$$

\vspace{3mm}
$$
\le \sum_{j_r=0}^{\infty}  \left(C_{j_r}(s,\tau)\right)^2
\sum_{j_r=0}^{\infty} \left(
\sum_{j_{r-2},\ldots, j_4, j_2=0}^{p}
C_{j_{r} j_{r-2} j_{r-2} \ldots j_4 j_4 j_2 j_2}(s,\tau)\right)^2\le
$$

\vspace{3mm}
$$
\le K_1 
\sum_{j_r=0}^{\infty} \left(
\sum_{j_{r-2},\ldots, j_4, j_2=0}^{p}
C_{j_{r} j_{r-2} j_{r-2} \ldots j_4 j_4 j_2 j_2}(s,\tau)\right)^2=
$$

\vspace{3mm}
$$
=K_1 
\sum_{j_r=0}^{\infty} \left(\int\limits_{\tau}^s \phi_{j_r}(u)
\sum_{j_{r-2},\ldots, j_4, j_2=0}^{p}
C_{j_{r-2} j_{r-2} \ldots j_4 j_4 j_2 j_2}(u,\tau)du\right)^2=
$$

\vspace{3mm}
$$
=K_1 
\int\limits_{\tau}^s \left(
\sum_{j_{r-2},\ldots, j_4, j_2=0}^{p}
C_{j_{r-2} j_{r-2} \ldots j_4 j_4 j_2 j_2}(u,\tau)\right)^2 du \le
$$

\vspace{3mm}
\begin{equation}
\label{febr2025b7}
\le K_1 K^2 
\int\limits_{\tau}^s du \le (T-t)K_1 K^2=K_2<\infty,
\end{equation}

\vspace{3mm}
\noindent
where constant $K_2$ does not depend on $p, s, \tau;$

\vspace{2mm}
$$
\left(\sum_{j_{r-2}=0}^{p} \sum_{j_{r}=0}^{p} C_{j_{r-2} j_{r} j_r}(s,\tau) 
\sum_{j_{r-4},\ldots, j_4, j_2=0}^{p} C_{j_{r-2} j_{r-4}j_{r-4}\ldots 
j_4 j_4 j_2 j_2}(s,\tau)\right)^2\le
$$

\vspace{3mm}
$$
\le \sum_{j_{r-2}=0}^{p} \left(\sum_{j_{r}=0}^{p} C_{j_{r-2} j_{r} j_r}(s,\tau)\right)^2 
\sum_{j_{r-2}=0}^{p}
\left(\sum_{j_{r-4},\ldots, j_4, j_2=0}^{p} C_{j_{r-2} j_{r-4}j_{r-4}\ldots 
j_4 j_4 j_2 j_2}(s,\tau)\right)^2\le
$$

\vspace{3mm}
$$
\le \sum_{j_{r-2}=0}^{\infty} \left(\sum_{j_{r}=0}^{p} C_{j_{r-2} j_{r} j_r}(s,\tau)\right)^2 
\sum_{j_{r-2}=0}^{\infty}
\left(\sum_{j_{r-4},\ldots, j_4, j_2=0}^{p} C_{j_{r-2} j_{r-4}j_{r-4}\ldots 
j_4 j_4 j_2 j_2}(s,\tau)\right)^2=
$$

\vspace{3mm}
$$
=\sum_{j_{r-2}=0}^{\infty} \left(
\int\limits_{\tau}^s \phi_{j_{r-2}}(u)
\sum_{j_{r}=0}^{p} C_{j_{r} j_r}(u,\tau) du\right)^2 \times
$$

\vspace{3mm}
$$
\times
\sum_{j_{r-2}=0}^{\infty}
\left(\int\limits_{\tau}^s \phi_{j_{r-2}}(u)
\sum_{j_{r-4},\ldots, j_4, j_2=0}^{p} C_{j_{r-4}j_{r-4}\ldots j_4 j_4 j_2 j_2}(u,\tau) du\right)^2=
$$

\vspace{3mm}
$$
=
\int\limits_{\tau}^s \left(
\sum_{j_{r}=0}^{p} C_{j_{r} j_r}(u,\tau)\right)^2 du \times
$$

\vspace{3mm}
\begin{equation}
\label{febr2025b8}
\times
\int\limits_{\tau}^s \left(
\sum_{j_{r-4},\ldots, j_4, j_2=0}^{p} C_{j_{r-4}j_{r-4}\ldots j_4 j_4 j_2 j_2}(u,\tau)\right)^2 du\le
K^4 (T-t)^2 = K_3 <\infty.
\end{equation}

\vspace{6mm}

Similarly, we get

\vspace{-1mm}
\begin{equation}
\label{febr2025b9}
\left(\sum_{j_{r-4}=0}^{p} \sum_{j_{r}, j_{r-2}=0}^{p} C_{j_{r-4} j_{r-2} j_{r-2} j_r j_r}(s,\tau) 
\sum_{j_{r-6},\ldots, j_4, j_2=0}^{p} C_{j_{r-4} j_{r-6}j_{r-6}\ldots j_4 j_4 j_2 j_2}(s,\tau)\right)^2
\le K_4 <
\infty,
\end{equation}

\vspace{2mm}
$$
\ldots
$$

\vspace{-2mm}
\begin{equation}
\label{febr2025b10}
\left(
\sum_{j_{4}=0}^{p} \sum_{j_r,j_{r-2},\ldots, j_6=0}^{p} 
C_{j_4 j_6 j_6\ldots j_{r-2}j_{r-2} j_r j_r}(s,\tau)
\sum_{j_{2}=0}^{p} C_{j_4 j_2 j_2}(s,\tau)\right)^2
\le K_4 <\infty,
\end{equation}
\begin{equation}
\label{febr2025b11}
\left(
\sum_{j_{2}=0}^{p} \sum_{j_r,j_{r-2},\ldots, j_4=0}^{p} 
C_{j_2 j_4 j_4\ldots j_{r-2}j_{r-2} j_r j_r}(s,\tau)\cdot  C_{j_2}(s,\tau)\right)^2
\le K_4 <\infty,
\end{equation}

\vspace{4mm}
\noindent
where constant $K_4$ does not depend on $p, s, \tau.$

Combining (\ref{febr2025b}), (\ref{febr2025b4})--(\ref{febr2025b6}),
(\ref{febr2025b7}), (\ref{febr2025b8}),
(\ref{febr2025b9})--(\ref{febr2025b11}),
we obtain (\ref{febr2025}).
The equality (\ref{09091}) is proved
for the case when the condition $(A)$
and the relation (\ref{copa5}) are satisfied
($\psi_1(\tau),\ldots,\psi_k(\tau)\equiv 1$).

\vspace{5mm}

\section{Further Development of the Approach
Based on Theorem~59 for the Case $\psi_1(\tau),\ldots, \psi_7(\tau)
\equiv 1$. Expansion of Iterated Stratonovich Stochastic Integrals
of Multiplicity 7 (The Cases of Legendre 
Polynomials and Trigonometric Functions)}

\vspace{5mm}

Unfortunately, the approach from the previous section 
can be generalized only partially to the case
when the condition $(A)$
and the relation (\ref{copa6}) are satisfied (see Sect.~30).
In particular, the mentioned approach
is applicable to the proof of inequality

\vspace{-1mm}
$$
\left|\sum\limits_{j_1,j_2,j_3=0}^p C_{j_3 j_2 j_1 j_3 j_2 j_1}(s,\tau)\right|\le K<\infty,
$$

\vspace{3mm}
\noindent
but is not applicable to the proof of inequality

\vspace{-1mm}
$$
\left|\sum\limits_{j_1,j_2,j_3=0}^p C_{j_2 j_3 j_3 j_1 j_2 j_1}(s,\tau)\right|\le K<\infty,
$$

\vspace{3mm}
\noindent
where $C_{j_k \ldots j_1}(s,\tau)$ is defined by (\ref{utoch1}),
constant $K$ does not depend on $p, s, \tau$
$(p\in \mathbb{N},\ t\le \tau< s\le T).$

In this section, we will restrict ourselves to the case
$k=7,$ $r=1,2,3$ and we will also assume that 
$\{\phi_j(x)\}_{j=0}^{\infty}$ is a complete orthonormal
system of Legendre polynomials or trigonometric functions
in the space $L_2([t, T])$.

Note that the condition (\ref{09091})
can be weakened. Namely, the constant $K^2$
can be replaced by the function $F$ such that
$\psi_{q_1}^2\ldots \psi_{q_{k-2r}}^2 F \in L_1([t, T]^{k-2r})$.
For the trigonometric case, we will prove (\ref{09091}) for $k=7,$ $r=1,2,3$.
For the polynomial case, we will prove a weakened version of
(\ref{09091}) for $k=7,$ $r=1,2,3$ (the constant $K$ and the above function 
$F$ will be used in the weakened version of 
(\ref{09091})).

Obviously, that the conditions
(\ref{2024december12})--(\ref{2024december11})
together with the following condition

\vspace{-1mm}
\begin{equation}
\label{cc123}
\left\vert \sum\limits_{j_1, j_2=0}^p C_{j_1}(s,\tau)C_{j_2}(\rho,v)
C_{j_1}(\theta,u)
C_{j_2}(\mu,w)\right\vert\le K
\end{equation}

\vspace{3mm}
\noindent
cover the case $k=7,$ $r=1,2$ (see (\ref{09091})),
where $p\in\mathbb{N},$ $t\le \tau < s \le T,$ $t\le u<\theta \le T,$
$t\le v<\rho \le T,$ $t\le w<\mu \le T,$ constant $K$ does not depend on 
$p, s, \tau, u, \theta, v, \rho, w, \mu$ (but only on $t, T$).
The inequality (\ref{cc123}) is easily verified
using (\ref{dsds14fffff}).

Now let us focus on the proof of (\ref{09091})
for the case $k=7$ and $r=3$. So, we need to prove that

\vspace{-1mm}
\begin{equation}
\label{march0001}
\left|\sum\limits_{j_{g_1},j_{g_3},j_{g_{5}}=0}^p
C_{j_{d_1} j_{d_1-1}j_{d_1-2}j_{d_1-3}j_{d_1-4}j_{d_1-5}}(s,\tau)
\biggl|_{j_{g_1}=j_{g_2},j_{g_3}=j_{g_4},j_{g_5}=j_{g_6}}\right|\le K<\infty,
\end{equation}

\vspace{2mm}
\begin{equation}
\label{march0002}
\left|\sum\limits_{j_{g_1},j_{g_3},j_{g_{5}}=0}^p
\bigl(C_{j_{d_2} j_{d_2-1}j_{d_2-2}j_{d_2-3}j_{d_2-4}}(s,\tau)
C_{j_{d_1}}(\theta,u)
\bigr)\biggl|_{j_{g_1}=j_{g_2},
j_{g_3}=j_{g_4},j_{g_5}=j_{g_6}}\right|\le K<\infty,
\end{equation}

\vspace{2mm}
\begin{equation}
\label{march0003}
\left|\sum\limits_{j_{g_1},j_{g_3},j_{g_{5}}=0}^p
\left(C_{j_{d_2} j_{d_2-1}j_{d_2-2}j_{d_2-3}}(s,\tau)
C_{j_{d_1}j_{d_1-1}}(\theta,u)\right)\biggl|_{j_{g_1}=j_{g_2},
j_{g_3}=j_{g_4},j_{g_5}=j_{g_6}}\right|\le K<\infty,
\end{equation}

\vspace{2mm}
\begin{equation}
\label{march0004}
\left|\sum\limits_{j_{g_1},j_{g_3},j_{g_{5}}=0}^p
\left(C_{j_{d_2} j_{d_2-1}j_{d_2-2}}(s,\tau)
C_{j_{d_1} j_{d_1-1}j_{d_1-2}}(\theta,u)\right)\biggl|_{j_{g_1}=j_{g_2},
j_{g_3}=j_{g_4},j_{g_5}=j_{g_6}}\right|\le K<\infty,
\end{equation}

\vspace{3mm}
\noindent
where $p\in\mathbb{N},$ $t\le \tau < s \le T,$ $t\le u<\theta \le T,$
constant $K$ does not depend on 
$p, s, \tau, u, \theta$ (but only on $t, T$) and may differ from line to line;
another notations are the same as in Sect.~30.

The inequalities (\ref{march0002})--(\ref{march0004})
are proved using the same technique as 
inequalities (\ref{2024december12})--(\ref{2024december11}) (see Sect.~31).
Here we will only prove as an example the following
special case of the inequality (\ref{march0003})
\begin{equation}
\label{march0005}
\left|\sum\limits_{j_1, j_2, j_3=0}^p C_{j_2 j_3 j_2 j_1}(s,\tau)C_{j_3 j_1}(\theta,u)\right|
\le K<\infty.
\end{equation}

\vspace{3mm}

Using the 
Cau\-chy--Bunyakovsky inequality as well as 
Fubini's Theorem, Parseval's equality and (\ref{2024december13}), we have

\vspace{-2mm}
$$
\left(\sum\limits_{j_1, j_2, j_3=0}^p C_{j_2 j_3 j_2 j_1}(s,\tau)C_{j_3 j_1}(\theta,u)\right)^2\le
$$

\vspace{2mm}
$$
\le
\sum\limits_{j_1,j_3=0}^p \left(\sum\limits_{j_2=0}^p C_{j_2 j_3 j_2 j_1}(s,\tau)\right)^2
\sum\limits_{j_1, j_3=0}^p C_{j_3 j_1}^2(\theta,u)\le
$$

\vspace{2mm}
$$
\le \sum\limits_{j_1,j_3=0}^{\infty} \left(\sum\limits_{j_2=0}^p 
\int\limits_{\tau}^s \phi_{j_2}(u)
\int\limits_{\tau}^u \phi_{j_3}(z)
\int\limits_{\tau}^{z} \phi_{j_2}(y)
\int\limits_{\tau}^{y} \phi_{j_1}(x)dx dy dz du\right)^2 \times
$$

\vspace{2mm}
$$
\times
\sum\limits_{j_1,j_3=0}^{\infty} C_{j_3 j_1}^2(\theta,u)=
$$

\vspace{2mm}
$$
=\sum\limits_{j_1,j_3=0}^{\infty} \left(\sum\limits_{j_2=0}^p 
\int\limits_{\tau}^s 
\phi_{j_3}(z)
\int\limits_{\tau}^{z} \phi_{j_2}(y)
\int\limits_{\tau}^{y} \phi_{j_1}(x)dx dy
\int\limits_z^s 
\phi_{j_2}(u)
du dz\right)^2 
\cdot \frac{(\theta-u)^2}{2}=
$$

\vspace{2mm}
$$
=\frac{(\theta-u)^2}{2}\sum\limits_{j_1,j_3=0}^{\infty} \left(\sum\limits_{j_2=0}^p 
\int\limits_{\tau}^s 
\phi_{j_3}(z)
\int\limits_{\tau}^{z} \phi_{j_1}(x) \int\limits_{x}^{z} \phi_{j_2}(y)
dy dx
\int\limits_z^s 
\phi_{j_2}(u)
du dz\right)^2 =
$$

\vspace{2mm}
$$
=\frac{(\theta-u)^2}{2}\sum\limits_{j_1,j_3=0}^{\infty} \left(
\int\limits_{\tau}^s 
\phi_{j_3}(z)
\int\limits_{\tau}^{z} \phi_{j_1}(x) \sum\limits_{j_2=0}^p  C_{j_2}(z,x) 
C_{j_2}(s,z) dx dz\right)^2 
=
$$

\vspace{2mm}
$$
=
\frac{(\theta-u)^2}{2}\int\limits_{\tau}^s 
\int\limits_{\tau}^{z} \left(\sum\limits_{j_2=0}^p  C_{j_2}(z,x) 
C_{j_2}(s,z)\right)^2 dx dz \le
$$

\vspace{2mm}
\begin{equation}
\label{marchto1}
\le K^2 \frac{(\theta-u)^2}{2}\frac{(s-\tau)^2}{2}\le K^2\frac{(T-t)^4}{4}=K_1.
\end{equation}

\vspace{5mm}
\noindent
The equality (\ref{march0005}) is proved.

The main difficulty is related to the proof
of the inequality (\ref{march0001}). Further, we prove (\ref{march0001})
for all 15 possible cases under the assumption
that 
$\{\phi_j(x)\}_{j=0}^{\infty}$ is a complete orthonormal
system of Legendre polynomials or trigonometric functions
in the space $L_2([t, T])$. As we noted above,
in some situations we will need a function $F\in L_1([t, T])$
instead of a constant $K^2$ for the polynomial case.

It is easy to see that (\ref{march0001}) reduces
to the following 15 inequalities

\vspace{-1mm}
\begin{equation}
\label{marsixsix8}
\left|\sum_{j_1, j_2, j_3=0}^{p}
C_{j_3 j_2 j_1 j_3 j_2 j_1}(s,\tau)\right|\le K<\infty,
\end{equation}
\begin{equation}
\label{marsixsix9}
\left|\sum_{j_1, j_2, j_3=0}^{p}
C_{j_1 j_3 j_2 j_3 j_2 j_1}(s,\tau)\right|\le K<\infty,
\end{equation}
\begin{equation}
\label{marsixsix10}
\left|\sum_{j_1, j_2, j_3=0}^{p}
C_{j_3 j_2 j_3 j_1 j_2 j_1}(s,\tau)\right|\le K<\infty,
\end{equation}
\begin{equation}
\label{marsixsix4}
\left|\sum_{j_1, j_2, j_3=0}^{p}
C_{j_1 j_2 j_3 j_3 j_2 j_1}(s,\tau)\right|\le K<\infty,
\end{equation}
\begin{equation}
\label{marsixsix14}
\left|\sum_{j_1, j_2, j_3=0}^{p}
C_{j_1 j_2 j_2 j_3 j_3 j_1}(s,\tau)\right|\le K<\infty,
\end{equation}
\begin{equation}
\label{marsixsix3}
\left|\sum_{j_1, j_2, j_3=0}^{p}
C_{j_3 j_3 j_2 j_2 j_1 j_1}(s,\tau)\right|\le K<\infty,
\end{equation}
\begin{equation}
\label{marsixsix7}
\left|\sum_{j_1, j_2, j_3=0}^{p}
C_{j_2 j_3 j_3 j_2 j_1 j_1}(s,\tau)\right|\le K<\infty,
\end{equation}
\begin{equation}
\label{marsixsix6}
\left|\sum_{j_1, j_2, j_3=0}^{p}
C_{j_3 j_2 j_3 j_2 j_1 j_1}(s,\tau)\right|\le K<\infty,
\end{equation}
\begin{equation}
\label{marsixsix1}
\left|\sum_{j_1, j_2, j_3=0}^{p}
C_{j_3 j_3 j_2 j_1 j_2 j_1}(s,\tau)\right|\le K<\infty,
\end{equation}
\begin{equation}
\label{marsixsix2}
\left|\sum_{j_1, j_2, j_3=0}^{p}
C_{j_3 j_3 j_1 j_2 j_2 j_1}(s,\tau)\right|\le K<\infty,
\end{equation}
\begin{equation}
\label{marsixsix5}
\left|\sum_{j_1, j_2, j_3=0}^{p}
C_{j_2 j_1 j_3 j_3 j_2 j_1}(s,\tau)\right|\le K<\infty,
\end{equation}
\begin{equation}
\label{marsixsix12}
\left|\sum_{j_1, j_2, j_3=0}^{p}
C_{j_3 j_1 j_2 j_3 j_2 j_1}(s,\tau)\right|\le K<\infty,
\end{equation}
\begin{equation}
\label{marsixsix11}
\left|\sum_{j_1, j_2, j_3=0}^{p}
C_{j_2 j_3 j_1 j_3 j_2 j_1}(s,\tau)\right|\le K<\infty,
\end{equation}
\begin{equation}
\label{marsixsix13}
\left|\sum_{j_1, j_2, j_3=0}^{p}
C_{j_3 j_1 j_3 j_2 j_2 j_1}(s,\tau)\right|\le K<\infty,
\end{equation}
\begin{equation}
\label{marsixsix15}
\left|\sum_{j_1, j_2, j_3=0}^{p}
C_{j_2 j_3 j_3 j_1 j_2 j_1}(s,\tau)\right|\le K<\infty,
\end{equation}

\vspace{3mm}
\noindent
where $p\in\mathbb{N},$ $t\le \tau < s \le T,$ 
constant $K$ does not depend on 
$p, s, \tau$ (but only on $t, T$) and may differ from line to line.

More precisely, the conditions (\ref{marsixsix8})--(\ref{marsixsix15})
need to be proved in two cases:
1.~$\tau=t,$\ \  2.~$s=T.$ Further, we will 
not carry out such a refinement if
some estimate from 
(\ref{marsixsix8})--(\ref{marsixsix15}) is true for all
$\tau, s\in [t, T]$ ($\tau<s$).
Looking ahead, we note that consideration
of Cases 1 and 2 will be required
only for some inequalities from (\ref{marsixsix8})--(\ref{marsixsix15})
for the polynomial case.

The relation (\ref{marsixsix3}) is a particular case of (\ref{febr2025}).
Let us prove the inequalities (\ref{marsixsix8})--(\ref{marsixsix14}),
(\ref{marsixsix7})--(\ref{marsixsix15}). 

\vspace{2mm}

{\bf Step~1.}\ First, we prove (\ref{marsixsix8})--(\ref{marsixsix14}),
(\ref{marsixsix5}) using special
symmetry properties of the
Fourier coefficients. 

By analogy with (\ref{sixsix40})
we obtain

\vspace{-1mm}
$$
C_{j_6 j_5 j_4 j_3 j_2 j_1}(s,\tau)+C_{j_1 j_2 j_3 j_4 j_5 j_6}(s,\tau)=
$$

$$
=
C_{j_6}(s,\tau)C_{j_5 j_4 j_3 j_2 j_1}(s,\tau)-C_{j_5 j_6}(s,\tau)C_{j_4 j_3 j_2 j_1}(s,\tau)+
$$

$$
+C_{j_4 j_5 j_6}(s,\tau)C_{j_3 j_2 j_1}(s,\tau)-C_{j_3 j_4 j_5 j_6}(s,\tau)C_{j_2 j_1}(s,\tau)+
$$

\begin{equation}
\label{sixsix40eee}
+
C_{j_2 j_3 j_4 j_5 j_6}(s,\tau)C_{j_1}(s,\tau).
\end{equation}

\vspace{3mm}

Using (\ref{sixsix40eee}), we get

$$
\sum_{j_1,j_2,j_3=0}^{p}
C_{j_3 j_2 j_1 j_3 j_2 j_1}(s,\tau)=
\frac{1}{2}\sum_{j_1,j_2,j_3=0}^{p}\biggl(
C_{j_3}(s,\tau)C_{j_2 j_1 j_3 j_2 j_1}(s,\tau)-\biggr.
$$

\vspace{3mm}
$$
-C_{j_2 j_3}(s,\tau)C_{j_1 j_3 j_2 j_1}(s,\tau)+
C_{j_1 j_2 j_3}(s,\tau)C_{j_3 j_2 j_1}(s,\tau)-
$$

\begin{equation}
\label{march0009}
\biggl.-C_{j_3 j_1 j_2 j_3}(s,\tau)C_{j_2 j_1}(s,\tau)+
C_{j_2 j_3 j_1 j_2 j_3}(s,\tau)C_{j_1}(s,\tau)\biggr),
\end{equation}

\vspace{5mm}
$$
\sum_{j_1,j_2,j_3=0}^{p}
C_{j_1 j_3 j_2 j_3 j_2 j_1}(s,\tau)=
\frac{1}{2}
\sum_{j_1,j_2,j_3=0}^{p}\biggl(
C_{j_1}(s,\tau)C_{j_3 j_2 j_3 j_2 j_1}(s,\tau)-\biggr.
$$

\vspace{3mm}
$$
-C_{j_3 j_1}(s,\tau)C_{j_2 j_3 j_2 j_1}(s,\tau)+
C_{j_2 j_3 j_1}(s,\tau)C_{j_3 j_2 j_1}(s,\tau)-
$$

\begin{equation}
\label{march00010}
\biggl.-C_{j_3 j_2 j_3 j_1}(s,\tau)C_{j_2 j_1}(s,\tau)+
C_{j_2 j_3 j_2 j_3 j_1}(s,\tau)C_{j_1}(s,\tau)\biggr),
\end{equation}

\vspace{5mm}

$$
\sum_{j_1,j_2,j_3=0}^{p}
C_{j_3 j_2 j_3 j_1 j_2 j_1}(s,\tau)=
\frac{1}{2}\sum_{j_1,j_2,j_3=0}^{p}\biggl(
C_{j_3}(s,\tau)C_{j_2 j_3 j_1 j_2 j_1}(s,\tau)-\biggr.
$$

\vspace{3mm}
$$
-C_{j_2 j_3}(s,\tau)C_{j_3 j_1 j_2 j_1}(s,\tau)+
C_{j_3 j_2 j_3}(s,\tau)C_{j_1 j_2 j_1}(s,\tau)-
$$

\begin{equation}
\label{march00011}
\biggl.-C_{j_1 j_3 j_2 j_3}(s,\tau)C_{j_2 j_1}(s,\tau)+
C_{j_2 j_1 j_3 j_2 j_3}(s,\tau)C_{j_1}(s,\tau)\biggr),
\end{equation}

\vspace{5mm}

$$
\sum_{j_1,j_2,j_3=0}^{p}
C_{j_1 j_2 j_3 j_3 j_2 j_1}(s,\tau)=
\frac{1}{2}
\sum_{j_1,j_2,j_3=0}^{p}\biggl(
C_{j_1}(s,\tau)C_{j_2 j_3 j_3 j_2 j_1}(s,\tau)-\biggr.
$$

\vspace{3mm}
$$
-C_{j_2 j_1}(s,\tau)C_{j_3 j_3 j_2 j_1}(s,\tau)+
\left(C_{j_3 j_2 j_1}(s,\tau)\right)^2-
$$

\begin{equation}
\label{march00012}
\biggl.-C_{j_3 j_3 j_2 j_1}(s,\tau)C_{j_2 j_1}(s,\tau)+
C_{j_2 j_3 j_3 j_2 j_1}(s,\tau)C_{j_1}(s,\tau)\biggr),
\end{equation}

\vspace{5mm}

$$
\sum_{j_1,j_2,j_3=0}^{p}
C_{j_1 j_3 j_3 j_2 j_2 j_1}(s,\tau)=
\frac{1}{2}
\sum_{j_1,j_2,j_3=0}^{p}
\biggl(
C_{j_1}(s,\tau)C_{j_3 j_3 j_2 j_2 j_1}(s,\tau)-\biggr.
$$

\vspace{3mm}
$$
-C_{j_3 j_1}(s,\tau)C_{j_3 j_2 j_2 j_1}(s,\tau)+
C_{j_3 j_3 j_1}(s,\tau)C_{j_2 j_2 j_1}(s,\tau)-
$$

\begin{equation}
\label{march00013}
\biggl.-C_{j_2 j_3 j_3 j_1}(s,\tau)C_{j_2 j_1}(s,\tau)+
C_{j_2 j_2 j_3 j_3 j_1}(s,\tau)C_{j_1}(s,\tau)\biggr),
\end{equation}

\vspace{5mm}
$$
\sum_{j_1,j_2,j_3=0}^{p}
C_{j_2 j_1 j_3 j_3 j_2 j_1}(s,\tau)=
\frac{1}{2}
\sum_{j_1,j_2,j_3=0}^{p}
\biggl(
C_{j_2}(s,\tau)C_{j_1 j_3 j_3 j_2 j_1}(s,\tau)-\biggr.
$$

\vspace{3mm}
$$
-C_{j_1 j_2}(s,\tau)C_{j_3 j_3 j_2 j_1}(s,\tau)+
C_{j_3 j_1 j_2}(s,\tau)C_{j_3 j_2 j_1}(s,\tau)-
$$

\begin{equation}
\label{march00014}
\biggl.-C_{j_3 j_3 j_1 j_2}(s,\tau)C_{j_2 j_1}(s,\tau)+
C_{j_2 j_3 j_3 j_1 j_2}(s,\tau)C_{j_1}(s,\tau)\biggr).
\end{equation}

\vspace{5mm}

Applying to the right-hand sides of (\ref{march0009})--(\ref{march00014})
the technique that led to the estimate (\ref{marchto1}), we obtain the inequalities
(\ref{marsixsix8})--(\ref{marsixsix14}), (\ref{marsixsix5}).

\vspace{2mm}

{\bf Step~2.}\ It is not difficult to see that

\vspace{-1mm}
\begin{equation}
\label{march00015}
\sum_{j_1, j_2, j_3=0}^{p}
C_{j_3 j_3 j_1 j_2 j_2 j_1}(s,\tau)=
\sum_{j_1, j_2, j_3=0}^{p}
C_{j_1 j_1 j_2 j_3 j_3 j_2}(s,\tau),
\end{equation}

\vspace{1mm}
\begin{equation}
\label{march00016}
\sum_{j_1, j_2, j_3=0}^{p}
C_{j_3 j_3 j_2 j_1 j_2 j_1}(s,\tau)=
\sum_{j_1, j_2, j_3=0}^{p}
C_{j_1 j_1 j_2 j_3 j_2 j_3}(s,\tau),
\end{equation}

\vspace{1mm}
\begin{equation}
\label{march00017}
\sum_{j_1, j_2, j_3=0}^{p}
C_{j_2 j_3 j_3 j_1 j_2 j_1}(s,\tau)=
\sum_{j_1, j_2, j_3=0}^{p}
C_{j_1 j_2 j_2 j_3 j_1 j_3}(s,\tau).
\end{equation}

\vspace{5mm}

Further, 
using (\ref{march00015})--(\ref{march00017}) and (\ref{sixsix40eee}), we get

$$
\sum_{j_1, j_2, j_3=0}^{p}
C_{j_2 j_3 j_3 j_2 j_1 j_1}(s,\tau)
+
\sum_{j_1, j_2, j_3=0}^{p}
C_{j_3 j_3 j_1 j_2 j_2 j_1}(s,\tau)
=
$$

$$
=\sum_{j_1, j_2, j_3=0}^{p}
C_{j_2 j_3 j_3 j_2 j_1 j_1}(s,\tau)
+
\sum_{j_1, j_2, j_3=0}^{p}
C_{j_1 j_1 j_2 j_3 j_3 j_2}(s,\tau)
=
$$

\vspace{3mm}
$$
=\sum_{j_1,j_2,j_3=0}^{p}
\biggl(
C_{j_2}(s,\tau)C_{j_3 j_3 j_2 j_1 j_1}(s,\tau)-\biggr.
$$

\vspace{3mm}
$$
-C_{j_3 j_2}(s,\tau)C_{j_3 j_2 j_1 j_1}(s,\tau)+
C_{j_3 j_3 j_2}(s,\tau)C_{j_2 j_1 j_1}(s,\tau)-
$$

\begin{equation}
\label{march00018}
\biggl.-C_{j_2 j_3 j_3 j_2}(s,\tau)C_{j_1 j_1}(s,\tau)+
C_{j_1 j_2 j_3 j_3 j_2}(s,\tau)C_{j_1}(s,\tau)\biggr),
\end{equation}

\vspace{6mm}

$$
\sum_{j_1, j_2, j_3=0}^{p}
C_{j_3 j_2 j_3 j_2 j_1 j_1}(s,\tau)+
\sum_{j_1, j_2, j_3=0}^{p}
C_{j_3 j_3 j_2 j_1 j_2 j_1}(s,\tau)=
$$

$$
=
\sum_{j_1, j_2, j_3=0}^{p}
C_{j_3 j_2 j_3 j_2 j_1 j_1}(s,\tau)+
\sum_{j_1, j_2, j_3=0}^{p}
C_{j_1 j_1 j_2 j_3 j_2 j_3}(s,\tau)=
$$

\vspace{3mm}
$$
=\sum_{j_1,j_2,j_3=0}^{p}
\biggl(
C_{j_3}(s,\tau)C_{j_2 j_3 j_2 j_1 j_1}(s,\tau)-\biggr.
$$

\vspace{3mm}
$$
-C_{j_2 j_3}(s,\tau)C_{j_3 j_2 j_1 j_1}(s,\tau)+
C_{j_3 j_2 j_3}(s,\tau)C_{j_2 j_1 j_1}(s,\tau)-
$$

\begin{equation}
\label{march00019}
\biggl.-C_{j_2 j_3 j_2 j_3}(s,\tau)C_{j_1 j_1}(s,\tau)+
C_{j_1 j_2 j_3 j_2 j_3}(s,\tau)C_{j_1}(s,\tau)\biggr),
\end{equation}

\vspace{6mm}

$$
\sum_{j_1, j_2, j_3=0}^{p}
C_{j_3 j_1 j_3 j_2 j_2 j_1}(s,\tau)+
\sum_{j_1, j_2, j_3=0}^{p}
C_{j_2 j_3 j_3 j_1 j_2 j_1}(s,\tau)=
$$

$$
=\sum_{j_1, j_2, j_3=0}^{p}
C_{j_3 j_1 j_3 j_2 j_2 j_1}(s,\tau)+
\sum_{j_1, j_2, j_3=0}^{p}
C_{j_1 j_2 j_2 j_3 j_1 j_3}(s,\tau)=
$$

\vspace{3mm}
$$
=\sum_{j_1,j_2,j_3=0}^{p}
\biggl(
C_{j_3}(s,\tau)C_{j_1 j_3 j_2 j_2 j_1}(s,\tau)-\biggr.
$$

\vspace{3mm}
$$
-C_{j_1 j_3}(s,\tau)C_{j_3 j_2 j_2 j_1}(s,\tau)+
C_{j_3 j_1 j_3}(s,\tau)C_{j_2 j_2 j_1}(s,\tau)-
$$

\begin{equation}
\label{march00020}
\biggl.-C_{j_2 j_3 j_1 j_3}(s,\tau)C_{j_2 j_1}(s,\tau)+
C_{j_2 j_2 j_3 j_1 j_3}(s,\tau)C_{j_1}(s,\tau)\biggr).
\end{equation}

\vspace{5mm}

Applying to the right-hand sides of (\ref{march00018})--(\ref{march00020})
the technique that led to the estimate (\ref{marchto1}), we obtain the inequalities

\vspace{-1mm}
\begin{equation}
\label{march00021}
\left|\sum_{j_1, j_2, j_3=0}^{p}
C_{j_2 j_3 j_3 j_2 j_1 j_1}(s,\tau)
+
\sum_{j_1, j_2, j_3=0}^{p}
C_{j_3 j_3 j_1 j_2 j_2 j_1}(s,\tau)\right|\le K <\infty,
\end{equation}

\begin{equation}
\label{march00022}
\left|\sum_{j_1, j_2, j_3=0}^{p}
C_{j_3 j_2 j_3 j_2 j_1 j_1}(s,\tau)+
\sum_{j_1, j_2, j_3=0}^{p}
C_{j_3 j_3 j_2 j_1 j_2 j_1}(s,\tau)\right|\le K <\infty,
\end{equation}

\begin{equation}
\label{march00023}
\left|\sum_{j_1, j_2, j_3=0}^{p}
C_{j_3 j_1 j_3 j_2 j_2 j_1}(s,\tau)+
\sum_{j_1, j_2, j_3=0}^{p}
C_{j_2 j_3 j_3 j_1 j_2 j_1}(s,\tau)
\right|\le K <\infty,
\end{equation}

\vspace{4mm}
\noindent
where $p\in\mathbb{N},$ $t\le \tau < s \le T,$ 
constant $K$ does not depend on 
$p, s, \tau$ (but only on $t, T$) and may differ from line to line.

Note that $\left|a\right|\le K_1+K$ follows from
$\left|b\right|\le K$ and $\left|a+b\right|\le K_1,$
where $a, b, K, K_1\in \mathbb{R}.$
Indeed, we have $\left|a\right|=\left|a+b-b\right|\le
\left|a+b\right|+\left|b\right|\le K_1+K$. Then
from (\ref{march00021})--(\ref{march00023})
it follows that if we prove (\ref{marsixsix1}), (\ref{marsixsix2}), (\ref{marsixsix15}), then 
(\ref{marsixsix6}), (\ref{marsixsix7}), (\ref{marsixsix13}) will be proved.
Thus, it remains to prove 
(\ref{marsixsix1}), (\ref{marsixsix2}), (\ref{marsixsix12}), (\ref{marsixsix11}),
(\ref{marsixsix15}).

\vspace{2mm}

{\bf Step~3.}\ Let us prove 
(\ref{marsixsix1}), (\ref{marsixsix2}), (\ref{marsixsix12}), (\ref{marsixsix11}),
(\ref{marsixsix15}).
Consider (\ref{marsixsix11}). 
Using the 
Cau\-chy--Bunyakov\-sky inequality as well as 
Fubini's Theorem, Parseval's equality, (\ref{after80xx}), (\ref{dsds14fffff})
and Lebesgue's Dominated Convergence Theorem, we have

$$
\left(\sum_{j_1, j_2, j_3=0}^{p}
C_{j_2 j_3 j_1 j_3 j_2 j_1}(s,\tau)\right)^2=
\left(\sum_{j_2=0}^{p} 1\cdot \sum_{j_1, j_3=0}^{p}
C_{j_2 j_3 j_1 j_3 j_2 j_1}(s,\tau)\right)^2\le
$$

\vspace{2mm}
$$
\le \sum_{j_2=0}^{p} 1^2 \cdot \sum_{j_2=0}^{p}\left(\sum_{j_1, j_3=0}^{p}
C_{j_2 j_3 j_1 j_3 j_2 j_1}(s,\tau)\right)^2=
$$

\vspace{2mm}
$$
=
(p+1)\sum_{j_2=0}^{p}\left(\sum_{j_1, j_3=0}^{p}
C_{j_2 j_3 j_1 j_3 j_2 j_1}(s,\tau)\right)^2=
$$

\vspace{2mm}
$$
=
(p+1)\sum_{j_2=0}^{p}\left(\sum_{j_1, j_3=0}^{p}
\int\limits_{\tau}^s \phi_{j_2}(t_6)
\int\limits_{\tau}^{t_6} \phi_{j_2}(t_2)
C_{j_1}(t_2,\tau)C_{j_3 j_1 j_3}(t_6,t_2) dt_2 dt_6
\right)^2\le
$$

\vspace{2mm}
$$
\le
(p+1)\sum_{j_2, j_2'=0}^{p}\left(\sum_{j_1, j_3=0}^{p}
\int\limits_{\tau}^s \phi_{j_2}(t_6)
\int\limits_{\tau}^{t_6} \phi_{j_2'}(t_2)
C_{j_1}(t_2,\tau)C_{j_3 j_1 j_3}(t_6,t_2) dt_2 dt_6
\right)^2\le
$$

\vspace{2mm}
$$
\le
(p+1)\sum_{j_2, j_2'=0}^{\infty}\left(
\int\limits_{\tau}^s \phi_{j_2}(t_6)
\int\limits_{\tau}^{t_6} \phi_{j_2'}(t_2)
\sum_{j_1, j_3=0}^{p}C_{j_1}(t_2,\tau)C_{j_3 j_1 j_3}(t_6,t_2) dt_2 dt_6
\right)^2=
$$

\vspace{2mm}
$$
=
(p+1)
\int\limits_{\tau}^s 
\int\limits_{\tau}^{t_6} 
\left(\sum_{j_1=0}^{p}C_{j_1}(t_2,\tau)\sum_{j_3=0}^{p}C_{j_3 j_1 j_3}(t_6,t_2)\right)^2 dt_2 dt_6
=
$$

\vspace{2mm}
$$
=
(p+1)
\int\limits_{\tau}^s 
\int\limits_{\tau}^{t_6} 
\left(\sum_{j_1=0}^{p}C_{j_1}(t_2,\tau)\sum_{j_3=p+1}^{\infty}C_{j_3 j_1 j_3}(t_6,t_2)\right)^2 dt_2 dt_6
\le
$$

\vspace{2mm}
$$
\le
(p+1)
\int\limits_{\tau}^s 
\int\limits_{\tau}^{t_6} 
\sum_{j_1=0}^{p}C_{j_1}^2(t_2,\tau)
\sum_{j_1=0}^{p}\left(\sum_{j_3=p+1}^{\infty}C_{j_3 j_1 j_3}(t_6,t_2)\right)^2 dt_2 dt_6\le
$$

\vspace{2mm}
$$
\le
(p+1)
\int\limits_{\tau}^s 
\int\limits_{\tau}^{t_6} 
\sum_{j_1=0}^{\infty}C_{j_1}^2(t_2,\tau)
\sum_{j_1=0}^{p}\left(\sum_{j_3=p+1}^{\infty}C_{j_3 j_1 j_3}(t_6,t_2)\right)^2 dt_2 dt_6=
$$

\vspace{2mm}
$$
\le
(p+1)
\int\limits_{\tau}^s 
\int\limits_{\tau}^{t_6} 
(t_2-\tau)
\sum_{j_1=0}^{p}\left(\sum_{j_3=p+1}^{\infty}C_{j_3 j_1 j_3}(t_6,t_2)\right)^2 dt_2 dt_6=
$$

\vspace{2mm}
$$
=
(p+1)
\int\limits_{\tau}^s 
\int\limits_{\tau}^{t_6} 
(t_2-\tau)
\sum_{j_1=0}^{p}\left(\sum_{j_3=p+1}^{\infty}
\int\limits_{t_2}^{t_6}\phi_{j_1}(\theta)
C_{j_3}(\theta,t_2)C_{j_3}(t_6,\theta)d\theta\right)^2 dt_2 dt_6=
$$

\vspace{2mm}
$$
=
(p+1)
\int\limits_{\tau}^s 
\int\limits_{\tau}^{t_6} 
(t_2-\tau)
\sum_{j_1=0}^{p}\left(
\int\limits_{t_2}^{t_6}\phi_{j_1}(\theta)\sum_{j_3=p+1}^{\infty}
C_{j_3}(\theta,t_2)C_{j_3}(t_6,\theta)d\theta\right)^2 dt_2 dt_6\le
$$

\vspace{2mm}
$$
\le
(p+1)
\int\limits_{\tau}^s 
\int\limits_{\tau}^{t_6} 
(t_2-\tau)
\sum_{j_1=0}^{\infty}\left(
\int\limits_{t_2}^{t_6}\phi_{j_1}(\theta)\sum_{j_3=p+1}^{\infty}
C_{j_3}(\theta,t_2)C_{j_3}(t_6,\theta)d\theta\right)^2 dt_2 dt_6=
$$

\vspace{2mm}
\begin{equation}
\label{march00025}
=
(p+1)
\int\limits_{\tau}^s 
\int\limits_{\tau}^{t_6} 
(t_2-\tau)
\int\limits_{t_2}^{t_6}\left(\sum_{j_3=p+1}^{\infty}
C_{j_3}(\theta,t_2)C_{j_3}(t_6,\theta)\right)^2d\theta dt_2 dt_6.
\end{equation}

\vspace{5mm}
                                    
For the trigonometric case (Fourier basis), we have the following obvious estimate

\vspace{-2mm}
\begin{equation}
\label{march00026}
\left|C_j(x,v)\right|=\left|
\int\limits_v^x\phi_{j}(\tau)d\tau
\right|< \frac{C}{j}\ \ \ (j>0),
\end{equation}

\vspace{3mm}
\noindent
where constant $C$ does not depend on $j, x, v.$

Recall that (see (\ref{obana}))
\begin{equation}
\label{march00027}
\sum\limits_{j=p+1}^{\infty}\frac{1}{j^2}
\le \int\limits_{p}^{\infty}\frac{dx}{x^2}=\frac{1}{p}.
\end{equation}

\vspace{3mm}

Combining (\ref{march00025})--(\ref{march00027}), we get

\vspace{-1mm}
$$
\left(\sum_{j_1, j_2, j_3=0}^{p}
C_{j_2 j_3 j_1 j_3 j_2 j_1}(s,\tau)\right)^2\le \frac{K_1(p+1)}{p^2}\le K^2,
$$

\vspace{3mm}
\noindent
where constants $K, K_1$ depend only on  $t, T.$
The inequality (\ref{marsixsix11}) is proved for the trigonometric case.

For the polynomial case, by analogy with (\ref{after1940})
we have 

\vspace{-1mm}
\begin{equation}
\label{march00028}
\left|C_j(x,v)\right|=\left|
\int\limits_v^x\phi_{j}(\tau)d\tau
\right| <
\frac{C}{j^{1-\varepsilon/2}}\Biggl(\frac{1}{(1-z^2(x))^{1/4-\varepsilon/4}}+
\frac{1}{(1-z^2(v))^{1/4-\varepsilon/4}}\Biggr),
\end{equation}

\vspace{3mm}
\noindent
where 
$j\in \mathbb{N},$ $z(x), z(v)\in (-1, 1)$ 
($z(x)$ is defined by (\ref{zz1})),
$x, v\in (t, T),$ $\varepsilon\in (0,1)$ is an arbitrary
small positive real number,
constant $C$ does not depend on $j$.

Recall that
(see (\ref{after1944}))
\begin{equation}
\label{march00029}
\sum\limits_{j=p+1}^{\infty}\frac{1}{j^{2-\varepsilon}}
\le \int\limits_{p}^{\infty}\frac{dx}{x^{2-\varepsilon}}=
\frac{1}{(1-\varepsilon)p^{1-\varepsilon}}.
\end{equation}

\vspace{3mm}

Combining (\ref{march00025}), (\ref{march00028}), (\ref{march00029})
($\varepsilon=1/4$), we obtain

\vspace{-1mm}
$$
\left(\sum_{j_1, j_2, j_3=0}^{p}
C_{j_2 j_3 j_1 j_3 j_2 j_1}(s,\tau)\right)^2\le \frac{K_1(p+1)}{p^{3/2}}\le K^2,
$$

\vspace{3mm}
\noindent
where constants $K, K_1$ depend only on  $t, T.$
The inequality (\ref{marsixsix11}) is proved for the polynomial case.

Let us prove (\ref{marsixsix12}). In complete analogy
with the proof of (\ref{marsixsix11}) we have

\vspace{-1mm}
$$
\left(\sum_{j_1, j_2, j_3=0}^{p}
C_{j_3 j_1 j_2 j_3 j_2 j_1}(s,\tau)\right)^2\le
$$

\vspace{2mm}
$$
\le
(p+1)
\int\limits_{\tau}^s (s-t_5)
\int\limits_{\tau}^{t_5} 
\int\limits_{t_1}^{t_5}\left(\sum_{j_2=p+1}^{\infty}
C_{j_2}(\theta,t_1)C_{j_2}(t_5,\theta)\right)^2d\theta dt_1 dt_5.
$$

\vspace{4mm}
\noindent
The further proof is the same as in the case of 
(\ref{marsixsix11}).
The inequality (\ref{marsixsix12}) is proved.

Let us prove (\ref{marsixsix15}). By analogy
with the proof of (\ref{marsixsix11}) (see (\ref{march00025}})) we get

\vspace{-1mm}
$$
\left(\sum_{j_1, j_2, j_3=0}^{p}
C_{j_2 j_3 j_3 j_1 j_2 j_1}(s,\tau)\right)^2\le
$$

\vspace{2mm}
\begin{equation}
\label{march00030}
\le
(p+1)
\int\limits_{\tau}^s (s-t_5)
\int\limits_{\tau}^{t_5} 
\int\limits_{\tau}^{t_4}\left(\sum_{j_1=p+1}^{\infty}
C_{j_1}(\theta,\tau)C_{j_1}(t_4,\theta)\right)^2 d\theta dt_4 dt_5.
\end{equation}

\vspace{4mm}
\noindent
The further proof for the trigonometric case is the same as for
the inequality (\ref{marsixsix11}).

Consider the polynomial case. In this case,
we note that it is actually
necessary to consider the following two cases of (\ref{march00030})
\begin{equation}
\label{march00040}
1.\ \tau=t,\ \ \ 2.\ s=T.
\end{equation}

\vspace{3mm}

For Case~1, the estimate (\ref{march00028}) is simplified
as follows (see (\ref{otit6000x}), (\ref{after5000}) and (\ref{after1940}))

\vspace{-1mm}
\begin{equation}
\label{march00031}
\left|C_j(x,t)\right|=\left|
\int\limits_t^x\phi_{j}(\tau)d\tau
\right| <
\frac{C}{j^{1-\varepsilon/2}}\frac{1}{(1-z^2(x))^{1/4-\varepsilon/4}},
\end{equation}

\vspace{3mm}
\noindent
where notations are the same as in (\ref{march00028}).

Combining (\ref{march00030}), (\ref{march00028}), (\ref{march00029}), (\ref{march00031})
($\varepsilon=1/4$), we obtain

\vspace{-1mm}
\begin{equation}
\label{march00042aa}
\left(\sum_{j_1, j_2, j_3=0}^{p}
C_{j_2 j_3 j_3 j_1 j_2 j_1}(s,t)\right)^2\le
\frac{K_1(p+1)}{p^{3/2}}\le K^2,
\end{equation}

\vspace{3mm}
\noindent
where constants $K, K_1$ depend only on  $t, T.$
The inequality (\ref{marsixsix15}) is proved for the polynomial case
(Case~1).

Consider Case~2.
Combining (\ref{march00030}), (\ref{march00028}), (\ref{march00029})
($\varepsilon=1/4$), we ob\-tain

\vspace{-1mm}
$$
\left(\sum_{j_1, j_2, j_3=0}^{p}
C_{j_2 j_3 j_3 j_1 j_2 j_1}(T,\tau)\right)^2\le
\frac{K_1(p+1)}{p^{3/2}}\frac{1}{(1-z^2(\tau))^{3/8}}\le 
$$

\vspace{1mm}
$$
\le 
\frac{K^2}{(1-z^2(\tau))^{3/8}} \stackrel{\sf def}{=} F(\tau),
$$

\vspace{4mm}
\noindent
where constants $K, K_1$ depend only on  $t, T$
and $F(\tau)\in L_1([t, T])$ (integrable majorant (see above in this section)).
The following weakened version of the inequality (\ref{marsixsix15})

\vspace{-1mm}
\begin{equation}
\label{march00042aaa}
\left(\sum_{j_1, j_2, j_3=0}^{p}
C_{j_2 j_3 j_3 j_1 j_2 j_1}(T,\tau)\right)^2\le F(\tau)
\end{equation}

\vspace{3mm}
\noindent
is proved for the polynomial case
(Case~2),
where

\vspace{-1mm}
$$
F(\tau)=\frac{K^2}{(1-z^2(\tau))^{3/8}}.
$$

\vspace{3mm}

Let us prove (\ref{marsixsix2}).
Using the 
Cau\-chy--Bunyakovsky inequality as well as 
Fubini's Theorem and Parseval's equality, we have

\vspace{-1mm}
$$
\left(\sum_{j_1, j_2, j_3=0}^{p}
C_{j_3 j_3 j_1 j_2 j_2 j_1}(s,\tau)\right)^2=
\left(\sum_{j_3=0}^{p} 1\cdot \sum_{j_1, j_2=0}^{p}
C_{j_3 j_3 j_1 j_2 j_2 j_1}(s,\tau)\right)^2\le
$$

\vspace{2mm}
$$
\le \sum_{j_3=0}^{p} 1^2 \cdot \sum_{j_3=0}^{p}\left(\sum_{j_1, j_2=0}^{p}
C_{j_3 j_3 j_1 j_2 j_2 j_1}(s,\tau)\right)^2=
$$

\vspace{2mm}
$$
=
(p+1)\sum_{j_3=0}^{p}\left(\sum_{j_1, j_2=0}^{p}
C_{j_3 j_3 j_1 j_2 j_2 j_1}(s,\tau)\right)^2=
$$

\vspace{2mm}
$$
=         
(p+1)\sum_{j_3=0}^{p}\left(\sum_{j_1, j_2=0}^{p}
\int\limits_{\tau}^s \phi_{j_3}(t_6)
\int\limits_{\tau}^{t_6} \phi_{j_3}(t_5)
C_{j_1 j_2 j_2 j_1}(t_5,\tau) dt_5 dt_6
\right)^2\le
$$

\vspace{2mm}
$$
\le
(p+1)\sum_{j_3, j_3'=0}^{p}\left(
\int\limits_{\tau}^s \phi_{j_3}(t_6)
\int\limits_{\tau}^{t_6} \phi_{j_3'}(t_5)
\sum_{j_1, j_2=0}^{p}C_{j_1 j_2 j_2 j_1}(t_5,\tau) dt_5 dt_6
\right)^2\le
$$

\vspace{2mm}
$$
\le
(p+1)\sum_{j_3, j_3'=0}^{\infty}\left(
\int\limits_{\tau}^s \phi_{j_3}(t_6)
\int\limits_{\tau}^{t_6} \phi_{j_3'}(t_5)
\sum_{j_1, j_2=0}^{p}C_{j_1 j_2 j_2 j_1}(t_5,\tau) dt_5 dt_6
\right)^2=
$$

\vspace{2mm}
$$
=
(p+1)
\int\limits_{\tau}^s 
\int\limits_{\tau}^{t_6}\left(
\sum_{j_1, j_2=0}^{p}C_{j_1 j_2 j_2 j_1}(t_5,\tau) \right)^2dt_5 dt_6
=
$$

\vspace{2mm}
$$
=
(p+1)
\int\limits_{\tau}^s 
\int\limits_{\tau}^{t_6}\left(
\sum_{j_2=0}^{p} 1\cdot \sum_{j_1=0}^{p}C_{j_1 j_2 j_2 j_1}(t_5,\tau) \right)^2dt_5 dt_6
\le
$$

\vspace{2mm}
$$
\le
(p+1)^2
\int\limits_{\tau}^s 
\int\limits_{\tau}^{t_6}
\sum_{j_2=0}^{p}\left(
\sum_{j_1=0}^{p}C_{j_1 j_2 j_2 j_1}(t_5,\tau)\right)^2 dt_5 dt_6
=
$$

\vspace{2mm}
$$
=
(p+1)^2
\int\limits_{\tau}^s 
\int\limits_{\tau}^{t_6}
\sum_{j_2=0}^{p}\left(
\sum_{j_1=0}^{p}
\int\limits_{\tau}^{t_5}\phi_{j_2}(t_3)
\int\limits_{\tau}^{t_3}\phi_{j_2}(t_2)
C_{j_1}(t_2,\tau)C_{j_1}(t_5,t_3)dt_2 dt_3\right)^2  \hspace{-1.5mm} dt_5 dt_6
\le
$$

\vspace{2mm}
$$
\le
(p+1)^2
\int\limits_{\tau}^s 
\int\limits_{\tau}^{t_6}
\sum_{j_2,j_2'=0}^{p}\left(
\int\limits_{\tau}^{t_5}\phi_{j_2}(t_3)\times \right.
$$

\vspace{2mm}
$$
\left. \times
\int\limits_{\tau}^{t_3}\phi_{j_2'}(t_2)
\sum_{j_1=0}^{p}C_{j_1}(t_2,\tau)C_{j_1}(t_5,t_3)dt_2 dt_3\right)^2 dt_5 dt_6
\le
$$

\vspace{2mm}
$$
\le
(p+1)^2
\int\limits_{\tau}^s 
\int\limits_{\tau}^{t_6}
\sum_{j_2,j_2'=0}^{\infty}\left(
\int\limits_{\tau}^{t_5}\phi_{j_2}(t_3)\times \right.
$$

\vspace{2mm}
$$
\left. \times
\int\limits_{\tau}^{t_3}\phi_{j_2'}(t_2)
\left(\sum_{j_1=0}^{\infty} - \sum_{j_1=p+1}^{\infty}\right)
C_{j_1}(t_2,\tau)C_{j_1}(t_5,t_3)dt_2 dt_3\right)^2 dt_5 dt_6
=
$$

\vspace{2mm}
\begin{equation}
\label{march00032}
=
(p+1)^2
\int\limits_{\tau}^s 
\int\limits_{\tau}^{t_6}
\int\limits_{\tau}^{t_5}
\int\limits_{\tau}^{t_3}
\left(\sum_{j_1=p+1}^{\infty}
C_{j_1}(t_2,\tau)C_{j_1}(t_5,t_3)\right)^2 dt_2 dt_3 dt_5 dt_6.
\end{equation}

\vspace{5mm}

Consider the trigonometric case.
Combining (\ref{march00032}), (\ref{march00026}), (\ref{march00027}}), we obtain

\vspace{-1mm}
$$
\left(\sum_{j_1, j_2, j_3=0}^{p}
C_{j_3 j_3 j_1 j_2 j_2 j_1}(s,\tau)\right)^2
\le \frac{K_1(p+1)^2}{p^2}\le K^2,
$$

\vspace{3mm}
\noindent
where constants $K, K_1$ depend only on  $t, T.$
The inequality (\ref{marsixsix2}) is proved for the trigonometric case.

Consider the polynomial case for two cases (\ref{march00040}). 
Let $\tau=t.$ The modification of the estimate
(\ref{march00028}) for $\varepsilon =0$ is as follows

\vspace{-1mm}
\begin{equation}
\label{march00042}
\left|C_j(x,v)\right|=\left|
\int\limits_v^x\phi_{j}(\tau)d\tau
\right| <
\frac{C}{j}\Biggl(\frac{1}{(1-z^2(x))^{1/4}}+
\frac{1}{(1-z^2(v))^{1/4}}\Biggr),
\end{equation}

\vspace{3mm}
\noindent
where 
$j\in \mathbb{N},$ $z(x), z(v)\in (-1, 1)$ 
($z(x)$ is defined by (\ref{zz1})),
$x, v\in (t, T),$
constant $C$ does not depend on $j$.

For $v=t$, the estimate (\ref{march00042}) is simplified
as follows (see (\ref{otit6000x}), (\ref{101xx}))

\vspace{-1mm}
\begin{equation}
\label{march00043}
\left|C_j(x,t)\right|=\left|
\int\limits_t^x\phi_{j}(\tau)d\tau
\right| <
\frac{C}{j(1-z^2(x))^{1/4}},
\end{equation}

\vspace{3mm}
\noindent
where notations are the same as in (\ref{march00042}).

Combining (\ref{march00032}), (\ref{march00042}), (\ref{march00043}), we get

\vspace{-1mm}
$$
\left(\sum_{j_1, j_2, j_3=0}^{p}
C_{j_3 j_3 j_1 j_2 j_2 j_1}(s,t)\right)^2
\le \frac{K_1(p+1)^2}{p^2}\le K^2,
$$

\vspace{3mm}
\noindent
where constants $K, K_1$ depend only on  $t, T.$
The inequality (\ref{marsixsix2}) is proved for the polynomial case 
($\tau=t$).

Now let $s=T.$
Combining (\ref{march00032}) and (\ref{march00042}), we ob\-tain

\vspace{-1mm}
$$
\left(\sum_{j_1, j_2, j_3=0}^{p}
C_{j_3 j_3 j_1 j_2 j_2 j_1}(T,\tau)\right)^2\le
\frac{K_1(p+1)^2}{p^2}\frac{1}{(1-z^2(\tau))^{1/2}}\le 
$$

\vspace{1mm}
$$
\le 
\frac{K^2}{(1-z^2(\tau))^{1/2}} \stackrel{\sf def}{=} F(\tau),
$$

\vspace{4mm}
\noindent
where constants $K, K_1$ depend only on  $t, T$
and $F(\tau)\in L_1([t, T])$ (integrable majorant (see above in this section)).
The following weakened version of the inequality (\ref{marsixsix2})

\vspace{-1mm}
\begin{equation}
\label{march00042aaaa}
\left(\sum_{j_1, j_2, j_3=0}^{p}
C_{j_3 j_3 j_1 j_2 j_2 j_1}(T,\tau)\right)^2\le F(\tau)
\end{equation}

\vspace{3mm}
\noindent
is proved for the polynomial case
($s=T$),
where

\vspace{-1mm}
$$
F(\tau)=\frac{K^2}{(1-z^2(\tau))^{1/2}}.
$$

\vspace{5mm}

Finally, we prove the inequality
(\ref{marsixsix1}). By analogy with (\ref{march00032}) we get

\vspace{-1mm}
$$
\left(\sum_{j_1, j_2, j_3=0}^{p}
C_{j_3 j_3 j_2 j_1 j_2 j_1}(s,\tau)\right)^2\le
$$

\vspace{2mm}
$$
\le
(p+1)\sum_{j_3=0}^{p}\left(\sum_{j_1, j_2=0}^{p}
C_{j_3 j_3 j_2 j_1 j_2 j_1}(s,\tau)\right)^2=
$$

\vspace{2mm}
$$
=         
(p+1)\sum_{j_3=0}^{p}\left(\sum_{j_1, j_2=0}^{p}
\int\limits_{\tau}^s \phi_{j_3}(t_6)
\int\limits_{\tau}^{t_6} \phi_{j_3}(t_5)
C_{j_2 j_1 j_2 j_1}(t_5,\tau) dt_5 dt_6
\right)^2\le
$$

\vspace{2mm}
$$
\le
(p+1)\sum_{j_3, j_3'=0}^{\infty}\left(
\int\limits_{\tau}^s \phi_{j_3}(t_6)
\int\limits_{\tau}^{t_6} \phi_{j_3'}(t_5)
\sum_{j_1, j_2=0}^{p}C_{j_2 j_1 j_2 j_1}(t_5,\tau) dt_5 dt_6
\right)^2=
$$

\vspace{2mm}
$$
=
(p+1)
\int\limits_{\tau}^s 
\int\limits_{\tau}^{t_6}\left(
\sum_{j_1, j_2=0}^{p}C_{j_2 j_1 j_2 j_1}(t_5,\tau) \right)^2dt_5 dt_6
=
$$

\vspace{2mm}
$$
\le
(p+1)^2
\int\limits_{\tau}^s 
\int\limits_{\tau}^{t_6}
\sum_{j_2=0}^{p}\left(
\sum_{j_1=0}^{p}C_{j_2 j_1 j_2 j_1}(t_5,\tau)\right)^2 dt_5 dt_6
=
$$

\vspace{2mm}
$$
=
(p+1)^2
\int\limits_{\tau}^s 
\int\limits_{\tau}^{t_6}
\sum_{j_2=0}^{p}\left(
\sum_{j_1=0}^{p}
\int\limits_{\tau}^{t_5}\phi_{j_2}(t_4)
\int\limits_{\tau}^{t_4}\phi_{j_2}(t_2)
C_{j_1}(t_2,\tau)C_{j_1}(t_4,t_2)dt_2 dt_4\right)^2  \hspace{-1.5mm} dt_5 dt_6
\le
$$

\vspace{2mm}
$$
\le
(p+1)^2
\int\limits_{\tau}^s 
\int\limits_{\tau}^{t_6}
\sum_{j_2,j_2'=0}^{\infty}\left(
\int\limits_{\tau}^{t_5}\phi_{j_2}(t_4)
\times \right.
$$

\vspace{2mm}
$$
\left. \times
\int\limits_{\tau}^{t_4}
\phi_{j_2'}(t_2)
\left(\sum_{j_1=0}^{\infty} - \sum_{j_1=p+1}^{\infty}\right)
C_{j_1}(t_2,\tau)C_{j_1}(t_4,t_2)dt_2 dt_4\right)^2 dt_5 dt_6
=
$$

\vspace{2mm}
\begin{equation}
\label{march00044}
=
(p+1)^2
\int\limits_{\tau}^s 
\int\limits_{\tau}^{t_6}
\int\limits_{\tau}^{t_5}
\int\limits_{\tau}^{t_4}
\left(\sum_{j_1=p+1}^{\infty}
C_{j_1}(t_2,\tau)C_{j_1}(t_4,t_2)\right)^2 dt_2 dt_4 dt_5 dt_6.
\end{equation}

\vspace{5mm}

The further proof of inequality (\ref{marsixsix1}) for the trigonometric case and the 
weakened analogue of inequality (\ref{marsixsix1}) for the polynomial case is 
completely analogous to the proof of (\ref{marsixsix15}) and its weakened 
analogue (see (\ref{march00030}), (\ref{march00042aa}), (\ref{march00042aaa})).

Thus, the following theorem is proved.

\vspace{2mm}

{\bf Theorem~62.}\ {\it Suppose that
$\{\phi_j(x)\}_{j=0}^{\infty}$ is a complete orthonormal
system of Legendre polynomials or tri\-go\-no\-met\-ric functions
in the space $L_2([t, T]).$
Then$,$ for the iterated Stra\-to\-no\-vich stochastic integral
of seventh multiplicity

$$
J^{*}[\psi^{(7)}]_{T,t}=
{\int\limits_t^{*}}^T
\ldots
{\int\limits_t^{*}}^{t_2}
d{\bf w}_{t_1}^{(i_1)}
\ldots d{\bf w}_{t_7}^{(i_7)}
$$

\vspace{1mm}
\noindent
the following 
expansion 

\vspace{-2mm}
$$
J^{*}[\psi^{(7)}]_{T,t}=
\hbox{\vtop{\offinterlineskip\halign{
\hfil#\hfil\cr
{\rm l.i.m.}\cr
$\stackrel{}{{}_{p\to \infty}}$\cr
}} }
\sum\limits_{j_1,\ldots, j_7=0}^{p}
C_{j_7 \ldots j_1}\zeta_{j_1}^{(i_1)}\ldots \zeta_{j_7}^{(i_7)}
$$

\vspace{4mm}
\noindent
that converges in the mean-square sense is valid, where 
$i_1,\ldots,i_7=0, 1,\ldots,m,$

\vspace{-1mm}
$$
C_{j_7\ldots j_1}=\int\limits_t^T
\phi_{j_7}(t_7)
\ldots
\int\limits_t^{t_2}
\phi_{j_1}(t_1)dt_1\ldots dt_7
$$

\vspace{2mm}
\noindent
and
$$
\zeta_{j}^{(i)}=
\int\limits_t^T \phi_{j}(\tau) d{\bf w}_{\tau}^{(i)}
$$ 

\vspace{2mm}
\noindent
are independent standard Gaussian random variables for various 
$i$ or $j$ {\rm (}in the case when $i\ne 0${\rm ),}
${\bf w}_{\tau}^{(i)}={\bf f}_{\tau}^{(i)}$ for
$i=1,\ldots,m$ and 
${\bf w}_{\tau}^{(0)}=\tau.$}

\vspace{5mm}

\section{Expansion of Iterated Stratonovich Stochastic Integrals
of Multiplicity 8 for the Case $\psi_1(\tau),\ldots, \psi_8(\tau)
\equiv 1$ (The Cases of Legendre 
Polynomials and Trigonometric Functions)}

\vspace{5mm}

This section is devoted to the following theorem.

\vspace{2mm}

{\bf Theorem~63.}\ {\it Suppose that
$\{\phi_j(x)\}_{j=0}^{\infty}$ is a complete orthonormal
system of Legendre polynomials or tri\-go\-no\-met\-ric functions
in the space $L_2([t, T]).$
Then$,$ for the iterated Stra\-to\-no\-vich stochastic integral
of eighth multiplicity 

$$
J^{*}[\psi^{(8)}]_{T,t}=
{\int\limits_t^{*}}^T
\ldots
{\int\limits_t^{*}}^{t_2}
d{\bf w}_{t_1}^{(i_1)}
\ldots d{\bf w}_{t_8}^{(i_8)}
$$

\vspace{1mm}
\noindent
the following 
expansion 

\vspace{-2mm}
$$
J^{*}[\psi^{(8)}]_{T,t}=
\hbox{\vtop{\offinterlineskip\halign{
\hfil#\hfil\cr
{\rm l.i.m.}\cr
$\stackrel{}{{}_{p\to \infty}}$\cr
}} }
\sum\limits_{j_1,\ldots, j_8=0}^{p}
C_{j_8 \ldots j_1}\zeta_{j_1}^{(i_1)}\ldots \zeta_{j_8}^{(i_8)}
$$

\vspace{4.5mm}
\noindent
that converges in the mean-square sense is valid, where 
$i_1,\ldots,i_8=0, 1,\ldots,m,$

$$
C_{j_8\ldots j_1}=\int\limits_t^T
\phi_{j_8}(t_8)
\ldots
\int\limits_t^{t_2}
\phi_{j_1}(t_1)dt_1\ldots dt_8
$$

\vspace{2mm}
\noindent
and
$$
\zeta_{j}^{(i)}=
\int\limits_t^T \phi_{j}(\tau) d{\bf w}_{\tau}^{(i)}
$$ 

\vspace{2mm}
\noindent
are independent standard Gaussian random variables for various 
$i$ or $j$ {\rm (}in the case when $i\ne 0${\rm ),}
${\bf w}_{\tau}^{(i)}={\bf f}_{\tau}^{(i)}$ for
$i=1,\ldots,m$ and 
${\bf w}_{\tau}^{(0)}=\tau.$}

\vspace{2mm}

{\bf Proof.}\ To prove the theorem, we need to check
the condition (\ref{09091}) (or its weakened version)
for the case $k=8>2r$, where $r=1, 2, 3$  (see Theorem~59).
Recall that the case
$k=2r$ is considered in Sect.~26 (see (\ref{july90000})).
Under the conditions of Theorem~63, this means that
$k=8=2r$, where $r=4$.
The relations
(\ref{2024december12})--(\ref{2024december11}), (\ref{cc123})
cover the case $k=8,$ $r=1,2$ (see (\ref{09091})).

Thus, it remains to consider the case $k=8,$ $r=3.$
The case $k=7,$ $r=3$ was considered in the previous section.
Here we will focus on the differences
between these two cases.

Since now $k=8,$ then along with inequalities
(\ref{march0001})--(\ref{march0004}), it is necessary to prove the following 
inequalities

\vspace{-1mm}
\begin{equation}
\label{march00045}
\left|\sum\limits_{j_{g_1},j_{g_3}. j_{g_5}=0}^p
\bigl(C_{j_{d_3} j_{d_3-1}j_{d_3-2}j_{d_3-3}}(s,\tau)
C_{j_{d_2}}(\theta,u) C_{j_{d_1}}(\rho,v)
\bigr)\biggl|_{j_{g_1}=j_{g_2},j_{g_3}=j_{g_4},j_{g_5}=j_{g_6}}\right|\le 
K<\infty,
\end{equation}

\vspace{2mm}
\begin{equation}
\label{march00046}
\left|\sum\limits_{j_{g_1},j_{g_3}. j_{g_5}=0}^p
\bigl(C_{j_{d_3} j_{d_3-1}j_{d_3-2}}(s,\tau)
C_{j_{d_2}j_{d_2-1}}(\theta,u) C_{j_{d_1}}(\rho,v)
\bigr)\biggl|_{j_{g_1}=j_{g_2},j_{g_3}=j_{g_4},j_{g_5}=j_{g_6}}\right|\le 
K<\infty,
\end{equation}

\vspace{2mm}
\begin{equation}
\label{march00047}
\left|\sum\limits_{j_{g_1},j_{g_3}. j_{g_5}=0}^p
\bigl(C_{j_{d_3} j_{d_3-1}}(s,\tau)
C_{j_{d_2}j_{d_2-1}}(\theta,u) C_{j_{d_1} j_{d_1-1}}(\rho,v)
\bigr)\biggl|_{j_{g_1}=j_{g_2},j_{g_3}=j_{g_4},j_{g_5}=j_{g_6}}\right|\le 
K<\infty,
\end{equation}

\vspace{5mm}
\noindent
where $p\in\mathbb{N},$ $t\le \tau < s \le T,$ $t\le u<\theta \le T,$
$t\le v<\rho \le T,$
constant $K$ does not depend on 
$p, s, \tau, \theta, u, \rho, v$ (but only on $t, T$) and may differ from line to line;
another notations are the same as in Sect.~30.

The inequalities (\ref{march00045})--(\ref{march00047})
are proved using the same technique as 
inequalities (\ref{2024december12})--(\ref{2024december11}) (see Sect.~31).
Here we will only prove as an example the following
special case of the inequality (\ref{march00047})
\begin{equation}
\label{march00048}
\left|\sum\limits_{j_1, j_2, j_3=0}^p C_{j_2 j_1}(s,\tau)C_{j_3 j_1}(\theta,u)
C_{j_2 j_3}(\rho,v)\right|
\le K<\infty.
\end{equation}

\vspace{3mm}

Using the 
Cau\-chy--Bunyakovsky inequality as well as 
Fubini's Theorem, Parseval's equality and (\ref{2024december13}), we have

\vspace{-1mm}
$$
\left(
\sum\limits_{j_1, j_2, j_3=0}^p C_{j_2 j_1}(s,\tau)C_{j_3 j_1}(\theta,u)
C_{j_2 j_3}(\rho,v)\right)^2=
$$

\vspace{2mm}
$$
=\left(\sum\limits_{j_2, j_3=0}^p  C_{j_2 j_3}(\rho,v)
\sum\limits_{j_1=0}^p C_{j_2 j_1}(s,\tau)C_{j_3 j_1}(\theta,u)
\right)^2\le
$$

\vspace{2mm}
$$
\le \sum\limits_{j_2, j_3=0}^p  C_{j_2 j_3}^2(\rho,v)
\sum\limits_{j_2, j_3=0}^p \left(\sum\limits_{j_1=0}^p C_{j_2 j_1}(s,\tau)C_{j_3 j_1}(\theta,u)
\right)^2\le
$$

\vspace{2mm}
$$
\le \sum\limits_{j_2, j_3=0}^{\infty}  C_{j_2 j_3}^2(\rho,v)
\sum\limits_{j_2, j_3=0}^{\infty} \left(\sum\limits_{j_1=0}^p C_{j_2 j_1}(s,\tau)C_{j_3 j_1}(\theta,u)
\right)^2=
$$

\vspace{2mm}
$$
=\frac{(\rho-v)^2}{2}
\sum\limits_{j_2, j_3=0}^{\infty} 
\left(\sum\limits_{j_1=0}^p \int\limits_{\tau}^s \phi_{j_2}(t_2)\int\limits_{\tau}^{t_2} \phi_{j_1}(t_1)
dt_1 dt_2
\int\limits_{u}^{\theta} \phi_{j_3}(t_4)\int\limits_{u}^{t_4} \phi_{j_1}(t_3)
dt_3 dt_4
\right)^2=
$$

\vspace{2mm}
$$
=\frac{(\rho-v)^2}{2}
\sum\limits_{j_2, j_3=0}^{\infty} \left(
\int\limits_{\tau}^s \int\limits_{u}^{\theta} \phi_{j_2}(t_2)\phi_{j_3}(t_4)\times
\right.
$$

\vspace{2mm}
$$
\left.\times
\sum\limits_{j_1=0}^p 
\int\limits_{\tau}^{t_2} \phi_{j_1}(t_1)
dt_1 \int\limits_{u}^{t_4} \phi_{j_1}(t_3)
dt_1 dt_3 dt_4 dt_2
\right)^2=
$$

\vspace{2mm}
$$
=\frac{(\rho-v)^2}{2}
\int\limits_{\tau}^s \int\limits_{u}^{\theta} 
\left(
\sum\limits_{j_1=0}^p  C_{j_1}(t_2,\tau) C_{j_1}(t_4,u)
\right)^2 dt_4 dt_2\le
$$

\vspace{2mm}
$$
\le
K_1^2 \frac{(\rho-v)^2}{2} (s-\tau)(\theta-u)\le
K_1^2 \frac{(T-t)^4}{2}=K.
$$

\vspace{5mm}
\noindent
The inequality (\ref{march00048}) is proved.

The inequalities (\ref{march0001})--(\ref{march0004})
for the case $k=8$ 
are proved similarly to the 
inequalities (\ref{march0001})--(\ref{march0004})
for the case $k=7$ (see Sect.~35).
There will be minor differences only when proving
(\ref{march0001}) for the case $k=8$ (polynomial case).
The above differences will be due to 
the fact that along with the two cases (\ref{march00040}) the following third
case 
$$
\tau, s \in (t, T)
$$

\vspace{3mm}
\noindent
will now appear when proving 
(\ref{marsixsix1}), (\ref{marsixsix2}), (\ref{marsixsix15}).

Using the technique that led to the estimates (\ref{march00042aaa}),
(\ref{march00042aaaa}),
we obtain for Case~3

\vspace{-1mm}
$$
\left(\sum_{j_1, j_2, j_3=0}^{p}
C_{j_3 j_3 j_2 j_1 j_2 j_1}(s,\tau)
\right)^2\le
\frac{K^2}{(1-z^2(\tau))^{3/8}}\stackrel{\sf def}
{=}F(\tau)\ \ \ (\hbox{for}\ (\ref{marsixsix1})),
$$

\vspace{3mm}
$$
\left(\sum_{j_1, j_2, j_3=0}^{p}
C_{j_3 j_3 j_1 j_2 j_2 j_1}(s,\tau)
\right)^2\le
\frac{K^2}{(1-z^2(\tau))^{1/2}}\stackrel{\sf def}
{=}F(\tau)\ \ \ (\hbox{for}\ (\ref{marsixsix2})),
$$

\vspace{3mm}
$$
\left(\sum_{j_1, j_2, j_3=0}^{p}
C_{j_2 j_3 j_3 j_1 j_2 j_1}(s,\tau)\right)^2\le 
\frac{K^2}{(1-z^2(\tau))^{3/8}}
\stackrel{\sf def}{=}F(\tau)\ \ \ (\hbox{for}\ (\ref{marsixsix15})),
$$

\vspace{5mm}
\noindent
where constant $K$ depends only on  $t, T$
and $F(\tau)\in L_1([t, T])$ (integrable majorant).
Theorem~63 is proved.

\vspace{5mm}

\section{Convergence of the Expansion (\ref{march000195})
to the Iterated Stratonovich Stochastic Integrals
in the Sense of Mathematical Expectation}

\vspace{5mm}

In the previous sections, we actually proved
that the value

\vspace{-1mm}
$$
\sum\limits_{j_1,\ldots,j_k=0}^{p}
C_{j_k \ldots j_1}\prod\limits_{l=1}^k \zeta_{j_l}^{(i_l)}
$$

\vspace{3mm}
\noindent
converges if $p\to\infty$ (under suitable conditions)
to the iterated Stratonovich stochastic integrals (\ref{str})
in the sense of mathematical expectation.
Let us explain this fact in more detail.

Suppose that $\psi_1(\tau),\ldots,\psi_k(\tau)$ ($k\in\mathbb{N}$) are
continuous functions on $[t, T]$ and consider Theorem~5.
First, let $k=2q+1,$ $q\in \mathbb{N}.$ We represent (w.~p.~1)
each stochastic integral $J[\psi^{(k)}]_{T,t}^{s_r,\ldots,s_1}$
from the right-hand side of (\ref{30.4})
using the transformation (\ref{febr15})
as a finite linear combination
of the iterated Ito stochastic integrals. Thus, we have (see (\ref{30.4}))

\vspace{-1mm}
\begin{equation}
\label{march000196}
{\sf M}\left\{J^{*}[\psi^{(k)}]_{T,t}\right\}=0,
\end{equation}

\vspace{3mm}
\noindent
where $J^{*}[\psi^{(k)}]_{T,t}$ is defined by
(\ref{str}). On the other hand,

\vspace{-1mm}
\begin{equation}
\label{march000197}
{\sf M}\left\{\sum\limits_{j_1,\ldots,j_k=0}^{p}
C_{j_k \ldots j_1}\prod\limits_{l=1}^k \zeta_{j_l}^{(i_l)}
\right\}=0,
\end{equation}

\vspace{3mm}
\noindent
since $\zeta_{j_l}^{(i_l)}$ has Gaussian distribution and
$k=2q+1,$ $q\in \mathbb{N}.$

Combining (\ref{march000196}) and (\ref{march000197}), we obtain

\vspace{-1mm}
\begin{equation}
\label{march000301}
\lim\limits_{p\to\infty}\left|{\sf M}\left\{
J^{*}[\psi^{(k)}]_{T,t}-
\sum\limits_{j_1,\ldots,j_k=0}^{p}
C_{j_k \ldots j_1}\prod\limits_{l=1}^k \zeta_{j_l}^{(i_l)}\right\}\right|=0.
\end{equation}

\vspace{3mm}

Now let $k=2q,$ $q\in \mathbb{N}.$ In this case, 
using the above reasoning, we get (see (\ref{30.4}))

$$
{\sf M}\left\{J^{*}[\psi^{(k)}]_{T,t}\right\}=
$$

\vspace{2mm}
$$
=
\frac{1}{2^q}{\bf 1}_{\{i_1=i_2\ne 0\}}
{\bf 1}_{\{i_3=i_4\ne 0\}}\ldots {\bf 1}_{\{i_{2q-1}=i_{2q}\ne 0\}}\times
$$

\vspace{2mm}
\begin{equation}
\label{march000199}
\times
\int\limits_t^T\psi_{2q}(t_{2q})\psi_{2q-1}(t_{2q}) \ldots 
\int\limits_t^{t_6}
\psi_4(t_4)\psi_3(t_4)
\int\limits_t^{t_4}
\psi_2(t_2)\psi_1(t_2)dt_2 dt_4\ldots  dt_{2q}.
\end{equation}

\vspace{4mm}

Recall that the multiple Wiener stochastic 
integral (\ref{mult11www}) has zero expectation.
Then, using (\ref{after8xxds1}), (\ref{july30016}) and (\ref{march000199}), we have

\vspace{-1mm}
$$
\lim\limits_{p\to\infty}
{\sf M}\left\{\sum\limits_{j_1,\ldots,j_k=0}^{p}
C_{j_k \ldots j_1}\prod\limits_{l=1}^k \zeta_{j_l}^{(i_l)}\right\}=
$$

\vspace{2mm}
$$
=
{\bf 1}_{\{i_1=i_2\ne 0\}}{\bf 1}_{\{i_3=i_4\ne 0\}}
\ldots {\bf 1}_{\{i_{2q-1}=i_{2q}\ne 0\}}
\lim\limits_{p\to\infty}
\sum_{j_q,j_{q-2},\ldots, j_2=0}^{p}
C_{j_q j_q j_{q-2} j_{q-2} \ldots j_2 j_2}=
$$

\vspace{2mm}
$$
=
\frac{1}{2^q}{\bf 1}_{\{i_1=i_2\ne 0\}}
{\bf 1}_{\{i_3=i_4\ne 0\}}\ldots {\bf 1}_{\{i_{2q-1}=i_{2q}\ne 0\}}\times
$$

\vspace{2mm}
$$
\times
\int\limits_t^T\psi_{2q}(t_{2q})\psi_{2q-1}(t_{2q}) \ldots 
\int\limits_t^{t_6}
\psi_4(t_4)\psi_3(t_4)
\int\limits_t^{t_4}
\psi_2(t_2)\psi_1(t_2)dt_2 dt_4\ldots  dt_{2q}=
$$

\vspace{2mm}
\begin{equation}
\label{march000300}
={\sf M}\left\{J^{*}[\psi^{(k)}]_{T,t}\right\}.
\end{equation}

\vspace{5mm}

Applying (\ref{march000300}), we obtain

$$
\lim\limits_{p\to\infty}
\left|{\sf M}\left\{
J^{*}[\psi^{(k)}]_{T,t}-\sum\limits_{j_1,\ldots,j_k=0}^{p}
C_{j_k \ldots j_1}\prod\limits_{l=1}^k \zeta_{j_l}^{(i_l)}\right\}\right|=
$$

\vspace{3mm}
$$
=\left|
{\sf M}\left\{
J^{*}[\psi^{(k)}]_{T,t}\right\}-
\lim\limits_{p\to\infty}
{\sf M}\left\{\sum\limits_{j_1,\ldots,j_k=0}^{p}
C_{j_k \ldots j_1}\prod\limits_{l=1}^k \zeta_{j_l}^{(i_l)}\right\}\right|=
$$

\vspace{2mm}
$$
=0.
$$

\vspace{4mm}

The equality (\ref{march000301}) is proved.

\vspace{5mm}

\section{Expansion of Iterated Stratonovich Stochastic Integrals of Arbitrary
Multiplicity. Proof of Hypothesis 8 for the Case $p_1=\ldots=p_k=p,$ 
$\psi_1(\tau),\ldots,\psi_k(\tau)\in C^1([t, T]),$ $k\in\mathbb{N}$}

\vspace{5mm}
                                   
{\bf Theorem~64} \cite{2018a}, \cite{19999999}.\ {\it
Let $\{\phi_j(x)\}_{j=0}^{\infty}$ be an arbitrary 
complete orthonormal system of functions in $L^2([t, T]),$
and let $\psi_1(\tau) \equiv \dots \equiv \psi_k(\tau) \equiv 1$. 
For any multiplicity $k \in \mathbb{N}$ 
and any indices $i_1, \dots, i_k \in \{1, \dots, d\},$
the partial sum 

\vspace{-1mm}
\begin{equation}
\label{ud11}
S_{T,t}^{(i_1, \dots, i_k)p}
=\sum_{j_1,\ldots,j_k=0}^{p}
C_{j_k\ldots j_1}
\prod_{l=1}^k \zeta_{j_l}^{(i_l)}
\end{equation}

\vspace{3mm}
\noindent
with equal truncation limits $p_1 = \dots = p_k = p$ 
converges in the mean-square sense 
to the iterated Stratonovich stochastic integral 
$J^{*}[\psi^{(k)}]_{T,t}^{(i_1\ldots i_k)}$ as $p \to \infty:$

\vspace{-1mm}
$$
\lim_{p \to \infty} 
{\sf M}\left\{
\left| S_{T,t}^{(i_1, \dots, i_k)p} 
- J^{*}[\psi^{(k)}]_{T,t}^{(i_1\ldots i_k)} \right|^2 
\right\} = 0,
$$

\vspace{3mm}
\noindent
where $J^{*}[\psi^{(k)}]_{T,t}^{(i_1\ldots i_k)}$ is defined by {\rm (\ref{str})}
for $i_1,\ldots,i_k=1,\ldots,d,$ 
$C_{j_k \ldots j_1}$ has the form {\rm (\ref{after3000}),} 
and
$$
\zeta_{j}^{(i)}=
\int\limits_t^T \phi_{j}(\tau) d{\bf w}_{\tau}^{(i)}\ \ \ (i=1,\ldots,d,\ j=0,1,2,\ldots)
$$ 

\vspace{2mm}
\noindent
are independent standard Gaussian random variables for various 
$i$ or $j.$ 
}

\vspace{2mm}

{\bf Proof.}\
Without loss of generality, 
by a linear time rescaling, 
we assume that $[t, T] = [0, 1]$, 
since all variation norms and estimates 
are invariant under linear time changes.

\vspace{2mm}

{\bf Step~1. Identification with the Karhunen--Lo\`eve lift.}
Let ${\bf w}_{\tau} = \bigl({\bf w}_{\tau}^{(1)}, \dots, {\bf w}_{\tau}^{(d)}\bigr)$, 
$\tau \in [0, 1]$, be a standard 
$d$-dimensional Wiener process with independent components and  ${\bf w}_0 = 0$ w.~p.~1. 
For each coordinate $\mathbf{w}^{(i)}$
($i = 1, \dots, d$), let $\mathcal{K}_i \subset L^2(\Omega)$ denote the
closed linear Gaussian subspace spanned by $\{\mathbf{w}_{\tau}^{(i)} : \tau \in [0,1]\}$.
Following Friz and Victoir (2010, Sect.~3.2) \cite{Friz}, the Cameron--Martin space
$\mathcal{H}_i$ associated with $\mathbf{w}^{(i)}$ is the Hilbert space of
paths $h : [0,1] \to \mathbb{R}$ of the form

\vspace{-1mm}
$$
h(\tau) = \mathsf{M}\{Z \mathbf{w}_{\tau}^{(i)}\}, \quad \tau \in [0, 1]
$$

\vspace{2mm}
\noindent
for some $Z \in \mathcal{K}_i$, equipped with the inner product

\vspace{-1mm}
$$
\langle h, g \rangle_{\mathcal{H}_i} := \mathsf{M}\{Z Z'\},
$$

\vspace{2mm}
\noindent
where $g(\tau) = \mathsf{M}\{Z' \mathbf{w}_{\tau}^{(i)}\}$ with $Z' \in \mathcal{K}_i$.

Since every element of $\mathcal{K}_i$ admits the Wiener integral representation

\vspace{-1mm}
$$
Z = \int\limits_0^1 f_1(s)\, d\mathbf{w}_s^{(i)}, \ \ \ 
Z' = \int\limits_0^1 f_2(s)\, d\mathbf{w}_s^{(i)}
$$

\vspace{2mm}
\noindent
for unique functions $f_1, f_2 \in L^2([0,1])$, the Ito isometry yields

$$
h(\tau) = \mathsf{M}\left\{\int\limits_0^1 f_1(s)\, d\mathbf{w}_s^{(i)}
\int\limits_0^{\tau} 1\, d\mathbf{w}_s^{(i)}\right\} = \int\limits_0^{\tau} f_1(s)\, ds,
$$

\vspace{2mm}
\noindent
and similarly 
$$
g(\tau) = \int\limits_0^{\tau} f_2(s) ds,
$$

\vspace{2mm}
\noindent
with inner product

\vspace{-1mm}
$$
\langle h, g \rangle_{\mathcal{H}_i}
= \mathsf{M}\left\{\int\limits_0^1 f_1(s)\, d\mathbf{w}_s^{(i)}
\int\limits_0^1 f_2(s)\, d\mathbf{w}_s^{(i)}\right\}
= \int\limits_0^1 f_1(s) f_2(s)ds.
$$

\vspace{2mm}

Defining the antiderivatives 

\vspace{-1mm}
$$
h_j(\tau) = \int\limits_0^{\tau} \phi_j(s)ds,\ \ \
\tau \in [0,\,1],\ \ \  j \ge 0,
$$

\vspace{2mm}
\noindent
we obtain 
$$
\langle h_j, h_m \rangle_{\mathcal{H}_i} 
= \int\limits_0^1 \phi_j(s)\phi_m(s)\,ds = \delta_{j,m}.
$$ 

\vspace{2mm}

Hence, $\{h_j\}_{j=0}^{\infty}$ forms 
an orthonormal basis of $\mathcal{H}_i$, 
and the vectors with $h_j$ in the $i$-th slot 
form an orthonormal basis of 
$$
\mathcal{H} = \bigoplus\limits_{i=1}^d \mathcal{H}_i.
$$

\vspace{3mm}

In Sect.~6 of Friz and Victoir (2010) \cite{Friz}, 
the general $L^2$-approximation of a Gaussian process $X_{\tau}$ 
for a subset $A \subset \mathbb{N}$ is defined by

\vspace{-1mm}
$$
X_{\tau}^A := {\sf M}\bigl\{ X_{\tau} \mid \mathcal{F}_A \bigr\} 
= \sum_{j \in A} Z_j h_j(\tau),\ \ \
\mathcal{F}_A = \sigma(Z_j : j \in A).
$$

\vspace{2mm}

Setting $X_{\tau} = {\bf w}_{\tau}^{(i)}$, $A = A_p := \{0, 1, \dots, p\}$, 
and 

\vspace{-1mm}
$$
Z_j = \zeta_j^{(i)} := \int\limits_0^1 \phi_j(s)d{\bf w}_{s}^{(i)},
$$ 

\vspace{2mm}
\noindent
this specializes to the absolutely continuous 
componentwise approximation (also see (\ref{um1xx}) for $t=0$)

\vspace{-1mm}
$$
{\bf w}_{\tau}^{(i)p}
:= {\sf M}\bigl\{ {\bf w}_{\tau}^{(i)} \mid \mathcal{F}_{A_p}^{(i)} \bigr\}
= \sum_{j=0}^{p} \zeta_j^{(i)} h_j(\tau),\ \ \ 
i = 1, \dots, d,
$$

\vspace{3mm}
\noindent
where $\mathcal{F}_{A_p}^{(i)} = \sigma\bigl(\zeta_j^{(i)} : 0 \le j \le p\bigr)$. 
We denote the resulting $d$-dimensional absolutely continuous approximation 
(having sample paths in $C^{1\text{-var}}([0,1],\mathbb{R}^d)$ w.~p.~1) by 

\vspace{-1mm}
$$
{\bf w}_{\tau}^{p} := \bigl({\bf w}_{\tau}^{(1)p}, \dots, {\bf w}_{\tau}^{(d)p}\bigr),\ \ \
\tau \in [0, 1].
$$

\vspace{2mm}

Next, for any path $x \in C^{1\text{-var}}([0,\,1], \mathbb{R}^d)$, 
its canonical step-$N$ lift $S_N(x)$ 
into the truncated tensor algebra 
$$
T^N(\mathbb{R}^d) = \bigoplus\limits_{m=0}^N (\mathbb{R}^d)^{\otimes m}
$$ 

\vspace{3mm}
\noindent
is given by \cite{Friz} (Friz and Victoir, 2010, Sect.~1.1):

$$
S_N(x)_{\tau} 
= 1 + \sum_{m=1}^{N}\
\int\limits_{0 < s_1 < \dots < s_m < \tau} 
dx_{s_1} \otimes \dots \otimes dx_{s_m}.
$$

\vspace{3mm}

Let $\{e_1, \dots, e_d\}$ be the standard orthonormal basis of $\mathbb{R}^d$ 
given by $e_r = (\delta_{r,1}, \dots, \delta_{r,d})^{\top}$ ($r = 1, \dots, d$). 
Then, in the $k$-th tensor power 
$(\mathbb{R}^d)^{\otimes k} = \underbrace{\mathbb{R}^d \otimes \dots \otimes \mathbb{R}^d}_{k \text{ times}}$ 
(of dimension $d^k$), 
the set 

\vspace{-1mm}
$$
\{e_{j_1} \otimes \dots \otimes e_{j_k}\}_{j_1, \dots, j_k = 1}^d
$$ 

\vspace{3mm}
\noindent
forms an orthonormal basis with respect to the canonical Euclidean structure, satisfying

\vspace{-1mm}
$$
\bigl\langle e_{j_1} \otimes \dots \otimes e_{j_k},\, e_{i_1} \otimes \dots \otimes e_{i_k} \bigr\rangle 
= \prod_{m=1}^k \delta_{j_m, i_m}.
$$

\vspace{2mm}

Substituting $N = k$, $x = {\bf w}^p$, and $\tau = 1$ (noting that 
$S_k({\bf w}^p)_{0,1} := \left(S_k({\bf w}^p)_0\right)^{-1} \otimes S_k({\bf w}^p)_1 = S_k({\bf w}^p)_1$ 
since $S_k({\bf w}^p)_0 = 1$), the lift $S_k({\bf w}^p)_{0,1}$ in the direct sum 
$$
T^k(\mathbb{R}^d) = \bigoplus\limits_{m=0}^k (\mathbb{R}^d)^{\otimes m}
$$ 

\vspace{2mm}
\noindent
is written componentwise as

\vspace{-1mm}
$$
\begin{aligned}
S_k({\bf w}^p)_{0,1} = 
\Biggl( 1,\; \int\limits_{0 < s_1 < 1} d{\bf w}^p_{s_1},\; 
&\int\limits_{0 < s_1 < s_2 < 1} d{\bf w}^p_{s_1} \otimes d{\bf w}^p_{s_2},\; \dots, \\
&\dots,\; \int\limits_{0 < s_1 < \dots < s_k < 1} d{\bf w}^p_{s_1} 
\otimes \dots \otimes d{\bf w}^p_{s_k} \Biggr).
\end{aligned}
$$

\vspace{3mm}

Applying the canonical projection $\pi_k : T^k(\mathbb{R}^d) \to (\mathbb{R}^d)^{\otimes k}$ 
onto the $k$-th tensor level, defined explicitly 
for any tuple 
$$
a = (a_0, a_1, \dots, a_k) \in \bigoplus\limits_{m=0}^k (\mathbb{R}^d)^{\otimes m}
$$ 

\vspace{3mm}
\noindent
by $\pi_k(a_0, a_1, \dots, a_k) = a_k$, we obtain the expansion

\vspace{-1mm}
$$
\pi_k\bigl(S_k({\bf w}^p)_{0,1}\bigr) 
= \int\limits_{0 < s_1 < \dots < s_k < 1} 
d{\bf w}^p_{s_1} \otimes \dots \otimes d{\bf w}^p_{s_k}
= 
$$

$$
=\sum_{j_1, \dots, j_k = 1}^{d} 
S^{(j_1, \dots, j_k)p} \, e_{j_1} \otimes \dots \otimes e_{j_k},
$$

\vspace{3mm}
\noindent
where $S^{(j_1, \dots, j_k)p}$ is defined by (\ref{ud11}) for $T=1,$ $t=0.$

Consequently, the partial sum $S^{(i_1, \dots, i_k)p}$ 
is the $(i_1, \dots, i_k)$-component of $\pi_k \bigl(S_k({\bf w}^p)_{0,1}\bigr)$:

\vspace{-1mm}
$$
S^{(i_1, \dots, i_k)p} 
= \int\limits_{0 < s_1 < \dots < s_k < 1} 
d{\bf w}^{(i_1)p}_{s_1} \cdots d{\bf w}^{(i_k)p}_{s_k} 
= 
$$

$$
=\bigl\langle \pi_k\bigl(S_k({\bf w}^p)_{0,1}\bigr),\, 
e_{i_1} \otimes \dots \otimes e_{i_k} \bigr\rangle.
$$

\vspace{3mm}

\textbf{Step 2. Convergence at level~$2$ in $p'$-variation.} 
Fix a rough path exponent $p' \in (2, 3)$ (in the sense of Lyons, 1998, 
Sect.~1.2 and~2.2 \cite{Lyons}), so that $\lfloor p' \rfloor = 2$. 
Recall that for any $N \ge 1$ and geometric $p'$-rough paths 
$\mathbf{x}, \mathbf{y} \in C_0^{0,p'\text{-var}}\bigl(\lbrack 0, 1 \rbrack, G^N(\mathbb{R}^d)\bigr)$ 
\cite{Friz} (Friz and Victoir, 2010, Sect.~1.1, p.~372), the $p'$-variation 
distance is defined by

\vspace{-1mm}
$$
d_{p'\text{-var};\lbrack 0, 1 \rbrack}(\mathbf{x}, \mathbf{y}) 
= \sup_{D=\{0=t_0<t_1<\dots<t_r=1\}} 
\left( \sum_{l=0}^{r-1} d(\mathbf{x}_{t_l, t_{l+1}}, \mathbf{y}_{t_l, t_{l+1}})^{p'} \right)^{1/p'},
$$

\vspace{3mm}
\noindent
where $d(g, h) = \|g^{-1} \otimes h\|$ denotes the Carnot--Carath\'eodory 
metric on $G^N(\mathbb{R}^d)$ and $\mathbf{x}_{s,\tau} := \mathbf{x}_s^{-1} \otimes \mathbf{x}_\tau$.

For the $d$-dimensional Wiener process $\mathbf{w} = (\mathbf{w}^{(1)}, \dots, \mathbf{w}^{(d)})$, 
the covariance of each coordinate $R(s, \tau) = \min(s, \tau)$ has finite 2D $1$-variation 
($\rho = 1$, so $p' > 2\rho = 2$; see Friz and Victoir, 2010, Sect.~3.1.1 \cite{Friz}). 
By Theorem~35(i), Definition~36, and Theorem~51 (eq.~(19)) in Friz and Victoir (2010) \cite{Friz} 
applied to $X = \mathbf{w}$ and its projection approximations $X^{A_p} = \mathbf{w}^p$ 
($A_p = \{0, \dots, p\}$), projecting onto level~$2$ via the canonical contraction 
$\pi_{0,2} : G^3(\mathbb{R}^d) \to G^2(\mathbb{R}^d)$ yields

\vspace{-1mm}
$$
\lim_{p \to \infty} \mathsf{M}\Bigl\{ d_{p'\text{-var};\lbrack 0, 1 \rbrack}\bigl(S_2(\mathbf{w}^p), \mathbf{W}^{(2)}\bigr)^q \Bigr\} = 0 
\qquad \forall q \in \lbrack 1, \infty),
$$

\vspace{3mm}
\noindent
where $\mathbf{W}^{(2)} := \pi_{0,2}(\mathbf{W}) 
\in C_0^{0,p'\text{-var}}\bigl(\lbrack 0, 1 \rbrack, G^2(\mathbb{R}^d)\bigr)$ 
and $\pi_{0,2}\bigl(S_3(\mathbf{w}^p)\bigr) = S_2(\mathbf{w}^p)$,
$\pi_{0,2}(1, g_1, g_2, g_3)= (1, g_1, g_2)$.

For any subinterval $\lbrack s, \tau \rbrack \subset \lbrack 0, 1 \rbrack$, choosing 
a partition $D$ containing $s$ and $\tau$ (e.g., $D = \{0, s, \tau, 1\}$) gives

\vspace{-1mm}
$$
d\bigl(S_2(\mathbf{w}^p)_{s,\tau}, \mathbf{W}^{(2)}_{s,\tau}\bigr) 
\le d_{p'\text{-var};\lbrack 0, 1 \rbrack}\bigl(S_2(\mathbf{w}^p), \mathbf{W}^{(2)}\bigr) 
\xrightarrow[p \to \infty]{L^q(\Omega)} 0 
\qquad \forall q \in \lbrack 1, \infty).
$$

\vspace{3mm}

By the equivalence of homogeneous norms on $G^2(\mathbb{R}^d)$ \cite{Friz} (Friz and Victoir, 
2010, Sect.~1.1), 
it follows that for every $\lbrack s, \tau \rbrack \subset \lbrack 0, 1 \rbrack$, 
the second-level coordinates of $S_2(\mathbf{w}^p)_{s,\tau}$ converge in $L^2(\Omega)$ 
to those of $\mathbf{W}^{(2)}_{s,\tau}$. 

On the other hand, for any complete orthonormal 
system $\{\phi_j(x)\}_{j=0}^\infty$ in $L^2(\lbrack 0, 1 \rbrack)$ and every 
subinterval $[s, \tau]\subset [0, 1]$, the second-order 
iterated Riemann--Stieltjes integrals of the series ap\-pro\-xi\-ma\-ti\-ons $\mathbf{w}^p$ 
(i.e. $S^{(i_1 i_2)p}_{\tau,s}$)
over $[s, \tau]$ $(0\le s<\tau\le 1)$ converge in $L^2(\Omega)$ to the $2$-fold iterated 
Stra\-to\-no\-vich stochastic integrals over $[s, \tau].$ 

Indeed, from Theorem~4 for 
$\psi_1(t_1)={\bf 1}_{\{s<t_1\}},$ $\psi_2(t_2)={\bf 1}_{\{s<t_2<\tau\}},$ $p_1=p_2=p,$
we obtain 

\vspace{-1mm}
\begin{equation}
\label{zor1q}
\int\limits_s^{\tau}
\int\limits_s^{t_2}
d{\bf w}_{t_1}^{(i_1)} d{\bf w}_{t_2}^{(i_2)}
=\hbox{\vtop{\offinterlineskip\halign{
\hfil#\hfil\cr
{\rm l.i.m.}\cr
$\stackrel{}{{}_{p\to \infty}}$\cr
}} }\sum_{j_1, j_2=0}^{p}
C_{j_2j_1}(\tau,s)\Biggl(\zeta_{j_1}^{(i_1)}\zeta_{j_2}^{(i_2)}
-{\bf 1}_{\{i_1=i_2\ne 0\}}
{\bf 1}_{\{j_1=j_2\}}\Biggr),
\end{equation}

\vspace{3mm}
\noindent
where $i_1,i_2=1,\ldots,d,$

\vspace{-1mm}
$$
C_{j_2j_1}(\tau,s)=
\int\limits_s^{\tau}\phi_{j_{2}}(t_{2})
\int\limits_s^{t_2}\phi_{j_{1}}(t_{1})dt_1 dt_2,
$$

$$
\zeta_{j}^{(i)}=
\int\limits_t^T \phi_{j}(s) d{\bf w}_s^{(i)} \ \ \ (i=1,\ldots,d,\ j=0,1,2,\ldots)
$$ 

\vspace{3mm}
\noindent
are independent
standard Gaussian random variables for various 
$i$ or $j$.

Moreover, using (\ref{novorigin5}) for 
$\psi_1(t_1)={\bf 1}_{\{s<t_1\}}$ and $\psi_2(t_2)={\bf 1}_{\{s<t_2<\tau\}},$ we get w.~p.~1

\vspace{-1mm}
$$
\int\limits_s^{\tau}
\int\limits_s^{t_2}
\circ\ d{\bf w}_{t_1}^{(i_1)}\circ d{\bf w}_{t_2}^{(i_2)}
=\int\limits_s^{\tau}
\int\limits_s^{t_2}
d{\bf w}_{t_1}^{(i_1)} d{\bf w}_{t_2}^{(i_2)}+
$$
\begin{equation}
\label{zor2q}
+
\frac{1}{2}{\bf 1}_{\{i_1=i_2\}}
\int\limits_s^{\tau}dt_1,
\end{equation}

\vspace{3mm}
\noindent
where 
$i_1, i_2=1,\ldots,d.$

Note that (see \cite{2018a}, Sect.~2.1.4), that the equality (\ref{novorigin20})
is valid for the case 
of an arbitrary complete orthonormal 
system of functions in $L_2([t, T])$ and
$\psi_1(\tau), \psi_2(\tau)\in L_2([t, T]).$
Furthermore, applying the above version of 
(\ref{novorigin20})
for
$\psi_1(t_1)={\bf 1}_{\{s<t_1\}},$ $\psi_2(t_2)={\bf 1}_{\{s<t_2<\tau\}},$ we have

\vspace{-1mm}
\begin{equation}
\label{zor3q}
\sum_{j=0}^{\infty}\int\limits_s^{\tau}
\phi_j(t_2)
\int\limits_s^{t_2} 
\phi_j(t_1)dt_1 dt_2=
\frac{1}{2}\int\limits_s^{\tau}d\tau.
\end{equation}

\vspace{3mm}

Combining (\ref{zor1q}), (\ref{zor2q}), and (\ref{zor3q}), we obtain

\vspace{-1mm}
$$
\int\limits_s^{\tau}
\int\limits_s^{t_2}
\circ\ d{\bf w}_{t_1}^{(i_1)}\circ d{\bf w}_{t_2}^{(i_2)}=
\hbox{\vtop{\offinterlineskip\halign{
\hfil#\hfil\cr
{\rm l.i.m.}\cr
$\stackrel{}{{}_{p\to \infty}}$\cr
}} }\sum_{j_1, j_2=0}^{p}
C_{j_2j_1}(\tau,s)\zeta_{j_1}^{(i_1)}\zeta_{j_2}^{(i_2)},
$$

\vspace{3mm}
\noindent
i.e. 
the second-order 
iterated Riemann--Stieltjes integrals 

\vspace{-1mm}
$$
S^{(i_1 i_2)p}_{\tau,s}=\sum_{j_1, j_2=0}^{p}
C_{j_2j_1}(\tau,s)\zeta_{j_1}^{(i_1)}\zeta_{j_2}^{(i_2)}
$$

\vspace{3mm}
\noindent
of the series approximations $\mathbf{w}^p$ 
over $[s, \tau]$ $(0\le s<\tau\le 1)$ converge in $L^2(\Omega)$ to the $2$-fold iterated 
Stratonovich stochastic integrals over $[s, \tau].$ 

By uniqueness of the 
$L^2(\Omega)$-limit and path continuity, $\mathbf{W}^{(2)}$ coincides w.~p.~$1$ with 
the canonical Stratonovich Brownian rough path of step~$2$ on $\lbrack 0, 1 \rbrack$.

\vspace{2mm}

\textbf{Step 3. Extension to arbitrary multiplicity $k$.} 
By Step~2 and Markov's inequality, 

\vspace{-1mm}
$$
d_{p'\text{-var};\lbrack 0, 1 \rbrack}\bigl(S_2(\mathbf{w}^p), 
\mathbf{W}^{(2)}\bigr) \xrightarrow{\mathsf{P}} 0
$$ 

\vspace{3mm}
\noindent
as $p \to \infty$. Assume $k \ge 3$ (the case $k \le 2$ is already covered by Step~2). 
Since $p' \in (2, 3)$, we have $\lfloor p' \rfloor = 2$.

By Lyons' Extension Theorem \cite{Lyons} (Lyons, 1998, Theorem~2.2.1), 
the extension map

\vspace{-1mm}
$$
\tilde{S}_k : C_0^{0,p'\text{-var}}\bigl(\lbrack 0, 1 \rbrack, G^2(\mathbb{R}^d)\bigr) 
\to C_0^{0,p'\text{-var}}\bigl(\lbrack 0, 1 \rbrack, G^k(\mathbb{R}^d)\bigr)
$$

\vspace{3mm}
\noindent
is unique and continuous with respect to $d_{p'\text{-var};\lbrack 0, 1 \rbrack}$, 
and satisfies $\tilde{S}_k\bigl(S_2(\mathbf{w}^p)\bigr) = S_k(\mathbf{w}^p)$. 
Since the family of $k$-fold iterated Stratonovich stochastic integrals of $\mathbf{w}$
up to order $k$ has the same first two levels as $\mathbf{W}^{(2)}$, satisfies
Chen's identity, obeys the classical shuffle product rule (hence takes values in
$G^k(\mathbb{R}^d)$), and has finite $p'$-variation w.~p.~1, it is uniquely
determined on all higher levels $3, \dots, k$ by its second-level rough path
$\mathbf{W}^{(2)}$ (via Lyons' Extension Theorem
\cite{Lyons} (Lyons, 1998, Theorem~2.2.1)). Consequently,
$\pi_k(\tilde{S}_k(\mathbf{W}^{(2)})_{0,1}) \in (\mathbb{R}^d)^{\otimes k}$
coincides w.~p.~1 with the tensor $J^*$ of $k$-fold iterated Stratonovich
stochastic integrals of $\mathbf{w}$ over $[0,1]$:

\vspace{-1mm}
$$
\pi_k\bigl(\tilde{S}_k(\mathbf{W}^{(2)})_{0,1}\bigr) 
= J^* 
= \sum_{j_1,\dots,j_k=1}^d J^{*(j_1,\dots,j_k)} e_{j_1} \otimes \dots \otimes e_{j_k},
$$

\vspace{3mm}
\noindent
so that $J^{*(i_1,\dots,i_k)} = \bigl\langle \pi_k\bigl(\tilde{S}_k(\mathbf{W}^{(2)})_{0,1}\bigr), 
e_{i_1} \otimes \dots \otimes e_{i_k} \bigr\rangle$, 
where $J^*[\psi^{(k)}]_{T,t}^{(i_1 \dots i_k)}$ for $\psi_1 \equiv \dots \equiv \psi_k \equiv 1$ 
and $T = 1$, $t = 0$ is denoted by $J^{*(i_1,\dots,i_k)}$.

Applying the continuity of the extension map $\tilde{S}_k$ with respect to 
$d_{p'\text{-var};\lbrack 0, 1 \rbrack}$ (via the Continuous Mapping Theorem \cite{Mann}), we obtain

\vspace{-1mm}
$$
d_{p'\text{-var};\lbrack 0, 1 \rbrack}\bigl(S_k(\mathbf{w}^p), \tilde{S}_k(\mathbf{W}^{(2)})\bigr) 
= d_{p'\text{-var};\lbrack 0, 1 \rbrack}
\Bigl(\tilde{S}_k\bigl(S_2(\mathbf{w}^p)\bigr), \tilde{S}_k\bigl(\mathbf{W}^{(2)}\bigr)\Bigr) 
\xrightarrow[p \to \infty]{\mathsf{P}} 0.
$$

\vspace{3mm}

To pass to individual tensor coordinates, recall from Sect.~1.1 of Friz and Victoir 
(2010) \cite{Friz} the equivalence of homogeneous norms on $G^k(\mathbb{R}^d)$ 
(where $\|g\| = d(1, g)$ is the Carnot--Carath\'eodory norm entering $d(g, h) = \|g^{-1} \otimes h\|$):

\vspace{-1mm}
$$
\frac{1}{K_k} \max_{1 \le m \le k} \lvert\pi_m(g)\rvert^{1/m} 
\le \|g\| 
\le K_k \max_{1 \le m \le k} \lvert\pi_m(g)\rvert^{1/m}
$$

\vspace{3mm}
\noindent
for some constant $K_k \ge 1$ and all $g \in G^k(\mathbb{R}^d)$. Choosing the coarse 
partition $D = \{0, 1\}$ in the definition of $d_{p'\text{-var};\lbrack 0, 1 \rbrack}$ for 
$\mathbf{x} = S_k(\mathbf{w}^p)$ and $\mathbf{y} = \tilde{S}_k(\mathbf{W}^{(2)})$ 
(since paths in $C_0^{0,p'\text{-var}}\bigl(\lbrack 0, 1 \rbrack, G^k(\mathbb{R}^d)\bigr)$ satisfy 
the initial condition $\mathbf{x}_0 = \mathbf{y}_0 = 1$, we have 
$\mathbf{x}_{0,1} = \mathbf{x}_0^{-1} \otimes \mathbf{x}_1 = 
\mathbf{x}_1$ and $\mathbf{y}_{0,1} = \mathbf{y}_1$), 
we obtain

\vspace{-1mm}
$$
d(\mathbf{x}_{0,1}, \mathbf{y}_{0,1}) = 
\|\mathbf{x}_{0,1}^{-1} \otimes \mathbf{y}_{0,1}\| 
\le d_{p'\text{-var};\lbrack 0, 1 \rbrack}(\mathbf{x}, \mathbf{y}).
$$

\vspace{3mm}

Using the algebraic identity $\mathbf{y}_{0,1} - \mathbf{x}_{0,1} = 
\mathbf{x}_{0,1} \otimes (\mathbf{x}_{0,1}^{-1} \otimes \mathbf{y}_{0,1} - 1)$ 
and the multiplication rule

\vspace{-1mm}
$$
\pi_k(\mathbf{a} \otimes \mathbf{b}) = \sum_{m=0}^k \pi_{k-m}(\mathbf{a}) \otimes \pi_m(\mathbf{b})
$$

\vspace{3mm}
\noindent
in $T^k(\mathbb{R}^d)$ \cite{Lyons} (Lyons, 1998, eq.~(2.2)), we obtain

\vspace{-1mm}
$$
\begin{aligned}
\pi_k(\mathbf{y}_{0,1}) - \pi_k(\mathbf{x}_{0,1}) 
&= \pi_k\Bigl(\mathbf{x}_{0,1} \otimes (\mathbf{x}_{0,1}^{-1} \otimes \mathbf{y}_{0,1} - 1)\Bigr) \\
&= \sum_{m=0}^k \pi_{k-m}(\mathbf{x}_{0,1}) \otimes \pi_m(\mathbf{x}_{0,1}^{-1} \otimes \mathbf{y}_{0,1} - 1) \\
&= \sum_{m=1}^k \pi_{k-m}(\mathbf{x}_{0,1}) \otimes \pi_m(\mathbf{x}_{0,1}^{-1} \otimes \mathbf{y}_{0,1}),
\end{aligned}
$$

\vspace{3mm}
\noindent
since $\pi_0(\mathbf{x}_{0,1}^{-1} \otimes \mathbf{y}_{0,1} - 1) = 1 - 1 = 0$ and 
$\pi_m(g - 1) = \pi_m(g)$ for all $m \ge 1$.

Applying $\lvert\pi_m(g)\rvert \le K_k^m \|g\|^m$ (which follows from the lower bound in 
the norm equivalence above) to $g = \mathbf{x}_{0,1}^{-1} \otimes \mathbf{y}_{0,1}$ 
(so that $\|g\| = d(\mathbf{x}_{0,1}, \mathbf{y}_{0,1})$) and using the triangle inequality 
together with the submultiplicativity 
$\lvert u \otimes v \rvert \le \lvert u \rvert\,\lvert v \rvert$ of the tensor norm, 
we deduce the explicit estimate

\vspace{-1mm}
$$
\bigl\lvert\pi_k(\mathbf{y}_{0,1}) - \pi_k(\mathbf{x}_{0,1})\bigr\rvert 
\le \sum_{m=1}^k K_k^m \lvert\pi_{k-m}(\mathbf{x}_{0,1})\rvert\, d(\mathbf{x}_{0,1}, \mathbf{y}_{0,1})^m 
\qquad \forall k \ge 1.
$$

\vspace{3mm}

Since $d_{p'\text{-var};\lbrack 0, 1 \rbrack}(\mathbf{x}, \mathbf{y}) = d_{p'\text{-var};\lbrack 0, 1 \rbrack}\bigl(S_k(\mathbf{w}^p), \tilde{S}_k(\mathbf{W}^{(2)})\bigr) \xrightarrow{\mathsf{P}} 0$, 
it follows that $d(\mathbf{x}_{0,1}, \mathbf{y}_{0,1}) \xrightarrow{\mathsf{P}} 0$ 
(and hence $d(\mathbf{x}_{0,1}, \mathbf{y}_{0,1})^m \xrightarrow{\mathsf{P}} 0$ for every $m = 1, \dots, k$). 
Since the estimate above holds at every level (with $\pi_0(\mathbf{x}_{0,1}) = 1$), 
it follows by induction on the tensor level that each lower-level projection 
$\pi_{k-m}(\mathbf{x}_{0,1})$ ($m \in \{1, \dots, k\}$) converges in probability 
to $\pi_{k-m}(\mathbf{y}_{0,1})$ and hence $\lvert\pi_{k-m}(\mathbf{x}_{0,1})\rvert$ is bounded 
in probability as $p \to \infty$. Since the product of a sequence bounded in 
probability and a sequence converging to $0$ in probability converges to $0$ 
in probability, each term in the sum and thus the entire right-hand side 
converges to $0$ in probability, so that $\lvert\pi_k(\mathbf{y}_{0,1}) - \pi_k(\mathbf{x}_{0,1})\rvert \xrightarrow{\mathsf{P}} 0$. 
Substituting back $\mathbf{x} = S_k(\mathbf{w}^p)$ and $\mathbf{y} = \tilde{S}_k(\mathbf{W}^{(2)})$, 
and noting that the $(i_1, \dots, i_k)$-coordinate of a tensor cannot exceed its 
Euclidean norm, we obtain for every multi-index $(i_1, \dots, i_k)$:

\vspace{-1mm}
$$
\bigl\lvert S^{(i_1,\dots,i_k)p} - J^{*(i_1,\dots,i_k)} \bigr\rvert 
\le \bigl\lvert\pi_k(\mathbf{y}_{0,1}) - \pi_k(\mathbf{x}_{0,1})\bigr\rvert 
= 
$$

$$
=\bigl\lvert\pi_k\bigl(\tilde{S}_k(\mathbf{W}^{(2)})_{0,1}\bigr) - 
\pi_k\bigl(S_k(\mathbf{w}^p)_{0,1}\bigr)\bigr\rvert 
\xrightarrow[p \to \infty]{\mathsf{P}} 0.
$$

\vspace{3mm}

\noindent\textbf{Step 4. Upgrading to mean-square convergence.}
By Proposition~21 in Friz and Victoir (2010) \cite{Friz}, 
the iterated integral $S^{(i_1, \dots, i_k)p}$ 
of the Gaussian process ${\bf w}^p$ (whose sample paths are absolutely continuous
and hence of finite $1$-variation w.~p.~1)
lies in the $k$-th (inhomogeneous) Wiener chaos $\mathcal{H}^{\le k}$. 
Since $\mathcal{H}^{\le k}$ is closed under convergence in probability, 
we also have $J^{*(i_1, \dots, i_k)} \in \mathcal{H}^{\le k}$, 
and thus the error 

\vspace{-1mm}
$$
Z_p := S^{(i_1, \dots, i_k)p} - J^{*(i_1, \dots, i_k)} 
\in \mathcal{H}^{\le k}.
$$

\vspace{3mm}

Recall Lemma~22 in Friz and Victoir (2010) \cite{Friz}
(hypercontractivity on the inhomogeneous Wiener chaos $\mathcal{H}^{\le k}$): 
for any $Z \in \mathcal{H}^{\le k}$ and $1 < p_0 < q < \infty$,

\vspace{-1mm}
$$ |Z|_{L^{p_0}(\Omega)} 
\le |Z|_{L^q(\Omega)} 
\le |Z|_{L^{p_0}(\Omega)} (k+1)(q-1)^{k/2} 
\max\bigl(1, (p_0 - 1)^{-k}\bigr).
$$

\vspace{3mm}

Furthermore, by Remark~24 in Friz and Victoir (2010) \cite{Friz}, 
on $\mathcal{H}^{\le k}$ convergence in probability 
is equivalent to convergence in $L^q(\Omega)$ 
for every $q \in (0, \infty)$:

\vspace{-1mm}
$$
Z_p \xrightarrow[p \to \infty]{{\sf P}} 0 
\quad \Longleftrightarrow \quad 
\lim_{p \to \infty} |Z_p|_{L^q(\Omega)} = 0 
\qquad \forall\, q \in (0, \infty).
$$

\vspace{3mm}

Substituting $Z_p = S^{(i_1, \dots, i_k)p} - J^{*(i_1, \dots, i_k)}$ 
and $q = 2$ into this equivalence, 
and using Step~3 ($Z_p \xrightarrow{{\sf P}} 0$), 
we finally obtain

\vspace{-1mm}
$$
\lim_{p \to \infty} |Z_p|_{L^2(\Omega)}^2 
= \lim_{p \to \infty} 
{\sf M}\left\{ 
\left| S^{(i_1, \dots, i_k)p} 
- J^{*(i_1, \dots, i_k)} \right|^2
\right\} = 0.
$$

\vspace{3mm}

Theorem~64 is proved. $\Box.$

\vspace{2mm}

{\bf Theorem~65\ (Extension to weighted iterated Stratonovich stochastic integrals)}\
\cite{2018a}, \cite{19999999}.\
{\it Let $k \ge 1,$ 
let $\{\phi_j(x)\}_{j=0}^{\infty}$ be an arbitrary 
complete orthonormal system of functions in $L^2([0, 1]),$
$(i_1, \dots, i_k) \in \{1, \dots, d\}^k,$
and let $\psi_1(\tau), \dots,$ $\psi_k(\tau) \in C^1([0,\,1])$ be deterministic 
continuously dif\-fe\-ren\-ti\-able weight functions. 
Define the weighted iterated stochastic integrals {\rm (}Riemann--Stieltjes and Stra\-to\-no\-vich$,$ 
res\-pec\-tive\-ly{\rm )}

\vspace{-1mm}
$$
S[\psi^{(k)}]^{(i_1, \dots, i_k)p}_{1,0}
:= \int\limits_{0 < t_1 < \dots < t_k < 1} 
\psi_1(t_1) \cdots \psi_k(t_k)\, 
d{\bf w}_{t_1}^{(i_1)p} \cdots d{\bf w}_{t_k}^{(i_k)p}, 
$$

$$
J^{*}[\psi^{(k)}]_{1,0}^{(i_1 \ldots i_k)}
:= \int\limits_{0 < t_1 < \dots < t_k < 1} 
\psi_1(t_1) \cdots \psi_k(t_k) 
\circ d{\bf w}_{t_1}^{(i_1)} \cdots \circ d{\bf w}_{t_k}^{(i_k)},
$$

\vspace{3mm}
\noindent
where 
$$
S[\psi^{(k)}]_{1,0}^{(i_1, \dots, i_k)p}
=\sum_{j_1,\ldots,j_k=0}^{p}
C_{j_k\ldots j_1}
\prod_{l=1}^k \zeta_{j_l}^{(i_l)};
$$

\vspace{3mm}
\noindent
another notations are the same as in Theorem~{\rm 64}. Then

\vspace{-1mm}
$$
\lim_{p \to \infty} 
{\sf M}\left\{
\left|S[\psi^{(k)}]^{(i_1, \dots, i_k)p}_{1,0}
- J^{*}[\psi^{(k)}]_{1,0}^{(i_1 \ldots i_k)} \right|^2 
\right\} = 0.
$$
}

\vspace{3mm}

{\bf Proof.} Recall that $\mathbf{w}^p = (\mathbf{w}^{(1)p}, \dots, \mathbf{w}^{(d)p})$ 
denotes the $p$-th partial sum expansion of the standard $d$-dimensional Wiener process 
$\mathbf{w}$ with independent components introduced in Theorem~64
(see (\ref{um1xx}) for $t=0$). 
For each $i_m\in\{1,\ldots,d\},$ $m \in \{1, \dots, k\},$ and $t \in \lbrack 0, 1 \rbrack$, 
define the weighted process

\vspace{-1mm}
$$
Y_t^{(i_m,m)} := \int\limits_0^t \psi_m(s)\, d\mathbf{w}_s^{(i_m)}
$$

\vspace{3mm}
\noindent
and its absolutely continuous approximation

\vspace{-1mm}
$$
Y_{t}^{(i_m,m)p} := \int\limits_0^t \psi_m(s)\, d\mathbf{w}_s^{(i_m)p},
$$

\vspace{3mm}
\noindent
where $d{\bf w}_s^{(i)p}$ $(i=1,\ldots,d)$ is defined by (\ref{um1xxx}).

Applying Ito's formula to $F(t, {\bf w}_t^{(i_m)}) = \psi_m(t){\bf w}_t^{(i_m)}$ 
(where $F(t,x)=\psi_m(t)x,$ $\partial^2 F/\partial x^2 \equiv 0$) 
with $\mathbf{w}_0^{(i_m)} = 0$, we obtain w.~p.~1

\vspace{-1mm}
$$
Y_t^{(i_m,m)} = \int\limits_0^t \psi_m(s)\, d\mathbf{w}_s^{(i_m)} 
= \psi_m(t)\mathbf{w}_t^{(i_m)} - \int\limits_0^t \psi_m'(s)\mathbf{w}_s^{(i_m)}\, ds,
$$

\vspace{3mm}
\noindent
and similarly, by classical integration by parts for the absolutely continuous path
$\mathbf{w}_s^{(i_m)p}$ (of finite $1$-variation w.~p.~1),

\vspace{-1mm}
$$
Y_{t}^{(i_m,m)p} = \psi_m(t)\mathbf{w}_t^{(i_m)p} - \int\limits_0^t \psi_m'(s)\mathbf{w}_s^{(i_m)p}\, ds\ \ \ 
\hbox{w.~p.~1}.
$$

\vspace{3mm}

Since $\psi_m(s)$ is deterministic and of bounded variation, the Ito and Stratonovich 
differentials coincide 
($\psi_m(t)\,d\mathbf{w}_t^{(i_m)} = \psi_m(t) \circ d\mathbf{w}_t^{(i_m)} = \circ\, dY_t^{(i_m,m)}$), 
and $dY_{t}^{(i_m,m)p} = \psi_m(t)d{\bf w}_t^{(i_m)p}$. 
Hence the weighted iterated 
integrals rewrite as unweighted iterated integrals of the $k$-dimensional processes 
$Y^p_t = (Y^{(i_1,1)p}_t,$ $\ldots, Y^{(i_k,k)p}_t)$ and
$Y_t = (Y^{(i_1,1)}_t, \ldots, Y^{(i_k,k)}_t),$ $i_1,\ldots,i_k=1,\ldots,d$:

\vspace{-1mm}
$$
S[\psi^{(k)}]^{(i_1, \dots, i_k)p}_{1,0}
= \int\limits_{0 < t_1 < \dots < t_k < 1} dY_{t_1}^{(i_1,1)p} \cdots dY_{t_k}^{(i_k,k)p}, 
$$

$$
J^{*}[\psi^{(k)}]_{1,0}^{(i_1 \ldots i_k)}
= \int\limits_{0 < t_1 < \dots < t_k < 1} \circ\, dY_{t_1}^{(i_1,1)} \cdots \circ dY_{t_k}^{(i_k,k)}.
$$

\vspace{3mm}

Note that the components $Y^{(i_m,m)}_t$ and $Y^{(i_l,l)}_t$ of
$Y_t = (Y^{(i_1,1)}_t,\ldots,Y^{(i_k,k)}_t)$
are independent when\-ever their first upper indices differ
($i_m \ne i_l$). Even though components with $i_m = i_l$
($m \ne l$) are correlated, the full covariance matrix
$R_Y$ has finite 2D $1$-variation ($\rho = 1$) on $[0,1]^2$ (see definition of 2D $1$-variation
in \cite{Friz}),
as shown below.

Now observe that $Y_t = (Y^{(i_1,1)}_t,\ldots,Y^{(i_k,k)}_t)$ is a centered $k$-dimensional 
Gaussian process whose covariance function $R_Y$ has finite 2D $1$-variation 
($\rho = 1$) on $[0, 1]^2$: indeed, by Ito's isometry, for any 
rectangle $[s_1, s_2]\times [t_1, t_2] \subset [0, 1]^2$, 
the rectangular covariance increment of $(Y^{(i_m,m)}_t, Y^{(i_l,l)}_t)$ vanishes if $i_m \ne i_l$ 
and equals 

\vspace{-1mm}
$$
\int\limits_{[s_1, s_2]\cap [t_1, t_2]} \psi_m(\tau)\psi_l(\tau) d\tau
$$ 

\vspace{3mm}
\noindent
if $i_m = i_l$, so summing over any grid partition of $[0, 1]^2$ gives 
a total variation bounded by 

\vspace{-1mm}
$$
\max\limits_{1 \le m \le k} \|\psi_m\|_\infty^2 < \infty.
$$

\vspace{3mm}

Next, let $\mathcal{F}_p := \sigma(\zeta_j^{(i)} : 0 \le j \le p,\, 1 \le i \le d)$ 
be the filtration generated by the first $p+1$ expansion coefficients of $\mathbf{w}_t$. 
Because $\{\phi_j(x)\}_{j=0}^\infty$ is a complete orthonormal basis of $L^2([0, 1])$, 
the coefficients $\zeta_j^{(i)}$ are independent standard Gaussian random variables 
and the $\sigma$-algebra 

\vspace{-1mm}
$$
\mathcal{F}_\infty := \sigma\left(\bigcup\limits_{p=0}^\infty \mathcal{F}_p\right)
$$ 

\vspace{3mm}
\noindent
coincides (up to ${\sf P}$-null sets)
with $\sigma(\mathbf{w}_t : t \in \lbrack 0, 1 \rbrack)$, so every $Y_t^{(i_m,m)}$ 
is $\mathcal{F}_\infty$-measurable. Moreover, since the tail series 

\vspace{-1mm}
$$
\sum_{j=p+1}^\infty \zeta_j^{(i)} \int\limits_0^t \phi_j(s)ds
$$ 

\vspace{3mm}
\noindent
is independent of $\mathcal{F}_p$ 
and has zero mean, we have w.~p.~1

\vspace{-1mm}
$$
\mathbf{w}_t^{(i)p} = {\sf M}\left\{\mathbf{w}_t^{(i)} \mid \mathcal{F}_p\right\}
$$ 

\vspace{3mm}
\noindent
for all $t \in \lbrack 0, 1 \rbrack$. By linearity of conditional expectation in the 
representations above, it follows that w.~p.~1

\vspace{-1mm}
$$
Y_{t}^{(i_m,m)p} = {\sf M}\left\{Y_t^{(i_m,m)} \mid \mathcal{F}_p \right\}
$$

\vspace{3mm}
\noindent
for all $m \in \{1, \dots, k\},$ $i_m\in\{1,\ldots,d\},$
and $t \in [0, 1]$.

In particular, $(Y_{t}^{(i_m,m)p})_{p \ge 0}$ is a closed centered Gaussian martingale with 
absolutely continuous sample paths (of finite $1$-variation w.~p.~1), 
and by L\'evy's Martingale Convergence Theorem, 

\vspace{-1mm}
$$
Y_{t}^{(i_m,m)p} \to Y_t^{(i_m,m)}
$$ 

\vspace{3mm}
\noindent
w.~p.~1 and in $L^2(\Omega)$ as $p \to \infty$ for every 
$t \in \lbrack 0, 1 \rbrack$.

Fix a rough path exponent $p' \in (2, 3)$ (so that $p' > 2\rho = 2$ and $\lfloor p' \rfloor = 2$). 
Applying Theorem~35(i), Definition~36, and Theorem~51 in Friz and Victoir (2010) \cite{Friz}
to the conditional expectation ap\-pro\-xi\-ma\-ti\-ons 
$Y_{t}^p = {\sf M}\left\{Y_t \mid \mathcal{F}_p \right\}$ 
and projecting onto level~$2$ via the canonical 
contraction $\pi_{0,2} : G^3(\mathbb{R}^k) \to G^2(\mathbb{R}^k)$, 
we obtain

\vspace{-1mm}
$$
\lim\limits_{p \to \infty} \mathsf{M}\Bigl\{ d_{p'\text{-var};\lbrack 0, 1 
\rbrack}\bigl(S_2(Y^p), \mathbf{Y}^{(2)}\bigr)^q \Bigr\} = 0 
\qquad \forall q \in \lbrack 1, \infty),
$$

\vspace{3mm}
\noindent
where $\mathbf{Y}^{(2)} := \pi_{0,2}(\mathbf{Y}) 
\in C_0^{0,p'\text{-var}}\bigl(\lbrack 0, 1 \rbrack, G^2(\mathbb{R}^k)\bigr)$. 

In particular, for every subinterval $\lbrack s, \tau \rbrack \subset \lbrack 0, 1 \rbrack$, 
the second-level coordinates of $S_2(Y^p)_{s,\tau}$ converge in $L^2(\Omega)$ to 
those of $\mathbf{Y}^{(2)}_{s,\tau}$. 

On the other hand, since the $(m, l)$-coordinate 
of $S_2(Y^p)_{s,\tau}$ is the weighted second-order iterated Riemann--Stieltjes integral 

\vspace{-1mm}
$$
\int\limits_{s < t_1 < t_2 < \tau} 
\psi_l(t_2)\psi_m(t_1) d\mathbf{w}_{t_1}^{(i_m)p} d\mathbf{w}_{t_2}^{(i_l)p}
$$ 

\vspace{3mm}
\noindent
with smooth weights $\psi_m(\tau), \psi_l(\tau) \in C^1(\lbrack 0, 1 \rbrack),$ then
we have
$L^2(\Omega)$-con\-ver\-gence of 

\vspace{-1mm}
$$
\int\limits_{s < t_1 < t_2 < \tau} 
\psi_l(t_2)\psi_m(t_1)  d\mathbf{w}_{t_1}^{(i_m)p} d\mathbf{w}_{t_2}^{(i_l)p}
$$ 

\vspace{3mm}
\noindent
to the $2$-fold iterated Stratonovich stochastic integral 

\vspace{-1mm}
$$
\int\limits_{s < t_1 < t_2 < \tau} \psi_l(t_2)\psi_m(t_1) 
\circ d\mathbf{w}_{t_1}^{(i_m)} \circ d\mathbf{w}_{t_2}^{(i_l)} 
= \int\limits_{s < t_1 < t_2 < \tau} \circ\, dY_{t_1}^{(i_m,m)} \circ dY_{t_2}^{(i_l,l)}.
$$ 

\vspace{3mm}

Indeed, from Theorem~4 for 
$\tilde \psi_m(t_1)=\psi_m(t_1){\bf 1}_{\{s<t_1\}}$ and 
$\tilde \psi_l(t_2)=\psi_l(t_2){\bf 1}_{\{s<t_2<\tau\}}$ 
($\psi_m(\tau),\psi_l(\tau) \in C^1([0, 1]$) (instead of $\psi_1(t_1), \psi_2(t_2)$),
$p_1=p_2=p,$
we obtain 

\vspace{-1mm}
$$
\int\limits_s^{\tau}\psi_l(t_2)
\int\limits_s^{t_2}\psi_m(t_1)
d{\bf w}_{t_1}^{(i_m)} d{\bf w}_{t_2}^{(i_l)}
=
$$

\begin{equation}
\label{zor1qx}
=\hbox{\vtop{\offinterlineskip\halign{
\hfil#\hfil\cr
{\rm l.i.m.}\cr
$\stackrel{}{{}_{p\to \infty}}$\cr
}} }\sum_{j_1, j_2=0}^{p}
C_{j_2j_1}(\tau,s)\Biggl(\zeta_{j_1}^{(i_m)}\zeta_{j_2}^{(i_l)}
-{\bf 1}_{\{i_m=i_l\ne 0\}}
{\bf 1}_{\{j_1=j_2\}}\Biggr),
\end{equation}

\vspace{3mm}
\noindent
where $i_m,i_l=1,\ldots,d,$
$$
C_{j_2j_1}(\tau,s)=
\int\limits_s^{\tau}\psi_l(t_2)\phi_{j_{2}}(t_{2})
\int\limits_s^{t_2}\psi_m(t_1)\phi_{j_{1}}(t_{1})dt_1 dt_2,
$$

$$
\zeta_{j}^{(i)}=
\int\limits_t^T \phi_{j}(s) d{\bf w}_s^{(i)} \ \ \ (i=1,\ldots,d,\ j=0,1,2,\ldots)
$$ 

\vspace{3mm}
\noindent
are independent
standard Gaussian random variables for various 
$i$ or $j$.

Moreover, using (\ref{novorigin5}) for 
$\tilde \psi_m(t_1)=\psi_m(t_1){\bf 1}_{\{s<t_1\}},$ 
$\tilde \psi_l(t_2)=\psi_l(t_2){\bf 1}_{\{s<t_2<\tau\}}$ 
($\psi_m(\tau), \psi_l(\tau) \in C^1([0, 1]$) (instead of $\psi_1(t_1), \psi_2(t_2)$), we get w.~p.~1

\vspace{-1mm}
$$
\int\limits_s^{\tau}\psi_l(t_2)
\int\limits_s^{t_2}\psi_m(t_1)
\circ\ d{\bf w}_{t_1}^{(i_m)}\circ d{\bf w}_{t_2}^{(i_l)}
=\int\limits_s^{\tau}\psi_l(t_2)
\int\limits_s^{t_2}\psi_m(t_1)
d{\bf w}_{t_1}^{(i_m)} d{\bf w}_{t_2}^{(i_l)}+
$$

\begin{equation}
\label{zor2qx}
+
\frac{1}{2}{\bf 1}_{\{i_m=i_l\}}
\int\limits_s^{\tau} \psi_l(t_1)\psi_m(t_1)dt_1,
\end{equation}

\vspace{3mm}
\noindent
where 
$i_m, i_l=1,\ldots,d.$

Recall (see \cite{2018a}, Sect.~2.1.4) that the equality (\ref{novorigin20})
is valid for the case 
of an arbitrary complete orthonormal 
system of functions in $L_2([t, T])$ and
$\psi_1(\tau), \psi_2(\tau)\in L_2([t, T]).$
Furthermore, applying the above version of 
(\ref{novorigin20})
for
$\tilde \psi_m(t_1)=\psi_m(t_1){\bf 1}_{\{s<t_1\}}$ and
$\tilde \psi_l(t_2)=\psi_l(t_2){\bf 1}_{\{s<t_2<\tau\}}$ 
($\psi_m(\tau), \psi_l(\tau) \in C^1([0, 1]$) (instead of $\psi_1(t_1), \psi_2(t_2)$), we have

\vspace{-1mm}
\begin{equation}
\label{zor3qx}
\sum_{j=0}^{\infty}\int\limits_s^{\tau}
\psi_l(t_2)\phi_j(t_2)
\int\limits_s^{t_2} 
\psi_m(t_1)\phi_j(t_1)dt_1 dt_2=
\frac{1}{2}\int\limits_s^{\tau} \psi_l(t_1)\psi_m(t_1) dt_1.
\end{equation}

\vspace{3mm}

Combining (\ref{zor1qx}), (\ref{zor2qx}), and (\ref{zor3qx}), we obtain

\vspace{-1mm}
$$
\int\limits_s^{\tau}\psi_l(t_2)
\int\limits_s^{t_2}\psi_m(t_1)
\circ\ d{\bf w}_{t_1}^{(i_m)}\circ d{\bf w}_{t_2}^{(i_l)}=
\hbox{\vtop{\offinterlineskip\halign{
\hfil#\hfil\cr
{\rm l.i.m.}\cr
$\stackrel{}{{}_{p\to \infty}}$\cr
}} }\sum_{j_1, j_2=0}^{p}
C_{j_2j_1}(\tau,s)\zeta_{j_1}^{(i_m)}\zeta_{j_2}^{(i_l)},
$$

\vspace{3mm}
\noindent
where
$$
\int\limits_{s < t_1 < t_2 < \tau} \psi_l(t_2) \psi_m(t_1)
\circ d\mathbf{w}_{t_1}^{(i_m)} \circ d\mathbf{w}_{t_2}^{(i_l)} 
= \int\limits_{s < t_1 < t_2 < \tau} \circ\, dY_{t_1}^{(i_m,m)} \circ dY_{t_2}^{(i_l,l)},
$$ 

\vspace{3mm}
\noindent
i.e. 
the second-order 
iterated Riemann--Stieltjes integrals 

\vspace{-1mm}
$$
S^{(i_1 i_2)p}_{\tau,s}=\sum_{j_1, j_2=0}^{p}
C_{j_2j_1}(\tau,s)\zeta_{j_1}^{(i_m)}\zeta_{j_2}^{(i_l)}
$$

\vspace{3mm}
\noindent
of the series approximations $\mathbf{w}^p$ 
over $[s, \tau]$ $(0\le s<\tau\le 1)$ converge in $L^2(\Omega)$ to the $2$-fold iterated 
Stratonovich stochastic integrals over $[s, \tau].$

By uniqueness of the $L^2(\Omega)$-limit and path continuity, $\mathbf{Y}^{(2)}$ 
coincides w.~p.~$1$ with the canonical Stratonovich rough path of $Y$ of step~$2$ 
on $\lbrack 0, 1 \rbrack$.

Finally, repeating Steps~3 and~4 of the proof of Theorem~64 verbatim for the unweighted 
iterated integrals of $Y^p_t$ and $Y_t$ in $\mathbb{R}^k$ with multi-index $(1, 2, \dots, k)$ 
(using the uniqueness and continuity of the Lyons extension $\tilde{S}_k$ in Lyons, 
1998, Theorem~2.2.1~\cite{Lyons}, and hypercontractivity on the inhomogeneous 
Wiener chaos $\mathcal{H}^{\le k}$ 
in Friz and Victoir, 2010, Lemma~22 and Remark~24 \cite{Friz})
yields

\vspace{-1mm}
$$
\lim_{p \to \infty} 
{\sf M}\left\{
\left|S[\psi^{(k)}]^{(i_1, \dots, i_k)p}_{1,0}
- J^{*}[\psi^{(k)}]_{1,0}^{(i_1 \ldots i_k)} \right|^2 
\right\} = 0.
$$

\vspace{3mm}
\noindent
Theorem~65 is proved. $\Box$

\vspace{2mm}

\textbf{Remark~1.}\ {\it Since the approximation errors

\vspace{-1mm}
$$
Z_p := S_{T,t}^{(i_1,\dots,i_k)p} - J^*[\psi^{(k)}]_{T,t}^{(i_1\dots i_k)}
$$

\vspace{3mm}
\noindent
in Theorem~{\rm 64} {\rm (}for $\psi_1(\tau) \equiv \dots \equiv \psi_k(\tau) \equiv 1${\rm )} and

\vspace{-1mm}
$$
Z_p[\psi^{(k)}] := S[\psi^{(k)}]_{1,0}^{(i_1,\dots,i_k)p}
- J^*[\psi^{(k)}]_{1,0}^{(i_1\dots i_k)}
$$

\vspace{3mm}
\noindent
in Theorem~{\rm 65} belong to the $k$-th inhomogeneous Wiener chaos
$\mathcal{H}^{\le k}$ and converge to zero in probability as $p \to \infty,$
Remark~{\rm 24} in Friz and Victoir {\rm (2010) \cite{Friz}} immediately implies that both
Theorems~{\rm 64} and~{\rm 65} remain valid for convergence in $L^q(\Omega)$ for every
$q \in (0, \infty)$ {\rm (}rather than only for $q = 2${\rm )}. Namely$,$ under the
assumptions of Theorems~{\rm 64} and~{\rm 65,} respectively$,$ we have

\vspace{-1mm}
$$
\lim_{p\to\infty} \mathsf{M}\left\{\left|
S_{T,t}^{(i_1,\dots,i_k)p} - J^*[\psi^{(k)}]_{T,t}^{(i_1\dots i_k)}
\right|^q\right\} = 0 \qquad \forall q \in (0, \infty)
$$

\vspace{3mm}
\noindent
and
$$
\lim_{p\to\infty} \mathsf{M}\left\{\left|
S[\psi^{(k)}]_{1,0}^{(i_1,\dots,i_k)p} - J^*[\psi^{(k)}]_{1,0}^{(i_1\dots i_k)}
\right|^q\right\} = 0 \qquad \forall q \in (0, \infty).
$$
}

\vspace{3mm}

{\bf Remark~2.}\ {\it Theorems {\rm 64} and {\rm 65} prove Hypothesis {\rm 2} 
$($for the case $p_1=\ldots=p_k=p,$ $i_1,\ldots,i_k=1,\ldots,d,$
and $\psi_1(\tau),\ldots,\psi_k(\tau)\in C^1([t, T]))$
on the mean-square convergence of expansions of iterated 
Stratonovich stochastic integrals of arbitrary multiplicity 
$(k \in \mathbb{N})$ with respect to components of a multidimensional 
Wiener process.
}

\end{document}